\documentclass[10pt,letterpaper]{report}
\usepackage{standalone}
\usepackage{import}
\usepackage{A2-customcommands}
\usepackage{xr} 
\graphicspath{{images/}}

\begin{document} 
\pagenumbering{arabic} 
\pagestyle{plain}


\pdfbookmark{Title Page}{title}
	\begin{center}
	\null\vspace{0.7in}
	\begin{spacing}{1.2}
		\bfseries\LARGE
		An Exposition on the \\
		Algebra and Computation \\
		of Persistent Homology

		\large (second version)
	\end{spacing}
	\vspace{0.5in}
	\begin{spacing}{1.2}
		An expository paper submitted to 
		
		{\bfseries OREGON STATE UNIVERSITY}

		in partial fulfillment of the requirements for the degree of 
		
		{\bfseries MASTER OF SCIENCE (MS) DEGREE IN MATHEMATICS}
	\end{spacing}
	\vspace{0.5in}
	\begin{spacing}{1.2}
		Author
		
		{\bfseries\large Jason Ranoa}
		\vspace{0.5in}

		Advisor

		{\bfseries\large Dr.\ Christine Escher}
		\vspace{0.5in}

		Initially submitted
		
		{\bfseries February 2024}

		with corrections submitted August 2024.
	\end{spacing}
	\end{center}
\newpage

\clearpage

	\begin{blockquote}
		\noindent
		{\bfseries\large An Exposition on the
		Algebra and Computation \\[3pt]
		of Persistent Homology.} 

		\noindent
		\textsc{Abstract.}
		We discuss the algebra behind the matrix reduction algorithm for persistent homology, as presented in the paper \textit{``Computing Persistent Homology''} by Afra Zomorodian and Gunnar Carlsson, in the lens of the more modern characterization of persistence modules as functors from a poset category to a category of vector spaces over a field adopted by authors such as Peter Bubenik, Frederik Chazal, and Ulrich Bauer.

		\spacer 

		\noindent 
		\textsc{Author's Note for the Second Version.} 

		This copy of the paper is a corrected version of the expository paper submitted in February 2024 before my final oral examination (defense). 
		Because of time constraints and other personal and health issues, 
			a non-negligible number of errors were not fixed or were not caught in time.
		Listed below are the major differences and changes:
		\begin{enumerate}[label={(\roman*)}]
			\item 
			\textbf{\fref{section:interval-decomposition-persmod}}
			and 
			\textbf{\fref{section:cat-equiv-graded-modules}}
			contain significant changes and corrections in the exposition and commentary.
			The definitions and results presented remain mostly the same, e.g.\ with minor changes in notation.

			\item 
			The entirety of \textbf{\fref{section:graded-mod-notation}} 
			is re-written.
			When I was first writing this section, I only had a tenuous grasp of graded module theory, and this was reflected in the definitions and explanations I provided.
			I was also missing a number of important definitions (e.g.\ graded submodule) and results (e.g.\ quotients of graded modules by graded submodules are graded). 

			\item 
			Introductions to \textbf{\fref{chapter:filtrations-and-pershoms}} 
			and to \textbf{\fref{chapter:matrix-calculation}} are added.

			\item 
			Errors in the commentary and examples of 
			\textbf{\fref{section:persistent-homology}} were fixed. 
			Additional commentary and images were included at several points.

			\item 
			\textbf{\fref{section:simplicial-persistent-homology}} 
			had several significant errors, 
			e.g.\ missing definitions for objects used in the next chapter, incorrect examples.
			Some crucial results, e.g.\ \fref{prop:simp-pers-hom-graded-to-persmod}, 
			were also not identified in the previous version. 

			\item 
			Section 4.3 \textit{``Graded Invariant Factor Decompositions and SNDs''} of the old version was divided into three sections for better exposition:
			\textbf{\fref{section:calculation-graded-ifds}}, 
			\textbf{\fref{section:matrix-reduction-of-graded-matrices}}, 
			and \textbf{\fref{section:snd-algorithm-in-the-graded-case}}. 
			The content is mostly unchanged, 
			with some additional commentary and minor changes to formatting.

			\item 
			\textbf{\fref{section:matrix-graded-chain-complex}}
			\textit{``Matrix Calculation of Homology of Graded Chain Complexes''} was missing in the old version because I forgot to uncomment the \texttt{import} line in my \LaTeX{} set-up and somehow missed this in the review. 
			This section contains a brief explanation about how the results for ungraded chain complexes extend nicely for graded chain complexes, and another example calculation of persistent homology.

			\item 
			The list of symbols in \textbf{\fref{appendix:miscellaneous}} is now grouped by topic. 
			Corrections made in the previous chapters are also reflected here. 
			Several errors and formatting issues were corrected in \textbf{\fref{appendix:matrix-theory}}.
			Some comments made in the previous chapters involving modules were collected and moved to \textbf{\fref{appendix:module-theory}}. 
		\end{enumerate}
		Please note that a significant part of the edits done in this version have not been reviewed or approved by my advisor (Dr.\ Christine Escher) as these were made after I had graduated.
		
		\vspace{\parskip}
		\hfill --- Jason Ranoa
	\end{blockquote}

	\clearpage

\renewcommand*\contentsname{Table of Contents}
\pdfbookmark{\contentsname}{toc}
\setcounter{tocdepth}{1}
\begin{spacing}{1.2} 
	\tableofcontents 
\end{spacing} 
\clearpage 


\addcontentsline{toc}{chapter}{Introduction}
\chapter*{Introduction}

\noindent
Topological data analysis (TDA) is a relatively new field of study that seeks real-world applications of the theory of algebraic topology.
One such application involves the identification of features on finite datasets (as finite metric subspaces of $\reals^N$) by investigating a family of topological spaces constructed using said dataset over some parameter.

For example, given the dataset $\mathbb{X} = \set{x_1, x_2, x_3, x_4} \subseteq \reals^2$ with $x_1 = (2,22)$, $x_2 = (6,15)$, $x_3 = (22,30)$, and $x_4 = (30,15)$, we can construct a family $\set{X_r: r \in \reals_{\geq}}$ of topological subspaces of $\reals^2$ by taking the union of $2$-disks centered at each $x_i$ with radius $r \in [0,\infty)$. This construction produces 6 subspaces distinct up to homeomorphism, as illustrated below:

\begin{center}\setlength{\tabcolsep}{1pt}
\newlength{\colballlength}\setlength{\colballlength}{0.16\linewidth}
\small
\begin{tabular}{ccc ccc}
	\includegraphics[width=\colballlength]{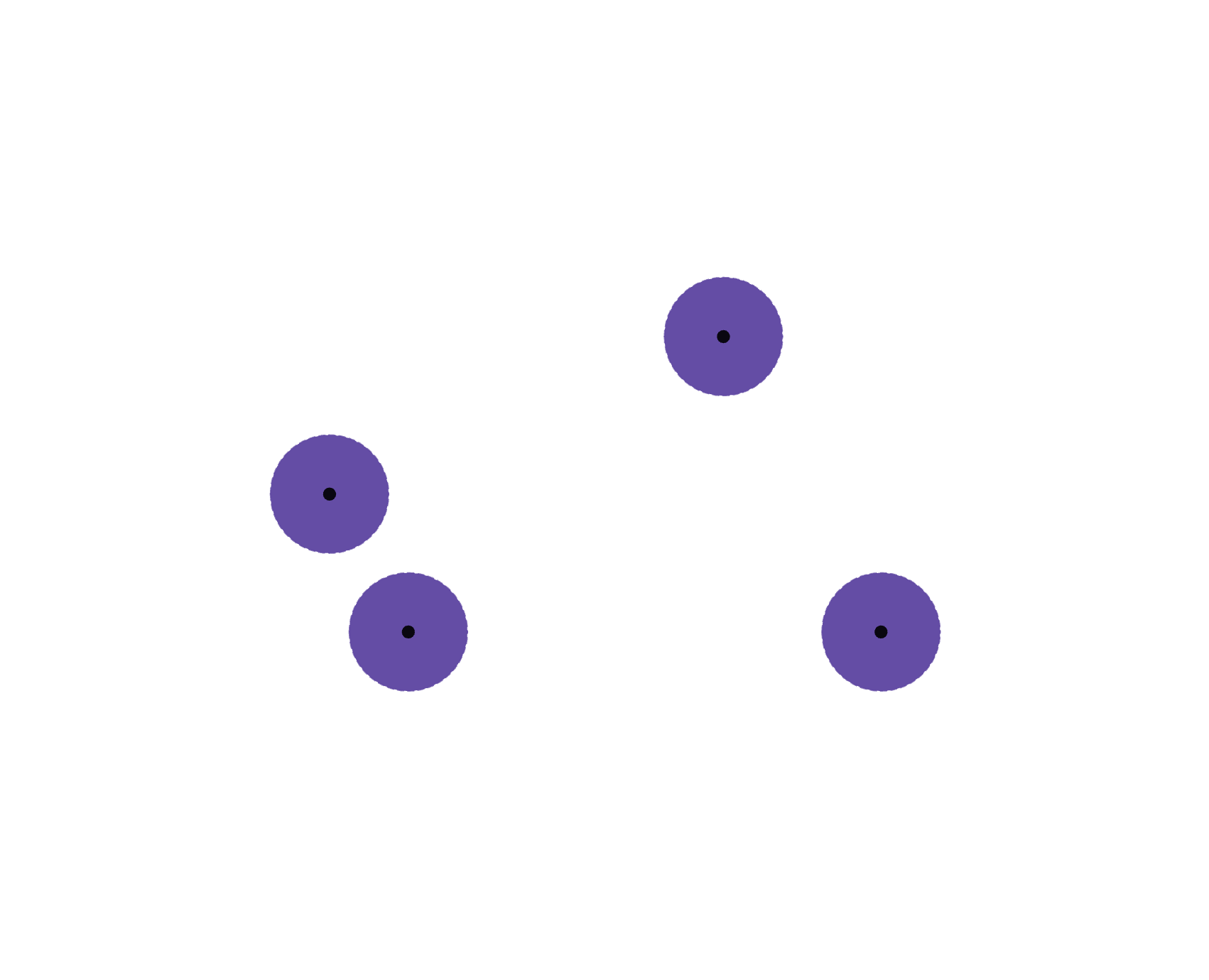}
	& \includegraphics[width=\colballlength]{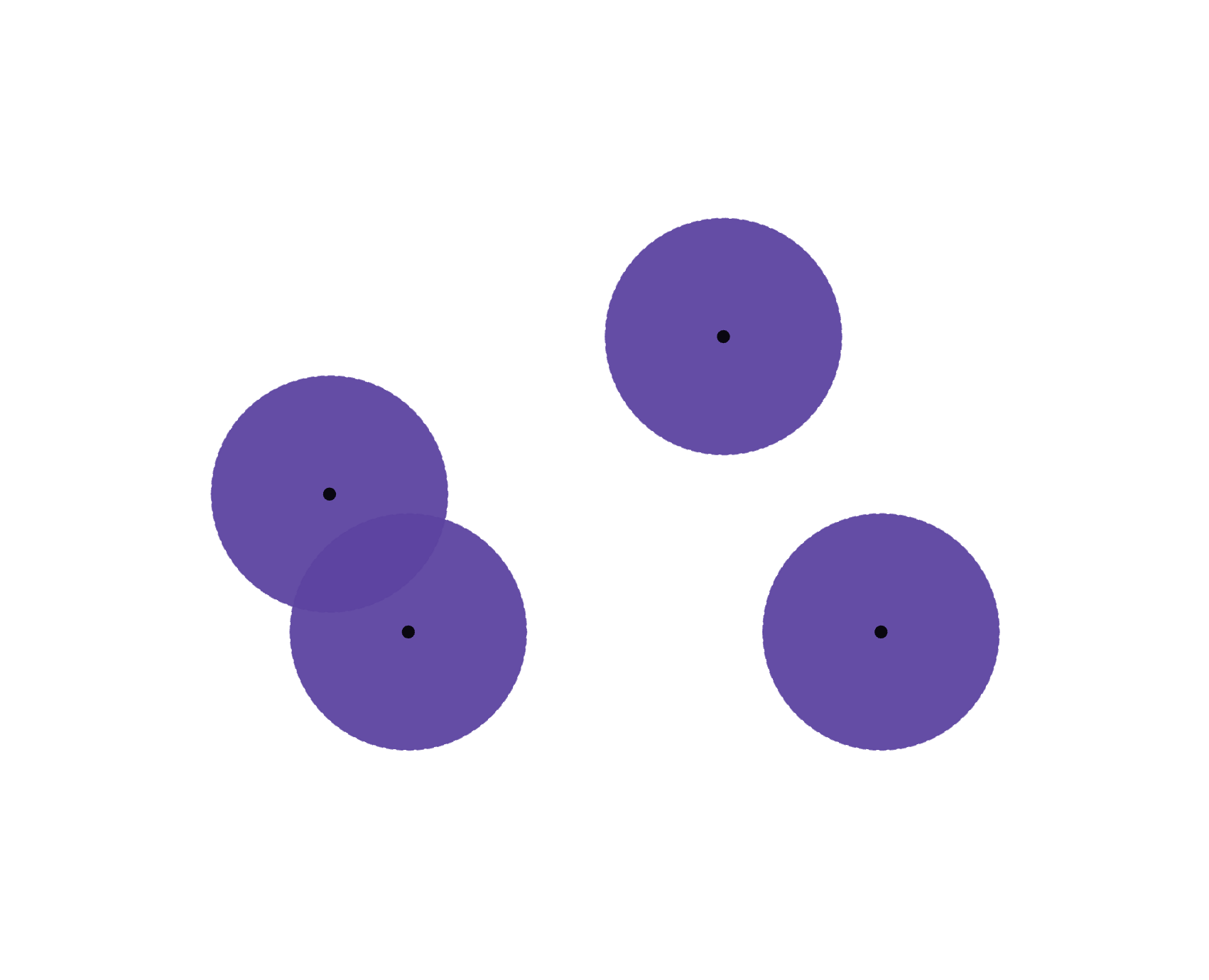}
	& \includegraphics[width=\colballlength]{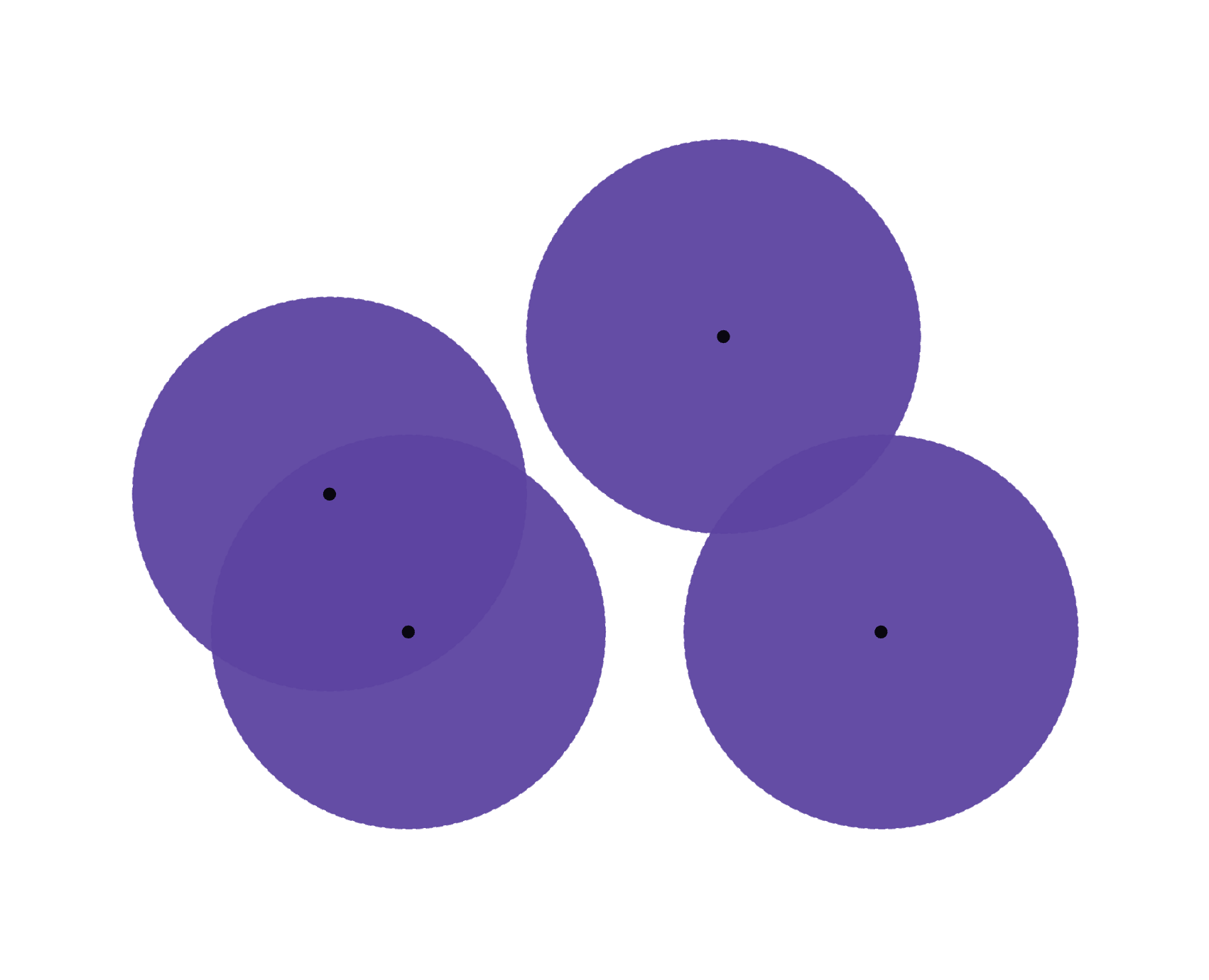}
	& \includegraphics[width=\colballlength]{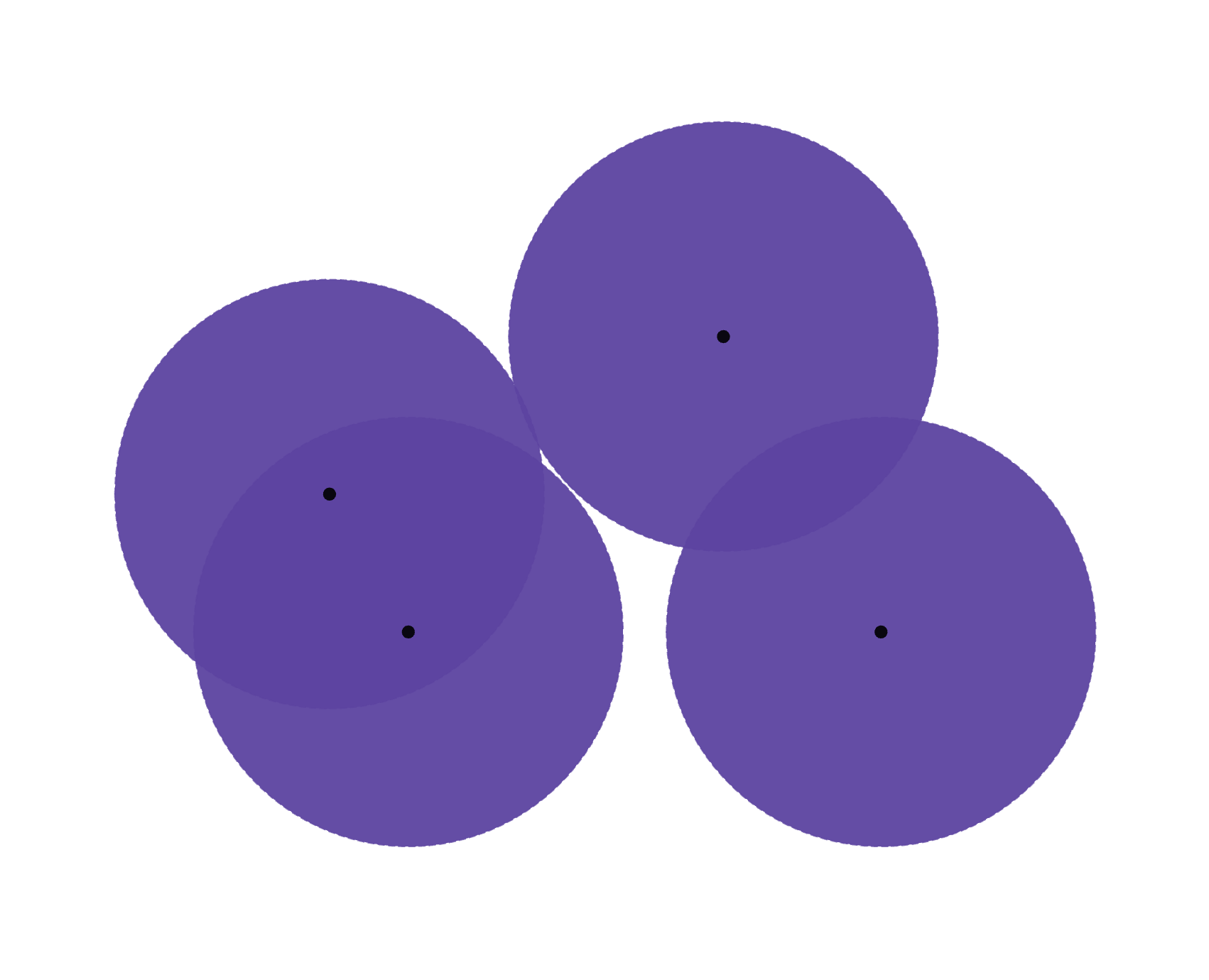}
	& \includegraphics[width=\colballlength]{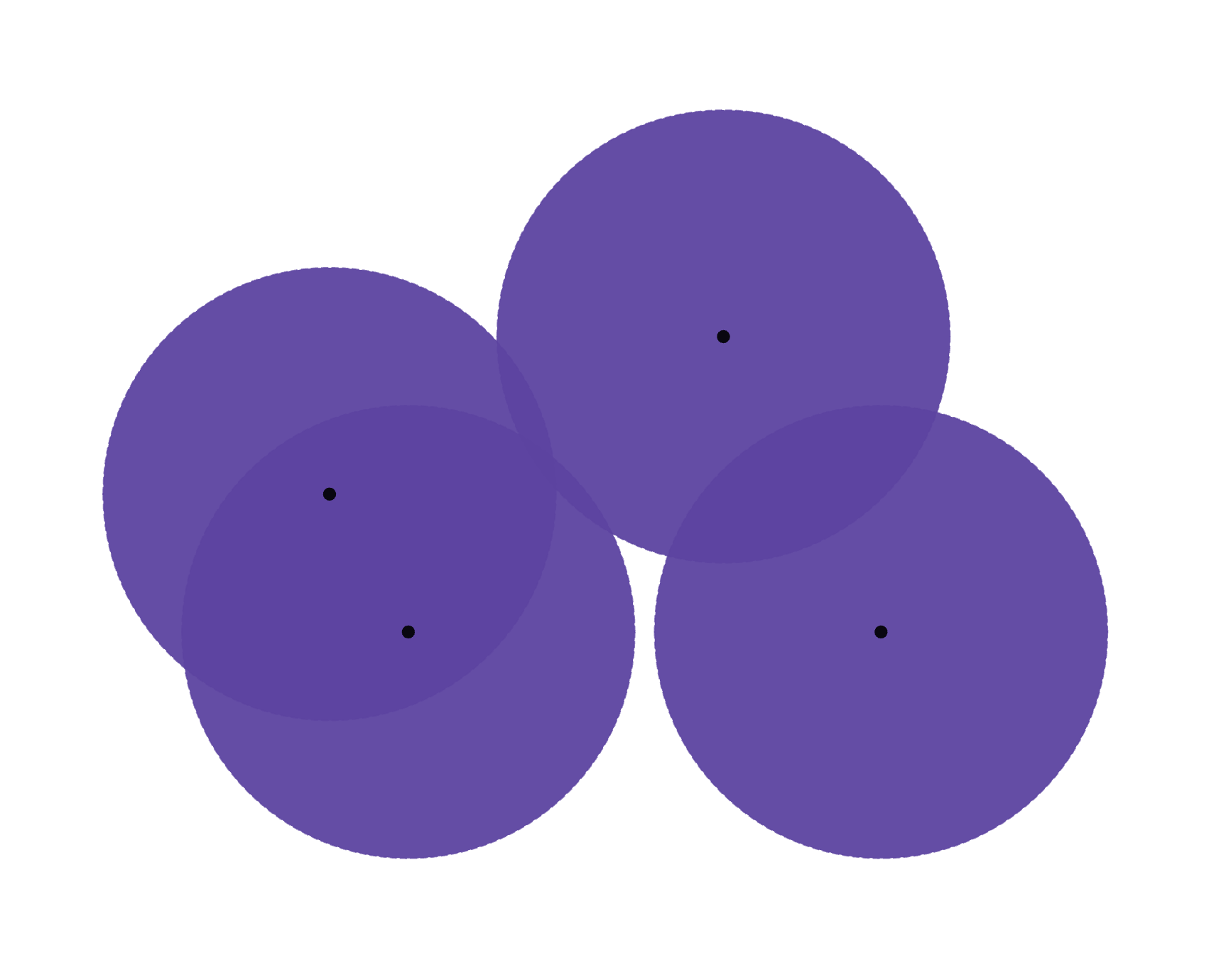}
	& \includegraphics[width=\colballlength]{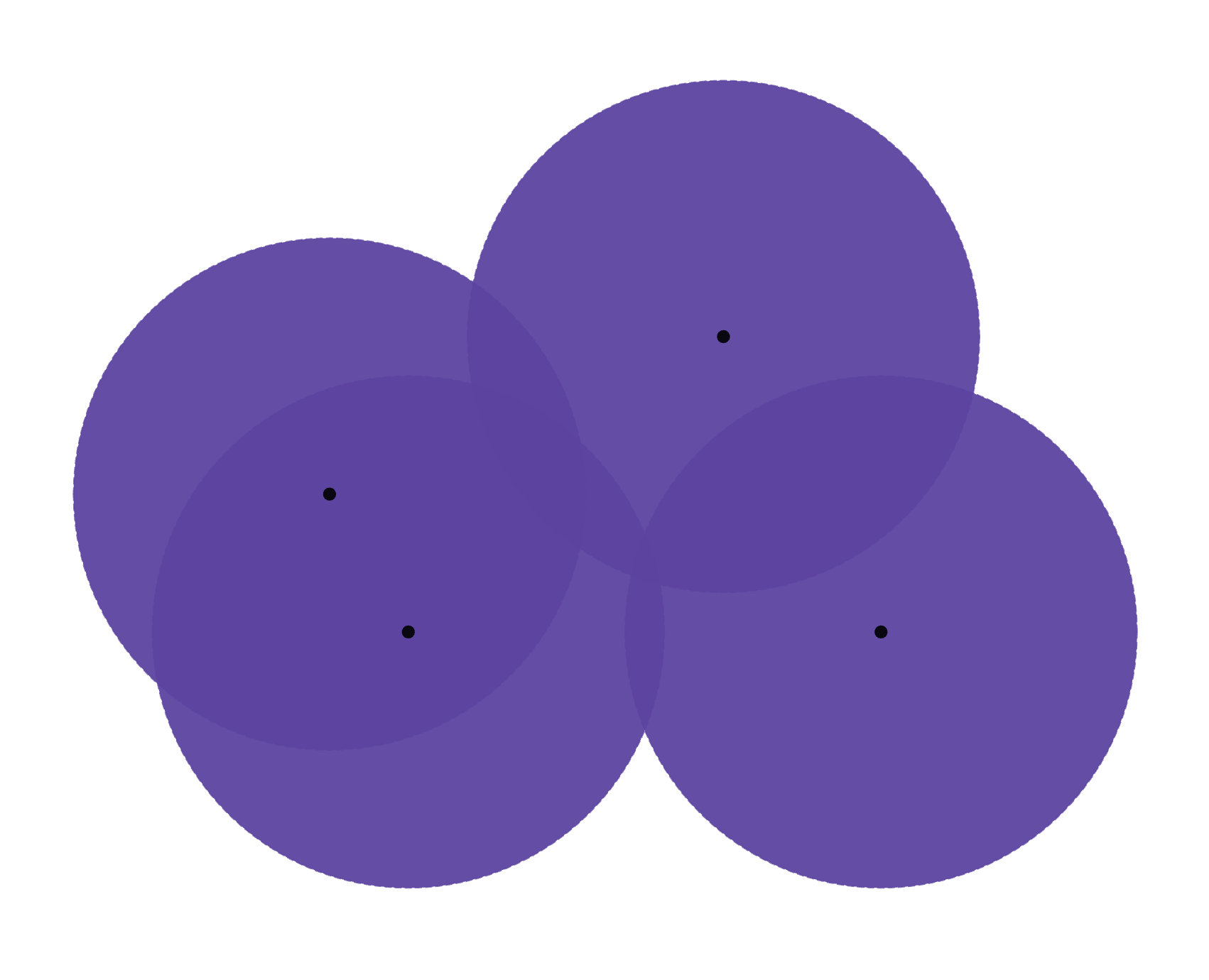}
	\\[-5pt] 
	$X_3$ & $X_6$ & $X_{10}$ & $X_{10.9}$ & $X_{11.5}$ & $X_{13}$
	\\
	$(r=3)$ & $(r=6)$ & $(r=10)$ & $(r=10.9)$ & $(r=11.5)$ & $(r=13)$
\end{tabular}
\vspace{0.5\baselineskip}
\end{center}

\noindent 
The collection of homology groups $H_n(X_r)$ with $r \in \reals_{\geq 0}$ and $n \in \nonnegints$ and the images of the maps $H_n(X_r) \to H_n(X_s)$ on homology induced by the inclusions $i^{s,r}: X_r \to X_s$ is what we call the \textit{persistent homology} of the family $\set{X_r}$.
This persistent homology then is used to determine certain characteristics of the dataset $\mathbb{X}$.
In the example above, we can determine that the points of $\mathbb{X}$ are near one another since the space $X_r$ merges into one path component as early as $r = 10.9$. 

In practice, the calculation of persistent homology happens at the level of simplicial complexes.
Following the example above,
	an abstract simplicial complex $C_r(\mathbb{X})$ 
	called the \textit{\v{C}ech complex} of $\mathbb{X}$ with parameter $r$ is constructed for each $r \in \reals_{\geq 0}$ using the following rule:
\begin{blockquote}
	For each subset $S$ of $\mathbb{X}$ with $n = \card(S) \geq 1$,
		$S$ is an $(n-1)$-simplex of $C_r(\mathbb{X})$
		if and only if the collection of the $2$-disks centered at each $x_i \in S$ with radius $r$ has non-trivial intersection as subspaces of $\reals^2$.
\end{blockquote}
There are seven distinct abstract simplicial complexes in the family $\set{C_r(\mathbb{X}) : r \in \reals_{\geq 0}}$ of \v{C}ech complexes, 
	as illustrated below.
For clarity, the elements $\set{x_1, x_2, x_3, x_4}$ of $\mathbb{X}$ are denoted as $\set{a,b,c,d}$ respectively when used as vertices of $C_r(\mathbb{X})$.

\begin{center}\setlength{\tabcolsep}{1pt}
\setlength{\colballlength}{0.135\linewidth}
\small\vspace{\baselineskip}
\begin{tabular}{ccc ccc c}
	\includegraphics[width=\colballlength]{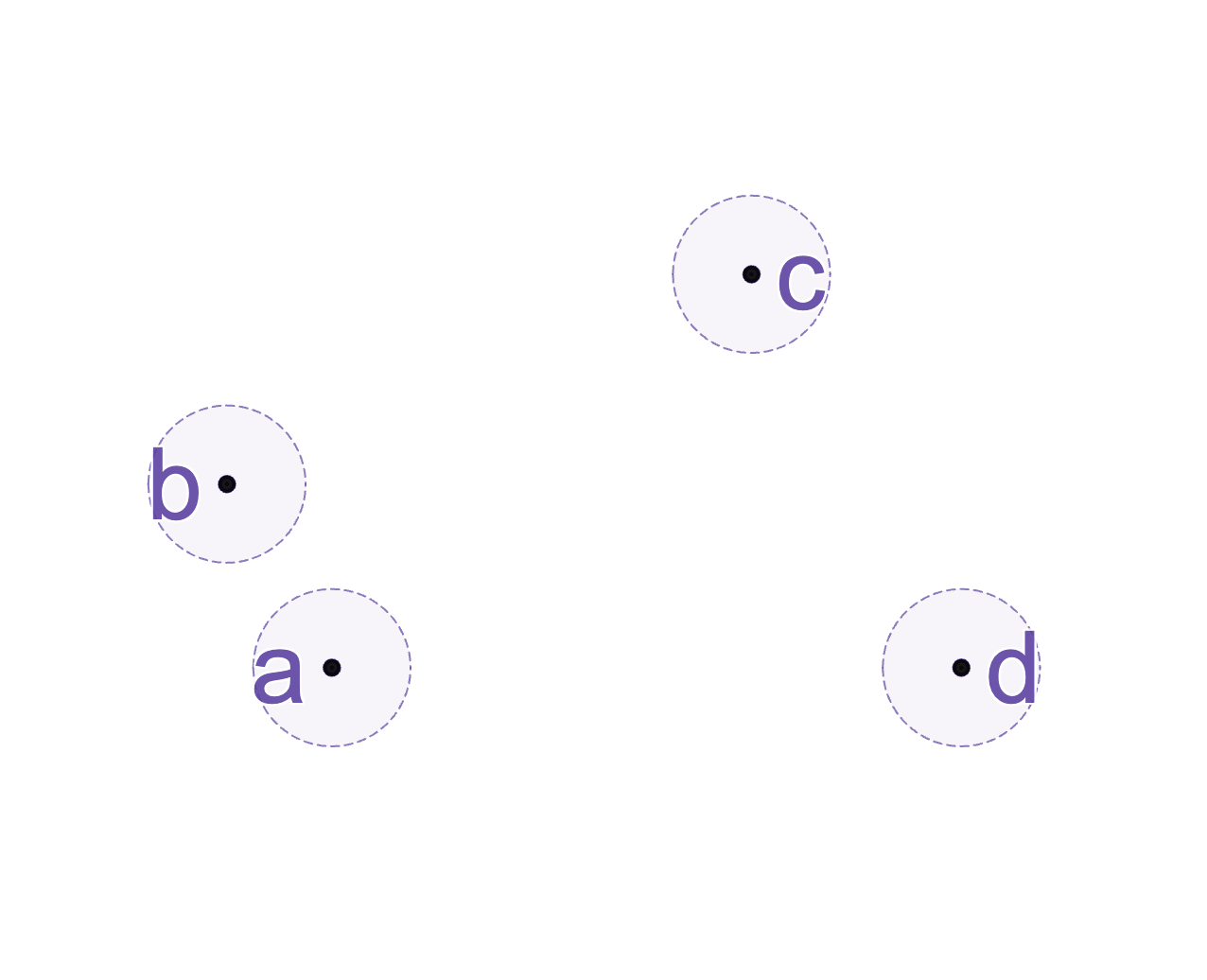}
	& \includegraphics[width=\colballlength]{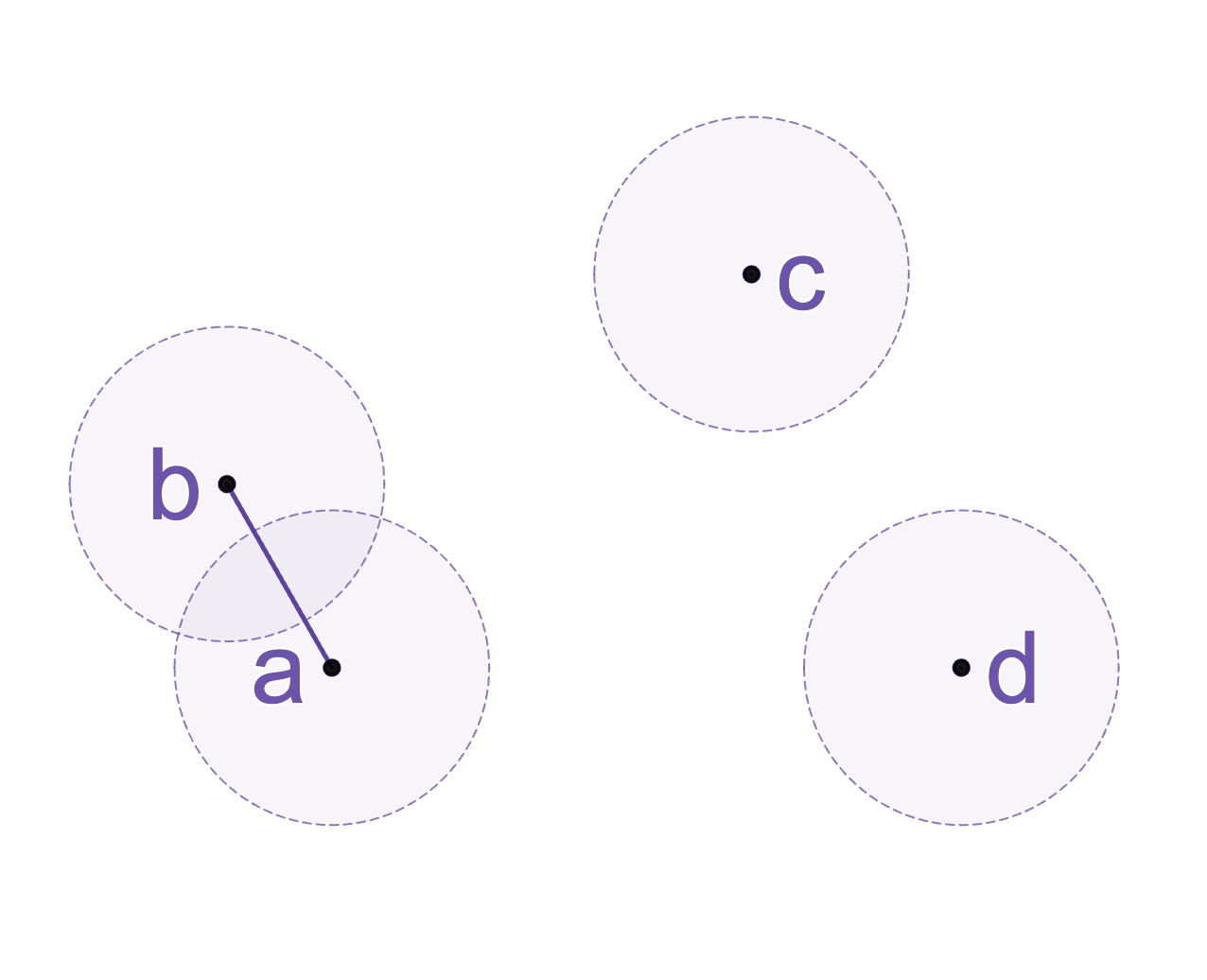}
	& \includegraphics[width=\colballlength]{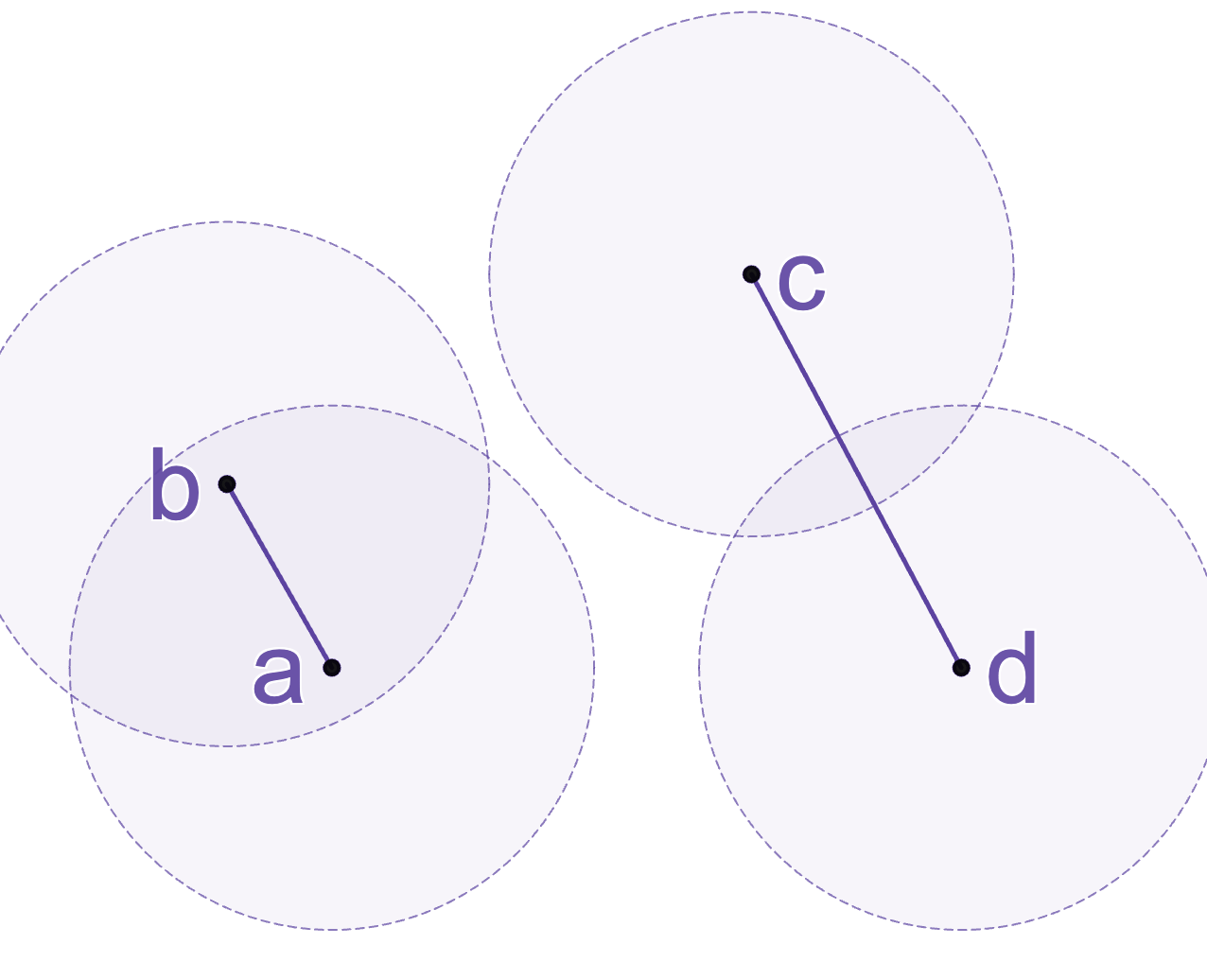}
	& \includegraphics[width=\colballlength]{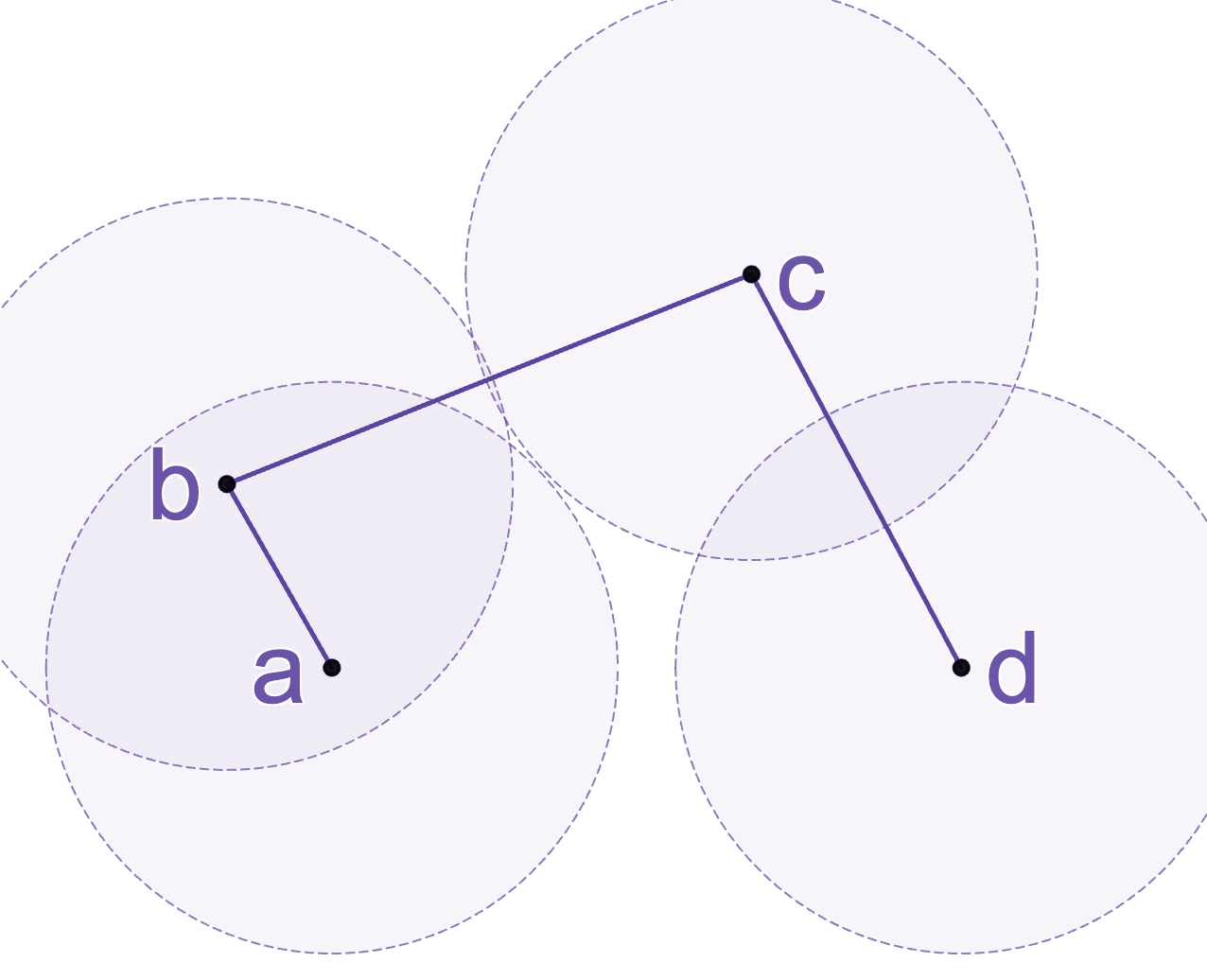}
	& \includegraphics[width=\colballlength]{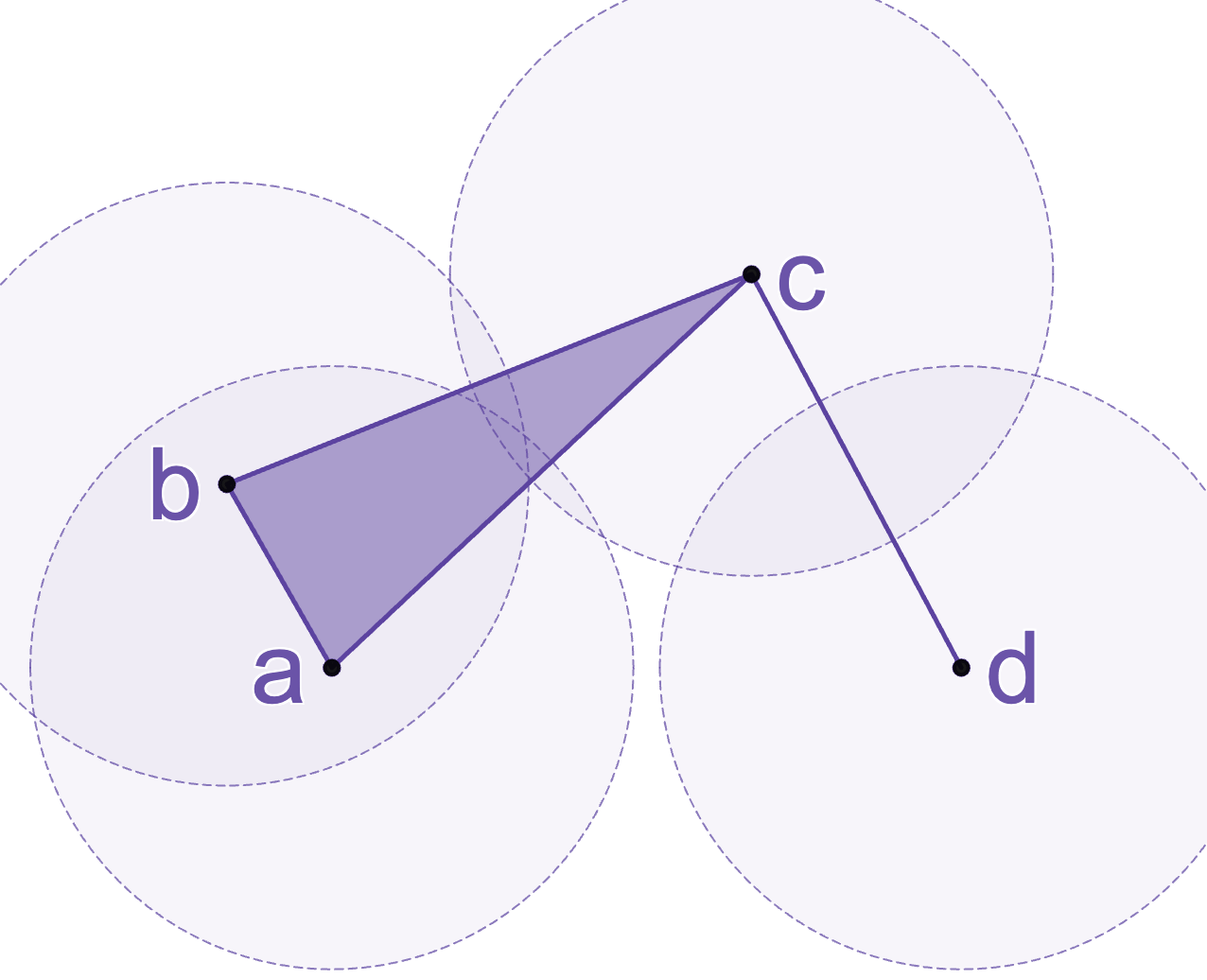}
	& \includegraphics[width=\colballlength]{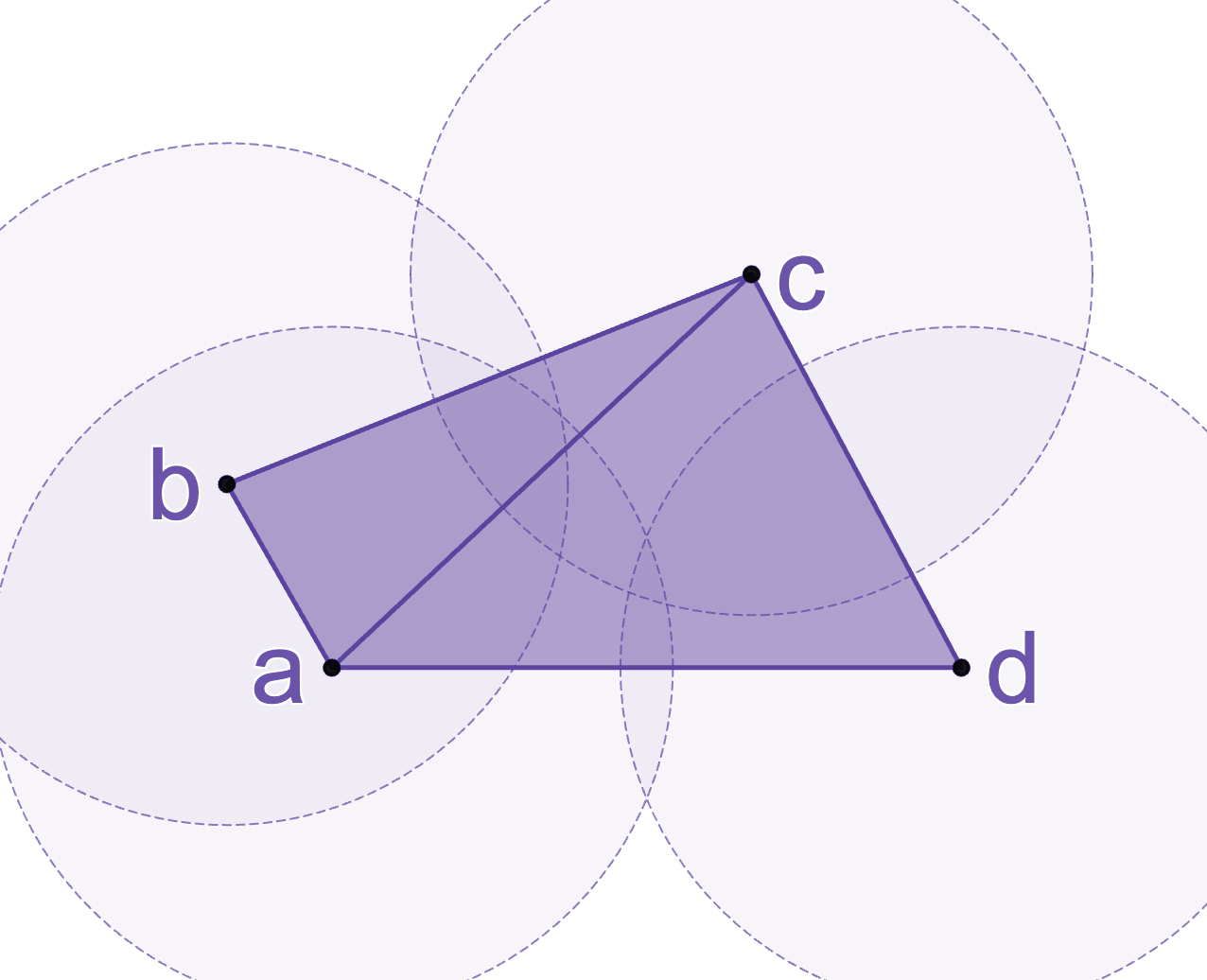}
	& \includegraphics[width=\colballlength]{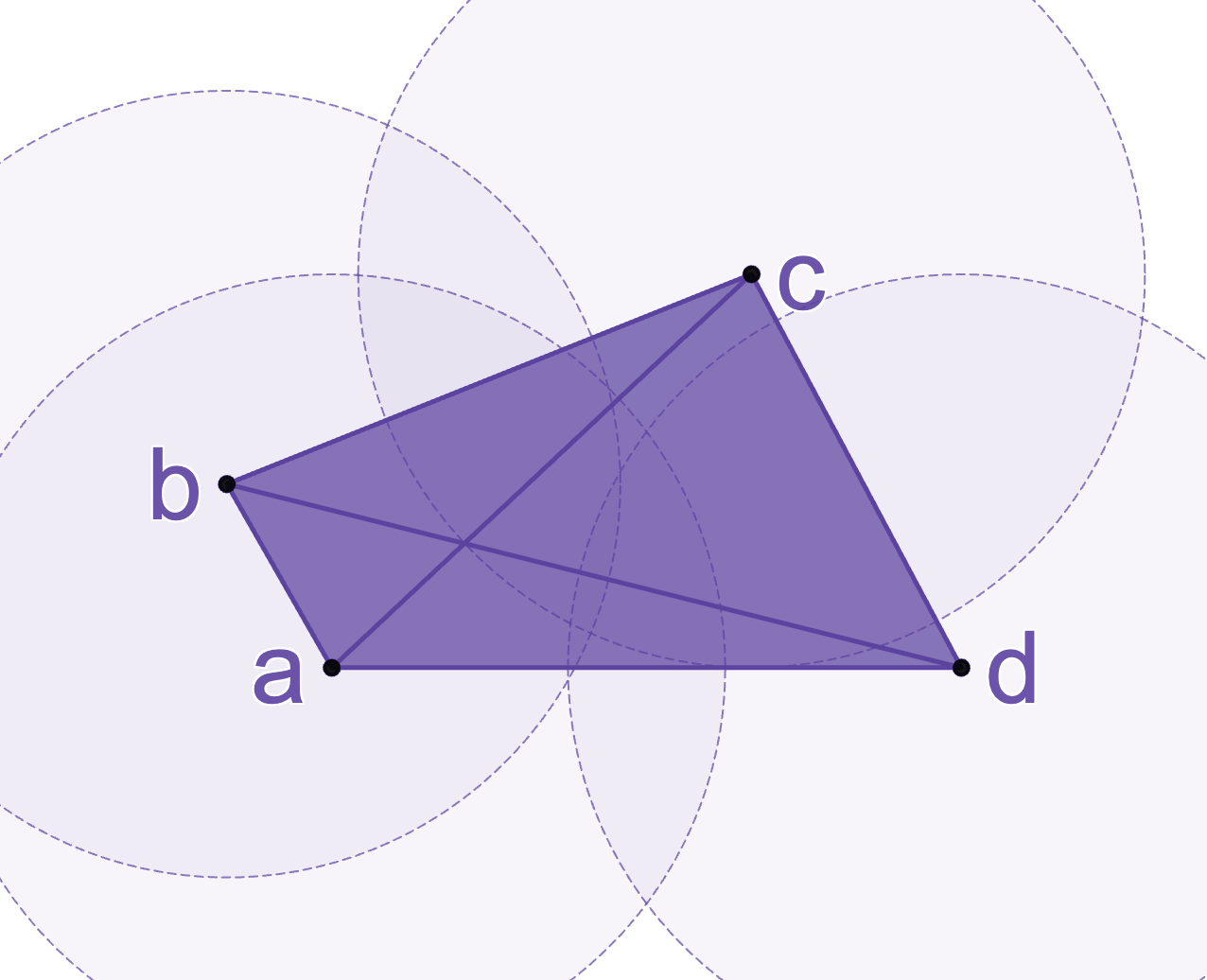}
	\\[-1pt] 
	$C_3(\mathbb{X})$ 
	& $C_{6}(\mathbb{X})$
	& $C_{10}(\mathbb{X})$
	& $C_{10.9}(\mathbb{X})$
	& $C_{11.5}(\mathbb{X})$
	& $C_{13}(\mathbb{X})$
	& $C_{15}(\mathbb{X})$
	\\[0pt] 
	added $a,b,c,d$
	& added $ab$
	& added $cd$
	& added $bc$
	& added $abc$
	& added $ad, acd$
	& full complex of $\mathbb{X}$ 
\end{tabular}\vspace{\baselineskip}
\end{center}

\noindent 
The collection $\set{ H_n(C_r(\mathbb{X})) : r \in \reals_{\geq 0}} $ of simplicial homology groups then determines the persistent homology of $\mathbb{X}$ by the \textit{Nerve Theorem},
	which primarily states that $H_n(X_r) \cong H_n(C_r(\mathbb{X}))$ as abelian groups
	for all $r \in \reals_{\geq 0}$ and $n \in \nonnegints$.

In this expository paper, we consider a more general case and discuss the persistent homology $H_n(\filt{K};\field)$ of $\nonnegints$-indexed filtrations $\filt{K} = \set{K_t}_{t \in \nonnegints}$ of abstract simplicial complexes (wherein $K_t \subseteq K_s$ whenever $t \leq s$).
More specifically, we talk about the algebraic background needed for the \textit{matrix reduction algorithm for persistent homology} discussed in the following paper:
\begin{center}
	\textit{Computing Persistent Homology} \cite{matrixalg:zomorodian} by Afra Zomorodian and Gunnar Carlsson.
\end{center}
The concepts discussed in this paper were then used as basis for a number of software packages such as Ripser \cite{ripser}, a \texttt{C++} package written by Ulrich Bauer that calculates the persistent homology of Vietoris-Rips filtrations, taking finite metric subspaces of $\mathbb{R}^N$ as input.

This expository paper talks about the concepts in \cite{matrixalg:zomorodian} in the perspective of category theory.
We have included an introductory-level discussion of category theory in
\fref{appendix:cat-theory}.
For general texts on category theory, we recommend 
\textit{Category Theory in Context} \cite{cattheory:rhiel} by Emily Riehl and \textit{Introduction to Homological Algebra} \cite{cattheory:rotman} by Joseph J.\ Rotman.

Much like how the homology groups $H_n(K;R)$ of a simplicial complex $K$ with coefficients in a principal ideal domain (PID) $R$ are encoded as $R$-modules,
	we characterize the persistent homology $H_n(\filt{K};\field)$ of a filtration $\filt{K}$ with coefficients in a field $\field$ as a \textit{persistence module} over $\field$.
This paper is organized as follows:
\begin{enumerate}[left=1in]
	\item[In \textbf{\fref{chapter:simplicial-homology}}:] 
	We provide necessary background information involving simplicial complexes and simplicial homology with coefficients in a PID $R$,
		as well as identify a number of results that will be useful in the characterization of persistent homology.

	\item[In \textbf{\fref{chapter:persistence-theory}}:] 
	We define persistence modules over a field $\field$ as a \textit{diagram},
		i.e.\ as functors $\posetN \to \catvectspace$,
	describe the interval decomposition of a persistence module (which is unique up to persistence isomorphism),
	and present a category equivalence between the category $\catpersmod$ of persistence modules and the category $\catgradedmod{\field}$ of graded $\field[x]$-modules.

	\item[In \textbf{\fref{chapter:filtrations-and-pershoms}}:] 
	We introduce simplicial filtrations as functors of the form $\posetN \to \catsimp$ and characterize the persistent homology of these filtrations as persistence modules.
	We also describe how simplicial homology can be extended to the case of persistence modules, which we call \textit{simplicial persistent homology}, and discuss the objects in $\catgradedmod{\field}$ brought about by applying the category equivalence between persistence modules and graded modules. 

	\item[In \textbf{\fref{chapter:matrix-calculation}}:] 
	We discuss how the invariant factor decomposition of the homology of free chain complexes of $R$-modules can be calculated using the Smith Normal Decomposition (SND) of matrices over $R$
	and how this calculation can be extended to the case of graded $\field[x]$-modules.
	We also discuss how this method of calculation, along with \textit{simplicial persistent homology}, can be used to find the interval decomposition of persistence modules corresponding to persistent homology.
\end{enumerate}

\HIDE{Given a family $\set{K_t}_{t \in \nonnegints}$ of simplicial complexes such that $K_t$ is a subcomplex of $K_s$ whenever $t \leq s$,
	we construct a filtration $\filt{K}: \posetN \to \catsimp$ 
	by $\filt{K}(t) = K_t$.
In \fref{section:persistent-homology}, we define the $n$\th persistent homology module $H_n(\filt{K};\field)$ of $\filt{K}$ with coefficients in a field $\field$ to be the persistence module represented by the following sequence:
\begin{equation*}
	H_n(\filt{K};\field) \quad:\quad \bigg[\quad
	H_n(K_0;\field) 
		\Xrightarrow{\quad i^{0}_\ast \quad}
	H_n(K_1;\field) 
		\Xrightarrow{\quad i^{1}_\ast \quad}
	H_n(K_2;\field) 
		\Xrightarrow{\quad i^{2}_\ast \quad} \cdots
	\quad\bigg]
\end{equation*}
where $i^{t}_\ast: H_n(K_t;\field) \to H_n(K_{t+1};\field)$ is the map on homology induced by the inclusion $i^t: K_t \to K_{t+1}$.
The calculation of the interval decomposition of $H_n(\filt{K};\field)$ as a persistence module is done as follows:
\begin{enumerate}[label={Step \arabic*.}, left=0.3in]
	\item 
	Following \fref{section:construction-of-persistent-homology},
	we construct the simplicial persistence complex $C_\ast(\filt{K};\field)$, which is a chain complex of persistence modules $C_n(\filt{K};\field)$ and persistence morphisms $\boundary^{\,\bullet}_n: C_n(\filt{K};\field) \to C_{n-1}(\filt{K};\field)$, as illustrated below:
	\begin{equation*}
		C_\ast(\filt{K};\field) \quad:\quad \bigg[\hspace{5pt}
			\cdots 
				\Xrightarrow{\quad}
			C_{n+1}(\filt{K};\field) 
				\Xrightarrow{\,\, \boundary^{\,\bullet}_{n+1} \,\,}
			C_n(\filt{K};\field) 
				\Xrightarrow{\,\, \boundary^{\,\bullet}_n \,\,}
			C_{n-1}(\filt{K};\field) 
				\Xrightarrow{\quad}
			\cdots
		\hspace{5pt}\bigg]
	\end{equation*}
	Like in the case of simplicial homology,
		the simplicial persistent homology 
		$H_n(\filt{K};\field)$ of $\filt{K}$ corresponds to the $n$\th homology of this chain complex, 
		i.e.\ 
		$H_n(\filt{K};\field) = 
			\ker(\boundary^{\,\bullet}_n) \bigmod 
		\im(\boundary^{\,\bullet}_{n+1})$ (as persistence modules over $\field$).

	\item 
	Using the functor $\togrmod: \catpersmod \to \catgradedmod{\field}$ discussed in \fref{section:cat-equiv-graded-modules},
	we construct a graded chain complex $C_\ast\graded(\filt{K};\field)$
	of graded $\field[x]$-modules $C_n\graded(\filt{K};\field) := \togrmod(C_n(\filt{K};\field))$ and graded differentials 
	$\boundary_n\graded := \togrmod(\boundary^{\,\bullet}_n)$, as illustrated below:
	\begin{equation*}
		C_\ast\graded(\filt{K};\field) \quad:\quad \bigg[\hspace{5pt}
			\cdots 
				\Xrightarrow{\quad}
			C_{n+1}\graded(\filt{K};\field) 
				\Xrightarrow{\,\, \boundary\graded_{n+1} \,\,}
			C_n\graded(\filt{K};\field) 
				\Xrightarrow{\,\, \boundary\graded_n \,\,}
			C_{n-1}\graded(\filt{K};\field) 
				\Xrightarrow{\quad}
			\cdots
		\hspace{5pt}\bigg]
	\end{equation*}

	\item 
	Following \fref{chapter:matrix-calculation},
		we calculate the graded invariant factor decomposition of 
		the $n$\th chain homology $H_n\graded(\filt{K};\field) =
		\ker(\boundary\graded_n) \bigmod 
		\im(\boundary\graded_{n+1})$
		of this graded chain complex
		using the SNDs of $[\boundary\graded_{n+1}]$ and $[\boundary\graded_n]$.
	This graded invariant factor decomposition is of the following form:
	\begin{equation*}
		H_n\graded(\filt{K};\field) \upgraded\cong
		\Sigma^{s_1}\!\paren{\frac{\field[x]}{(x^{t_1})}}
			\oplus \cdots \oplus 
			\Sigma^{s_r}\!\paren{\frac{\field[x]}{(x^{t_r})}}
			\oplus 
			\Sigma^{s_{r+1}}\field[x]
			\oplus \cdots \oplus 
			\Sigma^{s_{m}}\field[x]
	\end{equation*}
	where $\upgraded\cong$ refers to an isomorphism in the category $\catgradedmod{\field}$ of graded $\field[x]$-modules and graded homomorphisms.
	
	Note that the term \textit{matrix reduction algorithm} generally refers to this step in the computation,
		wherein optimizations are applied for faster and more efficient calculation.
	For example, the Ripser package \cite{ripser} simulates this calculation on matrices over $\field = \ints_p$, as opposed to matrices over $\ints_p[x]$, and uses several results to justify skipping specific steps in the reduction (e.g.\ clearing).

	\item 
	Using the functor $\topersmod: \catgradedmod{\field} \to \catpersmod$ discussed in \fref{section:cat-equiv-graded-modules},
	we transform this graded invariant factor decomposition into an interval decomposition of $H_n(\filt{K};\field)$ as follows:
	\begin{align*}
		H_n(\filt{K};\field) 
		&\,\uppersmod\cong\,
		\topersmod\Biggl(
			\Sigma^{s_1}\!\paren{\frac{\field[x]}{(x^{t_1})}}
			\oplus \cdots \oplus 
			\Sigma^{s_r}\!\paren{\frac{\field[x]}{(x^{t_r})}}
			\oplus 
			\Sigma^{s_{r+1}}\field[x]
			\oplus \cdots \oplus 
			\Sigma^{s_{m}}\field[x]
		\Biggr) \\
		&\,\uppersmod\cong\,
		\intmod{[s_1, s_1+t_1)}
		\oplus \cdots \oplus 
		\intmod{[s_r, s_r+t_r)}
		\oplus 
		\intmod{[s_{r+1},\infty)}
		\oplus \cdots \oplus 
		\intmod{[s_{m},\infty)}
	\end{align*}
	where $\uppersmod\cong$ refers to an isomorphism in the category $\catpersmod$ of persistence modules over $\field$.
\end{enumerate}
In summary, the computation of the persistent homology is done at the level of graded $\field[x]$-modules and is allowed by the category equivalence between $\catpersmod$ and $\catgradedmod{\field}$.

\clearpage}


\chapter{Simplicial Homology}
\label{chapter:simplicial-homology} 
 
In this chapter, we review a number of definitions and results involving the simplicial homology of simplicial complexes.
In the field of persistent homology theory, 
	the term \textit{simplicial complex} is conventionally used to refer to an \textit{abstract simplicial complex}

Note that this chapter is not supposed to be a rigorous or thorough treatment of simplicial homology theory but instead covers constructions that are immediately relevant in our discussion of persistent homology.
As such, proofs for most of the propositions and theorems in this chapter are not provided in this paper but can be found in standard algebraic topology textbooks.
This chapter is divided into the following sections:
\begin{enumerate}[chapterdecomposition, listparindent=\parindent]
	\item \textbf{Simplicial Complexes and Geometric Realizations} 

	\noindent
	We discuss abstract simplicial complexes, which are representations of certain topological spaces, and review some terminology and results.

	Note that, after this section, we often drop the modifier \textit{abstract} when referring to an \textit{abstract simplicial complex}, following convention.
	
	\item \textbf{Simplicial Homology with Coefficients in a PID}

	\noindent
	We provide definitions for the constructions involved in the simplicial homology of (abstract) simplicial complexes and identify homology as an invariant of the homeomorphism type of topological spaces.

	\item \textbf{Functorial Constructions in Simplicial Homology}
	
	\noindent
	We define the category $\catsimp$ of (abstract) simplicial complexes and simplicial maps and discuss why the constructions discussed in the previous section correspond to functors $\catsimp \to \catmod{R}$ and $\catsimp \to \catchaincomplex{\catmod{R}}$ where 
		$\catmod{R}$ refers to the category of $R$-modules and $R$-module homomorphisms 
		and $\catchaincomplex{\catmod{R}}$ refers to the corresponding chain complex category.

	This functorial perspective of simplicial homology will be relevant in \fref{chapter:filtrations-and-pershoms} when we discuss persistent homology.

	

\end{enumerate}
Listed below are the main references used in this chapter (in order of decreasing relevance):
\begin{enumerate}
	\item 
	\textit{An Introduction to Algebraic Topology} \cite{algtopo:rotman} by Joseph J. Rotman.

	\item 
	\textit{Algebraic Topology} \cite{algtopo:hatcher} by Allen Hatcher.

	\item 
	\textit{Elements of Algebraic Topology}~\cite{algtopo:munkres} by James Munkres.
\end{enumerate}

	\clearpage

\section{Simplicial Complexes and Geometric Realizations}
\label{section:simplicial-complexes}

In this section, we define a combinatorial representation of topological spaces called a \textit{simplicial complex}.
There are generally two types of simplicial complexes used in algebraic topology:
\begin{enumerate}
	\item Abstract Simplicial Complexes (given in \fref{defn:simp-complex-abstract}), which are combinatorial representations.
	We have found that, in the field of persistent homology, the term \textit{simplicial complex} conventionally refers to this type of simplicial complex. We follow this convention in this expository paper.

	\item Geometric Simplicial Complexes (given in \fref{defn:simp-complex-geometric}), which are representations of topological spaces as a collection of subsets of $\reals^N$ for sufficiently high $N$.
	Note that, in the more general field of algebraic topology, the term \textit{simplicial complex} usually refers to this type of simplicial complex, as seen in \cite{algtopo:hatcher} and \cite{algtopo:munkres}.
\end{enumerate}
For clarity, we do not drop the modifier \textit{abstract} and \textit{geometric} when discussing simplicial complexes in this section. 
However, under certain finiteness conditions, there is a correspondence between these two types of simplicial complexes.

We start with a definition for abstract simplicial complexes, taken from \cite[p141]{cattheory:rotman}. 

\begin{definition}\label{defn:simp-complex-abstract}
	An \textbf{abstract simplicial complex} $K$ is a collection of nonempty subsets of some set $V$ such that 
		the following properties are satisfied:
	\begin{enumerate}
		\item For all $v \in V$, $\set{v} \in K$, i.e.\ all vertices are in $K$. 
		\item 
		For all sets $\tau \in K$, if $\sigma$ is non-empty subset of $\tau$, then $\sigma \in K$.
	\end{enumerate}
	The set $V$ is called the \textbf{vertex set} of $K$, denoted $\Vertex(K) := V$, and an element $v \in V$ is called a \textbf{vertex} of $K$.
	If $V$ is finite, $K$ is called a \textbf{finite simplicial complex}.

	A \textbf{simplex} $\sigma$ is a set in $K$.
	If $\card(\sigma) = n+1$, then $\sigma \in K$ is called an $n$-\textbf{simplex} of $K$ with \textbf{dimension} $n \in \nonnegints$ denoted by $\dim(\sigma) = n$.
	A \textbf{face} of $\sigma$ is another simplex $\tau$ of $K$ such that $\tau$ is a subset of $\sigma$, denoted $\tau \subseteq \sigma$.
	If $\dim(\tau) = \dim(\sigma) - 1$, $\tau$ is called a \textbf{facet} of $\sigma$.

	The \textbf{dimension}  of $K$, denoted $\dim(K)$, is the maximum dimension of its simplices. In this case, we say $K$ is $n$-\textbf{dimensional} and write $\dim(K) = n$.
	A simplicial complex $L$ is called a \textbf{subcomplex} of $K$ if every simplex $\tau \in L$ is in $K$.
\end{definition}
\remarks{
	\item 
	Note that, in the context of calculating simplicial homology (as in \fref{section:simplicial-homology}) and simplicial persistent homology (later in \fref{section:construction-of-persistent-homology}),
	we often assume that the simplicial complex $K$ is finite. 
	While we make an effort to identify this finiteness condition when relevant, 
		some papers (e.g. \cite{matrixalg:zomorodian}) implicitly assume that simplicial complexes are finite.

	\item 
	Outside of this section, we often denote $n$-simplices using string notation as opposed to set notation for brevity, e.g.\ we write $abc$ to refer to the $2$-simplex $\set{a,b,c}$.
	We justify this convention later in \fref{section:simplicial-homology} under \fref{remark:shorthand-for-oriented-simplices} in the context of simplicial homology.
}

Note that $\card(A)$ refers to the cardinality of some set $A$.
In this paper, we have chosen not to use the conventional notation of $|A|$ for cardinality since the operator $|-|$ is used differently in the context of simplicial complexes.
Note that for simplices $\sigma$, $\card(\sigma) = \dim(\sigma) + 1$.
We state some remarks about the vertex set $V := \Vertex(K)$ of a simplicial complex $K$:
\begin{enumerate}
	\item 
	The vertex set $V$ is often not explicitly identified in a specification of $K$ since $V$ can be identified by 
	\begin{equation*}
	\begin{array}{ccc}
		V = \bigcup_{\sigma \in K} \sigma
		&\quad\text{ and }\quad& 
		V = \bigcup\bigl\{ \sigma \in K : \card(\sigma)=1 \bigr\} 
		\\[5pt]
		\text{i.e.\ as the union of all simplices of $K$}
		&& 
		\text{i.e.\ as the union of all $0$-simplices of $K$}
	\end{array}
	\end{equation*}
	\vspace{-\baselineskip}

	\item 
	Observe that $K$ can be seen as a subset of the power set $2^V$ of $V$. Therefore, if $V$ is finite, then $2^V$ is finite and $K \subseteq 2^V$ is also finite as a set.

	As a sidenote, the notation $2^V$ for the power set of $V$ comes from the fact that for each set $P \in 2^V$, the set $P$ can be seen as a collection of choices, one for each $v \in V$, whether either $v$ is included in $P$ or not included in $P$ (i.e.\ two possible choices for each $v \in V$).
\end{enumerate}
While the vertex set $V := \Vertex(K)$ of an abstract simplicial complex can consist of points in $\reals^N$ for sufficiently high $N$, 
	the elements of $V$ are often defined and interpreted to be indeterminates.
Consequently, the vertex set $V$ and the simplicial complex $K$ generally do not have inherently topological or geometric characteristics.
Hence, the modifier \textit{abstract} in the name \textit{abstract simplicial complexes}.
As a sidenote, some references call abstract simplicial complexes as \textit{combinatorial simplicial complexes} for the same reason.
We provide an example of an abstract simplicial complex below.

\begin{example}\label{ex:abs-simp-one}
	Let $V = \set{a,b,c,d}$ be a set of indeterminates.
	The collection $K$, as given below, is an abstract simplicial complex:
	\begin{equation*}
		K = \set{\,
			\begin{aligned}
				& \set{a}, \set{b}, \set{c}, \set{d}, \\
				& \set{a,b}, \set{a,c}, \set{b,c}, \set{a,d}, \set{b,d}, \\
				& \set{a,b,c}, \set{a,b,d}
			\end{aligned}
		\,}
	\end{equation*}
	Observe that the non-empty subsets of $\set{a,b,c}$, i.e.\ the sets $\set{a,b}, \set{b,c}, \set{a,c}$ and $\set{a},\set{b},\set{c}$ are all in $K$. The same applies to those of $\set{a,b,d}$.
\end{example}

Abstract simplicial complexes can be related to each other using simplicial maps. We provide a definition below, taken from \cite[p141]{algtopo:rotman}.

\begin{definition}\label{defn:simplicial-maps}
	Let $K$ and $L$ be simplicial complexes.
	\begin{enumerate}
		\item 
		A function $f: K \to L$ is called a \textbf{simplicial map} if for all $n$-simplices $\sigma = \set{v_0, \ldots, v_n}$ of $K$,
		$f(\sigma) = \bigl\{ f(v_0), \ldots, f(v_n) \bigr\}$ 
		is a simplex of $L$ (not necessarily of dimension $n$).

		\item 
		If there exists a simplicial map $f: K \to L$ such that 
			$\sigma$ is an $n$-simplex of $K$ if and only if $f(\sigma)$ is an $n$-simplex of $L$,
		then we call $f$ a \textbf{simplicial isomorphism} and say $K$ is \textbf{isomorphic} to $L$.
		Note that $f: K \to L$ is a simplicial isomorphism if and only if $f: K \to L$ is a bijection (between sets of sets).
	\end{enumerate}
	Note that a simplicial map $f: K \to L$ is determined by its restriction $f: \Vertex(K) \to \Vertex(L)$ into the vertex sets of $K$ and $L$.
\end{definition}

Observe that a simplicial isomorphism is a correspondence between the vertex sets of two simplicial complexes. 
As such, we can interpret a simplicial isomorphism to be a \textit{relabeling} or \textit{renaming} of the vertices of a simplicial complex.
We provide an example of this \textit{renaming} below:

\begin{example}
	Let $K$ be as given in \fref{ex:abs-simp-one} and define the simplicial complex $L$ using the simplicial map 
	$f: \Vertex(K) \to \Vertex(L)$ given by $a \mapsto v_1$, $b \mapsto v_2$, $c \mapsto v_3$, and $d \mapsto v_4$.
	Then, $L$ is as follows:
	\begin{equation*}
		K = \set{\,
			\begin{aligned}
				& \set{a}, \set{b}, \set{c}, \set{d}, \\
				& \set{a,b}, \set{a,c}, \set{b,c}, \set{a,d}, \set{b,d}, \\
				& \set{a,b,c}, \set{a,b,d}
			\end{aligned}
		\,}
		\,\,\Xrightarrow{\quad f \quad}\,\,
		L = \set{\,
		\begin{aligned}
			& \set{v_1}, \set{v_2}, \set{v_3}, \set{v_4}, \\
			& \set{v_1,v_2}, \set{v_1,v_3}, \set{v_2,v_3}, \set{v_1,v_4}, \set{v_2,v_4}, \\
			& \set{v_1,v_2,v_3}, \set{v_1,v_2,v_4}
		\end{aligned}
	\,}
	\end{equation*}
	We can also define $L$ first as above
	and define the simplicial map $f: K \to L$ by $a \mapsto v_1$, $b \mapsto v_2$, $c \mapsto v_3$, and $d \mapsto v_4$, i.e.\ writing a correspondence between the vertices of $K$ to those of $L$.
	In this case, $f$ defines a simplicial isomorphism.
\end{example}

To gain topological information from abstract simplicial complexes, 
	we relate them to geometric simplicial complexes.
As suggested by the modifier \textit{geometric}, these are composed of geometric spaces like lines, triangles, and their $n$-dimensional analogs and are, therefore, more suited for topology-related interpretations, e.g.\ as spaces in $\mathbb{R}^N$.
Before we provide a definition for geometric $n$-simplices, 
	we list some notation and terminology for affine spaces and convex hull 
	below, taken from \cite[Chapter 2]{algtopo:rotman}.
\begin{enumerate}
	\item 
	A subset $A \subseteq \reals^N$ is an \textbf{affine space} if for every pair of distinct points $x,y \in A$, the line passing through $x$ and $y$ is contained in $A$. Equivalently,
		$A \subseteq \reals^N$ is an affine space if there exists $a \in \reals^N$ and a vector subspace $V$ of $\reals^N$ such that 
		$A = a + V = \set{a + v : v \in V}$, i.e.\ an affine space is a translated vector subspace of $\reals^N$.

	\item 
	Given a set $X \subseteq \reals^N$, the \textbf{affine hull} $\aff(X)$ of $X$ is the smallest affine space in $\reals^N$ that contains $X$.
	If $X = \set{x_0, x_1, \ldots, x_n}$, $\aff(X)$ is given by 
	\begin{align*}
		\aff(X) 
		&= \set{
			\sum_{i=0}^n t_i x_i \in \reals^N 
			:
			\sum_{i=0}^n t_i = 1
		}
		= \set{
			x_0 + \sum_{i=1}^n t_i (x_i - x_0)
			:
			t_1, \ldots, t_n \in \reals
		}
		\\
		&= x_0 + \spanv\set{x_1 - x_0, \ldots, x_n - x_0}
	\end{align*}
	We call $\aff(X)$ the affine hull \textbf{spanned} by $X$.

	\item 
	A set $X = \set{x_0, x_1, \ldots, x_n}$ in $\reals^N$ is called \textbf{affinely independent} if the only subset $Y$ of $X$ such that $\aff(Y) = \aff(X)$ is $X$ itself.
	Equivalently, $X$ is affinely independent if the set $\set{x_1 - x_0, \ldots, x_n - x_0}$ is linearly independent.
	$X$ is called \textbf{affinely dependent} if $X$ is not affinely independent.

	\item 
	Let $A \subseteq \reals^N$ be an affine space such that $A = a+V$ for some $a \in \reals^N$ and subspace $V$ of $\reals^N$ with dimension $n$ (as a vector subspace).
	Then, the \textbf{dimension} of $A$ is $n$, denoted $\dim(A) = n$, and we call $A$ an $n$-\textbf{plane}.
	If $A$ is spanned by an affinely independent set $\set{x_0, \ldots, x_n}$, then 
	$\dim(A) = n$.

	\item 
	The \textbf{convex hull} $\conv(X)$ of a set $X \subseteq \reals^N$ is the smallest convex set in $\reals^N$ containing $X$.
	If $X = \set{x_0, \ldots, x_n}$, then $\conv(X)$ is given as follows:
	\begin{equation*}
		\conv(X) 
		= \conv\bigl(\set{x_0, \ldots, x_n}\bigr) 
		= \set{
			\sum_{i=0}^n t_i x_i \in \reals^N 
			:
			\sum_{i=0}^n t_i = 1 
			\text{ with }
			t_i \geq 0
		}
	\end{equation*}
	Observe that $\conv(X) \subseteq \aff(X)$.
\end{enumerate}
Note that there are various properties listed above that need to be proven, e.g.\ an affine space is translated vector subspace and $\dim(A)$ is well-defined. 
We refer to \cite[Chapter 2]{algtopo:rotman} for such properties.
We provide an example below involving affine spaces and the affine independence/dependence of sets.

\begin{example}
	The sets $A = \set{a,b,c} \subset \mathbb{R}^2$ and $P = \set{p,q,r,w} \subset \mathbb{R}^3$, as illustrated below, are affinely dependent sets. 
	The affine dependence of $A$ can be shown by the non-empty intersection of the affine hull of $\set{a,b}$ (a line) with $\set{c}$ and that of $P$ by the non-empty intersection of the affine hull of $\set{p,q,r}$ (the $yz$-plane) and $\set{w}$.
	\begin{center}
		\includegraphics[width=0.8\linewidth]{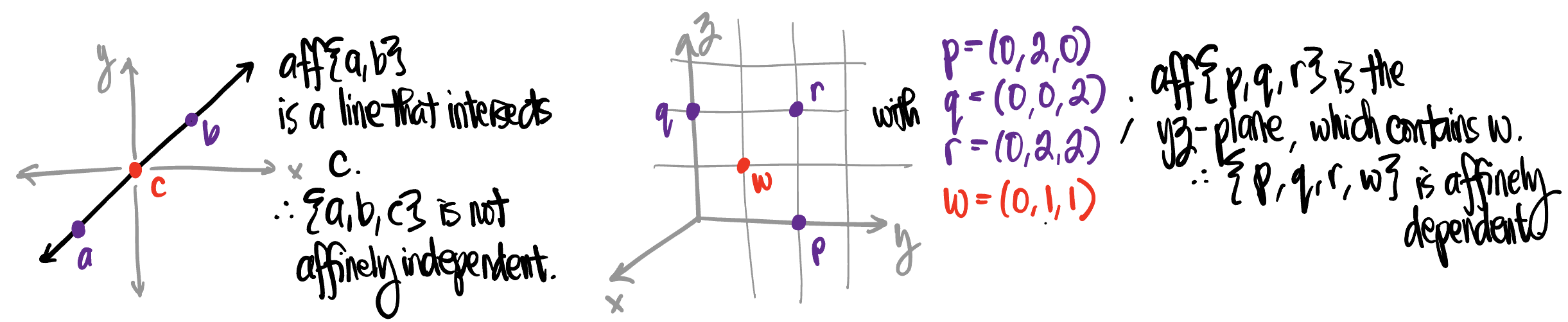}
	\end{center}
	Symbolically, $\aff\set{A} = \aff\set{a,b}$ 
	and $\aff\set{P} = \aff\set{p,q,r}$ despite $\set{A} \neq \set{a,b}$ and $\set{P} \neq \set{p,q,r}$.
	Note that sets $\set{a,b} = A \setminus \set{c}$ and $\set{p,q,r} = P \setminus \set{w}$ are both affinely independent.

\end{example}

Having these definition for affine planes and convex hulls allow a nice characterization of geometric $n$-simplices, as stated below.

\begin{definition}\label{defn:geometric-simplex}
	A \textbf{geometric $n$-simplex} $\sigma = [x_0, \ldots, x_n]$ in $\reals^N$ is the convex hull $\sigma := \conv(x_0, \ldots, x_n)$ of some affinely independent set $X = \set{x_0, \ldots, x_n}$ of points in $\reals^N$. 
	The elements $x_0, \ldots, x_n$ are called the \textbf{vertices} of $\sigma$ and we say that $\sigma$ is \textbf{spanned} by $X$.
	The \textbf{dimension} $\dim(\sigma)$ of $\sigma$ is given by $\dim(\sigma) := n$.
	When $n$ is arbitrary or unambiguous, we may refer to $\sigma$ as a \textbf{geometric simplex}.

	A simplex $\tau$ spanned by some $Y \subseteq X$ is called a \textbf{face} of $\sigma$. When $\dim(\tau) = \dim(\sigma) - 1$, $\tau$ is called a \textbf{facet} of $\sigma$.
	The \textbf{(geometric) boundary} of $\sigma$ is the collection of all its facets.
\end{definition}
\remark{
	Observe that the terminology for geometric simplices is similar to those of (abstract) simplices. This choice is intentional.
}

We provide an example below.

\begin{example}
	Let $X_1 = \set{(0,0)}$, $X_2 = \set{(-1,1), (1, -1)}$, and $X_3 = \set{(0,1), (-1,-1), (1,-1)}$ be collections of points in $\mathbb{R}^2$. All $X_1, X_2, X_3$ are affinely independent. Illustrated below are the simplices spanned by $X_1, X_2, X_3$ respectively with selected points expressed as linear combinations of the corresponding vertex sets.
	\begin{center}
		\includegraphics[width=0.75\linewidth]{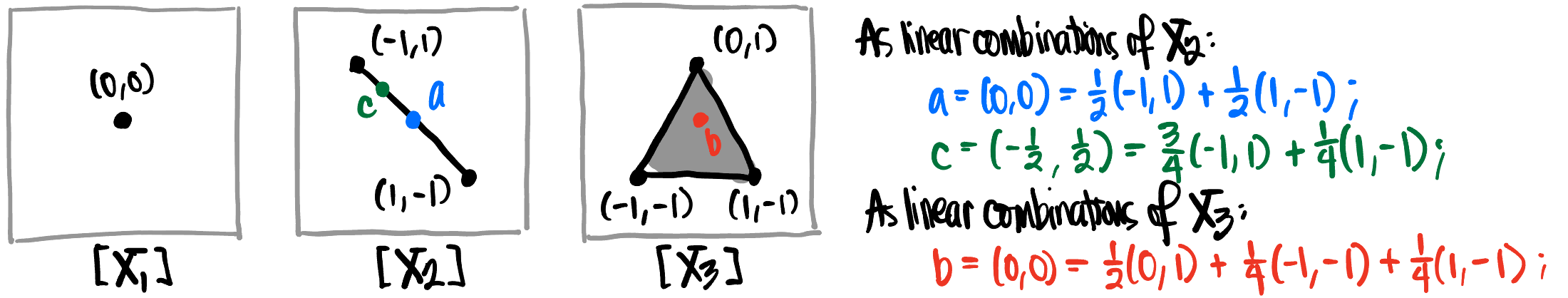}
	\end{center}
\end{example}

As a consequence of the affine independence requirement of the geometric simplices, it can be proven that each geometric $n$-simplex is a homeomorphic copy of every other geometric $n$-simplex, not necessarily living in the same ambient space $\reals^N$. See \cite[Theorem 2.10, Exercise 2.2.6]{algtopo:rotman}.
Intuitively, this means that there is a distinct geometric shape expected of an $n$-simplex for each $n \in \nonnegints$. 
For example:
\begin{enumerate}
	\item Any geometric $0$-simplex is a point.
	\item Any geometric $1$-simplex is a line, homeomorphic to the $1$-disk $D^1 = [-1,1] \subseteq \reals$.
	\item Any geometric $2$-simplex is a filled-in triangle, homeomorphic to the $2$-disk $D^2$ (a circle with its interior).
	\item Any geometric $3$-simplex is a solid tetrahedron, homeomorphic to the $3$-disk $D^3$ (a solid sphere).
\end{enumerate}
Thus, there is usually no reference to the ambient space $\mathbb{R}^N$ and we assume that the ambient dimension $N$ is large enough, similarly to how we refer to other topological constructs like $S^n$ (the $n$-sphere). 
More generally, an $n$-simplex is homeomorphic to the $n$-disk $D^n$ for all $n \in \nonnegints$ (with $D^0$ taken to be point $\reals^0 = \set{0}$).
This relationship also extends to the geometric boundary of $n$-simplices for $n \geq 1$.
For example:
\begin{enumerate}
	\item The boundary of a geometric $1$-simplex $[x_0, x_1]$ is the set $\set{[x_1], [x_2]}$ of $0$-simplices, which is homeomorphic to the $1$-sphere $S^0 = \set{-1,1} \subseteq \reals$.
	Note that the topological boundary of $D^1 = [-1,1]$ as a topological subspace of $\reals$ is given by $\boundary(D^1) = S^0$.

	\item 
	The boundary of a geometric $2$-simplex $[x_0, x_1, x_2]$ is the set 
	$\set{[x_0, x_1], [x_1, x_2], [x_0, x_2]}$ of $1$-simplices.
	Observe that the union of these $1$-simplices form a loop, which is homeomorphic to the $1$-sphere $S^1$, the topological boundary of the $2$-disk $D^2$ as a topological subspace of $\reals^2$,
	as illustrated below:
	\begin{center}
		\includegraphics[width=0.55\linewidth]{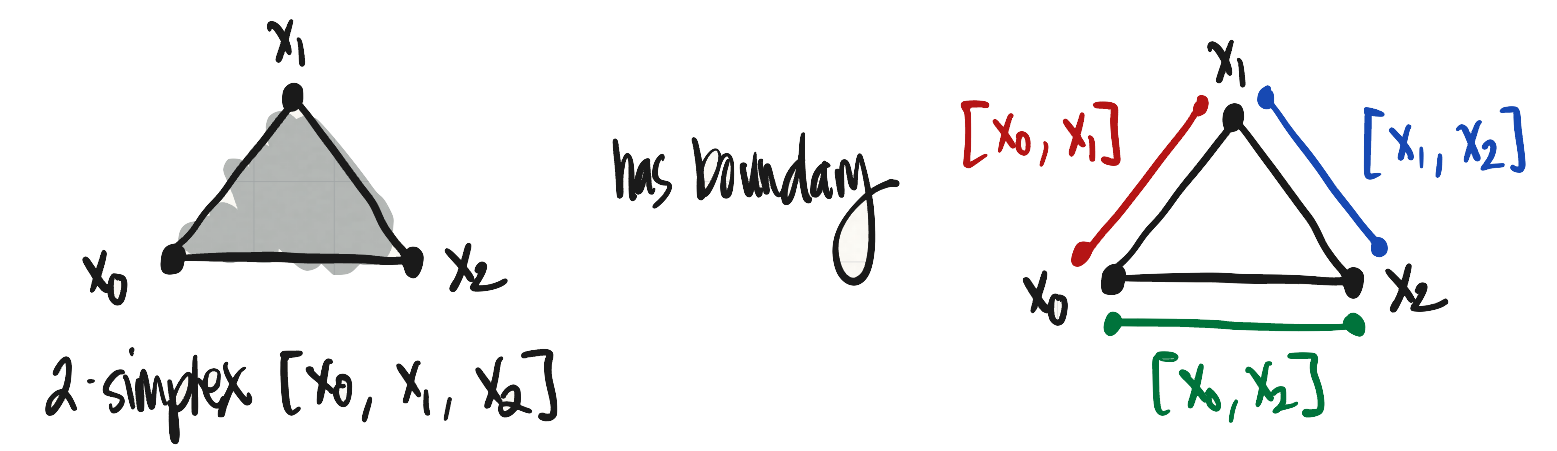}
	\end{center}
\end{enumerate}
That is, the union of the boundary of a geometric $n$-simplex is homeomorphic to the topological boundary $\boundary{D^n} = S^{n-1}$ of the $n$-disk as a subset of $\reals^n$
where $S^n$ refers to the $n$-sphere $S^n = \set{x \in \reals^{n+1} : |x| = 1}$.
Once we have established the relationship between geometric simplicial complexes and abstract simplicial complexes, this notion of boundary extends to the simplices of abstract simplicial complexes as well.
We see this later in \fref{section:simplicial-homology} in the context of simplicial homology.

A geometric simplicial complex is then a collection of these geometric simplices with certain properties. 
We state this in more detail below.

\begin{definition}\label{defn:simp-complex-geometric}
	A \textbf{geometric simplicial complex} $K$ is a collection of simplices in $\reals^N$ such that 
	\begin{enumerate}
		\item 
		If $\sigma \in K$, every face of $\sigma$ also belongs in $K$.

		\item 
		If $\sigma, \tau \in K$, then $\sigma \cap \tau$ is either empty or a common face of $\sigma$ and of $\tau$.
	\end{enumerate}
	We write $\Vertex(K)$ to denote the \textbf{vertex set} of $K$ given by the union of all geometric $0$-simplices of $K$.
	The Euclidean space $\reals^N$ is called the \textbf{ambient space} of $K$.
	A simplicial complex $L$ is called a \textbf{subcomplex} of $K$ if $\Vertex(L) \subseteq \Vertex(K)$.
	If $K$ is a finite set, then $K$ is called a \textbf{finite geometric simplicial complex}.

	The \textbf{underlying space} or \textbf{polytope} $|K|$ of $K$ is the topological subspace of $\reals^N$ given by $|K| := \bigcup_{\sigma \in K}\sigma$, i.e.\ the union of all simplices of $K$.
	A topological space $X$ is called a \textbf{polyhedron} if there exists a geometric simplicial complex $K$ and a homeomorphism $h: |K| \to X$.
	In this case, the pair $(K,h)$ is called a \textbf{triangulation} of $X$.
\end{definition}
\remark{
	In this expository paper, we somewhat abuse notation for triangulations of topological spaces by not specifying the homeomorphism $h: |K| \to X$, i.e.\ we often say that $K$ is a triangulation of $X$.
}

We want to emphasize that geometric simplicial complexes are a specific type of representation of topological spaces and that both conditions, as given in \fref{defn:simp-complex-geometric}, have to be satisfied.
In particular, an arbitrary collection of simplices $\set{\sigma_i}_{i \in I}$ does not make a simplicial complex.
We provide some examples below.

\begin{example}
	The first two collections of simplices listed below correspond to geometric simplicial complexes but the third does not.
	\begin{enumerate}
		\item 
		The solid tetrahedron in $\mathbb{R}^3$ as a geometric simplicial complex.
		\begin{center}
			\includegraphics[width=0.6\linewidth]{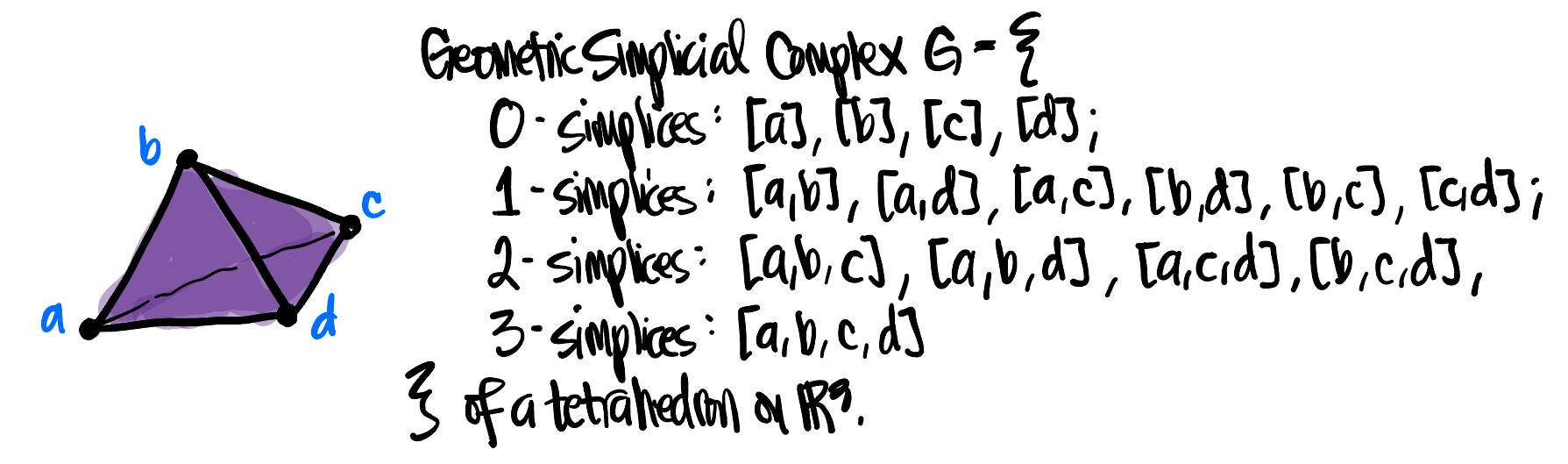}
		\end{center}
		Note that $a,b,c,d$ are points in $\reals^3$ with unidentified coordinates and are not indeterminates.

		\item 
		A geometric simplicial complex in $\reals^2$.
		\begin{center}
			\includegraphics[width=0.6\linewidth]{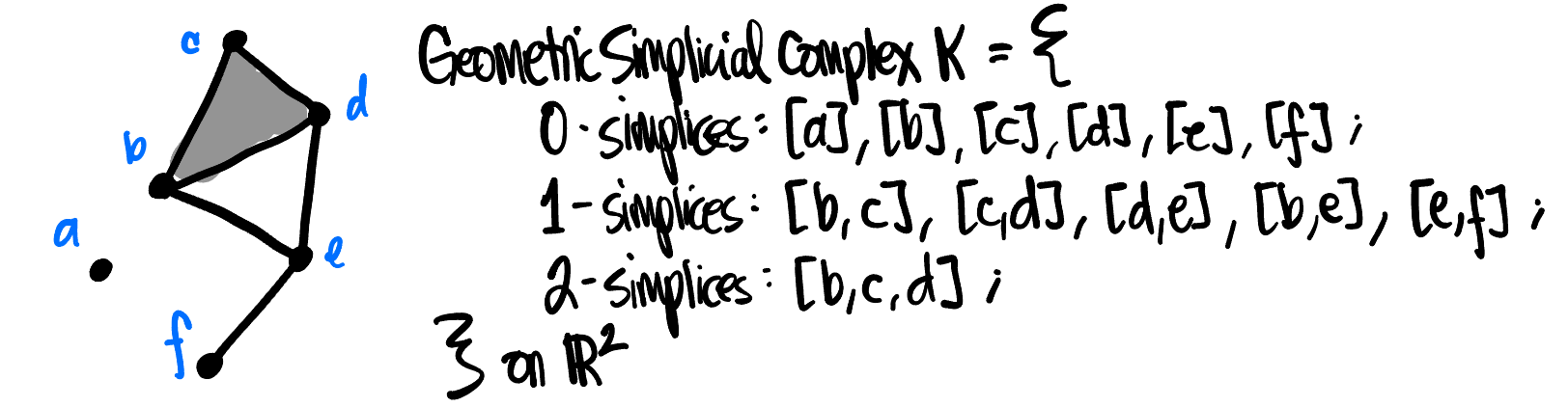}
		\end{center}
		As with the previous item, $a,b,c,d,e,f$ are points in $\reals^2$.

		\item 
		The collection $H$, as illustrated below, is not a geometric simplicial complex since it violates condition (ii) of \fref{defn:simp-complex-geometric}. 
		One such pair that violates condition (ii) is the pair $[a,c]$ and $[c,d]$ since their intersection, $[c,d]$, is not a face of $[a,c]$. 
		\begin{center}
			\includegraphics[width=0.6\linewidth]{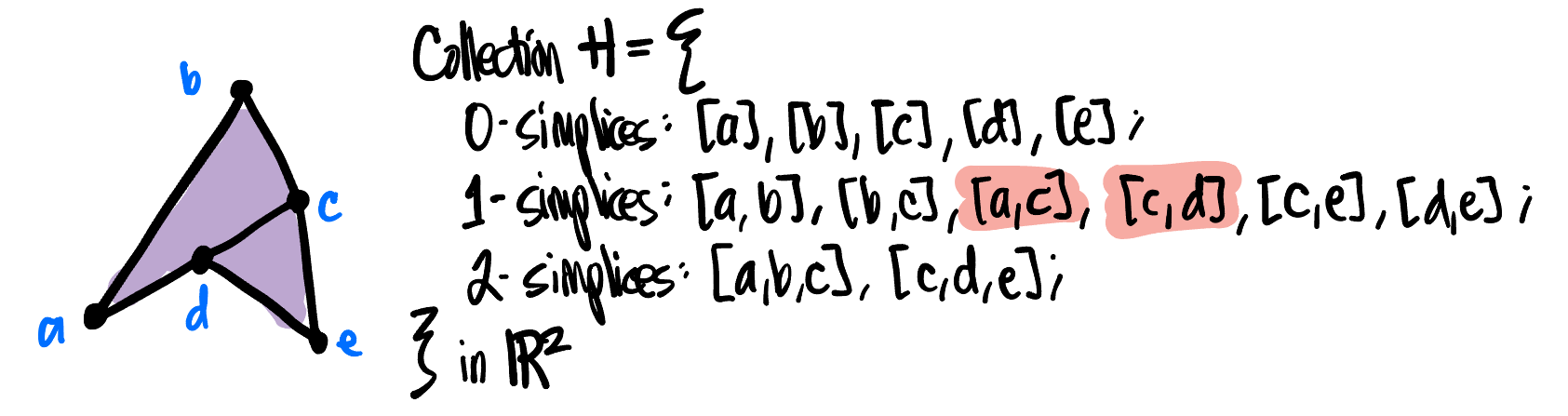}
		\end{center}
		Observe that simply adding the simplex $[a,b,c,d]$ to $H$ does not fix the issue. 
		The set $\set{a,b,c,d}$ is not an affinely independent set and, therefore, cannot generate a simplex. That is, the simplex $[a,b,c,d]$ is not defined.
	\end{enumerate}
\end{example}

\noindent
Abstract simplicial complexes and geometric simplicial complexes are related using the following concept:

\begin{definition}
	The \textbf{vertex scheme} $K_V$ of a geometric simplicial complex $K_G$ is the abstract simplicial complex given by $K_V := \set{ A \subseteq V : A \text{ spans some simplex in } K_G }$ with $\Vertex(K_V) = \Vertex(K_G)$ (as sets).
	Given a simplicial complex $K_A$, we say that geometric simplicial complex $K_G$ is a \textbf{geometric realization} of $K_A$ if the vertex scheme $K_V$ of $K_G$ is isomorphic to $K_A$.
\end{definition}

Observe that, while we can always generate an abstract simplicial complex from a geometric simplicial complex, the converse may not be true.
However, if the abstract simplicial complex in question has finite dimension, there always exist a corresponding geometric simplicial complex.
We state this as a theorem below.

\begin{theorem}\label{thm:simplicial-complex-uniqueness}
	\textbf{Existence and Uniqueness (up to Homeomorphism) of Geometric Realizations}

	Every abstract simplicial complex $K$ with dimension $\dim(K) = d$ has a geometric realization in $\reals^{2d+1}$.
	Furthermore, the polytopes of the geometric realizations of $K$ are unique up to homeomorphism.
\end{theorem}
\remark{
	We refer to \cite[Theorem 1.6.1]{algtopo:matousek} for proof of the existence claim and 
	\cite[Proposition 1.5.4]{algtopo:matousek} for the claim of uniqueness up to homeomorphism.

	The proof of the existence claim relies on constructing a geometric simplicial complex $K_G$ using the \textit{moment curve}.
	Label the vertices of $K$ by $\Vertex(K) = \set{v_i}_{i \in I}$.
	By assumption, $\dim(\sigma) \leq d$ and $\card(\sigma) \leq d+1$ for all simplices $\sigma \in K$.
	We can determine a lower bound on the dimension $N$ of the ambient space $\reals^N$ needed for $K_G$ by considering an upper bound on $\card(\sigma \cup \tau)$ 
	with $\sigma$ and $\tau$ distinct simplices of $K$:
	\begin{equation*}
		\card(\sigma \cup \tau) \leq \card(\sigma) + \card(\tau) \leq (d+1) + (d+1) = 2d+2
	\end{equation*}
	We define $N = 2d+1$ and consider the moment curve $\gamma(t): \reals \to \reals^{2d+1}$ given by $\gamma(t) = (t, t^2, \ldots, t^{2d+1})$.
	The moment curve has the property that for any set $G \subseteq \gamma(\reals)$ with $\card(G) \geq 2d+2$, any subset of $G$ of cardinality $2d+2$ is affinely independent (i.e.\ $G$ is in general position).
	We construct a vertex set $X = \set{x_i}_{i \in I}$ of $K_G$ by assigning a unique $t_i \in \reals$ for each vertex $v_i \in \Vertex(K)$ and defining $x_i := \gamma(t_i)$.
	Then, we define $K_G$ by the condition: $\set{v_{i_0}, \ldots, v_{i_n}}$ is an $n$-simplex of $K$ if and only if $[x_{i_0}, \ldots, x_{i_n}]$ is a geometric $n$-simplex of $K_G$.
}

Observe that this result tells us that every abstract simplicial complex of finite dimension determines a unique topological space by a geometric simplicial complex and that, assuming the dimension is finite, we can go back and forth between representations as convenient.
We provide an example below.

\begin{example}
	The figure below illustrates a number of representations of the $2$-disk $D^2$.
	\begin{center}
		\includegraphics[width=0.95\linewidth]{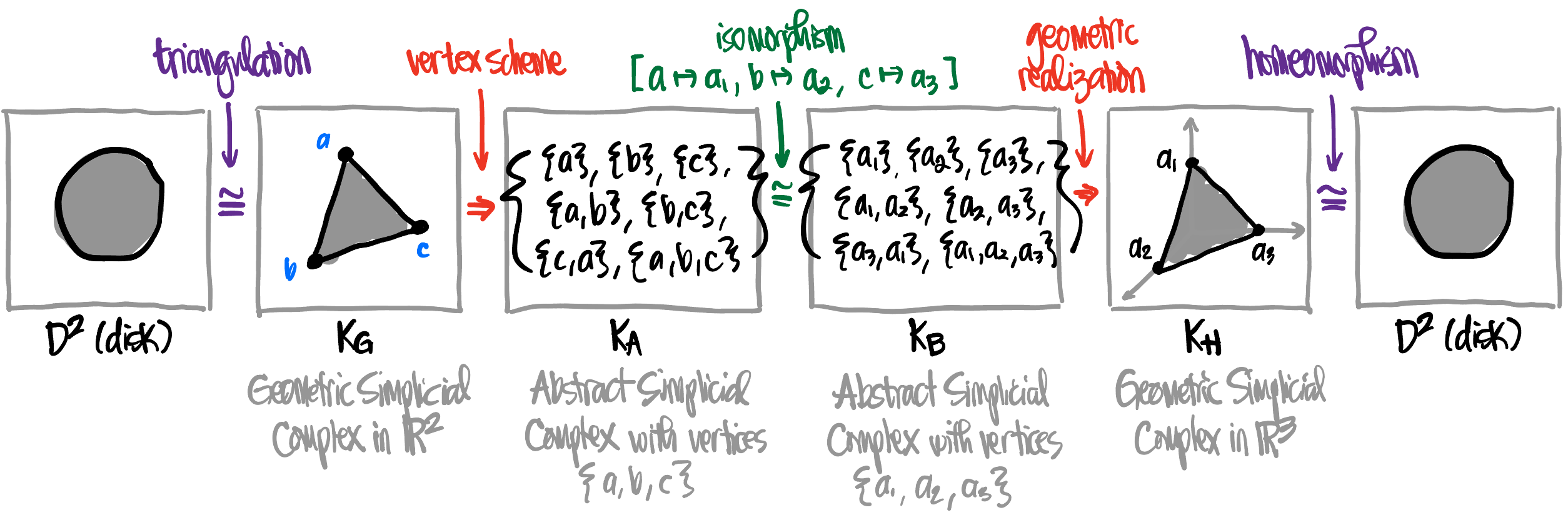}
	\end{center}
	The correspondence between these representations are identified below:
	\begin{enumerate}
		\item 
		The geometric simplicial complex $K_G$ in $\reals^2$ with $\Vertex(K_G) = \set{a,b,c}$ is a triangulation of $D^2$ with vertex scheme given by the abstract simplicial complex $K_A$.

		\item 
		We relate $K_A$ to another abstract simplicial complex $K_B$ with $\Vertex(K_B) = \set{a_1, a_2, a_3}$ using a simplicial isomorphism.
		We can interpret this to be renaming of the vertices $a,b,c$ of $K_A$ to $a_1, a_2, a_3$ of $K_B$ respectively.

		\item 
		$K_B$ has a geometric realization given by the geometric simplicial complex $K_H$ in $\reals^3$ with $\Vertex(K_H) = \set{e_1, e_2, e_3}$ (the standard basis vectors on $\reals^3$).
		Observe that $K_H$ is a triangulation of $D^2$ either by inspection of $K_H$ or by application of the uniqueness result of \fref{thm:simplicial-complex-uniqueness} on $|K_G| \cong D^2$.
	\end{enumerate}
	Observe that by \fref{thm:simplicial-complex-uniqueness}, $K_A$ with $d=2$ is guaranteed to have a geometric realization in $\mathbb{R}^{2d+1} = \mathbb{R}^5$. 
	Note that the theorem does not state that $2d+1$ is the minimal dimension and, therefore, this example does not contradict said theorem.
\end{example}

	\clearpage


\section{Simplicial Homology with Coefficients in a PID}
\label{section:simplicial-homology}

In this section, we discuss the homology of simplicial complexes with coefficients in a PID $R$.
Note that the examples presented in this paper mostly consider cases wherein the coefficient ring is $R = \ints$ (in this case, the homology group is an abelian group) or $R = \field$ is a field (in this case, the homology group is a $\field$-vector space).


The construction of simplicial homology involves the creation of a chain complex of $R$-modules from an abstract simplicial complex.
We use the following definition for chain complexes, taken from \cite[p317]{algtopo:rotman}.

\begin{blockquote}
	A \textit{chain complex} $C_\ast = (C_n, \boundary_n)_{n \in \ints}$ of $R$-modules is an $\ints$-indexed collection of $R$-modules $C_n$ and $R$-module homomorphisms $\boundary_n: C_n \to C_{n-1}$ such that for all $n \in \ints$, $\boundary_n \circ\ \boundary_{n+1} = 0$.
		The $R$-modules $C_n$ are generally called \textit{chain groups} and the homomorphisms $\boundary_n$ are called \textit{differentials} or \textit{connecting homomorphisms}.
\end{blockquote}

\noindent 
We follow the convention in \cite{algtopo:rotman} and use an asterisk $(\ast)$ as the ``placeholder'' for the index of a chain complex.
	In contrast, other references, e.g.\ \cite{cattheory:rhiel,cattheory:rotman}, use a bullet $(\bullet)$ and write $C_\bullet = (C_n, \boundary_n)_{n \in \ints}$.
	We reserve the use of bullets for persistence modules, introduced later in \fref{chapter:persistence-theory}.
	
We start with the definition for the chain groups in simplicial homology, adapted from
\cite[Chapter 7]{algtopo:rotman}.

\begin{definition}\label{defn:chain-groups}
	Let $K$ be a simplicial complex and $R$ a PID.
	For each $n \in \ints$ with $n \geq 0$,
		define the $n$\th \textbf{simplicial chain group $C_n(K; R)$ of $K$ with coefficients in $R$}
	to be the $R$-module with the following presentation:
	\begin{enumerate}[labelwidth=2in, leftmargin=1in]
		\item[Generators:]
			All $(n+1)$-tuples $(v_0, v_1, \ldots, v_n)$ such that $v_i \in \Vertex(K)$ and $\set{v_0, \ldots, v_n}$ is an $n$-simplex of $K$.
		\item[Relations:]
			For $n=0$, none.
			For $n \geq 1$,
			$(v_0, v_1, \ldots, v_n) - (\sgn\pi)(
				v_{\pi(0)}, v_{\pi(1)}, \ldots, v_{\pi(n)}
			)$ for all permutations $\pi: [n] \to [n]$ with 
			$[n] = \set{0, 1, \ldots, n}$ where  
			$\sgn\pi = \pm 1$ refers to the parity of $\pi$. 
	\end{enumerate}
	For $n \in \ints$ with $n \leq -1$, define $C_n(K;R)$ to be the trivial $R$-module.
	An $n$-\textbf{chain} is an element of $C_n(K;R)$.
	If $R = \ints$,
		we often write $C_n(K) := C_n(K;\ints)$ and call $C_n(K)$ the $n$\th \textbf{simplicial chain group of $K$}, i.e.\ without reference to the ring $\ints$.
\end{definition}
\remarks{
	\item 
	For brevity, we usually use the term ``chain group'' when referring to a simplicial chain group if the relation to a simplicial complex is clear from context, e.g.\ we say $C_n(K;R)$ is a chain group of $K$ if it has been stated that $K$ is a simplicial complex. 

	\item 
	The $(n+1)$-tuple $(v_0, v_1, \ldots, v_n)$ described above can be seen as a total order (alternatively, linear order) on the $n$-simplex $\sigma := \set{v_0, v_1, \ldots, v_n}$.
	In the statement above, it just so happens that the indexing $i$ on the vertices $v_i$ seems to follow the total order on $[n] = \set{0, \ldots, n} \subseteq \nonnegints$.
	Note that there cannot be repeated vertices on $(v_0, \ldots, v_n)$. Otherwise, $\card(\sigma) < n$ and $\sigma$ cannot be an $n$-simplex.
}

In the context of simplicial homology, the index $n \in \ints$ is usually called the \textit{dimension} since each $n$-simplex corresponds to an $n$-dimensional object.
We follow his convention throughout this paper to distinguish $n$ from other indices.
Note that the dimension being possibly negative for simplicial chain groups is not an issue since $C_n(K;R)$ for $n < 0$ is trivial.
In fact, some references like \cite{algtopo:hatcher,algtopo:munkres} do not define such groups since chain complexes in said references are not typically discussed outside a topological perspective.

Observe that \fref{defn:chain-groups} defines simplicial chain groups as quotient modules and that the elements of $C_n(K;R)$ are cosets of formal sums of orderings of $n$-simplices of $K$.
Relative to the conventional coset notation, elements of $C_n(K;R)$ should be written as sums of terms of the form $\bigl[(v_0, \ldots, v_n)\bigr]$ where $(v_0, \ldots, v_n)$ is a specific ordering of the $n$-simplex $\set{v_0, \ldots, v_n}$.
To avoid the abundance of grouping symbols and for brevity, we introduce alternative notation below.

\begin{bigremark}\label{remark:shorthand-for-oriented-simplices}
	For each $n$-simplex $\sigma = \set{v_0, v_1, \ldots, v_n}$ of some simplicial complex $K$,
	define 
	the string $v_0 v_1 \cdots v_n$ 
	and the symbol
	$[v_0, \ldots, v_n]$
	to refer to the coset $\bigl[(v_0, \ldots, v_n)\bigr] \in C_n(K;R)$,
	specifically with coset representative $(v_0, \ldots, v_n)$.
\end{bigremark}

Note that we prefer using the string notation outside this section since we will be defining simplicial homology groups later.
That is, we prefer writing $v_0 v_1 \cdots v_n$ for elements of $C_n(K;R)$ so that we write $[v_0 v_1 \cdots v_n]$, as opposed to $\bigl[ [v_0, v_1, \ldots, v_n] \bigr]$, for the elements of the simplicial homology group $H_n(K;R)$ (defined later in \fref{defn:simplicial-chain-complex}).
The notation $[v_0, \ldots, v_n]$ seems to be taken from that of geometric simplices, as stated in \fref{defn:geometric-simplex}.

Next, we state a result about simplicial chain groups being free modules, taken from \cite[Lemma 7.10]{algtopo:rotman}.

\begin{proposition}\label{prop:chain-groups-are-free}
	Let $K$ be a simplicial complex and $R$ a PID.
	For all $n \in \ints$, $C_n(K;R)$ is a free $R$-module.
\end{proposition}
\begin{proof}
	For $n \leq -1$, $C_n(K;R)$ is trivial and therefore free.
	Assume $n \geq 0$.
	Let $G_n$ be the free $R$-module generated by $(n+1)$-tuples $(v_0, \ldots, v_n)$ such that $\set{v_0, \ldots, v_n}$ is an $n$-simplex of $K$.
	Let $S_n$ be the (free) submodule of $G_n$ generated by elements of the form
	$
		(v_0, \ldots, v_n) - (\sgn\pi)(v_{\pi(0)}, \ldots v_{\pi(n)})
	$
	with $\set{v_0, \ldots, v_n}$ an $n$-simplex of $K$ and $\pi: [n] \to [n]$ a permutation. 
	By definition, $C_n(K;R) = G_n \bigmod S_n$.

	For each $n$-simplex $\sigma = \set{v_0, \ldots, v_n}$ of $K$:
	define $G_{n,\sigma}$ to be the free $R$-module generated by all $(n+1)$-tuples of $\sigma$ and 
		$S_{n,\sigma}$ to be the (free) submodule of $G_{n,\sigma}$ generated by elements of the form 
	$
		(v_0, \ldots, v_n) - (\sgn\pi)(v_{\pi(0)}, \ldots v_{\pi(n)})
	$. Then, 
	\begin{equation*}
		G_n = \bigoplus_{\mathclap{\substack{\sigma \in K \\ \dim(\sigma) = n}}}
			G_{n,\sigma} 
		\quad,\quad 
		S_n = \bigoplus_{\mathclap{\substack{\sigma \in K \\ \dim(\sigma) = n}}}
			S_{n,\sigma} 
		\quad\text{ and }\quad
		C_n(K;R) = \frac{G_n}{S_n} =
		\bigoplus_{\mathclap{\substack{\sigma \in K \\ \dim(\sigma) = n}}}
			\left( \frac{G_{n,\sigma}}{S_{n,\sigma}} \right)
	\end{equation*}
	Therefore, it suffices to check that each $(G_{n,\sigma} \bigmod S_{n,\sigma})$ is a free $R$-module.

	Let $\sigma = \set{v_0, \ldots, v_n}$ be an $n$-simplex of $K$.
	Note that by labelling the vertices of $\sigma$, we are implicitly imposing a total order on $\sigma$.
	Let $(\sigma) := (v_0, \ldots, v_n)$.
	Let $\pi: [n] \to [n]$ be a permutation and let $(\sigma_\pi) \in G_{n,\sigma}$ refer to the tuple $(\sigma_\pi) := (v_{\pi(0)}, \ldots, v_{\pi(n)})$.
	If $\pi$ is an identity permutation, 
	then $(\sigma) - (\sgn\pi)(\sigma_\pi) = (\sigma) - (1)(\sigma) = 0$.
	Otherwise, $(\sigma) - (\sgn\pi)(\sigma_\pi) \neq 0$.

	Let $B_{n,\sigma}$ be a set of elements in $G_{n,\sigma}$ such that $(\sigma) - (\sgn\pi)(\sigma_\pi) \in B_{n,\sigma}$ if and only if $\pi:[n] \to [n]$ is not the identity permutation on $[n]$.
	Then, $B_{n,\sigma}$ is a basis of $S_{n,\sigma}$ 
		and $\set{(\sigma)} \cup B_{n,\sigma}$ is a basis of $G_{n,\sigma}$.
	Therefore, the quotient $G_{n,\sigma} \bigmod S_{n,\sigma}$ is a free $R$-module with basis $\set{[\sigma]}$ with $[\sigma] := (\sigma) + S_{n,\sigma}$.

	Therefore, $C_n(K;R)$ is a free $R$-module.
	A basis $B_n$ of $C_n(K;R)$ can be generated from the union of the bases $\set{[\sigma]}$ for each summand $G_{n,\sigma} \bigmod S_{n,\sigma}$.
\end{proof}
\remark{
	In the proof above, we consider the direct sums to correspond to \textit{internal direct sums}, wherein the elements of the direct sum need not be tuples. For example, $R\ket{a} \oplus R\ket{b} = \set{r_1a + r_2b : r_1,r_2 \in R}$ (internal), as opposed to 
	$R\ket{a} \oplus R\ket{b} = \set{(r_1a, r_2b) : r_1, r_2 \in R}$ (external).
	However, since these two notions uniquely determine an $R$-module up to $R$-module isomorphism, these are used somewhat interchangeably. 
}

We want to emphasize that in most texts including \cite{algtopo:hatcher} and \cite{algtopo:rotman},
	simplicial chain groups are typically defined on 
	\textit{oriented} simplicial complexes.
Below, we provide definitions for these notions, relative to our definition of chain groups as stated in \fref{defn:chain-groups}. 

\begin{definition}\label{defn:simplicial-complex-orientation}
	Given an $n$-simplex $\sigma = \set{v_0, \ldots, v_n}$ of $K$,
	the coset $[v_0, \ldots, v_n] := \big[(v_0, \ldots, v_n)\big] \in C_n(K;R)$ is called an \textbf{oriented} $n$-\textbf{simplex}.
	An \textbf{orientation} on a simplicial complex $K$ is a total order  on the vertex set $\Vertex(K)$ of $K$.
	An \textbf{oriented simplicial complex} is a simplicial complex equipped with an orientation.
\end{definition}
\remark{
	Recall that a total order $\leq$ on a finite set  $A = \set{a_0, \ldots, a_n}$ corresponds to a linear ordering of the same set, 
	i.e.\ an ordered sequence $(a_0, a_1, \ldots, a_n)$ of elements of $A$ defines a unique total order $a_0 < a_1 < \cdots < a_n$ of $A$.
	Following that convention, 
		we usually equip a orientation to a simplicial complex $K$ by 
		using a tuple $(v_0, \ldots, v_N)$ corresponding to a total order on $\Vertex(K) = \set{v_0, \ldots, v_N}$.
	Note that if the vertices of $K$ are indexed with a totally ordered set (e.g.\ some subset of $\ints$), the total order on $\Vertex(K)$ is often defined to correspond to the total order on said indexing set.
}

Observe that \fref{defn:chain-groups} does not require that the simplicial complex be equipped with an orientation beforehand.
Here, we interpret orientation on the simplices of $K$ (not of the simplicial complex $K$) to be more of a consequence of the definition of chain groups.
In particular, the relation $(v_0, \ldots, v_n) - (\sgn\pi)(v_{\pi(0)}, \ldots, v_{\pi(n)}) \in C_n(K;R)$ corresponding to an $n$-simplex $\sigma = \set{v_0, \ldots, v_n}$ of $K$ with $n \geq 1$ in \fref{defn:chain-groups}
corresponds to the following equivalence relation $\sim$ on the set of all total orders (written as tuples) on $\sigma = \set{v_0, \ldots, v_n}$:
\begin{equation*}
\begin{array}{r c l c l}
	(v_0, \ldots, v_n) &\sim& \phantom{-}(v_{\pi(0)}, \ldots, v_{\pi(n)})
	&\text{ if }& \pi:[n] \to [n] \text{ is an even permutation } 
	\\
	(v_0, \ldots, v_n) &\sim& -(v_{\pi(0)}, \ldots, v_{\pi(n)})
	&\text{ if }& \pi:[n] \to [n] \text{ is an odd permutation } 
\end{array}
\end{equation*}
Note that parity (i.e.\ $\sgn\pi$) is not defined for permutations $\pi: \set{0} \to \set{0}$ on one element.
The relation $\sim$ produces two equivalence classes for each $n$-simplex $\sigma$ of $K$ with $n \geq 1$, which are often interpreted to correspond to \textit{orientations} of the $n$-simplex $\sigma$ as a \textit{geometric} object (hence the name).
We list some of these interpretations below:
\begin{enumerate}[left=0.7in]
	\item[For $n=1$:]
	The geometric simplex corresponding to a $1$-simplex $\set{a,b}$ is a line segment.
	Then, an orientation on $\set{a,b}$ is interpreted to be a choice in direction of said line segment, as illustrated below:
	\begin{center}
		\includegraphics[width=0.4\linewidth]{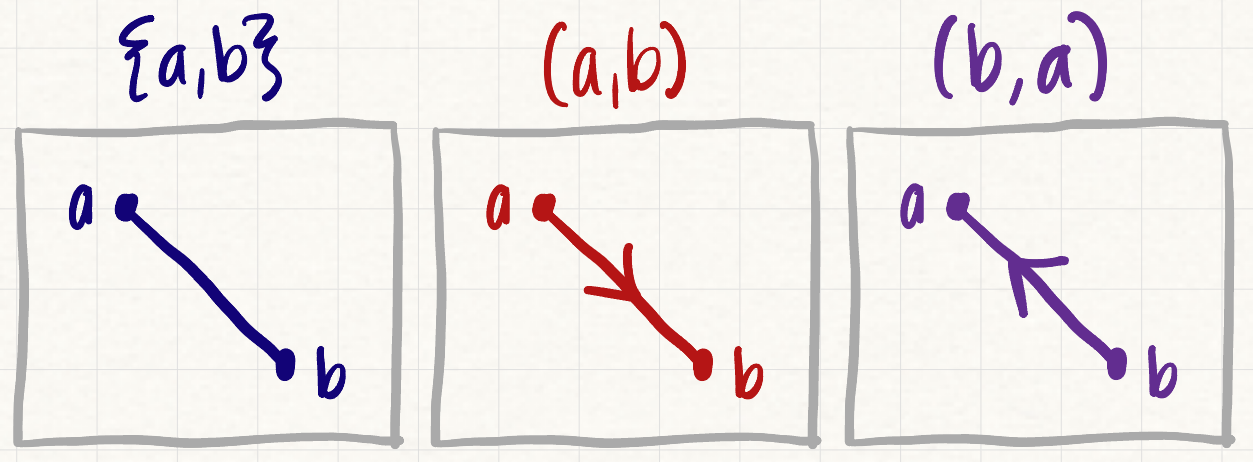}
	\end{center}

	\item[For $n=2$:]
	The geometric simplex corresponding to a $2$-simplex $\set{a,b,c}$ is a triangle (more specifically, the convex hull of a triangle).
	An orientation on this simplex is often interpreted to be a choice in the direction of rotation about an axis normal to the affine 2-plane spanned by $\set{a,b,c}$.
	The $6=3!$ possible total orders on $\set{a,b,c}$ partition into two orientations, as illustrated below:
	\begin{center}
		\includegraphics[width=0.55\linewidth]{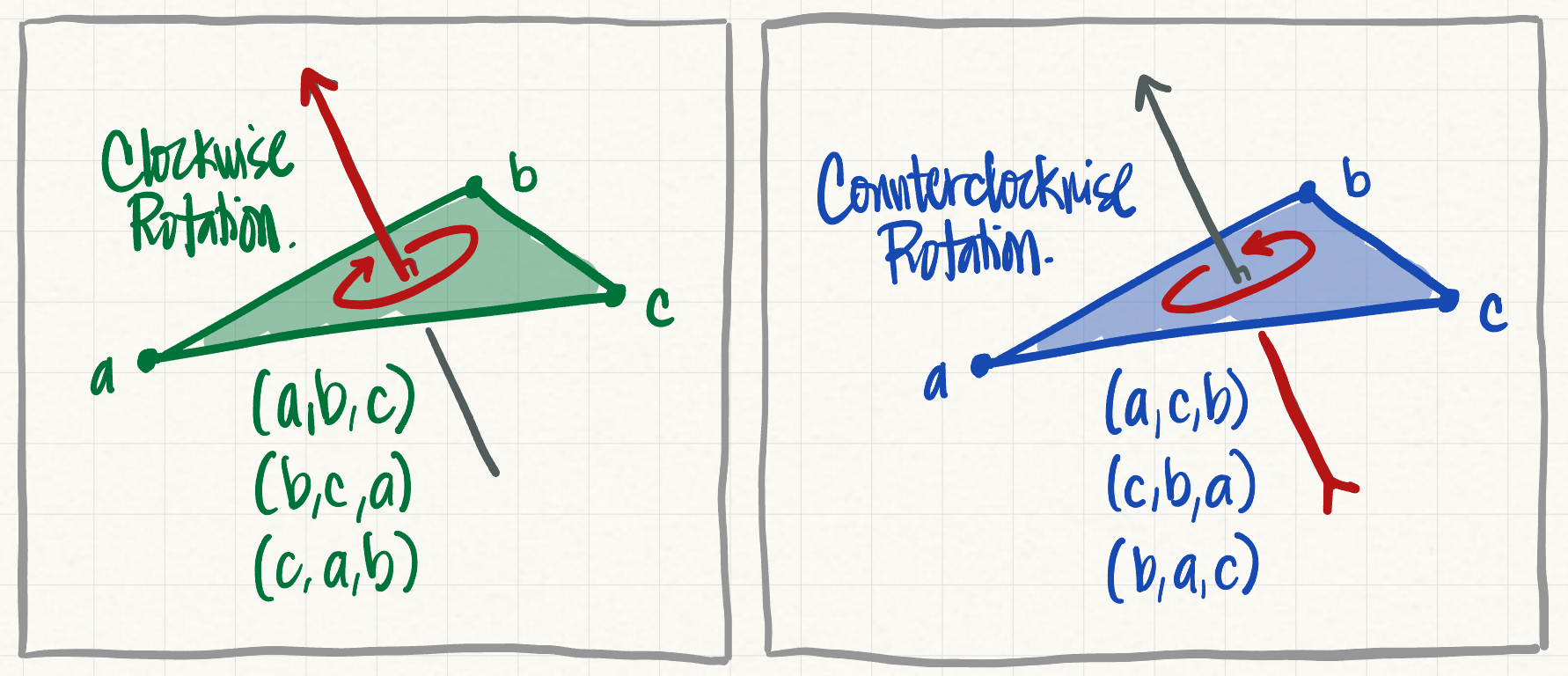}
	\end{center}
	Note that which rotation is considered ``counterclockwise'' is determined by the right hand rule on the \textit{chosen} direction of the normal vector.
	This rotation can also be interpreted as a choice of direction of the vector normal to the place, denoted in the illustration above by the red part of the axis. 
\end{enumerate}
Then, multiplication of $[v_0, \ldots, v_n]$ by $-1 \in R$ is interpreted to be a reversal in orientation.
Unfortunately, this interpretation loses effectiveness in higher dimensions ($n \geq 3$) as it becomes harder to describe or visualize the orientation of $n$-dimensional objects.
It does, however, provide a justification for the definition of boundary maps.

Instead, we interpret an orientation on $K$ to correspond to a choice of basis on each of the chain groups $C_n(K;R)$ of $K$.
Since each $n$-simplex $\sigma = \set{v_0, \ldots, v_n}$ of a simplicial complex $K$ is a subset of the vertex set $\Vertex(K)$ of $K$,
	an orientation on $K$ imposes a total order on $\sigma$ by restriction.
Consequently, an orientation on a simplicial complex $K$ uniquely identifies an $(n+1)$-tuple for each $n$-simplex of $K$.
This motivates the following definition:

\begin{definition}\label{defn:standard-basis-on-chain-groups}
	Let $K$ be a simplicial complex and $R$ a PID.
	Let $\leq'$ be a total order on $\Vertex(K)$, i.e.\ an orientation on $K$.
	For each $n \in \ints$ with $n \geq 0$,
		let the \textbf{standard ordered basis $\basis{K}_n$} of $C_n(K;R)$
		\textbf{induced by the orientation on $K$} be the following set, ordered lexicographically with respect to $\leq'$:
	\begin{equation*}
		\basis{K}_n := \biggl\{
			[v_0, \ldots, v_n] 
			: 
			\set{v_0, \ldots, v_n} \text{ is an $n$-simplex of $K$ and }
			v_0 \leq' \cdots \leq' v_n
		\biggr\}
	\end{equation*}
	with $[v_0, \ldots, v_n]$ as described in \fref{remark:shorthand-for-oriented-simplices}, i.e.\ 
	the symbol $[v_0, \ldots, v_n]$ and the string $v_0 v_1 \cdots v_n$ both refer to the coset in $C_n(K;R)$ specifically with coset representative $(v_0, \ldots, v_n)$.
\end{definition}
\remark{
	Observe that there is an immediate one-to-one correspondence between the elements of $\basis{K}_n$ and the $n$-simplices of $K$.
	Earlier in the proof of \fref{prop:chain-groups-are-free},
		$C_n(K;R)$ is shown to be a direct sum of free $R$-modules 
		$G_{n,\sigma} \bigmod S_{n,\sigma}$,
		one for each $n$-simplex $\sigma$ of $K$ with basis determined by some arbitrarily chosen total order (expressed as a tuple) on $\sigma$.
	$\basis{K}_n$ simply specifies which total order of $\sigma$ is chosen for each summand $G_{n,\sigma} \bigmod S_{n,\sigma}$.
	That is, the claim of $\basis{K}_n$ being a basis of $C_n(K;R)$ is given by the same proof as \fref{prop:chain-groups-are-free}.
}

Note that the lexicographic order on $\basis{K}_n$ is not entirely critical for the construction of chain groups.
For example, \cite{algtopo:rotman} defines an orientation on $K$ more generally to be a partial order on $K$ that restricts to a total order on $\sigma$ for each simplex $\sigma \in K$.
In this case, we are not guaranteed that every vertex $v_k$ and $w_k$ of $K$ are comparable and lexicographic order may not be well-defined.
Then, $\basis{K}_n$ can only be considered a \textit{standard basis}, as opposed to a standard \textit{ordered} basis.
However, since we will be calculating simplicial homology using matrices in \fref{chapter:matrix-calculation}, having a ``default'' order on $\basis{K}_n$ becomes convenient for exposition.

Since we have defined an orientation on a simplicial complex $K$ to be a total order $\leq'$ on $\Vertex(K)$,
	a lexicographic order on $\basis{K}_n$ is well-defined.
More specifically, we refer to a lexicographic order on the tuples $(v_0, \ldots, v_n)$ of the coset representative of $C_n(K;R)$.
Given two $n$-simplices $\set{v_0, \ldots, v_n}$ and $\set{w_0, \ldots, w_n}$ of $K$:
\begin{equation*}
	[v_0, \ldots, v_n] \leq [w_0, \ldots, w_n] 
	\qquad\text{ if and only if }\qquad
	\left\{\begin{array}{c}
		\text{ there exists } k \in [n] \text{ such that } \\
		v_i = w_i \text{ for all } i = 0, \ldots, k-1 \\
		\text{ and } v_k \leq' w_k
	\end{array}\right\}
\end{equation*}
In practice, this means that we compare the $k$\th elements of $(v_0, \ldots, v_n)$ and of $(w_0, \ldots, w_n)$ in increasing index $k \in \set{0, \ldots, n}$ until the vertices $v_k$ and $w_k$ are different. Then, whichever coset comes first in the order of $\basis{K}_n$ is determined by whether $v_k$ or $w_k$ comes first in the order of $\Vertex(K)$ (i.e.\ the orientation on $K$).

Below, we provide an example wherein we identify the chain groups of a simplicial complex $K$ along with the standard bases on said chain groups.

\begin{example}
	Let $K$ be a simplicial complex on the full simplex of $V = \set{a,b,c,d}$, i.e.\ the $n$-simplices of $K$ are exactly the subsets of $V$ of cardinality $n+1$ and $K$ has the geometric realization of a solid tetrahedron.

	Equip $K$ with the orientation by $V = (a,b,c,d)$. Observe that the total order on $\Vertex(K) = V$ is denoted using a tuple.
	Given below are the standard ordered bases $\basis{K}_n$ on the chain groups $C_n(K) = C_n(K;\ints)$ of the oriented simplicial complex $K$:
	\begin{equation*}
		\begin{aligned}
			\basis{K}_0 &= \Big( [a],[b],[c],[d] \Big)	\\
			\basis{K}_1 &= \Big( [a,b], [a,c], [a,d], [b,c], [b,d], [c,d] \Big)
		\end{aligned}
		\qquad\quad 
		\begin{aligned}
			\basis{K}_2 &= \Big( [a,b,c], [a,b,d], [a,c,d], [b,c,d] \Big) 
		\\
		\basis{K}_3 &= \Big( [a,b,c,d] \Big)
		\end{aligned}
	\end{equation*}
	with $\basis{K} = \emptyset$ for $n \in \ints$ with $n \neq 0,1,2,3$.
	Observe that for $n=2$, the oriented $n$-simplex $[a,b,c]$ comes before $[a,b,d]$ since the first two entries match and we have $c \leq d$ as elements of $V$ on the third entry.
	Similarly, $[a,c,d]$ is listed before $[b,c,d]$ since on the first entry, we have $a \leq b$ on the first entry.

	Following \fref{remark:shorthand-for-oriented-simplices},
	we may write $[v_0, \ldots, v_n]$ as a string $v_0 \cdots v_n$ of vertices for brevity.
	Then, we can describe the $n$\th chain groups $C_n(K)$ of $K$ as follows:
	\begin{equation*}
		C_n(K) = \begin{cases}
			\ints\ket{a,b,c,d} &\text{ if } n = 0 \\
			\ints\ket{ab, ac, ad, bc, bd, cd} 
				&\text{ if } n = 1 \\
			\ints\ket{abc, abd, acd, bcd} 
				&\text{ if } n = 2 \\
			\ints\ket{abcd} &\text{ if } n = 3 \\
			0 &\text{ otherwise }
		\end{cases}
	\end{equation*}
	Examples of elements of $C_2(K)$ include $\sigma_1 = 2abc - abd$ and $\sigma_2 = 3acd + 7abc - 10bcd$.
\end{example}

\spacer

Next, we provide a definition for the boundary map related to the simplicial chain groups.

\begin{definition}\label{defn:boundary-map}
	Let $K$ be a simplicial complex and $R$ a PID.
	For each $n \in \ints$ with $n \geq 0$,
		define the $n$\th \textbf{boundary map} or the $n$\th \textbf{boundary homomorphism} $\boundary_n: C_{n}(K;R) \to C_{n-1}(K;R)$ to be the $R$-module homomorphism given by
	\begin{equation}\label{eqn:simplicial-boundary-map}
		\boundary_n\biggl(
			[v_0, \ldots, v_n]
		\biggr) 
		:= \sum_{i=0}^n (-1)^i [v_0, \ldots, \hat{v}_i, \ldots, v_n] 
		= \sum_{i=0}^n (-1)^i [v_0, \ldots, v_{i-1}, v_{i+1}, \ldots, v_n]
	\end{equation}
	where $\hat{v}_i$ indicates the removal of the vertex $v_i$ in the ordering $(v_0, \ldots, v_n)$.
	For $n < 0$, $\boundary_n: C_n(K;R) \to C_{n-1}(K;R)$ can only be the trivial homomorphism since $C_n(K;R) = 0$. 

	An $n$-\textbf{cycle} is an element of $\ker(\boundary)_n \subseteq C_n(K;R)$
	and an $n$-\textbf{boundary} is an element of $\im(\boundary_{n+1}) \subseteq C_{n+1}(K;R)$.
\end{definition}
\remark{
	It can be verified that Equation\ \ref{eqn:simplicial-boundary-map} is well-defined for any choice of coset representative of $\sigma \in C_n(K;R)$ (relative to \fref{defn:chain-groups}) and produces a well-defined homomorphism.
	In particular, for any oriented $n$-simplex $[\sigma] = [v_0, \ldots, v_n]$ with an arbitrary coset representative $(v_0, \ldots, v_n)$ and for any permutation $\pi: [n] \to [n]$,
	it can verified that
	\begin{align*}
		\boundary_n\bigl( [v_0, \ldots, v_n] \bigr)
		:=
		\sum_{i=0}^n (-1)^i [v_0, \ldots, \hat{v}_i, \ldots, v_n]
		&=
		(\sgn\pi)\sum_{i=0}^n (-1)^i [
			v_{\pi(0)}, \ldots, \hat{v}_{\pi(i)}, \ldots, v_{\pi(n)}
		]
		\\[-2pt] 
		&=:
		(\sgn\pi)\,
		\boundary_n\bigl( [v_{\pi(0)}, \ldots, v_{\pi(n)}] \bigr).
	\end{align*}	
	with the first and second equality being given by Equation (1).
}

\HIDE{Given a simplicial complex $K$, the chain groups $C_n(K;R)$ of $K$ induce the notion of orientation on each $n$-simplex $\sigma = \set{v_0, \ldots, v_n}$ of $K$ with $n \geq 1$.
The boundary map $\boundary_n$ can be understood as a specification for how the orientation of an $n$-simplex $\sigma$ induces orientation on each of its facets (i.e.\ a $(n-1)$-simplex contained in $\sigma$).
This specification also allows for the ``cancellation'' of facets in a sum of oriented $n$-simplices that matches with geometric intuition.
We describe this for $n=1$ and $n=2$ below.
For convenience, let $R = \ints$.
\begin{enumerate}[left=0.7in]
	\item[For $n=1$:]
	The orientation on $\set{a,b}$ corresponds to a direction on the line segment between $a$ and $b$.
	The boundary map $\boundary_1: C_1(K;\ints) \to C_0(K;\ints)$ then assigns the origin vertex an ``orientation'' of $(-1)$ and the end vertex with $(+1)$, as illustrated below:
	\centernote{
		Doggo.
	}
	Observe that the alternating signs $(-1)^i$ in the definition of $\boundary_n$ are there to account for the choice of basis of $C_n(K;R)$.

	A path along the $1$-simplices of $K$ can be represented using a sum $\rho \in C_1(K;\ints)$ of oriented $1$-simplices (with the chosen orientation matching the direction of the path).
	The calculation by the $1$\st boundary map $\boundary_1$ then produces two vertices, as we would expect geometrically.
	We illustrate this below:
	\centernote{
		Blah.
	}
	Then, given a loop along the $1$-simplices of $K$ represented by $\theta \in C_1(K;\ints)$, we would expect $\boundary_1(\theta) = 0$. We illustrate this below:
	\centernote{
		Loopy-loop.
	}
	In fact, the $1$-cycles $\theta \in \ker(\boundary_1)$ of $C_1(K;\ints)$ correspond exactly to collections of these loops.
	Furthermore, the $1$-boundaries $\beta \in \im(\boundary_2)$ of $C_1(K;\ints)$ correspond to loops that bound some area in $K$, as illustrated below:
	\centernote{
		bound-dead.
	}
	As a sidenote, observe that a $1$-cycle forms a shape homeomorphic to a $1$-sphere $S^1$ (i.e.\ a circle) and that $1$-boundaries function similarly to the topological boundary of $S^1 \cong \boundary(D^2)$ of the $2$-disk $D^2$ as a set in $\reals^2$.
	The boundaries of $1$-chains also corresponds to a collection of $0$-spheres $S^0$ (i.e.\ $S^0 = \set{-1,1}$) which are topological boundaries of the $1$-disk $D^1$ (i.e.\ an interval $[-1,1]$) as a set in $\reals$.
	
	\item[For $n=2$:]
	An orientation on a $2$-simplex $\set{a,b,c}$ can be interpreted as a choice of rotation on the triangle by $\set{a,b,c}$.
	The boundary map $\boundary_2: C_2(K;\ints) \to C_1(K;\ints)$ then induces orientations on the $1$-simplices of $\set{a,b,c}$ that follow this rotation, as illustrated below:
	\centernote{
		time to twirl
	}
	The boundary map also allows the cancellation of edges on $2$-simplices that share a facet (assuming these simplices are appropriately oriented), as illustrated below.
	\centernote{
		time to twirl part 2.
	}
	Then, $2$-cycle would correspond to a collection of oriented $2$-simplices (again, appropriately oriented) whose union is homeomorphic to a $2$-sphere $S^2$ (or a collection of $2$-spheres).
	We illustrate this below.
	\centernote{
		time to twirl, ball edition.
	}
	Similarly, a $2$-boundary is a $2$-cycle that bounds a volume covered by $K$. That is, a $2$-boundary is similarly to the topological boundary $\boundary(D^2) = S^2$ of the $2$-disk $D^2$ (i.e.\ a solid ball) as a set in $\reals^3$.
\end{enumerate}
Note that, in both cases, we expect $n$-boundaries to also be $n$-cycles.
This is confirmed by the following result.}

Next, we state the property that determines that the collection of chain groups and boundary maps form a chain complex.

\begin{proposition}\label{prop:boundary-squared-zero}
	For each $n \in \ints$, $\boundary_{n}\boundary_{n+1} = 0$.
	That is, $\im(\boundary_{n+1}) \subseteq \ker(\boundary_{n})$.
\end{proposition}
\remark{
	This can be proven by direct calculation on a coset $[v_0, \ldots, v_{n+1}] \in C_{n+1}(K;R)$ with representative $(v_0, \ldots, v_{n+1})$.
	A similar argument is presented under the proof of \cite[Theorem 4.6]{algtopo:rotman}.
}

Observe that this property on $\set{\boundary_n}$ allows us to take the quotient of the $R$-module $\ker(\boundary_n)$ by its submodule $\im(\boundary_{n+1})$.
This quotient is exactly the simplicial homology of a simplicial complex. We state this in more detail below.

\begin{definition}\label{defn:simplicial-chain-complex}
	The \textbf{simplicial chain complex} of a simplicial complex $K$ with coefficients in a PID $R$ is the chain complex 
	$C_\ast(K;R) := \bigl( C_n(K;R), \boundary_n \bigr)_{n \in \ints}$
	of simplicial chain groups $C_n(K;R)$ and boundary maps $\boundary_n: C_n(K;R) \to C_{n-1}(K;R)$, illustrated as the following sequence:
	\begin{equation*}
		\cdots 			\Xrightarrow{\,\boundary_{n+2}}
		C_{n+1}(K;R) 	\Xrightarrow{\,\boundary_{n+1}}
		C_{n}(K;R)		\Xrightarrow{\,\boundary_{n}}
		C_{n-1}(K;R)	\Xrightarrow{\,\boundary_{n-1}}
		\cdots 
	\end{equation*}
	For each $n \in \ints$, define the $n$\th \textbf{simplicial homology group} $H_n(K;R)$ 
	and the $n$\th \textbf{Betti number} $\betti_n(K;R)$ of $K$ with coefficients in $R$ as follows:
	\begin{equation*}
		H_n(K;R) := \frac{\ker(\boundary_n)}{\im(\boundary_{n+1})}
		\qquad\text{ and }\qquad
		\beta_n(K;R) := \rank\bigl( H_n(K;R) \bigr)
	\end{equation*}
	If $R = \ints$, we write $C_\ast(K) := C_\ast(K;\ints)$ and $H_n(K) := H_n(K;\ints)$.
\end{definition}
\remark{
	We follow the convention by \cite{algtopo:rotman} and use $(\ast)$ as the ``placeholder'' of the index $n \in \ints$ of the chain complex $C_\ast = (C_n, \boundary_n)$ with chain groups $C_n$ and differentials $\boundary_n: C_n \to C_{n-1}$.
}

Observe that if the simplicial complex $K$ is a finite simplicial complex, then the $n$\th chain group $C_n(K;R)$ must also be finitely generated for all $n \in \ints$.
In this case, $H_n(K;R)$ is a finitely-generated module over a PID $R$ and the Structure Theorem (\fref{thm:structure-theorem}) for finitely-generated modules over $R$ applies.
We discuss this and a method of calculating the homology of free chain complexes using matrices over $R$ later in \fref{chapter:matrix-calculation}.

One of the key characteristics of simplicial homology is that it is an invariant of the homeomorphism type of topological spaces.
We state this in a theorem below.

\begin{theorem}\label{thm:homology-invariant-under-homeo}
	Let $K$ and $L$ be simplicial complexes and let $R$ be a PID.
	If the geometric realizations of $K$ and $L$ are homeomorphic as topological spaces, 
		then $H_n(K;R) \cong H_n(L;R)$ as $R$-modules.
\end{theorem}
\remark{
	We refer to \cite[Theorem 7.13]{algtopo:rotman} for a proof.
	This theorem actually generalizes to the case of \textit{homotopy equivalence}, which is an equivalence weaker than a homeomorphism. 
}


Note that there are other ways to define the homology of a topological space.
With simplicial homology, we can define the homology group $H_n(X;R)$ of a \textit{topological space} $X$ to be the homology group $H_n(K;R)$ of a triangulation $K$ of $X$.
Observe that this is well-defined by \fref{thm:homology-invariant-under-homeo}, i.e.\ $H_n(K;R)$ is determined up to isomorphism regardless of the triangulation $K$.
Note that this definition requires that a triangulation of $X$ exist in the first place.
Homology can also be defined using other representations of spaces, e.g.\ singular homology of topological spaces, 
	simplicial homology of $\Delta$-complexes,
	and cellular homology of CW-complexes.
The key point here is that all these homology theories must produce isomorphic homology groups as with simplicial homology.
We refer to the discussion in \cite[Theorem 2.27]{algtopo:hatcher} about singular homology and \cite[Theorem 2.35]{algtopo:hatcher} about cellular homology.

We can take advantage of this by comparing results from other homology theories against results from simplicial homology.
This becomes particularly useful when the triangulations of a relatively simple space are cumbersome or unwieldy.
We provide an example of this below.

\begin{example}
	Let $K$ be the abstract simplicial complex given as follows, with simplices written as strings (following \fref{remark:shorthand-for-oriented-simplices})
	for brevity and with vertex indices always taken modulo $3$.

	\vspace{\parskip}
	\begin{center}
	\begin{tabular}{CCC}
		\text{nine $0$-simplices of $K$} &:&
			\Vertex(K) = \set{ a_0, a_1, a_2, b_0, b_1, b_2, c_0, c_1, c_2 }
		\\[5pt]
		\text{twenty-seven $1$-simplices of $K$} &:&
			\left\{\begin{array}{ccc}
				a_i a_{i+1} 
				& b_i b_{i+1} 
				& c_i c_{i+1} 
				\\
				a_{i} b_{i}
				& b_i c_i
				& c_i a_i 
				\\
				a_i b_{i+1} 
				& b_i c_{i+1} 
				& c_i a_{i+1}
			\end{array}
			\,\text{ for } i=0,1,2 \,
			\right\}\vspace{5pt}
		\\[5pt]
		\text{eighteen $2$-simplices of $K$} &:&
			\left\{\begin{array}{ccc}
				a_i a_{i+1} b_{i+1}, 
				& b_i b_{i+1} c_{i+1},
				& c_i c_{i+1} a_{i+1},
				\\
				a_i b_i b_{i+1},
				& b_i c_i, c_{i+1},
				& c_i a_i, a_{i+1}
			\end{array}
			\,\text{ for } i=0,1,2 \,
			\right\}
	\end{tabular}
	\end{center}
	
	\vspace{\parskip}
	\noindent
	We can interpret $K$ to be a geometric simplicial complex in $\reals^3$ by defining the vertices of $K$ as follows:
	\begin{gather*}
		\text{For each } i \in \set{0,1,2}: \quad
		a_i := T\left(
			\frac{2\pi}{3}i, \frac{\pi}{3}
		\right)
		,\quad
		b_i := T\left(
			\frac{2\pi}{3}i, -\frac{\pi}{3}
		\right)
		,\quad
		c_i := T\left(
			\frac{2\pi}{3}i, -\pi
		\right),
		\\[3pt]
		\text{ with }
		T(u,v) = \left(\,\begin{matrix}
			\cos(u)\bigl( R + r\cos(v) \bigr) \\
			\sin(u)\bigl( R + r\cos(v) \bigr) \\
			r\sin(v)
		\end{matrix}\,\right)
		\text{ with } r=1 \text{ and } R = 5.
	\end{gather*}
	We illustrate $K$ as a geometric simplicial complex below and claim that there exists a triangulation $|K| \to T^2$ to the torus $T^2$ parametrized by $T(u,v)$ with $(u,v) \in \reals\times\reals$.

	\vspace{5pt}
	\begin{center}\begin{minipage}[t]{0.3\linewidth}\centering
		The Geometric \\
		Simplicial Complex $K$: \\[2pt]
		\includegraphics[width=0.9\linewidth]{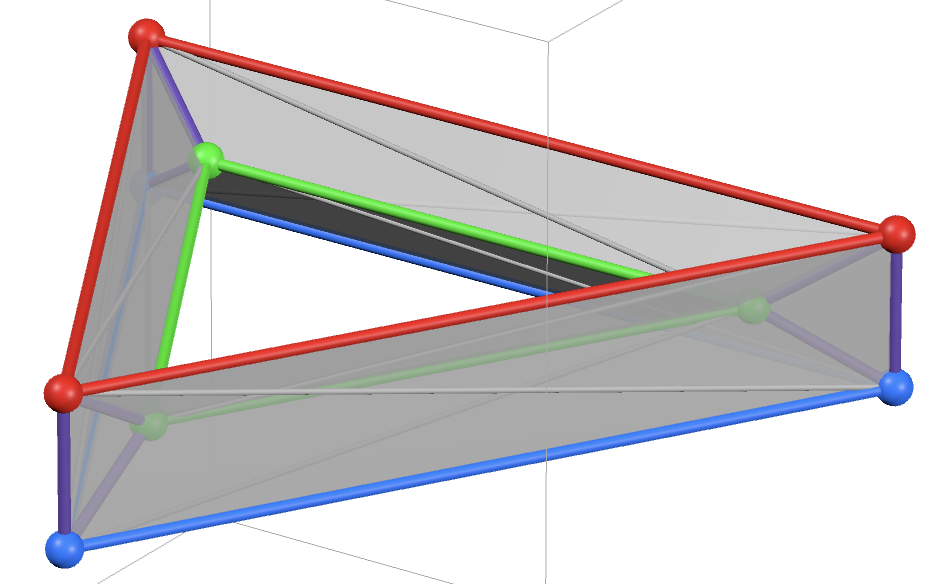}

		\vspace{10pt}
		The Torus $T^2$ \\
		by $T(u,v)$, $(u,v) \in \reals \times \reals$: \\
		\includegraphics[width=0.9\linewidth]{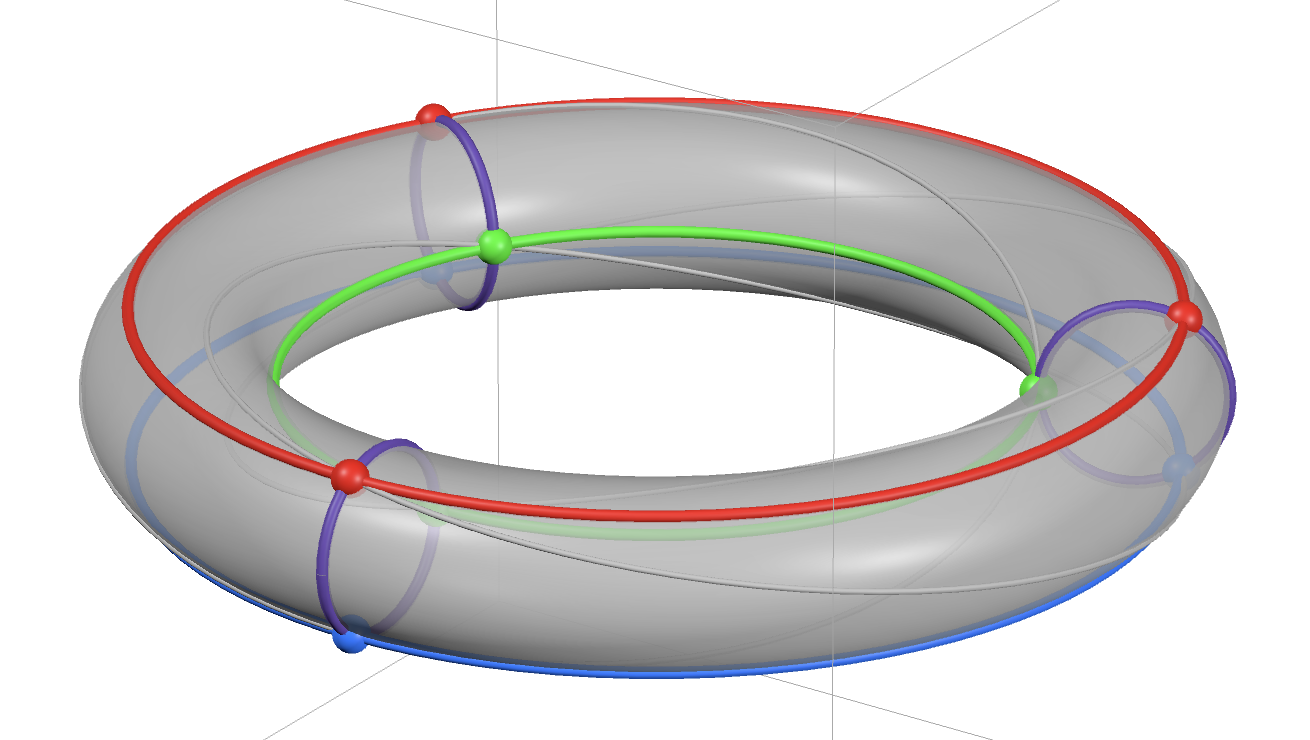}
	\end{minipage}\begin{minipage}[t]{0.7\linewidth}
		We use the following color scheme for the simplices of $K$ and their corresponding images in $T^2$ under the triangulation $|K| \to T^2$.
	
		\begin{enumerate}[itemsep=5pt, topsep=5pt, labelsep=5pt, leftmargin=50pt]
			\item[\redtag:]
			In $K$, 
			the vertices $a_i$ and the three $1$-simplices $a_i a_{i+1}$ for $i \in \set{0,1,2}$. 
			In $T^2$, the major circle given by 
			$T(u, \frac{\pi}{3})$ with $u \in \reals$.
			
			\item[\bluetag:]
			In $K$, the vertices $b_i$ and the three $1$-simplices $b_ib_{i+1}$ for $i \in \set{0,1,2}$.
			In $T^2$, the major circle given by 
			$T(u, -\frac{\pi}{3})$ with $u \in \reals$.
	
			\item[\greentag:]
			In $K$, the vertices $c_i$ and the three $1$-simplices $c_i c_{i+1}$ for $i \in \set{0,1,2}$.
			In $T^2$, the major circle given by 
			$T(u,-\pi)$ with $u \in \reals$.
	
			\item[\purpletag:]
			In $K$, the nine $1$-simplices $a_ib_i$, $b_ic_i$, and $c_i a_i$ for $i = 0,1,2$.
			In $T^2$, the three minor circles connecting the vertices $a_i, b_i, c_i$ given by 
			$T(\frac{2\pi}{3}i, v)$ for $v \in \reals$
			for each $i = 0,1,2$.
	
			\item[\colortag{gray}{gray}{15}:]
			In $K$, all eighteen $2$-simplices (triangles) and 
			the nine $1$-simplices $a_i b_{i+1}$, $b_i c_{i+1}$, and $c_i a_{i+1}$ for $i \in \set{0,1,2}$. \\[5pt]
			For each $i \in \set{0,1,2}$,
			the thin gray loop in $T^2$ that connects the vertices $a_i, b_{i+1}, c_{i+2}$ are given by 
			$
				T(t, \frac{\pi}{3} + \frac{2\pi}{3}i - t)
			$ with $t \in \reals$.
		\end{enumerate}
	\end{minipage}\end{center}

	\vspace{\parskip}\noindent
	Observe that $K$ has $9+27+18=54$ simplices and calculating its homology groups by hand can be very cumbersome to do by hand, e.g.\ $C_1(K;\ints)$ has $27$ basis elements.
	Since $K$ is a triangulation of $T^2$, we can use \fref{thm:homology-invariant-under-homeo} and the known homology groups of $T^2$ to determine $H_n(K;\ints)$ up to $\ints$-module isomorphism as follows:
	\begin{equation*}
		H_n(K_T;\ints) \cong 
		H_n(T^2;\ints) = H_n(T^2) \cong \begin{cases}
		\ints 				&\text{ if } n=0 \\
		\ints \oplus \ints 	&\text{ if } n=1 \\
		\ints 				&\text{ if } n=2 \\
		0 					&\text{ if } n \geq 3 \\
		\end{cases}
	\end{equation*}
	For constrast, \cite[Example 2.3, p106]{algtopo:hatcher} uses delta-complexes to calculate $H_n(T^2)$ and the chain groups corresponding to $C_0(T^2), C_1(T^2), C_2(T^2)$ have $1$, $3$, and $2$ basis elements respectively.
	As a sidenote, it has been proven that the minimal triangulation of $T^2$ consists of fourteen $2$-simplices (triangles).
\end{example}

	\clearpage


\section{Functorial Constructions in Simplicial Homology}
\label{section:functors-in-simplicial-homology}

In this section, 
we discuss the functorial nature of the construction of the simplicial chain complex $C_\ast(K;R) = \bigl(C_n(K;R), \boundary_n\bigr)_{n \in \ints}$ for a simplicial complex $K$, as discussed in \fref{section:simplicial-homology}.
For reference, we briefly discuss categories and functors in \fref{appendix:cat-theory}. 
We consider the following categories in this section:
\begin{enumerate}
	\item 
	The category $\catmod{R}$ of modules over a PID $R$ and $R$-module homomorphisms.

	\item 
	The category $\catchaincomplex{\catmod{R}}$ of chain complexes of $R$-modules and of chain maps.
\end{enumerate}
Note that the category $\catchaincomplex{\catmod{R}}$ has an accompanying \textit{chain homology functor} 
$H_n(-): \catchaincomplex{\catmod{R}} \to \catmod{R}$
with the following assignments:
Let $C_\ast = (C_n, \boundary_n)_{n \in \ints}$ and $A_\ast = (A_n, \alpha_n)_{n \in \ints}$ be chain modules with $R$-modules $C_n$ and $A_n$ and differentials $\boundary_n: C_n \to C_{n-1}$ and $\alpha_n: A_n \to A_{n-1}$.
\begin{enumerate}
	\item 
	$H_n(-)$ maps $C_\ast$ to its $n$\th homology group $H_n(C_\ast)$ defined by $H_n(C_\ast) = \ker(\boundary_n) \bigmod \im(\boundary_{n+1})$.
	Note that this is well-defined since $\im(\boundary_{n+1}) \subseteq \ker(\boundary_{n})$ by the definition of chain complexes.

	\item 
	Given a chain map $f_\ast: C_\ast \to A_\ast$ with $f_\ast = \set{f_n : C_n \to A_n}_{n \in \ints}$,
	$H_n(-)$ maps $f_\ast$ to the $R$-module homomorphism 
	$H_n(C_\ast) \to H_n(A_\ast)$ induced by the map $f_n: C_n \to A_n$ and the quotient (or cokernel) construction.
	Note that the map on $H_n(C_\ast) \to H_n(A_\ast)$ is well-defined since $f_{n-1} \boundary_{n} = \alpha_{n} \circ f_n$ by definition of a chain map. 
\end{enumerate}
We claim that simplicial homology, as discussed in \fref{section:simplicial-homology}, corresponds to a composition of functors: one that creates the simplicial chain complex $C_\ast(K;R)$ of a simplicial complex $K$ and the $n$\th homology functor on $\catchaincomplex{\catmod{R}}$.
Recall that a functor $F: \catname{C} \to \catname{D}$ between categories $\catname{C}$ and $\catname{D}$ is an assignment of objects and morphisms of $\catname{C}$ to those of $\catname{D}$ such that the following properties are satisfied:
\begin{enumerate}[label=(F\arabic*), ref=F\arabic*]
	\item \label{functor-domain-codomain}
	For each morphism $f: X \to Y$ in $\catname{C}$,
	the functor maps $f$ to a morphism with domain $F(X)$ and codomain $F(Y)$, i.e. $f \mapsto F(f) : F(X) \to F(Y)$.

	\item \label{functor-identity}
	The functor respects \textit{identity} of objects and of morphisms: For all objects $X$ in $\catname{C}$ with identity morphism $\id_X: X \to X$, the morphism $F(\id_X): F(X) \to F(X)$ is exactly the identity of $F(X)$ in $\catname{D}$,
	i.e.\ $F(\id_X) = \id_{F(X)}$.

	\item \label{functor-composition}
	The functor respects \textit{composition} of morphisms:
	For all composable pairs of morphisms 
		$f: X \to Y$ and $g: Y \to Z$,
		$F(g \circ f) = F(g) \circ F(f)$.
\end{enumerate}
In this section, we argue that the construction of the $n$\th chain group $C_n(K;R)$ of a simplicial complex $K$ and of the simplicial chain complex $C_\ast(K;R)$ correspond to the functors by defining the morphism assignments for both and showing that the properties above are satisfied.
To start, we state a result involving the formation of the category of simplicial complexes.

\begin{theorem}\label{thm:cat-simp-theorem}
	Abstract simplicial complexes and simplicial maps form a well-defined category.
\end{theorem}
\remark{
	This theorem is taken from \cite[Theorem 7.7]{algtopo:rotman}, 
	which states ``a routine check'' as proof. 
	We believe this refers to checking that the axioms listed in the definition of category, i.e.\ as listed \fref{defn:category}, are satisfied.
	Since simplicial maps can be seen as functions between sets of sets, 
		the composition law and the identity morphism designed to each abstract simplicial complex follow those on the category $\catname{Set}$ of sets and functions.
	Then, it suffices to check that composition is unital with identity morphisms and is associative, i.e.\ \fref{defn:category}(iv,v).
}

\begin{definition}\label{defn:cat-simp}
	Denote the category given in \fref{thm:cat-simp-theorem} by $\catsimp$ and call it the \textbf{category of simplicial complexes and simplicial maps}, or the \textbf{category of simplicial complexes} for convenience.
\end{definition}


We want to show that the construction of the chain groups in \fref{defn:chain-groups} correspond to a functor $\catsimp \to \catmod{R}$ with $\catmod{R}$ being the category of $R$-modules and $R$-module homomorphisms.
\fref{defn:chain-groups} determines the object assignment.
Below, we provide a definition for the corresponding morphism assignment, taken from~\cite{algtopo:rotman}.

\begin{definition}\label{defn:induced-maps-on-chain-groups}
	Let $f: K \to L$ be a simplicial map between simplicial complexes $K$ and $L$.
	For each $n \in \nonnegints$, 
	define the \textbf{homomorphism $f_{n,\hash}: C_n(K;R) \to C_n(L;R)$ on the $n$\th chain groups induced by $f$} as follows:
	\begin{equation*}
		f_{n,\hash}\bigl([\sigma]\bigr)
		= f_{n,\hash}\bigl([v_0, v_1, \ldots, v_n] \bigr)
		:= \begin{cases}
			\bigl[ f(v_{0}), f(v_{1}), \ldots, f(v_{n}) \bigr]
			&\text{ if $f(\sigma)$ is an $n$-simplex in $L$} \\
			0 &\text{ otherwise }
		\end{cases}
	\end{equation*}
	where $(v_0, \ldots, v_n)$ is the coset representative of $[\sigma] \in C_n(K;R)$ with $\sigma = \set{v_0, \ldots, v_n}$
	and $(f(v_0), \ldots, f(v_n))$ that of $f_{n,\hash}([\sigma]) \in C_n(K;R)$ of the $n$-simplex $f(\sigma) = \set{f(v_0), \ldots, f(v_n)}$.
\end{definition}
\remark{
	For brevity, we often suppress the dimension $n$ in $f_{n,\hash}$, i.e.\ we write $f_\hash$ for $f_{n,\hash}$ 
}

Note that this construction is defined on the coset representatives of $C_n(K;R)$.
Consequently, if $L$ is equipped with an orientation and the standard basis $L[n]$ by \fref{defn:standard-basis-on-chain-groups} is used,
appropriate sign changes might be needed for $[ f(v_{0}),\ldots, f(v_{n}) ] \in C_n(L;R)$, i.e.\ 
\begin{equation*}
	\bigl[ f(v_{0}),\ldots, f(v_{n}) \bigr] 
	= (\sgn\pi)\bigl[
		f(v_{\pi(0)}), \ldots, f(v_{\pi(n)})
	\bigr] 
\end{equation*}
where $\pi: [n] \to [n]$ permutes the vertices of $f(\sigma)$ such that $\bigl(f(v_{\pi(0)}), \ldots, f(v_{\pi(n)})\bigr)$ is ordered with respect to the orientation on $L$.

Observe that the definition $f_{n,\hash}: C_n(K;R) \to C_n(L;R)$ above already satisfies property \ref{functor-domain-codomain}.
Below, we state that this also satisfies the other two properties required for functors.

\begin{proposition}\label{prop:chain-group-functorial}
	Let $R$ be a PID and let $n \in \ints$.
	\begin{enumerate}
		\item 
		Let $K$ be a simplicial complex with identity simplicial map $\id_K: K \to K$.
		Then, $(\id_K)_\hash = \id_{C_n(K;R)}$.

		\item 
		Let $f: K_1 \to K_2$ and $g: K_2 \to K_3$ be simplicial maps on simplicial complexes $K_1$, $K_2$, and $K_3$.
		Then, $(g \circ f)_\hash = g_\hash \circ f_\hash$.
	\end{enumerate}
\end{proposition}
\remark{
	This can be proven by calculation on some arbitrary element $[\sigma] = [v_0, \ldots, v_n]$ of $C_n(K;R)$ or $C_n(K_1;R)$ with $\sigma = \set{v_0, \ldots, v_n}$ an $n$-simplex of $K$.
}

Since \fref{defn:induced-maps-on-chain-groups} satisfies properties 
\ref{functor-domain-codomain}, 
\ref{functor-identity},
and \ref{functor-composition}
for functors, we can define the construction of the $n$\th chain group $C_n(K;R)$ as a functor, as given below.

\begin{definition}\label{defn:simplicial-chain-group-functor}
	For each $n \in \ints$,
	define \textbf{$n$\th simplicial chain group functor} $C_n(-;R): \catsimp \to \catmod{R}$ 
	with coefficients in a PID $R$ as follows:
	\begin{enumerate}
		\item 
		The object assignment maps a simplicial complex $K$ to the $n$\th chain group $C_n(K;R)$ as given by \fref{defn:chain-groups}.
		
		\item 
		The morphism assignment maps a simplicial map $f: K \to L$ between simplicial complexes $K$ and $L$ to map $f_\hash: C_n(K;R) \to C_n(L;R)$ given by \fref{defn:induced-maps-on-chain-groups}.
	\end{enumerate}
	If $R = \ints$, we write $C_n(-) := C_n(-;\ints)$.
\end{definition}

\HIDE{
	Next, we state a minor result that will be useful later in \fref{section:construction-of-persistent-homology} in the context of persistent homology.

	\begin{lemma}\label{lemma:map-on-chain-groups-by-inclusion}
		Let $L$ be subcomplex of a simplicial complex $K$ with inclusion map $i: L \to K$.
		For all $n \in \ints$:
			the map $i_\hash: C_n(L;R) \to C_n(K;R)$
			is determined by the identity map on $C_n(K;R)$,
		i.e.\ for all $[\sigma] \in C_n(L;R)$, 
		\begin{equation*}
			i_\hash([\sigma]) = \id_{C_n(L;R)}([\sigma])
			= \id_{C_n(K;R)}([\sigma])
		\end{equation*}
	\end{lemma}
	\begin{proof}
		Let $[\sigma] = [v_0, \ldots, v_n] \in C_n(L;R)$ with $\sigma = \set{v_0, \ldots, v_n}$. 
		Since $L \subseteq K$, $\sigma$ is an $n$-simplex in $K$ and 
		\begin{equation*}
			i_\hash([v_0, \ldots, v_n]) 
			\overset{{\substack{\text{Prop} \\ \text{\ref{prop:chain-group-functorial}}\\[-12pt]\strut}}}{=} 
			[v_0, \ldots, v_n]
			= 
			\id_{C_n(K;R)}\bigl( [v_0, \ldots, v_n] \bigr)
		\end{equation*}
		Note the last equality is valid since the order of vertices in the coset representative $(v_0, \ldots, v_n)$ is unchanged.
	\end{proof}
}

To construct the functor $\catsimp \to \catchaincomplex{\catmod{R}}$, we claim that we can simply collect all maps $f_{n,\hash}: C_n(K;R) \to C_n(L;R)$ on the chain groups induced by a simplicial map $f: K \to L$ to define a chain map between the simplicial chain complexes $C_\ast(K;R)$ and $C_\ast(L;R)$.
The following result allows us to do this:

\begin{proposition}
	\label{prop:induced-maps-on-chains-respect-boundary}
	Let $f: K \to L$ be a simplicial map between simplicial complexes $K$ and $L$.
	For all $n \in \ints$,
	\begin{equation*}
		f_{n-1,\hash} \circ \boundary_n^K = \boundary_n^L \circ f_{n,\hash}
	\end{equation*}
	where $\boundary_n^K: C_n(K;R) \to C_{n-1}(K;R)$ 
	and $\boundary_n^L: C_n(K;L) \to C_{n-1}(K;L)$ refer to the $n$\th simplicial boundary map of $K$ and $L$ respectively,
	i.e.\ the following diagram commutes:
	\vspace{-0.5\baselineskip}
	\begin{displaymath}
	\begin{tikzcd}[column sep=4em]
		C_n(K;R) \arrow[r, "\boundary_n^K"] 
					\arrow[d, "f_{n,\hash}" swap]
			& C_{n-1}(K;R) 
					\arrow[d, "f_{n-1,\hash}"]
		\\
		C_n(L;R) \arrow[r, "\boundary_n^L"] & C_{n-1}(L;R)
	\end{tikzcd}
	\end{displaymath}  
	That is, the collection $\set{f_{n,\hash}}_{n \in \ints}$ determines a chain map $\set{f_{n,\hash}}: C_\ast(K;R) \to C_\ast(L;R)$.
\end{proposition}
\remark{
	This can be proven by direct calculation on an arbitrary element $[v_0, \ldots, v_n] \in C_n(K;R)$ with coset representative $(v_0, \ldots, v_n)$ on an $n$-simplex $\set{v_0, \ldots, v_n}$ of $K$.
	Note that boundary maps, as given in \fref{defn:boundary-map}, are well-defined on any choice of coset representative.
}

Observe that the chain map given by $\set{f_{n,\hash}}_{n \in \ints}$, as denoted above, already satisfies the domain and codomain property of functors as stated in Property \ref{functor-domain-codomain}.
The fact that $\set{f_{n,\hash}}_{n \in \ints}$ respects identity as in Property \ref{functor-identity} 
and respects composition as in Property \ref{functor-composition} follows from those properties being satisfied by $f_{n,\hash}$ for each $n \in \ints$ separately.
Therefore, we can define a functor on $\catsimp \to \catchaincomplex{\catmod{R}}$ as follows: 

\begin{definition}\label{defn:simplicial-chain-complex-functor}
	Define the $n$\th \textbf{simplicial chain complex functor} 
	$C_\ast(-;R): \catsimp \to \catchaincomplex{\catmod{R}}$
	with coefficients in a PID $R$ as follows:
	\begin{enumerate}
		\item $C_\ast(-;R)$ maps a simplicial complex $K$ to the simplicial chain complex $C_\ast(K;R) = (C_n(K;R), \boundary_n)_{n \in \ints}$ as given by \fref{defn:simplicial-chain-complex}.
		
		\item 
		$C_\ast(-;R)$ maps a simplicial map $f: K \to L$ between simplicial complexes $K$ and $L$ to the collection $\set{f_{n,\hash}}_{n \in \ints}$ of maps $f_{n,\hash}: C_n(K;R) \to C_n(L;R)$ as given by \fref{defn:induced-maps-on-chain-groups}.
	\end{enumerate}
	If $R = \ints$, we write $C_\ast(-) := C_\ast(-;\ints)$ for brevity.
\end{definition}

Bringing this all together, we have the following proposition:

\begin{proposition}\label{prop:decomposition-simplicial-homology}
	Let $K$ be a simplicial complex and $R$ a PID.
	For all $n \in \ints$:
	\begin{equation*}
		H_n(K;R) = \Bigl(
			H_n \circ C_\ast(-;R)
		\Bigr)(K)
		= \Bigl(
			H_n(C_\ast(K;R))
		\Bigr)
	\end{equation*}
	where $H_n(K;R)$ is the simplicial chain group of $K$ by 
	\fref{defn:simplicial-chain-complex},
	$H_n: \catchaincomplex{\catmod{R}} \to \catmod{R}$ is the chain homology functor,
	and $C_\ast(-;R)$ is the simplicial chain complex functor by 
	\fref{defn:simplicial-chain-complex-functor}.
\end{proposition}
\begin{proof}
	The object assignment by $H_n(-) \circ C_\ast(-;R)$ is exactly as described in \fref{defn:simplicial-chain-complex} for $H_n(K;R)$.
	Note that this proposition disregards the morphism assignment of the relevant functors.
\end{proof}

For convenience, we define the composition of the simplicial chain complex functor and the homology functor to be the simplicial chain group functor.

\begin{definition}\label{defn:simplicial-homology-functor}
	For each $n \in \ints$, define the $n$\th \textbf{simplicial homology functor} 
	$H_n(-;R): \catsimp \to \catmod{R}$ with coefficients in a PID $R$ as the following composition of functors:
	\begin{equation*}
		H_n(-;R) := H_n \circ C_n(-;R)
	\end{equation*}
	where $H_n: \catchaincomplex{\catmod{R}} \to \catmod{R}$ refers to the chain homology functor.
	Given a simplicial map $f: K \to L$, 
		let $f_\ast: H_n(K;R) \to H_n(L;R)$ be the map on homology induced by application of the functor $H_n(-;R)$.
\end{definition}

Below, we identify a result that will be useful later in \fref{section:construction-of-persistent-homology} in the context of persistent homology.

\begin{lemma}\label{lemma:boundary-on-chain-groups-by-inclusion}
	Let $L$ be a subcomplex of some simplicial complex $K$ and let $i: L \to K$ be the corresponding inclusion map.
	For each $n \in \ints$:
	\begin{enumerate}
		\item 
		The $n$\th chain group $C_n(L;R)$ of $L$ is a submodule of $C_n(K;R)$ and 
		for all $[\sigma] \in C_n(L;R)$,
		$i_{n,\hash}\bigl([\sigma]\bigr) = \id_{C_n(K;R)}\bigl([\sigma]\bigr) = [\sigma]$
			where $\id_{C_n(K;R)}$ refers to the identity map on $C_n(K;R)$.

		\item 
		For all $[\sigma] \in C_n(L;R)$,
		$\boundary_n^L\bigl([\sigma]\bigr) = \boundary_n^K\bigl([\sigma]\bigr)$
			where $\boundary_n^L: C_n(L;R) \to C_{n-1}(L;R)$ refers to the boundary map on $L$ and $\boundary_n^K: C_n(K;R) \to C_{n-1}(K;R)$ that on $K$.
	\end{enumerate}
\end{lemma}
\begin{proof}
	Let $n \in \ints$.
	If $n < 0$, then $C_n(L;R) = 0$, $C_n(K;R) = 0$, and the proposition is trivially satisfied.
	Assume $n \geq 0$.
	Since $i: L \to K$ is an inclusion map, $i(\sigma) = \id_K(\sigma)$ for all $n$-simplices $\sigma \in L$.
	By \fref{prop:chain-group-functorial}, $i_{n,\hash} = \id(C_n(K;R)) = (\id_K)_\hash$. By \fref{prop:induced-maps-on-chains-respect-boundary}, 
	for $[\sigma] \in C_n(L;R)$:
	\begin{equation*}
		(i_{n,\hash} \circ \boundary_n^L)([\sigma])
		= (\boundary_{n}^K \circ i_{n,\hash})([\sigma]) 
		= \boundary_n^K([\sigma])
	\end{equation*}
	as desired.
\end{proof}

	\clearpage



\onlyifstandalone{\setcounter{chapter}{1}}
\chapter{Introduction to Persistence Theory}
\label{chapter:persistence-theory} 

The study of persistence modules or \textit{persistence theory} first came about in the study of persistent homology.
Herbert Edelsbrunner, David Letscher, and Afra Zomorodian in the paper \textit{Topological Persistence and Simplification}~\cite{pershom:edelsbrunner-homsimp} had first defined persistence in terms of homology classes of some collection of spaces.
In particular, persistence theory is interested in the following quantity:
\begin{equation*}
	\rank\!\Big(
		H_n(X_t; R) \Xrightarrow{i_\ast} H_n(X_s; R)
	\Big)
\end{equation*}
where $X_t$ and $X_s$ are two topological spaces such that $X_t \subseteq X_s$ and $R$ is some ring.
Here, the term \textit{persistence} refers to how homology classes in $H_n(X_t; R)$ map to $H_n(X_s; R)$ under the homomorphism $i_\ast$ induced by the inclusion $i: X_t \to X_s$, where we say a homology class $[\sigma] \in H_n(X_t; R)$ \textit{persists} if $i_\ast([\sigma]) \neq 0$.

Persistence theory was later generalized.
For example, 
Afra Zomorodian and Gunnar Carlsson in~\cite{matrixalg:zomorodian} (published in February 2005) then introduced a definition for persistence modules in terms of $R$-modules and considered persistent homology to be a specific example of a persistence module.
In particular, a persistence module was defined to be a collection $\set{M_t}_{t \in \nonnegints}$ of $R$-modules together with a collection of homomorphisms $\set{\phi_t: M_t \to M_{t+1}}_{t \in \nonnegints}$. 
Relative to this characterization, persistence theory studies the following quantity:
\begin{equation*}
	\rank\!\Big(
		M_t \Xrightarrow{\phi_t} M_{t+1}
	\Big).
\end{equation*}
However, definitions for certain constructions involving persistence modules may seem arbitrary or unmotivated given this collection definition. 
For example, the direct sum between persistence modules $\set{M_t}_{t \in \nonnegints}$ and $\set{N_t}_{t \in \nonnegints}$ is defined to be done pointwise as follows:
\begin{equation*}
	\set{M_t}_{t \in \nonnegints} 
	\oplus_\text{Pers} 
	\set{N_t}_{t \in \nonnegints}
	:= 
	\Big\{
		M_t \oplus_\text{Mod} N_t
	\Big\}_{t \in \nonnegints}
\end{equation*}
where $\oplus_\text{Pers}$ refers to the direct sum of persistence modules and $\oplus_\text{Mod}$ to that of $R$-modules.
This notion of direct sum between persistence modules is then related to that between graded $R[x]$-modules, without much justification as to how the definition for $\oplus_\text{Pers}$ is consistent with the categorical definition (so as to be consistent with that of $\oplus_\text{Mod}$).

Later, Bubenik and Scott in the paper \textit{Categorification of Persistence Modules}~\cite{persmod:bubenik-categorification} (presented in May 2014) introduced a new definition for persistence modules in terms of category theory: one that covers the ideas of Zomorodian and Carlsson, provides a more robust foundation to the theory, and is generalizable to concepts such as zigzag persistence and multiparameter persistence.
This new definition in terms of functors is often used in the current literature for persistence theory.

In this chapter,
we discuss persistence theory relative to this characterization of persistence modules as functors,
and present results involving persistence modules relevant to concepts discussed in \cite{matrixalg:zomorodian}.
This chapter is organized as follows:
\begin{enumerate}[chapterdecomposition]
	\item \textbf{Persistence Modules as Functors.}

			We define persistence modules over a field $\field$ as functors of the form $\posetN \to \catvectspace$
			and introduce a type of persistence module called a \textit{finite-type} persistence module.

	\item \textbf{The Category of Persistence Modules.}

			We discuss the category $\catpersmod$ of persistence modules over $\field$, defined to be the category of functors of the form $\posetN \to \catvectspace$,
			We also provide definitions for several algebraic constructions involving persistence modules,
			e.g.\ isomorphisms, direct sums, and chain complexes,
			arising from this characterization as a functor category and briefly discuss how these definitions are consistent with their corresponding categorical definitions.

	\item \textbf{Interval Decompositions of Persistence Modules.}

			We introduce a specific type of direct sum decomposition of a persistence module called an \textit{interval decomposition}, which can be proven to be unique up to persistence isomorphism,
			and discuss the notion of a \textit{persistence barcode}.

	\item \textbf{The Category of Graded Modules over Polynomial Rings}
	
			We discuss the category $\catgradedmod{R}$
			of graded modules over the polynomial ring $R[x]$ with $R$ a PID 
			and review relevant definitions, terminology, and results involving these graded modules.
			We also introduce \textit{graded invariant factor decompositions} for graded $\field[x]$-modules,
				which correspond to invariant factor decompositions that respect the graded structure,
			and present the \textit{Graded Structure Theorem}.

	\item \textbf{The Equivalence between Persistence Modules and Graded Modules}


			We present an isomorphism of categories between $\catpersmod$ and $\catgradedmod{\field}$ 
			and discuss how interval decompositions of persistence modules over $\field$ correspond to graded invariant factor decompositions of graded $\field[x]$-modules.
			We also talk about the correspondence between the algebraic constructions in $\catpersmod$ and those of $\catgradedmod{\field}$ resulting from this isomorphism of categories.
			
\end{enumerate}
The majority of the definitions and results we present below involving persistence modules are taken or adapted from the following papers, listed in increasing order of their initial publication dates:
\begin{enumerate}
	\item \textit{The Categorification of Persistence Modules}~\cite{persmod:bubenik-categorification} by Bubenik and Scott. 

	\item \textit{The Structure and Stability of Persistence Modules}~\cite{persmod:chazal-structure} by Chazal, de Silva, Glisse, and Oudot.

	\item \textit{The Observable Structure of Persistence Modules}~\cite{persmod:chazal-observable} by Chazal, Crawley-Boevey, and de Silva.

	\item \textit{Homological Algebra for Persistence Modules}~\cite{persmod:bubenik-homoloalg} by Bubenik and \Milicevic.
\end{enumerate}
Note that most of the concepts in this chapter are discussed in terms of category theory. 
Some introductory category theory definitions and results are presented in \fref{appendix:cat-theory}. 
For a more detailed treatment of category theory, we recommend reading
\textit{Category Theory in Context} by Emily Riehl \cite{cattheory:rhiel} 
and \textit{Introduction to Homological Algebra} by Joseph Rotman \cite{cattheory:rotman}. 

\clearpage


\section{Persistence Modules as Functors}
\label{section:persistence-modules-as-functors}

\noindent 
We start this section by providing a functor definition for persistence modules, adapted from~\cite[Section 1.3]{persmod:chazal-observable}. 
For reference, $\posetN$ refers to the category induced by the partially ordered set (poset) of nonnegative integers $\nonnegints = \set{n \in \ints: n \geq 0}$ under the $\leq$ relation (see \fref{defn:poset-cat}) and $\catvectspace$ refers to the category of vector spaces over a scalar field $\field$ and linear maps. 
 
\begin{definition}\label{defn:persmod} 
	A \textbf{persistence module} $(\persmod{V}, \alpha_\bullet)$ \textbf{over a field} $\field$ is a functor of the form $\posetN \to \catvectspace$.
	For convenience, we may write $\persmod{V}$ to refer to $(\persmod{V}, \alpha_\bullet)$.
	When the field $\field$ is arbitrary or unambiguous, we may refer to $\persmod{V}$ as a \textbf{persistence module}.  
	
	We identify the following terminology for certain features of persistence modules:
	\begin{enumerate}
		\item 
		For each $t \in \nonnegints$,
			define $V_t := (V_\bullet, \alpha_\bullet)(t)$,
			i.e.\ $V_t$ is the vector space over $\field$ obtained by evaluating the functor $(V_\bullet, \alpha_\bullet)$ on the object $t$ of $\posetN$.
		The vector spaces of $(\persmod{V}, \alpha_\bullet)$ generally refer to the collection $\set{V_t : t \in \nonnegints}$ of vector spaces.

		\item 
		A \textbf{structure map} of $(\persmod{V}, \alpha_\bullet)$ refers to a linear map $\alpha_{s,t}: V_t \to V_s$ by 
		$\alpha_{s,t} := (\persmod{V}, \alpha_\bullet)(t \to s)$ for $t,s \in \nonnegints$ with $t \leq s$,
			i.e.\ $\alpha_{s,t}$ is obtained by evaluating the 	functor $(V_\bullet, \alpha_\bullet)$ on the morphism $t \to s$ of $\posetN$.
		For brevity, we may write 
			$\alpha_t: V_t \to V_{t+1}$ to refer to
			$\alpha_{t} := \alpha_{t+1,t} = \persmod{V}(t \to t+1)$.
		
		\item 
		We call $\posetN$ the \textbf{indexing category} of $\persmod{V}$ and $\nonnegints$ the \textbf{indexing set} of $\persmod{V}$. 
		We may also refer to the index $t \in \nonnegints$ as the \textbf{parameter} or \textbf{scale} of the vector space $V_t$ in the persistence module $\persmod{V}$.
	\end{enumerate}
\end{definition}
\remarks{
	\item 
	The location of the bullet $(\bullet)$ in $(V_\bullet, \alpha_\bullet)$ determine the location of the indices of the vector spaces and structure maps in the notation, 
	e.g.\ if $(W^\bullet, \gamma^\bullet)$ were a persistence module, we denote its vector spaces as $W^t$ and its structure maps by $\gamma^{s,t}: W^t \to W^s$.
	Unlike in the case of chain complexes and cochain complexes, 
	we use the same definition for persistence modules regardless of the location of the bullet $(\bullet)$ (i.e.\ as a subscript or superscript).

	\item 
	In this expository paper, 
	we use an asterisk $(\ast)$ as the ``placeholder'' for the index $n \in \ints$ for a chain complex $C_\ast = (C_n, \boundary_n)_{n \in \ints}$, following \cite{algtopo:rotman}, 
	and we reserve the use of bullets $(\bullet)$ for persistence modules.
	This distinction will be important after \fref{defn:category-of-persistence-chains}, where we introduce chain complexes of persistence modules.
	
	\item 
	For the case of structure maps, the indices $t,s \in \nonnegints$ of $\alpha_{s,t}$ are written in right-to-left order, following the notation for function composition, so that compositions of structure maps are written like $\alpha_{s,r} \circ \alpha_{r,t} = \alpha_{s,t}$ (with equality given later by \fref{lemma:persmod-functor-props}).
	Note that, in this paper, when we say $\alpha_{s,t}$ is a structure map, then it is implied that $t,s \in \nonnegints$ and $t \leq s$.
}

We would like to point out that the definition given in~\cite[Section 1.3]{persmod:chazal-observable} applies to a more general family of persistence modules since it allows for different indexing categories for the domain category $\poset(I, \leq)$ and different categories such as $\catmod{R}$ for the codomain category.
That is, a persistence module $\persmod{V}$ is defined as a functor $\persmod{V}: \poset(I, \leq) \to \catmod{R}$ where $(I, \leq)$ is some partially ordered set and $R$ is some (commutative) ring.
Observe that if $R$ is a field, then an $R$-module is an $R$-vector space by definition, i.e.\ $\catmod{\field}$ and $\catvectspace$ are the same category. 
We bring attention to this since homology is typically introduced using $\ints$ coefficients. Therefore, it may be natural to consider functors of the form $\posetN \to \catmod{\ints}$, so that we can consider a persistence module where each index in $\nonnegints$ corresponds to the $n$\th homology group with coefficients in $\ints$ of some space for some fixed $n \in \nonnegints$.
In this paper, we restrict persistence modules to codomain category $\catvectspace$ since $\ints[x]$ is not a principal ideal domain (PID) but $\field[x]$ for any field $\field$ is, which allows us to use the Graded Structure Theorem (\fref{thm:graded-structure-theorem}) later in \fref{section:cat-equiv-graded-modules}.


We would also like to emphasize that since we introduce persistence theory to discuss persistent homology more concretely, 
we will primary use $\field=\rationals$ and $\field=\ints_p$ for our examples. 
We discuss this choice in more detail later 
in \fref{remark:persistent-hom-why-intsprime} under \fref{section:persistent-homology}.
Below, we provide an example of a persistence module over $\rationals$.

\begin{example}\label{ex:first-persmod-example}
	Let $(\persmod{A}, \alpha_\bullet)$ be a persistence module over $\rationals$ with the vector spaces $A_t$ given as follows:
	\begin{equation*}
		A_t = \begin{cases}
			\rationalsket{a_1, a_2} 		&\text{ if } t = 0	\\
			\rationalsket{b} 				&\text{ if } t = 1 	\\
			\rationalsket{c_1, c_2, c_3}	&\text{ if } t = 2	\\
			\rationalsket{d}				&\text{ if } t \geq 3
		\end{cases}
	\end{equation*}
	and the structure maps $\alpha_{s,t}: A_t \to A_s$ of $\persmod{A}$ are defined as follows:
	\begin{equation*}
		\begin{aligned}[t]
			\alpha_{1,0}: A_0 &\to A_1 \\
				a_1 &\mapsto 2b \\
				a_2 &\mapsto -3b
			\\ \null
			\\[0.5\baselineskip]
			\alpha_{2,0}: A_0 &\to A_2 \\
				a_1 &\mapsto 6c_1 + 4c_2 + 2c_3 \\
				a_2 &\mapsto -9c_1 - 6c_2 - 3c_3
		\end{aligned}
		\qquad
		\begin{aligned}[t]
			\alpha_{2,1}: A_1 &\to A_2 \\
				b &\mapsto 3c_1 + 2c_2 + c_3
			\\ \null
			\\ \null
			\\[0.5\baselineskip]
			\alpha_{t,1}: A_1 &\to A_t \text{ with } t \geq 3 \\
				b &\mapsto 5d
		\end{aligned}
		\qquad\qquad
		\begin{aligned}[t]
			\alpha_{t,2}: A_2 &\to A_3 
				\text{ with } t \geq 3 \\
				c_1 &\mapsto d \\
				c_2 &\mapsto 0 \\
				c_3 &\mapsto 2d
			\\[0.5\baselineskip]
			\alpha_{t,0}: A_0 &\to A_t 
				\text{ with } t \geq 3 \\
				a_1 &\mapsto 10d \\
				a_2 &\mapsto -15d
		\end{aligned}
	\end{equation*}
	For all $t \in \nonnegints$, the structure map $\alpha_{t,t}: A_t \to A_t$ is the identity map on $A_t$. 
	For all $t,s \geq 3$, 
		the structure map $\alpha_{s,t}: A_t \to A_s$ is the identity map on $A_t = \rationalsket{d}$, i.e.\ $d \mapsto d$. 
	
	Observe that the specification for the structure maps above satisfy the functorial property of composition.
	By definition of $\posetN$ (as given in \fref{defn:poset-cat}), there is exactly one morphism $t \to s$ for any $t,s \in \nonnegints$ with $t \leq s$.
	Therefore, the composition $(t \to r) \circ (r \to s)$ of morphisms in $\posetN$ with $t \leq r \leq s$ must be the morphism $t \to s$.
	We talk about this in more detail below in \fref{lemma:persmod-functor-props}.

	For example, any composition of structure maps resulting in a domain of $A_0 = \rationalsket{a_1, a_2}$ and a codomain of $A_3 = \rationalsket{d}$ must exactly be the structure map $\alpha_{3,0}: V_0 \to V_3$ (i.e.\ domain, codomain, and assignments/evaluation must be the same).
	We list some of these compositions below, evaluated on $2a_1 \in A_0$.
	\begin{align*}
		\alpha_{3,0}(2a_1) 
			&= 2(10d) 
			&= 20d 		
		\\
		(\alpha_{3,2} \circ \alpha_{2,0})(2a_1) 
			&= \alpha_{3,2}(
				12c_1 + 8c_2 + 4c_3
			) = 12d + 8d 
			&= 20d 
		\\
		(\alpha_{3,1} \circ \alpha_{1,0})(2a_1)
			&= \alpha_{3,1}(4b)
			&= 20d
		\\
		(\alpha_{3,2} \circ \alpha_{2,1} \circ \alpha_{1,0})(2a_1)
			&= (\alpha_{3,2} \circ \alpha_{2,1})(2b)
			= \alpha_{3,2}(6c_1 + 4c_2 + c_3)
			= 12d + 8d 
			&= 20d
		\\
		(\alpha_{3,0} \circ \alpha_{0,0})(2a_1)
			&= \alpha_{3,0}(2a_1) 
			&= 20d
		\\
		(\alpha_{3,2} \circ \alpha_{2,2} \circ \alpha_{2,0})(2a_1) 
			&= (\alpha_{3,2} \circ \alpha_{2,2})(
				12c_1 + 8c_2 + 4c_3
			) 
			= \alpha_{3,2}(
				12c_1 + 8c_2 + 4c_3
			)
			= 12d + 8d 
			&= 20d 
	\end{align*}
\end{example}

Defining persistence modules as functors is a succinct way to impose additional algebraic structure to the construction. 
As mentioned in the introduction, 
\cite[Definition 3.2]{matrixalg:zomorodian} defines a persistence module as a collection of $R$-modules $\set{V^t}$ and homomorphisms $\set{\phi^t: V_t \to V_{t+1}}$, both indexed by $t \in \nonnegints$. 
Observe that the collection definition does not explicitly state, for example, that there exists a linear map $V_2 \to V_4$.
While it may be somewhat natural to assume the linear map $\alpha_{4,2}: V_2 \to V_4$ is given by the composition $\alpha_{4,2} = \alpha_{3} \circ \alpha_2$, 
having to state it as a separate condition can be somewhat confusing or cumbersome and can also feel arbitrary.
In comparison, by assuming that said collection corresponds to a functor, we get a number of properties as a consequence of the functor definition. We list some below.

\begin{lemma}\label{lemma:persmod-functor-props}
	Let $(\persmod{V}, \alpha_\bullet)$ be a persistence module. 
	\begin{enumerate}
		\item For all $t \in \nonnegints$, the structure map $\alpha_{t,t}: V_t \to V_t$ is the identity map on $V_t$.

		\item For all $t,r,s \in \nonnegints$ with $t \leq r \leq s$, $\alpha_{s,t} = \alpha_{s,r} \circ \alpha_{r,t}$. 
		That is, the following diagram commutes:
		\begin{center}
			\includegraphics[height=0.7in]{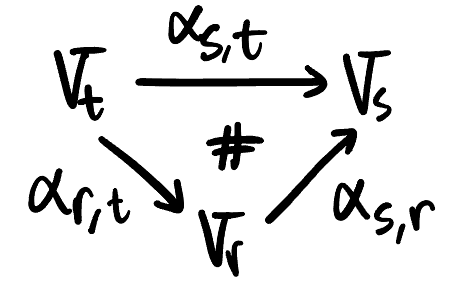}
		\end{center} 

		\item For all $t,s \in \nonnegints$ with $t \leq s$, there is exactly one linear map $V_t \to V_s$. 
	\end{enumerate}
\end{lemma}
\begin{proof}
	Let $n \in \nonnegints$. By definition of $\posetN$, the unique morphism $n \to n$ must be the identity morphism on $n$. By definition, a functor must pass $\persmod{V}(n \to n)$ to the identity map on $V_n$. Therefore, the structure map $\persmod{V}(n \to n) = \alpha_{n,n} = \id_{V_n}$.

	Let $t,s,r \in \nonnegints$ such that $t \leq r \leq s$. By definition of $\posetN$, we have the equality $(s \leftarrow r) \circ (r \leftarrow t) = (s \leftarrow t)$ as morphisms of $\posetN$. Since functors must respect composition (indicated below by $\star$), we have 
	\begin{equation*}
		\alpha_{s,t} =
		\persmod{V}(s \leftarrow t)
		=
		\persmod{V}\big(
			(s \leftarrow r) \circ (r \leftarrow t)
		\big) 
		\overset{\star}{=} 
		\persmod{V}(s \leftarrow r) \circ \persmod{V}(r \leftarrow t)
		= 
		\alpha_{s,r} \circ \alpha_{r,t}
	\end{equation*}
	Claim (\textbf{iii}) is a re-statement of Claim (\textbf{ii}), wherein all linear maps with domain $V_t$ and domain $V_s$ must be equal.
\end{proof}

\noindent 
Note that \fref{lemma:persmod-functor-props} allows us to represent persistence modules as sequences and use sequences to refer to a persistence module unambiguously. 
In particular, we can describe a persistence module $\persmod{V}$ as the following sequence:
\begin{equation*}
	V_0 \Xrightarrow{\alpha_0} V_1
		\Xrightarrow{\alpha_1} V_2
		\Xrightarrow{\alpha_2} V_3
		\Xrightarrow{\alpha_3} V_4
		\Xrightarrow{\alpha_4} \cdots
\end{equation*}
with $\alpha_t$ being the linear map $\alpha_t: V_t \to V_{t+1}$ for all $t \in \nonnegints$. 
As a sidenote, this sequence representation applies to any diagram
with $\posetN$ as the domain category.
Note that we use the term \textit{diagram} to refer to functors of the form $\poset(I, \leq) \to \catname{C}$, as described in~\fref{defn:cat-diagrams}.


Conversely, \fref{lemma:persmod-functor-props} also implies that there is enough information in the sequence (i.e.\ the collection of vector spaces $\set{V_n}$ and morphisms $\set{\alpha_{t}: V_t \to V_{t+1}}$) to determine a persistence module $\persmod{V}$. 
In particular, for a given persistence module $\persmod{V}$, the structure maps of the form $\alpha_{s,t}: V_t \to V_s$ with $s \neq t+1$ are uniquely determined by the set of $\set{\alpha_t}_{t \in \nonnegints}$ of linear maps.
Therefore, the collection definition given in~\cite[Definition 3.2]{matrixalg:zomorodian} is compatible with the functor definition given in~\fref{defn:persmod}. 
We state this as a lemma below.

\begin{proposition}\label{prop:persmod-collection}
	Let $\set{W_t}_{t \in \nonnegints}$ and $\set{\phi: W_t \to W_{t+1}}_{t \in \nonnegints}$ be $\nonnegints$-indexed collections of $\field$-vector spaces and linear maps respectively.
	Let $(\persmod{V}, \alpha_\bullet)$ be a persistence module over $\field$ such that $V_t = W_t$ for all $t \in \nonnegints$ and the structure maps of $(\persmod{V}, \alpha_\bullet)$ of the form $\alpha_{t}: V_t \to V_{t+1}$ are given by $\alpha_{t} = \phi_t$.
	Then, there is exactly one choice for the remaining structure maps of $(\persmod{V}, \alpha_\bullet)$, as listed below:
	\begin{enumerate}
		\item For all $t \in \nonnegints$, the structure map $\alpha_{t,t}: V_t \to V_t$ must be the identity map $\id_{V_t}: V_t \to V_t$ of $V_t$.

		\item 
		For all $t \in \nonnegints$ and $k \geq 2$,
			the structure map $\alpha_{t+k,t}: V_t \to V_{t+k}$
			is given by the following composition:
		\begin{equation*}
			\alpha_{s,t} = \alpha_{s-1} \circ \alpha_{s-2} \circ \cdots 
				\circ \alpha_{t+1} \circ \alpha_{t}
		\end{equation*}
	\end{enumerate}
	That is, the pair of collections $\set{W_t}$ and $\set{\phi_t}$ unambiguously determines a persistence module $(W_\bullet, \phi_\bullet)$.
\end{proposition}
\begin{proof}
	This follows immediately from the definition of $\posetN$ and that of functors. The proof is similar to that of  \fref{lemma:persmod-functor-props}.
\end{proof}

\noindent 
Note that by \fref{lemma:persmod-functor-props} and \fref{prop:persmod-collection},
we get a correspondence between the set of persistence modules $(\persmod{V}, \alpha_\bullet)$ and the set of pairs of collections $(\set{V_t}, \set{\alpha_t: V_t \to V_{t+1}})$.
For convenience, we often define persistence modules using \fref{prop:persmod-collection}.
As an example, we describe the persistence module $(\persmod{A}, \alpha_\bullet)$ given in~\fref{ex:first-persmod-example} using a collection of vector spaces and of linear maps. We also describe the sequence corresponding to $(\persmod{A}, \alpha_\bullet)$.

\begin{example}
	Let $(\persmod{B}, \beta_\bullet)$ be a persistence module over $\rationals$ with vector spaces $B_t$ given as follows:
	\begin{equation*}
		B_t = \begin{cases}
			\rationalsket{a_1, a_2} 		&\text{ if } t = 0	\\
			\rationalsket{b} 				&\text{ if } t = 1 	\\
			\rationalsket{c_1, c_2, c_3}	&\text{ if } t = 2	\\
			\rationalsket{d}				&\text{ if } t \geq 3
		\end{cases}
	\end{equation*}
	Let the structure maps $\beta_{t}: B_t \to B_{t+1}$ of $\persmod{B}$ be as given below:
	\begin{equation*}
		\begin{aligned}[t]
			\beta_0: B_0 &\to B_1 	\\
					 	a_1 &\mapsto 2b 	\\
					 	a_2 &\mapsto -3b
		\end{aligned}
		\qquad\quad
		\begin{aligned}[t]
			\beta_1: B_1 &\to B_2 	\\
					 	b 	&\mapsto 3c_1 + 2c_2 + c_3
		\end{aligned}
		\qquad\quad
		\begin{aligned}[t]
			\beta_2: B_2 &\to B_3 	\\
						c_1 &\mapsto d \\
						c_2 &\mapsto 0 \\
						c_3 &\mapsto 2d
		\end{aligned}
		\qquad\quad
		\begin{aligned}[t]
			\beta_t: B_t &\to B_{t+1} \text{ for } t \geq 3 \\
						d &\mapsto d
		\end{aligned}
	\end{equation*}
	Use \fref{prop:persmod-collection} to define the remaining structure maps of $\persmod{B}$.
	Then, $\persmod{B} = \persmod{A}$ (as functors) where $\persmod{A}$ is the persistence module defined in~\fref{ex:first-persmod-example}.
	Observe that the specification for $\persmod{B}$ is much shorter that that for $\persmod{A}$.
	Also, $\persmod{A}$ corresponds to the following sequence:
	\begin{equation*}
		\equalsupto{A_0}{\rationalsket{a_1, a_2}} 
			\Xrightarrow{\,\alpha_0\,}
		\equalsupto{A_1}{\rationalsket{b}} 
			\Xrightarrow{\,\alpha_1\,}
		\equalsupto{A_2}{\rationalsket{c_1, c_2, c_3}} 
			\Xrightarrow{\,\alpha_2\,}
		\equalsupto{A_3}{\rationalsket{d}}  
			\Xrightarrow{\,\alpha_3\,}
		\equalsupto{A_4}{\rationalsket{d}}  
			\Xrightarrow{\,\alpha_4\,}
			\cdots
	\end{equation*}
\end{example}


Observe that~\fref{defn:persmod} allows for persistence modules that have non-finite characteristics.
For example, some of the vector spaces $V_t$ of a persistence module $\persmod{V}$ may be infinite-dimensional. 
We can also have a persistence module $(\persmod{W}, \gamma_\bullet)$ wherein none of the structure maps $\gamma_{t}: W_t \to W_{t+1}$ are isomorphisms or can be made into isomorphisms by a change of basis on $W_t$. 
In this case, we are considering an infinite number of vector spaces.
Persistence modules of this nature do present some problems.
Indeed, the results presented in~\cite{matrixalg:zomorodian} only consider a specific type of persistence module, one that is finite in specific ways.
We use a definition of such, adapted from~\cite[Definition 3.3]{matrixalg:zomorodian}.

\begin{statement}{Definition}\label{defn:persmod-constant-on-interval}
	A persistence module $(\persmod{V}, \alpha_\bullet)$ is \textbf{constant on an interval $I \subseteq \nonnegints$} if
	for all $t,s \in I$ with $t \leq s$, 
		the structure map $\alpha_{s,t}: V_t \to V_s$ is 
		a vector space isomorphism.
	A persistence module $\persmod{V}$ is called \textbf{finite-type} 
	if all of its vector spaces $V_t$ are finite-dimensional and $\persmod{V}$ is constant on $[N, \infty)$ for some $N \in \nonnegints$.
\end{statement}
\begin{miniremark}
	Since persistence theory is a relatively new field, 
	there exist terms that are commonly used in the literature but are defined differently depending on the author(s).
	The term \textit{finite-type} is one of these.

	Later in \fref{prop:finite-type-has-interval},
		using the definition above for \textit{finite-type}, 
		we state that interval decompositions 
		(defined in \fref{defn:interval-decomposition})
		exist for finite-type persistence modules.
	In contrast, 
		\cite[Definition 4.1]{persmod:bubenik-categorification} 
		and 
		\cite{persmod:chazal-structure}
	define the term \textit{finite-type} to refer to persistence modules for which an interval decomposition exists.
	The term \textit{tame} in \cite{persmod:bubenik-categorification} serves the same function as \textit{finite-type} in \cite{matrixalg:zomorodian}, i.e.\ it serves as a characterization of persistence modules that have interval decompositions. 
	However, the term \textit{tame} also has its problems, as discussed in \cite[p5]{persmod:chazal-structure}.
\end{miniremark}


All of the examples we will present in this expository involve finite-type persistence modules (not including the one example given below for comparison).
As we will discuss in \fref{section:persistent-homology}, under the assumption that the simplicial complex in question is finite, the persistence modules generated by said simplicial complex will necessarily be of finite-type. 
This assumption is justified since most practical applications of persistent homology, particularly those involving calculation, 
have simplicial complexes generated from finite datasets.
The usual construction involves letting the dataset $V$ be the vertex set of the to-be constructed simplicial complex $K$. Since $K$ is a subset of the power set $2^V$ and power sets of finite sets are finite, $K$ must also be finite.
This explains why some of the introductory literature on persistent homology describe persistence modules as finite sequences, 
i.e.\ persistence modules are described to be finite sequences as given below:
\begin{equation*}
	{V_0} 
		\Xrightarrow{\alpha_0} V_1 
		\Xrightarrow{\alpha_1} V_2 
		\Xrightarrow{\alpha_2} 
		\cdots 
		\Xrightarrow{\alpha_{N-1}} 
	V_N
\end{equation*}
Using our definition, we would interpret the sequence above to correspond to a persistence module $\persmod{V}$ that is constant on $[N, \infty)$ and assume that $V_n = V_N$ for all $n \in [N+1, \infty)$.
As an example, we look at the persistence module given earlier in~\fref{ex:first-persmod-example}.

\begin{example}
	Let $(\persmod{A},\alpha_\bullet)$ be as defined in~\fref{ex:first-persmod-example}.
	Then, $\persmod{A}$ is a finite-type persistence module and is constant on $[3,\infty)$. 
	Also, $\persmod{A}$ corresponds to the following sequence:
	\begin{equation*}
		\equalsupto{A_0}{\rationalsket{a_1, a_2}} 
			\Xrightarrow{\alpha_0}
		\equalsupto{A_1}{\rationalsket{b}} 
			\Xrightarrow{\alpha_1}
		\equalsupto{A_2}{\rationalsket{c_1, c_2, c_3}} 
			\Xrightarrow{\alpha_2}
		\equalsupto{A_3}{\rationalsket{d}}  
	\end{equation*}
\end{example}

\noindent
For comparison, we provide an example below of a persistence module that is not finite-type. 

\begin{example}
	Let $A = \set{a_n}_{n \in \nonnegints}$ be a set of indeterminates. 
	Define the persistence module $(\persmod{F}, \eta_\bullet)$ over $\rationals$ using \fref{prop:persmod-collection} as follows:
	\begin{enumerate}
		\item 
		For each $t \in \nonnegints$, $F_t := \rationalsket{a_0, a_1, \ldots, a_t}$.

		\item 
		For each $t \in \nonnegints$, define the structure map $\eta_t: F_t \to F_{t+1}$ by $\eta_t(a_i) = a_i$ for all $i \in \set{0, \ldots, t}$.
	\end{enumerate}
	Then, the following sequence represents the persistence module $(\persmod{F}, \eta_\bullet)$:
	\begin{equation*}
		\equalsupto{F_0}{\rationalsket{a_0}}
			\Xrightarrow[a_0 \mapsto a_0]{\eta_0}
		\equalsupto{F_1}{\rationalsket{a_0, a_1}}
			\Xrightarrow[
				a_i \mapsto a_i, i=0,1
			]{\eta_1}
		\equalsupto{F_2}{\rationalsket{a_0, a_1, a_2}}
			\Xrightarrow[
				a_i \mapsto a_i, i=0,1,2
			]{\eta_2}
		\equalsupto{F_3}{\rationalsket{a_0, a_1, a_2, a_3}}
			\Xrightarrow{\eta_3}
			\cdots
	\end{equation*}
	While all vector spaces $F_t$ are finite-dimensional, $\persmod{V}$ fails to be constant on the interval $[N, \infty)$ for any $N \in \nonnegints$.
	We identify two ways we can see this:
	\begin{enumerate}
		\item For all $t \in \nonnegints$, $\eta_t: F_t \to F_{t+1}$ cannot be an isomorphism since $\eta_t\inv(\set{a_{t+1}}) = \emptyset$, i.e.\ $\eta_t$ is not surjective.
		\item Any two vector spaces $F_t$ and $F_s$ of $\persmod{V}$ with $t \neq s$ cannot be isomorphic since $\dim(F_t) = t+1 \neq s+1 = \dim(F_s)$.
	\end{enumerate}
	Therefore, $\persmod{F}$ is \textbf{not} a finite-type persistence module.
\end{example}

 \clearpage

\section{The Category of Persistence Modules}
\label{section:cat-persistence-modules}

One of the key goals of persistence theory is to characterize persistence modules that admit a special decomposition called an interval decomposition.
Note that we use the term \textit{(direct sum) decomposition} to refer to an isomorphism to a direct sum, similar to the case of vector spaces, modules, and of chain complexes.
However, to understand decompositions, we first need to define the notions of isomorphisms and of direct sums between persistence modules.

In persistence theory, these algebraic constructions are defined by forming the category $\catpersmod$ of persistence modules over $\field$ as a \textit{functor category} with \textit{natural transformations} as morphisms.
Below, we provide definitions specific for the category of persistence modules, {adapted} from \cite[Definition 1.3]{persmod:chazal-observable}.

\begin{statement}{Definition}\label{defn:persmod-cat}
	The \textbf{category $\catpersmod$ of persistence modules over a field $\field$} 
	is the category of functors of the form $\posetN \to \catvectspace$ and natural transformations, i.e.\ $\catpersmod$ consists of the following:
	\begin{enumerate}
		\item The objects in $\catpersmod$ are persistence modules over $\field$, as given in \fref{defn:persmod}.

		\item 
		The morphisms in $\catpersmod$ are persistence morphisms, defined as follows:

		A \textbf{persistence morphism} $\phi_\bullet: (\persmod{V},\alpha_\bullet) \to (\persmod{W},\gamma_\bullet)$ between persistence modules 
		$(\persmod{V},\alpha_\bullet)$ and $(\persmod{W}, \gamma_\bullet)$ over $\field$ 
		is a collection of linear maps $\phi_\bullet = (\phi_t: V_t \to W_t)_{t \in \nonnegints}$ such that for all $t,s \in \nonnegints$ with $t \leq s$, 
		the composition relation 
		$\gamma_{s,t} \circ \phi_t = \phi_s \circ \alpha_{s,t}$
		is satisfied,
		i.e.\ 
		the following diagram commutes:
		\vspace{-5pt}
		\begin{displaymath}
		\begin{tikzcd}
			V_t \arrow[r, "\alpha_{s,t}"] 
						\arrow[d, "\phi_t" swap]
				& V_s
						\arrow[d, "\phi_s"]
			\\
			W_t \arrow[r, "\gamma_{s,t}"] & W_s
		\end{tikzcd}\vspace{-0.5\baselineskip}
		\end{displaymath} 
	\end{enumerate}
\end{statement}
\remarks{
	\item 
	In this paper, writing the symbol $\phi_\bullet$ to represent a persistence morphism
		determines that $\phi_t$ refers to the corresponding linear map $\phi_t: V_t \to W_t$ for each $t \in \nonnegints$,
		i.e.\ the bullet $(\bullet)$ in $\phi_\bullet$ is replaced by the index $t \in \nonnegints$, much like in the case of persistence modules as remarked under \fref{defn:persmod}.

	\item 
	By \cite[Section 1.7, p44]{cattheory:rhiel} and \cite[Example 1.19(i), p27]{cattheory:rotman},
	for any two categories $\catname{C}$ and $\catname{D}$, 
	there exists a category
	denoted $\catname{D}^\catname{C}$ 
	consisting of functors of the form $\catname{C} \to \catname{D}$ as objects and natural transformations between functors. 
	This category is called a \textit{functor category}.
	We use this result for the claim of $\catpersmod$ being a well-defined category.

	As a sidenote, \cite{persmod:bubenik-categorification} denotes $\catpersmod$ as $\catvectspace^{\raisebox{-2.5pt}{$\scriptstyle (\nonnegints, \leq)$}}$ where $(\nonnegints, \leq)$ refers to the poset category $\posetN$ (see remarks under \fref{defn:poset-cat} for the notation of poset categories).
}


We provide an example of a persistence morphism below. Note that not all collections of linear maps $\set{\phi_t: V_t \to W_t}_{t \in \nonnegints}$ produce a persistence morphism $\phi_\bullet: \persmod{V} \to \persmod{W}$. 
We need to check that the composition relation, i.e.\ the commuting squares and diagrams, is satisfied.
\begin{example}\label{ex:first-pers-morphism}
	Let $(\persmod{V}, \alpha_\bullet)$ and 
	$(\persmod{W}, \gamma_\bullet)$ be persistence modules over $\rationals$ 
	with the following $\rationals$-vector spaces defined over indeterminates $A,B,C$ for $V_t$ and $a,b,c$ for $W_t$:
	\begin{equation*}
		V_t = \begin{cases}
			0 &\text{ if } t=0 \\
			\rationals\ket{A} 
				&\text{ if } t=1 \\
			\rationals\ket{A,B}
				&\text{ if } t=2 \\
			\rationals\ket{A,B,C}
				&\text{ if } t \geq 3
		\end{cases}
		\qquad\text{ and }\qquad 
		W_t = \begin{cases}
			\rationals\ket{a,b} 
				&\text{ if } t=0 \\
			\rationals\ket{a,b,c} 
				&\text{ if } t \geq 1
		\end{cases}
	\end{equation*}
	The structure maps $\alpha_t: V_t \to V_{t+1}$ of $\persmod{V}$ 
	and $\gamma_t: W_t \to W_{t+1}$ of $\persmod{W}$ are determined by the identity maps on 
		$\rationals\ket{A,B,C}$ and 
		on $\rationals\ket{a,b,c}$ respectively, i.e.\ 
	$\alpha_t(X) = \id_{\rationals\ket{A,B,C}}(X)$ for all $X \in V_t$ and 
	$\gamma_t(x) = \id_{\rationals\ket{a,b,c}}(x)$ for all $x \in W_t$.

	Define a family $\set{\phi_t}_{t \in \nonnegints}$ of linear maps $\phi_t: V_t \to W_t$ as follows:
	Let $\Phi: \rationals\ket{A,B,C} \to \rationals\ket{a,b,c}$ be given by
		$A \mapsto b-a$,
		$B \mapsto a-c$,
		and $C \mapsto c-b$.
	Define $\phi_t: V_t \to W_t$ by 
		$\phi_t(X) := \Phi(X)$ for all $X \in V_t$,
		i.e.\ we have the following assignments for $\phi_t$:
	\begin{equation*}\setlength{\arraycolsep}{0.5\arraycolsep}
	\begin{array}{rclcl}
		\phi_t(A) &=& b-a &\text{ for }& t \geq 1 \\
		\phi_t(B) &=& a-c &\text{ for }& t \geq 2 \\
		\phi_t(C) &=& c-b &\text{ for }& t \geq 3
	\end{array}
	\end{equation*}
	Note that $\phi_0: V_0 \to W_0$ is necessarily the trivial map since the domain $V_0$ is the trivial vector space.
	We can confirm that the family $\set{\phi_t}_{t \in \nonnegints}$ satisfies the commutativity requirement as follows:
	for all $X \in V_t$ with $t \geq 1$:
	\begin{equation*}
		(\phi_s \circ \alpha_{s,t})(X)
		= \Bigl( \Phi \circ \id_{\rationals\ket{A,B,C}} \Bigr)(X)
		= \Phi(X) 
		= \Bigl( \id_{\rationals\ket{a,b,c}} \mathrel{\circ} \Phi \Bigr)(X)
		= (\gamma_{s,t} \circ \phi_t)(X) 
	\end{equation*}
	This is illustrated in the diagram below:
	\vspace{-5pt} 
	\begin{displaymath}
	\begin{tikzcd}[row sep=normal, column sep=large]
		{} &[-25pt] \color{black}\scriptstyle (t=0) 
		& \color{black}\scriptstyle (t=1) 
		& \color{black}\scriptstyle (t=2) 
		& \color{black}\scriptstyle (t=3) 
		\\[-20pt]
		(V_\bullet, \alpha_\bullet) :
		& 0 
			\arrow[r,"\alpha_0"]
			\arrow[d, "\phi_0" swap]
		& \rationals\ket{A} 
			\arrow[r, "\alpha_1"]
			\arrow[d, "\phi_1" swap,
			"\color{black}\substack{
				A \,\mapsto\, b-a
			}"]
		& \rationals\ket{A,B} 
			\arrow[r, "\alpha_2"]
			\arrow[d, "\phi_2" swap,
			"\color{black}\substack{
				A \,\mapsto\, b-a \\[1pt]
				B \,\mapsto\, a-c
			}"]
		& \rationals\ket{A,B,C}
			\arrow[r, "\alpha_3"]
			\arrow[d, "\phi_3" swap,
			"\color{black}\substack{
				A \,\mapsto\, b-a \\[1pt]
				B \,\mapsto\, a-c \\[1pt]
				C \,\mapsto\, c-b
			}"]
		& \cdots
		\\[15pt] 
		(W_\bullet, \gamma_\bullet) :
		& \rationals\ket{a,b} 
			\arrow[r, "\gamma_0" swap]
		& \rationals\ket{a,b,c} 
			\arrow[r, "\gamma_1" swap]
		& \rationals\ket{a,b,c} 
			\arrow[r, "\gamma_2" swap]
		& \rationals\ket{a,b,c}
			\arrow[r, "\gamma_3" swap]
		& \cdots
	\end{tikzcd}
	\end{displaymath}  
	Therefore, the set $\set{\phi_t}_{t \in \nonnegints}$ determines a persistence morphism $\phi_\bullet: (V_\bullet, \alpha_\bullet) \to (W_\bullet, \gamma_\bullet)$.
\end{example} 

Observe that, given two persistence morphisms $\phi_\bullet: (A_\bullet, \alpha_\bullet) \to (B_\bullet, \beta_\bullet)$ and 
$\psi_\bullet: (B_\bullet, \beta_\bullet) \to (C_\bullet, \gamma_\bullet)$
between persistence modules 
$(A_\bullet, \alpha_\bullet)$,
$(B_\bullet, \beta_\bullet)$,
and $(C_\bullet, \gamma_\bullet)$,
	the collection of compositions 
	$\set{\psi_t \circ \phi_t}_{t \in \nonnegints}$
	also determines a persistence morphism
	${\psi_\bullet} \circ {\phi_\bullet}: A_\bullet \to C_\bullet$.
Relative to the commuting squares, 
	the squares of $\phi_\bullet$ and $\psi_\bullet$ together produce the squares for the composition, as illustrated below:
\begin{displaymath}
	\begin{array}{c c c c c}
		\mathclap{\color{black}\left(\,\substack{
			\text{commuting square of} \\[2pt]
			\phi_\bullet \,:\, (A_\bullet, \alpha_\bullet) \to (B_\bullet,\beta_\bullet)
		}\,\right)}
		&& 
		\mathclap{\color{black}\left(\,\substack{
			\text{commuting square of} \\[2pt]
			\psi_\bullet \,:\, (B_\bullet, \beta_\bullet) \to (C_\bullet,\gamma_\bullet)
		}\,\right)}
		&&
		\mathclap{\color{black}\left(\,\substack{
			\text{commuting square for} \\[2pt]
			\text{the composition }\psi_\bullet \circ\, \phi_\bullet
		}\,\right)}
		\\[5pt]
		\begin{tikzcd}
			A_t \arrow[r, "\alpha_{s,t}"]
				\arrow[d, "\phi_t" swap]
				\arrow[rd, "\hash" description,draw=none]
			& A_s 
				\arrow[d, "\phi_s" ]
			\\
			B_t \arrow[r, "\beta_{s,t}" swap]
			& B_s
		\end{tikzcd}
		& 
		\text{ with }
		&
		\begin{tikzcd}[row sep=normal, column sep=normal]
			B_t \arrow[r, "\beta_{s,t}"]
				\arrow[d, "\psi_t" swap]
				\arrow[rd, "\hash" description,draw=none]
			& B_s 
				\arrow[d, "\psi_s" ]
			\\
			C_t \arrow[r, "\gamma_{s,t}" swap]
			& C_s
		\end{tikzcd}
		&
		\text{ implies } 
		& 
		\begin{tikzcd}[row sep=normal, column sep=normal]
			A_t \arrow[r, "\alpha_{s,t}"]
				\arrow[d, "\psi_t \circ \phi_s" swap]
				\arrow[rd, "\hash" description,draw=none]
			& A_s 
				\arrow[d, "\psi_s \circ \phi_s" ]
			\\
			C_t \arrow[r, "\gamma_{s,t}" swap]
			& C_s
		\end{tikzcd}
	\end{array}
\end{displaymath}
Therefore, the composition $\psi_\bullet \mathrel{\circ} \phi_\bullet: (A_\bullet, \alpha_\bullet) \to (C_\bullet, \gamma_\bullet)$ of persistence morphisms is also a persistence morphism.
Note that this is part of the proof of $\catpersmod$ being a well-defined category.

\spacer 

Since the category $\catpersmod$ of persistence modules is defined to be a functor category of functors of the form $\posetN \to \catvectspace$,
	there is a natural extension of the algebraic constructions in $\catvectspace$ to the case of $\catpersmod$.
In particular, we can define constructions in $\catpersmod$ by pointwise evaluation for each $t \in \nonnegints$ using the corresponding constructions in $\catvectspace$.
Note that we need to be careful here since the constructions in $\catpersmod$ have to respect the structure of the morphisms.
We discuss this in more detail for the following constructions in $\catpersmod$:
	isomorphism relations, direct sum operations, subobject relations, kernel, image, and cokernel of morphisms, and chain complexes.

We want to point out that this extension process (for lack of a better term) is not unique to persistence modules.
More generally, this seems to be a common property of functor categories with functors of the form  $\poset{(I,\leq)} \to \catname{A}$ where $\catname{A}$ is some \textit{abelian} category.
Note that $\catmod{R}$ over a PID $R$ and $\catvectspace$ over a field $\field$ are examples of abelian categories.
Abelian categories are outside the scope of this paper but for those interested, this is discussed in more detail in \cite[Section 5.5: Proposition 5.93 and Corollary 5.94]{cattheory:rotman}
and in \cite[Appendix A.4]{hom-algebra:weibel}.
This might explain why most of the literature we have read on persistence theory do not explicitly describe these constructions, 
	e.g.\ \cite{persmod:chazal-observable,persmod:bubenik-homoloalg} wherein a reference to $\catpersmod$ being an \textit{abelian} or \textit{Grothendieck} category seems to suffice.

Note that an equality between persistence modules in $\catpersmod$ corresponds to an equality between functors.
That is, two persistence modules $(\persmod{V}, \alpha_\bullet)$ and $(\persmod{W},\gamma_\bullet)$ are equal if 
$V_t = W_t$ as $\field$-vector spaces for all $t \in \nonnegints$ and $\alpha_{s,t} = \gamma_{s,t}$ as linear maps for all $t,s \in \nonnegints$ with $t \leq s$. 
We start with a definition for isomorphisms between persistence modules.

\begin{definition}\label{defn:persistence-isomorphism}
	A \textbf{persistence isomorphism} $\phi_\bullet: V_\bullet \to W_\bullet$ between persistence modules $(V_\bullet, \alpha_\bullet)$ and $(W_\bullet, \gamma_\bullet)$ over $\field$ is a persistence morphism $\phi_\bullet: V_\bullet \to W_\bullet$ with $\phi_\bullet = (\phi_t: V_t \to W_t)_{t \in \nonnegints}$
	such that for all $t \in \nonnegints$,
		the linear map $\phi_t: V_t \to W_t$ is an $\field$-vector space isomorphism.
	If such a persistence isomorphism exists,
		we say that $V_\bullet$ and $W_\bullet$ are \textbf{isomorphic} (as persistence modules) and write 
		$V_\bullet \cong W_\bullet$.
\end{definition}
\remark{
	If needed, we may write $\uppersmod\cong$ instead of $\cong$ to emphasize that the isomorphism relation is of the category $\catpersmod$.
	We use this notation sometimes in \fref{chapter:matrix-calculation} where we talk about isomorphisms between $R$-modules, between graded $\field[x]$-modules, and between persistence modules in the same context.
}
\HIDE{\begin{miniremark}
	Technically speaking, the definition for isomorphism above is a consequence of the definition for persistence morphism.
	In general, a morphism $\phi: C \to D$ in a category $\catname{C}$ is an \textit{isomorphism} if there exists another morphism $\eta: D \to C$ such that $\eta \circ \phi = $ and $\phi \circ \eta = \mathbf{1}_D$
\end{miniremark}}

We want to emphasize that a persistence isomorphism has to be given by a persistence morphism and that the structure maps of the persistence modules cannot be ignored.
That is, given two persistence modules $(V_\bullet, \alpha_\bullet)$ and $(W_\bullet, \gamma_\bullet)$,
a collection $\set{\phi_t}_{t \in \nonnegints}$ of $\field$-vector space isomorphisms $\phi_t: V_t \to W_t$ generally does not make a persistence isomorphism.
Below, we provide an example of a collection $\set{\phi_t}_{t \in \nonnegints}$ of vector spaces isomorphisms that is not a persistence morphism, and another collection $\set{\psi_t}_{t \in \nonnegints}$ that is a persistence isomorphism.

\begin{example}
	Define the persistence modules $(V_\bullet, \alpha_\bullet)$ 
	and $(W_\bullet, \gamma_\bullet)$ over $\rationals$ as follows:
	For all $t \in \nonnegints$, define $V_t := \rationals\ket{a,b}$ and $W_t := \rationals\ket{x,y}$ with indeterminates $a$, $b$, $x$, and $y$.
	For $t=0$, define the structure map $\alpha_0: V_0 \to V_1$ by $\alpha_0(a) := a+b$ and $\alpha_0(b) := b$.
	For $t \geq 1$, let $\alpha_t: V_t \to V_{t+1}$ be the identity map on $\rationals\ket{a,b}$.
	For all $t \in \nonnegints$, let $\gamma_t : W_t \to W_{t+1}$ be the identity map on $\rationals\ket{x,y}$.
	\begin{enumerate}[label={Part\,(\alph*).}, left=0.2in]
		\item 
		For each $t \in \nonnegints$, let $\phi_t: V_t \to W_t$ be given by $a \mapsto x$ and $b \mapsto y$.
		Since $V_t = \rationals\ket{a,b}$ and $W_t = \rationals\ket{x,y}$ for all $t \in \nonnegints$, each $\phi_t$ is a vector space isomorphism with obvious inverse.
		However, the collection $\set{\phi_t}_{t \in \nonnegints}$ is \textbf{not} a persistence morphism since 
		\begin{equation*}
			x+y 
			= \bigl( \phi_1 \circ \alpha_{0} \bigr)(a)
			= \bigl( \phi_1 \circ \alpha_{1,0} \bigr)(a) 
			\neq 
			\bigl( \gamma_{1,0} \circ \phi_0 \bigr)(a) 
			= \bigl( \gamma_{0} \circ \phi_0 \bigr)(a) 
			= x
		\end{equation*}
		That is, the following diagram does not commute:
		\vspace{-5pt}
		\begin{displaymath}
		\begin{tikzcd}[column sep=huge, row sep=large]
			\mathllap{V_0 =\ } \rationals\ket{a,b} 
				\arrow[r, "\alpha_{1,0}", "\color{black}\substack{
					a \,\mapsto\, a+b ,\, b \,\mapsto\, b
				}" swap]
				\arrow[d, "\phi_0", "\color{black}\substack{
					a \,\mapsto\, x \\[2pt] b \,\mapsto\, y
				}" swap]
			&[25pt] \rationals\ket{a,b} \mathrlap{\ = V_1}
				\arrow[d, "\phi_1" swap, "\color{black}\substack{
					a \,\mapsto\, x \\[2pt] b \,\mapsto\, y
				}"]
			\\
			\mathllap{W_0 =\ } \rationals\ket{x,y}
				\arrow[r, "\gamma_{1,0}", "\color{black}\substack{
					x \,\mapsto\, x ,\, y \,\mapsto\, y
				}" swap]
			& \rationals\ket{x,y} \mathrlap{\ = W_1}
		\end{tikzcd}\vspace{-5pt}
		\end{displaymath}
		Therefore, $\set{\phi_t}_{t \in \nonnegints}$ is \textbf{not} a persistence isomorphism.

		\item 
		It turns out that there exists a persistence isomorphism between $(V_\bullet, \alpha_\bullet)$ and $(W_\bullet, \gamma_\bullet)$.
		This becomes clearer if we perform a change of basis on $V_t$ and $W_t$ for $t \in \nonnegints$.
		Observe that $\rationals\ket{a,b} = \rationals\ket{c,b}$ with $c := a+b$, i.e.\ $\set{c,b}$ is a basis of $\rationals\ket{a,b}$.
		Then, we can express the persistence module $(V_\bullet,\alpha_\bullet)$ as follows. Note that $\alpha_t = \id_{\rationals\ket{a,b}}$ for $t \geq 1$.
		\vspace{-5pt}
		\begin{displaymath}
		\begin{tikzcd}[column sep=huge]
			{} &[-45pt] \color{black}\scriptstyle (t=0) 
			&[10pt] \color{black}\scriptstyle (t=1) 
			& \color{black}\scriptstyle (t=2) 
			& \color{black}\scriptstyle (t=3) 
			\\[-20pt]
			(V_\bullet,\alpha_\bullet) :
			& \rationals\ket{a,b} 
				\arrow[r, "\alpha_0", "\substack{
					a \,\mapsto\, a+b \,=:\, c \\[2pt] 
					b \,\mapsto\, b
				}" swap]
			& \rationals\ket{c,b}
				\arrow[r, "\alpha_1", "\substack{
					c \,\mapsto\, c \\[2pt]
					b \,\mapsto\, b 
				}" swap]
			& \rationals\ket{c,b}
				\arrow[r, "\alpha_2", "\substack{
					c \,\mapsto\, c \\[2pt]
					b \,\mapsto\, b 
				}" swap]
			& \rationals\ket{c,b}
				\arrow[r] &[-30pt] \cdots
		\end{tikzcd}
		\end{displaymath}
		Similarly, observe that 
		$\rationals\ket{x,y} = \rationals\ket{z,y}$ with $z := x+y$, i.e.\ $\set{z,y}$ is a basis for $\rationals\ket{x,y}$.
		We apply this change of basis on $W_t$ for all $t \in \nonnegints$.
		Since $\gamma_t = \id_{\rationals\ket{x,y}}$ for all $t \in \nonnegints$, the persistence module $(W_\bullet, \gamma_\bullet)$ can be illustrated as follows:
		\vspace{-5pt}
		\begin{displaymath}
		\begin{tikzcd}[column sep=huge]
			{} &[-45pt] \color{black}\scriptstyle (t=0) 
			&[10pt] \color{black}\scriptstyle (t=1) 
			& \color{black}\scriptstyle (t=2) 
			& \color{black}\scriptstyle (t=3) 
			\\[-20pt]
			(W_\bullet,\gamma_\bullet) :
			& \rationals\ket{z,y} 
				\arrow[r, "\gamma_0", "\substack{
					z \,\mapsto\, z \\[2pt] 
					y \,\mapsto\, y
				}" swap]
			& \rationals\ket{z,y}
				\arrow[r, "\gamma_1", "\substack{
					z \,\mapsto\, z \\[2pt] 
					y \,\mapsto\, y
				}" swap]
			& \rationals\ket{z,y}
				\arrow[r, "\gamma_2", "\substack{
					z \,\mapsto\, z \\[2pt] 
					y \,\mapsto\, y
				}" swap]
			& \rationals\ket{z,y}
				\arrow[r] &[-30pt] \cdots
		\end{tikzcd}
		\end{displaymath}
		For $t=0$, define $\psi_0: V_0 \to W_0$ by $a \mapsto z$ and $b \mapsto y$.
		For $t \geq 1$, define $\psi_t: V_t \to W_t$ by $c \mapsto z$ and $b \mapsto y$.
		We claim that $\psi_\bullet = \set{\psi_t}_{t \in \nonnegints}$ determines a persistence morphism $\psi_\bullet: (V_\bullet, \alpha_\bullet) \to (W_\bullet, \gamma_\bullet)$ since these two diagrams commute (with $t \geq 1$ for the diagram on the right hand side):
		\vspace{-5pt}
		\begin{displaymath}
		\begin{tikzcd}[column sep=huge, row sep=large]
			\mathllap{V_0 =\ } \rationals\ket{a,b} 
				\arrow[r, "\alpha_{1,0}", "\color{black}\substack{
					a \,\mapsto\, c ,\, b \,\mapsto\, b
				}" swap]
				\arrow[d, "\psi_0", "\color{black}\substack{
					a \,\mapsto\, z \\[2pt] b \,\mapsto\, y
				}" swap]
			&[15pt] \rationals\ket{c,b} \mathrlap{\ = V_1}
				\arrow[d, "\psi_1" swap, "\color{black}\substack{
					c \,\mapsto\, z \\[2pt] b \,\mapsto\, y
				}"]
			\\
			\mathllap{W_0 =\ } \rationals\ket{z,y}
				\arrow[r, "\gamma_{1,0}", "\color{black}\substack{
					z \,\mapsto\, z ,\, y \,\mapsto\, y
				}" swap]
			& \rationals\ket{z,y} \mathrlap{\ = W_1}
		\end{tikzcd}
		\qquad\text{ and }\quad 
		\begin{tikzcd}[column sep=huge, row sep=large]
			\mathllap{V_t =\ } \rationals\ket{c,b} 
				\arrow[r, "\alpha_{t}", "\color{black}\substack{
					\text{(identity)}
				}" swap]
				\arrow[d, "\psi_t", "\color{black}\substack{
					\text{(identity)}
				}" swap]
			& \rationals\ket{c,b} \mathrlap{\ = V_{t+1}}
				\arrow[d, "\psi_{t+1}" swap, "\color{black}\substack{
					\text{(identity)}
				}"]
			\\
			\mathllap{W_t =\ } \rationals\ket{z,y}
				\arrow[r, "\gamma_{t}", "\color{black}\substack{
					\text{(identity)}
				}" swap]
			& \rationals\ket{z,y} \mathrlap{\ = W_{t+1}}
		\end{tikzcd}\vspace{-5pt}
		\end{displaymath}

		Observe that, unlike the case for $\set{\phi_t}_{t \in \nonnegints}$, 
		$\psi_1 \circ \alpha_0$ and $\gamma_0 \circ \psi_0$ agree at $a \in V_0$ for $\set{\psi_t}_{t \in \nonnegints}$:
		\begin{equation*}
			\bigl( \psi_1 \circ \alpha_0 \bigr)(a)
			= \psi_1(c) = z 
			= \gamma_0(z)
			= \bigl( \gamma_0 \circ \psi_0 \bigr)(a)
		\end{equation*}
		Therefore, $\psi_\bullet = \set{\psi_t}_{t \in \nonnegints}$ is a persistence morphism 
		$\psi_\bullet: V_\bullet \to W_\bullet$.
		Since each linear map $\psi_t: V_t \to W_t$ is a linear isomorphism, $\psi_\bullet$ is a persistence isomorphism and $(V_\bullet, \alpha_\bullet) \cong (W_\bullet, \gamma_\bullet)$.
	\end{enumerate}
\end{example}


\spacer
Next, we provide a characterization of the direct sum of persistence modules.

\begin{definition}\label{defn:persmod-directsum}
	The \textbf{direct sum} $(\persmod{V}, \alpha_\bullet) \oplus (\persmod{W}, \beta_\bullet) =: (\persmod{U}, \gamma_\bullet)$ 
	of two persistence modules $(\persmod{V}, \alpha_\bullet)$ and $(\persmod{W}, \beta_\bullet)$ over a field $\field$
	is the persistence module over $\field$ with vector spaces given by $U_t := V_t \oplus W_t$ (i.e.\ a direct sum of vector spaces) for all $t \in \nonnegints$ and structure maps $\gamma_{s,t} = \alpha_{s,t} \oplus \beta_{s,t}$ (i.e.\ the unique linear map induced by the direct sum) for all $t,s \in \nonnegints$ with $t \leq s$, i.e.\ 
	\begin{equation*}\setlength{\arraycolsep}{0.5\arraycolsep}
	\begin{array}{cccc}
		\gamma_{s,t} :& 
		U_t := V_t \oplus W_t 
		&\to
		&V_s \oplus W_s =: U_s \\[1pt]
		& (v_t, w_t) &\mapsto 
		&\bigl(\alpha_{s,t}(v_t), \beta_{s,t}(w_t)\bigr)
	\end{array}
	\end{equation*}
\end{definition}


Observe that this definition extends to \textit{finite} direct sums of persistence modules.
Since direct sums of persistence modules are given by those of vector spaces, properties of the direct sum of vector spaces extend to the case of persistence modules.
We identify some of these below:
\begin{enumerate}
	\item Direct sums of persistence modules are \textit{commutative}, i.e.\ 
	$
		(\persmod{V}, \alpha_\bullet) \oplus (\persmod{W}, \beta_\bullet) 
		\cong 
		(\persmod{W}, \beta_\bullet) \oplus 
		(\persmod{V}, \alpha_\bullet) 
	$.

	\item 
	Finite direct sums of persistence modules are \textit{associative}, i.e.\ 
	\begin{equation*}
		\Bigl( 
			(\persmod{V}, \alpha_\bullet) \oplus (\persmod{W}, \beta_\bullet)
		\Bigr)
		\oplus (\persmod{U}, \gamma_\bullet)
		\cong 
			(\persmod{V}, \alpha_\bullet) \oplus 
		\Bigl( 
				(\persmod{W}, \beta_\bullet)
				\oplus (\persmod{U}, \gamma_\bullet)
		\Bigr)
	\end{equation*}
	Therefore, a statement such as $\persmod{V} \oplus \persmod{W} \oplus \persmod{U}$ is unambiguous (up to persistence isomorphism).

	\item 
	The zero persistence module $\persmod{0}$ with trivial vector spaces and trivial structure maps is the identity of the direct sum operation in $\catpersmod$,
	i.e.\ for any persistence module $(\persmod{V}, \alpha_\bullet)$, $(\persmod{V},\alpha_\bullet) \oplus \persmod{0} \cong (\persmod{V},\alpha_\bullet)$.

	\item 
	The distinction between \textit{internal} direct sums and \textit{external} direct sums of vector spaces also extend to the case of persistence modules, 
	That is, there is no difference between the two constructions up to vector space isomorphism when there are a finite number of non-trivial summands:
	\begin{enumerate}
		\item 
		An internal direct sum $A = A_1 \oplus A_2$ on a vector space $A$ is defined if $A_1$ and $A_2$ are both subspaces of $A$ and that $A_1 \cap A_2 = \set{0}$.
		Here, the elements of $A$ are generally not of the form $(a_1, a_2)$ with $a_1 \in A_1$ and $a_2 \in A_2$.
		For example, $A = \reals^2$, $A_1 = \vspan\set{(1,1)}$, and $A_2 = \vspan\set{(1,-1)}$.

		\item 
		An external direct sum $A := A_1 \oplus A_2$ of vector spaces $A_1$ and $A_2$ refers to direct/Cartesian product $A_1 \times A_2$ with a vector space structure induced by those of $A_1$ and of $A_2$.
		Here, the spaces $A_1$ and $A_2$ are interpreted to be distinct vector spaces and 
			the element of $A$ are exactly of the form $(a_1, a_2)$ with $a_1 \in A_1$ and $a_2 \in A_2$.
		Note that the direct sum and the direct product of a finite collection of vector spaces produce isomorphic vector spaces.
	\end{enumerate}
	For the case of persistence modules, 
		the collection of vector space isomorphisms between these two constructions for each index $t \in \nonnegints$ 
		forms a persistence isomorphism.

	In this paper, 
		we prefer to interpret direct sums as internal direct sums if possible and avoid denoting the elements of the direct sum as tuples.
\end{enumerate}
We provide an example of a direct sum of two persistence modules below. 

\begin{example}
	Let the persistence modules 
	$(\persmod{A},\alpha_\bullet)$ 
	and $(B_\bullet, \beta_\bullet)$ over $\ints_2$ be given as follows, with indeterminates $a_i$ for each $i \in \set{1,2,3}$ and $b_j$ for each $j \in \set{1,2,3}$:
	\begin{align*}
		A_t &= \begin{cases}
			\ints_2\ket{a_1, a_2} 		&\text{ if } t=0 \\
			\ints_2\ket{a_1, a_2, a_3}	&\text{ if } t=1 \\
			\ints_2\ket{a_2, a_3}		&\text{ if } t \geq 2
		\end{cases}
		&\hspace{-2pt}\text{ with }\hspace{7pt}
		&\left\{\begin{array}{r !{:} l !{\text{by}} lll}
			\alpha_0 & A_0 \to A_1 
				& a_1 \mapsto a_1, 
				& a_2 \mapsto a_2 \\
			\alpha_1 & A_1 \to A_2 
				& a_1 \mapsto 0, 
				& a_2 \mapsto a_2, 
				& a_3 \mapsto a_3
			\\
			\alpha_t & A_t \to A_{t+1} 
				& 
				& a_2 \mapsto a_2, 
				& a_3 \mapsto a_3
				\hspace{5pt}\text{ for } t \geq 2
		\end{array}\right\}
		\\[2pt]
		B_t &= \begin{cases}
			\ints_2\ket{b_1} 		&\text{ if } t=0 \\
			\ints_2\ket{b_1,b_2}	&\text{ if } t=1 \\
			\ints_2\ket{b_1,b_2,b_3}		&\text{ if } t \geq 2
		\end{cases}
		&\hspace{-2pt}\text{ with }\hspace{7pt}
		&\left\{\begin{array}{r !{:} l !{\text{by}} lll}
			\beta_0 & B_0 \to B_1 
				& b_1 \mapsto b_1 \\
			\beta_1 & B_1 \to B_2 
				& b_1 \mapsto b_1, 
				& b_2 \mapsto b_2
			\\
			\beta_t & B_t \to B_{t+1} 
				& b_1 \mapsto b_1, 
				& b_2 \mapsto b_2
				& b_3 \mapsto b_3
				\hspace{5pt}\text{ for } t \geq 2
		\end{array}\right\}
	\end{align*}
	These are illustrated below:
	\begin{center}
		\includegraphics[width=0.75\linewidth]{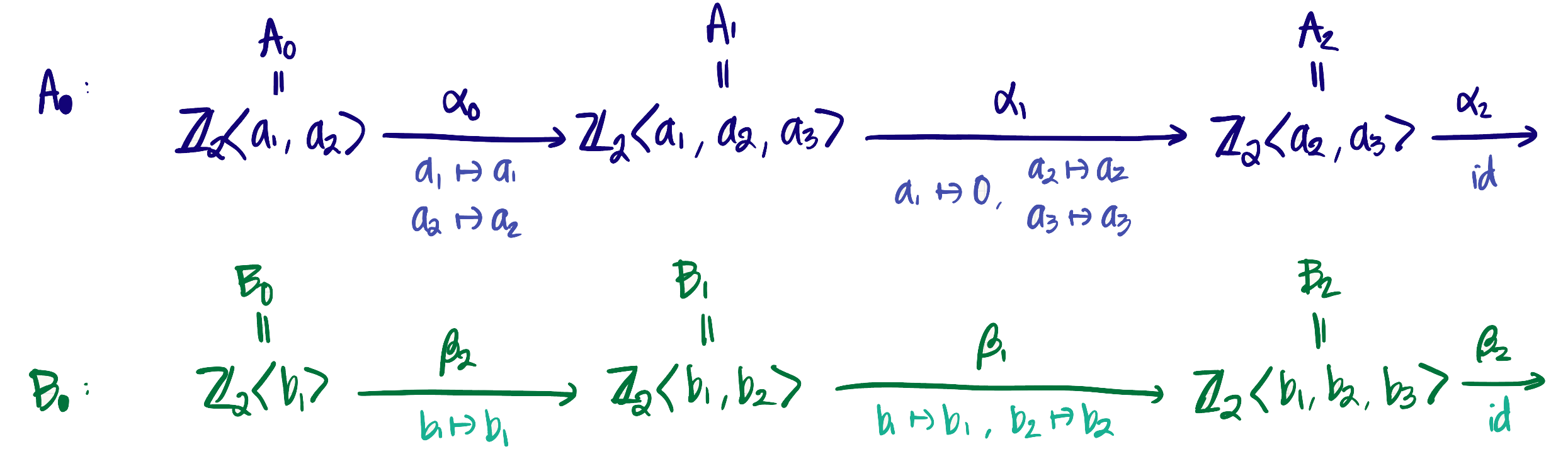}
	\end{center}
	Let $(Y_\bullet, \gamma_\bullet) :=
	(\persmod{A},\alpha_\bullet) \oplus (B_\bullet, \beta_\bullet)
	$.
	The vector spaces $Y_t$ of $(Y_\bullet, \gamma_\bullet)$ are given as follows:
	\begin{equation*}
		Y_t = \left\{\rule[0pt]{0pt}{26pt}\!\begin{array}{lcl c ll}
			\ints_2\ket{a_1,a_2} 
				&\oplus& \ints_2\ket{b_1}
				&\cong& \ints_2\ket{a_1,a_2,b_1}
				&\text{ if } t=0 \\
			\ints_2\ket{a_1, a_2, a_3} 
				&\oplus& \ints_2\ket{b_1,b_2}
				&\cong& \ints_2\ket{a_1,a_2,a_3,b_1,b_2}
				&\text{ if } t=1 \\
			\ints_2\ket{a_2,a_3} 
				&\oplus& \ints_2\ket{b_1,b_2,b_3} 
				&\cong& \ints_2\ket{a_2,a_3,b_1,b_2,b_3}
				&\text{ if } t \geq 2
		\end{array}\right.
	\end{equation*}
	Observe that the isomorphism relations above are valid since $a_i$'s and $b_j$'s are defined to be indeterminates 
	and that $\ints_2\ket{a_1, a_2, a_3} \cap \ints_2\ket{b_1,b_2,b_3} = \set{0}$.

	The structure maps $\gamma_t: Y_t \to Y_{t+1}$ can be described by collecting the assignments of $\alpha_t$ on the basis elements $a_1,a_2,a_3$ and those of $\beta_t$ on $b_1,b_2,b_3$ (whichever is included in the direct sum).
	For example, $\gamma_1: Y_1 \to Y_1$ has the following assignments:
	\begin{equation*}
		\gamma_1: Y_1 \to Y_2 
		\quad\text{ is given by }\quad
		\left\{\rule[0pt]{0pt}{16pt}\begin{array}{lll}
			a_1 \mapsto 0, 	
				& a_2 \mapsto a_2, 	
				& a_3 \mapsto a_3, \\
			b_1 \mapsto b_1,
				& b_2 \mapsto b_2,
		\end{array}\right\}
	\end{equation*}
	Since $\alpha_t = \id_{A_2}$ and $\beta_t = \id_{B_2}$ for all $t \geq 2$, 
		$(A_\bullet, \alpha_\bullet)$ and $(B_\bullet, \beta_\bullet)$ are both constant on $[2,\infty)$
		and the same applies to the direct sum $(Y_\bullet, \gamma_\bullet)$.
\end{example}


For an arbitrary category $\catname{C}$,
	the terms \textit{subobject}, \textit{kernel}, \textit{cokernel}, and \textit{image} are generally defined to be morphisms of $\catname{C}$ satisfying certain properties (as opposed to objects).
For an arbitrary category $\catname{C}$,
	a subobject of an object $x$ in $\catname{C}$ refers to an injective morphism $y \rightarrowtail x$ with codomain $x$ \cite[Definition 4.6.8]{cattheory:rhiel}.
For some of the categories that are more accessible to the introductory learner,
	subobjects $y \rightarrowtail x$ of an object $x$ (assuming the morphism agrees with the identity $\id_x: x \to x$, i.e.\ $y \hookrightarrow x$) are characterized by the domain $y$ 
	and are labeled using a term specific to the category.
We list some examples below.
\begin{enumerate}
	\item 
	In the category $\catvectspace$ of vector spaces over a field $\field$,
		subobjects correspond to \textit{(vector) subspaces}
	and each vector subspace $W$ of a vector space $V$ has a corresponding inclusion map $W \hookrightarrow V$.
	Since $W$ has to be a vector space (i.e.\ an object of $\catvectspace$),
	the inclusion $W \hookrightarrow V$ must be a linear map, i.e.\ a morphism in $\catvectspace$.

	\item 
	In the category $\catsimp$ of (abstract) simplicial complexes and simplicial maps (as given in \fref{defn:cat-simp}),
		subobjects coincide with (simplicial) subcomplexes.
	Given a subcomplex $L$ of a simplicial complex $K$,
		the set-wise inclusion map $L \hookrightarrow K$ is also a simplicial map (as described in \fref{defn:simplicial-maps}).

	\item 
	In the category $\cattop$ of topological spaces,
		subobjects are (topological) subspaces.
	Given any subset $Y$ of a topological space $X$,
		$Y$ can be equipped with the subspace topology of $X$.
	The set-wise inclusion map $Y \hookrightarrow X$
		is then made a continuous map relative to this topology on $Y$.
\end{enumerate}
Below, we provide a corresponding definition for subobjects on the category $\catpersmod$ of persistence modules.

\begin{definition}\label{defn:persistence-submodule}
	Fix a field $\field$.
	Let $(V_\bullet, \alpha_\bullet)$ be a persistence module over a field $\field$.
	A persistence module $(W_\bullet, \gamma_\bullet)$ over $\field$ is called a \textbf{persistence submodule} or \textbf{submodule} of $(V_\bullet, \alpha_\bullet)$ 
		if the following are true:
	\begin{enumerate}
		\item 
		$W_t$ is a vector subspace of $V_t$ for all $t \in \nonnegints$.

		\item 
		The structure maps $\gamma_{s,t}: W_t \to W_s$ of 
			$(W_\bullet, \gamma_\bullet)$
			are exactly the structure maps of 
			$\alpha_{s,t}: V_t \to V_s$ of 
			$(W_\bullet, \gamma_\bullet)$ with domain and codomain restricted to $W_t$ and $W_s$ respectively.
	\end{enumerate}
	That is, 
	the collection $\set{i_t}_{t \in \nonnegints}$ of inclusion maps $i_t: W_t \hookrightarrow V_t$ forms a persistence morphism $(W_\bullet,\gamma_\bullet) \hookrightarrow (V_\bullet, \alpha_\bullet)$.
	In this case, we write $(W_\bullet, \gamma_\bullet) \subseteq (V_\bullet, \alpha_\bullet)$ (or $W_\bullet \subseteq V_\bullet$ for convenience).
\end{definition}

This definition of subobjects for persistence modules allows characterizations of kernels, cokernels, and images of persistence morphisms like those of linear maps in $\catvectspace$, i.e.\ as objects in the category as opposed to morphisms.
In particular,
	we extend the set-wise definition of kernels, images, and cokernels in $\catvectspace$ to the case of $\catpersmod$.
The composition relation on the structure maps of persistence modules (i.e.\ the commuting squares condition) guarantee that the resulting collection of vector subspaces form a persistence submodule.
We state this in more detail below.

\begin{statement}{Definition}\label{defn:pers-morphisms-features}
	Let $\phi_\bullet: (\persmod{V},\alpha_\bullet) \to (\persmod{W},\gamma_\bullet)$ be a persistence morphism between two persistence modules $(\persmod{V},\alpha_\bullet)$ and $(\persmod{W},\gamma_\bullet)$ over a field $\field$ with $\phi_\bullet = \set{\phi_t: V_t \to W_t}_{t \in \nonnegints}$.
	\begin{enumerate}
		\item 
		The \textbf{kernel} $\ker(\phi_\bullet) =: (K_\bullet, \kappa_\bullet)$ of $\phi_\bullet$ 
		is the submodule of $(V_\bullet, \alpha_\bullet)$ with 
			$K_t = \ker(\phi_t) \subseteq V_t$ for all $t \in \nonnegints$.

		\item 
		The \textbf{image} $\image(\phi_\bullet) =: (B_\bullet, \beta_\bullet)$ of $\phi_\bullet$ is the submodule of $(W_\bullet, \gamma_\bullet)$
		with 
			$B_t = \im(\phi_t) \subseteq W_t$ for all $t \in \nonnegints$.

		\item 
		The \textbf{cokernel} $\coker(\phi_\bullet) =: (C_\bullet, \lambda_\bullet)$ of $\phi_\bullet$ is the persistence module 
		with vector spaces 
		given by $C_t := \coker(\phi_t) = W_t \bigmod \im(\phi_t)$ for all $t \in \nonnegints$
		and structure maps
			$\lambda_{s,t}: C_t \to C_s$ being the linear maps induced by $\gamma_{s,t}: W_t \to W_s$ and the cokernel/quotient construction on $\catvectspace$,
			i.e.\ $\eta_{s,t}$ maps $[w] \to [\gamma_{s,t}(w)]$ with $w \in W_t$,
			for all $t,s\in \nonnegints$ with $t \leq s$.
	\end{enumerate}
\end{statement}
\remark{
	The symbols $K$, $B$, and $C$ for the kernel, image, and cokernel respectively are used here for convenience, i.e.\ we do not typically use these letters for the kernel, image, and cokernel respectively.
}

Observe that the vector spaces for the kernel, image, and cokernel of persistence morphisms are all well-defined since the kernel, image, and cokernel of linear maps are well-defined.
For the kernel and image,
	the question is whether the resulting structure maps by \fref{defn:persistence-submodule} as persistence submodules are well-defined.
We can show this using a diagram chase  
using a diagram chase on the commuting squares of $(V_\bullet, \alpha_\bullet)$ and $(W_\bullet, \gamma_\bullet)$ on the indices $t$ and $s$, as illustrated below.
\begin{center}
	\includegraphics[width=0.65\linewidth]{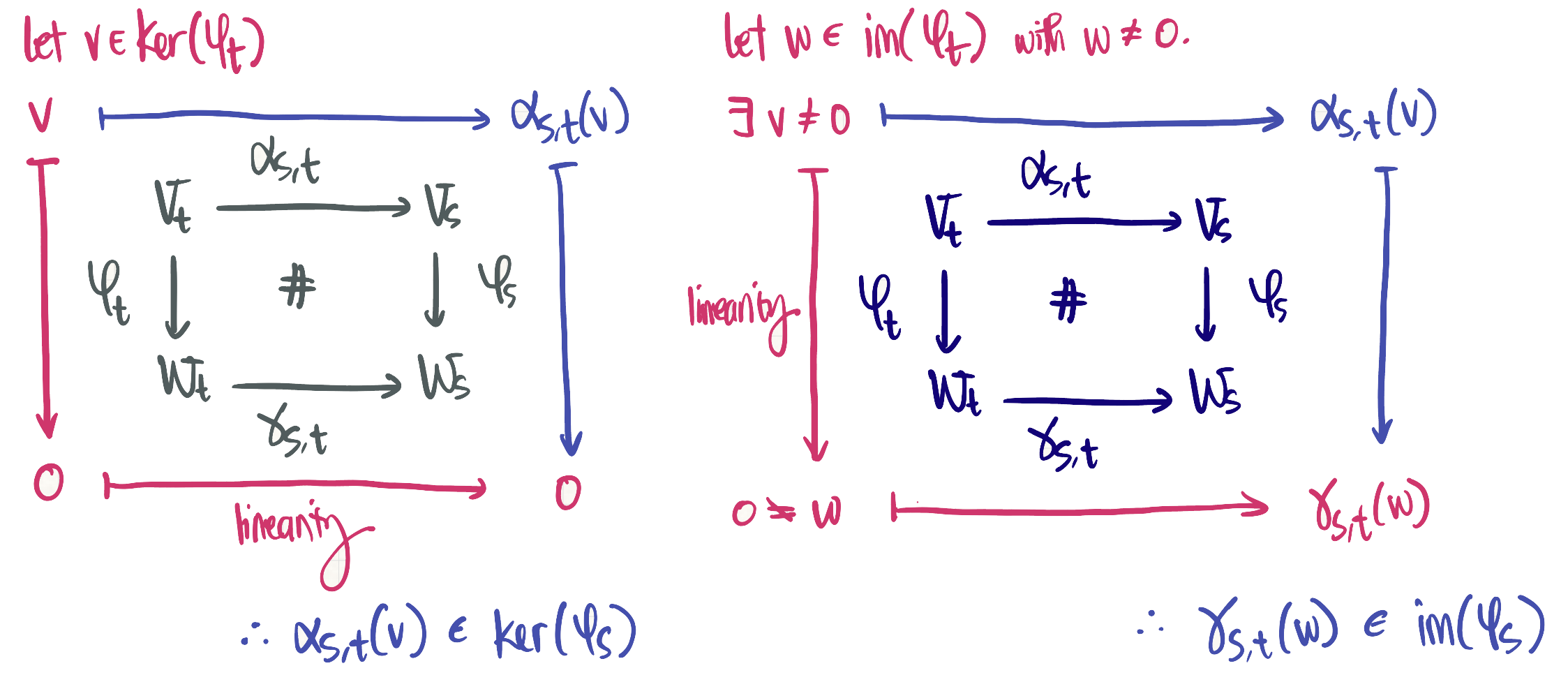}
\end{center}
The diagram chase on the left implies that 
	$\alpha_{s,t}\bigl( \ker(\phi_t) \bigr) \subseteq \ker(\phi_s) \subseteq V_s$ 
	and that the restriction of the domain $V_t$ to $\ker(\phi_t)$ and the codomain $V_s$ to $\ker(\phi_s)$ 
	of the structure map $\alpha_{s,t}: V_t \to V_s$ is well-defined. 

The diagram chase on the right implies that 
	$\gamma_{s,t}\bigl( \im(\phi_t) \bigr) \subseteq \im(\phi_s) \subseteq W_s$.
This tells us that the restriction of the domain and codomain of the structure map $\gamma_{s,t}: W_t \to W_s$ to $\im(\phi_t)$ and $\im(\phi_s)$ respectively is well-defined.
This also tells us that the structure maps of $\coker(\phi_\bullet)$ are well-defined since
\begin{equation*}
	\gamma_{s,t}\bigl( w + \im(\phi_t) \bigr) 
	= 
	\gamma_{s,t}(w) + \im(\phi_s) \in W_s
\end{equation*}
and that, for any coset representative $w \in \im(\phi_t)$ of the trivial element $[w] = 0 \in W_t \bigmod \im(\phi_t)$,
	$\lambda_{s,t}([w]) = [\gamma_{s,t}(w)] = 0 \in W_s \bigmod \im(\phi_s)$.

We provide an example below where the kernel, image, and cokernel of a persistence morphism are identified.

\begin{example}
	Let $(\persmod{A}, \alpha_\bullet)$ and $(\persmod{X}, \chi_\bullet)$ be persistence modules over $\rationals$ with vector spaces $A_n$ and $X_n$ defined as follows:
	\begin{equation*}
		A_t = \begin{cases}
			\rationalsket{a} 	&\text{ if } t = 0 \\
			\rationalsket{a,b}	&\text{ if } t = 1 \\
			\rationalsket{a,b,c,d}		&\text{ if } t \geq 2
		\end{cases}
		\qquad\text{ and }\qquad 
		X_t = \begin{cases}
			\rationalsket{x,y} 	&\text{ if } t = 0,1 \\
			\rationalsket{x,y,z}	&\text{ if } t = 2 \\
			\rationalsket{x,y,z,w}		&\text{ if } t \geq 3
		\end{cases}
	\end{equation*}
	where the structure maps $\alpha_{s,t}: A_t \to A_s$
	of $\persmod{A}$ and 
	$\chi_{s,t}: X_t \to X_s$
	of $\persmod{X}$ are restrictions of the identity maps on $\rationalsket{a,b,c,d}$ and $\rationalsket{x,y,z,w}$ respectively.
	Define the map $\Phi: \rationals\ket{a,b,c,d} \to \rationals\ket{x,y,z,w}$ by 
	\begin{equation*}
		[\Phi] = \hspace{3pt}
		\begin{NiceArray}{>{\color{black}}{c} !{\,\,} cccc}
			\RowStyle{\color{black}}
			{} & a & b & c & d \\
			{x} & 1 & -1 & 1 & 0 \\
			{y} & 0 & 0 & 1 & 1 \\
			{z} & 0 & 0 & 1 & 1 \\
			{w} & 0 & 0 & 0 & 0
		\CodeAfter 
			\SubMatrix({2-2}{5-5})[left-xshift=3pt, right-xshift=3pt]
		\end{NiceArray}
		\hspace{9pt}\qquad\text{ or equivalently }\qquad
		\left\{\begin{array}{cc r@{\,}r@{\,}r}
			a &\mapsto& x \\
			b &\mapsto& -x 	\\
			c &\mapsto& x&+\,y&+\,z \\
			d &\mapsto& &y&+\,z	\\
		\end{array}\right.
	\end{equation*}
	For each $t \in \nonnegints$, define $\phi_t: A_t \to X_t$ by $\phi_t(f) = \Phi(f)$ for all $f \in A_t$.
	Observe that the codomains of $\phi_t$ with $t=0,1,2,3$ are well-defined.
	We claim that 
		$\phi_\bullet = \set{\phi_t}_{t \in \nonnegints}$ is a persistence morphism $\phi_\bullet: (A_\bullet, \alpha_\bullet) \to (X_\bullet, \chi_\bullet)$.

	The vector spaces of the persistence module $\ker(\phi_\bullet)$ by $\ker(\phi_\bullet)(t) := \ker(\phi_t)$ and of the persistence module $\im(\phi_\bullet)$ by $\im(\phi_\bullet)(t) := \im(\phi_t)$
	are as follows:
	\begin{equation*}
		\ker(\phi_t) = \begin{cases}
			0 					&\text{ if } t=0 \\
			\rationals\ket{a+b}	&\text{ if } t=1 \\
			\rationals\ket{a+b, a-c+d}
					&\text{ if } t \geq 2 \\
		\end{cases}
		\qquad\text{ and }\qquad 
		\im(\phi_t) = \begin{cases}
			\rationals\ket{x} 	&\text{ if } t = 0,1 \\
			\rationals\ket{x,y+z} &\text{ if } t \geq 2
		\end{cases}
	\end{equation*}
	Since the vector subspace relations $A_t \subseteq A_s$ and $X_t \to X_s$ for $t \leq s$ are preserved by $\phi_t$ and $\phi_s$, 
		the structure maps $\ker(\phi_\bullet)(t \to s)$ and $\im(\phi_\bullet)(t \to s)$, which are inclusion maps as restrictions of identity maps, are both well-defined.
	Observe that $\ker(\phi_\bullet)$ is a persistence submodule of $(A_\bullet, \alpha_\bullet)$ and $\im(\phi_\bullet)$ is a persistence submodule of $(X_\bullet, \chi_\bullet)$.

	The vector spaces of $\coker(\phi_\bullet)$ by 
		$\coker(\phi_\bullet)(t) = \coker(\phi_t)$ can be described as follows:
	\begin{equation*}
		\coker(\phi_t) = \frac{X_t}{\im(\phi_t)} =
		\left\{\rule[0pt]{0pt}{26pt}\!\begin{array}{
			r@{\,}c@{\,\,}l cl cl l}
			\rationals\ket{x,y} &\bigmod& \rationals\ket{x}
				&\cong& \rationals\ket{[y]}
				&\cong& \rationals
				&\text{ if } t=0 \\
			\rationals\ket{x,y} &\bigmod& \rationals\ket{x}
				&\cong& \rationals\ket{[y]}
				&\cong& \rationals
				&\text{ if } t=1 \\
			\rationals\ket{x,y,z} &\bigmod& \rationals\ket{x,y+z}
				&\cong& \rationals\ket{[z]}
				&\cong& \rationals
				&\text{ if } t=2 \\
			\rationals\ket{x,y,z,w}&\bigmod& \rationals\ket{x,y+z}
				&\cong& \rationals\ket{[z],[w]}
				&\cong& \rationals^2
				&\text{ if } t \geq 3 \\
		\end{array}\right.
	\end{equation*} 
	Since vector spaces have no torsion, the ranks of the $\coker(\phi_t)$ for all $t \in \nonnegints$ is determined by the rank of $X_t$ minus the rank of $\im(\phi_t)$.
	The structure maps of $\coker(\phi_\bullet)$ are illustrated by the following sequence:
	\begin{displaymath}
	\begin{tikzcd}
			\coker(\phi_0)		\arrow[d,equal]
			& \coker(\phi_1)	\arrow[d,equal]
			&[45pt] \coker(\phi_2)	\arrow[d,equal]
			&[40pt] \coker(\phi_3)	\arrow[d,equal]
			&[45pt]
		\\[-10pt]
		\rationals\ket{[y]} \arrow[r, "\substack{
				[y] \,\mapsto\, [y]
			}"]
		& \rationals\ket{[y]} 
			\arrow[r, "\substack{
				[y] \,\mapsto\, 0
			}", "\substack{
				\text{ $[y]$ becomes trivial } \\ 
				\text{ on addition of $y+z$ } \\
				\text{ to $\im(\phi_1)$ }
			}" swap]
		& \rationals\ket{z} 
			\arrow[r, "\substack{
				[z] \,\mapsto\, [z]
			}", "\substack{
				\text{ $[z]$ is unaffected } \\ 
				\text{ by addition of $w$ } \\
				\text{ to $X_3$ }
			}" swap]
		& \rationals\ket{z,w} 
			\arrow[r, "\substack{
				\text{structure maps} \\ 
				\text{are identity maps} \\
				\text{for $t\geq 3$}
			}" swap]
		& \cdots
	\end{tikzcd}
	\end{displaymath}
\end{example}


Note that defining the cokernel of morphisms in an abelian category defines the quotient operation on said category.
For example, in the category of $\catmod{R}$ of modules over a PID $R$, the quotient module $X \bigmod Y$ of an $R$-module $X$ by its submodule $Y$ is exactly the cokernel of the inclusion map $Y \hookrightarrow X$, i.e.\ 
\begin{equation*}
	X \bigmod Y = X \bigmod \im(Y \hookrightarrow X) = \coker(Y \hookrightarrow X)
\end{equation*}
Consequently, \fref{defn:pers-morphisms-features} also tells us that the quotient operation on a persistence module by its submodule is done pointwise for each $t \in \nonnegints$.
Since we have definitions for persistence submodule relations and quotients of persistence modules,
	we can form chain complexes of persistence modules.

More generally, given any abelian category $\catname{A}$,
there is a corresponding category $\catchaincomplex{\catname{A}}$ of chain complexes in $\catname{A}$, defined very similarly as in the case of chain complexes of $R$-modules.
For those interested, we refer to \cite[Section 5.5]{cattheory:rotman} and \cite[Chapter 1]{hom-algebra:weibel} for a more general and detailed discussion.
We state what this means specifically for $\catpersmod$ below.

\begin{definition}\label{defn:category-of-persistence-chains}
	The \textbf{category of persistence complexes over a field $\field$}, denoted $\catchaincomplex{\catpersmod}$, refers to the category of chain complexes on $\catpersmod$, with objects and morphism described below:
	\begin{enumerate}
		\item 
		The objects of $\catchaincomplex{\catpersmod}$ are persistence complexes, defined as follows:

		A \textbf{persistence complex} $(V^\bullet_\ast, \alpha^\bullet_\ast, \boundary_\ast^\sbullet)$ over $\field$ is a $\ints$-indexed collection of persistence modules 
		$(V^\bullet_n, \alpha^\bullet_n)$ over $\field$
		and persistence morphisms 
			$\boundary_n\updot: (V^\bullet_n, \alpha^\bullet_n) \to (V^\bullet_{n-1}, \alpha^\bullet_{n-1})$
		such that 
			$\boundary_n\updot \circ \boundary_{n-1}\updot = 0^\bullet$ 
			for all $n \in \ints$,
			where $0^\bullet$ denotes the zero persistence morphism.
		Illustrated below is $(V^\bullet_\ast, \alpha^\bullet_\ast, \boundary^\bullet_\ast)$ as a filtered sequence of persistence modules and persistence morphisms:
		\begin{displaymath}
		\begin{tikzcd}[column sep=huge]
			\cdots \arrow[r, "\boundary_{n+2}^\sbullet"]
			& V_{n+1}^\bullet \arrow[r, "\boundary_{n+1}^\sbullet"]
			& V_{n}^\bullet \arrow[r, "\boundary_{n}^\sbullet"]
			& V_{n-1}^\bullet \arrow[r, "\boundary_{n-1}^\sbullet"]
			& \cdots
		\end{tikzcd}
		\end{displaymath}
		For convenience, we may suppress the structure maps of $V^\bullet$ and write $(V^\bullet_\ast, \boundary^\sbullet_\ast) := (V^\bullet_\ast, \alpha^\bullet_\ast, \boundary_\ast^\sbullet)$.

		\item 
		The morphisms of $\catchaincomplex{\catpersmod}$ are persistence chain morphisms, defined as follows:

		A \textbf{persistence chain morphism} 
		$
			f^{\,\bullet}_\ast: (V^\bullet_\ast, \boundary^{\,\bullet}_\ast)
				\to (W^\bullet_\ast, \delta^{\,\bullet}_\ast)
		$
		between two persistence complexes 
		$(V^\bullet_\ast, \boundary^{\,\bullet}_\ast)$
		and 
		$(W^\bullet_\ast, \delta^{\,\bullet}_\ast)$ over $\field$
		is a $\ints$-indexed collection of persistence morphisms
		$f^\sbullet_n: V^\bullet_n \to W^\bullet_n$ 
		such that for all $n \in \ints$,
			$f^\sbullet_{n-1} \circ \boundary_n^\sbullet 
			= \delta_n^\sbullet \circ f^\sbullet_{n}$,
		i.e.\ we have the following commutative squares:
		\vspace{-5pt}
		\begin{equation*}
			\begin{tikzcd}[column sep=large, row sep=normal]
				V_n^\bullet 
					\arrow[r, "\boundary_n^\sbullet"]
					\arrow[d, "f_n^\sbullet", swap]
					\arrow[rd, "\hash" description,draw=none]
					& V_{n-1}^\bullet 
					\arrow[d, "f_{n-1}^\sbullet"]
					\\
				W_n^\bullet 
					\arrow[r, "\delta_n^\sbullet", swap]
					& W_{n-1}^\bullet
			\end{tikzcd}
		\end{equation*}
	\end{enumerate}
\end{definition}

\vspace{-5pt}
Observe that the notation for the persistence complex $(V^\bullet_\ast, \boundary^\sbullet_\ast)$ has two ``placeholder'' indices: an index $t \in \nonnegints$ denoted by the bullet $(\bullet)$ and an index $n \in \ints$ denoted by the asterisk $(\ast)$.
Since each $n \in \ints$ identifies a persistence module $(V^\bullet_n, \alpha^\bullet_n)$ 
and each $t \in \nonnegints$ identifies an $\field$-vector space $V^{t}_n$,
	a persistence complex $(V^\bullet_\ast, \boundary^\sbullet_\ast)$ corresponds to the following commutative grid of vector spaces:

	\null\vspace{-1em}
	\begin{displaymath}
		\begin{tikzcd}[column sep=6em, row sep=normal, inner sep=10pt]
			{} 	&[-70pt]& V^\bullet_{n+1} 
				& V^\bullet_{n} 
				& V^\bullet_{n-1} 
			\\[-20pt]
				{} && \rotatebox{90}{:}
				& \rotatebox{90}{:}
				& \rotatebox{90}{:}
			\\[-25pt]
			{} 	&& \vdots 
				& \vdots 
				& \vdots
			\\
			V_\ast^{t+1} \quad&:\quad
			\cdots 
				\arrow[r, "\boundary^{\hspace{1pt}t+1}_{n+2}"]
				& V^{t+1}_{n+1}
					\arrow[r, "\boundary^{\hspace{1pt}t+1}_{n+1}"]
					\arrow[u, "\alpha^{t+1}_{n+1}", swap]
				& V^{t+1}_{n} 
					\arrow[r, "\boundary^{\hspace{1pt}t+1}_{n}"]
					\arrow[u, "\alpha^{t+1}_{n}", swap]
				& V^{t+1}_{n-1} 
					\arrow[r, "\boundary^{\hspace{1pt}t+1}_{n-1}"]
					\arrow[u, "\alpha^{t+1}_{n-1}", swap]
				& \cdots 
			\\
			V_\ast^{t} \quad&:\quad
			\cdots \arrow[r, "\boundary^{\hspace{1pt}t}_{n+2}"]
				& V^{t}_{n+1} 
					\arrow[r, "\boundary^{\hspace{1pt}t}_{n+1}"]
					\arrow[u, "\alpha^{t}_{n+1}", swap]
				& V^{t}_{n} 
					\arrow[r, "\boundary^{\hspace{1pt}t}_{n}"]
					\arrow[u, "\alpha^{t}_{n}", swap]
				& V^{t}_{n-1} 
					\arrow[r, "\boundary^{\hspace{1pt}t}_{n-1}"]
					\arrow[u, "\alpha^{t}_{n-1}", swap]
				& \cdots 
			\\
			V_\ast^{t-1} \quad&:\quad
			\cdots \arrow[r, "\boundary^{\hspace{1pt}t-1}_{n+2}"]
				& V^{t-1}_{n+1} 
					\arrow[r, "\boundary^{\hspace{1pt}t-1}_{n+1}"]
					\arrow[u, "\alpha^{t-1}_{n+1}", swap]
				& V^{t-1}_{n} 
					\arrow[r, "\boundary^{\hspace{1pt}t-1}_{n}"]
					\arrow[u, "\alpha^{t-1}_{n}", swap]
				& V^{t-1}_{n-1} 
					\arrow[r, "\boundary^{\hspace{1pt}t-1}_{n-1}"]
					\arrow[u, "\alpha^{t-1}_{n-1}", swap]
				& \cdots 
			\\
			{} 	&& \vdots \arrow[u, "\alpha^{t-2}_{n+1}", swap]
				& \vdots \arrow[u, "\alpha^{t-2}_{n}", swap]
				& \vdots \arrow[u, "\alpha^{t-2}_{n-1}", swap]
		\end{tikzcd}
	\end{displaymath}

Note that each column represents a single persistence module and each row represents a chain complex of vector spaces and linear maps.
Persistence complexes will be significant later in \fref{section:construction-of-persistent-homology} of \fref{chapter:filtrations-and-pershoms} where we extend the construction of simplicial homology of simplicial complexes to the case of persistence modules.

The category $\catchaincomplex{\catname{A}}$ of chain complexes on an abelian category $\catname{A}$ also brings with it a family of chain homology functors $H_n(-): \catchaincomplex{\catname{A}} \to \catname{A}$, one for each $n \in \ints$. 
We state a corresponding definition specific to the category $\catpersmod$ below.

\begin{definition}
	For each $n \in \ints$, 
	the $n$\th \textbf{chain homology functor} 
	$H_n: \catchaincomplex{\catpersmod} \to \catpersmod$ 
	on $\catchaincomplex{\catpersmod}$
	sends a persistence complex $(V_\ast^\bullet, \boundary_\ast^\sbullet)$ to its $n$\th \textbf{chain homology} $H_n(V_\ast^\bullet, \boundary_\ast^\sbullet)$, which is the persistence module given by 
	\begin{equation*}
		H_n(V_\ast^\bullet, \boundary_\ast^\sbullet) 
		= \frac{ \ker(\boundary_n^\sbullet) }{ \im(\boundary_{n+1}^\sbullet) }
		= \coker\Bigl(
			\im(\boundary_{n+1}^\sbullet) 
			\hookrightarrow
			\ker(\boundary_n^\sbullet)
		\Bigr)
	\end{equation*}
\end{definition}
 \clearpage

\section{Interval Decompositions of Persistence Modules} 
\label{section:interval-decomposition-persmod}

Now that we have definitions for isomorphisms and direct sums between persistence modules, 
it is possible for us to talk about decompositions of persistence modules.
Note that we use the term \textit{decomposition}, or more specifically \textit{direct sum decomposition}, in the same way we would for vector spaces and modules. 
That is, a decomposition of a persistence module consists of a (finite) direct sum of other persistence modules such that the direct sum is isomorphic to the original persistence module.

In persistence theory, we are interested in a unique form of decomposition called \textit{interval decomposition}. 
In this section, we define the notions of interval modules and interval decompositions and provide examples for each.
To start, we provide a definition for interval modules, adapted from~\cite[Section 1.4]{persmod:chazal-structure}.

\begin{definition}\label{defn:persistence-interval-modules}
	Let $\field$ be a field.
	The $J$-\textbf{interval module $\intmod{J}$ over $\field$} for some interval $J \subseteq \nonnegints$
	is the persistence module $\intmod{J} = (\intmod{J}, i_\bullet^J): \posetN \to \catvectspace$ with vector spaces $\intmod[t]{J}$ and structure maps $i_{s,t}^J: \intmod[t]{J} \to \intmod[s]{J}$ given as follows:
	\begin{equation*}
		\intmod[t]{J} = \begin{cases}
			\field 			&\text{ if } t \in J \\
			0 				&\text{ if } t \not\in J
		\end{cases}
		\quad\text{ for all } t \in \nonnegints 
		\quad\text{ and }\quad 
		i_{s,t}^J
		 = \begin{cases}
			\id_{\field} 			&\text{ if } t,s \in J \\
			0 				&\text{ otherwise } 
		\end{cases}
		\quad\text{ for all } t,s \in \nonnegints \text{ with } t \leq s.
	\end{equation*}
	If the interval $J$ and the field $\field$ are arbitrary, we may drop references to both and refer to $\intmod{J}$ as an \textbf{interval module}.
	If $J = [a,b)$, we may write $\intmod{[a,b)}$ to refer to $\intmod{J}$. Similarly, we may write $\intmod{[a, \infty)}$ if $J = [a, \infty)$.
\end{definition}
\remark{
	The condition that $i_{s,t}^J = \id_{\field}$ if $t,s \in J$ is stated as $i_{s,t} = 1$ in \cite[Section 1.4]{persmod:chazal-structure}.
	In this case, $i_{s,t}^J = 1$ refers to the linear map $\field \to \field$ by $k \mapsto 1 \mathrel{\cdot} k = k$ for all $k \in \field$.
}

Since we are only using $\nonnegints$ as the indexing set for persistence modules, we can characterize all intervals in $\nonnegints$ as pairs of values.
In particular, any interval in $\nonnegints$ can be represented using exactly one of the two representations below:
\begin{equation*}
	[a,b) = \set{n \in \nonnegints : a \leq n < b}
	\qquad\text{ or }\qquad 
	[a, \infty) = \set{n \in \nonnegints : a \leq n}
\end{equation*}
Observe that, assuming we allow $b$ in $[a,b)$ to have $\infty$ as a value (i.e.\ $b \in \nonnegints \cup \set{\infty}$), every interval in $\nonnegints$ can be unambiguously represented by the two endpoints $a$ and $b$. 
Thus, in some papers (e.g.~\cite{matrixalg:zomorodian,ripser}), the interval $[a,b)$ is represented using an ordered pair $(a,b)$. 

Since we will be working at the level of persistence isomorphisms,
	it will be helpful to have a characterization of persistence modules that are isomorphic to interval modules.
We provide such a characterization below in terms of ranks.

\begin{lemma}\label{lemma:int-mod-isom}
	A persistence module $(V_\bullet, \alpha_\bullet)$ over $\field$ is isomorphic to the interval module $\intmod{J}$ for some interval $J \subseteq \nonnegints$ 
	if and only if the vector spaces $V_t$ and structure maps $\alpha_{s,t}: V_t \to V_s$ satisfy the following:
	\begin{equation}
		\rank(V_t) = \begin{cases}
			1 			&\text{ if } t \in J \\
			0 			&\text{ if } t \not\in J
		\end{cases}
		\qquad\text{ and }\qquad 
		\rank(\alpha_{s,t}) = 
		\begin{cases}
			1		&\text{ if } t,s \in J \\
			0		&\text{ otherwise } 
		\end{cases}
		\label{eqn:int-mod-isom}
	\end{equation}
	where $\rank(-)$ refers to the rank of $\field$-vector spaces and of $\field$-linear maps.
\end{lemma}
\begin{proof}
	Let $(V_\bullet, \alpha_\bullet)$ be a persistence module over $\field$.

	For the forward direction, 
		assume that $V_\bullet \cong \intmod{J}$ for some interval $J \subseteq \nonnegints$
		with $\intmod{J} = (\intmod{J}, i_\bullet^J)$.
	Then, there exists a persistence isomorphism 
		$\phi_\bullet: V_\bullet \to \intmod{J}$
		with $\phi_\bullet = 
		\set{\phi_t: V_t \to W_t}_{t \in \nonnegints}$.
	By \fref{defn:persistence-isomorphism}, $\phi_t: V_t \to W_t$ is an $\field$-vector space isomorphism for each $t \in \nonnegints$.
	Since linear isomorphisms preserve the ranks of vector spaces, we have the following for each $t \in \nonnegints$:
	\begin{equation*}
		\rank(V_t)
		= \rank\bigl( \phi_t(V_t) \bigr)
		= \rank\bigl( \intmod[t]{J} \bigr)
		= \left\{\begin{array}{c@{\ }c@{\ }c l}
			\rank(\field) &=& 1 &\text{ if } t \in J \\[2pt]
			\rank(0) &=& 0 &\text{ otherwise }
		\end{array}\right.
	\end{equation*}
	Similarly, composition with linear isomorphisms preserve the ranks of linear maps. Then, for all $t,s \in \nonnegints$ with $t \leq s$:
	\begin{equation*}
		\rank(\alpha_{s,t})
		= \rank(\phi_s \circ \alpha_{s,t} \circ \phi_t)
		= \rank(i_{s,t})
		= \left\{\begin{array}{c@{\ }c@{\ }c l}
			\rank(\id_{\field}) &=& 1 &\text{ if } t,s \in J \\[2pt]
			\rank(0) &=& 0 &\text{ otherwise }
		\end{array}\right.
	\end{equation*}
	Therefore, $V_\bullet$ satisfies \ref{eqn:int-mod-isom} for the interval $J \subseteq \nonnegints$.

	For the backwards direction, assume that $(V_\bullet, \alpha_\bullet)$ satisfies \ref{eqn:int-mod-isom} for some interval $J \subseteq \nonnegints$. 
	To construct the persistence isomorphism $\phi_\bullet: V_\bullet \to \intmod{J}$ with 
	$\intmod{J} = (\intmod{J}, i_\bullet^J)$,
		we build the linear maps $\phi_t: V_t \to \intmod[t]{J}$ inductively.
	For clarity, let $1_t \in \field$ refer to the multiplicative identity of $\field$, specifically as the vector space of $\intmod{J}$ at index $t \in J$, 
		i.e.\ $i_{s,t}^J(1_t) = \id_{\field}(1_t) = 1_s$ for all $t,s \in J$ with $t \leq s$.
	\begin{enumerate}
		\item 
		Let $a = \min(J)$, which exists by the well-ordering principle on $\nonnegints$.
		That is, $J = [a,b)$ for some $b \in \nonnegints$ or $J = [a,\infty)$.
		Choose a nonzero $\sigma_a \in V_a$, which exists by assumption of $\rank(V_a) = 1$.
		Observe that $\set{\sigma_a}$ is a basis of $V_a$.
		Define the linear map $\phi_a: V_a \to \intmod[a]{J} = \field$ by $\sigma_a \mapsto 1_a$. 

		\item 
		Let $t \in J$ with $t \neq a$.
		Let $\sigma_t \in V_t$ be such that 
		$\alpha_{t,a}(\sigma_a) = \sigma_t$.
		By assumption, $\rank(V_t) = 1$ and $\rank(\alpha_{a,t}) = 1$.
		Then, $\sigma_t \neq 0$ and $\set{\sigma_t}$ is a basis for $V_t$.
		Define the linear map 
			$\phi_t: V_t \to \intmod[t]{J} = \field$ 
			by $\sigma_t \mapsto 1_t$.

		\item 
		Let $t \in \nonnegints$ with $t \not\in J$.
		By assumption, $\rank(V_t) = 0$ and $V_t$ is the trivial vector space.
		Then, $\phi_t: V_t \to \intmod[t]{J} = 0$ can only be the trivial map.
	\end{enumerate}
	Let $t,s \in \nonnegints$ with $t \leq s$.
	By \fref{defn:persmod-cat}(ii),
		we need to show that the following composition relation is satisfied: 
		$\phi_s \circ \alpha_{s,t} = i_{s,t}^J \circ \phi_t$
		where $\alpha_{s,t}: V_t \to V_s$ is a structure map of $V_\bullet$ 
		and $i_{s,t}^J: \intmod[t]{J} \to \intmod[s]{J}$ is that of $\intmod{J}$.
	We examine three cases:
	\begin{enumerate}
		\item 
		Assume $t,s \in J$. 
		Then, 
			$V_t$ and $V_s$ are non-trivial $\field$-vector spaces with bases $\set{\sigma_t}$ and $\set{\sigma_s}$ respectively.
		By \fref{lemma:persmod-functor-props}(ii),
			$\sigma_s =
			\alpha_{s,a}(\sigma_a)
			= \bigl( 
				\alpha_{s,t} \circ \alpha_{t,a}
			\bigr)(\sigma_a) 
			= \alpha_{s,t}(\sigma_t)$.
		Then, 
		\begin{equation*}
			\bigl( \phi_s \circ \alpha_{s,t} \bigr)(\sigma_t)
			= 
			\phi_s(\sigma_s)
			= 1_s 
			= i_{s,t}^J(1_t) 
			= \bigl( i_{s,t}^J \circ \phi_t \bigr)(\sigma_t)
		\end{equation*}
		Therefore, the relation $\phi_s \circ \alpha_{s,t} = i_{s,t}^J \circ \phi_t$ is satisfied.

		\item 
		Assume $t \not\in J$.
		Then, $V_t$ is the trivial vector space by assumption of $\rank(V_t) = 0$ and the relation $\phi_s \circ \alpha_{s,t} = i_{s,t}^J \circ \phi_t$ is trivially satisfied.

		\item 
		Assume $s \not\in J$.
		Then, $\intmod[s]{J} = 0$ and the maps $\phi_s: V_s \to \intmod[s]{J}$ and $i_{s,t}^J: \intmod[t]{J} \to \intmod[s]{J}$ are necessarily zero.
		Then, the composition relation  
		$\phi_s \circ \alpha_{s,t} = i_{s,t}^J \circ \phi_t = 0$ is trivially satisfied 
	\end{enumerate}
	Therefore, $\phi_\bullet: V_\bullet \to \intmod{J}$ is a persistence isomorphism and $V_\bullet \cong \intmod{J}$.
\end{proof}

As discussed in \fref{section:persistence-modules-as-functors},
	the structure maps of a persistence module $(V_\bullet, \alpha_\bullet)$ are uniquely determined by the collection $\set{\alpha_t: V_t \to V_{t+1}}_{t \in \nonnegints}$ of linear maps.
Let the interval $J \subseteq \nonnegints$ be given by $J = [a_k, b_k)$ for some $a_k \in \nonnegints$ and $b_k \in [a_k,\infty) \cup \set{\infty}$.
The condition \ref{eqn:int-mod-isom} in \fref{lemma:int-mod-isom} can equivalently be stated as follows:
\begin{equation}
	\rank(V_t) = \begin{cases}
		1 			&\text{ if } t \in [a,b) \\
		0 			&\text{ otherwise }
	\end{cases}
	\qquad\text{ and }\qquad 
	\rank(\alpha_{t}) = 
	\begin{cases}
		1		&\text{ if } t \in [a,b-1) \\
		0		&\text{ otherwise } 
	\end{cases}
	\label{eqn:easy-int-mod-isomorphism}
\end{equation}
where $[a,b-1)$ is interpreted to be $[a,\infty)$ if $b= \infty$.
Since the structure map $\alpha_{b-1}: V_{b-1} \to V_b$ has trivial codomain $V_b = 0$, it must also be trivial.
\fref{lemma:int-mod-isom} also implies that if $V_\bullet \cong \intmod{[a,b)}$, then any nonzero element $\sigma_a \in V_a$ will satisfy $\alpha_{t,a}(\sigma) \neq 0$ for all $t \in [a,b-1)$.
We illustrate this below relative to the sequence representation of $V_\bullet$, with the additional assumption that $b \in \nonnegints$:
\begin{displaymath}
\begin{tikzcd}[column sep=3em]
	0 			\arrow[r, equals]
	&[-23pt] V_{a-1} 	\arrow[r, "\alpha_{a-1}"] 
	& V_a 		\arrow[r, "\,\alpha_a\,"]
	& V_{a+1} 	\arrow[r, "\,\alpha_{a+1}\,"]
	& V_{a+2}	\arrow[r]
	&[-15pt] \cdots 	\arrow[r]
	&[-15pt] V_{b-1} 	\arrow[r, "\,\alpha_{b-1}\,"]
	& V_b 		\arrow[r, equals]
	&[-23pt] 0
	\\[-13pt] 
	&& \sigma_a 		\arrow[r, mapsto]
	& v_{a+1} 	\arrow[r, mapsto] 
			\arrow[equal, "/" marking]{d}
	& v_{a+2} 	\arrow[r, mapsto] 
			\arrow[equal, "/" marking]{d}
	& \cdots 	\arrow[r, mapsto]
	& v_{b-1}	\arrow[r, mapsto]
			\arrow[equal, "/" marking]{d}
	& 0	
	\\[-10pt]
	&&& 0 & 0 && 0
\end{tikzcd}\vspace{0\baselineskip} 
\end{displaymath} 
where $v_t := \alpha_{t,a}(\sigma_a)$ for $t \geq a+1$.
Observe that for all $t \in [a,b)$,
	$\set{v_t}$ is a basis for $V_t$ since $v_t \neq 0$ and $\rank(V_t) = 1$ by assumption.
For $t \geq b$,
	$v_t = 0$ and $V_t$ is the trivial vector space.
We give an example of a persistence module that is isomorphic to an interval module below.	

\begin{example}\label{ex:interval-mod-basis-change}
	Let $(Q_\bullet, \gamma_\bullet)$ be a persistence module over $\rationals$
	with vector spaces $Q_t$ and structure maps $\gamma_t: Q_t \to Q_{t+1}$ given as follows:
	\begin{equation*}
		\begin{aligned}[c]
			Q_t = \begin{cases}
				\rationalsket{a} 	&\text{ if } t = 3 \\
				\rationalsket{b}	&\text{ if } t = 4 \\
				\rationalsket{c}	&\text{ if } t = 5 \\
				0 					&\text{ if } t \not\in \set{3,4,5}
			\end{cases}
		\end{aligned}
		\quad\text{ and }\quad 
		\begin{array}{l @{\ }c@{\ } r@{\ }c@{\ }l}
			\gamma_3: Q_3 \to Q_4 
				&\text{ by }& 
				a &\mapsto& 2b 
			\\[2pt] 
			\gamma_4: Q_4 \to Q_5 
				&\text{ by }& 
				b &\mapsto& \frac{1}{3}c
			\\[2pt] 
			\gamma_5: Q_5 \to Q_6 
				&\text{ by }& 
				c &\mapsto& 0
		\end{array}
	\end{equation*}
	Note that for all $t \not\in \set{3,4,5}$, the structure map $\gamma_{t}: Q_t \to Q_{t+1}$ is necessarily the zero map since $Q_t = 0$.

	By \fref{lemma:int-mod-isom},
		$Q_\bullet \cong \intmod{J}$ 
		with $J = \set{3,4,5} = [3,6)$.
	We can also determine a corresponding persistence isomorphism $\phi_\bullet: Q_\bullet \to \intmod{[3,6)}$ using the same arguments presented in the proof of \fref{lemma:int-mod-isom}.
	Choose $a \in Q_3 = \rationals\ket{a}$.
	Then, the images of $a$ under the structure maps of $Q_\bullet$ is illustrated as follows:
	\vspace{0pt}
	\begin{displaymath}
	\begin{tikzcd}[column sep=3em]
		0 				\arrow[r, equals]
		&[-23pt] Q_2 	\arrow[r, "\gamma_{2}"] 
		& Q_3 			\arrow[r, "\gamma_3"]
		& Q_4 			\arrow[r, "\gamma_3"]
		& Q_5 			\arrow[r, "\gamma_3"]
		& Q_6 	\arrow[r, equals]
		&[-23pt] 0
		\\[-18pt] 
		& 
		& a \arrow[r, mapsto]
		& 2b \arrow[r, mapsto] 
		& {\textstyle\frac{2}{3}c} \arrow[r, mapsto]
		& 0
	\end{tikzcd}\vspace{0\baselineskip} 
	\end{displaymath}
	Observe that $\set{2b}$ and $\set{\frac{2}{3}c}$ are bases of $Q_4 = \rationals\ket{b}$ and $Q_5 = \rationals\ket{c}$ respectively.
	The linear maps $\phi_t: Q_t \to \intmod[t]{[3,6)}$ of the persistence isomorphism $\phi_\bullet$ are then given as follows:
	\begin{equation*}
		\begin{aligned}
			\phi_3: Q_3 = \rationals\ket{a}
				&\to \intmod[3]{[3,6)} = \rationals
			\\ 
			a &\mapsto 1
		\end{aligned}
		\quad,\quad 
		\begin{aligned}
			\phi_4: Q_4 = \rationals\ket{b}
				&\to \intmod[4]{[3,6)} = \rationals
			\\ 
			2b &\mapsto 1
		\end{aligned}
		\quad,\quad 
		\begin{aligned}
			\phi_5: Q_5 = \rationals\ket{c}
				&\to \intmod[5]{[3,6)} = \rationals
			\\ 
			{\textstyle\frac{2}{3}c} &\mapsto 1
		\end{aligned} 
	\end{equation*}
	Note that for all $t \not\in [3,6)$,
		the linear map $\phi_t: Q_t \to \intmod[t]{[3,6)}$ is the trivial map since $Q_t = 0$ and $\intmod[t]{[3,6)} = 0$.

	Let $\intmod{[3,6)} = (\intmod{[3,6)}, i_\bullet)$, i.e.\ 
		let $i_{s,t}: \intmod[t]{[3,6)} \to \intmod[s]{[3,6)}$ denote the structure maps of $\intmod{[3,6)}$.
	Since $\phi_\bullet$ is a persistence morphism, 
		the commutativity relation 
		$\phi_s \circ \gamma_{s,t} = i_{s,t} \circ \phi_t$
		is satisfied for all $t,s \in \nonnegints$ with $t \leq s$.
	This can be visualized using the following diagram, where highlighted in \redtag are elements of $\redmath{Q_t}$ and in \bluetag are those of $\bluemath{\intmod[t]{[3,6)}}$.

	\begin{displaymath}
	\begin{tikzcd}[column sep=-5pt]
		\mathllap{0=\ } Q_2
			\arrow[r, "\gamma_2"]
			\arrow[d, "\phi_2", swap]
		&[5em] Q_3 
			\arrow[rr, "\gamma_3"]
			\arrow[d, "\phi_3" swap]
			& \raisebox{-18pt}{$\smash{\redmath{a}}$}
			\arrow[rr, mapsto, yshift=-15pt, \redcolorname]
			\arrow[d, mapsto, shorten=7pt, yshift=3pt, \purplecolorname]
		&[3.5em] Q_4 
			\arrow[rr, "\gamma_4"]
			\arrow[d, "\phi_4" swap]
			& \raisebox{-15pt}{$\smash{\redmath{2b}}$}
			\arrow[rr, mapsto, yshift=-15pt, \redcolorname]
			\arrow[d, mapsto, shorten=7pt, yshift=3pt, \purplecolorname]
		&[3.5em] Q_5 
			\arrow[rr, "\gamma_5"]
			\arrow[d, "\phi_5" swap]
			& \raisebox{-15pt}{$\smash{\redmath{\textstyle\frac{2}{3}c}}$}
			\arrow[rr, mapsto, yshift=-15pt, \redcolorname]
			\arrow[d, mapsto, shorten=8.5pt, yshift=2pt, \purplecolorname]
		&[3.5em] Q_6\mathrlap{\ =0}
			\arrow[d, "\phi_6"]
			& \raisebox{-15pt}{$\smash{\redmath{0}}$}
			\arrow[d, mapsto, shorten=7pt, yshift=2.5pt, \purplecolorname]
		\\[15pt] 
		\mathllap{0=\ } \intmod[2]{[3,6)}
			\arrow[r, "0", swap]
		& \intmod[3]{[3,6)} 
			\arrow[rr, "\id_{\rationals}" swap]
			& \raisebox{20pt}{$\smash{\bluemath{1}}$}
			\arrow[rr, mapsto, yshift=20pt, \bluecolorname]
		& \intmod[4]{[3,6)} 
			\arrow[rr, "\id_{\rationals}" swap]
			& \raisebox{20pt}{$\smash{\bluemath{1}}$}
			\arrow[rr, mapsto, yshift=20pt, \bluecolorname]
		& \intmod[5]{[3,6)} 
			\arrow[rr, "0" swap]
			& \raisebox{20pt}{$\smash{\bluemath{1}}$}
			\arrow[rr, mapsto, yshift=20pt, \bluecolorname]
		& \intmod[6]{[3,6)}\mathrlap{\ =0}
			& \raisebox{17pt}{$\smash{\bluemath{0}}$}
	\end{tikzcd}
	\end{displaymath}
\end{example}

\spacer 

Next, we provide a definition for interval decompositions, adapted from \cite[Section 1.5]{persmod:chazal-structure}.

\begin{definition}\label{defn:interval-decomposition}
	An \textbf{interval decomposition} of a persistence module $(V_\bullet, \alpha_\bullet)$ over a field $\field$ is a finite direct sum $\bigoplus_k \intmod{J_k}$
	of interval modules where $\set{J_k}_{k=1}^m$ is some multiset of intervals in $\nonnegints$ such that 
	\begin{equation*}
		\persmod{V} 
			\cong 
			\bigoplus_{k=1}^m \, \intmod{J_k}
			\cong 
				\Big( \intmod{J_1} \Big)
				\oplus 
				\Big( \intmod{J_2} \Big)
				\oplus 
				\cdots
				\oplus 
				\Big( \intmod{J_m} \Big)
	\end{equation*}
	We say that $V_\bullet$ \textbf{admits an interval decomposition} if there exists an interval decomposition of $V_\bullet$.
\end{definition}
\remark{
	Some authors refer to $\set{J_k}$ as a collection instead of a multiset. 
	While we may conventionally use the term \textit{collection} to refer to a set or some set-like object, 
		we will use the term \textit{collection} to refer to multiset
		whenever we are talking about interval decompositions.
	For example, the collection $\set{[1,4), [1,4)}$ corresponds to the interval decomposition $\intmod{[1,4)} \oplus \intmod{[1,4)}$.
}

Observe that the definition of the term ``interval decomposition'' does not require that a persistence isomorphism $\phi_\bullet: V_\bullet \to \bigoplus_k \intmod{J_k}$ be explicitly given, only that at least one exists.
This convention is justified by the following uniqueness theorem.

\begin{theorem}\label{thm:uniqueness-interval-mods}
	Let $(V_\bullet, \alpha_\bullet)$ be a persistence module over $\field$.
	Assume that there exists two interval decompositions for $V_\bullet$ given as follows:
	\begin{equation*}
		\persmod{V} 
			\cong \bigoplus_{k \in K} \intmod{J_k}
		\qquad\text{ and }\qquad 
		\persmod{V}
			\cong \bigoplus_{m \in M} \intmod{L_m}
	\end{equation*}
	where $\set{J_k : k \in K}$ and $\set{L_m : m \in M}$ are two multisets of intervals in $\nonnegints$.
	Then,
	there exists a bijection $\pi: K \to M$ between the indexing sets $K$ and $M$ such that for all $k \in K$, $J_k = L_{\pi(k)}$ as intervals in $\nonnegints$.
\end{theorem}
\remark{
	A proof is available under \cite[Theorem 1.3]{persmod:chazal-structure}.
}

The existence of the bijection $\pi: K \to M$, as denoted in the theorem above, 
implies that the intervals of $\set{J_k}$ are exactly the intervals of $\set{L_m}$.
If the multisets $\set{J_k}$ and $\set{L_m}$ are linearly ordered,
	then the bijection $\pi: K \to M$ corresponds to a permutation or re-ordering of the intervals in $\set{J_k}$.
This explains why some authors describe interval decompositions as being unique up to permutation.

Furthermore, this implies that the multiset $\set{J_k}$ of intervals that determine the interval decomposition $\bigoplus_{k} \intmod{J_k}$ is unique up to persistence isomorphism.
This uniqueness result explains why we can use the article \textit{the} when talking about interval decompositions, 
	e.g.\ we talk of \textit{the} interval decomposition of a persistence module.
Note that the persistence isomorphism $V_\bullet \to \bigoplus_{k \in K} \intmod{J_k}$ behind the isomorphism relation $V_\bullet \cong \bigoplus_{k \in K} \intmod{J_k}$ is not generally unique.
This motivates the following terminology.

\begin{definition}\label{defn:persistence-barcode}
	The \textbf{persistence barcode} $\barcode(V_\bullet)$
	of a persistence module $(V_\bullet, \alpha_\bullet)$ is the multiset $\set{J_k}_{k=1}^m$ of intervals in $\nonnegints$ such that $\bigoplus_{k=1}^m \intmod{J_k}$ is an interval decomposition of $V_\bullet$, i.e.\ 
		$V_\bullet \cong \bigoplus_{k=1}^m \intmod{J_k}$.
	Note that $\barcode(V_\bullet)$ does not exist if $V_\bullet$ does not have an interval decomposition.
\end{definition}

The Structure Theorem (\fref{thm:structure-theorem}) for finitely-generated modules over a PID $R$ states that the invariant factors of a finitely-generated $R$-module are uniquely determined up to $R$-module isomorphism.
An equivalent way of saying this is that the invariant factors of an $R$-module (if an invariant factor decomposition exists) is an invariant of the isomorphism type of $R$-modules.

In the same vein, the persistence barcode (if it exists) is an invariant of the isomorphism type of persistence modules.
Furthermore,
	the persistence barcode is a concise characterization of the ranks of the vector spaces and structure maps of a persistence module.
We state this in more detail below.

\begin{proposition}\label{prop:interval-decomposition-and-ranks}
	Let $(V_\bullet, \alpha_\bullet)$ be a persistence module over $\field$.
	Assume $V_\bullet$ admits an interval decomposition 
	and let $\barcode(V_\bullet)$ be given by the multiset $\barcode(V_\bullet) = \set{J_k}_{k=1}^m$ of intervals.
	Then, 
	\begin{align}
		\rank(V_t) &= \card\Bigl\{
			J \in \barcode(V_\bullet) : t \in J
		\Bigr\} 
		\label{eqn:interval-decomposition-vector-space}
		\\
		\rank(\alpha_{s,t}) 
		&= \card\Bigl\{
			J \in \barcode(V_\bullet) : [t,s] = [t,s+1) \subseteq J
		\Bigr\} 
		\label{eqn:interval-decomposition-linear-map}
	\end{align}
	where $\rank(-)$ refers to the rank of $\field$-vector spaces and $\field$-linear maps
	and $\card(-)$ refers to the number of elements in a set, i.e.\ set cardinality.
\end{proposition}
\begin{proof}
	Assume that $V_\bullet \cong \bigoplus_{k=1}^m \intmod{J_k}$ for some multiset $\set{J_k}_{k=1}^m$ of intervals in $\nonnegints$.
	Then, there exists a persistence isomorphism $\phi_\bullet: V_\bullet \to \bigoplus_{k} \intmod{J_k}$
	with linear isomorphisms $\phi_t: V_t \to \bigoplus_{k=1}^m \intmod[t]{J_k}$ for $t \in \nonnegints$.

	\noindent
	For \ref{eqn:interval-decomposition-vector-space}:
	Let $t \in \nonnegints$.
	By definition of interval module,
		we have the following for each $ \in \set{1, \ldots, m}$,
	\begin{equation*}
		\rank\Bigl( \intmod[t]{J_k}
		\!\Bigr) = \left\{
		\!\!\begin{array}{c @{\ }c@{\ } c l}
			\rank(\field)
				&=& 1
				&\text{ if }
				t \in J_k 
			\\[2pt]
			\rank(0)
				&=& 0 
				&\text{ otherwise }
		\end{array}\right.
	\end{equation*}
	Since $\phi_t: V_t \to W_t$ is a linear isomorphism and by properties of the direct sum of vector spaces:
	\begin{equation*}
		\rank(V_t) 
		= \rank\bigl( \phi_t(V_t) \bigr)
		= \rank\Biggl(\,
			\bigoplus_{k=1}^m \intmod[t]{J_k}
		\!\Biggr)
		= \sum_{k=1}^m \,\rank\Bigl(
			\intmod[t]{J_k}
		\!\Bigr)
		= \card\Bigl\{
			J \in \barcode(V_\bullet) : t \in J
		\Bigr\}
	\end{equation*}
	Therefore, \fref{eqn:interval-decomposition-vector-space} is satisfied.

	\noindent
	For \ref{eqn:interval-decomposition-linear-map}:
		Let $t,s \in \nonnegints$ with $t \leq s$.
		For each $k \in \set{1, \ldots, m}$,
			let $\gamma_{s,t}^{\,k}$ refer to the structure map $\gamma_{s,t}^{\,k}: \intmod[t]{J_k} \to \intmod[s]{J_k}$
			of the interval module $\intmod{J_k}$.
		Note that $t,s \in J_k$ if and only if $[t,s] \subseteq J_k$.
		Therefore,
		\begin{equation*}
			\rank(\gamma_{s,t}^{\,k})
			= \left\{\begin{array}{c @{\ }c@{\ } cl}
				\rank(\id_{\field})
					&=& 1
					&\text{ if }
					[t,s] \subseteq J_k 
				\\[2pt]
				\rank(0)
					&=& 0 
					&\text{ otherwise }
			\end{array}\right.
		\end{equation*}
		Let $\gamma_{s,t}: \bigoplus_{k=1}^m \intmod[t]{J_k} 
		\to \bigoplus_{k=1}^m \intmod[s]{J_k}$ be the structure map of $\bigoplus_{k=1}^m \intmod{J_k}$ on $t \to s$.
		By definition of the direct sum of persistence modules,
			$\gamma_{s,t} = \bigoplus_{k=1}^m \gamma_{s,t}^{\,k}$
			where $\bigoplus_{k=1}^m \gamma_{s,t}^{\,k}$ is the linear map induced by the direct sum operation between $\field$-vector spaces.
		Since $\phi_\bullet$ must be a persistence morphism,
			the composition relation $\gamma_{s,t} \circ \phi_t = \phi_s \circ \alpha_{s,t}$ is satisfied.
		Therefore, 
		\begin{align*}
			\rank(\alpha_{s,t})
			&= \rank\Bigl(
				(\phi_s)\inv \circ \gamma_{s,t} \circ \phi_t
			\Bigr)
			= \rank(\gamma_{s,t})
			= \rank\Biggl(\,\,
				\bigoplus_{k=1}^m \gamma_{s,t}^{\,k}
			\Biggr)
			\\
			&= \sum_{k=1}^m \,\rank\bigl(
				\gamma_{s,t}^{\,k} \bigr)
			= \card\Bigl\{
				J \in \barcode(V_\bullet) : [t,s] = [t,s+1) \subseteq J
			\Bigr\}
		\end{align*}
		Therefore, \fref{eqn:interval-decomposition-linear-map} is also satisfied.
\end{proof}

It can be proven that finite-type persistence modules admit interval decompositions.
We state this later in \fref{prop:finite-type-has-interval} under \fref{section:cat-equiv-graded-modules} and prove it using an equivalence of categories.

The calculation of the interval decomposition of a given finite-type persistence module $(V_\bullet, \alpha_\bullet)$ requires that both the vector spaces $V_t$ and the structure maps $\alpha_{s,t}: V_t \to V_s$ for all $t \in \nonnegints$ be considered, for lack of a better term, \textit{simultaneously} for all indices $t \in \nonnegints$. 
In practice, this means that the decomposition of the vector spaces $V_t$ must be \textit{compatible} with the structure maps, in that the collection of vector space isomorphisms involving each $V_t$ must form a persistence morphism.
We give an example of this below.

\begin{example}\label{ex:easy-int-decomp}
	Let ($\persmod{V}, \gamma_\bullet)$ be a persistence module over $\ints_5$ 
	with vector spaces given as follows:
	\begin{equation*}
		V_t = \begin{cases}
			\intsket[5]{x_1} 			&\text{ if } t \in [0,13) \\
			\intsket[5]{y_1, y_2}		&\text{ if } t \in [13, 21) \\
			\intsket[5]{z_1, z_2, z_3}	&\text{ if } t \in [21, \infty)
		\end{cases}
	\end{equation*}
	For $t \neq 12$ and $t \neq 20$, let $\gamma_t: V_t \to V_{t+1}$ be the identity map on $V_t$, which is well-defined since $V_t = V_{t+1}$.
	Define the structure map $\gamma_t: V_t \to V_{t+1}$ of $V_\bullet$ for $t=12$ and $t=20$ as follows:
	\begin{equation*}
		\begin{aligned}[t]
			\gamma_{12}: V_{12} &\to V_{13} \\
						x_1 &\mapsto y_2
		\end{aligned}
		\qquad\qquad
		\begin{aligned}[t]
			\gamma_{20}: V_{20} &\to V_{21} \\
						y_1 &\mapsto z_1 \\
						y_2 &\mapsto z_3
		\end{aligned}
	\end{equation*}
	To determine the interval decomposition of $C_\bullet$, 
	we need to find a decomposition of $C_\bullet$ into persistence submodules such each submodule is isomorphic to an interval module.
	For each $t \in \nonnegints$, 
		$A_t$, $B_t$, and $C_t$ be vector subspaces of $V_t$ given as follows:
	\begin{equation*}
		A_t = \begin{cases}
			\ints_5\ket{x_1} &\text{ if } t \in [0,13) \\
			\ints_5\ket{y_2} &\text{ if } t \in [13, 21) \\
			\ints_5\ket{z_3} &\text{ if } t \in [21, \infty)
		\end{cases}
		\qquad 
		B_t = \begin{cases}
			0 &\text{ if } t \in [0,13) \\
			\ints_5\ket{y_1} &\text{ if } t \in [13, 21) \\
			\ints_5\ket{z_1} &\text{ if } t \in [21, \infty)
		\end{cases}
		\qquad 
		C_t = \begin{cases}
			0 &\text{ if } t \in [0,13) \\
			0 &\text{ if } t \in [13, 21) \\
			\ints_5\ket{z_2} &\text{ if } t \in [21, \infty)
		\end{cases}
	\end{equation*}
	Observe that, for all $t \in \nonnegints$, $A_t$ is a $\ints_5$-vector subspace of $V_t$ and that the structure map $\gamma_t: V_t \to V_{t+1}$ of $V_\bullet$ satisfies $\gamma_t(A_t) \subseteq A_{t+1} \subseteq V_{t+1}$.
	Therefore, the collection $\set{A_t}_{t \in \nonnegints}$ of vector spaces determines a persistence submodule $A_\bullet$ of $V_\bullet$ with $A_\bullet(t) = A_t$ and 
	the structure map $A_t \to A_{t+1}$ defined by restricting the domain and codomain of the structure map $\gamma_t: V_t \to V_{t+1}$ of $V_\bullet$ for all $t \in \nonnegints$.
	Similarly, the collections $\set{B_t}_{t \in \nonnegints}$ and $\set{C_t}_{t \in \nonnegints}$ of $\ints_5$-vector spaces determine persistence submodules $B_\bullet$ and $C_\bullet$ of $V_\bullet$ respectively.

	Since for each $t \in \nonnegints$, $V_t \cong A_t \oplus B_t \oplus C_t$ as $\ints_5$-vector spaces (with the direct sum interpreted as an internal direct sum of vector spaces), 
		we have the following decomposition of $V_\bullet$ as a persistence module:
	\begin{equation*}
		V_\bullet \cong A_\bullet \oplus B_\bullet \oplus C_\bullet 
	\end{equation*}
	with the persistence isomorphism $i_\bullet: A_\bullet \oplus B_\bullet \oplus C_\bullet \to V_\bullet$ given by the inclusion maps $i_t: A_t \oplus B_t \oplus C_t \hookrightarrow V_t$ for all $t \in \nonnegints$.
	This decomposition of $V_\bullet$ is illustrated in the following diagram:
	\begin{displaymath}
	\begin{tikzcd}[column sep=14pt, every label/.append style = {font = \small}]
		{} &[-15pt]{}&[-15pt]
				(t=0) &  \cdots &  (t=12) &
				(t=13) &  \cdots &  (t=20) & 
				(t=21) &  \cdots
		\\[-15pt]
		V_\bullet &:&
			\ints_5\ket{x_1}\arrow[r]
			& \cdots \arrow[r]
			& \ints_5\ket{x_1}\arrow[r]
			& \ints_5\ket{y_1, y_2}	\arrow[r]
			& \cdots \arrow[r]
			& \ints_5\ket{z_1, z_2, z_3}\arrow[r]
			& \ints_5\ket{z_1, z_2, z_3}\arrow[r]
			& \cdots
		\\[-20pt]
		\sideequals 
			&& \sideequals && \sideequals & \sideequals && \sideequals & \sideequals
		\\[-20pt]
		\bluemath{A_\bullet} &:& 
			\bluemath{\ints_5\ket{x_1}}\arrow[r]
			& \bluemath{\cdots\hspace{-1pt}}\arrow[r]
			& \bluemath{\ints_5\ket{x_1}}\arrow[r]
			& \bluemath{\ints_5\ket{y_2}}\arrow[r]
			& \bluemath{\cdots\hspace{-1pt}}\arrow[r]
			& \bluemath{\ints_5\ket{y_2}}\arrow[r]
			& \bluemath{\ints_5\ket{z_3}}\arrow[r]
			& \bluemath{\cdots\hspace{-1pt}}
		\\[-20pt]
		\oplus 
			&& && & \oplus && \oplus & \oplus
		\\[-20pt]
		\greenmath{B_\bullet} &:& 
			&& 
			& \greenmath{\ints_5\ket{y_1}}\arrow[r]
			& \greenmath{\cdots\hspace{-1pt}}\arrow[r]
			& \greenmath{\ints_5\ket{y_1}}\arrow[r]
			& \greenmath{\ints_5\ket{z_1}}\arrow[r]
			& \greenmath{\cdots\hspace{-1pt}}
		\\[-20pt]
		\oplus 
			&&  && & && & \oplus
		\\[-20pt]
		\redmath{C_\bullet} &:& 
			&& &&&
			& \redmath{\ints_5\ket{z_2}}\arrow[r]
			& \redmath{\cdots\hspace{-1pt}}
	\end{tikzcd}
	\end{displaymath}
	From the diagram above, we can conclude that the persistence submodules $A_\bullet$, $B_\bullet$, and $C_\bullet$ of $V_\bullet$ are isomorphic to interval modules as follows:
	\begin{enumerate}[left=1in]
		\item[For $A_\bullet \cong \intmod{[0,\infty)}$\,:] 
		Let $A_\bullet = (A_\bullet, \alpha_\bullet)$, i.e.\ denote the structure maps of $A_\bullet$ by $\alpha_t: A_t \to A_{t+1}$.
		Recall that $\alpha_t(a) = \gamma_t(a)$ for all $a \in A_t$ and for all $t \in \nonnegints$.
		Observe that for all $t \in \nonnegints = [0,\infty)$, 
			$\rank(A_t) = 1$ and $\rank(\alpha_t) = 1$.
			Then, $A_\bullet \cong \intmod{[0,\infty)}$ by \fref{lemma:int-mod-isom}, relative to characterization by \fref{eqn:easy-int-mod-isomorphism}.
		
		Alternatively, we can construct the persistence isomorphism $\phi^A_\bullet: A_\bullet \to \intmod{[0,\infty)}$ 
		explicitly as follows:
		Denote the multiplicative identity of $\ints_5$ as $1$.
		Then, for each $t \in \nonnegints$, define the linear map 
			$\phi^A_t: A_t \to \intmod[t]{[0,\infty)}$ as follows:
		\begin{enumerate}[align=parleft, labelwidth=0.9in, leftmargin=1.3in]
			\item[If $t \in [0,13)$:]
				Let $\phi^A_t: A_t = \ints_5\ket{x_1} \to \intmod[t]{[0,\infty)} = \ints_5$ be given by $x_1 \mapsto 1$.

			\item[If $t \in [13,21)$:]
				Let $\phi^A_t: A_t = \ints_5\ket{y_2} \to \intmod[t]{[0,\infty)} = \ints_5$ be given by $y_2 \mapsto 1$.

			\item[If $t \in [21,\infty)$:]
				Let $\phi^A_t: A_t = \ints_5\ket{z_3} \to \intmod[t]{[0,\infty)}$ be given by $z_3 \mapsto 1$.
		\end{enumerate}
		It should be straightforward to check that the linear maps $\phi^A_t: A_t \to \intmod[t]{[0,\infty)}$ commute with the structure maps $\alpha_t: A_t \to A_{t+1}$ of $A_\bullet$ and those of $\intmod{[0,\infty)}$.
		Therefore, $\phi^A_\bullet$ by $\phi^A_\bullet = \set{\phi_t}_{t \in \nonnegints}$ is a well-defined persistence morphism. Since each $\phi^A_t$ is a vector space isomorphism, $\phi^A_\bullet$ is a persistence isomorphism.

		\item[For $B_\bullet \cong \intmod{[13,\infty)}$\,:]  
		Let $B_\bullet = (B_\bullet, \beta_\bullet)$, i.e.\ denote the structure maps of $B_\bullet$ by $\beta_t: B_t \to B_{t+1}$.
		Observe that if $t \in [13,\infty)$, $\rank(B_t) = 1$ and $\rank(\beta_t) = 1$. If $t \in [0,12)$, $B_t$ is trivial and $\rank(B_t) = 0$ and $\rank(\beta_t) = 0$.
		Then, $B_\bullet \cong \intmod{[13,\infty)}$ by \fref{lemma:int-mod-isom}.

		The persistence isomorphism $\phi^B_\bullet: B_\bullet \to \intmod{[13,\infty)}$ can be constructed similarly as with the case of $\phi^A_\bullet$, i.e.\ 
		\begin{enumerate}[align=parleft, labelwidth=0.9in, leftmargin=1.3in]
			\item[If $t \in [0,13)$:]
				$B_t$ is trivial and $\phi^B_t: B_t \to \intmod[t]{[13,\infty)} = 0$ can only be the trivial map. 

			\item[If $t \in [13,21)$:]
				The map $\phi^B_t: B_t = \ints_5\ket{y_1} \to \intmod[t]{[13,\infty)} = \ints_5$ is given by $y_1 \mapsto 1$.
			
			\item[If $t \in [21,\infty)$:]
				The map $\phi^B_t: B_t = \ints_5\ket{z_1} \to \intmod[t]{[13,\infty)} = \ints_5$ is given by 
				$z_1 \mapsto 1$.
		\end{enumerate}

		\item[For $B_\bullet \cong \intmod{[13,\infty)}$\,:] 
			Let $C_\bullet = (C_\bullet, \psi_\bullet)$.
			Note that if $t \in [21,\infty)$, $\rank(C_t) = 1$ and $\rank(\psi_t: C_t \to C_{t+1}) = 1$.
			If $t \in [0,21)$, $\rank(C_t) = 0$ and $\rank(\psi_t: C_t \to C_{t+1}) = 0$.
			By \fref{lemma:int-mod-isom}, $C_\bullet \cong \intmod{[21,\infty)}$.

			The persistence isomorphism $\phi^C_\bullet: C_\bullet \to \intmod{[21,\infty)}$ is given as follows:
			If $t \in [21,\infty)$, $\phi^C_t: C_t = \ints_5\ket{z_2} \to \intmod[t]{[21,\infty)} = \ints_5$ is given by 
			$z_2 \mapsto 1$. If $t \in [0,21)$, $C_t = 0$ and $\phi^C_t$ can only be the trivial map.
	\end{enumerate}
	Therefore, we have the following decomposition for $V_\bullet$:
	\begin{equation*}
		V_\bullet 
		\cong 
		A_\bullet \oplus B_\bullet \oplus C_\bullet 
		\cong 
		\intmod{[0,\infty)} 
		\oplus \intmod{[13,\infty)} \oplus \intmod{[21,\infty)}
	\end{equation*}
	Consequently, the persistence barcode $\barcode(V_\bullet)$ of $V_\bullet$ is given by $\barcode(V_\bullet) 
	= \bigl\{
		[0,\infty), [13,\infty), [21,\infty)
	\bigr\}$.
\end{example}

Observe that, in the above example, the persistence barcode of $V_\bullet$ corresponds to a collection of bases for each $V_t$ that is compatible with the structure maps.
In particular, \fref{ex:easy-int-decomp} is written such that elements of the bases $\set{x_1, x_2, x_3}$, $\set{y_1, y_2}$, and $\set{z_1, z_2, z_3}$ for $V_t$ with $t \in [0,13)$, $t \in [13,21)$, and $t \in [21,\infty)$ respectively are mapped either to $1 \in \ints_5$ or $0 \in \ints_5$ in the given interval decomposition for $V_\bullet$.

Unfortunately, this nice correspondence does not generally apply to arbitrary finite-type persistence modules. 
In the example below, we define the vector spaces of a persistence module using some set of indeterminates but show an interval decomposition of said persistence modules that requires a change of basis for the vector spaces involved.

\begin{example}\label{ex:int-mod-intsthree}
	Let $(W_\bullet, \gamma_\bullet)$ be a persistence module over $\ints_3$ with vector spaces given as follows:
	\begin{equation*}
		W_t = \begin{cases}
			\intsket[3]{a_1} 	
				&\text{ if } t=0 \\
			\intsket[3]{b_1, b_2}
				&\text{ if } t=1 \\
			\intsket[3]{c_1, c_2}
				&\text{ if } t=2 \\			
			\intsket[3]{d_1}
				&\text{ if } t \geq 3 \\
		\end{cases}
	\end{equation*}
	Let the structure maps of $W_\bullet$ of the form $\gamma_t: W_t \to W_{t+1}$ for $t \in \nonnegints$ be given by 
	\begin{align*}
		\gamma_0: W_0 = \intsket[3]{a_1} &\to W_1 &
		\gamma_1: W_1 &\to W_2 &
		\gamma_2: W_2 &\to W_3 &
		\gamma_t: W_t &\to W_{t+1} \text{ for all } t \geq 3
		\\
		a_1 &\mapsto 2b_1 + b_2 &
		b_1 &\mapsto c_1 + c_2 &
		c_1 &\mapsto 0 &
		d_1 &\mapsto d_1
		\\
		&&
		b_2 &\mapsto c_1 + 2c_2 &
		c_2 &\mapsto d_1
	\end{align*}
	To determine the interval decomposition of $W_\bullet$, we must first find a decomposition of $W_\bullet$ into persistence submodules such that each submodule is isomorphic to an interval module.
	For each $t \in \nonnegints$, define the $\ints_3$-vector spaces $U_t$ and $V_t$ as follows:
	\begin{equation*}
		A_t = \begin{cases}
			\ints_3\ket{a_1} &\text{ if } t = 0 \\
			\ints_3\ket{2b_1 + b_2} &\text{ if } t = 1 \\
			\ints_3\ket{c_2} &\text{ if } t=2 \\
			\ints_3\ket{d_1} &\text{ if } t \geq 3
		\end{cases}
		\qquad\text{ and }\qquad
		B_t = \begin{cases}
			0 &\text{ if } t = 0 \\
			\ints_3\ket{b_1 + b_2} &\text{ if } t = 1 \\
			\ints_3\ket{2c_1} &\text{ if } t=2 \\
			0 &\text{ if } t \geq 3
		\end{cases}
	\end{equation*}
	For the collections $\set{A_t}_{t \in \nonnegints}$ and $\set{B_t}_{t \in \nonnegints}$ to determine persistence submodules $A_\bullet = (A_\bullet, \alpha_\bullet)$ and $B_\bullet = (B_\bullet, \beta_\bullet)$ of $W_\bullet$ respectively, 
		we must show that linear maps $\alpha_t: A_t \to A_{t+1}$ and $\beta_t: B_t \to B_{t+1}$ obtained by appropriately restricting the domain and codomain of $\gamma_{t}: V_t \to V_{t+1}$ are well-defined.
	This is shown by the following calculations:
	\begin{enumerate}[left=1in]
		\item[For $\set{A_t}_{t \in \nonnegints}$\,:] 
		It suffices to consider the images of the basis elements $a_1$, $2b_1+b_2$, and $c_2$ of $A_0$, $A_1$, and $A_2$ respectively under the appropriate structure map $\gamma_{t}$ of $W_\bullet$:
		\begin{equation*}\def\arraystretch{1.2}
		\begin{array}{c @{\ }c@{\ } c @{\ }c@{\ } c}
			\gamma_0(a_1)
			&=& \cdots 
			&=& 2b_1 + b_2
			\\
			\gamma_1(2b_1+b_2)
			&=& 2(c_1 + c_2) + (c_1 + 2c_2) 
				= 2c_1 + 2c_2 + c_1 + 2c_2
			&=& c_2
			\\ 
			\gamma_2(c_2)
			&=& \cdots 
			&=& d_1 
		\end{array}
		\end{equation*}
		Note that for all $t \geq 3$, $\gamma_3$ is the identity map on $V_t = \ints_3\ket{d_1}$.
		Then, for all $t \in \nonnegints$, $\gamma_t(A_t) \subseteq A_{t+1}$ 
		and the linear map $\alpha_t:A_t \to A_{t+1}$ obtained by restricting $\gamma_t: W_t \to W_{t+1}$ is well-defined.
		Therefore, $A_\bullet = (A_\bullet,\alpha_\bullet)$ by $A_\bullet(t) = A_t$ and $A_\bullet(t \to t+1) = \alpha_t$ is a persistence submodule of $W_\bullet$.

		\item[For $\set{B_t}_{t \in \nonnegints}$\,:] 
		Note that for $t \neq 1,2$, $B_t$ is trivial.
		It suffices to consider the image of $b_1 + b_2 \in B_1$ and $2c_1 \in B_2$ under the appropriate structure map $\gamma_t$ of $W_\bullet$:
		\begin{equation*}\def\arraystretch{1.2}
		\begin{array}{c @{\ }c@{\ } c @{\ }c@{\ } c}
			\gamma_1(b_1+b_2)
			&=& (c_1 + c_2) + (c_1 + 2c_2) = 2c_1 + 3c_2
			&=& 2c_1
			\\
			\gamma_2(2c_1)
			&=& 2(0)
			&=& 0
		\end{array}
		\end{equation*}
		Then, $\gamma_1(B_1) \subseteq \ints_3\ket{2c_1} = B_2$ and 
		$\gamma_2(B_2) \subseteq \ints_3\ket{0} = B_3$,
		and the linear maps $\beta_1: B_1 \to B_2$ and $\beta_2: B_1 \to B_2$ obtained by restricting $\gamma_1: W_1 \to W_2$ and $\gamma_2: W_2 \to W_3$ are well-defined.
		For $t \neq 1,2$, define $\beta_t: B_t \to B_{t+1}$ to be the trivial map.
		Therefore, $B_\bullet = (B_\bullet, \beta_\bullet)$ by $B_\bullet(t) = B_t$ and $B_\bullet(t \to t+1) = \beta_t$ is a persistence submodule of $W_\bullet$.
	\end{enumerate}
	Observe that for $t=0$ and for $t \geq 3$,
		$A_t \oplus B_t \cong A_t = V_t$ since $B_t$ is trivial.
	For the case of $t=2$,
		we have that $\set{2b_1 + b_2, b_1 + b_2}$ is a basis for $W_2 = \ints_3\ket{b_1, b_2}$ by the following calculation:
	\begin{equation*}\def\arraystretch{1.2}
	\begin{array}{c @{\ }c@{\ } c @{\ }c@{\ } c  @{\ }c@{\ } c}
		(1)(2b_1 + b_2) + (2)(b_1 + b_2)
		&=& 2b_1 + b_2 + 2b_1 + 2b_2 
		&=& 4b_1 + 3b_2
		&=& b_1
		\\
		(2)(2b_1 + b_2) + (2)(b_1 + b_2)
		&=& 4b_1 + 2b_2 + 2b_1 + 2b_2 
		&=& 6b_1 + 4b_2 
		&=& b_2
	\end{array}
	\end{equation*}
	Similarly, $\set{2c_1, c_2}$ is a basis for $W_3 = \ints_3\ket{c_1, c_2}$ since $2(2c_1) = 4c_1 = c_1$.
	Therefore, $W_t = A_t \oplus B_t$ for all $t \in \nonnegints$ (taken as an internal direct sum) and we have the following persistence isomorphism relation:
	\begin{equation*}
		W_\bullet \cong A_\bullet \oplus B_\bullet
	\end{equation*}
	This is illustrated in the diagram below:
	\begin{displaymath}
	\begin{tikzcd}[column sep=14pt, every label/.append style = {font = \small}]
		{} &[-15pt]{}&[-15pt]
				(t=0) & (t=1) &
				(t=2) & (t=3) & \cdots
		\\[-15pt]
		W_\bullet &:
			& \ints_5\ket{a_1}\arrow[r]
			& \ints_5\ket{b_1, b_2}\arrow[r]
			& \ints_5\ket{c_1, c_2}	\arrow[r]
			& \ints_5\ket{d_1}\arrow[r]
			& \cdots
		\\[-20pt]
		\sideequals 
			&& \sideequals & \sideequals & \sideequals & \sideequals
		\\[-20pt]
		\bluemath{A_\bullet} &:
			& \bluemath{\ints_5\ket{a_1}}\arrow[r]
			& \bluemath{\ints_5\ket{2b_1+b_2}}\arrow[r]
			& \bluemath{\ints_5\ket{c_2}}\arrow[r]
			& \bluemath{\ints_5\ket{d_1}}\arrow[r]
			& \bluemath{\cdots\hspace{-1pt}}
		\\[-20pt]
		\oplus 
			&& \oplus & \oplus & \oplus & \oplus 
		\\[-20pt]
		\greenmath{B_\bullet} &:
			& 0 \arrow[r]
			& \greenmath{\ints_3\ket{b_1+b_2}}\arrow[r]
			& \greenmath{\ints_3\ket{2c_1}}\arrow[r]
			& 0 \arrow[r]
			& \cdots
	\end{tikzcd}
	\end{displaymath}
	The submodules $A_\bullet$ and $B_\bullet$ of $W_\bullet$ are also isomorphic to interval modules, as described below:
	\begin{enumerate}[left=1in]
		\item[For $A_\bullet \cong \intmod{[0,\infty)}$\,:] 
		The persistence isomorphism $\phi^A_\bullet: A_\bullet \to \intmod{[0,\infty)}$ is given as follows:
		\begin{enumerate}[align=parleft, labelwidth=0.7in, leftmargin=1in]
			\item[For $t=0$:]
				The map $\phi^A_0: A_0
				\to \intmod[0]{[0,\infty)} = \ints_3$ is given by 
				$a_1 \mapsto 1 \in \ints_3$.

			\item[For $t=1$:]
				The map $\phi^A_1: A_1 \to \intmod[1]{[0,\infty)} = \ints_3$ is given by 
				$2b_1 + b_2 \mapsto 1 \in \ints_3$.

			\item[For $t=2$:]
				The map $\phi^A_2: A_2 \to \intmod[2]{[0,\infty)} = \ints_3$ is given by 
				$c_2 \mapsto 1 \in \ints_3$.

			\item[For $t \geq 2$:]
				The map $\phi^A_t: A_t \to \intmod[t]{[0,\infty)} = \ints_3$ is given by 
				$d_1 \mapsto 1 \in \ints_3$.
		\end{enumerate}
		We can confirm that the maps $\phi^A_t: A_t \to \intmod[t]{[0,\infty)}$ commutes with the structure maps $\alpha_{t}: A_t \to A_{t+1}$ of $A_\bullet$ and those of $\intmod{[0,\infty)}$.
		We can also conclude that $A_\bullet \cong \intmod{[0,\infty)}$ using \fref{lemma:int-mod-isom} since $\rank(A_t)=1$ and $\rank(\alpha_t: A_t \to A_{t+1}) = 1$ for all $t \in \nonnegints$.

		\item[For $B_\bullet \cong \intmod{[1,3)}$\,:] 
		The persistence isomorphism $\phi^B_\bullet: B_\bullet \to \intmod{[1,3)}$ is given as follows:
		\begin{enumerate}[align=parleft, labelwidth=0.7in, leftmargin=1in]
			\item[For $t=0$:]
				$B_0$ is trivial and $\phi^B_0: B_0 \to \intmod{[1,3)} = 0$ is necessarily the trivial map.

			\item[For $t=1$:]
				The map $\phi^B_1: B_1 \to \intmod[1]{[1,3)} = \ints_3$ is given by 
				$b_1 + b_2 \mapsto 1 \in \ints_3$.

			\item[For $t=2$:]
				The map $\phi^B_2: B_2 \to \intmod[2]{[1,3)} = \ints_3$ is given by 
				$2c_1 \mapsto 1 \in \ints_3$.

			\item[For $t \geq 2$:]
				$B_t$ and $\intmod[t]{[1,3)}$ are both trivial and 
				$\phi^B_t$ is necessarily the trivial map.
		\end{enumerate}
		Observe that the maps $\set{\phi^B_t}_{t \in \nonnegints}$ commute with the structure maps of $B_\bullet$ and those of $\intmod{[1,3)}$. In particular, we have the following for the case of $t=1$ and $t=2$:
		\begin{equation*}\def\arraystretch{1.2}
		\begin{array}{r @{\ }c@{\ } c @{\ }c@{\ } l}
			\Bigl(
				\phi^B_2 \mathrel{\circ} \beta_1
			\Bigr)(b_1 + b_2)
			&=& \phi^B_2(2c_1)
			= 1
			= \id_{\ints_3}(1)
			&=& \Bigl(
				\id_{\ints_3} \mathrel{\circ} \phi^B_1
			\Bigr)(b_1 + b_2)
			\\[5pt]
			\Bigl(
				\phi^B_3 \mathrel{\circ} \beta_2
			\Bigr)(2c_1)
			&=& \phi^B_3(0)
			= 0
			= \id_{\ints_3}(0)
			&=& \Bigl(
				\id_{\ints_3} \mathrel{\circ} \phi^B_2
			\Bigr)(2c_1)
		\end{array}
		\end{equation*} 
		Alternatively, we can use \fref{lemma:int-mod-isom} (relative to its characterization by \fref{eqn:easy-int-mod-isomorphism}) to determine that $B_\bullet \cong \intmod{[1,3)}$ since we have the following:
		\begin{equation*}
			\rank(B_t) = \begin{cases}
				1 &\text{ if } t \in [1,3) \\
				0 &\text{ otherwise }
			\end{cases}
			\qquad\text{ and }\qquad 
			\rank(\beta_t) = \begin{cases}
				1 &\text{ if } t \in [1,3-1) = 1 \\
				0 &\text{ otherwise }
			\end{cases}
		\end{equation*}
	\end{enumerate}
	Therefore, we have the following persistence isomorphism relation:
	\begin{equation*}
		W_\bullet \cong A_\bullet \oplus B_\bullet 
		\cong \intmod{[0,\infty)} \oplus \intmod{[1,3)}
	\end{equation*}
	and the persistence barcode of $W_\bullet$ is given by 
	$\barcode(W_\bullet) = \bigl\{
		[0,\infty), [1,3)
	\bigr\}$.
\end{example}

In \fref{chapter:matrix-calculation}, 
we discuss a method of calculating the interval decomposition of finite-type persistence modules that are the $n$\th chain homology of a chain complex of finite-type persistence modules.
 \clearpage

\section{The Category of Graded Modules over Polynomial Rings}\label{section:graded-mod-notation}

Let $R$ be a commutative ring with identity $1_R \in R$.
One of the key results in persistence theory is the equivalence between the category $\catpersmod$ of persistence modules and the category $\catgradedmod{R}$ of graded $R[x]$-modules, discussed in \fref{section:cat-equiv-graded-modules}.

However, we have found that definitions and descriptions for graded modules, particularly those involving the \textit{category} of graded modules, seem to be sparse in introductory-level abstract algebra texts, 
and that rigorous definitions often only appear in more specialized fields such as commutative algebra and homological algebra.
To avoid confusion, we identify a number of definitions, notation, and results involving graded $R[x]$-modules in this section that are used throughout the paper.
We use the following texts as our primary references.
\begin{enumerate}
	\item 
	\textit{Graded Syzygies} \cite{algebra:peeva} by Irena Peeva.

	\item 
	\textit{Methods of Graded Rings} \cite{algebra:nastasescu} by Constantin \Nastasescu{} and Freddy Van Oystaeyen.
\end{enumerate}
To start, most of the literature regarding graded modules and persistence modules use the term \textit{action} to refer to the scalar multiplication operation on a module.
We state a definition of \textit{action} below.

\begin{definition}\label{defn:action-on-module}
	The \textbf{action} of $R$ on an $R$-module $M$ is a biadditive group action $\cdot: R \times M \to M$ that defines the scalar multiplication operation on $M$.
\end{definition}
\noindent 
Note that the requirements that the \textit{action} $\cdot: R \times M \to M$ be a group action and be biadditive satisfy the usual conditions for scalar multiplication on an $R$-module. We state this in more detail below:
\begin{enumerate}
	\item 
	The condition that $\cdot: R \times M \to M$ be a \textit{group action} requires that $\cdot$ must satisfy the identity and compatibility axioms with respect to the multiplication operation $\cdot_R: R \times R \to R$ on $R$ (as a ring), i.e.\ for all $m \in M$ and $r,s \in R$,
	\begin{equation*}
		1_R \cdot m = m 
		\qquad\text{ and }\qquad 
		r \cdot (s \cdot m) = (r \cdot_R s) \cdot m
	\end{equation*}

	\item 
	The condition that $\cdot: R \times M \to M$ be \textit{biadditive} requires that $\cdot$ be linear on the first and second arguments
	with respect to addition $+_R: R \times R \to R$ on $R$ (as a ring) and to addition $+_M: M \times M \to M$ on $M$ (as an abelian group), i.e.\ for all $r,s \in R$ and $m,n \in M$,
	\begin{align*}
		(r +_R s) \cdot m &= (r \cdot m) +_M (s \cdot m) \\
		r \cdot (m +_M n) &= (r \cdot m) +_M (r \cdot n)
	\end{align*}
	The two equalities above are sometimes called the distributivity properties of scalar multiplication (of modules).
\end{enumerate}
Next, we provide definitions for the standard grading on $R[x]$ taken from \cite[Section 1]{algebra:peeva}.

\begin{definition}\label{defn:graded-polynomial-ring}
	The \textbf{standard grading} on the polynomial ring 
	$R[x]$ is given by polynomial degree.
	\begin{enumerate}
		\item 
		A nonzero element $f \in R[x]$ is \textbf{homogeneous} if $f = ax^t$ for some nonzero $a \in R$ and $t \in \nonnegints$, i.e.\ $f$ is a monomial.
		In this case, we say that $f = ax^t$ is \textbf{homogeneous of degree $t$} and write $\degh(f) = t$.

		\item 
		For all $t \in \nonnegints$ with $t \geq 1$, define $Rx^t$ as the subring of $R[x]$ consisting of all monomials of degree $t$ and the zero polynomial, i.e.\ $Rx^t = \set{rx^t : r \in R}$.
		Define $Rx^0 := R$. 

		\item 
		Under the standard grading, the \textbf{homogeneous component of $R[x]$ of degree $t \in \nonnegints$} is the subring $Rx^t$ of $R[x]$.
		Observe that
		$
			R[x] = \bigoplus_{t \in \nonnegints} Rx^t 
		$ as abelian groups with $\oplus$ taken as an internal direct sum.
	\end{enumerate}
	In this paper: when $R[x]$ is viewed as a graded ring, we always equip $R[x]$ with the standard grading.
\end{definition}
\remark{
	For a general definition for graded rings (i.e.\ not specific to $R[x]$ under the standard grading), see \cite{algebra:peeva,algebra:nastasescu}. 
	Note that the grading on a ring $R$ can be defined relative to any abelian group $G$.
	Here, the degrees of homogeneous elements take values in $G$ and say that $R$ is $G$-graded.

	In our case, we consider the standard grading on $R[x]$ to be an $\nonnegints$-grading on $R[x]$ since we only expect the degrees of homogeneous elements of $R[x]$ to take values in $\nonnegints$.
	We point this out since, in some references, e.g.\ \cite{algebra:webb,algebra:marley,algebra:nastasescu}, 
	$R[x]$ is described to be $\ints$-graded with the homogeneous components of $R[x]$ of degree $t \in \ints$ with $t < 0$ being trivial abelian groups. 
	Here, $R[x]$ is often seen as a subring of the Laurent polynomial ring $R[x, x^{-1}]$.
}
\HIDE{
\remarks{
	\item

	\item 
	We want to emphasize that we have not described/defined the degree of $0 \in R[x]$ as an element of a graded ring, i.e.\ $\degh(0)$ is undefined, following \cite{algebra:nastasescu}. 
	We state this explicitly since definitions for the degree of $0 \in R[x]$ vary depending on the context/reference. 

	For example, viewing $R[x]$ as a ring (not necessarily graded), 
		\cite{algebra:dummit} uses $\deg(0) := 0$ 
		but \cite[p73]{algebra:adkins} uses $\deg(0) := -\infty$.
	Viewing $R[x]$ as a graded ring, \cite{algebra:peeva} defines $\degh(0)$ to be arbitrary since $0$ is an element of $Rx^t$ (the homogeneous component of degree $t$) for all $t \in \nonnegints$. 
	Note that we get around this issue with $\degh(0)$ by only considering nonozero elements for relations involving $\degh(-)$.
} }

With $R[x]$ equipped with the standard grading, observe that we denote the degree of $f \in R[x]$ by $\degh(f)$, as opposed to $\deg(f)$.
Note that, conventionally, $\deg(f)$ denotes the degree of a polynomial $f \in R[x]$, disregarding the equipped grading on $R[x]$. 
More specifically, for all nonzero $f \in R[x]$, we have that:
\begin{equation*}
	\degh(f) = \begin{cases}
		\deg(f) 			&\text{ if $f$ is homogeneous, i.e.\ $f$ is a monomial } \\
		\text{undefined} &\text{ otherwise }
	\end{cases}
\end{equation*}
We identify two advantages to this approach:
\begin{enumerate}
	\item 
	Writing $\degh(f)$ emphasizes that $f \in R[x]$ is to be considered as an element of a graded ring.
	This additional notation suggests that $\degh(f)$ and $\deg(f)$, while related or similar, are two distinct properties.

	We have found this observation to be important when first learning graded module theory especially since, relative to a more general definition of graded ring, $R[x]$ may equipped with a grading other than the standard grading.
	For example, we may equip $R[x]$ with the trivial grading wherein all elements of $R[x]$ are defined to be homogeneous of degree $0 \in \nonnegints$. In this case, we have that $\degh(3x^2) = 0$ despite $\deg(3x^2) = 2$. 
	While we do not consider cases like these in this paper, it is helpful to know this distinction for other applications.

	\item 
	Writing a relation in terms of $\degh(-)$ emphasizes that the relation is generally only valid when the arguments are homogeneous elements of a graded ring.
	\HIDE{For example, consider the relation $\deg(f) = \deg(g)$ for nonzero $f,g \in \rationals[x]$. 
	Given below are some polynomials that satisfy this relation:
	\begin{equation*}
		\begin{array}{r !{=} l} 
			\deg(2x) & \deg(x) \\
			\deg(x^2+2x+1) & \deg(3x - 2x^2) \\
			\deg(3x^3 + 2) & \deg(x^3)
		\end{array}
	\end{equation*}
	Observe that if we replace $\deg(-)$ with $\degh(-)$, i.e.\ we consider the relation $\degh(f) = \degh(g)$ for nonzero $f,g \in \rationals[x]$, only the first equality $\degh(2x) = \degh(x)$ is valid. That is, 
	\begin{equation*}
		\begin{array}{r !{=} c@{\quad}c@{\quad} c !{=} l}
			\degh(2x) & 1
				&=& 1 & \degh(x) \\
			\degh(x^2 + 2x+1) & \text{undefined} 
				&\neq& \text{undefined} & \degh(3x - 2x^2) \\
			\degh(3x^3 + 2) & \text{undefined}
				&\neq& 3 & \degh(x^3)
		\end{array}
	\end{equation*}}
	This will be relevant later in \fref{section:calculation-graded-ifds} where we require and expect nonzero elements to be homogeneous.
\end{enumerate}

\spacer 

The notion of graded rings gives rise to the notion of graded modules.
In this paper, since we only consider the graded structure of modules over $R[x]$ equipped with the standard grading, we use a definition that is specific to this type of modules.
We state this below.

\begin{definition}\label{defn:graded-module}
	An $R[x]$-module $M$ is a \textbf{graded $R[x]$-module} or \textbf{graded module}
	if there exists an $\nonnegints$-indexed family $\set{M_t}_{t \in \nonnegints}$ of additive subgroups $M_t$ of $M$ such that two conditions are satisfied:
	\begin{enumerate}
		\item 
		The underlying abelian group $M$ decomposes into $M = \bigoplus_{t \in \nonnegints} M_t$ with equality and (internal) direct sum taken to be of abelian groups (seen as additive groups). 

		\item For all $s,t \in \nonnegints$, $Rx^s \cdot M_t \subseteq M_{s+t}$, i.e.\ 
		the action of $R[x]$ on $M$ respects the grading on $R[x]$ and on $M$.
	\end{enumerate}
	For all $t \in \nonnegints$, we call the additive subgroup $M_t \subseteq M$ the \textbf{homogeneous component of $M$ of degree $t$}.
	We call the direct sum $\bigoplus_{t \in \nonnegints} M_t$ the \textbf{homogeneous decomposition of $M$}.
	An nonzero element $f \in M$ is \textbf{homogeneous} if there exists $t \in \nonnegints$ such that $f \in M_t$. In this case, we say that $f$ is \textbf{homogeneous of degree $t$} and write $\degh(f) = t$.
	A subset $A \subseteq M$ is \textbf{homogeneous} is each element $a \in A$ is homogeneous (not necessarily of the same degree).
	A $R[x]$-submodule $N$ of $M$ is called a \textbf{graded submodule} if $N$ is also a graded $R[x]$-module.
\end{definition}
\remarks{
	\item 
	The statement $Rx^s \cdot M_t$ is specific to the ring being $R[x]$ equipped with the standard grading. 
	This condition can be written more generally as $R_s M_t \subseteq M_{s+t}$ where $R_s$ is the homogeneous component of degree $s \in \nonnegints$ of the graded ring $R = \bigoplus_{t \in \nonnegints} R_t$.

	\item 
	The remarks about $G$-grading of rings below \fref{defn:graded-polynomial-ring} also apply to graded modules. 
	More specifically, if the graded ring is considered to be $G$-graded, then the graded modules over said ring is also $G$-graded.
	Since we consider $R[x]$ to be $\nonnegints$-graded in this paper, this means that our graded $R[x]$-modules are also $\nonnegints$-graded.

	\HIDE{
	Note that, given a $G$-graded module $M$, \cite{algebra:peeva,algebra:nastasescu,persmod:bubenik-homoloalg} denotes the homogeneous component of $M$ of degree $g \in G$ by $M_g$, i.e.\ attaching the degree $g$ as a subscript of $M$.
	We try not to follow this notation since, as seen later in \fref{chapter:filtrations-and-pershoms}, we construct graded modules on objects that already have a few subscripts in their notation.

	Later in \fref{remark:assumption-on-graded-modules}, we state that we can assume that the homogeneous component of a graded $\field[x]$-module $M$ is given by $M_t x^t$ instead of $M_t$, i.e.\ 
	the abelian group $M$ decomposes into $M = \bigoplus_{t \in \nonnegints} M_t x^t$,
	so that the degree $t$ is explicitly identified by the presence of $x^t$.
	}
}

By the definition above, the homogeneous decomposition of a graded $R[x]$-module $M$ refers to a direct sum of abelian groups. 
In practice, we usually define graded $R[x]$-modules at the level of $R$-modules. That is, we define a graded $R[x]$-module using two collections of objects:
\begin{enumerate}
	\item 
	A family $\set{M_t}_{t \in \nonnegints}$ of $R$-modules 
	such that $M_t \cap M_s = \set{0}$ whenever $t \neq s$, i.e.\ each $R$-module $M_t$ is seen as having nonzero elements distinct from the other $R$-modules in said family.
	Note that we usually force the condition $M_t \cap M_s = \set{0}$ by using the family $\set{M_t x^t}_{t \in \nonnegints}$ constructed from $\set{M_t}$.

	We then define $M$ to be the (internal) direct sum $M = \bigoplus_{t \in \nonnegints} M_t$ of $R$-modules.
	This induces an $R$-module structure on $M$.
	Observe that the nonzero elements in one summand $M_t$ cannot be generated using sums of $R$-multiples of elements from other summands, i.e.\ $M_s$ with $t \neq s$.
	If $R$ is a field $\field$,
		then $M$ is a (graded) $\field$-vector space. 
	As a sidenote, 
	a \textit{graded $\field$-vector space} is an $\field$-vector space (i.e.\ $\field$-module) that satisfies Condition (i) of \fref{defn:graded-module}.
	Note that all graded modules over $\field[x]$ are graded $\field$-vector spaces.
	However, not all graded $\field$-vector spaces form graded modules, e.g.\ the action of $\field[x]$ may not be defined.
	
	\item 
	A collection of assignments of the form $x \cdot m_t := m_{t+1}$,
	defining the product between $x \in R[x]$ of degree $1$ and each element $m_t \in M_t$ of degree $t \in \nonnegints$ to some element $m_{t+1} \in M_{t+1}$ of degree $t+1$.

	These assignments are then extended linearly 
	over the action of $R$ on each $M_t$, the ring operation on $R[x]$, and the addition operation on $M$ as an $R$-module.
	If the action of $R[x]$ is well-defined, then $M$ becomes an $R[x]$-module.
\end{enumerate}
The condition that $M_t \cap M_s = \set{0}$ whenever $t \neq s$ satisfies the condition that $M$ is a direct sum $\bigoplus_{t \in \nonnegints} M_t$, i.e.\ each element $m \in M$ can be represented \textit{uniquely} as a sum of homogeneous elements.
The set of assignments in the form $x \cdot m_t := m_{t+1}$ forces the action of $R[x]$ on $M$ to satisfy $Rx \cdot M_t \subseteq M_{t+1}$. 
If the assignments result in a well-defined action of $R[x]$ on $M$, then the condition that $Rx^s \cdot M_t \subseteq M_{s+t}$ is satisfied.
Therefore, the resulting $R[x]$-module $M$ is also a graded $R[x]$-module.
We provide an example of this construction below.

\begin{example}\label{ex:linear-extension-graded-mod}
	For each $t \in \nonnegints$, 
	define the $\ints$-module $M_t$ by $M_t := \ints\ket{a,b}x^t = \set{k_1 a x^t + k_2 b x^t : k_1, k_2 \in \ints}$
	where $\ints\ket{a,b}$ is the free $\ints$-module with basis $\set{a,b}$.
	Then, the direct sum below results in an $\ints$-module:
	\begin{equation*}
		M = \bigoplus_{t \in \nonnegints} M_t 
		= \bigoplus_{t \in \nonnegints} \ints\ket{a,b}x^t 
		= \Biggl\{\,\,
			\sum_{t=0}^\infty \Bigl( k_t^a ax^t + k_t^b bx^t \Bigr)
			:\,
			\text{only finitely many $k_t^a, k_t^b \in \ints$ are nonzero }
		\Biggr\}
	\end{equation*}
	Define the action of $\ints[x]$ on $M$, which we denote as 
	$\star: \ints[x] \times M \to M$ in this example for clarity,
	by linearly extending the following set of assignments on $M$:
	\begin{equation*}
		x \mathrel{\star} a x^{t} := 2a x^{t+1} 
		\quad\text{ and }\quad 
		x \mathrel{\star} b x^{t} := 3b x^{t+1}
		\qquad\text{ for all } t \in \nonnegints
	\end{equation*}
	Here, linearly extending means that we use the compatibility and biadditivity axioms, as discussed under \fref{defn:action-on-module}, 
	to calculate assignments of $\star: \ints[x] \times M \to M$ not of the form $x \star ax^t$ and $x \star bx^t$
	Below, we provide examples of this linear extension.
	\begin{longtable}{R@{\ }C@{\ }L @{\hspace{12pt}} l}
		\rlap{In the equality $\redmath{=\rule{0pt}{5.5pt}\hspace{-2.5pt}}$, we extend the
		assignments $x \mathrel{\star} ax^t$ and $x \mathrel{\star} bx^t$ over \ldots}\hfill\null
		\\[5pt]
		\hline\\[-3pt]
		x \cdot (ax^{3} + bx^{5}) 
			&\mathrel{\redmath{=\rule{0pt}{5.5pt}\hspace{-2.5pt}}}&
			(x \mathrel{\star} ax^{3}) + (x \mathrel{\star} bx^{5})
			=
			2a x^{4} + 3b x^{6}
			& \ldots\ addition $+$ on $M$ as a $\ints$-module.
		\\[5pt] 
		x \cdot 3a x^{11} 
			&\mathrel{\redmath{=\rule{0pt}{5.5pt}\hspace{-2.5pt}}}&
				3 \cdot (x \mathrel{\star} ax^{11}) 
			= 3 \cdot (2ax^{12}) = 6ax^{12}
			& \ldots\ 
			the action $\cdot$ of $\ints$ on $M_{11} = \rationals\ket{a,b}x^{11}$.
		\\[5pt]
		x^2 \mathrel{\star} ax^3 
			= (x \cdot x) \star ax^3
			&\mathrel{\redmath{=\rule{0pt}{5.5pt}\hspace{-2.5pt}}}&
			x \mathrel{\star} (x \mathrel{\star} ax^3)
			= x \mathrel{\star} (2a x^{4})
			= 4ax^5
			& \ldots\ ring operation $\cdot$ on $\ints[x]$.
		\\[5pt]
		(x + x^{2}) \mathrel{\star} b
			&\mathrel{\redmath{=\rule{0pt}{5.5pt}\hspace{-2.5pt}}}&
			(x \star bx^0) + (x^2 \star bx^{0})
			= 3bx + 9bx^2
			& \ldots\ addition $+$ on $\ints[x]$ as a ring.
	\end{longtable}
	\noindent 
	We claim that this describes a well-defined action of $\ints[x]$ on $M$, which makes $M$ a $\ints[x]$-module.
	We can verify that the decomposition $M = \bigoplus_{t \in \nonnegints} \ints\ket{a,b}x^t$, along with the action $\star$ on $M$, satisfies both conditions in \fref{defn:graded-module}.
	This determines $M$ to be a graded $\ints[x]$-module 
	with homogeneous component of degree $t \in \nonnegints$ given by $\ints\ket{a,b}x^t$.
	We list some examples of elements of $M$ and determine their homogeneity or degree:
	\begin{enumerate}
		\item $m_1 = 2ax + 3bx^2$ is not a homogeneous element of $M$, i.e.\ $\degh(2ax + 3bx^2)$ is undefined. 
		Note that it may feel natural to extend $\deg(-): R[x] \to \nonnegints$ using powers of $x$, e.g.\ 
		we may consider that $\deg(2ax + 3bx^2) = 2$.
		However, 
		the distinction between $\degh(-)$ and $\deg(-)$ in the case of $R[x]$ also holds in the case of $M$.

		\item $m_2 = 4ax^5 - 3bx^5 = (4a-3b)x^5$ is a homogeneous element of $M$ with $\degh(m_2) = 5$.

		\item $m_3 = -7a$ is a homogeneous element with $\degh(m_3) = 0$.
	\end{enumerate}
	Note that, by construction, the homogeneous elements of $M$ are exactly monomials of the form $f x^t$ with nonzero $f \in \ints\ket{a,b}$ and $t \in \nonnegints$.
\end{example}

We want to point out that the elements of an $R[x]$-module $M$ are generally not polynomials in $x$ (unlike the example above).
Consequently, if the $R[x]$-module $M$ is graded, the degrees of the homogeneous elements of $M$ are not generally immediately apparent from the notation for said elements.

It would be very helpful (calculation-wise) if we can, without loss of generality, assume that the notation of the elements of a graded $R[x]$-module also carry information about the degree.
To formalize this, we need a notion of similarity or equivalence between graded $R[x]$-modules.
First, we provide a definition for the family of homomorphisms between graded modules, which in turn induces a definition for isomorphisms between graded modules.

\begin{definition}\label{defn:graded-homomorphism}
	Let $\phi: M \to N$ be an $R[x]$-module homomorphism between graded $R[x]$-modules $M$ and $N$. 
	Let $M_t$ and $N_t$ refer to the homogeneous component of $M$ and $N$ of degree $t \in \nonnegints$ respectively.
	\begin{enumerate}
		\item 
		We call $\phi: M \to N$ a \textbf{graded $R[x]$-module homomorphism} or \textbf{graded homomorphism} 
		if for all $t \in \nonnegints$, 
			$\phi(M_t) \subseteq N_t$.
		That is, $\phi$ sends homogeneous elements of $M$ to either zero or homogeneous elements of $N$ of the same degree.

		\item 
		If $\phi: M \to N$ is both an $R[x]$-module isomorphism and a graded homomorphism, then $\phi: M \to N$ is called a \textbf{graded $R[x]$-module isomorphism} or \textbf{graded isomorphism} and we say that $M$ and $N$ are \textbf{graded isomorphic}, denoted $M \cong N$. 
	\end{enumerate}
\end{definition}

Graded modules and graded homomorphisms give rise to a \textit{category} of graded modules.
We state this as a theorem below.

\begin{theorem}\label{thm:category-of-graded-modules}
	Let $R$ be a PID. Graded $R[x]$-modules and graded $R[x]$-module homomorphisms form a well-defined category, denoted $\catgradedmod{R}$ and called the \textbf{category of graded $R[x]$-modules}.
\end{theorem}
\remark{
	We refer to \cite[Section 2.2]{algebra:nastasescu} for the claim of $\catgradedmod{R}$ being a well-defined category.
	Note that the category of graded modules can be more generally defined using grading over any abelian group $G$ and over some $G$-graded ring $R$.

	For our case, the symbol $\catgradedmod{R}$ refers specifically to the category of $\nonnegints$-graded modules over $R[x]$ equipped with the standard grading.
	Note that while it is possible to define $\catgradedmod{R}$ to be of $\ints$-graded modules (following remarks under \fref{defn:graded-polynomial-ring}), 
	the category equivalence presented in \fref{section:cat-equiv-graded-modules} between $\catpersmod$ and $\catgradedmod{\field}$ is not valid for said definition.
}

Observe that all graded $R[x]$-modules and all graded $R[x]$-module homomorphisms are $R[x]$-modules and $R[x]$-module homomorphisms respectively.
Equivalently, we can also say that $\catgradedmod{R}$ is a \textit{subcategory} of $\catmod{R[x]}$
or that there exists a forgetful functor $\catgradedmod{R} \to \catmod{R[x]}$ (see \cite[Section 2.5]{algebra:nastasescu}).

Later in \fref{section:calculation-graded-ifds},
	we use this subcategory relationship to investigate an algorithm involving $\field[x]$-modules 
	and see if it can be extended to the case of graded $\field[x]$-modules.
The main hurdle here is that the graded structure of the graded modules involved may not be preserved across multiple steps, which may result in an isomorphism that is only true at the level of $R[x]$-modules and not at the level of graded $R[x]$-modules.
The distinction between isomorphisms in $\catmod{R[x]}$ and graded isomorphisms in $\catgradedmod{R[x]}$ becomes relevant when we use the category equivalence discussed in \fref{section:cat-equiv-graded-modules}.

In \fref{chapter:matrix-calculation}, we will talk about different levels of isomorphisms involving graded $R[x]$-modules. 
To avoid ambiguity, we identify some notation below. 

\begin{statement}{Remark}\label{remark:shorthand-for-graded-relations}
	When given an equality $(=)$ or isomorphism relation $(\cong)$ involving $R[x]$-modules, the following symbols denote the category in which the relation applies.
	\begin{enumerate}
		\item 
		The symbol \raisebox{0.5pt}{\scalebox{0.8}{$\catname{Ab}$}}
		in $\upabelian=$ or $\upabelian\cong$
		denotes a relation either in the category $\catabgroup$ of abelian groups, $\catmod{R}$ of $R$-modules, or $\catvectspace$ of $\field$-vector spaces.
		In these cases, only the additive structure and (if defined) the action of $R$ on the modules involved are respected by the relations.

		In this paper, we often use this shorthand to describe the homogeneous decomposition of graded $\field[x]$-modules, wherein an accompanying collection of assignments in the form $x \cdot m x^t$ is expected. 

		\item 
		The symbol \raisebox{0.5pt}{\scalebox{0.8}{$\catname{Mod}$}} in $\upmod=$ or $\upmod\cong$ denotes a relation in the category 
		$\catmod{R[x]}$ of $R[x]$-modules, 
		i.e.\ the relation respects the action of $R[x]$ on the $R[x]$-modules involved (and not just $R$) but it may not respect the graded structure (if such exists).

		\item 
		The symbol \raisebox{0.5pt}{\scalebox{0.8}{$\catname{GrMod}$}} in $\upgraded=$ or $\upgraded\cong$ denotes a relation in the category $\catgradedmod{R}$ of graded $R[x]$-modules, i.e.\ the relation involves graded $R[x]$-modules and 
		is given by a graded isomorphism.
	\end{enumerate}
\end{statement}

Note that the symbols 
	\raisebox{0.5pt}{\scalebox{0.8}{$\catname{Ab}$}}, 
	\raisebox{0.5pt}{\scalebox{0.8}{$\catname{Mod}$}}, 
	and \raisebox{0.5pt}{\scalebox{0.8}{$\catname{GrMod}$}}
can be seen as increasing restrictions on the relation. 
For example, the relation $M \upgraded= N$ implies $M \upmod= N$ which implies $M \upabelian= N$.
The converse is not generally true.

\spacer 

Next, we state the result that allows us to assume without loss of generality that the notation for the elements of a graded $R[x]$-module carry information about their homogeneity or degree.

\begin{proposition}\label{prop:assumption-on-graded-modules}
	Let $M$ be a graded $R[x]$-module and denote its homogeneous component of degree $t \in \nonnegints$ by $M_t \subseteq M$.
	Then, there exists a graded $R[x]$-module $M^\prime$ such that the homogeneous component of $M^\prime$ of degree $t \in \nonnegints$ is given by $M_t x^t$
	and $M^\prime$ is graded isomorphic to $M$.
\end{proposition}
\begin{proof}
	Assume $M$ is a graded $R[x]$-module.
	By definition, there exists a family $\set{M_t}_{t \in \nonnegints}$ of additive subgroups $M_t$ of $M$ such that
	\begin{equation*}
		M \upabelian= \bigoplus_{t \in \nonnegints} M_t 
		\quad\text{ and }\quad 
		Rx^s \cdot M_t \subseteq M_{t+s}
	\end{equation*}
	where $\cdot$ denotes the action of $R[x]$ on $M$.
	The homogeneous component of $M$ of degree $t \in \nonnegints$ is given by $M_t$.
	By restriction of the action of $R[x]$ on $M$, each additive subgroup $M_t$ is also an $R$-module.

	For each $t \in \nonnegints$,
		define the $R$-module $M_t x^t$ by the canonical isomorphism $\eta_t: M_t \to M_t x^t$ by $m_t \mapsto m_t x^t$ for all $m_t \in M_t$.
	Since $M_t x^t \cap M_s x^s = \set{0}$ whenever $t \neq s$,
		$\bigoplus_{t \in \nonnegints} M_t x^t$ is a well-defined (internal) direct sum of $R$-modules.
	Define the $R$-module $M^\prime$ by 
	\begin{equation*}
		M^\prime := \bigoplus\nolimits_{t \in \nonnegints} M_t x^t
	\end{equation*}	
	Then, for all $m \in M^\prime$, 
		$m$ decomposes uniquely into 
		$\sum_{t \in \nonnegints} m_t x^t$ with $m_t \in M_t$ for all $t \in \nonnegints$ 
		such that only finitely many $m_t$'s are nonzero. 

	Let $\eta := \bigoplus_{t \in \nonnegints} \eta_t$ be the $R$-module homomorphism $\eta: M \to M^\prime$
	induced by the direct sum operation between $R$-modules.
	Observe that $\eta$ must be an $R$-module isomorphism.
	Let $\star: R[x] \times M \to M$ 
	be given by 
	$x^s \star m := \eta\bigl( x^s \cdot \eta\inv(m) \bigr)$ for all $m \in M^\prime$ and $s \in \nonnegints$.
	More specifically, if $m \in M'$ decomposes into $m = \sum_{t \in \nonnegints} m_t x^t$ with $m_t \in M_t$ for all $t \in \nonnegints$, then 
	\begin{equation*}
		x^s \star m 
		= x^s \star \sum_{t \in \nonnegints} m_t x^t
		= \sum_{t \in \nonnegints} 
			\eta\bigl( x^s \cdot \eta\inv(m_t x^t) \bigr)
		= \sum_{t \in \nonnegints} ( x^s \cdot m_t ) x^{t+s}
	\end{equation*}
	Since $Rx^s \cdot M_t \subseteq M_{t+s}$, $(x^s \cdot m_t) \in M_{t+s} x^{t+1}$.
	Since $\eta$ is an $R$-module isomorphism,
		$\star$ forms a well-defined action of $R[x]$ on $M^\prime$ and makes $M^\prime$ an $R[x]$-module.
	
	Let $t,s \in \nonnegints$.
	Since $Rx^s \star M_t x^t = (Rx^s \cdot M_t) x^{t+s} \subseteq (M_{t+s})x^{t+1} = M_{t+s}x^{t+s}$,
		$M^\prime$ is a graded $R[x]$-module 
		with homogeneous component of degree $t \in \nonnegints$ given by $M_t x^t$.
	Since $\eta(M_t) = M_t x^t$, $\eta$ is a graded $R[x]$-module homomorphism.
	Since $\eta$ is an $R[x]$-module isomorphism,
		$\eta$ is a graded isomorphism.
	Therefore, $M$ and $M^\prime$ are graded isomorphic.
\end{proof}

Observe that the fact that $M \upabelian= \bigoplus_{t \in \nonnegints} M_t$ is a direct sum of $R$-modules is critical here since it implies that each nonzero homogeneous element of $M$ has a unique degree.
Then, using the proposition above, we can assume the following:
 
\begin{bigremark}\label{remark:assumption-on-graded-modules}
	Let $M$ be a graded $R[x]$-module.
	In this paper, assume without loss of generality
	that there exists a family $\set{M_t}_{t \in \nonnegints}$ of $R$-modules such that 
	\begin{equation*}
		M \upabelian= \bigoplus_{t \in \nonnegints} M_t x^t 
		\qquad\text{ and }\qquad 
		Rx^s \cdot M_t \subseteq M_{s+t} \text{ for all } t,s \in \nonnegints 
	\end{equation*}
	and the homogeneous component of $M$ of degree $t \in \nonnegints$ is exactly $M_t x^t$.
	Therefore, a nonzero element $m \in M$ is homogeneous of degree $t \in \nonnegints$ if and only if $m = m_t x^t$ for some $m_t \in M_t$.
\end{bigremark}

Note that some $R[x]$-module isomorphisms, if between two graded $R[x]$-modules, can be made into graded $R[x]$-module isomorphisms by applying appropriate shifting in grading. 
Below, we provide a definition for this grading shift.

\begin{definition}\label{defn:upwards-shift-graded-module}
	Let $M$ be a graded $R[x]$-module and assume, following \fref{remark:assumption-on-graded-modules}, that the homogeneous component of $M$ of degree $t$ is given by $M_t x^t$ for some $R$-module $M_t$ for all $t \in \nonnegints$.
	For each $k \in \nonnegints$, 
	define the $k$-\textbf{upwards shift} $\Sigma^k M$ of 
	a graded $R[x]$-module $M$ 
	to be the graded $R[x]$-module given as follows.
	\begin{enumerate}
		\item 
		The underlying abelian group of $\Sigma^k M$ is given by $
			\Sigma^k M 
			:= \bigoplus_{t \in \nonnegints} M_t x^{t+k}
			= \bigoplus_{t = k}^\infty M_{t-k} x^t
		$.

		\item 
		The action of $R[x]$ on $\Sigma^k M$ satisfies 
		$x \cdot ax^{t+k} = bx^{t+k+1}$ whenever 
		$x \cdot_M ax^t = bx^{t+1}$ with $ax^t \in M_t x^t \subseteq M$, $bx^{t+1} \in M_{t+1}x^{t+1}$ 
		where $\cdot_M$ refers to the action of $R[x]$ on $M$.
	\end{enumerate}
	The homogeneous component of $\Sigma^k M$ of degree $t \in \nonnegints$ is given by $M_{t-k} x^t$ if $t \geq k$ and trivial if $t < k$.
	The index $k \in \nonnegints$ is called the \textbf{grading shift} or \textbf{shift} of $\Sigma^k M$.
\end{definition}
\remark{
	In \cite{algebra:peeva, algebra:nastasescu}, the $k$-upwards shift $\Sigma^k M$ of $M$ is denoted as $M(-k)$ since the homogeneous component of $\Sigma^k M = M(-k)$ of degree $t \geq k$ is exactly the homogeneous component of $M$ of degree $t-k$.
	The notation $\Sigma^k M$ follows that of \cite{matrixalg:zomorodian}. 
}

One key application of upwards shifts involves the isomorphism between $R[x]$ and $\Sigma^k R[x]$.
Let $\phi: R[x] \to \Sigma^k R[x]$ be given by $r x^t \mapsto r x^{t+k}$.
We claim that $\phi$ is an $R[x]$-module homomorphism with inverse $rx^{t} \mapsto rx^{t-k}$. Note that $\Sigma^k R[x]$ has no elements of the form $r x^{t}$ with $t < k$.
Observe that if $k \geq 1$, $\phi$ sends the homogeneous elements of $R[x]$ of degree $t \in \nonnegints$ to the homogeneous component of $\Sigma^k R[x]$ of degree $t+k \neq t$.
Therefore, $\phi$ is not a graded isomorphism.
In other words, we have that 
\begin{equation*}
	\text{If $k \geq 1$: }\qquad
	R[x] \upmod\cong \Sigma^k R[x] 
	\qquad\text{ but }\qquad 
	R[x] \upgraded{\not\cong} \Sigma^k R[x] 
\end{equation*}
Note that this affects how the Structure Theorem on finitely-generated $\field[x]$-modules in in the category $\catmod{\field[x]}$
translates to the category $\catgradedmod{\field}$ of graded $\field[x]$-modules.

\spacer 

Next, we identify a number of useful results involving algebraic constructions in $\catgradedmod{R}$, particularly those that extend, in a specific sense, those of the ungraded category $\catmod{R[x]}$ of $R[x]$-modules.

\begin{proposition}\label{prop:graded-submodule-hom-generators}
	An $R[x]$-submodule $N$ of a graded $R[x]$-module $M$ is a graded $R[x]$-module if and only if 
		there exists a homogeneous system of generators for $N$,
		i.e.\ there exists a homogeneous set of elements in $M$ that generates $N$.
\end{proposition}
\remark{
	The first steps of the proof can be found under \cite[Proposition 2.1]{algebra:marley}.
	This proposition is also listed as \cite[Exercise 2.8]{algebra:peeva}.
}
\negativespacer
\begin{proposition}
	\label{prop:kernels-and-images-are-graded}
	Kernels and images of graded homomorphisms are graded, i.e.\ 

	If $\phi: M \to N$ is a graded $R[x]$-module homomorphism 
	between graded $R[x]$-modules $M$ and $N$, then 
		$\ker(M)$ is graded submodule of $M$ and $\im(M)$ is a graded submodule of $N$.
\end{proposition}
\remark{
	A proof is available under \cite[Proposition 2.9, p8-9]{algebra:peeva}.
}
\negativespacer
\begin{proposition}\label{prop:cokernels-of-graded-modules}
	Quotients of graded modules over graded submodules are graded, i.e.\ 

	Let $M$ be a graded $R[x]$-module and $N$ a graded submodule of $M$. Following \fref{remark:assumption-on-graded-modules}, let the homogeneous component of $M$ and $N$ of degree $t \in \nonnegints$ be given by $M_t x^t$ and $N_t x^t$ with $R$-modules $M_t$ and $N_t$ respectively.
	Then, $M \bigmod N$ is a graded $R[x]$-module with homogeneous component of degree $t \in \nonnegints$ given by $(M_t \bigmod N_t) x^t$.
\end{proposition}
\remark{
	Note that $M \bigmod N = \coker(N \hookrightarrow M)$ and the inclusion map $N \hookrightarrow M$ is a graded homomorphism.
	\cite[Proposition 2.9, p8-9]{algebra:peeva} states that cokernels of graded homomorphisms are also graded along with a proof.
}
\negativespacer
\begin{proposition}\label{prop:direct-sums-are-graded}
	Direct sums of graded modules are graded, i.e.\ 
	
	Let $M$ and $N$ be graded $R[x]$-modules. Following \fref{remark:assumption-on-graded-modules}, let the homogeneous component of $M$ and $N$ of degree $t \in \nonnegints$ be given by $M_t x^t$ and $N_t x^t$ with $R$-modules $M_t$ and $N_t$ respectively.
	Then, $M \oplus N$ is a graded $R[x]$-module with homogeneous component of degree $t \in \nonnegints$ given by $(M_t \oplus N_t)x^t$.
\end{proposition}
\remark{
	This is discussed in \cite[p20]{algebra:nastasescu}
	and listed as \cite[Exercise 1.4]{algebra:marley}.
	\cite{algebra:peeva} does not seem to explicitly state this but assumes this in other results such as \cite[Exercise 2.9, Theorem 2.10]{algebra:peeva}.
}

\spacer 

Next, we state a result involving upward shifts in grading and quotients.

\begin{proposition}\label{prop:shift-distribute-over-quotient}
	Let $M$ be a graded $R[x]$-module and let $N$ be a graded submodule of $M$.
	For all $s \in \nonnegints$,
	\begin{equation*}
		\Sigma^s \Big( M \bigmod N \Big) \upgraded\cong 
		\Big( \Sigma^s M \Big) 
			\bigmod \Big( \Sigma^s N \Big)
	\end{equation*}
	That is, the upwards shifts distribute over quotients.
\end{proposition}
\begin{proof}
	Since $N$ is a graded submodule of $M$, then $\Sigma^s N$ must also be a graded submodule of $\Sigma^s M$ and the quotient $\Sigma^s M \bigmod \Sigma^s N$ is well-defined.
	Denote the homogeneous component of degree $t \in \nonnegints$ of a graded module by addition of the subscript $t$.
	For all $t < s$, $(\Sigma^s M)_t = 0$, $(\Sigma^s N)_t = 0$, and $(\Sigma^s (M \bigmod N))_t = 0$.
	For all $t \geq s$, we have that 
	\begin{equation*}
		\Bigl(
			\Sigma^s \bigl( M \bigmod N \bigr)
		\Bigr)_t
		=
		(M \bigmod N)_{t-s} 
		= M_{t-s} \bigmod N_{t-s}
		= (\Sigma^s M)_t \bigmod (\Sigma^s N)_t
		= \Bigl(
			(\Sigma^s M) \bigmod (\Sigma^s N)
		\Bigr)_t
	\end{equation*}
	by \fref{defn:upwards-shift-graded-module} (for $\Sigma^k M$) and by \fref{prop:shift-distribute-over-quotient} (involving quotients).
\end{proof}

\spacer 

Recall that $\field[x]$ is a PID for any field $\field$.
The Structure Theorem (\fref{thm:structure-theorem}) on $\catmod{\field[x]}$ guarantees the existence and uniqueness 
of \textit{invariant factor decompositions} for finitely generated $\field[x]$-modules up to $\field[x]$-module isomorphism.
In \fref{section:calculation-graded-ifds},
	we state a corresponding result in $\catgradedmod{\field}$ which we call the \textit{The Graded Structure Theorem} (\fref{thm:graded-structure-theorem})
	for \textit{graded invariant factor decompositions}
	(given in \fref{defn:graded-invariant-factor-decomposition}) for finitely generated graded modules.
For convenience, statements of 
\fref{thm:graded-structure-theorem} 
and \fref{defn:graded-invariant-factor-decomposition}
are included below.

\begingroup
\vspace{\statementspacing}\setlength{\parindent}{0in}\noindent 
\textbf{\fref{thm:graded-structure-theorem}. The Graded Structure Theorem.}
 
	Let $M$ be a finitely generated graded $\field[x]$-module over $\field[x]$ for some field $\field$.
	Then, 
	there exists a finite direct sum of shifted cyclic graded submodules of $\field[x]$ that is graded isomorphic to $M$ as follows:
	\begin{equation*}
		M \,\upgraded\cong\,
			\mathlarger\Sigma^{s_1}\!\paren{\frac{\field[x]}{(x^{t_1})}}
			\oplus \cdots \oplus 
			\mathlarger\Sigma^{s_r}\!\paren{\frac{\field[x]}{(x^{t_r})}}
			\oplus 
			\Sigma^{s_{r+1}}\field[x]
			\oplus \cdots \oplus 
			\Sigma^{s_{m}}\field[x]
	\end{equation*}
	with indices $s_1, \ldots, s_r, \ldots, s_m \in \nonnegints$ and 
	non-zero, non-unit $x^{t_1}, x^{t_2}, \ldots, x^{t_r} \in \field[x]$
	such that the divisibility relation 
	$x^{t_1} \divides x^{t_2} \divides \cdots \divides x^{t_r}$ is satisfied.
	Furthermore, the collection 
	$
		\bigl\{
			(x^{t_1}, s_1), \ldots, (x^{t_r}, s_r),
			(0, s_{r+1}), \ldots, (0, s_m)
		\bigr\}
	$
	is uniquely determined by $M$ up to graded isomorphism.
\vspace{\statementspacing}
\endgroup
\remarks{
	\item 
	We use the term \textit{cyclic} similarly in the case of $R[x]$-modules, i.e.\ 
	a cyclic graded module is a graded module that can be generated by a single element. 
	It can be verified that cyclic graded modules can be generated by a single \textit{homogeneous} element as well.

	\item 
	Note that the ideal $(x^t) = \field[x] \cdot x^t$ is a graded submodule of $\field[x]$ since it is generated by a homogeneous element $x^t$.
	Therefore, the quotient $\field[x] \bigmod (x^t)$ results in a graded module.
	In this case, we call $(x^t)$ a \textit{graded ideal} of $\field[x]$.
}


\begingroup\setlength{\parindent}{0in}\noindent 
\textbf{\fref{defn:graded-invariant-factor-decomposition}.}
	Let $M$ be a finitely generated $\field[x]$-module and let the following direct sum decomposition of $M$ be as denoted in the Graded Structure Theorem 
	(\fref{thm:graded-structure-theorem}):
	\begin{equation*}
		M \,\upgraded\cong\,
			\mathlarger\Sigma^{s_1}\!\paren{\frac{\field[x]}{(x^{t_1})}}
			\oplus \cdots \oplus 
			\mathlarger\Sigma^{s_r}\!\paren{\frac{\field[x]}{(x^{t_r})}}
			\oplus 
			\Sigma^{s_{r+1}}\field[x]
			\oplus \cdots \oplus 
			\Sigma^{s_{m}}\field[x]
	\end{equation*}
	This decomposition is called the \textbf{graded invariant factor decomposition} of $M$.
	The \textbf{invariant factors} of $M$ are given by 
		$x^{t_1}, \ldots, x^{t_r} \in \field[x]$ (which are non-zero and non-unit)
	and the \textbf{grading shifts} of $M$ by $s_1, \ldots, s_m \in \nonnegints$.
\vspace{\statementspacing}\endgroup

The notion of grading on $R[x]$-modules and $R[x]$-module homomorphisms also extend to the case of chain complexes. We provide relevant definitions below.

\begin{definition}
	A chain complex $C_\ast = (C_n, \boundary_n)_{n \in \ints}$ of $R[x]$-modules is a \textbf{graded chain complex} if for all $n \in \ints$, $C_n$ is a graded $R[x]$-module and $\boundary_n: C_n \to C_{n-1}$ is a graded homomorphism.
	A chain map $\phi_\ast: C_\ast \to A_\ast$ with $\phi_\ast = (\phi_n: C_n \to A_n)_{n \in \ints}$ between graded chain complexes $C_\ast$ and $A_\ast$ is a \textbf{graded chain map} if each homomorphism $\phi_n$ is a graded homomorphism.
\end{definition}

The collection of graded chain complexes and graded chain maps also form a chain complex category and a corresponding chain homology functor. We identify some notation below.

\begin{definition}\label{defn:cat-graded-chain-complexes}
	Let $\catchaincomplex{\catgradedmod{R}}$ denote the \textbf{category of graded chain complexes} of graded $R[x]$-modules and graded chain maps.
	To each $n \in \ints$, let $H_n(-)$ refer to the $n$\th \textbf{chain homology functor} $H_n(-): \catchaincomplex{\catgradedmod{R}} \to \catgradedmod{R}$.
\end{definition}
\remark{
	For the claim of $\catchaincomplex{\catgradedmod{R}}$ being well-defined, we refer to \cite[Section 2.2]{algebra:nastasescu} which claims that $\catgradedmod{R}$ is an \textit{abelian category}, and to \cite[Section 5.5]{cattheory:rotman} which states that each abelian category has a corresponding category of chain complexes and family of homology functors.
	Note that \fref{prop:cokernels-of-graded-modules} states that the $n$\th chain homology of a graded chain complex, as the cokernel of a graded homomorphism, must be a graded module.
}

Note that, much like how $\catgradedmod{R}$ can be seen as a subcategory of $\catmod{R[x]}$, 
	the category $\catchaincomplex{\catgradedmod{R}}$ of graded chain complexes can be seen as a subcategory of 
	the category $\catchaincomplex{\catmod{R[x]}}$ of chain complexes of $R[x]$-modules.

\spacer 

Lastly, the additional structure required for the grading for graded $R[x]$-modules also brings about graded versions of certain properties of ungraded $R[x]$-modules.
For example, \cite[p21]{algebra:nastasescu} defines a graded $\field[x]$-module to be \textit{graded-free} or \textit{gr-free} if it has a homogeneous basis, and provides 
a counterexample in which a graded module over some graded ring (other than $R[x]$) with a non-homogeneous basis is not graded-free.
For clarity, we explicitly make the following remark.

\begin{bigremark}
	When we state that a graded $R[x]$-module $M$ is a \textbf{free graded} module,
		we mean that $M$ is both a graded module and a free $R[x]$-module, i.e.\ $M$ has a basis that is not necessarily homogeneous.
	In this paper, if $M$ has a homogeneous basis, we do not use the term \textit{graded-free} for $M$ to avoid ambiguity.
\end{bigremark}

 \clearpage

\section{The Equivalence between Persistence Modules and Graded Modules}\label{section:cat-equiv-graded-modules}

In this section, we discuss the isomorphism of categories between the category $\catpersmod$ of persistence modules over $\field$ (see \fref{defn:persmod-cat}) and 
the category $\catgradedmod{\field}$ of $\nonnegints$-graded $\field[x]$-modules (see \fref{thm:category-of-graded-modules}).
We also discuss how this isomorphism allows us to correspond the algebraic constructions in $\catpersmod$ (e.g.\ direct sums, persistence isomorphisms, chain complexes) to those in $\catgradedmod{\field}$.

Before we begin, we want to emphasize that the results presented in this section rely heavily on concepts involving category theory and homological algebra, and proofs for said results are not provided in the following discussion.
Since we mainly need to \textit{use} the results in this section for the calculation of persistent homology, 
	it should suffice to know the definitions involving categories and functors presented in \fref{appendix:cat-theory} 
	and the following definition for an   
	\textit{isomorphism of categories}:
\begin{blockquote}
	\noindent 
	Given two categories $\catname{C}$ and $\catname{D}$ with identity functors $\id_{\catname{C}}: \catname{C} \to \catname{C}$ and $\id_{\catname{D}}: \catname{D} \to \catname{D}$ respectively,
	an \textbf{isomorphism of categories} or \textbf{category isomorphism} between $\catname{C}$ and $\catname{D}$ refers to a pair of functors $F: \catname{C} \to \catname{D}$ and $G: \catname{D} \to \catname{C}$ such that 
	$G \circ F = \id_{\catname{C}}$ and 
	$F \circ G = \id_{\catname{D}}$ \cite[p20]{cattheory:rhiel}.
\end{blockquote}
\noindent 
Roughly speaking, the functors $F$ and $G$ describe a bijective correspondence between the objects and morphisms of $\catname{C}$ and those of $\catname{D}$ respectively.

\spacer 
\remark{
	An \textit{equivalence of categories} or a \textit{category equivalence} between categories $\catname{C}$ and $\catname{D}$ refers to a pair of functors $F: \catname{C} \to \catname{D}$ and $G: \catname{D} \to \catname{C}$ such that there exists a \textit{natural isomorphism} between the functors $G \circ F$ and $\id_{\catname{C}}$ and another between 
	$F \circ G$ and $\id_{\catname{D}}$.
	The existence of these natural isomorphisms is usually denoted by 
	$G \circ F \cong \id_{\catname{C}}$ and $F \circ G \cong \id_{\catname{D}}$.

	We point this out since \cite[Theorem 3.1]{matrixalg:zomorodian} claims the existence of an \textit{equivalence of categories} between $\catpersmod$ and $\catgradedmod{\field}$.
	For our discussion, it should suffice to know that the existence of an isomorphism of categories is a stronger requirement than that of an equivalence of categories, 
		i.e.\ a category isomorphism implies a category equivalence but the converse is not generally true,
	and that the results we present below also apply for the weaker case of category equivalence.
}

In this section, the isomorphism of categories between $\catpersmod$ and $\catgradedmod{\field}$ is given by 
the functors 
	$\togrmod: \catpersmod \to \catgradedmod{\field}$ (described in \fref{defn:togrmod})
	and 
	$\topersmod: \catgradedmod{\field} \to \catpersmod$ (described in \fref{defn:topersmod}).
Then, 
	we present the claim that $\togrmod$ and $\topersmod$ indeed determine an isomorphism of categories in \fref{thm:cat-isom-persmod-grmod}.
Finally, the propositions following \fref{thm:cat-isom-persmod-grmod} describe the correspondence between algebraic constructions.

We begin with a description of $\togrmod: \catpersmod \to \catgradedmod{\field}$ below, adapted from \cite[p8]{persmod:bubenik-homoloalg}.

\begin{definition}\label{defn:togrmod}
	Fix a field $\field$.
	Define $\togrmod: \catpersmod \to \catgradedmod{\field}$ to be the following assignment of the objects and morphisms of $\catpersmod$ to those of $\catgradedmod{\field}$ respectively.
	\begin{enumerate}
		\item 
		\textit{Object Assignment:}
		Let $(V_\bullet, \alpha_\bullet)$ be a persistence module over $\field$.
		Define $\togrmod(V_\bullet, \alpha_\bullet) = \togrmod(V_\bullet)$ 
		to be the graded $\field[x]$-module given as follows:
		\begin{equation*}
			\togrmod(V_\bullet) 
			\upabelian= 
			\bigoplus_{t \in \nonnegints} V_t x^t 
			= \set{\,\,
				\sum_{t=0}^\infty v_t x^t \,:\, v_t \in V_t \text{ for all } t \in \nonnegints 
			\,}
		\end{equation*}
		where only finitely many of the $v_t$'s are nonzero for each element of $\togrmod(V_\bullet)$.
		Let the action of $\field[x]$ on $\togrmod(V_\bullet)$ be given by $x^{s} \cdot v_t x^t := \alpha_{t+s,t}(v_t) x^{s+t}$.
		Then, for each $t \in \nonnegints$, 
			the homogeneous component of $\togrmod(V_\bullet, \alpha_\bullet)$ of degree $t$ is given exactly by $V_t x^t$.

		\item 
		\textit{Morphism Assignment:}
		Let $\phi_\bullet: (V_\bullet, \alpha_\bullet) \to (W_\bullet, \gamma_\bullet)$ 
		be a persistence morphism between persistence modules over $\field$
		with $\phi_\bullet = (\phi_t: V_t \to W_t)_{t \in \nonnegints}$.
		Define $\togrmod(\phi_\bullet)$ to be the graded $\field[x]$-homomorphism 
		$\togrmod(\phi_\bullet): \togrmod(V_\bullet) \to \togrmod(W_\bullet)$ given by 
		\begin{equation*}
			\sum_{t=0}^\infty v_t x^t 
			\,\longmapsto\,
			\sum_{t=0}^\infty \phi_t(v_t) x^t 
		\end{equation*}
	\end{enumerate}
\end{definition}
\remarks{
	\item 
	The notation $\togrmod$ is not standard or convention in persistence literature.
	We introduced this notation since we had found that explicitly identifying the ``conversion'' of persistence modules to graded modules to be helpful in understanding the theory.
	
	On a less serious note, 
		we chose uppercase gamma $(\mathbf\Gamma)$ since 
		it looks like the letter \textbf{T}.
	In this case, the notation $\togrmod$ references the conversion of a persistence module \textbf{to} a graded module in $\catgradedmod{\field}$.
	For contrast, \cite{matrixalg:zomorodian} denotes the functor corresponding to $\togrmod$ simply as $\alpha$ (lowercase alpha), i.e.\ they would write $\alpha(V_\bullet)$ instead of $\togrmod(V_\bullet)$.

	\item 
	The statement of $\togrmod$ above follows \fref{remark:assumption-on-graded-modules}, where we explicitly include the power $x^t$ on the notation for the elements of graded modules to identify the degree $t \in \nonnegints$.
	For contrast, 
		\cite{matrixalg:zomorodian} and \cite{persmod:bubenik-homoloalg} defines $\togrmod(V_\bullet)$ as 
	$\togrmod(V_\bullet) := \bigoplus_{t \in \nonnegints} V_t$, i.e.\ the summands are $V_t$ instead of $V_t x^t$.
	In this case, the direct sum of $\field$-vector spaces is to be interpreted as an \textit{external direct sum} and the elements of $\togrmod(V_\bullet)$ are $\nonnegints$-indexed tuples $(v_0, v_1, \ldots)$ with the entry at index $t \in \nonnegints$ being an element of $V_t$.
}

The assignment $\togrmod$ as denoted above determines a functor, as stated below.

\begin{proposition}\label{prop:togrmod-is-a-functor}
	Fix a field $\field$. 
	The object and morphism assignment $\togrmod$ by \fref{defn:togrmod} determines a functor 
	$\togrmod: \catpersmod \to \catgradedmod{\field}$. 
\end{proposition}
\remark{
	A brief discussion of the proof is available in \cite[p8]{persmod:bubenik-homoloalg},
		which relies on the commutativity requirement on the structure maps of persistence modules (see \fref{lemma:persmod-functor-props})
		and on the linear maps of persistence morphisms (see \fref{defn:persmod-cat}(ii)).
}

One way to prove the proposition above involves verifying that the assignments are well-defined (in that, they indeed produce graded modules and graded homomorphisms)
and that the functorial axioms are satisfied (see \fref{defn:functor}(i-iv)).
In the following discussion, we roughly explain why the object and morphism assignment by $\togrmod$ produces graded $\field[x]$-modules and graded $\field[x]$-module homomorphisms respectively.

\spacer 

\noindent 
\textbf{About the Object Assignment of $\togrmod(-)$:}
The construction of the graded $\field[x]$-module by $\togrmod(-)$, as stated in \fref{defn:togrmod}(i), can be understood to happen in three stages:
Let $V_\bullet = (V_\bullet, \alpha_\bullet)$ be a persistence module over some field $\field$.
\begin{enumerate}
	\item 
	The statement $\togrmod(V_\bullet) \upabelian= \bigoplus_{t \in \nonnegints} V_t x^t$ first defines $\togrmod(V_\bullet)$ as a (graded) $\field$-vector space.

	For each $t \in \nonnegints$,
		$V_t$ is an $\field$-vector space by assumption of $(V_\bullet,\alpha_\bullet)$ being a persistence module 
		and $V_t x^t$ is a $\field$-vector space given by the canonical linear isomorphism $V_t \to V_t x^t$ by $v \mapsto vx^t$ for all $v \in V_t$.
	Then, 
		the direct sum $\bigoplus_{t=0}^\infty V_t x^t$ of $\field$-vector spaces induces an action of $\field$ on the set $\bigoplus_{t=0}^\infty V_t x^t$ (as discussed in \fref{remark:assumption-on-graded-modules}),
	making $\bigoplus_{t=0}^\infty V_t x^t$ an $\field$-module, i.e.\ a $\field$-vector space.

	Note that the characterization of the elements of $\togrmod(V_\bullet)$ as $\sum_{t=0}^\infty v_t x^t$ refers to how we prefer to interpret the elements of graded modules as formal sums of powers of $x^t$ with coefficients in $V_t$, per \fref{remark:assumption-on-graded-modules}.

	\item 
	The assignments $x^s \cdot v_t x^t = \alpha_{t+s,t}(v_t) x^{t+s}$ determine the action of $\field[x]$ on 
	$\togrmod(V_\bullet)$ and makes $\togrmod(V_\bullet)$
	a $\field[x]$-module, not necessarily graded at this point.

	For clarity, 
	let $\star: \field[x] \times M \to M$ with $M := \togrmod(V_\bullet)$ be given by $x^s \star v_t x^t := \alpha_{t+s,t}(v_t) x^{t+s}$,
		i.e.\ we denote the supposed action of $\field[x]$ on $M$ by $\star$.
	Since $\alpha_{t+s,t}(V_t) \subseteq V_{t+s}$ by definition of $V_\bullet$ as a persistence module,
		$\alpha_{t+s,t}(v_t) x^{t+s} \in M$.

	Recall that, by \fref{defn:action-on-module}, an action of $\field[x]$ on $\togrmod(V_\bullet)$ must be a biadditive group action.
	Since the restriction of $\star$ onto $\field \times M \to M$ describes the action of $\field$ on $M$ and $M$ has been shown to be a $\field$-vector space,
		it should suffice to show that for all $t,s,r \in \nonnegints$ and for all $v_t \in V_t$,
	\begin{equation*}
		x^r \star (x^s \star v_t x^t) = (x^r \cdot x^s) \star v_t x^t
	\end{equation*}
	where $\cdot$ refers to the multiplication on $\field[x]$.
	By \fref{lemma:persmod-functor-props},
		we have the commutativity relation 
		$\alpha_{t+s+r,t} = \alpha_{t+s+r,t+s} \circ \alpha_{t+s,t}$ (the structure maps of $V_\bullet$ commute with each other).
	Then, 
	\begin{align*}
		x^r \star (x^s \star v_t x^t)
		&= x^r \star \alpha_{t+s,t}(v_t) \, x^{t+s} 
		= \bigl(
			\alpha_{t+s+r, t+s} \circ \alpha_{t+s, t}
		\bigr)(v_t) \, x^{t+s+r}
		\\
		&= \alpha_{t+s+r, t}(v_t) \, x^{t+s+r}
		= x^{s+r} \star v_t x^t
		= (x^r \cdot x^s) \star v_t x^t
	\end{align*}
	Therefore, $\star$ is a well-defined action of $\field[x]$ on $M$ and $M = \togrmod(V_\bullet)$ is a $\field[x]$-module.

	\item 
	Finally, the specified domain and codomain of the structure maps $\alpha_{s,t}: V_t \to V_s$ of $V_\bullet$ for each $t \in \nonnegints$ makes $M = \togrmod(V_\bullet)$ a \textit{graded} $\field[x]$-module.
	More specifically, 
		the relation $\alpha_{t+s,t}(V_t) \subseteq V_{t+s}$ 
		implies that the action $\star$ of $\field[x]$ on $M$ satisfies the following relation for all $t,s \in \nonnegints$:
	\begin{equation*}
		\field x^{s} \star V_t x^t 
		= \alpha_{t+s,t}(V_t) \, x^{t+s}
		\subseteq V_{t+s} x^{t+s}
	\end{equation*}
	This, together with the decomposition $M = \bigoplus_{t \in \nonnegints} V_t x^t$,
		satisfies the conditions of \fref{defn:graded-module} and determines that $\togrmod(V_\bullet)$ is a \textit{graded} $\field[x]$-module with homogeneous component of degree $t \in \nonnegints$ given exactly by $V_t x^t$.
\end{enumerate}
Recall that, by \fref{prop:persmod-collection}, the structure maps $\alpha_{s,t}: V_t \to V_s$ of $V_\bullet$ are determined uniquely by the collection $\set{\alpha_t: V_t \to V_{t+1}}_{t \in \nonnegints}$.
Similarly,
	we can unambiguously describe the action $\star$ of $\field[x]$ on $\togrmod(V_\bullet)$ using assignments of the following form:
\begin{equation*}
	x \star v_t x^t = \alpha_t(v_t) x^{t+1}
	\quad\text{ for all }\quad 
	v_t x^t \in V_t x^t \subseteq \togrmod(V_\bullet) 
	\quad\text{ with }\quad
	\alpha_t: V_t \to V_{t+1}
\end{equation*}
As a sidenote, \cite{matrixalg:zomorodian} uses this simpler characterization to define the object assignment of the functor  $\catpersmod \to \catgradedmod{\field}$.
We give an example of the object assignment of $\togrmod$ in action below.

\begin{example}
	Let $(V_\bullet, \alpha_\bullet)$ be a persistence module over $\rationals$ with vector spaces given by $V_t = \rationals\ket{a}$ 
	and linear maps $\alpha_t: V_t \to V_{t+1}$ given by $a \mapsto a$ for all $t \in \nonnegints$.
	Then, the graded $\field[x]$-module $\togrmod(V_\bullet)$ is given by 
	\begin{equation*}
		\togrmod(V_\bullet) \upabelian= 
		\bigoplus_{t \in \nonnegints} V_t x^t 
		= \bigoplus_{t \in \nonnegints} \rationals\ket{a} x^t 
		= \rationals[x]\ket{a}
	\end{equation*}
	with the action of $\rationals[x]$ on $\togrmod(V_\bullet)$ given by $x \cdot ax^t = \alpha_t(a)x^{t+1} = ax^{t+1}$.
	Note that $\rationals[x]\ket{a}$ is the free $\field[x]$-module generated by the indeterminate $a$ with $a$ interpreted to have degree $0$.
	Since $x \cdot ax^t = ax^{t+1}$,
		it can verified that the map $\togrmod(V_\bullet) \to \rationals[x]$ by $ax^t \mapsto x^t \in \field[x]$
		is a graded isomorphism.
	That is, $\togrmod(V_\bullet)$ is graded isomorphic to $\rationals[x]$.
\end{example}

The graded isomorphism in the example above is relatively straightforward to determine since 
	the chosen bases for the vector spaces $\set{V_t}_{t \in \nonnegints}$ are ``compatible'' with the structure maps and the ring operation in $\rationals[x]$.
Here, we interpret a ``compatible'' choice to mean the following:
	given a persistence module $(V_\bullet, \alpha_\bullet)$, the structure map $\alpha_t : V_t \to V_{t+1}$ sends a basis element $a_t$ of $V_t$ to another basis element $b_{t+1}$ of $V_{t+1}$ so that the action of $\field[x]$ on $\togrmod(V_\bullet)$ is given by $x \cdot a_t x^t \mapsto b_{t+1} x^{t+1}$.
Loosely speaking, the bases are chosen such that multiplication by $x$ on an element $vx^t$ of $\togrmod(V_\bullet)$ changes the basis element and adds $1$ to the power of $x^t$.
We provide an example below where the immediate choice for bases is not ``compatible'' in this same sense.

\begin{example}
	Let $(W_\bullet, \gamma_\bullet)$ be a persistence module over $\rationals$ with vector spaces given by $W_t = \rationals\ket{a}$ and linear maps $\gamma_t: W_t \to W_{t+1}$ by $a \mapsto 2a$ for all $t \in \nonnegints$.
	Define $M$ to be the graded $\rationals[x]$-module given by $M := \togrmod(W_\bullet)$.
	Then, $M$ has the following underlying abelian group:
	\begin{equation*}
		\togrmod(W_\bullet) 
		\upabelian= 
		\bigoplus_{t \in \nonnegints} W_t x^t 
		= \bigoplus_{t \in \nonnegints} \rationals\ket{a} x^t 
		= \rationals[x]\ket{a}
	\end{equation*}
	For clarity, let $\star: \rationals[x] \times M \to M$ denote the action of $\rationals[x]$ on $M$
	and let $\cdot$ refer to the usual scalar multiplication on $\rationals[x]$ given by $x \cdot kx^t = kx^{t+1}$ for all $k \in \field$ and $t \in \nonnegints$.
	Then, $\star$ is given by $x \star \gamma_t(a)x^t = 2ax^{t+1}$ and more generally by $x^s \star ax^t = 2^s ax^{t+1}$ for all $t,s \in \nonnegints$.

	Define a map $f: M \to \rationals[x]$ by $ax^t \mapsto x^{t}$ for all $t \in \nonnegints$.
	Note that $f$ is only well-defined as a linear map with $M$ and $\rationals[x]$ being seen as $\rationals$-vector spaces.
	In particular, the assignment $ax^t \mapsto x^{t}$ fails to commute with scalar multiplication over $\rationals[x]$,
	i.e.\ $x \cdot f(ax^t) \neq f(x \star ax^t)$ for all $t \in \nonnegints$ since $x \cdot f(ax^t) = x \cdot x^t = x^{t+1}$ but $f(x \star ax^t) = f(2ax^{t+1}) = 2x^{t+1}$.

	In contrast, define the map $g: M \to \rationals[x]$ by $ax^t \mapsto 2^{-t} x^t$.
	Since $\set{ax^t}$ and $\set{2^{-t}x^t}$ both serve as bases for $\rationals\ket{a}x^t$ and $\rationals x^t$ respectively, 
	$g$ is well-defined as a linear map between $\rationals$-vector spaces.
	Observe that for all $t \in \nonnegints$,
		$g(x \star ax^t) = x \cdot g(ax^t)$ since 
	\begin{align*}
		g(x \star ax^t) &= g(2ax^{t+1}) = 2^{-(t+1)} \cdot 2x^{t+1} = 2^{-t}x^{t+1}
		\\
		x \cdot g(ax^t) &= x \cdot 2^{-t}x^t = 2^{-t} x^{t+1}
	\end{align*}
	Since the action $\star$ of $\field[x]$ commutes with $g$,
		$g$ is also a $\rationals[x]$-module homomorphism.
	Since for all $t \in \nonnegints$, $g$ satisfies 
		$\rationals x \star \rationals\ket{a} x^t \subseteq \rationals\ket{a} x^{t+1}$,
		$g$ is also a graded $\field[x]$-module homomorphism.
	Since $g$ has an obvious inverse by $x^t \mapsto 2^t a x^t$,
		$M = \togrmod(W_\bullet)$ and $\rationals[x]$ are graded isomorphic with graded isomorphism given by $g: M \to \rationals[x]$.
	Note that the graded isomorphism by $g: M \to \rationals[x]$ may be easier to see if we choose the basis element $b_t := 2^t a$ for the vector space $W_t$ for each $t \in \nonnegints$. 
\end{example}

\spacer 
\noindent 
\textbf{About the Morphism Assignment of $\togrmod(-)$:}
The construction of the graded $\field[x]$-module homomorphism by $\togrmod$ as stated in \fref{defn:togrmod}(ii) can also be understood to happen in three stages.
Let $(V_\bullet, \alpha_\bullet)$ and $(W_\bullet, \gamma_\bullet)$ be persistence modules over $\field$.
Let $\phi_\bullet: V_\bullet \to W_\bullet$ be a persistence morphism with $\phi_\bullet = (\phi_t: V_t \to W_t)_{t \in \nonnegints}$.
\begin{enumerate}
	\item 
	At the level of $\field$-modules, i.e.\ $\field$-vector spaces:

	For each $t \in \nonnegints$, $\phi_t: V_t \to W_t$ must be a well-defined linear map between $\field$-vector spaces $V_t$ and $W_t$ by definition of persistence morphism.
	The canonical isomorphisms $V_t \to V_t x^t$ and $W_t \to W_t x^t$ given by $v \mapsto v x^t$ for $v \in V_t$ and $w \mapsto w x^t$ for $w \in W$ respectively
		define the linear map $\phi\graded_t: V_t x^t \to W_t x^t$ given by $v_t x^t \to \phi_t(v_t) x^t$.
	Then, the (internal) direct sum operation on $\field$-vector spaces induces $\phi\graded: \togrmod(V_\bullet) \to \togrmod(W_\bullet)$ 
	to be a linear map $\phi\graded: \bigoplus_{t \in \nonnegints} V_t x^t \to \bigoplus_{t \in \nonnegints} W_t x^t$ 
		between (graded) $\field$-vector spaces 
	given by $\phi\graded = \bigoplus_{t \in \nonnegints} \phi\graded_t$.

	\item 
	At the level of $\field[x]$-modules, not necessarily graded at this point:
	
	For $\phi\graded$ to be a well-defined $\field[x]$-module homomorphism,
		it must be shown that the action of $\field[x]$ on both $\togrmod(V_\bullet)$ and $\togrmod(W_\bullet)$ commute with $\phi\graded$.
	Since $\phi\graded$ is a well-defined linear map, 
		it should suffice to check that the following is true:
		\begin{equation*}
			\phi\graded(x^s \star_V v) = x^s \star_W \phi\graded(v)
			\quad\text{ for all }
			t,s \in \nonnegints \text{ and }
			v \in \togrmod(V_\bullet)
		\end{equation*}
		where $\star_V$ and $\star_W$ refer to the action of $\field[x]$ on $\togrmod(V_\bullet)$ and $\togrmod(W_\bullet)$ respectively.
	This condition can be equivalently stated as follows:
	\begin{equation*}
		\bigl( \phi_{t+s} \circ \alpha_{t+s, t} \bigr)(v_t)
		= 
		\bigl( \gamma_{t+s, t} \circ \phi_{t} \bigr)(v_t)
		\quad\text{ for all }
		t,s \in \nonnegints \text{ and for all } v_t \in V_t
	\end{equation*}
	Observe that, by \fref{defn:persmod-cat}(ii), a persistence morphism $\phi_\bullet$ must form commuting squares with the structure maps of $V_\bullet$ and of $W_\bullet$, 
		i.e.\ 
			$\phi_{t+s} \circ \alpha_{t+s,t} = \gamma_{t+s,s} \circ \phi_t$ for all $t,s \in \nonnegints$.
	Therefore, $\phi\graded: \togrmod(V_\bullet) \to \togrmod(W_\bullet)$ is a well-defined $\field[x]$-module homomorphism.

	\item 
	At the level of graded $\field[x]$-modules:

	Since $\phi\graded(V_t x^t) \subseteq W_t x^t$ for all $t \in \nonnegints$
	and $V_t x^t$ and $W_t x^t$ are exactly the homogeneous components of $\togrmod(V_\bullet)$ and $\togrmod(W_\bullet)$ of degree $t \in \nonnegints$ respectively,
		$\phi\graded$ is also a well-defined \textit{graded} $\field[x]$-module homomorphism.
\end{enumerate}
We provide an example of the morphism assignment of $\togrmod(-)$ below.

\begin{example}
	Let $(V_\bullet, \alpha_\bullet)$ and $(W_\bullet, \alpha_\bullet)$ be persistence modules over $\rationals$ with vector spaces given as follows with distinct indeterminates $a$, $b$, and $ab$:
	\begin{equation*}
		V_t = \begin{cases}
			0 	&\text{ if } t < 5 \\
			\rationals\ket{ab} &\text{ if } t \geq 5
		\end{cases} \\
		\qquad\text{ and }\qquad 
		W_t = \begin{cases}
			0 &\text{ if } t < 2 \\
			\rationals\ket{a,b} &\text{ if } t \geq 3
		\end{cases}
	\end{equation*}
	Define $\alpha_{t}: V_t \to V_{t+1}$ by $ab \mapsto ab$ for all $t \geq 5$ and $\gamma_t: W_t \to W_{t+1}$ by $a \mapsto a$ and $b \mapsto b$ for all $t \geq 3$.
	Let $V\graded$ and $W\graded$ be graded $\rationals[x]$-modules given by $V\graded := \togrmod(V_\bullet)$ and $W\graded := \togrmod(W_\bullet)$.
	We claim that $V\graded$ and $W\graded$ are graded isomorphic to the following:
	\begin{equation*}
		V\graded \upgraded\cong
			\Sigma^5 \rationals[x]\ket{ab}
			= \rationals[x]\ket{abx^5}
		\qquad\text{ and }\qquad 
		W\graded \upgraded\cong 
			\Sigma^3 \rationals[x]\ket{a,b}
			= \rationals[x]\ket{ax^3, bx^3}
	\end{equation*}
	where $\Sigma^k$ refers to the upward shift in grading (see \fref{defn:upwards-shift-graded-module}).
	That is, the relation is valid at the level of graded $\rationals[x]$-module homomorphisms and that the action of $\rationals[x]$ on $V\graded$ and $W\graded$ is given by $x \cdot abx^t = abx^t$ for all $t \geq 5$ and $x \cdot ax^t = ax^{t+1}, x \cdot bx^t = bx^{t+1}$ for all $t \geq 3$ respectively.

	For each $t \in \nonnegints$, define $\phi_t: V_t \to W_t$ to be the restriction of the map $\Phi: \rationals\ket{ab} \to \rationalsket{a,b}$ given by $ab \mapsto b-a$.
	We claim that $\phi_\bullet = (\phi_t)_{t \in \nonnegints}$ defines a persistence morphism $\phi_\bullet: V_\bullet \to W_\bullet$.
	Let $\phi\graded: V\graded \to W\graded$ be the graded $\rationals[x]$-module homomorphism 
	given by $\phi\graded := \togrmod(\phi_\bullet)$.
	Then, for all $t \in \nonnegints$ with $t \geq 5$:
	\begin{equation*}
		\phi\graded(abx^t) = \phi_t(ab)x^t = (b-a)x^t
	\end{equation*}
	Observe that $\phi\graded$ sends homogeneous elements of $V\graded$ of degree $t \geq 5$ to homogeneous elements of $W\graded$ of the same degree.
	Note that $V\graded$ has no homogeneous elements of degree $t < 5$. 
\end{example}

\spacer 

Next, we provide a description of $\topersmod: \catgradedmod{\field} \to \catpersmod$ below. 

\begin{definition}\label{defn:topersmod}
	Fix a field $\field$.
	Define $\topersmod: \catgradedmod{\field} \to \catpersmod$ 
	to be the following assignment of the objects and morphisms of $\catgradedmod{\field}$ to those of $\catpersmod$ respectively.
	\begin{enumerate}
		\item 
		\textit{Object Assignment:}
		Let $M$ be a graded $\field[x]$-module and let its homogeneous component be given by $M_t x^t$ for some $\field$-vector space $M_t$ for all $t \in \nonnegints$.
		Then, $M \upabelian= \bigoplus_{t \in \nonnegints} M_t x^t$.

		Define $\topersmod(M) =: (V_\bullet, \alpha_\bullet)$ to be the persistence module over $\field$ 
		such that for all $t \in \nonnegints$, the vector space $V_t$ is given by $V_t := M_t$,
		and for all $t,s \in \nonnegints$ with $t \leq s$,
			the structure map $\alpha_{s,t}: V_t \to V_s$
			is given by $\alpha_{s,t}(v_t) = v_{s}$ whenever 
			$x^{s-t} \cdot v_t x^t = v_s x^s$ with $v_s \in M_s$.

		\item 
		\textit{Morphism Assignment:}
		Let $f: M \to N$ be a graded $\field[x]$-module homomorphism between graded $\field[x]$-modules $M$ and $N$.
		Let the homogeneous component of $M$ and $N$ be given by $M_t x^t$ and $N_t x^t$ for some $\field$-vector spaces $M_t$ and $N_t$ for all $t \in \nonnegints$ respectively.

		Define $\topersmod(f) =: \phi_\bullet$ to be the persistence morphism $\phi_\bullet: V_\bullet \to W_\bullet$ between persistence morphisms $(V_\bullet, \alpha_\bullet) := \topersmod(M)$ and $(W_\bullet, \gamma_\bullet) := \topersmod(M)$ 
		with $\phi_t: V_t \to W_t$ given by 
		$\phi(v_t) := w_t$ whenever $f(v_t x^t) = w_t x^t$ for all $v_t \in V_t$ and for all $t \in \nonnegints$.
	\end{enumerate}
\end{definition}
\remarks{
	\item 
	As is the case of $\togrmod$, the notation $\topersmod$ is not standard or convention in persistence literature. 
	In fact, 
		notation for the functor $\catgradedmod{\field} \to \catpersmod$ in the category isomorphism was not explicitly identified in \cite{matrixalg:zomorodian}.

	\item 
	The graded modules and graded homomorphisms follow \fref{remark:assumption-on-graded-modules}, where $x^t$ is added to the notation of the elements of graded modules to help identify the degree of homogeneous elements.
}

\begin{proposition}\label{prop:topersmod-is-a-functor}
	Fix a field $\field$. The object and morphism assignment $\topersmod$ by \fref{defn:topersmod} determines a functor $\topersmod: \catpersmod \to \catgradedmod{\field}$.
\end{proposition}
\remark{
	A brief discussion of the proof is available in \cite[p8]{persmod:bubenik-homoloalg}.
}

Note that the arguments given for $\togrmod$ producing well-defined graded modules and graded homomorphisms, when presented in reverse order, 
	also tells us that application of $\topersmod$ results in well-defined persistence modules and persistence morphisms.
We provide an example of $\topersmod$ in action below.

\begin{example}
	Let $M$ be a graded $\rationals[x]$-module given as follows:
	\begin{equation*}
		M \upgraded= \Sigma^4\biggl(
			\rationals[x] \bigmod \rationals[x]\ket{x^5}
		\biggr)
	\end{equation*}
	Since $(x^5) = \rationals[x]\ket{x^5}$ is generated by a homogeneous element of $\rationals[x]$, 
	$M$ must be a graded $\rationals[x]$-module by \fref{prop:graded-submodule-hom-generators}.
	Let $\set{M_t}_{t \in \nonnegints}$ be the family of $\rationals$-vector spaces such that $M \upabelian= \bigoplus_{t \in \nonnegints} M_t x^t$
	and $M_t x^t$ is the homogeneous component of $M$ of degree $t \in \nonnegints$ for all $t \in \nonnegints$.
	Then, $M_t$ is given as follows:
	\begin{equation*}
		M_t x^t = \begin{cases}
			0 &\text{ if } t < 4 \\
			\rationals x^t &\text{ if } t \in \set{4, 5, \ldots, 8} = [4,9) \\
			0 &\text{ if } t \geq 9
		\end{cases}
	\end{equation*}
	That is, $M$ consists of $\rationals$-linear combinations of $x^4, x^5, x^6, x^7$ and $x^8$.
	The action of $\rationals[x]$ on $M$ is also given by $x \cdot x^t = x^{t+1}$ for all $t \in \set{4, \ldots, 8}$ and by 
	$x \cdot x^8 = 0$. 

	Let $(V_\bullet, \alpha_\bullet)$ be the persistence module over $\rationals$ given by $(V_\bullet, \alpha_\bullet) := \topersmod(M)$.
	Then, the vector spaces of $V_\bullet$ are given as follows:
	\begin{equation*}
		V_t = \begin{cases}
			0 &\text{ if } t < 4 \\
			\rationals &\text{ if } t \in [4,9) \\
			0 &\text{ if } t \geq 9
		\end{cases}
	\end{equation*}
	For $t \in \set{4, 5, \ldots, 7}$, the structure map $\alpha_t: V_t \to V_{t+1}$ is given by $V_t \ni 1 \mapsto 1 \in V_{t+1}$, i.e.\ $\alpha_t = \id_{\rationals}$.
	For $t = 8$, the structure map $\alpha_8: V_8 \to V_9$ sends every vector in $V_8$ to the zero vector in $V_9 = 0$.
	For $t \not\in [4,9)$, $\alpha_t: V_t \to V_{t+1}$ is trivial since $V_t = 0$.
	Equivalently, we have that for all $t,s \in \nonnegints$ with $t \leq s$:
	\begin{equation*}
		\alpha_{s,t} = \begin{cases}
			\id_{\rationals} &\text{ if } t \in [4,9) \\
			0 &\text{ otherwise }
		\end{cases}
	\end{equation*}
	We claim that $V_\bullet$ is isomorphic to the interval module $\intmod{[4,9)}$ over $\rationals$ as persistence modules.

	As a sidenote, 
	we can consider the elements of $M$ to be $\rationals$-linear combinations of $x^t$ with $t \in [4,9)$ since $(x^5)$ being a graded submodule of $\rationals[x]$ means that we can consider the homogeneous component of degree $t \in \nonnegints$ separately.
	More specifically: Let $N = (x^5)$. The homogeneous component $N_t x^t$ of $N$ is given by $N_t x^t = \rationals x^t$ if $t \geq 5$ and $N_t x^t = 0$ otherwise.
	Then, following \fref{prop:direct-sums-are-graded},
	\begin{equation*}
		\frac{\rationals[x]}{(x^5)}
		\upabelian= 
		\bigoplus_{t \in \nonnegints} 
		\left( \frac{\rationals x^t}{ N_t x^t } \right)
		= 
		\bigoplus_{t \in \nonnegints} 
		\left( \frac{\rationals }{ N_t } \right) x^t 
		= 
		\Biggl(
			\bigoplus_{t = 0}^{4} 
			\Bigl( \rationals \bigmod 0 \Bigr) x^t 
		\Biggr)
		\oplus 
		\Biggl(
			\bigoplus_{t = 5}^\infty 
			\Bigl( \rationals \bigmod \rationals \Bigr) x^t
		\Biggr)
		= 
		\bigoplus_{t = 0}^{4} \rationals x^t
	\end{equation*}
	Note that we write $\upabelian=$ since the direct sums only consider the action of $\field$ on the quotient module, not the action of $\field[x]$.
	Alternatively, since $(x^5)$ is a graded submodule, 
		the quotient $\rationals[x] \bigmod \rationals\ket{x^5}$ makes all $\rationals$-multiples of powers $x^t$ with $t \geq 5$ trivial
		and each coset of $\rationals[x] \bigmod \rationals\ket{x^5}$ is represented uniquely by a $\rationals$-linear combination of $x^t$ with $t \in \set{1,\ldots,4}$.
\end{example}

\spacer 

Finally, we state our theorem involving the isomorphism of categories between that of $\nonnegints$-indexed persistence modules over $\field$ and $\nonnegints$-graded $\field[x]$-modules.

\begin{theorem}\label{thm:cat-isom-persmod-grmod}
	Fix a field $\field$.
	The pair $\togrmod$ and $\topersmod$ determine an isomorphism of categories between the category $\catpersmod$ of persistence modules over $\field$ and the category $\catgradedmod{\field}$ on $\nonnegints$-graded $\field[x]$-modules, i.e.\ 
	\begin{equation*}
		\topersmod \circ \togrmod = \id_{\catpersmod} 
		\quad\text{ and }\quad 
		\togrmod \circ \topersmod = \id_{\catgradedmod{\field}}
	\end{equation*}
	where $\id_{\catpersmod}$ and $\id_{\catgradedmod{\field}}$ denote the identity functors on $\catpersmod$ and $\catgradedmod{\field}$ respectively.
	Note that this also implies that $\togrmod$ and $\topersmod$ form an equivalence of categories.
\end{theorem}
\remarks{
	\item 
	This theorem is a special case of \cite[Theorem 2.21]{persmod:bubenik-homoloalg}
	where $\catpersmod$
	and $\catgradedmod{\field}$ 
	are denoted as 
	$\catname{Mod}_R^{\catname{P}}$ and 
	$Gr^{\catname{P}}\text{-}\catname{Mod}_{R[U_0]}$ respectively
	with $\catname{P} = \posetN$,
	$R = \field$, and $\field[x] = R[U_0]$. 

	\item 
	A weaker version of this theorem is presented in \cite[Theorem 3.1]{matrixalg:zomorodian}, which only claims a category equivalence between the subcategory of $\catpersmod$ limited to \textit{finite-type} persistence modules over $\field$ and the subcategory of $\catgradedmod{\field}$ limited to \textit{finitely-generated} graded $\field[x]$-modules.
	\cite{matrixalg:zomorodian} claims that the Artin-Rees theorem in commutative algebra is sufficient for the proof.

	Alternatively, \cite{persmod:corbet} provides a detailed proof of \cite[Theorem 3.1]{matrixalg:zomorodian} (which it calls the ZC Representation Theorem) without using the Artin-Rees theorem.
	A sketch of an alternative proof that uses the Artin-Rees theorem is provided in \cite[Appendix C]{persmod:corbet}.
}

\spacer 
\noindent 
\textbf{About the Correspondence between Algebraic Constructions.}

\fref{thm:cat-isom-persmod-grmod} is particularly significant in this paper since it allows us to correspond algebraic constructions in $\catpersmod$, i.e.\ involving persistence modules, to those in $\catgradedmod{\field}$, i.e.\ involving graded modules.
In the second half of this section, we state the propositions that tell us which properties of the original persistence complex are preserved as we go from the category of persistence modules, to that of graded modules, then back to that of persistence modules.

We begin with a statement involving isomorphisms in $\catpersmod$ and $\catgradedmod{\field}$.

\begin{proposition}\label{prop:catequiv-preserves-isom}
	The functors $\togrmod(-)$ and $\topersmod(-)$ preserve
	and reflect isomorphisms, i.e.\ 
	\begin{enumerate}
		\item 
		A persistence morphism $\phi_\bullet: V_\bullet \to W_\bullet$ between persistence modules $V_\bullet$ and $W_\bullet$ is a persistence isomorphism if and only if 
		$\togrmod(\phi_\bullet): \togrmod(V_\bullet) \to \togrmod(W_\bullet)$ is a graded $\field[x]$-module isomorphism.

		\item  
		A graded $\field[x]$-module homomorphism $f: M \to N$ between graded $\field[x]$-modules $M$ and $N$ is a graded isomorphism if and only if 
		$\topersmod(f): \topersmod(M) \to \topersmod(N)$ is a persistence isomorphism.
	\end{enumerate}
\end{proposition}
\begin{proof}
	For part (i): For the forward direction, assume $\phi_\bullet: V_\bullet \to W_\bullet$ is a persistence isomorphism.
	By \cite[Lemma 1.3.8]{cattheory:rhiel}, which states that all functors preserve isomorphisms,
		$\togrmod(\phi_\bullet): \togrmod(V_\bullet) \to \togrmod(W_\bullet)$ must be a graded $\field[x]$-module isomorphism.
	For the converse, assume that $\togrmod(\phi_\bullet)$ is a graded isomorphism.
	By \fref{thm:cat-isom-persmod-grmod},
		$\topersmod \circ \togrmod = \id_{\catpersmod}$
		and 
	\begin{equation*}
		\bigl( \topersmod \circ \togrmod \bigr)(V_\bullet) = V_\bullet, 
		\quad 
		\bigl( \topersmod \circ \togrmod \bigr)(W_\bullet) = W_\bullet,
		\quad\text{ and }\quad
		\bigl( \topersmod \circ \togrmod \bigr)(\phi_\bullet)
		= \phi_\bullet
		.
	\end{equation*}
	By \cite[Lemma 1.3.8]{cattheory:rhiel} on $\topersmod(-)$,
		$\phi_\bullet$ must be a persistence isomorphism.
	A similar argument applies for part (ii).
\end{proof}
\remark{
	This results holds even if $\togrmod(-)$ and $\topersmod(-)$ only form an equivalence of categories since \cite[Theorem 1.5.9]{cattheory:rhiel} states that any functor defining an equivalence of categories is full and faithful and \cite[Exercise 1.5.iv]{cattheory:rhiel} states that full and faithful functors reflect isomorphisms.
}

Since a graded invariant factor decomposition of a graded module is given by a graded isomorphism,
	\fref{prop:catequiv-preserves-isom} tells us that the application of $\topersmod(-)$ on such a decomposition
	will produce a persistence isomorphism.
Moreover, we can show that this graded decomposition corresponds to an interval decomposition.
Below, we describe the relationship between interval modules and the cyclic summands of graded invariant factor decompositions.

\begin{lemma}\label{lemma:intmods-gradedmods-corr}
	Fix a field $\field$.
	The interval persistence modules correspond to shifted cyclic graded $\field[x]$-modules,
		i.e.\ for all $t,s \in \nonnegints$,
	\begin{equation*}
		\intmod{[s,\infty)} = \topersmod\bigg(
			\Sigma^{s} \field[x]
		\bigg)
		\qquad\text{ and }\qquad
		\intmod{[s, s + t)} = \topersmod\bigg(
			\Sigma^{s} \Bigl( \field[x] \bigmod (x^{t}) \Bigr)
		\bigg)
	\end{equation*}
	where 
	$\intmod{J}$
	denotes the interval modules over the interval $J \subseteq \nonnegints$
	(see 
	\fref{defn:persistence-interval-modules}),
	$\Sigma^s(-)$ denotes an upwards $s$-shift in grading (see \fref{defn:upwards-shift-graded-module})
	and $(x^t) = \field[x]\ket{x^t}$ is the graded $\field[x]$-module generated by $x^t$.
\end{lemma}
\begin{proof}
	Let $t,s \in \nonnegints$.
	First, we want to show that $\intmod{[s,\infty)} \cong \topersmod(\Sigma^s \field[x])$.
	Note that $\Sigma^s \field[x] = \field[x]\ket{x^s}$ is a graded $\field[x]$-module with 
	homogeneous component given by $\field x^r$ if $r \geq 5$ and trivial otherwise.
	The action of $\field[x]$ on $\Sigma^s \field[x]$ is given by $x \cdot x^r = x^{r+1}$ for all $r \geq 5$.
	Equivalently, $x^{q-r} \cdot x^r = x^{q}$ for all $r,q \in [s,\infty)$ with $r \leq q$.
	Let $(V_\bullet, \alpha_\bullet) := \topersmod\bigl( \Sigma^s \field[x] \bigr)$.
	For all $r,q \in \nonnegints$ with $r \leq q$,
	the vector space $V_r$ and the structure map $\alpha_{q,r}: V_r \to V_q$ are given as follows:
	\begin{equation*}
		V_r = \begin{cases}
			\field &\text{ if } r \in [s,\infty) \\
			0 &\text{ otherwise }
		\end{cases}
		\quad\text{ and }\quad 
		\alpha_{q,r} = \begin{cases}
			\id_{\field} &\text{ if } r,q \in [s,\infty) \\
			0 & \text{ otherwise }
		\end{cases}
	\end{equation*}
	By \fref{defn:persistence-interval-modules},
		$V_\bullet = \intmod{[s,\infty)}$.

	Next, we want to show that 
	$\intmod{[s,s+t)} = \topersmod(M)$
	with $M := \Sigma^s \bigl( \field[x] \bigmod (x^t) \bigr)$.
	By \fref{prop:cokernels-of-graded-modules}, 
		the homogeneous component of $\field[x] \bigmod (x^t)$ of degree $r \in \nonnegints$ is given by the quotient of the homogeneous component of $\field[x]$ of degree $r$ by that of $(x^t) = \field[x]\ket{x^t}$.
	That is, the homogeneous component $M_r x^r$ of $M$ of degree $r \in \nonnegints$ is as follows:
	\begin{equation*}
		M_r x^r = \left\{\begin{array}{
			r@{\,}l@{\,}l c ccl
		}
			&0& &\cong& 0
			&\text{ if }& r \in [0,s) 
			\\
			\field x^r &\bigmod& 0 &\cong& \field x^r 
				&\text{ if }& r \in [s,s+t) 
			\\
			\field x^r &\bigmod& \field x^r &\cong& 0
				&\text{ if }& r \in [s+t,\infty)
		\end{array}\right.
	\end{equation*}
	The action of $\field[x]$ on $M$ is given by $x \cdot x^r = x^{r+1}$ for $r \in [s, s+t-1)$ and by $x \cdot x^{s+t-1} = 0$.
	Equivalently, for all $r,q \in \nonnegints$ with $r \leq q$,
		$x^{q-r} \cdot x^r = x^q$ if $r,q \in [s,s+t)$
		and 
		$x^{q-r} \cdot x^r = 0$ otherwise.
	Let $(V_\bullet, \alpha_\bullet) := \topersmod(M)$.
	Then, the vector spaces $V_r$ and structure maps $\alpha_{q,r}: V_r \to V_{q}$ of $V_\bullet$ are given as follows:
	\begin{equation*}
		V_r = \begin{cases}
			\field &\text{ if } r \in [s,s+t) \\
			0 &\text{ otherwise }
		\end{cases}
		\quad\text{ and }\quad 
		\alpha_{q,r} = \begin{cases}
			\id_{\field} &\text{ if } r,q \in [s,s+t)
			\text{ with } r \leq q \\
			0 &\text{ otherwise }
		\end{cases}
	\end{equation*}
	By \fref{defn:persistence-interval-modules},
		$V_\bullet = \intmod{[s,s+t)}$.
\end{proof}

Note that we can also determine the homogeneous component of degree $r \in \nonnegints$ of $M := \Sigma^s \bigl( \field[x] \bigmod (x^t) \bigr)$ by using \fref{prop:shift-distribute-over-quotient} and distributing the upwards shift $\Sigma^s(-)$ across the quotient, i.e.\
\begin{equation*}
	\Sigma^s \left( \frac{ \field[x] }{ (x^t) } \right)
	\upgraded\cong\,
	\frac{ \Sigma^s \field[x] }{ \Sigma^s (x^t) }
	= \frac{ \field[x]\ket{x^s} }{ \field[x]\ket{x^{t+s}} }
	= \Bigl\{
		k x^r \in \field[x] : k \in \field \text{ and } r \in [s,s+t)
	\Bigr\}
\end{equation*}
More rigorously, any coset of $M$ as a quotient module can be represented uniquely by $\field$-linear combinations of powers $x^r$ with $r \in [s,s+t)$.
Following \fref{remark:assumption-on-graded-modules}, the only non-trivial homogeneous components of $M$ are those of degree $r \in [s,s+t)$ given exactly by $\field x^r$.

To state the correspondence between graded invariant factor decompositions and interval decompositions of persistence modules, we need the following result involving finite direct sums.

\newcommand{\sideP}[1]{\sublabel{P}{#1}}
\newcommand{\sideG}[1]{\sublabel{G}{#1}}

\begin{proposition}\label{prop:distribute-over-direct-sums}
	The functors $\togrmod(-)$ and $\topersmod(-)$ distribute over finite direct sums, i.e.\ 
	given persistence modules $V_\bullet$ and $W_\bullet$ and graded $\field[x]$-modules $M$ and $N$,
		\begin{equation*}
		\begin{array}{c@{\ }c@{\ }c l}
			\togrmod\paren{
				\persmod{V} \sideP\oplus \persmod{W}
			} 
			&\upgraded\cong& 
			\togrmod(\persmod{V}) \sideG\oplus \togrmod(\persmod{W})
			\\[2pt]
			\topersmod\paren{
					M \sideG\oplus N
			} &\uppersmod\cong& 
				\topersmod(N) \sideP\oplus \topersmod(N)
		\end{array}
		\end{equation*}
		where $\sideG\oplus$ refers to a direct sum of graded $\field[x]$-modules 
		and $\sideP\oplus$ refers to that of persistence modules.
\end{proposition}
\remark{
	This can be seen as a consequence $\togrmod(-)$ and $\topersmod(-)$ forming an equivalence of categories by \fref{thm:cat-isom-persmod-grmod}.
	We are unable to provide a proof of this (or a reference for such) at this moment but we have outlined two possible arguments below.
	\begin{enumerate}
		\item 
		We can use \cite[Lemma 3.3.6]{cattheory:rhiel}, which roughly states that an equivalence of categories preserves all limits and colimits,
			as described in \cite[Chapter 3]{cattheory:rhiel}.
		Then, as argued in the proof of 
			\cite[Proposition 4.5.10]{cattheory:rhiel},
		direct sums and kernels in abelian categories are finite limits.
		Since $\togrmod(-)$ and $\topersmod(-)$ form an equivalence of categories by \fref{thm:cat-isom-persmod-grmod},
			these must both preserve finite direct sums.

		\item 
		Alternatively, we can use \cite[Proposition 4.4.5]{cattheory:rhiel}, which roughly states that any equivalence of categories determine an adjoint equivalence (also see \cite[Footnote 37, p30]{cattheory:rhiel}),
		and \cite[Corollary 4.5.11]{cattheory:rhiel}, which then implies that these functors are exact and additive.

		Note that \cite{cattheory:rhiel} defines an additive functor to be a functor that preserves direct sums.
		In contrast, \cite[Corollary 5.88]{algtopo:rotman} states this to be a property of additive functors.
		A more detailed discussion on properties preserved by exact and additive functors can also be found in \cite{cattheory:bailie} with the above result given as part of \cite[Theorem 27]{cattheory:bailie}.

	\end{enumerate}
	We are unsure if the above arguments apply for arbitrary direct sums. However, this is not an issue for us since arbitrary direct sums of persistence modules are not needed for the matrix reduction algorithm for persistent homology.
}

Below, we state the correspondence between interval decompositions of persistence modules and graded invariant factor decompositions of graded modules.

\begin{statement}{Corollary}\label{cor:interval-decomp-from-structure-theorem}
	Let $(V_\bullet, \alpha_\bullet)$ be a persistence module.
	Assume that $M := \togrmod(V_\bullet)$ admits the following graded invariant factor decomposition, 
		as described in the Graded Structure Theorem (\fref{thm:graded-structure-theorem}):
	\begin{equation*}
		\togrmod(V_\bullet) = M 
		\,\upgraded\cong\, 
		\Sigma^{s_1}\!\left( \frac{\field[x]}{(x^{t_1})} \right)
		\oplus \cdots \oplus
		\Sigma^{s_r}\!\left( \frac{\field[x]}{(x^{t_r})} \right)
		\oplus
		\Sigma^{s_{r+1}} \field[x] 
		\oplus \cdots \oplus 
		\Sigma^{s_m} \field[x]
	\end{equation*}
	with invariant factors $\set{x^{t_1}, \ldots, x^{t_r}}$
	and grading shifts $\set{s_1, \ldots, s_r, \ldots, s_m}$.
	Then, $V_\bullet$ admits the following interval decomposition, as defined in \fref{defn:interval-decomposition},
	\begin{equation*}
		V_\bullet 
		= \topersmod(M) 
		\,\uppersmod\cong\, 
			\intmod{[s_1, s_1 + t_1)}
			\,\oplus \cdots \oplus\, 
			\intmod{[s_r, s_r + t_r)}
			\,\oplus\, 
			\intmod{[s_{r+1},\infty)}
			\,\oplus \cdots \oplus\, 
			\intmod{[s_{m}, \infty]}
	\end{equation*}
	with persistence barcode given by 
	$\barcode(V_\bullet) = 
	\bigl\{ [s_1,s_1 + t_1), \ldots, [s_r, s_r+t_r] \bigr\} \cup \bigl\{ [s_{r+1},\infty), \ldots, [s_m, \infty) \bigr\}$.
\end{statement}
\begin{proof}
	\fref{prop:catequiv-preserves-isom} implies that 
		the graded isomorphism of $M$ to its graded invariant factor decomposition determines a persistence isomorphism by application of $\topersmod(-)$.
	Then, we distribute $\topersmod(-)$ over the graded decomposition using \fref{prop:distribute-over-direct-sums}
	and apply \fref{lemma:intmods-gradedmods-corr} to each cyclic summand.
\end{proof}

The Graded Structure Theorem (\fref{thm:graded-structure-theorem}) guarantees the existence of graded invariant factor decompositions for finitely generated graded $\field[x]$-modules.
This, along with \fref{cor:interval-decomp-from-structure-theorem}, gives us a benchmark for the existence of interval decompositions.
We state this in more detail below.

\begin{proposition}\label{prop:finite-type-has-interval}
	Let $(V_\bullet,\alpha)$ be a persistence module over $\field$.
	If $V_\bullet$ is a finite-type persistence module,
	then $\togrmod(V_\bullet)$ is a finitely generated graded $\field[x]$-module and $V_\bullet$ admits an interval decomposition.
\end{proposition}
\begin{proof}
	Assume $(V_\bullet, \alpha_\bullet)$ is a finite-type persistence module.
	By \fref{defn:persmod-constant-on-interval},
		the vector space $V_t$ is finite dimensional for all $t \in \nonnegints$ 
		and $V_\bullet$ is constant on $[N,\infty)$ for some $N \in \nonnegints$.
	For each $t \in [0,N]$, let $\basis{B}_t$ be a basis of $V_t$. 
	Let $\basis{B} := \bigcup_{t=0}^N \basis{B}_t$.
	Note that $\basis{B}$ is a finite set.

	Let $M := \togrmod(V_\bullet)$.
	By definition of $\togrmod(-)$,
		the homogeneous component of $M$ of degree $t \in \nonnegints$ is given exactly by $V_t x^t$
		and the action of $\field[x]$ on $M$ is given by 
		$x^s \cdot v_t x^t = \alpha_{s+t, t}(v_t) x^{t+s}$ with $\alpha_{s+t, t}: V_t \to V_{s+t}$ for all $v_t \in V_t$ and $t \in \nonnegints$.
	Let $t \in [0,N]$.
		Note that for all $b \in \basis{B}_t$,
			$b x^t \in V_t x^t \subseteq M$ by construction.
		Define $\basis{B}_t x^t := \set{ bx^t \in M : b \in \basis{B}_t }$
		and $\basis{B}\graded := \bigcup_{t=0}^N \basis{B}_t x^t$. 
	We want to show that $\basis{B}\graded$ generates $M$, i.e.\ 
	$M = \set{ f(x) \cdot bx^t : f(x) \in \field[x], bx^t \in \basis{B}\graded}$. 

	Since $M \upabelian= \bigoplus_{t \in \nonnegints} V_t x^t$,
		$m = \sum_{t \in \nonnegints} v_t x^t$
		for a unique set $\set{v_t}_{t \in \nonnegints}$ of elements $v_t \in V_t$, only finitely many of which are nonzero.
	Then, it suffices to check if each homogeneous component $v_t x^t$ is generated by $\basis{B}\graded$.
	Let $m = v_t x^t \in V_t x^t \subseteq M$ for some $t \in \nonnegints$. 
	\begin{enumerate}
		\item 
		Assume $t \in [0,N]$.
		Since $\basis{B}_t$ is a basis of $V_t$ by assumption,
			$\basis{B}_t x^t$ must also be a basis of $V_t x^t$, viewing $V_t x^t$ as an $\field$-vector space.
		Therefore, there exists an $\field$-linear combination in $\basis{B}_t x^t \subseteq \basis{B}\graded$ that equals $v_t x^t$.

		\item 
		Assume $t \in [N+1,\infty)$.
		By assumption of $V_\bullet$ being constant on $[N,\infty)$, the structure map $\alpha_{t,N}: V_N \to V_t$ must be an $\field$-vector space isomorphism.
		Let $v_N \in V_N$ such that $\alpha_{t,N}(v_N) = v_t$.
		Then, $v_N x^N \in M$ and 
			$x^{t-N} \cdot v_N x^N = \alpha_{t,N}(v_N) x^t = v_t x^t$.
	\end{enumerate}
	Therefore, any element $m \in M$ is an $\field[x]$-linear combination on $\basis{B}\graded$.
	Since $M \upabelian= \bigoplus_{t \in \nonnegints} V_t x^t$ by definition,
		all $\field[x]$-linear combinations on $\basis{B}\graded$ must be on $M$.
	Then, 
		\begin{equation*}
			M = 
			\field[x] \basis{B}\graded 
			=
			\Biggl\{\,
				\sum_{t=0}^{\infty} f_t(x) \cdot b_tx^t : 
					f_t(x) \in \field[x] 
					\text{ finitely many of which are nonzero and }
					b_tx^t \in \basis{B}\graded
			\,\Biggr\}
		\end{equation*}
	and $\basis{B}$ is a finite system of generators for $M$.
	Therefore, 
	$M$ is a finitely generated graded $\field[x]$-module
	and admits a graded invariant factor decomposition by the Graded Structure Theorem (\fref{thm:graded-structure-theorem}).
	By \fref{cor:interval-decomp-from-structure-theorem},
		there exists an interval decomposition for $V_\bullet$.
\end{proof}
\remark{
	If we consider functors of the form $\posetN \to \catmod{\ints}$ to be persistence modules
	and let $\catpersmod[\ints]$ be the corresponding functor category,
	we can modify \fref{thm:cat-isom-persmod-grmod} to state an isomorphism of categories between $\catpersmod[\ints]$ and the category $\catgradedmod{\ints}$ of $\nonnegints$-graded $\ints[x]$-modules.
	However, while the graded $\ints[x]$-module constructed by applying $\togrmod(-)$ to some finite-type persistence module may be finitely generated, 
		the Graded Structure Theorem still would not apply to said graded $\ints[x]$-module since $\ints[x]$ is not a PID.
	Therefore, we cannot use \fref{prop:finite-type-has-interval} for the existence of interval decompositions for persistence modules of the form $\posetN \to \catmod{\ints}$.
}

\spacer 

Finally, we state how the chain homology of persistence modules and that of graded modules interact with the equivalence of categories given by $\togrmod(-)$ and $\topersmod(-)$ below.

\begin{proposition}\label{prop:catequiv-preserves-chains}
	The functors $\togrmod$ and $\topersmod$ preserve chain complexes and commute with the homology functor,
	i.e.\ let $H_n^{\suplabel{Pers}}(-): \catchaincomplex{\catpersmod} \to \catpersmod$ denote the $n$\th homology functor on persistence complexes 
	and let $H_n^{\suplabel{Gr}}(-): \catchaincomplex{\catgradedmod{\field}} \to \catgradedmod{\field}$ denote that on graded chain complexes.
	Then,
	\begin{enumerate}
		\item 
		Let $(V^\bullet_\ast, \boundary^\sbullet_\ast) = (V^\bullet_{n}, \boundary^\sbullet_{n})_{n \in \ints}$ be a persistence complex with 
			persistence modules $V^\bullet_{n}$ over $\field$
			and persistence morphisms 
			$\boundary^\sbullet_{n}: V^\bullet_{n} \to V^\bullet_{n-1}$.
		For all $n \in \ints$,
		there exists a graded isomorphism such that 
		\begin{equation*}
			\togrmod\Bigl(
				H_n^{\suplabel{Pers}}\bigl( V^\bullet_\ast, \boundary^\sbullet_\ast \bigr)
			\Bigr)
			\upgraded\cong
			H_n^{\suplabel{Gr}}(M_\ast, d_\ast)
		\end{equation*}
		where $(M_\ast, d_\ast) = (M_n, d_n)_{n \in \ints}$ is the graded chain complex generated by component-wise application of $\togrmod(-)$
		with $M_n := \togrmod(V^\bullet_n)$ and $d_n: M_n \to M_{n-1}$ by 
		$d_n := \togrmod(\boundary^\sbullet_n)$ for all $n \in \ints$.

		\item 
		Let $(M_\ast, d_\ast) = (M_n, d_n)_{n \in \ints}$ be a graded chain complex with graded $\field[x]$-modules $M_n$ and graded homomorphisms $d_n: M_n \to M_{n-1}$.
		For all $n \in \ints$, there exists a persistence isomorphism such that 
		\begin{equation*}
			\topersmod\Bigl(
				H_n^{\suplabel{Gr}}( M_\ast, d_\ast )
			\Bigr)
			\uppersmod\cong
			H_n^{\suplabel{Pers}}\bigl(
				V^\bullet_\ast, \boundary^\sbullet_\ast
			\bigr)
		\end{equation*}
		where $(V^\bullet_\ast, \boundary^\sbullet_\ast) = 
		(V_n^\bullet, \boundary_n^\sbullet)_{n \in \ints}$ is the persistence complex generated by component-wise application of $\topersmod$,
		with $V_n^\bullet := \topersmod(M_n)$ and 
			$\boundary_n^\sbullet := \topersmod(d_n)$
			for all $n \in \ints$. 
	\end{enumerate}
\end{proposition}
\remark{
	We are unable to provide a rigorous proof (or a reference for such) at this moment.
	However, if we can claim that $\togrmod(-)$ and $\topersmod(-)$ are exact functors, 
		following the remarks under \fref{prop:distribute-over-direct-sums}, 
	then the above result is stated as part of \cite[Theorem 27]{cattheory:bailie}.
	A similar result is listed as \cite[Exercise 6.8, p339]{cattheory:rotman}, which claims that exact additive functors between categories of modules over different rings commute with homology.
	We believe this applies more generally to exact additive functors between abelian categories, as discussed on the introduction of \cite[Chapter 6]{algtopo:rotman}.
}

We use the proposition above later in \fref{section:simplicial-persistent-homology} in the context of simplicial homology. \clearpage

\onlyifstandalone{\input{"../BX. Reference Page.tex"}}

\onlyifstandalone{\setcounter{chapter}{2}} 
\chapter{Filtrations and Persistent Homology}
\label{chapter:filtrations-and-pershoms}

Persistent homology theory is interested in the characterization, calculation, and representation of the persistent homology of filtrations of topological spaces.


Generally speaking, a \textit{filtration} refers to any 
	collection $\set{C_t : t \in \Lambda}$
 of objects, indexed over some set $\Lambda$, along with a set of subobject relations. 
The indexing set $\Lambda$ is usually equipped with a partial order $\leq$ that determines the expected set of subobject relations on $\set{C_t: t \in \Lambda}$, i.e.\ the subobject relation $C_t \subseteq C_s$ is present if and only if $t \leq s$ in the poset $(\Lambda, \leq)$.
Then, the persistent homology of a filtration $\set{X_t: t \in \Lambda}$ of topological spaces $X_t$ refers to some characterization of the following collections over all dimensions $n \in \nonnegints$:
\begin{equation*}
	\left\{\,\,\begin{gathered}
		\text{homology groups} \\
		H_n(X_t) : t \in \Lambda
	\end{gathered}\,\,\right\}
	\quad\text{ and }\quad
	\left\{\,\,\begin{gathered}
		\text{maps on homology induced by inclusions} \\
		H_n(X_t) \to H_n(X_s) : t,s \in \Lambda \text{ with } t \leq s
	\end{gathered}\,\,\right\}
\end{equation*}
In practice, the homology groups are taken with coefficients in a field $\field$ (usually $\field = \ints_p$ for prime $p$), and the calculation of persistent homology is done at the level of ranks, i.e.\ we want to find the following quantities:
\begin{equation*}
	\biggl\{\,\,
		\rank\bigl( H_n(X_t; \field) \bigr) : t \in \Lambda 
	\,\,\biggr\}
	\quad\text{ and }\quad 
	\biggl\{\,\, 
		\rank\Bigl(
			H_n(X_t; \field) \Xrightarrow{i_\ast} H_n(X_s; \field)
		\Bigr) : t,s \in \Lambda \text{ with } t \leq s
	\,\,\biggr\}
\end{equation*}
These are then represented succinctly using a multiset of intervals in $\Lambda$ (as a poset) called the \textit{persistence barcode} of the filtration, which is different albeit very similar to the persistence barcode of a persistence module.

In this chapter, we discuss the key ideas and constructions presented in the paper 
\textit{Computing Persistent Homology} \cite{matrixalg:zomorodian} by Afra Zomorodian and Gunnar Carlsson.
In particular, we restrict our attention to 
$\nonnegints$-indexed filtrations $K_\bullet := \set{K_t}_{t \in \nonnegints}$ of (abstract) simplicial complexes $K_t$ such that $K_t \subseteq K_{t+1}$ for all $t \in \nonnegints$ and study the following sequence of homology groups and induced maps:
\begin{equation*}
	H_n(K_0; \field) \Xrightarrow{\quad i_0^\ast \quad} 
	H_n(K_1; \field) \Xrightarrow{\quad i_1^\ast \quad} 
	H_n(K_2; \field) \Xrightarrow{\quad i_2^\ast \quad} 
	H_n(K_3; \field) \Xrightarrow{\quad i_3^\ast \quad} 
	\cdots
\end{equation*}
where $i_t: K_t \to K_{t+1}$ denotes the inclusion map for all $t \in \nonnegints$.
We discuss these notions relative to the characterization of persistence modules presented in \fref{chapter:persistence-theory},
i.e.\ as functors of the form $\posetN \to \catvectspace$.
This chapter is structured as follows:

\begin{enumerate}[chapterdecomposition]
	\item 
	We characterize simplicial filtrations, i.e.\ filtrations of simplicial complexes, as functors of the form $K_\bullet: \posetN \to \catsimp$
	where $\catsimp$ denotes the category of (abstract) simplicial complexes and simplicial maps.
	We also discuss the notion of \textit{finite-type filtrations}.	

	\item 
	We formalize the notion of the persistent homology of a simplicial filtration $\filt{K}$ with coefficients in a field $\field$ 
	by defining a persistence module over $\field$ called the $n$\th persistent homology module $H_n(K_\bullet;\field): \posetN \to \catvectspace$ for each dimension $n \in \ints$.

	We also identify a number of relevant terminology and interpretations involving $H_n(K_\bullet;\field)$, and show that the interval decomposition of $H_n(K_\bullet;\field)$ exists assuming $K_\bullet$ is a finite-type filtration.
	
	\item 
	We extend the simplicial chain complex of (abstract) simplicial complexes, as discussed in \fref{section:simplicial-homology},
	to the case of persistence modules 
	and construct a chain complex of persistence modules called 
	the simplicial persistence complex $C_\ast(K_\bullet;\field) = (C_n(K_\bullet;\field), \boundary_n^\sbullet)_{n \in \ints}$.
	The $n$\th chain homology of this persistence complex is then shown to be isomorphic to $H_n(\filt{K};\field)$.
 
	We also discuss how the isomorphism of categories between $\catpersmod$ and $\catgradedmod{\field}$ discussed in \fref{section:cat-equiv-graded-modules}
	allows us to calculate the $n$\th chain homology of $C_\ast(K_\bullet;\field)$ at the level of graded $\field[x]$-modules.
	This will serve as the basis for the \textit{matrix reduction algorithm for persistent homology}, later discussed in \fref{chapter:matrix-calculation}.
\end{enumerate}

\clearpage

\section{Filtrations of Simplicial Complexes}
\label{section:filtrations}

In this section, we discuss specific $\nonnegints$-indexed collections of simplicial complexes $\set{K_t}_{t \in \nonnegints}$ called \textit{filtrations} and provide a characterization of these as functors of the form $\posetN \to \catsimp$.

\begin{definition}\label{defn:filtration}
	A \textbf{simplicial filtration} $K_\bullet$ of a simplicial complex $K$ 
	is a functor $\filt{K}: \posetN \to \catsimp$ with the following properties:
	\begin{enumerate}
		\item 
		For all $t \in \nonnegints$, 
			$\filt{K}(t)$ is a subcomplex of $K$.
		For brevity, 
		we often write $K_t := \filt{K}(t)$,
			i.e.\ the bullet is replaced with $t \in \nonnegints$.
		The \textbf{index} $t \in \nonnegints$ of the simplicial complex $K_t$ in $\filt{K}$ is sometimes called the \textbf{scale} or \textbf{parameter} of $K_t$ in $\filt{K}$.
		
		\item 
		For all $t,s \in \nonnegints$ with $t \leq s$, 
			$K_t$ is a subcomplex of $K_s$
			and 
			$\filt{K}(t \to s): K_t \to K_s$ is exactly the inclusion map $K_t \hookrightarrow K_s$.
		Let $i^{s,t}: K_t \to K_s$ and $i^{t}: K_t \to K_{t+1}$
		denote the inclusion maps 
		$i^{s,t} := K_\bullet(t \to s)$ and $i^{t} := K_\bullet(t \to t+1)$ respectively.

		\item 
		For each simplex $\sigma \in K$,
			there must exist $t \in \nonnegints$ such that $\sigma$ is a simplex of $K_t$.
	\end{enumerate} 
	When the relation to a simplicial complex $K$ is clear from context, we may say $K_\bullet$ is a filtration of $K$ for brevity.
\end{definition}
\remark{
	Condition (iii) above determines that each simplicial filtration $\filt{K}$ corresponds to a unique simplicial complex, i.e.\ 
	$K = \bigcup_{t \in \nonnegints} K_\bullet(t) = \bigcup_{t \in \nonnegints} K_t$ if $K_\bullet$ is a filtration of $K$.
	It serves a similar purpose as condition \fref{defn:simp-complex-abstract}(i) for simplicial complexes.
}

Let $\filt{K}$ be a filtration of a simplicial complex $K$.
Observe that the subcomplex relation $K_t \subseteq K_s$ for all $t,s \in \nonnegints$ with $t \leq s$ stated in \fref{defn:filtration}(ii) implies that for all $t \in \nonnegints$, $K_t$ is a subcomplex of $K_{t+1}$.
Therefore, the simplicial filtration $\filt{K}$ determines the following nested sequence of simplicial complexes:
\begin{equation*}
	K_0 \,\subseteq\,
	K_1 \,\subseteq\,
	K_2 \,\subseteq\,
	K_3 \,\subseteq\,
	K_4 \,\subseteq\, \cdots
\end{equation*}
We can also show that increasing nested sequences of simplicial complexes determine simplicial filtrations, as we expect it should.
We state this in a proposition below.

\begin{proposition}
	Let $\set{K_t}_{t \in \nonnegints}$ be a collection of simplicial complexes such that $K_t \subseteq K_{t+1}$ for all $t \in \nonnegints$.
	Then, $\set{K_t}_{t \in \nonnegints}$ determines a filtration $\filt{K}$ of the simplicial complex $K := \bigcup_{t \in \nonnegints} K_t$ by $\filt{K}(t) := K_t$.
\end{proposition}
\begin{proof}
	Note that the union of simplicial complexes is a simplicial complex. Therefore, $K$ is well-defined.
	For each $t \in \nonnegints$, let $K_\bullet(t) := K_t$. Observe that $K_t$ is a subcomplex of $K$.
	For all $t,s \in \nonnegints$ with $t \leq s$,
		define $K_\bullet(t \to s)$ to be the inclusion map $K_t \hookrightarrow K_s$, which exists since $K_t \subseteq K_s$ by assumption.

	Let $t,r,s \in \nonnegints$ such that $t \leq r \leq s$.
	Since the composition $(K_t \hookrightarrow K_s) \circ (K_s \hookrightarrow K_r)$ of inclusion maps is exactly the inclusion $K_t \hookrightarrow K_r$, the composition axiom as stated in \fref{defn:functor}(iii) is satisfied.
	Since the inclusion map $K_t \hookrightarrow K_t$ is exactly the identity map, the identity axiom as stated in \fref{defn:functor}(iv) is satisfied.
	Therefore, $K_\bullet$ is a well-defined functor.
	The conditions of \fref{defn:filtration} are also all satisfied.
\end{proof}

The subcomplex relations also allow us to characterize a simplicial filtration $\filt{K}$ by identifying the simplicial complex $K_0$ at index $0 \in \nonnegints$ and specifying the collection $\set{\sigma_{t,i}}$ of simplices that is appended to $K_{t-1}$ to form $K_t = K_{t-1} \cup \set{\sigma_{t,i}}$ for each $t \geq 1$.
Note that each $K_t$ is required to be a simplicial complex by definition and sequentially appending an arbitrary set of simplices to some simplicial complex $L$ generally does not result in a simplicial filtration.
We provide an example of a simplicial filtration below.

\begin{example}\label{ex:baby-filtration}
	Let the simplicial complex $K$ with vertex set $\Vertex(K) = \set{a,b,c,d}$ and let the filtration $\filt{K}$ on $K$ be given by the following illustrations:
	\begin{center}
		\includegraphics[height=1in]{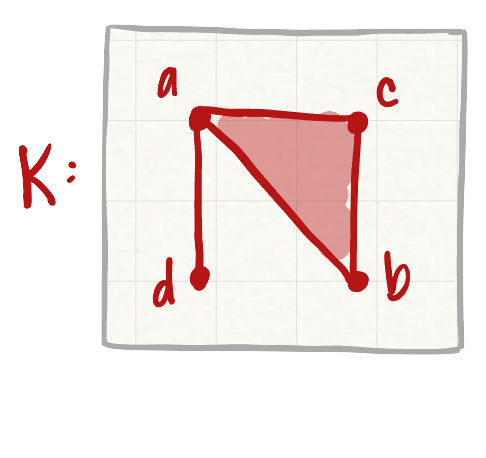}
		\quad 
		\includegraphics[height=1in]{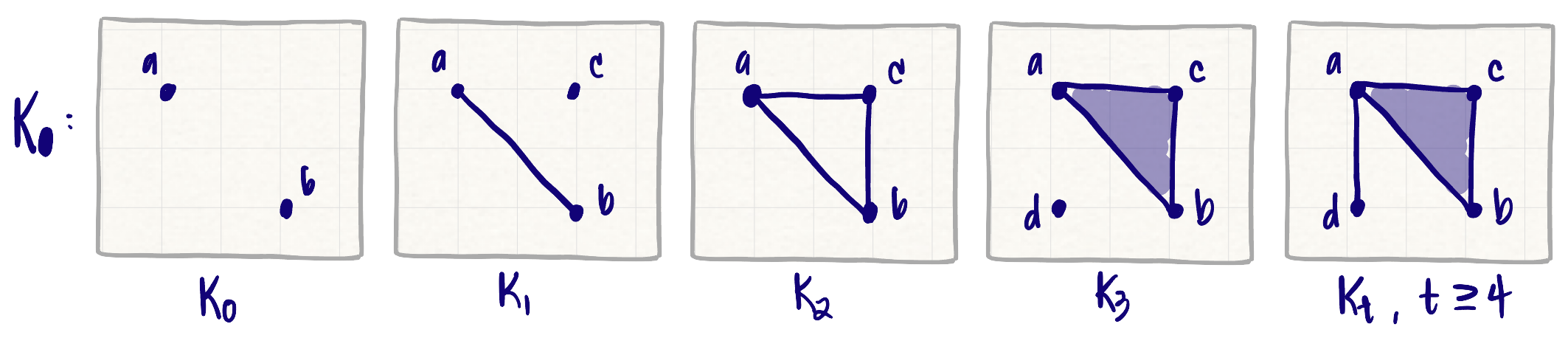}
	\end{center}
	The simplicial complexes $K_t$ of $\filt{K}$ can be described abstractly as follows, with the simplices of $K_t$ written as strings of vertices, following the remarks under \fref{defn:simp-complex-abstract}.
	\begin{equation*}
		K_t = \begin{cases}
			\set{a,b}			
				&\text{ if } t = 0 \\
			\set{a,b,c, ab}		
				&\text{ if } t = 1 \\
			\set{a,b,c, ab,ac,bc}
				&\text{ if } t = 2 \\
			\set{a,b,c,d, ab,ac,bc, abc}
					&\text{ if } t = 3 \\
			\set{a,b,c,d, ab,ac,bc,ad, abc}
					&\text{ if } t \geq 4 \\
		\end{cases}
	\end{equation*}
	Note that the abstract description of $K_t$ can be cumbersome to work with by hand, e.g.\ confirming that each element corresponds to a simplicial complex and that the collection corresponds to a filtration can be become tedious even for a relatively small number of simplices.

	In the illustration below, we describe the filtration $\filt{K}$ by specifying which simplices are added as we increase the index $t \in \nonnegints$.
	In particular, for each $t \in \nonnegints$, the simplices colored in \redtag are the simplices in $K_t$ that are not present in $K_{t-1}$, with $K_{-1}$ interpreted to be $\emptyset$.
	These simplices are also listed in \redtag below the illustration of each $K_t$.
	\begin{center}
		\includegraphics[height=1.1in]{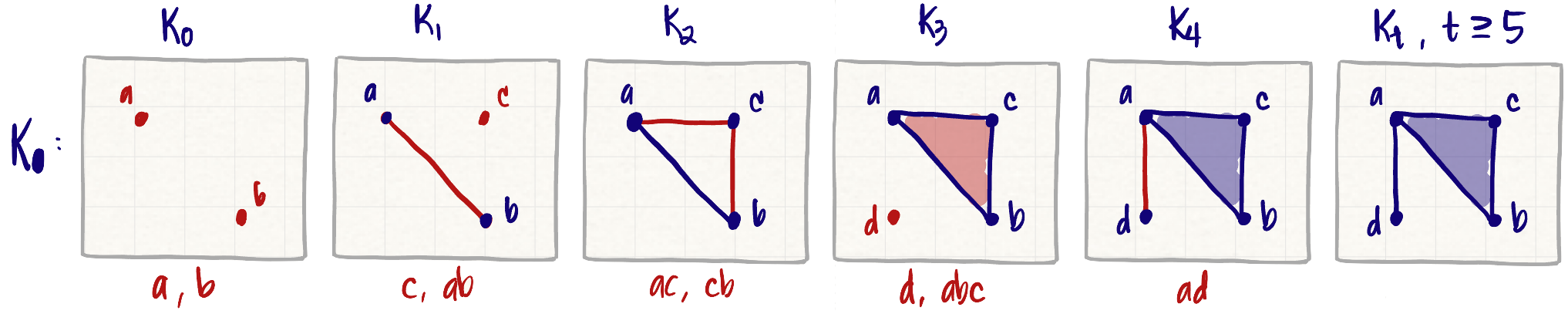}
	\end{center}
	Observe that there are no simplices added to $K_t$ for $t \geq 5$.
\end{example}

We would like to point out that the subcomplex relations on the simplicial complexes of a filtration $\filt{K}$ can be seen as a total order on the set $\set{K_t}_{t \in \nonnegints}$ induced by the total order on $\nonnegints$.
That is, we can define a total order $\leq_c$ on the set of simplicial complexes by defining $K \leq_c L$ if and only if $K$ is a subcomplex of $L$.
Since functors with $\posetN$ as the domain category can generally be represented using a sequence with arrows (see discussion on poset categories in~\fref{appendix:cat-theory}),
we can interpret an arrow $K_t \to K_s$ to refer to the subcomplex relation $K_t \subseteq K_s$, i.e.\ 
\begin{equation*}
	K_0 \Xrightarrow{i^0}
	K_1 \Xrightarrow{i^1}
	K_2 \Xrightarrow{i^2}
	K_3 \Xrightarrow{i^3}
	\cdots
	\qquad\text{ corresponds to }\qquad\quad
	K_0 \,\subseteq 
	K_1 \,\subseteq 
	K_2 \,\subseteq 
	K_3 \,\subseteq 
	\cdots
\end{equation*}
One consequence of this is the use of colloquial language when talking about filtrations. 
For example,
	\textit{going up} a filtration usually implies going from a simplicial complex $K_t$ from some $t \in \nonnegints$ and then considering a simplicial complex $K_s$ at a higher index $t < s$.

Note that in most introductory literature for persistent homology, simplicial filtrations are generally represented using \textit{finite} nested sequences.
That is, filtrations are sometimes described to be finite collections $\set{K_0, K_1, \ldots, K_T}$ of simplicial complexes such that 
\begin{equation*}
	K_0 \,\subseteq\,
	K_1 \,\subseteq\,
	K_2 \,\subseteq\,
	\cdots \,\subseteq\,
	K_T
\end{equation*}
This motivates the following terminology. 
	
\begin{definition}
	Let $K_\bullet$ be a simplicial filtration.
	We say that $\filt{K}$ is \textbf{constant on an interval} $I \subseteq \nonnegints$ if for all $t,s \in I$ with $t \leq s$, 
		$K_t = K_s$.
	We say that $\filt{K}$ is \textbf{finite-type} if $K_t$ is a finite simplicial complex for all $t \in \nonnegints$, and there exists $T \in \nonnegints$ such that $\filt{K}$ is constant on $[T,\infty)$.
\end{definition}

Note that these definitions mimic those of \fref{defn:persmod-constant-on-interval}, which defines the same terms for the case of persistence modules.
Observe that if $\filt{K}$ is a finite-type filtration,
	then there can only be finitely many distinct simplices in $\filt{K}$.
As we will see in \fref{section:construction-of-persistent-homology}, 
	this finite-type condition on filtrations implies that the persistence modules constructed for simplicial persistent homology are also of finite-type.
We provide an example of a finite-type filtration below.

\begin{example}
	The simplicial filtration $\filt{K}$ defined in \fref{ex:baby-filtration} is a finite-type filtration that is constant on $[4,\infty)$.
	In other words, knowing that $\filt{K}$ is constant on $[4,\infty)$,
		we can determine $\filt{K}$ using only the following illustrations of simplicial complexes:
	\begin{center}
		\includegraphics[width=0.65\linewidth]{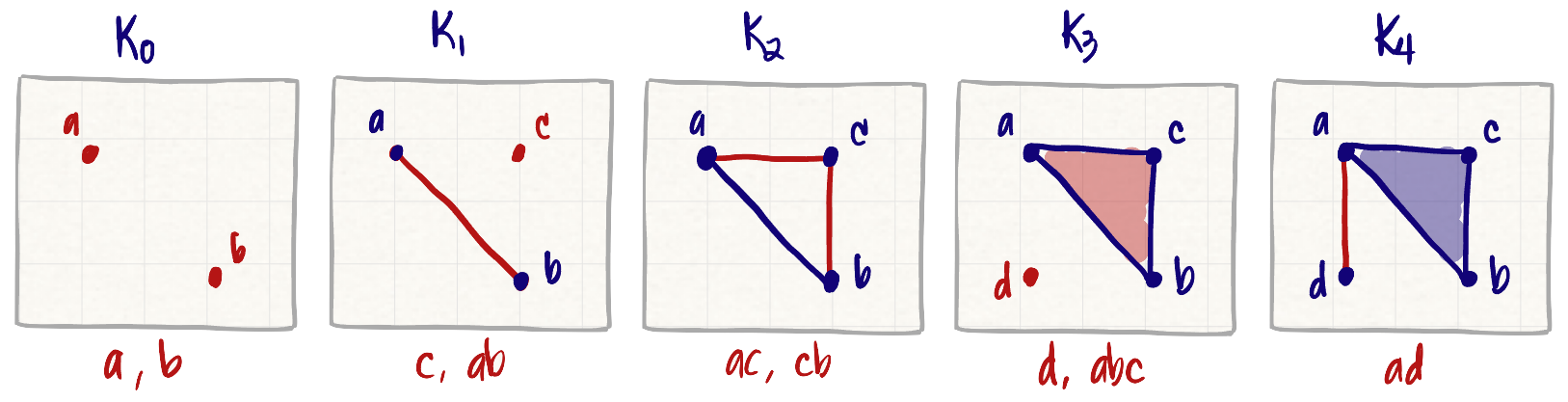}
	\end{center}
	The assumption that $\filt{K}$ is constant on $[4,\infty)$ determines that for all $t \geq 5$, $K_t = K_4$.
\end{example}

In practice, simplicial filtrations are usually constructed to be filtrations on \textit{finite} simplicial complexes.
Relative to \fref{defn:filtration}, we assume that a filtration described this way is constant on $[T,\infty)$.
We provide an example of this below. 

\begin{example}
	Given below is a copy of \cite[Figure 1]{matrixalg:zomorodian},
	which describes a nested sequence of simplicial complexes 
	$K_0 \subseteq K_1 \subseteq K_2 \subseteq \cdots \subseteq K_5$.
	Under each illustration of a simplicial complex, 
		the index $t \in \set{0, 1, \ldots, 5}$ in $K_t$ is denoted on the bottom-left,
		and the list of simplices present in $K_t$ but not on $K_{t-1}$ is on the bottom-right (with $K_{-1} := \emptyset$).
	The $0$-simplices in said list for each $K_t$ are drawn with \redtagged{light\,red} shaded circles,
	the $1$-simplices with \textit{dashed lines},
	and the $2$-simplices with \redtagged{light\,red} shaded triangles.
	\vspace{2pt}
	\begin{center}
		\includegraphics[height=1.2in]{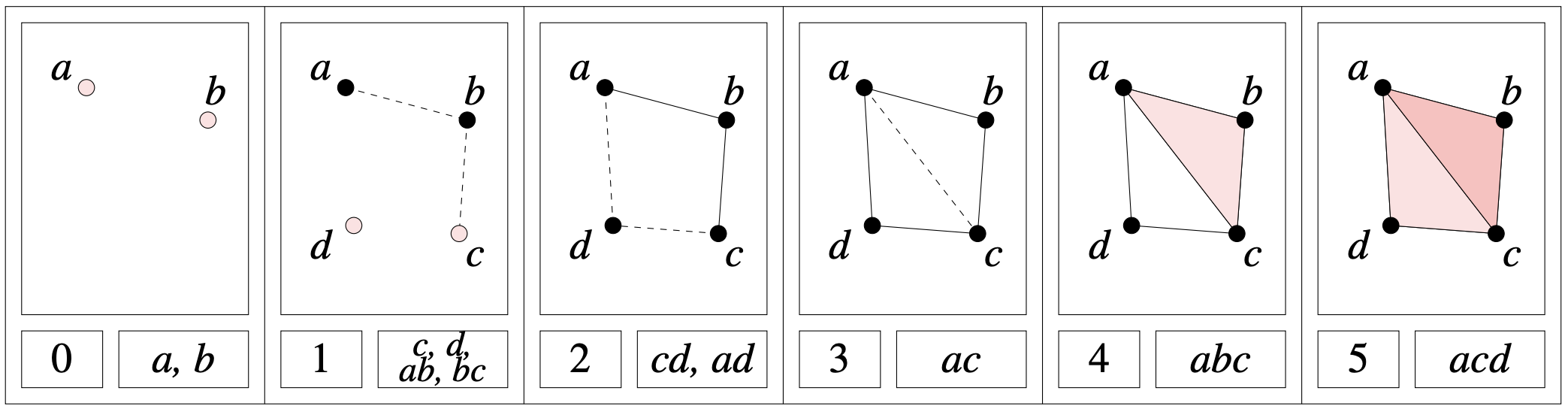}
	\end{center}
	Relative to \fref{defn:filtration}, 
	this determines a filtration $\filt{K}$ on the simplicial complex $K := K_5$ 
	with $K_\bullet(t) = K_t$ illustrated below for all $t \in \nonnegints$:
	\begin{center}
		\baselineCenter{\includegraphics[height=1.1in]{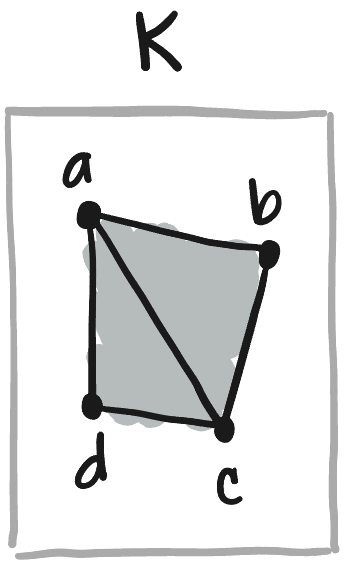}}
		\quad\quad 
		\baselineCenter{\includegraphics[height=1.2in]{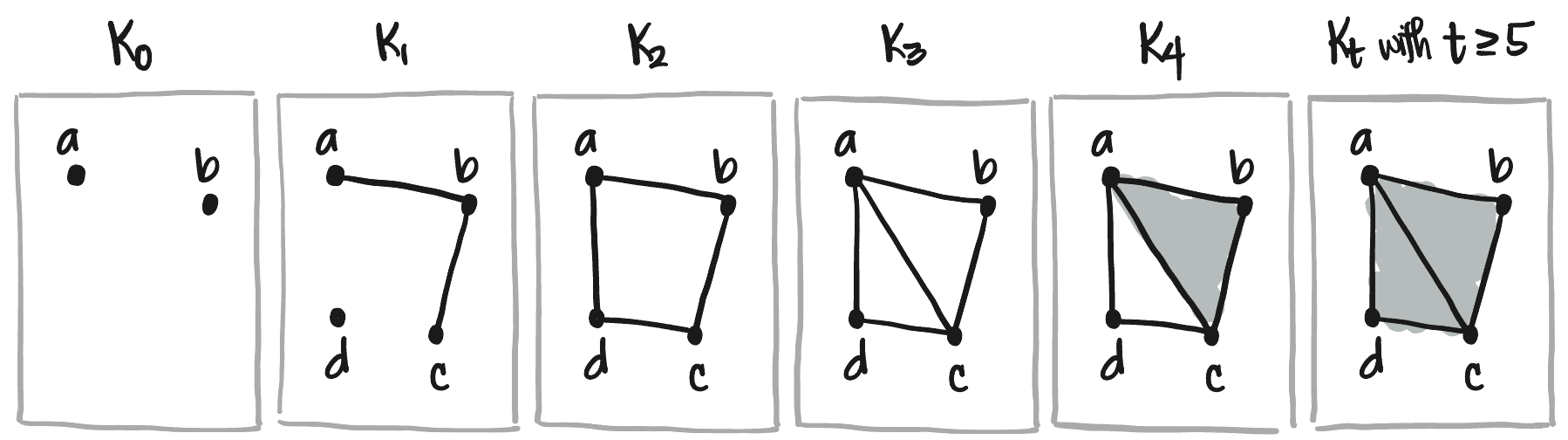}}
	\end{center}
	\vspace{5pt}
	Observe that $\filt{K}$ is constant on $[5,\infty)$.
\end{example}

Filtrations on finite simplicial complexes are necessarily of finite-type.
We state this in more detail below.

\begin{lemma}\label{lemma:filtration-constant}
	Let $\filt{K}$ be a filtration on a simplicial complex $K$.
	If $K$ is a finite simplicial complex,
		then $\filt{K}$ is a finite-type filtration and 
		there exists some $T \in \nonnegints$ such that $\filt{K}$ is constant on $[T,\infty)$.
\end{lemma}
\begin{proof}
	Assume that $K$ is a finite simplicial complex
	and let $\filt{K}$ be some filtration on $K$.
	For all $t \in \nonnegints$,
		$\filt{K}(t) = K_t$ must be a finite simplicial complex since $K_t \subseteq K$ by \fref{defn:filtration}(i) and $K$ is finite by assumption.

	Assume, for the sake of contradiction, that there does not exist $T \in \nonnegints$ such that $\filt{K}$ is constant on $[T,\infty)$.
	Let $r_0 = 0$. 
	For each $t \geq 1$, let $r_t \in \nonnegints$ be such that $r_{t-1} \leq r_t$ and $K_\bullet(r_{t-1}) \neq K_\bullet(r_t)$.
	Note that $r_t$ exists since $\filt{K}$ cannot be constant on $[r_{t-1},\infty)$ by assumption.
	Since $K_\bullet(r_{t-1}) \subseteq K_\bullet(r_t)$ by \fref{defn:filtration}(ii), 
		there exists a simplex $\sigma_t \in K$ such that 
		$\sigma_t \in K_\bullet(r_t)$ and $\sigma_t \not\in K_\bullet(r_{t-1})$.
	This inductive process constructs a collection $\set{\sigma_0, \sigma_1, \sigma_2, \ldots}$ of simplices of $K$ such that each $\sigma_t$ is distinct.
	Since this process can continue indefinitely by assumption of the non-existence of $T \in \nonnegints$,
		we can construct an infinite collection $\set{\sigma_t}_{t=0}^\infty$ of distinct simplices of $K$.
	Since $\bigcup_{t=0}^\infty \sigma_t \subseteq K$,
		$K$ must be infinite.
	This contradicts the assumption that $K$ is finite.
\end{proof}


\HIDE{
	
Next, since it can be proven that addition of exactly one simplex changes the homology of a simplicial complex, i.e.\ an addition of an $n$-simplex $\sigma$ to some simplicial complex $K$ (assuming $K \cup \set{\sigma}$ is also a simplicial complex) does exactly one of two things:
it creates a free component in $H_n(K;R)$ or it destroys a free component in $H_{n-1}(K;R)$.
In the context of persistent homology, filtrations that build simplicial complexes one simplex at a time have properties that can be used to optimize the matrix reduction algorithm.
We provide a definition of this kind of filtration below, adapted from~\cite{ripser}.

\begin{definition}\label{defn:filtration-simplexwise-refinement}
	Let $\filt{K}$ be a filtration on some simplicial complex $K$.
	\begin{enumerate}
		\item We say that $\filt{K}$ is a \textbf{simplexwise filtration} if $K_0 = \set{v}$ for some vertex $v \in K$ and for all $t \in \nonnegints$, $K_{t+1} \setminus K_t = \set{\sigma}$ for some simplex $\sigma \in K$.

		\item A simplex-wise filtration $\filt{S}$
		is called a \textbf{simplexwise refinement} of $\filt{K}$ if there exists an increasing map $r: \nonnegints \to \nonnegints$ such that for all $t \in \nonnegints$, $K_t = S_{r(t)}$.
		In this case, we call $r$ a \textbf{reindexing map} from $\filt{K}$ to $\filt{S}$.
	\end{enumerate}
\end{definition}

Note that simplexwise refinements are generally not unique.
In practice, simplexwise refinements are generated by \textit{expanding} a filtration.
For our purposes, this means building a simplexwise filtration inductively in increasing index $t \in \nonnegints$.
In this process, if there exists $t$ such that $\card(K_{t+1} \setminus K_t) \geq 2$, i.e.\ more than one simplex is appended to $K_t$,
there is a choice of 
which simplices in $K_{t+1} \setminus K_t$ to append first in the constructed simplexwise filtration.
We provide an example of this below.

\begin{example}
	Let $\filt{K}$ be the filtration on $K$ as defined in~\fref{ex:baby-filtration}. An illustration of $\filt{K}$ is copied below:
	\begin{center}
		\includegraphics[width=0.7\linewidth]{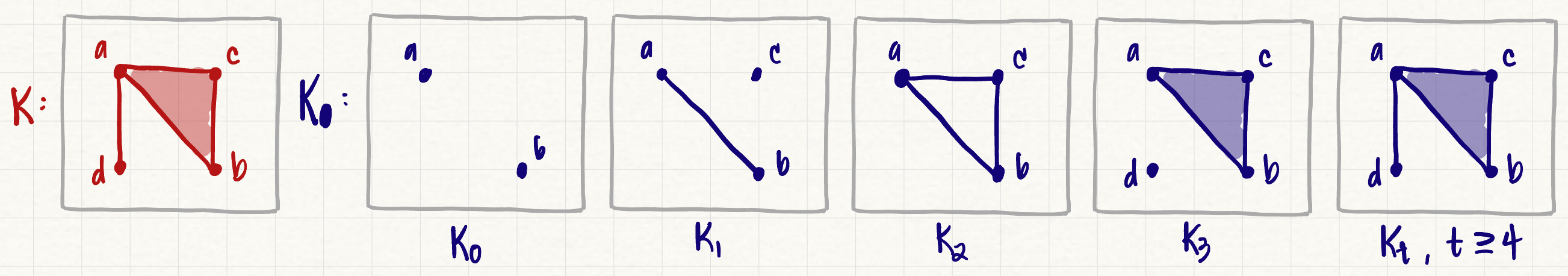}
	\end{center}
	We have constructed a simplexwise refinement $\filt{X}$ of $\filt{K}$, illustrated below:
	\begin{center}
		\includegraphics[width=0.6\linewidth]{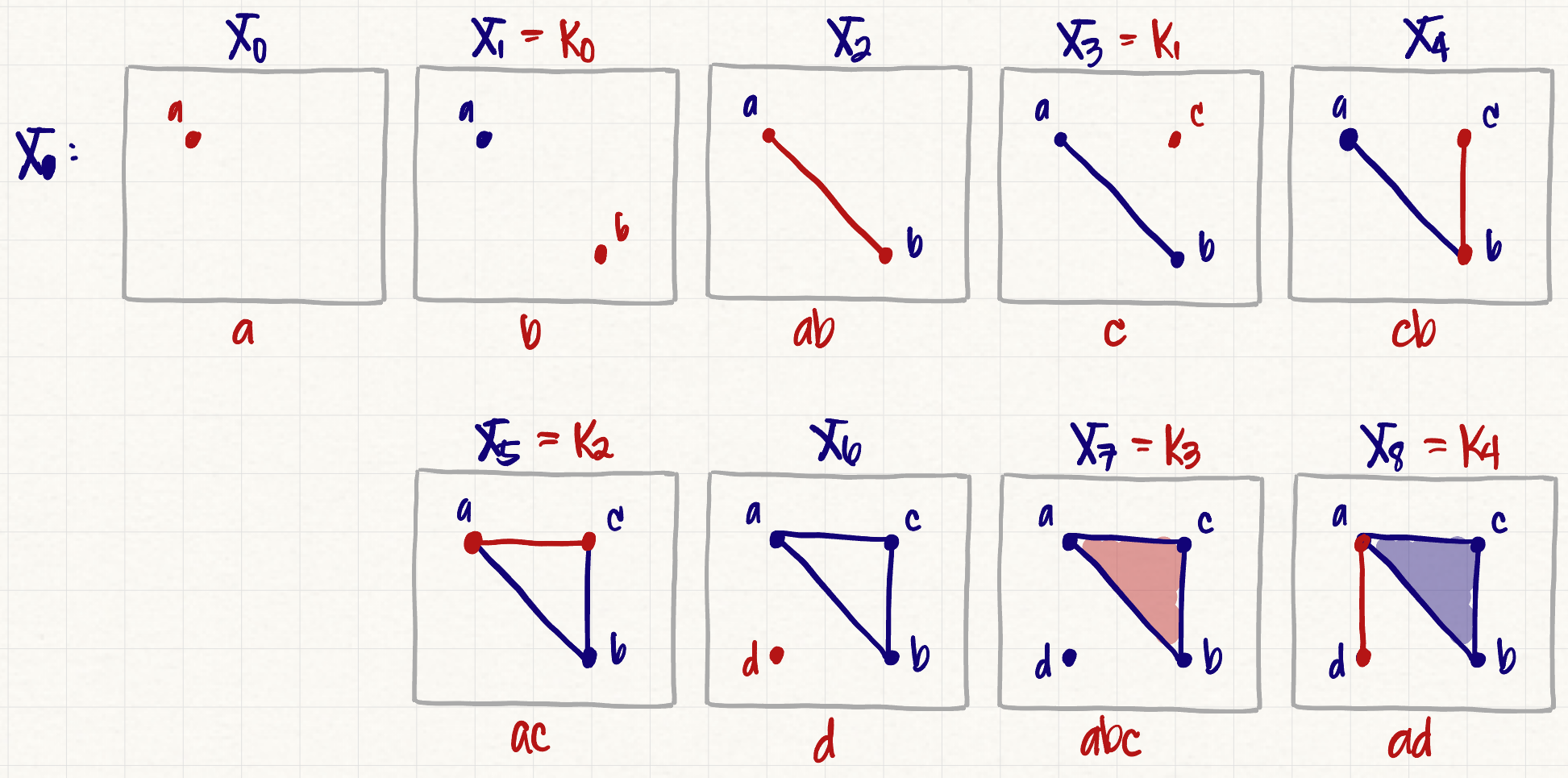}
	\end{center}
	Note that the simplices colored in \textit{red} are the simplices append at each index of the filtration.
	
	To start, $\filt{X}$ is constructed by defining an ordering of the simplices of $K_0$ (generally in increasing dimension to ensure that simplicial complexes are well-defined). In this case, we assigned the $0$-simplex $a$ to $X_0$ and assigned the $0$-simplex $b$ to the next space $X_1$. Next, we consider $K_1 \setminus K_0 = \set{c,ab}$. 
	We arbitrarily chose to assign $ab$ first to $X_2$ and $c$ to $X_3$.
	Note that defining $X_2 = X_1 \cup \set{c}$ and $X_3 = X_2 \cup \set{ab}$ would have worked as well.
	
	This process continues until all simplices of $K$ are added.
	Observe that with this method, any constructed simplex-wise refinement must be constant on $[\card(K)-1, \infty)$.
	Then, we can define the reindexing map $r: \nonnegints \to \nonnegints$ from $\filt{K}$ by $r(t) = r(\card(K_t) - 1)$ for all $t \in \nonnegints$.
\end{example} }

As a final remark,
we want to point out that, like with persistence modules, the definition of a filtration can be generalized to allow for any poset category to be the domain category. That is, we can define filtrations to be functors of the form $\poset(P, \leq) \to \catsimp$ for any poset $(P, \leq)$.

\HIDE{In most applications, $R = \reals$ or $R$ is some finite set in $\reals$.
For this expository paper, we justify the restriction of the indexing set to $\nonnegints$ since the focus of this paper is on the theory behind the matrix reduction algorithm for persistent homology.
In particular, filtrations of the form $\poset(R, \leq) \to \catsimp$ have to be represented as filtrations of the $\posetN \to \catsimp$, as defined \fref{defn:filtration}, in order for the matrix reduction algorithm to apply.
Note that this representation has to be, in some sense, lossless and this is generally a strong condition.

}\clearpage

\section{The Persistent Homology of Filtrations}
\label{section:persistent-homology}

In this section, we define the persistent homology of filtrations in terms of persistence modules
and discuss the interval decompositions of persistence modules in the context of persistent homology.
We start with a definition of persistent homology of filtrations, adapted from \cite{ripser}.

\begin{statement}{Definition}\label{defn:pershom}
	\label{defn:n-persistent-homology}
	Let $n \in \ints$.
	The $n$\th \textbf{persistent homology module} $H_n(\filt{K}; \field)$ of a filtration $\filt{K}$ \textbf{with coefficients in $\field$}
	is the persistence module $H_n(\filt{K}; \field): \posetN \to \catvectspace$ given by the following functor composition:
	\begin{equation*}
		H_n(\filt{K}; \field) := H_n(-; \field) \circ \filt{K} 
	\end{equation*}
	where $H_n(-; \field): \catsimp \to \catvectspace$ denotes the $n$\th simplicial homology functor with $\field$ coefficients.
	The $n$\th \textbf{persistence barcode $\barcode_n(\filt{K}; \field)$} of $\filt{K}$ \textbf{with coefficients in $\field$} is
	the persistence barcode of $H_n(\filt{K;\field})$ as a persistence module, i.e.\ 
	$\barcode_n(\filt{K};\field) := \barcode(H_n(\filt{K};\field))$.
\end{statement}

We have some remarks regarding the notation and terminology of the vector spaces and structure maps of $H_n(K_\bullet;\field)$, taken as a persistence module, relative to our characterization of persistence modules by \fref{defn:persmod}.
\begin{enumerate}
	\item 
	Observe that $H_n(K_\bullet;\field) = H_n(K_t;\field)$ for all $t \in \nonnegints$. 
	This is consistent with the convention of replacing the bullet $(\bullet)$ of a persistence module $V_\bullet = (V_\bullet,\alpha_\bullet)$ to denote the vector space $V_t$ at index $t \in \nonnegints$.
	We may also call the index $t \in \nonnegints$, in the context of $H_n(K_\bullet;\field)$ as a persistent homology module, as the \textit{scale} or \textit{parameter} of $H_n(K_t;\field)$.

	\item 
	The structure maps of a filtration $K_\bullet$, as stated in \fref{defn:filtration}, are denoted as $i^{s,t}: K_t \to K_s$.
	Following the notation for the induced maps on homology,
		the structure maps of $H_n(K_\bullet;\field)$ are denoted as $i^{s,t}_\ast: H_n(K_t;\field) \to H_n(K_s;\field)$.
	Relative to the notation $V_\bullet = (V_\bullet, \alpha_\bullet)$ for an arbitrary persistence module,
		$H_n(\filt{K}; \field) = \bigl( 
			H_n(\filt{K}; \field) , i^\sbullet_\ast
		\bigr)$.
\end{enumerate}
Since we may be dealing with homology classes in $H_n(K_t;\field)$ of $H_n(K_\bullet;\field)$ across different indices $t \in \nonnegints$, we identify some alternative notation below.

\begin{bigremark}
	We modify the coset notation for the homology classes to include the index $t \in \nonnegints$ as a subscript when discussed in the context of $H_n(K_\bullet;\field)$,
	i.e.\ we write $[\sigma]_t$ to refer to the homology class $[\sigma] \in H_n(K_t;\field)$ as an element of the vector space of $H_n(K_\bullet;\field)$ specifically at index $t \in \nonnegints$.

	Note that $\sigma$ in $[\sigma]_t \in H_n(K_t;\field)$
	denotes a cycle representative $\sigma \in \ker(\boundary_n^{\w t})$, with 
	$\boundary_n^{\w t}: C_n(K_t;\field) \to C_{n-1}(K_t;\field)$ referring to the $n$\th simplicial boundary map of $K_t$.
\end{bigremark}

Since we defined persistent homology modules using functor composition, 
	we can interpret the persistent homology module $H_n(\filt{K}; \field)$ to be calculated 
	using simplicial homology first at each index $t \in \nonnegints$,
	and 
second by assembling the resulting homology groups into a persistence module, 
i.e.\ we have two separate operations, illustrated as follows:
\vspace{-5pt}
\begin{equation*}
\begin{tikzcd}[column sep=7em]
	\posetN 
		\arrow[r, "\text{filtration } K_\bullet"]
		\arrow[r, mapsto, yshift=-1.8em]
	&[0.5em] \catsimp 
		\arrow[r, "\text{simplicial homology } H_n(-;\field)"]
		\arrow[r, mapsto, yshift=-1.8em]
	&[8em] \quad\catvectspace\quad
	\\[-1.8em]
	t
	& K_t 
	& H_n(K_t;\field)
	\\[-2.1em] 
	\mathclap{\text{\footnotesize(scale or index)} }
	& \mathclap{\text{\footnotesize(simplicial complex)}} 
	& \mathclap{\text{\footnotesize(homology group)}}
\end{tikzcd}\vspace{-5pt}
\end{equation*} 
As such, the terminology for simplicial complexes and simplicial homology also extends to the homology groups of a persistent homology module at some specific index $t \in \nonnegints$.

Below, we provide an example of a persistent homology module of a filtration.

\begin{example}\label{ex:pershom-one}
	Let the simplicial complex $K$ and the filtration $\filt{K}$ on $K$ be as illustrated below:\vspace{0pt}
	\begin{center}
		\baselineCenter{\includegraphics[height=1.05in]{zomfig1/zom-simp.png}}
		\quad
		\baselineCenter{\includegraphics[height=1.1in]{zomfig1/zom1-cleaned.png}}
	\end{center}
	\vspace{3pt}
	We consider the $n$\th persistent homology module $H_n(\filt{K};\rationals)$ of $K_\bullet$ with $\rationals$ coefficients in all dimensions $n \in \nonnegints$.

	\vspace{0.5\baselineskip}\noindent
	\textbf{Part A. Pointwise Calculation of Homology Groups with Rational Coefficients}

	We consider 
		the chain groups $C_n(K_t;\rationals)$ 
		and boundary maps $\boundary_n^{\,t}: C_n(K_t;\rationals) \to C_{n-1}(K_t;\rationals)$ 
	of $K_t$ at each scale $t \in \nonnegints$ separately.
	Orient each simplicial complex $K_t$ with the vertex order $\Vertex(K) = (a,b,c,d)$ restricted on $\Vertex(K_t)$.
	Illustrated below is the induced orientation on the $1$-simplices of $K$ and those of $K_t$ for all $t \in \nonnegints$:
	\begin{center}
		\baselineCenter{\includegraphics[height=1.1in]{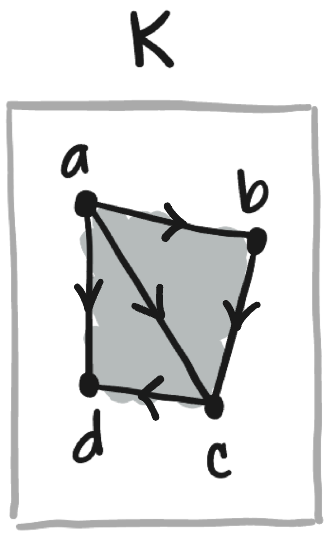}}
		\quad
		\baselineCenter{\includegraphics[height=1.1in]{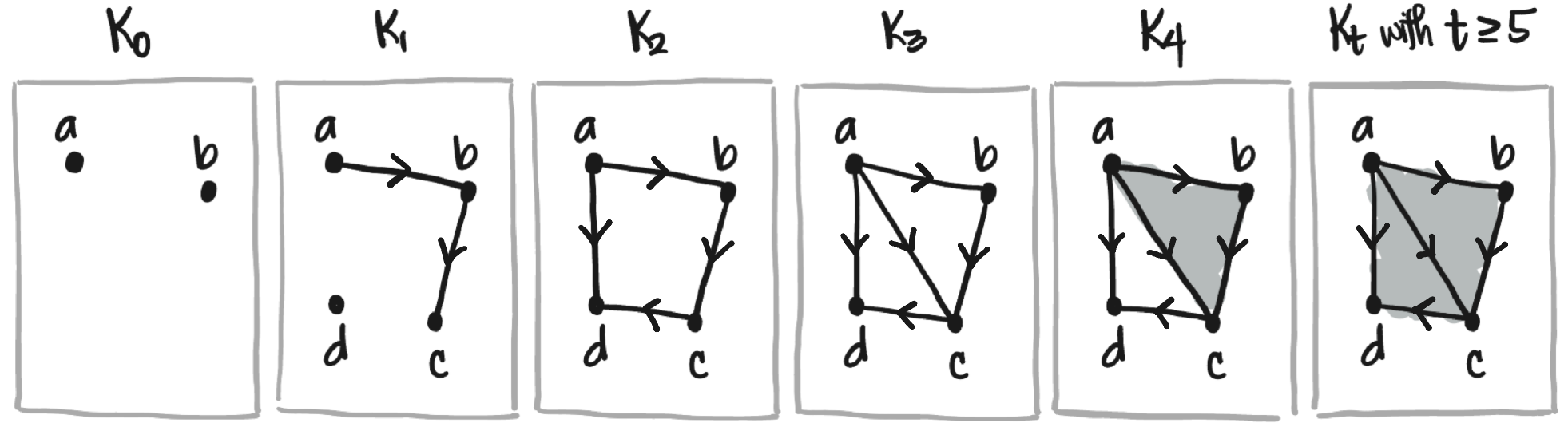}}
	\end{center}
	For all $n \in \nonnegints$,
		the $n$\th homology group $H_n(K_t;\rationals)$ of $K_t$ with $\rationals$ coefficients is then calculated as usual with
	\begin{equation*}
		H_n(K_t;\rationals) = \ker\bigl(
			\boundary_n^{\w t}
		\bigr) \bigmod \im\bigl(
			\boundary_{n+1}^{\,t}
		\bigr)
	\end{equation*} 
	Observe that $H_n(K_t;\rationals)$ is a $\rationals$-vector space for all $n \in \nonnegints$.
	Since $C_n(K_t;\rationals) = 0$ for all $n \geq 3$,
		$\ker(\boundary_n^{\w t}) = 0$ and 
		$H_n(K_t;\rationals) = 0$ for all $n \geq 2$.
	Then, $H_0(K_t;\rationals)$ and $H_1(K_t;\rationals)$ for all $t \in \nonnegints$ are as follows:
	\begin{equation*}
		H_0(K_t; \rationals) \cong \begin{cases}
			\rationals\ket{[a]_0, [b]_0} 		
				&\text{ if } t = 0 		\\
			\rationals\ket{[a]_1, [d]_1} 		
				&\text{ if } t = 1 		\\
			\rationals\ket{[a]_t} 		
				&\text{ if } t \geq 2
		\end{cases}
		\quad\text{ and }\quad 
		H_1(K_t; \rationals) = \begin{cases}
		0 		
			&\text{ if } t = 0,1 		\\
		\rationals\ket{[ab + bc + cd - ad]_2} 		
			&\text{ if } t = 2 		\\
		\rationals\ket{[ab + bc - ac]_3, [ac + cd - ad]_3} 		
			&\text{ if } t = 3 \\
		\rationals\ket{[ac + cd - ad]_4} 
			&\text{ if } t = 4 \\
		0 &\text{ if } t \geq 5
		\end{cases}
	\end{equation*}

	\noindent
	\textbf{Part B. Interpretation of Homology Groups with Rational Coefficients}

	The interpretation of the homology classes of $H_0(K_t;\ints)$ as corresponding to path components, see \cite[Proposition 2.6, p109]{algtopo:hatcher}, also applies to homology in $\rationals$ coefficients.
	Illustrated below are the path components of $K_t$ relating to basis elements of $H_0(K_t;\rationals)$ for all $t \in \nonnegints$.
	\vspace{5pt}
	\begin{center}
		\begin{minipage}[t]{0.63\linewidth}\null\vspace{-0.7\baselineskip}
			\includegraphics[width=\linewidth]{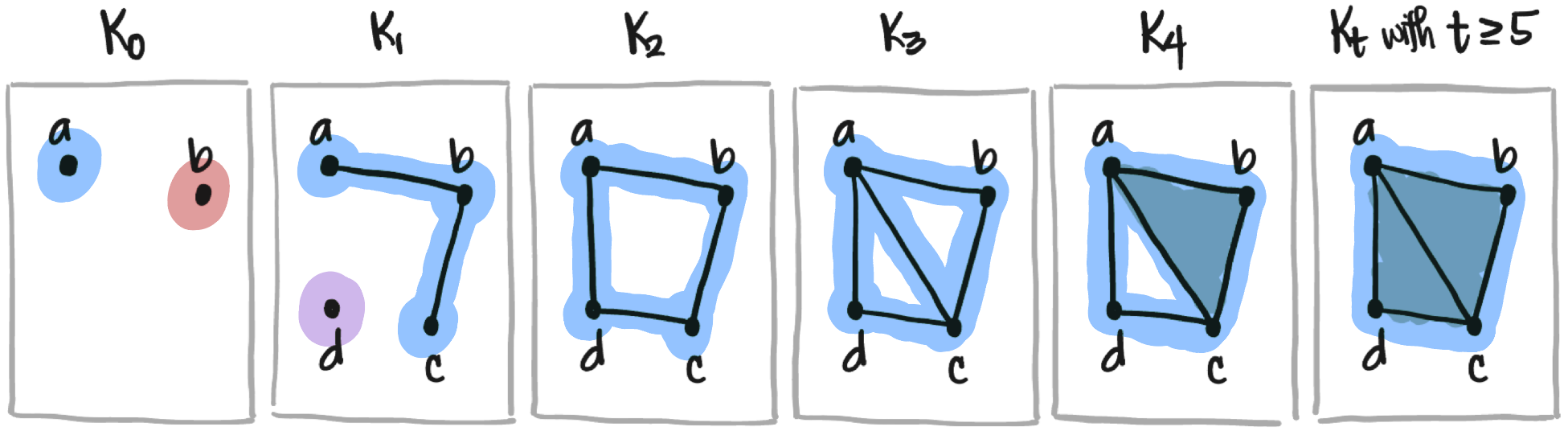}\null
		\end{minipage}\hspace{0.01\linewidth}\begin{minipage}[t]{0.34\linewidth}\raggedright 
			Color Scheme for Path Components: \\[2pt]
			\begin{enumerate}[topsep=0pt, leftmargin=55pt, labelsep=0pt]
				\item[\bluetag: ] 
				that of $[a]_t \in H_0(K_t;\rationals)$ at all scales $t \in \nonnegints$
	
				\item[\redtag: ]
				that of $[b]_0 \in H_0(K_0;\rationals)$ only at scale $t=0$
	
				\item[\purpletag: ]
				that of $[d]_1 \in H_0(K_1;\rationals)$ only at scale $t=1$
			\end{enumerate}
		\end{minipage}%
	\end{center}
	\vspace{5pt}

	The cycle representatives of $H_1(K_t;\ints)$ correspond to directed loops in $K_t$, consisting of $\ints$-multiples of oriented $1$-simplices.
	The same interpretation also applies to $H_1(K_t;\rationals)$ since we can scale any cycle representative such that it consists of $\ints$-multiples of oriented $1$-simplices.
	Illustrated below are the loops in $K_t$ with $t=2,3,4$ corresponding to homology classes in $H_1(K_t;\rationals)$:
	\vspace{5pt}
	\begin{center}
		\begin{minipage}[t]{0.32\linewidth}\null\vspace{-0.8\baselineskip}
			\includegraphics[width=\linewidth]{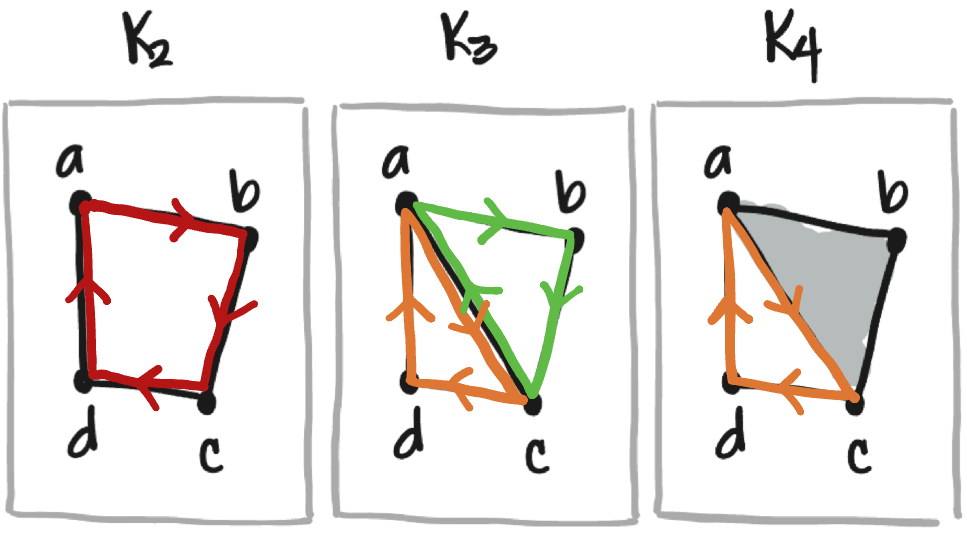}\null
		\end{minipage}\hspace{0.01\linewidth}\begin{minipage}[t]{0.65\linewidth}\raggedright 
			Color Scheme for Directed Loops: \\[2pt]
			\begin{enumerate}[topsep=0pt, leftmargin=47pt, labelsep=0pt, itemsep=3pt] 
				\item[\redtag: ] 
				$\redmath{\alpha := ab + bc + cd - ad}$, representing $[\alpha]_t \in H_1(K_t;\rationals)$ for $t \geq 2$. 
				Note that $[\alpha]_t = 0$ if $t \geq 5$.

				\item[\greentag: ] 
				$\greenmath{\tau := ab + bc - ac}$,
				representing $[\tau]_t \in H_1(K_t;\rationals)$ for $t \geq 3$. 
				Note that $[\tau]_t = 0$ if $t \geq 4$.

				\item[\orangetag: ]
				$\orangemath{\sigma := ac + cd - ad}$,
				representing $[\sigma]_t \in H_1(K_t;\rationals)$ for $t \geq 3$.
				Note that $[\sigma]_t = 0$ if $t \geq 5$
				and that $[\alpha]_3 = [\tau]_3 + [\sigma]_3$.
			\end{enumerate}
		\end{minipage}%
	\end{center}
	\vspace{5pt}
	Negative multiples of oriented $1$-simplices denote a reversal in direction. For example, $ad \in C_1(K_2;\ints)$ is a line from vertex $a$ to vertex $d$ but $-ad \in C_1(K_2;\ints)$ is a line from $d$ to $a$.
	The addition of an oriented $1$-simplex to an oppositely oriented copy of itself can be interpreted to be an annihilation of said simplex, e.g.\ $ad + (-ad) = 0$.
	
	\vspace{0.5\baselineskip}\noindent 
	\textbf{Part C. Formation of Persistent Homology Modules}

	Let $n \in \nonnegints$.
	The $n$\th persistent homology module $H_n(K_\bullet;\rationals)$ 
		has vector spaces given by $H_n(K_\bullet;\rationals)(t) = H_n(K_t;\rationals)$ for all $t \in \nonnegints$.
	The structure maps of $H_n(K_\bullet;\rationals)$ are the maps on homology induced by inclusions $i^{s,t}: K_t \to K_s$ for all $t,s \in \nonnegints$ with $t \leq s$.
	For now, it suffices to know that these structure maps exist.
		
	Observe that $H_n(\filt{K};\rationals)$ is the trivial persistence module over $\rationals$ for $n \geq 2$ since $H_n(K_t;\rationals)=0$ for all $t \in \nonnegints$.
	We talk more about $H_0(K_\bullet;\rationals)$ and $H_1(K_\bullet;\rationals)$ in the following examples in this section.		
\end{example}


\begin{bigremark}\label{remark:persistent-hom-why-intsprime}
	The preferred coefficient field for persistent homology is $\field = \ints_p$ with prime $p$.
	The majority of the examples in this paper use $\field = \rationals$ for the coefficient field for convenience.
\end{bigremark} 

\HIDE{
The reason for this choice is that the focus of this expository paper is on the theory behind the matrix reduction algorithm for persistent homology and said algorithm is typically done using $\ints_p$ coefficients (as opposed to $\reals$ or $\rationals$ coefficients).
We attribute this convention to the well-known \textit{instability} of algorithms on matrices over $\reals$, which is rooted in the fact that real numbers cannot be stored exactly by computers.
Here, we use the term \textit{instability} to refer to how small errors from arithmetic operations can propagate and lead to massive errors in the output.
For example, according to the \href{https://docs.python.org/3/tutorial/floatingpoint.html}{Official Python3 Documentation}, 
the numbers $1$ and $10$ have exact binary floating point representations
but the number $0.1$ does not, assuming that a finite number of bits are used,
since the decimal expansion of $0.1$ in base-2 is infinite.
In this case, $0.1$ needs to be approximated to the nearest floating point number and therefore, $0.1$ cannot be exactly stored.
Because of this, we get 
results such as $\texttt{1/10} \neq \texttt{0.1}$ and $\texttt{1/10 - 0.1} \neq \texttt{0}$, relative to equality between floating point numbers.
Issues like this leads to the well-known instability of Gaussian elimination of matrices over $\reals$.

	As a sidenote, while elements $\frac{p}{q} \in \rationals$ can be stored exactly using a pair $(p,q)$ of integers $p$ and $q$, the finite nature of storage means not all numbers in $\ints$ can be represented exactly. For example, the usual unsigned 16-bit representation for integers can only store integers in the interval $[0, 2^{16}-1]$. 
	A rational number such as $1 / (2^{32} + 4)$ cannot be stored exactly in this manner.
	Consequently, this representation is not generally closed under addition and multiplication and may cause some issues with matrix reduction.
}
	
The reason for $\ints_p$ being the preferred coefficient field is that numbers in $\ints_p$ can be represented exactly using a finite number of bits. 
For example, the Ripser package described in \cite{ripser} stores numbers in $\ints_p$ as unsigned 16-bit integers and allows for calculation of persistent homology in $\ints_p$ coefficients with $p \leq 2^{16} - 1$. 
Furthermore, a number in $\ints_2$ can be represented using one bit and arithmetic operations in $\ints_2$ can be done as bit operations, which are extremely fast computationally.
This explains why Cohen-Steiner, et al.\ in~\cite{matrixalg:cohen-steiner} and other authors (particularly those in the field of computer science) often define persistent homology with coefficients in $\field = \ints_2$, i.e.\ their persistence modules are of the form $\posetN \to \catvectspace[\ints_2]$.

However, operations in $\ints_p$ can be somewhat cumbersome to do by hand, which became an issue when we were creating examples for this expository paper. 
As a compromise, we use the field $\field = \rationals$ for most of our examples since the relation between
homology groups over $\rationals$ and those of $\ints_p$ is relatively straightforward.
For reference, see \cite[Corollary 3A.6]{algtopo:hatcher}.

\spacer

The persistence barcode of a persistence module, as discussed in \fref{chapter:persistence-theory}, characterizes the ranks of the vector spaces and the ranks of the structure maps of the persistence module. The same applies for persistent homology modules, assuming interval decompositions and persistence barcodes exist.
We state an existence result below.

\begin{proposition}\label{prop:filtrations-finite-type-has-interval}
	Let $\filt{K}$ be a filtration on a finite simplicial complex $K$ and fix a field $\field$.
	For each $n \in \nonnegints$,
		the $n$\th persistent homology module $H_n(\filt{K};\field)$ is of finite-type and therefore admits an interval decomposition.
\end{proposition}
\begin{proof}
	By \fref{lemma:filtration-constant}, 
	$K_t$ is a finite simplicial complex for all $t \in \nonnegints$ and 
	there exists $T \in \nonnegints$ such that $\filt{K}$ is constant on $[T,\infty)$.
	For $t \in [0,T]$,
		$H_n(K_t;\field)$ must be finitely generated since 
		$K_t$ is a finite simplicial complex.
	For all $t \in [T,\infty)$,
		$H_n(K_T;\field) = H_n(K_t;\field)$ since $K_T = K_t$.
	Therefore, $H_n(\filt{K};\field)$ is a finite-type persistence module that is constant on $[T,\infty)$.
	By \fref{prop:finite-type-has-interval},
		$H_n(K_\bullet;\field)$ admits an interval decomposition.
\end{proof}

The calculation of this interval decomposition is the goal of the matrix reduction algorithm for persistent homology.
Below, we identify terminology used in \cite{ripser} and other persistent homology literature involving these interval decompositions.

\begin{statement}{Definition}
	Let $\filt{K}$ be a filtration on a finite simplicial complex $K$.
	Let the persistence barcode $\barcode_n(\filt{K};\field)$ of the $n$\th persistent homology module $H_n(\filt{K};\field)$ of $K_\bullet$ with coefficients in a field $\field$
	be given as follows:
	\begin{equation*}
		\barcode_n(\filt{K}; \field)
		= 
		\barcode\bigl( H_n(\filt{K}; \field) \bigr)
		= 
		\Big\{ 
			[b_j, d_j)
		\Big\}_{j=1}^{r}
		\cup 
		\Big\{ 
			[a_i, \infty)
		\Big\}_{i=1}^{m}
	\end{equation*}
	where $\set{[b_j, d_j)}_{j=1}^r$ and $\set{[a_i, \infty)}_{i=1}^m$ both consist of intervals in $\nonnegints$ with $b_j, d_j, a_i \in \nonnegints$.
	\begin{enumerate}
		\item 
		The pair $(b_j, d_j)$ of indices from the interval $[b_j, d_j)$
		is called an \textbf{index persistence pair} of $\filt{K}$.
		In this case, $b_j$ is called a \textbf{birth index} and $d_j$ is called its corresponding \textbf{death index}.

		\item 
		The index $a_i$ in the interval $[a_i, \infty)$ is called an 
		\textbf{essential birth index} or \textbf{essential index} of $\filt{K}$.
	\end{enumerate}
\end{statement}

By \fref{thm:uniqueness-interval-mods},
	persistence barcodes of persistence modules are unique up to persistence isomorphism.
Therefore, the collection of index persistence pairs and essential indices from $H_n(\filt{K};\field)$ can be considered an invariant of the isomorphism type of persistence modules.

\spacer

As briefly discussed in \fref{section:filtrations}, 
it seems to be common practice to assume the index $t \in \nonnegints$ represents some time value.
Roughly speaking, 
	the terms \textit{birth index} and \textit{death index} correlate to how homology classes are created (alternatively, born) or become trivial (alternatively, die or get destroyed) 
	as we increase the time parameter $t \in \nonnegints$.
The term \textit{persistence} then corresponds to how long a homology class that is born at some index $t_B \in \nonnegints$ remains non-trivial (alternatively, lives) as we increase the index $t \geq t_B$.

\HIDE{
	We typically talk about birth, death, and persistence relative to the cycle representatives of homology classes.
	We state a relevant corollary below.

	\begin{corollary}
		Let $\filt{K}$ be a simplicial filtration and fix a field $\field$.
		For each index $t \in \nonnegints$
		and for all $n \in \ints$,
			let $C_n(K_t;\field)$ and $\boundary_n^{\,t}: C_n(K_t;\field) \to C_{n-1}(K_t;\field)$ denote the $n$\th chain group and $n$\th boundary map on $K_t$.
		Then, 
		for all $t,s \in \nonnegints$ with $t \leq s$,
		\begin{equation*}\tag{i}
			i^{s,t}\Bigl(
				\ker\bigl( \boundary_n^{\,t} \bigr) 
			\Bigr)
			\subseteq 
			\ker\bigl( \boundary_n^{\,s} \bigr)
		\end{equation*}
		Therefore, the structure map $i^{s,t}_\ast: H_n(K_t;\field) \to H_n(K_s;\field)$ can be calculated at the level of coset representatives, i.e.\ for all 
		homology classes $[\sigma]_t \in H_n(K_t;\rationals)$ in $K_t$ with cycle representative $\sigma \in \ker(\boundary_n^{\,t})$,
		\begin{equation*}\tag{ii}
			i^{s,t}(\sigma) \in \ker(\boundary_n^{\,s})
			\quad\text{ and }\quad
			i^{s,t}_\ast\bigl(
				[\sigma]_t
			\bigr) = 
			\bigl[
				i^{s,t}(\sigma)
			\bigr]_s \in H_n(K_s;\rationals)
		\end{equation*}
	\end{corollary}
	\begin{proof}
		For (i),
		apply \fref{prop:induced-maps-on-chains-respect-boundary} to the inclusion $i^{s,t}: K_t \to K_s$.
		For (ii), the map $i^{s,t}: H_n(K_t;\field) \to H_n(K_s;\field)$ is the linear map induced by the cokernel on 
		$i^{s,t}_\hash: C_n(K_t;\field) \to C_n(K_s;\field)$.
	\end{proof}
}

Then, the phrase ``{as we increase the time parameter}'' refers to the application of the structure map $i^{s,t}_\ast: H_n(K_t;\field) \to H_n(K_s;\field)$.
We can also colloquially say that we are going up the filtration $K_\bullet$ from time $t$ to time $s \geq t$.
To avoid ambiguity, we provide definitions for the terms such as \textit{lives}, \textit{dies}, and \textit{persists} below.

\begin{statement}{Definition}\label{defn:pershom-intuition}
	Let $\filt{K}$ be a simplicial filtration.
	Let $[\sigma]_t \in H_n(K_t;\field)$ be a non-trivial homology class with cycle representative 
		$\sigma \in \ker(\boundary_n^{\,t}) \subseteq C_n(K_t;\field)$.
	Let $s \in \nonnegints$ such that $t \leq s$.
	\begin{enumerate}
		\item 
		We say that $[\sigma]_t$ \textbf{lives} or \textbf{persists in $K_s$} if 
		$i^{s,t}_\ast([\sigma]_t) \neq 0$.
		
		\item 		
		If $[\tau]_t \in H_n(K_t;\field)$ is a non-trivial homology class such that $[\tau]_t \neq [\sigma]_t$
		and $i_\ast^{s,t}([\sigma]_t) = i_\ast^{s,t}([\tau]_t) \neq 0$,
		we say that $[\tau]_t$ and $[\sigma]_t$ \textbf{merge in $K_s$}.

		\item 
		We say that $[\sigma]_t$ is \textbf{born in $K_B$} or \textbf{at index $B \in \nonnegints$} if 
		$$
			B = \max\Bigl\{
				b \in \nonnegints 
				: b \leq t 
				\text{ and }
					[\sigma]_t \not\in \image\bigl( i_\ast^{t,b-1} \bigr)
			\Bigr\}
		$$
		where $\image\bigl(i_\ast^{t, -1}\bigr)$ is taken to be the trivial vector space.

		\item 
		We say that $[\sigma]_t$ \textbf{dies} or \textbf{is destroyed in $K_D$} or \textbf{at index $D \in \nonnegints$} if 
		$$
			D = \min\Bigl\{
				d \in \nonnegints 
				: t \leq d \text{ and }
				i_\ast^{d,t}([\sigma]_t) = 0
			\Bigr\}
		$$
		If such a minimum does not exist, then 
		we say that $[\sigma]_t$ \textbf{does not die in the filtration}.
		Note that if $[\sigma]_t$ is destroyed in $K_D$, then $i^{s,t}([\sigma]_t) 
		= \bigl( i^{s,D} \circ i^{D,t} \bigr)([\sigma]_t)
		= i^{s,D}(0)
		= 0$ for any $s \geq D$.

		\item 
		We say that $[\sigma]_t$ has \textbf{persistence} or \textbf{lifespan} $D-B \in \nonnegints$ if $[\sigma]_t$ is born at index $B \in \nonnegints$ and dies at index $D \in \nonnegints$.
		If $[\sigma]_t$ does not die in $\filt{K}$, then persistence is interpreted to be $\infty$.
	\end{enumerate}
\end{statement}

We provide an example of these terms in use below, relative to the $0$\th persistent homology module.

\begin{example}
	Let $\filt{K}$ be as given  in \fref{ex:pershom-one}.
	For convenience, we copied the illustration involving $H_0(K_t;\rationals)$ and the path components of $K_t$ at each $t \in \nonnegints$ below:
	\vspace{5pt}
	\begin{center}
		\begin{minipage}[t]{0.63\linewidth}\null\vspace{-0.7\baselineskip}
			\includegraphics[width=\linewidth]{zomfig1/zom1-path-components.png}\null
		\end{minipage}\hspace{0.01\linewidth}\begin{minipage}[t]{0.34\linewidth}\raggedright 
			Color Scheme for Path Components: \\[2pt]
			\begin{enumerate}[topsep=0pt, leftmargin=55pt, labelsep=0pt]
				\item[\bluetag: ] 
				that of $[a]_t \in H_0(K_t;\rationals)$ at all scales $t \in \nonnegints$
	
				\item[\redtag: ]
				that of $[b]_0 \in H_0(K_0;\rationals)$ only at scale $t=0$
	
				\item[\purpletag: ]
				that of $[d]_1 \in H_0(K_1;\rationals)$ only at scale $t=1$
			\end{enumerate}
		\end{minipage}%
	\end{center}
	\vspace{5pt}
	We consider some $0$\th homology classes in $H_0(K_t;\rationals)$ for varying $t \in \nonnegints$ below.
	\begin{enumerate}
		\item 
		The homology classes $[a]_0, [b]_0 \in H_0(K_0;\rationals)$ persist to $K_2$ since 
		$i^{2,0}\bigl( [a]_0 \bigr)
		= [a]_2 \neq 0$
		and 
		$i^{2,0}\bigl( [b]_0 \bigr)
		= [b]_2 \neq 0$.
		Moreover, $[a]_0$ and $[b]_0$ live or persist to $K_t$ for all $t \geq 1$.  
		At index $t=0$, \bluetagged{$[a]_0$\!} is highlighted in \bluetag and \redtagged{$[b]_0$\!} in \redtag in the illustration above. 

		\item 
		The homology classes $[a]_0, [b]_0 \in H_0(K_0;\rationals)$ merge in $K_1$ since 
		$i^{1,0}([a]_0) = [a]_1 = [b]_1 = i^{1,0}([b]_0)$.
		That is, the $0$-cycles $a$ and $b$ represent distinct homology classes in $K_0$ but represent the same homology class in $K_1$. 
		Furthermore, for all $t \geq 1$,
			\bluetagged{$[a]_t$} and \bluetagged{$[b]_t$} represent the same path component, highlighted in \bluetag in the illustration above.

		\item 
		The homology class $[a]_3 \in H_0(K_3;\rationals)$ is born in $K_0$ since $\im(i^{3,0}_\ast) = \rationals\ket{ [a]_3 }$.
		Similarly, the homology class $[d]_4 \in H_0(K_3; \rationals)$ is born in $K_0$ since $[d]_4 = [a]_4 = i^{4,0}([a]_0)$.
		These homology classes are highlighted in \bluetag in the illustration above.

		\item 
		The homology class \purpletagged{$[d]_1 \in H_0(K_1;\rationals)$}\!, specifically at index $t=1$ and highlighted in \purpletag above, is born at index $t=1$ since $\im(i^{1,0}_\ast) = \rationals\ket{[a]_1}$ and $[d]_1 \neq [a]_1$.
		Note that for all $t \geq 2$, \bluetagged{$[d]_t$} corresponds to the path component highlighted in \bluetag\!.
	\end{enumerate}
	Observe that none of the homology classes listed above die in the filtration.
\end{example}

Let $\filt{K}$ be a filtration on a simplicial complex $K$.
Observe that, in dimension $n=0$,
	a $0$\th homology class in $H_0(K_t;\field)$ for any $t \in \nonnegints$ cannot die in a filtration since path components do not disappear as we increase the time parameter.
If $K$ is path connected, 
	then all of these $0$\th homology classes merge into one homology class at sufficiently high $T \in \nonnegints$.
Persistence becomes more interesting in dimensions $n \geq 1$ since it becomes possible for homology classes to become trivial as we increase the time parameter.
We provide an example below on dimension $n=1$.

\begin{example}\label{ex:zom1-h2-loops-interpretation}
	Let $\filt{K}$ be as given in \fref{ex:pershom-one}.
	For convenience, we copied the illustration involving $H_0(K_t;\rationals)$ and directed loops in $K_t$ below.
	Note that if $t \not\in \set{2,3,4}$, then $H_1(K_t;\field)=0$.
	\vspace{5pt}
	\begin{center}
		\begin{minipage}[t]{0.32\linewidth}\null\vspace{-0.8\baselineskip}
			\includegraphics[width=\linewidth]{zomfig1/zom1-loops.png}\null
		\end{minipage}\hspace{0.01\linewidth}\begin{minipage}[t]{0.65\linewidth}\raggedright 
			Color Scheme for Directed Loops: \\[2pt]
			\begin{enumerate}[topsep=0pt, leftmargin=47pt, labelsep=0pt, itemsep=3pt] 
				\item[\redtag: ] 
				$\redmath{\alpha := ab + bc + cd - ad}$, representing $[\alpha]_t \in H_1(K_t;\rationals)$ for $t \geq 2$. 
				Note that $[\alpha]_t = 0$ if $t \geq 5$.

				\item[\greentag: ] 
				$\greenmath{\tau := ab + bc - ac}$,
				representing $[\tau]_t \in H_1(K_t;\rationals)$ for $t \geq 3$. 
				Note that $[\tau]_t = 0$ if $t \geq 4$.

				\item[\orangetag: ]
				$\orangemath{\sigma := ac + cd - ad}$,
				representing $[\sigma]_t \in H_1(K_t;\rationals)$ for $t \geq 3$.
				Note that $[\sigma]_t = 0$ if $t \geq 5$
				and that $[\alpha]_3 = [\tau]_3 + [\sigma]_3$.
			\end{enumerate}
		\end{minipage}%
	\end{center}
	\vspace{5pt}
	We consider some $1$\st homology classes in $H_1(K_t;\rationals)$ with $t \in \set{2,3,4}$ below.
	\begin{enumerate}
		\item 
		$[\alpha]_2 = [ab+bc+cd-ad]_2 \in H_1(K_2;\rationals)$ at index $t=2$ is born at index $t=2$, dies at index $t=5$, and has a lifespan of $3$.
		Illustrated below is the image of $[\alpha]_2$ on $H_n(K_t;\rationals)$ for $t=3,4$.

		\vspace{2pt}
		\begin{center}
		\begin{minipage}[t]{0.45\linewidth}\vspace{-0.9\baselineskip}
			\includegraphics[width=\linewidth]{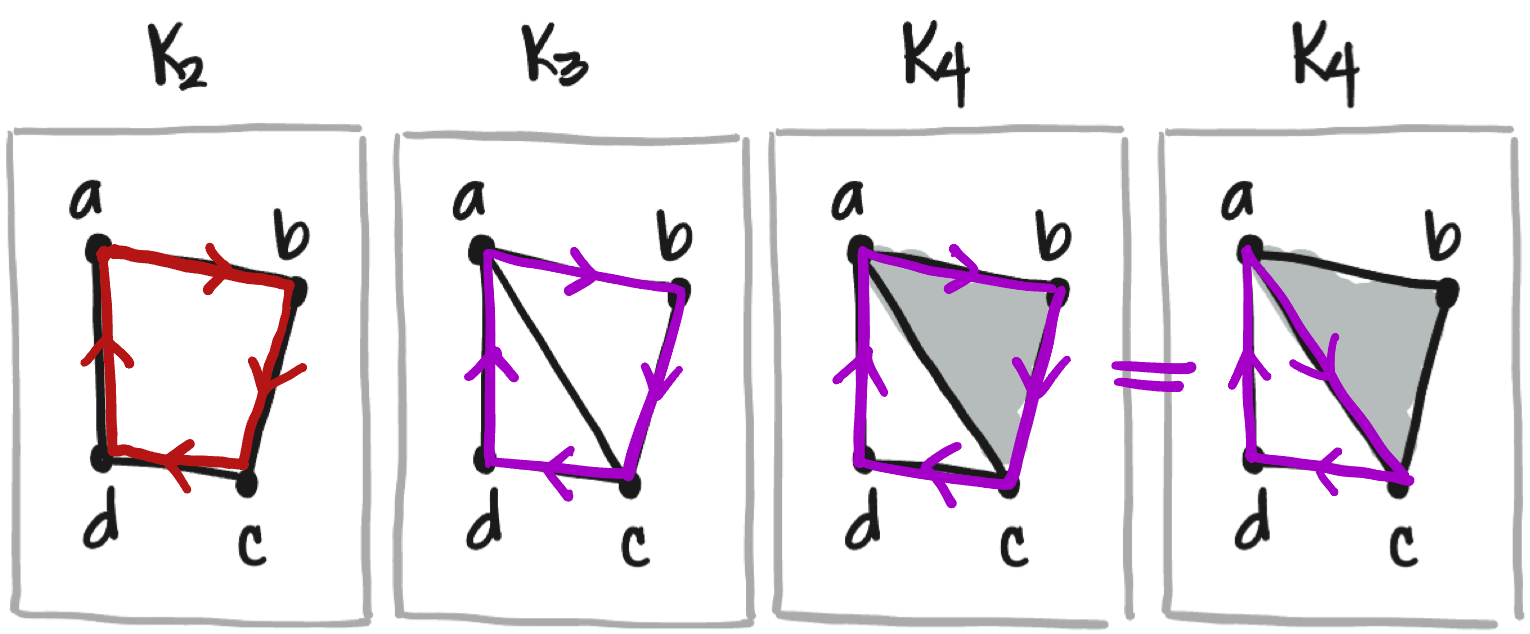}\null
		\end{minipage}\hspace{0.01\linewidth}\begin{minipage}[t]{0.5\linewidth}\raggedright 
			Color Scheme: \\[2pt]
			\begin{enumerate}[topsep=0pt, leftmargin=55pt, labelsep=0pt, itemsep=3pt]
				\item[\redtag: ] $\alpha$ in $K_2$, representing $[\alpha]_2$.
				\item[\purpletag: ] $\alpha$ in $K_3$ and $K_4$, representing $[\alpha]_3$ and $[\alpha]_4$ respectively.
			\end{enumerate}
		\end{minipage}%
		\end{center}
		\vspace{5pt}
		
		\item 
		$[\alpha]_3 = [ab+bc+cd-ad]_3 \in H_1(K_3;\rationals)$ at index $t=3$ is born at index $t=2$ since $i^{3,2}_\ast([\alpha]_2) = [\alpha]_3$.

		\item 
		$[\alpha]_4 = [ab+bc+cd-ad]_4 \in H_1(K_4;\rationals)$ is born at index $t=3$ since $i^{4,3}([\sigma]_3) = [\sigma]_4 = [\alpha]_4$ with $\sigma=ac+cd-ad$, as illustrated below:

		\vspace{2pt}
		\begin{center}
		\begin{minipage}[t]{0.35\linewidth}\vspace{-0.9\baselineskip}
			\includegraphics[width=\linewidth]{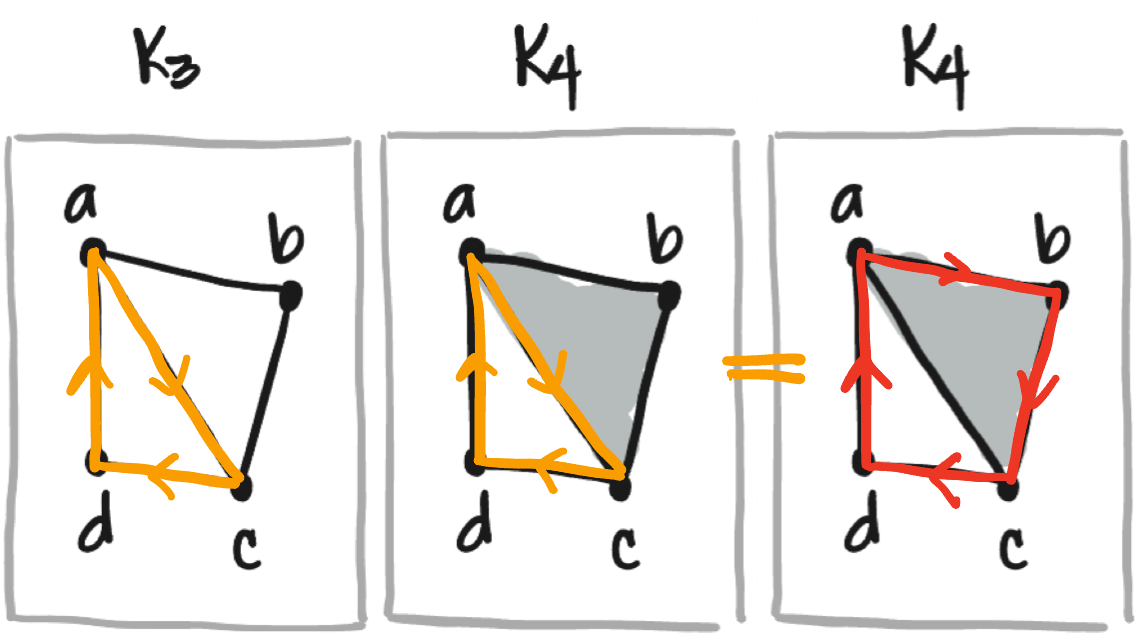}\null
		\end{minipage}\hspace{0.01\linewidth}\begin{minipage}[t]{0.55\linewidth}\raggedright 
			Color Scheme: \\[2pt]
			\begin{enumerate}[topsep=0pt, leftmargin=55pt, labelsep=0pt, itemsep=3pt]
				\item[\redtag: ] $\alpha$ in $K_4$, representing $[\alpha]_4$.

				\item[\orangetag: ] (middle) $\sigma \in K_4$, representing $[\sigma]_4 = [\alpha]_4$ 
				\\[2pt]
				(left) $\sigma$ in $K_3$, representing $[\sigma]_3$.
			\end{enumerate}
		\end{minipage}%
		\end{center}
		\vspace{5pt}

		\item 
		$[\tau]_3 = [ab+bc-ac]_3 \in H_1(K_3;\rationals)$ at index $t=3$ is born at the same index $t=3$.
		Since the $2$-simplex $abc$ is added to $K_4$,
			$\boundary_2^4(abc) = bc - ac + ab = \tau$
			and $\tau$ becomes a $1$-boundary in $K_4$.
		Therefore, $i^{4,3}([\tau]_3) = [\tau]_4 = 0$
		and $[\tau]_3$ dies in $K_4$.

		\item 
		$[\sigma]_3 \in H_1(K_3;\rationals)$ at index $t=3$ is born at index $t=3$, dies at index $t=5$, and has persistence $5-3=2$.
	\end{enumerate}
\end{example}

\spacer

Let $\filt{K}$ be a filtration on a finite simplicial complex $K$ and assume that $\filt{K}$ is constant on $[T,\infty)$. 
Note that $T \in \nonnegints$ exists by \fref{lemma:filtration-constant} and that $K = K_T$.
The persistence barcode $\barcode_n(K_\bullet;\field)$ of a filtration $K_\bullet$ with $\field$ coefficients determines the interval decomposition $H_n(K_\bullet;\field)$ as a persistence module.
Let $\barcode_n(K_\bullet;\field)$ be given as follows
\begin{equation*}
	\barcode_n(K_\bullet;\field) 
	= \barcode\bigl( H_n(K_\bullet;\field)\bigr)
	= \Bigl\{ 
		[s_1, s_1 + t_1),
		\ldots,
		[s_r, s_r + t_r),
		[s_{r+1}, \infty),
		\ldots, 
		[s_m, \infty]
	\Bigr\}
\end{equation*}
for some $s_1, \ldots, s_m \in \nonnegints$ and $t_1, \ldots, t_r \in \nonnegints$.
The intervals $[s_j,s_j+t_j)$ and $[s_j,\infty)$ in $\barcode_n(K_\bullet;\field)$ are denoted following \fref{cor:interval-decomp-from-structure-theorem},
where we state the correspondence between graded invariant factor decompositions in $\catgradedmod{\field}$ and interval decompositions in $\catpersmod$.
Therefore, there exists a persistence isomorphism as follows:
\begin{equation*}
	\phi_\bullet: H_n(K_\bullet;\field) 
	\to 
	\Bigl( 
		\intmod{[s_1, s_1+t_1)} 
		\oplus \cdots \oplus
		\intmod{[s_r, s_r + t_r)} 
		\oplus 
		\intmod{[s_{r+1}, \infty)} 
		\oplus \cdots \oplus 
		\intmod{[s_{m}, \infty)}
	\Bigr)
\end{equation*}
This persistence isomorphism identifies a collection of homology classes on varying indices $t \in \nonnegints$ by:
\begin{equation*}
	[\sigma_k]_{s_k} \in (\phi_{s_k})\inv 
	\Bigl(
		\intmod{J_k}(s_k)
	\Bigr)
	\subseteq H_n(K_{s_k}; \field)
	\quad\text{ with }\quad 
	J_k := \begin{cases}
		[s_k, s_k + t_k) &\text{ if } k \in \set{1, \ldots, r} \\
		[s_k,\infty) &\text{ if } k \in \set{r+1, \ldots, m}
	\end{cases}
\end{equation*}
Note that $(\phi_{s_k})\inv 
\bigl( \intmod{J_k}(s_k) \bigr)$ must be a one-dimensional $\field$-vector space and there exists a non-trivial homology class $[\sigma_k]_{s_k}$ for each $k \in \set{1, \ldots, m}$.
The birth, death, and persistence of $[\sigma_k]_{s_k}$ as a homology class in $H_0(K_\bullet;\field)$
correspond to the indices in the intervals $[s_k, s_k+t_k)$ or $[s_k, \infty)$.
In particular:
\begin{enumerate}
	\item 
	For $k \in \set{1, \ldots, r}$,
		the homology class $[\sigma_k]_{s_k}$ from $\intmod{J_k} = \intmod{[s_k,s_k + t_k)}$
	is born at index $s_k$, dies at index $s_k + t_k$, and has persistence $t_k$.
	Furthermore, 
		the behavior of the structure maps match the interval $[s_k, s_k + t_k)$, i.e.\ 
	\begin{equation*}
		i^{t,s_k}_\ast\bigl( [\sigma_k]_{s_k} \bigr)
		= \begin{cases}
			[\sigma_k]_t \neq 0 &\text{ if } t \in [s_k, s_k+t_k) \\
			0 &\text{ if } t \geq t_k
		\end{cases}
	\end{equation*}
	Observe that the bounded interval $[s_k,s_k+t_k)$ corresponds to homology classes that become trivial in $K=K_T$, i.e.\ $[\sigma_k]_T = i^{T,s_k}([\sigma_k]_{s_k}) = 0$ as an element of $H_n(K_T;\field) = H_n(K_T;\field)$.

	\item 
	For $k \in \set{r+1, \ldots, m}$,
		the homology class $[\sigma_k]_{s_k}$ from $\intmod{J_k} = \intmod{[s_k,\infty)}$ 
	is born at index $s_k$ and does not die in the filtration, i.e.\ for all $t \geq s_k$, 
	$
		i^{t,s_k}_\ast\bigl( [\sigma_k]_{s_k} \bigr)
		= [\sigma_k]_t \neq 0
	$.

	Note that the unbounded intervals $[s_k,\infty)$ correspond to homology classes that do not die in $K=K_T$,
		i.e.\ $[\sigma_k]_T = i^{T,s_k}([\sigma_k]_{s_k}) \neq 0$ 
		in $H_n(K;\field) = H_n(K_T;\field)$ for all $k \in \set{r+1, \ldots, m}$.
	This might explain why $s_k$ in the interval $[s_k,\infty)$ is called an \textit{essential} birth index of $K_\bullet$.
\end{enumerate}
The problem here is that such a collection $\set{[\sigma_1]_{s_1}, \ldots, [\sigma_m]_{s_m}}$ of homology classes
	is not generally immediately obvious by examination of the homology groups $H_n(K_t;\field)$ at every $t \in \nonnegints$ since the structure maps, i.e.\ the maps on homology induced by inclusion, cannot be ignored.

The matrix reduction algorithm for persistent homology, discussed in \fref{chapter:matrix-calculation},
	not only determines the persistence barcode of $H_n(\filt{K};\field)$ but also determines possible cycle representatives corresponding to each interval in said barcode.
In the two examples below, we interpret results that will be calculated in \fref{chapter:matrix-calculation}.

\begin{example}\label{ex:zom1-zero-pershom-interval-decomp}
	Let $\filt{K}$ be given as in \fref{ex:pershom-one}.
	Given below is a summary of the results of the calculation involving $H_0(\filt{K};\rationals)$ started in \fref{ex:zom1-gradedboundary1-permutation} and finished in \fref{ex:zom1-gradedboundary1-finalcalc}:
	\begin{longtable}{ccccc}
		Interval Module 
		& Interval
		& Birth Index
		& Cycle Representative 
		& Persistence 
		\\[5pt]
		$\intmod{[0,1)}$ & $J_1 = [0,1)$
		& $s_1=0$
		& $\sigma_1 = b-a$
		& $t_1 = 1-0 = 1$
		\\[2pt] 
		$\intmod{[1,2)}$ & $J_2 = [1,2)$
		& $s_2=1$
		& $\sigma_2 = d-a$
		& $t_2=2-1 = 1$
		\\[2pt]
		$\intmod{[0,\infty)}$ & $J_3 = [0,\infty)$
		& $s_3=0$ 
		& $\sigma_3 = a$
		& $t_3=\infty$
	\end{longtable}
	\vspace{-\parskip}\noindent
	The interval decomposition of $H_0(K_\bullet;\rationals)$ is then calculated to be 
	\begin{equation*}
		H_0(\filt{K}; \rationals) 
		\,\uppersmod\cong\,
		\intmod{[0, \infty)} \oplus 
		\intmod{[0,1)} \oplus 
		\intmod{[1,2)}
	\end{equation*}
	Observe that each interval module $\intmod{J_k}$ corresponds to a homology class that is born at index $s_k$ and dies at index $s_k + t_k$.
	However, these homology classes do not correspond directly to the path components of each $K_t$.
	For comparison, we copied the illustration in \fref{ex:pershom-one} involving $H_0(K_\bullet;\rationals)$:
	\vspace{5pt}
	\begin{center}
		\begin{minipage}[t]{0.63\linewidth}\null\vspace{-0.7\baselineskip}
			\includegraphics[width=\linewidth]{zomfig1/zom1-path-components.png}\null
		\end{minipage}\hspace{0.01\linewidth}\begin{minipage}[t]{0.34\linewidth}\raggedright 
			Color Scheme for Path Components: \\[2pt]
			\begin{enumerate}[topsep=0pt, leftmargin=55pt, labelsep=0pt]
				\item[\bluetag: ] 
				that of $[a]_t \in H_0(K_t;\rationals)$ at all scales $t \in \nonnegints$
	
				\item[\redtag: ]
				that of $[b]_0 \in H_0(K_0;\rationals)$ only at scale $t=0$
	
				\item[\purpletag: ]
				that of $[d]_1 \in H_0(K_1;\rationals)$ only at scale $t=1$
			\end{enumerate}
		\end{minipage}%
	\end{center}
	\vspace{5pt}
	Observe that $\intmod{[0,\infty)}$ corresponds to the homology classes $[a]_t$ with $t \in \nonnegints$ and has a direct relation to the path component of $K_t$ containing the $0$-cycle $a$.
	However, the same does not apply for $\intmod{[0,1)}$ and $\intmod{[1,2)}$, which correspond to the homology classes 
	$[b-a]_0$ and $[d-a]_1$ respectively.
\end{example}

\begin{example}\label{ex:zom1-one-pershom-interval-decomp}
	Let $\filt{K}$ be given as in \fref{ex:pershom-one}.
	A summary of the results calculated in \fref{ex:zom1-h1-matrix-calculation} involving the interval decomposition of $H_1(K_\bullet;\rationals)$
	is given below:
	\begin{longtable}{ccccc}
		Interval Module 
		& Interval
		& Birth Index
		& Cycle Representative 
		& Persistence 
		\\[5pt]
		$\intmod{[2,5)}$ & $J_1 = [2,5)$
		& $s_1=2$
		& $\beta_1 = ab+bc-ad+cd$
		& $t_1 = 5-2 = 3$
		\\[2pt] 
		$\intmod{[3,4)}$ & $J_2 = [3,4)$
		& $s_2=3$
		& $\beta_2 = -ab -bc + ac$
		& $t_2=4-3=1$
	\end{longtable}
	\vspace{-\parskip}\noindent
	The interval decomposition of $H_1(K_\bullet;\rationals)$ is given by $H_1(K_\bullet;\rationals) 
	= \intmod{[2,5)} \oplus \intmod{[3,4)}$.
	For convenience,
		we copied the illustration in \fref{ex:pershom-one} describing the orientation on the $1$-simplices of $K_t$ at each $t \in \nonnegints$ induced by the vertex order $\Vertex(K) = (a,b,c,d)$ below.
	\begin{center}
		\baselineCenter{\includegraphics[height=1.1in]{zomfig1/zom1-simp-oriented.png}}
		\quad
		\baselineCenter{\includegraphics[height=1.1in]{zomfig1/zom1-filt-oriented.png}}
	\end{center}
	\vspace{5pt}\noindent
	The first interval module $\intmod{[2,5)}$ corresponds to the 
	homology classes $[\beta_1]_t$ with $t \in [2,5)$.
	Observe that the $1$-cycle $\beta_1 = ab+bc-ad+cd$, illustrated in \redtag below for $t \in [2,5)$,
		first appears at index $s_1 = 2$.
	The $2$-cycle $abc + acd$, that makes $\beta_1$ a $1$-boundary and the homology class $[\beta_1]_t$ trivial, 
	first appears at index $s_1 + t_1 = 2+3=5$.
	\begin{center}
		\includegraphics[width=4in]{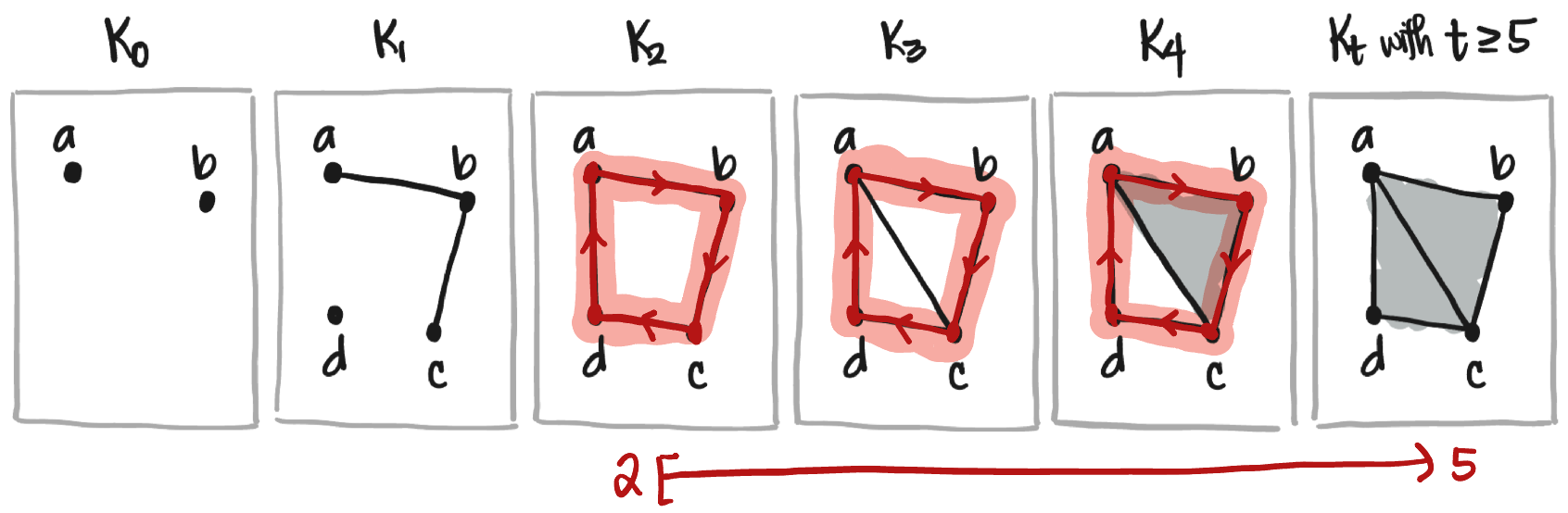}
	\end{center}
	The second interval module $\intmod{[3,4)}$
	corresponds to the homology classes $[\beta_2]_t$ with $t \in [3,4)$.
	Observe that the $1$-cycle $\beta_2 = -ab-bc+ac$, illustrated below in \purpletag for $t \in [3,4)$, first appears at $s_2 = 2$.
	The $2$-cycle $-abc$, which makes $\beta_2$ a $1$-boundary and the homology class $[\beta_2]_t$ trivial, first appears at index $s_2 + t_2 = 4$.
	\begin{center}
		\includegraphics[width=4in]{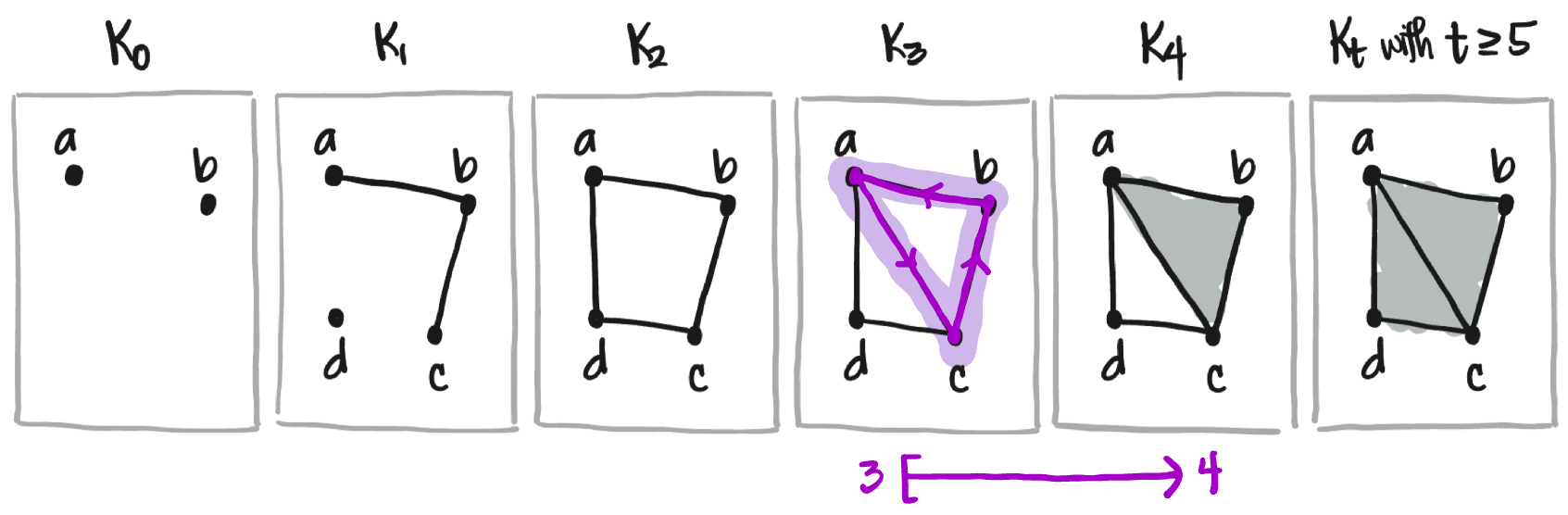}
	\end{center}
	Note that the $1$-cycle $ac + cd - ad$ is not represented in either interval module for $H_1(K_\bullet;\rationals)$, despite being a valid cycle representative for homology classes in $H_1(K_t;\rationals)$ for $t \geq 3$.
\end{example}

For the sake of completion, we include an alternate characterization for persistent homology below, one that is prevalent in persistent homology literature, e.g.\ in \cite{matrixalg:zomorodian,matrixalg:cohen-steiner}.
While this perspective is not as useful as the functor definition of persistent homology (\fref{defn:pershom}) for rigorously discussing the origins of the matrix reduction algorithm (as presented in \fref{chapter:matrix-calculation}),
it does offer some more insight to the significance of persistence barcodes.
The following definitions are taken from~\cite[Section 2.6]{matrixalg:zomorodian}. 

\begin{statement}{Definition}\label{defn:p-persistent-bullshit}
	Let $\filt{K}$ be a simplicial filtration and let $\field$ be a field.
	The $p$-\textbf{persistent} $n$\th \textbf{homology group} $H_n(K_t; p; \field)$ \textbf{with coefficients in $\field$} of the simplicial complex $K_t$ in $\filt{K}$ is the $\field$-vector space given by 
	\begin{equation*}
		H_n(K_t; p; \field) 
		= i_\ast^{t+p, t}\Bigl( H_n(K_t; \field) \Bigr)
		\subseteq 
		H_n(K_{t+p}; \field)
	\end{equation*}
	and the rank of $H_n(K_t;p;\field)$ is called the $p$-\textbf{persistent} $n$\th \textbf{Betti number} $\betti_n(K_t; p; \field)$ of $K_t$ \textbf{with coefficients in $\field$}.
\end{statement}

Observe that the characterization of persistent homology as a persistence module $H_n(\filt{K};\field)$ by \fref{defn:pershom} accounts for the $p$-persistent homology groups $H_n(K_t; p; \field)$ for all $p, t \in \nonnegints$
and is a more concise way to describe how all of these homology groups are related. 
Furthermore, the persistence barcode $\barcode_n(\filt{K};\field)$ of a filtration $\filt{K}$ is a concise representation of the ranks of all $p$-persistent homology groups.
In other words, the persistence barcode encodes the evolution of the Betti numbers of $K_t$ across the filtration $\filt{K}$.
We state this in more detail below.

\begin{corollary}\label{cor:betti-numbers-pers-barcode}
	Let $\filt{K}$ be a filtration on a finite simplicial complex $K$.
	For all $n \in \nonnegints$,
		the $n$\th Betti numbers $\betti_n(K_t; \field)$
		of the simplicial complex $K_t$
		and the $n$\th $p$-persistent Betti numbers $\beta_n(K_t;p;\field)$ of $K_t$ in $\filt{K}$ 
	are determined by the persistence barcode 
	$\barcode_n(K_\bullet;\field)$ as follows:
	\begin{align*}
		\beta_n(K_t; \field) &= \card\!\Big\{
			J \in \barcode_n(\filt{K}; \field)
			\,:\,
			t \in J
		\Big\}
		\\
		\beta_n(K_t; p; \field) &= \card\!\Big\{
			J \in \barcode_n(\filt{K}; \field)
			\,:\,
			[t, t+p] = [t,t+p+1) \subseteq J
		\Big\}
	\end{align*}
\end{corollary}
\begin{proof}
	By \fref{prop:filtrations-finite-type-has-interval},
		an interval decomposition for $H_n(K_\bullet; \field)$ exists.
	Applying \fref{prop:interval-decomposition-and-ranks} on the $\field$-vector space $H_n(K_t;\field)$ for all $t \in \nonnegints$, we have that
	\begin{equation*}
		\rank\bigl( H_n(K_t);\field \bigr)
		= \beta_n(K_t;\field) = \card\!\Big\{
			J \in \barcode_n(\filt{K}; \field)
			\,:\,
			t \in J
		\Big\}.
	\end{equation*}
	Applying \fref{prop:interval-decomposition-and-ranks} on structure maps $i^{t+p,t}_\ast$ of $H_n(K_t;\field)$ for all $t,p \in \nonnegints$, we have that
	\begin{align*}
		\rank\bigl( i^{t+p,t}_\ast \bigr)
		&= \rank\Bigl( 
			i^{t+p,t}_\ast\bigl( H_n(K_t;\field) \bigr) 
		\Bigr)
		= \rank\Bigl( H_n(K_t;p;\field) \Bigr)
		= \beta_n(K_t;p;\field)
		\\
		&= \card\!\Big\{
			J \in \barcode_n(\filt{K}; \field)
			\,:\,
			[t, t+p] = [t,t+p+1) \subseteq J
		\Big\}.
	\end{align*}\par\vspace{-41pt}
\end{proof}

\begin{example}
	Let $\filt{K}$ be given as in \fref{ex:pershom-one}. 
	For convenience, the illustration of the simplicial complex $K$ and the filtration $\filt{K}$ is copied below:
	\begin{center}
		\baselineCenter{\includegraphics[height=1.05in]{zomfig1/zom-simp.png}}
		\quad
		\baselineCenter{\includegraphics[height=1.1in]{zomfig1/zom1-cleaned.png}}
	\end{center}
	\vspace{3pt}
	Following the discussion on \fref{ex:zom1-zero-pershom-interval-decomp},
		the interval decomposition 
		$H_0(K_\bullet;\rationals) \cong \intmod{[0,1)} \oplus \intmod{[1,2)} \oplus \intmod{[0,\infty)}$ implies that 
	\begin{equation*}
		\barcode_0(K_\bullet;\rationals) = \Bigl\{
			[0,\infty), [0,1), [1,2)
		\Bigr\}
	\end{equation*}
	Below, we compare the $p$-persistent $0$\th homology group $H_0(K_t;p;\rationals)$ and the $p$-persistent $0$\th Betti number 
	$\betti_0(K_t;p;\rationals)$ for selected $p,t \in \nonnegints$ and compare these against the results expected by \fref{cor:betti-numbers-pers-barcode}.
	\begin{enumerate}
		\item 
		Let $t = 1$.
		Since $K_t = K_1$ has two path components,
			we expect $\betti_0(K_1;\rationals) = 2$.
		There are $2$ intervals in $\barcode_0(K_\bullet;\rationals)$ that contain $t=1$: $[0,\infty), [1,2) \in \barcode_0(K_\bullet;\rationals)$.

		\item 
		Let $t=0$ and $p=1$.
		Since $K_t = K_0$ has three path components and $K_{t+p} = K_1$ has two path components, 
			we expect that $H_0(K_0;1;\rationals) \cong \rationals^2$
			and $\betti_0(K_0;1\rationals) = 2$.
		Observe that there are $2$ intervals in $\barcode_0(K_\bullet;\rationals)$ containing the interval $[t,t+p] = [0,1]$: $[0,\infty), [1,2) \in \barcode_0(K_\bullet;\rationals)$.

		\item 
		Let $t=1$ and $p=3$.
		Since $K_t = K_1$ has two path components and $K_{t+p} = K_4$ has one path component,
			we expect that $H_0(K_1;3;\rationals) \cong \rationals$ and $\betti_0(K_1;3;\rationals) = 3$.
		Observe that there is only $1$ interval in 
		$\barcode_0(K_\bullet;\rationals)$ containing $[t,t+p] = [1,4]$: $[0,\infty) \in \barcode_0(K_\bullet;\rationals)$.
	\end{enumerate}
\end{example}\clearpage

\section{Simplicial Persistent Homology}
\label{section:construction-of-persistent-homology}
\label{section:simplicial-persistent-homology}

In the previous section,
we interpreted the $n$\th persistent homology module $H_n(K_\bullet;\field)$ of a simplicial filtration $K_\bullet$ as done pointwise, with each index $t \in \nonnegints$ considered separately. 
In this section, we consider an approach wherein we consider all indices $t \in \nonnegints$ of $K_\bullet$ simultaneously.
In particular, 
	we extend the notions of simplicial chain groups, boundary maps, and chain complexes to the case of \textit{persistence modules}.
For brevity, we will refer to this extension as \textit{simplicial persistent homology}.

We want to emphasize that the persistence modules and persistence complexes we introduce in this section are not exactly new. Rather, they come as a natural consequence of the functorial nature of simplicial homology.
Earlier in \fref{section:functors-in-simplicial-homology}, 
we discussed how the calculation of the simplicial homology can be seen as a composition of functors, as illustrated below for the simplicial complex $K_t$ of $K_\bullet$:
\begin{equation*}
\begin{tikzcd}[column sep=7em]
	\hspace{5pt}\catsimp\hspace{5pt}
		\arrow[r, "C_\ast(-;R)"]
		\arrow[r, mapsto, yshift=-1.75em]
	& \hspace{0.7in}
		\catchaincomplex{\catvectspace}
		\hspace{0.7in}
		\arrow[r, "H_n(-)"]
		\arrow[r, mapsto, yshift=-1.75em]
	& \hspace{10pt}\catvectspace\hspace{10pt}
	\\[-1.8em] 
	K_t
		& C_\ast(K_t;\field) = \bigl( C_n(K_t;R), \boundary_n \bigr)_{n \in \ints}
		& H_n(K_t;\field) 
	\\[-2.0em] 
	\mathclap{\text{(simplicial complex)}} 
		& \mathclap{\text{(simplicial chain complex)}} 
		& \mathclap{\text{(simplicial homology)}}
\end{tikzcd}	
\end{equation*} 
In simplicial persistent homology, we leave the parameter $t \in \nonnegints$ in $K_\bullet$ arbitrary. We present a loose visualization of this below, with the parts highlighted in \redtag to be discussed in this section.
\newcommand{\functorcatthing}[1]{\ensuremath\catname{[}#1\catname{]}}
\begin{equation*}
\begin{tikzcd}
	\functorcatthing{ \posetN, \catsimp }
		\arrow[r]
		\arrow[r, mapsto, yshift=-3.4em]
	& \redmath{\begin{gathered}
			\functorcatthing{ \posetN, \catchaincomplex{\catvectspace} } \\
			= \catchaincomplex{\catpersmod} 
		\end{gathered}}
		\arrow[r]
		\arrow[r, mapsto, yshift=-3.4em]
	& \functorcatthing{ \posetN, \catvectspace } = \catpersmod
	\\[-1.9em] 
	K_\bullet 
	& \redmath{C_\ast(K_\bullet;\field)} 
	& H_n(K_\bullet;\field) 
	\\[-1.9em] 
	\text{(filtration)}
	& \redmath{\text{(simplicial persistence complex?)}} 
	& \text{(persistent homology module)}
\end{tikzcd}	
\end{equation*}
As we construct new objects in $\catpersmod$, 
	we also investigate the corresponding graded structures in $\catgradedmod{\field}$ resulting from the application of $\togrmod: \catpersmod \to \catgradedmod{\field}$ to these objects 
	and discuss which properties are preserved under this category equivalence.

\spacer 

We start with a definition of \textit{filtered chain modules} of simplicial filtrations, which are simplicial chain groups extended to the case of persistence modules. 
Note that this definition is adapted from \cite{ripser}. 

\begin{statement}{Definition}
	\label{defn:filtered-n-chain-module}
	Let $\filt{K}$ be a filtration on a simplicial complex $K$ and fix a field $\field$.
	For each $n \in \ints$,
		define the $n$\th \textbf{filtered chain module $C_n(\filt{K}; \field)$ of $\filt{K}$ with coefficients in $\field$} 
		to be the persistence module over $\field$ given by the following functor composition:
	\begin{equation*}
		C_n(\filt{K}; \field) := C_n(-; \field) \circ \filt{K}
	\end{equation*}
	where $C_n(-; \field): \catsimp \to \catvectspace$ refers to the $n$\th chain group functor with coefficients in $\field$. 
\end{statement}
\remarks{
	\item 
	$C_n(K_\bullet;\field)$ is well-defined since the codomain category of $K_\bullet: \posetN \to \catsimp$ and the domain category of $C_n(-;\field): \catsimp \to \catvectspace$ match. 
	Note that $K_\bullet$ is a functor by decomposition and $C_n(-;\field)$ is a well-defined functor as discussed above \fref{defn:simplicial-chain-group-functor}. 

	\item 
	The modifier ``\textit{filtered}'' implies that the chain groups are those from a simplicial \textit{filtration}. 
	We decided to use the term ``\textit{module}'', instead of calling $C_n(K_\bullet;\field)$ a filtered chain \textit{group},
	to emphasize that $C_n(K_\bullet;\field)$ is to be interpreted mainly as a persistence \textit{module}.
}

Observe that $C_n(K_\bullet;\field)(t) = C_n(K_t;\field)$ for all $t \in \nonnegints$ and the notation involving the vector spaces of $C_n(K_\bullet;\field)$ should work as expected.
For example, passing the parameter $t \in \nonnegints$ to $C_n(K_\bullet;\field)$ returns the simplicial chain group of $K_t$ of $K_\bullet$.
This is also consistent with the convention of writing $V_t := V_\bullet(t)$ for an arbitrary persistence module $V_\bullet$, as introduced in \fref{defn:persmod}.

As given in \fref{defn:filtration}, 
	the inclusion maps of $K_\bullet$ are denoted as $i^{s,t}: K_t \to K_s$, with the shorthand of $i^{t} = i^{t+1,t}: K_t \to K_{t+1}$.
Since the structure maps of $C_n(K_\bullet;\field)$ must be the maps on the simplicial chain groups induced by inclusions, 
	we can extend the hash $(\hash)$ notation identified in \fref{defn:induced-maps-on-chain-groups}
	and denote the structure maps of $C_n(K_\bullet;\field)$ as follows:
\begin{equation*}
\begin{array}{l ccc}
	&	i_\hash^{s,t} = C_n(\filt{K}; \field)(t \to s)
		& \text{ with }
		& i_\hash^{s,t}: C_n(K_t;\field) \to C_n(K_s;\field) 
	\\[5pt]
	\text{ and } 
		& i_\hash^{t} = C_n(\filt{K};\field)(t \to t+1)
		& \text{ with } 
		& i_\hash^{t}: C_n(K_t;\field) \to C_n(K_{t+1};\field)
\end{array}
\end{equation*}
We provide an example of a filtered chain module below.

\begin{example}\label{ex:chains-one}
	Let $\filt{K}$ and $K$ be given as in \fref{ex:pershom-one}
	and orient $K$ with the vertex order $\Vertex(K) = (a,b,c,d)$.
	For convenience, an illustration of $K$ and $K_\bullet$ (without orientation) is copied below:
	\begin{center}
		\baselineCenter{\includegraphics[height=1.05in]{zomfig1/zom-simp.png}}
		\quad
		\baselineCenter{\includegraphics[height=1.1in]{zomfig1/zom1-cleaned.png}}
	\end{center}
	\vspace{3pt}
	The $0$\th and $1$\st filtered chain modules $C_0(\filt{K}; \rationals)$ 
	and $C_0(\filt{K}; \rationals)$ of $\filt{K}$ with rational coefficients
	have the following vector spaces, described relative to their respective standard bases:
	\begin{equation*}
		C_0(K_t; \rationals) = \begin{cases}
			\rationalsket{a,b}			&\text{ if } t = 0 \\
			\rationalsket{a,b,c,d}		&\text{ if } t \geq 1
		\end{cases}
		\qquad\text{ and }\qquad 
		C_1(K_t; \rationals) = \begin{cases}
			0								&\text{ if } t = 0 \\
			\rationalsket{ab, bc}			&\text{ if } t = 1 \\
			\rationalsket{ab, bc, cd, ad}	&\text{ if } t = 2 \\
			\rationalsket{ab, bc, cd, ad, ac}			&\text{ if } t \geq 3 \\
		\end{cases}
	\end{equation*}
	Below, we list some examples of $0$-chains and $1$-chains at selected values of $t \in \nonnegints$:
	\begin{enumerate}
		\item 
			The $0$-chain $2a+b$ is an element of $C_0(K_t; \rationals)$ for all $t \in \nonnegints$. 
			In contrast, the chain $a + 2b + c$ is not an element of $C_0(K_0; \rationals)$ since the simplex $c$ is not in $K_0$.
		\item 
			The $1$-chain $ab + bc - ac$ first appears in $C_1(K_3; \rationals)$, in that $t=3$ is the minimal index for which $ab+bc-ac \in C_1(K_t; \rationals)$.
			Observe that $ab+bc-ac \in C_1(K_s; \rationals)$ for any $s \geq 3$.
		\item 
			The $1$-chain $2ab - bc$ as an element of $C_1(K_1; \rationals)$ is mapped to $i^{s,1}_\hash(2ab - bc) = 2ab - bc \in C_1(K_s; \rationals)$ for any $s \geq 1$.
	\end{enumerate}
\end{example}

Given a filtration $K_\bullet$ of a simplicial complex $K$,
	the vector spaces $C_n(K_t;\field)$ of the filtered chain module $C_n(K_\bullet;\field)$ can all be considered subspaces of the simplicial chain group $C_n(K;\field)$ of $K$.
Additionally, the structure maps of $C_n(K_\bullet;\field)$ are all inclusions and agree with the identity map on $C_n(K;\field)$.
We state this in more detail below.

\begin{corollary}\label{cor:chain-module-inclusion}
	Let $\filt{K}$ be a filtration of some simplicial complex $K$.
	For all $t,s \in \nonnegints$ with $t \leq s$,
		$C_n(K_t;\field)$ is a vector subspace of $C_n(K;\field)$ 
		and the structure map 
		$i^{s,t}_\hash: C_n(K_t; \field) \to C_n(K_s; \field)$ 
		satisfies
		$i^{s,t}_\hash(\sigma) = \id_{C_n(K;\field)}(\sigma) = \sigma$
		for all $\sigma \in C_n(K_t;\field)$.
\end{corollary}
\begin{proof}
	Let $t,s \in \nonnegints$ with $t \leq s$.
	By definition of $K_\bullet$, $K_t \subseteq K_s \subseteq K$.
	By \fref{lemma:boundary-on-chain-groups-by-inclusion}(i),
		$C_n(K_t;\field) \subseteq C_n(K_s;\field) \subseteq C_n(K;\field)$ 
		and for all $n$-chains $\sigma \in C_n(K_t;\field)$,
	$
		i^{s,t}_\hash(\sigma)
		= \id_{C_n(K_s;\field)}(\sigma) 
		= i^{[s]}_\hash(\sigma)
		= \id_{C_n(K;\field)}(\sigma)
		= \sigma
	$
	where $i^{[s]}$ refers to the inclusion map $i^{[s]}: K_s \to K$.
\end{proof}

The corollary above does not take any imposed orientation on either $K$ or on any $K_t$ of $K_\bullet$ into account.
However, for convenient calculation, 
	we usually set an orientation on $K$ and let all simplicial complexes $K_t$ of $K_\bullet$ inherit said orientation by restriction.
	In this case, the standard basis for any $C_n(K_t;\field)$ at index $t \in \nonnegints$ (relative to the inherited orientation) must be a subset of that of $C_n(K;\field)$.
	This property, along with the corollary above, makes the graded module obtained by applying $\togrmod$ to the filtered chain modules have a relatively uncomplicated structure.
Below, we identify notation for said graded module.

\begin{definition}\label{defn:graded-chain-module}
	For each $n \in \ints$, 
		the $n$\th \textbf{graded chain module} of a simplicial filtration $K_\bullet$ {with coefficients in a field} $\field$ 
		is the graded $\field[x]$-module given by 
	$
		C_n\graded(K_\bullet;\field) := \togrmod(
			C_n\graded(K_\bullet;\field)
		)
	$.
	We call an element of $C_n(K_\bullet;\field)$ a \textbf{filtered} $n$-\textbf{chain}.
\end{definition}
\remark{
	The term ``graded chain module'' is not used in most of the literature for persistent homology theory. 
	We introduced this term in this paper to emphasize 
	the difference between $C_n(K_\bullet;\field)$ as a persistence module and 
	$C_n\graded(K_\bullet;\field)$ as a graded module. The symbol $C_n\graded(K_\bullet;\field)$ is added for brevity.
}

Following the description of $\togrmod$ given in \fref{defn:togrmod}, 
	we have the following set description for the graded chain module $C_n\graded(K_\bullet;\field)$:
	\begin{equation*}
		C_n\graded(K_\bullet;\field) 
		\upabelian= \bigoplus_{t \in \nonnegints} C_n(K_t;\field) x^t
	\end{equation*}
Then, a filtered $n$-chain $\sigma(x) \in C_n\graded(K_\bullet;\field)$ is an $\field[x]$-formal sum $\sigma(x) = \sum_{t \in \nonnegints} \sigma_t x^t$ of $n$-chains such that $\sigma_t \in C_n(K_t;\field)$, i.e.\ an $n$-chain at index $t \in \nonnegints$, for all $t \in \nonnegints$.
Note that the direct sum characterization also implies that only finitely many of $\sigma_t$ can be non-trivial.

Observe that this formal sum notation for filtered $n$-chains makes identifying the homogeneity and degree of elements in $C_n\graded(K_\bullet;\field)$ straightforward.
This becomes useful since we prefer dealing with homogeneous elements of graded chain modules, as we will see later in \fref{chapter:matrix-calculation}.
To avoid losing this benefit, we usually avoid using the indeterminate $x$ of $\field[x]$ as a vertex of the simplicial complex $K$.

Before we proceed with an example, 
we have identified two useful properties of $C_n\graded(K_\bullet;\field)$ resulting from \fref{cor:chain-module-inclusion} below.

\begin{lemma}\label{lemma:chain-group-graded-mod-char}
	Let $\filt{K}$ be a filtration of some simplicial complex $K$.
	Let $n \in \ints$.

	\begin{enumerate}
		\item 
		For all $t \in \nonnegints$,
		$C_n(K_t;\field) x^t$ is a vector subspace of $C_n(K;\field) x^t$. 
		Moreover, 
		if a filtered $n$-chain is of the form $\sigma x^t \in C_n\graded(K_\bullet;\field)$ with $\sigma \in C_n(K;\field)$ and $t \in \nonnegints$,
		then $\sigma \in C_n(K_t;\field)$.

		\item 
		The action of $\field[x]$ on $C_n\graded(K_\bullet;\field)$ satisfies 
		$x^{s} \cdot \sigma x^t = \sigma x^{t+s}$ for all $\sigma x^t \in C_n\graded(K_\bullet;\field)$ and for all $t,s \in \nonnegints$.
	\end{enumerate}
\end{lemma}
\begin{proof}
	For \textbf{(i)}:
		$C_n(K_t;\field) x^t \subseteq C_n(K;\field) x^t$ since 
		$C_n(K_t;\field)$ is a vector subspace of $C_n(K;\field)$ by \fref{cor:chain-module-inclusion}.
	Assume there exists $\sigma x^t \in C_n\graded(K_\bullet;\field)$ with $\sigma \in C_n(K;\field)$ and $t \in \nonnegints$.
		We can assume, without loss of generality, $x^t \not\in C_n(K_s;\field)$ for all $s \in \nonnegints$.
		Therefore, $\sigma x^t$ can only be present in $C_n(K_t;\field)x^t$, the homogeneous component of $C_n\graded(K_\bullet;\field)$ of degree $t$.
		Therefore, $\sigma \in C_n(K_t;\field)$.
		
	For \textbf{(ii)}: 
		Let $t,s \in \nonnegints$ and let $\sigma x^t \in C_n\graded(K_\bullet;\field)$. Then, $\sigma \in C_n(K_t;\field)x^t$.
		By definition of $\togrmod$, 
			$x^s \cdot \sigma x^t = i_\hash^{t+s,s}(\sigma) x^{t+s}$.
		By \fref{cor:chain-module-inclusion}, 
			$i_\hash^{t+s,s}(\sigma) x^{t+s} = \sigma x^{t+s}$.
\end{proof}

We provide an example of a filtered chain module viewed as a graded $\rationals[x]$-module below.

\begin{example}\label{ex:zom1-filtered-chains}
	Let $K$ and $\filt{K}$ be as given in \fref{ex:pershom-one}, illustrations of which are copied below for convenience.
	Orient $K$ and each $K_t$ by $(a,b,c,d)$, with the vertex set restricted when appropriate.
	\begin{center}
		\baselineCenter{\includegraphics[height=1.05in]{zomfig1/zom-simp.png}}
		\quad
		\baselineCenter{\includegraphics[height=1.1in]{zomfig1/zom1-cleaned.png}}
	\end{center}
	\vspace{3pt}
	\noindent 
	A description of the $0$\th graded chain module $C_0\graded(K_\bullet;\rationals)$ of $K_\bullet$ is given below.
	Note that the direct sums below are to be interpreted as internal direct sums of $\rationals$-vector spaces.
	\begin{align*}
		C_0\graded(K_\bullet;\rationals) 
		&= 
		\,\togrmod\bigl( C_0(K_\bullet;\rationals) \bigr)
		\upabelian= \bigoplus_{t \in \nonnegints} C_0(K_t;\rationals) x^t 
		= \rationals\ket{a,b} 
		\oplus \left(\, \bigoplus_{t=1}^\infty \rationals\ket{a,b,c,d}x^t \right)
		\\ 
		&= \biggl\{ 
			a \cdot f_1(x) + b \cdot f_2(x) + c \cdot x f_3(x) + d \cdot x f_4(x)
			\,:\, 
			f_i(x) \in \rationals[x] \text{ for } i=1,2,3,4
		\biggr\}
	\end{align*}
	Using \fref{lemma:chain-group-graded-mod-char}, 
	we know that the action of $\rationals[x]$ on $C_0\graded(K_\bullet;\field)$ 
	satisfies the following:
	\begin{equation*}
	\begin{array}{c c c l}
		x \cdot a x^t = a x^{t+1} 
			& \text{ and } 
			& x \cdot b x^t = b x^{t+1} 
			& \text{ for all } t \in \nonnegints 
		\\[5pt]
		x \cdot c x^t = c x^{t+1} 
			& \text{ and } 
			& x \cdot d x^t = d x^{t+1} 
			& \text{ for all } t \geq 1 
	\end{array}
	\end{equation*}
	Listed below are some filtered $0$-chains in $C_0\graded(K_\bullet;\rationals)$, along with some comments.
	\begin{enumerate}
		\item 
		The filtered $0$-chain $\sigma_1(x) := 2ax + 3bx^3$ corresponds to the pair 
		of
			the $0$-chain $2a \in C_0(K_1; \rationals)$ at index $1$
		and 
			the $0$-chain $3b \in C_0(K_3; \rationals)$ at index $3$.
		Note that $\sigma_1(x)$ is not a homogeneous element since 
		$\degh(2ax) = 1 \neq 3 = \degh(3bx^3)$.
		In this case, $\degh(\sigma_1(x))$ is undefined.

		\item 
		The element $c$ is not a filtered $0$-chain of $K_\bullet$ since 
		$c \not\in C_0(K_0;\rationals) = \rationals\ket{a,b}$.
		In contrast, for all $t \geq 1$, $cx^t \in C_0\graded(K_\bullet;\rationals)$ since $c$ is a vertex in $K_t$.

		\item 
		The filtered $0$-chain $\sigma_2(x) := (a + 2b - c)x^4$ 
		corresponds to the $1$-chain $a + 2b - c$ as an element of $C_0(K_4; \rationals)$ at index $4$
		and $\degh(\sigma_2(x)) = 4$, i.e.\ $\sigma_2(x)$ is homogeneous of degree $4$.

		\item 
		The filtered $0$-chain $a \in C_0\graded(K_\bullet;\rationals)$ generates the $\rationals[x]$-submodule $\rationals[x]\ket{a}$
		consisting of elements of the form $ax^t \in C_0\graded(K_\bullet;\rationals)$ with $t \in \nonnegints$.
		Observe that $\rationals[x]\ket{a}$ is a graded submodule of $C_0\graded(K_\bullet;\rationals)$.
	\end{enumerate}
\end{example}

Recall that simplicial chain groups $C_n(K;R)$ with coefficients in a PID $R$ are free $R$-modules.
We have a similar result for graded chain modules.

\begin{proposition}
	Let $\filt{K}$ be a simplicial filtration.
	For all $n \in \nonnegints$, $C_n\graded(\filt{K};\field)$ is free.
\end{proposition}
\begin{proof}
	Let $n \in \ints$. If $n \geq -1$, then $C_n(K_\bullet;\field)$ is trivial and therefore free.
	Assume $n \geq 0$. 
	To determine that $C_n\graded(\filt{K};\field)$ is free, it suffices to show that it is torsion-free as a $\field[x]$-module.
	Let $\sigma x^t \in C_n\graded(K_\bullet;\field)$. 
	Then, $\sigma \in C_n(K_t;\field)$. 
	Note that $C_n\graded(\filt{K};\rationals)$ is a $\field$-vector space and is torsion-free.
	If $x^s \cdot \sigma x^t = 0$ for some $s \in \nonnegints$, 
		then $i^{t+s,s}_\hash(\sigma) = 0$, contradicting \fref{lemma:chain-group-graded-mod-char}.
	Therefore, $C_n\graded(\filt{K};\field)$ is torsion-free and, therefore, free.
\end{proof}

Earlier in \fref{defn:standard-basis-on-chain-groups},
	we identified a natural choice for basis on the chain groups of a simplicial complex based on its equipped orientation.
Similarly, if $K_\bullet$ is a filtration on a simplicial complex $K$,
	an orientation on $K$ induces a natural choice of basis on its graded chain modules.
	We state this as a proposition below.

\begin{proposition}\label{prop:basis-graded-chain-module}
	Let $K_\bullet$ be a filtration of an oriented finite simplicial complex $K$.
	For each $n \in \nonnegints$,
		$C_n\graded(K_\bullet;\field)$ is free with the following homogeneous basis:
	\begin{equation}
		\basis{K}_n\graded := \Bigl\{ 
			\sigma_1 x^{t_1}, \sigma_2 x^{t_2}, \ldots, \sigma_m x^{t_m}
		\Bigr\}
		\quad\text{ with }\quad 
		t_i := \min\Bigl\{t \in \nonnegints : \sigma_i \in C_n(K_t;\field) \Bigr\} 
		\text{ for all } i \in \set{1, \ldots, m}
	\end{equation}
	where $K[n] = \set{\sigma_1, \ldots, \sigma_m}$ is the standard basis of $C_n(K_t;\field)$ induced by the orientation on $K$ (see Defn.\,\ref{defn:standard-basis-on-chain-groups}).
\end{proposition}
\begin{proof}
	Fix $n \in \nonnegints$. 
	Since $K$ is a finite simplicial complex, $K[n]$ is finite with $m := \card(K[n]) < \infty$
	and we can label the oriented $n$-simplices in $K[n]$ by $\set{\sigma_1, \ldots, \sigma_n}$.
	Let $i \in \set{1, \ldots, m}$.
	By \fref{defn:filtration}(iii), there must exist some $t \in \nonnegints$ such that $\sigma_i$ corresponds to an $n$-simplex in $K$
	and $\sigma_i \in C_n(K_t;\field)$.
	Since $\nonnegints$ is bounded below, 
		a minimal $t_i \in \nonnegints$ for $\sigma_i$ such that 
		$\sigma_i \in C_n(K_{t_i};\field)$ must exist.
	Then, $\sigma_i x^{t_i} \in C_n(K_{t_i};\field) x^{t_i} \subseteq C_n\graded(K_\bullet;\field)$ is homogeneous with degree $t_i$.
	Therefore, $\basis{K}_n\graded$, as given above, is well-defined as a set.
	
	Assume, without loss of generality, that $t_1 \leq t_2 \leq \cdots \leq t_m$, i.e.\ $\basis{K}_n\graded$ is indexed in order of increasing degree. 
	We need to show that $\basis{K}_n\graded$ generates $C_n\graded(K_\bullet;\field)$ and that $\basis{K}_n\graded$ is $\field[x]$-linearly independent.

	Fix $t \in \nonnegints$ and consider the vector subspace $C_n(K_t;\field) x^t$ of $C_n\graded(K_\bullet;\field)$.
	Let $r \in \set{1, \ldots, m}$ be maximum such that $t_r \leq t$.
	By minimality of the $t_i$'s, 
		$\sigma_i$ corresponds to an $n$-simplex in $K_t$ if and only if $i \in \set{1, \ldots, p}$.
	Let $K_t$ inherit the orientation on $K$ by restriction and 
	let $K_t[n]$ be the standard basis of $C_n(K_t;\field)$ induced by this new orientation.
	Then, $K_t[n] = \set{\sigma_1, \ldots, \sigma_p}$ and 
	the following set is a basis for $C_n(K_t;\field)$:
	\begin{equation}\label{eqn:basis-per-degree}
		\basis{B}(t) := \Bigl\{
			\sigma_1 x^{t}, \sigma_2 x^t, \ldots, \sigma_p x^{t}
		\Bigr\}
		= 
		\Bigl\{
			x^{t-t_1} \cdot \sigma_1 x^{t_1}, 
			x^{t-t-2} \cdot \sigma_2 x^{t_2},
			\ldots, 
			x^{t-t_p} \cdot \sigma_r x^{t_p}
		\Bigr\}
	\end{equation}
	Since $t-t_i \geq 0$ for each $i \in \set{1, \ldots, r}$,
		each element in $C_n(K_t;\field)$ is equal to a $\field[x]$-linear combination in $\basis{K}_n\graded$.
	Since $t \in \nonnegints$ is arbitrary, $\basis{K}_n\graded$ generates $C_n\graded(K_\bullet;\field)$. 

	Consider the $\field[x]$-linear dependence relation
	$f(x) := \sum_{i=1}^m f_i \cdot \sigma_i x^{t_i} = 0$ 
	for $\basis{K}_n\graded$
	with $f_1, \ldots, f_m \in \field[x]$.
	Recall that $C_n\graded(K_\bullet;\field) = \bigoplus_{t \in \nonnegints} C_n(K_t;\field) x^t$.
	For each $i \in \set{1, \ldots, m}$,
		let $\set{k_{i,t}}_{t \in \nonnegints}$ be such that
	$
		f_i(x) \cdot \sigma_i x^{t_i} = \sum_{t=0}^\infty k_{i,t} \cdot \sigma_i x^t
	$
	with $k_{i,t} \in \field$ for all $t \in \nonnegints$.
	Observe that, for all $i \in \set{1,\ldots,m}$,
		only finitely many $\set{k_{i,t}}_{t \in \nonnegints}$ are nonzero and 
		$k_{i,t} = 0$ if $t < t_i$.
	Then, $f(x)$ decomposes into
	\begin{equation*}
		f(x) = \sum_{i=1}^m f_i \cdot \sigma_i x^{t_i}
		= \sum_{i=1}^m \left(
			\sum_{t=0}^\infty k_{i,t} \cdot \sigma_i x^t
		\right)
		= \sum_{t=0}^\infty \left( \sum_{i=1}^m k_{i,t} \cdot \sigma_i x^t \right)
	\end{equation*}
	and the homogeneous component of $f(x)$ of degree $t \in \nonnegints$ is given by $\sum_{i=1}^m k_{i,t} \cdot \sigma_i x^t$.
	To solve linear dependence relation $f(x)=0$,
		it suffices to consider $\sum_{i=1}^m k_{i,t} \cdot \sigma_i x^{t}=0$ separately for each $t \in \nonnegints$.
	Let $t \in \nonnegints$.
	Observe that $\sum_{i=1}^m k_{i,t} \cdot \sigma_i x^{t} \in C_n(K_t;\field) x^t$ is a $\field$-linear combination in $\basis{B}(t)$, with $\basis{B}(t)$ as defined in Equation \ref{eqn:basis-per-degree}.
	Since $\basis{B}(t)$ is a basis for $C_n(K_t;\field) x^t$, 
		$k_{1,t} = k_{2,t} = \cdots = k_{m,t} = 0$.
	Since $t$ is arbitrary, $k_{i,t} = 0$ for all $t \in \nonnegints$.
	Then, for each $i \in \set{1, \ldots, m}$,
	\begin{equation*}
		f_i 
			= \sum_{t=0}^\infty k_{i,t} x^{t-t_i}
			= \sum_{t=t_i}^\infty k_{i,t} x^{t-t_i}
			= \sum_{t=t_i}^\infty 0 \cdot x^{t-t_i}
			= 0
	\end{equation*}
	Therefore, $f_1 = f_2 = \cdots = f_n$ and 
	$\basis{K}_n\graded$ is $\field[x]$-linearly independent.

	Since $\basis{K}_n\graded$ is $\field[x]$-linearly independent, generates $C_n\graded(K_\bullet;\field)$, 
	and consists of homogeneous elements,
	$\basis{K}_n\graded$ is a homogeneous $\field[x]$-basis for $C_n\graded(K_\bullet;\field)$. 
\end{proof}

In the example calculations presented later in \fref{chapter:matrix-calculation},
we often use the basis $\basis{K}_n\graded$ as a starting point.
For convenience, we provide a name for this basis below.

\begin{definition}\label{defn:standard-basis-graded-chain}
	Let $K_\bullet$ be a filtration of an oriented finite simplicial complex $K$.
	For each $n \in \nonnegints$, 
		let $\basis{K}_n\graded$, as denoted in \fref{prop:basis-graded-chain-module} and ordered first in increasing degree then by lexicographic order of the vertices of $K$,
		be the \textbf{standard ordered basis}
		of $C_n\graded(K_\bullet;\rationals)$
		\textbf{induced by the orientation on $K$}.
\end{definition}
\remark{
	Other authors may use a different ordering
	e.g.\ \cite{matrixalg:zomorodian} orders $\basis{K}_n\graded$ in \textit{decreasing} degree instead.
	We introduced the notion for standard \textit{ordered} basis here for convenience, as we will deal with coordinate matrices later in \fref{chapter:matrix-calculation}.
}

In the example below, we identify the basis $\basis{K}_n\graded$ corresponding to the given graded chain modules.

\begin{example}
	Let $K$ and $\filt{K}$ be as given in \fref{ex:pershom-one}, illustrations of which are copied below for convenience, and orient $K$ by the vertex order $(a,b,c,d)$.
	\begin{center}
		\baselineCenter{\includegraphics[height=1.05in]{zomfig1/zom-simp.png}}
		\quad
		\baselineCenter{\includegraphics[height=1.1in]{zomfig1/zom1-cleaned.png}}
	\end{center}
	\vspace{3pt}
	\noindent 
	Provided below are the ordered bases $\basis{K}_n$ for $C_n\graded(K_\bullet;\rationals)$ for $n=0,1,2$, as described in \fref{prop:basis-graded-chain-module}.
	\begin{align*}
		\basis{K}_0\graded = (a,b,cx,dx), \qquad
		\basis{K}_1\graded = (abx, bcx, adx^2, cdx^2, acx^3), \qquad 
		\basis{K}_2\graded = (abc x^4, acd x^5)
	\end{align*}
	Therefore, we can describe $C_n\graded(K_\bullet;\rationals)$ for $n=0,1,2$ as follows:
	\begin{align*}
		C_0\graded(K_\bullet;\rationals)
		= \rationals[x]\Bigl\lket 
			a,b,cx,dx
		\Bigr\rket ,
		\quad 
		C_1 = \rationals[x]\Bigl\lket 
			abx, bcx, adx^2, cdx^2, acx^3
		\Bigr\rket ,
		\quad 
		C_2 = 
		\rationals[x]\Bigl\lket 
			abc x^4, acd x^5
		\Bigr\rket
	\end{align*}
	Standard ordered bases allow us to represent filtered $n$-chains using coordinate matrices and vectors. 
	Following the ordering of $\basis{K}_0\graded$ and 
	$\basis{K}_1\graded$ denoted above, 
		we listed some examples of these coordinate vectors below.
	\vspace{-5pt}
	\begin{longtable}{
		L @{\quad$\implies$\quad} L
	} 
	\sigma_1 = 2ax + 3bx^3 = (a)(2x) + (b)(2x^3)
		& [\sigma_1] = \hspace{3pt}
		\begin{NiceArray}{>{\color{gray}}c @{\hspace{13pt}}c}
			a & 2x \\
			b & 3x^3 \\
			cx & 0 \\
			dx & 0
		\CodeAfter
			\SubMatrix({1-2}{4-2})[left-xshift=3pt, right-xshift=3pt]
		\end{NiceArray}
	\\[22pt]
	\begin{aligned}
		\sigma_2 &= (ab+bc-ac)x^4 \\
		&= (abx)(x^3) + (bcx)(x) + (acx^3)(-x)
	\end{aligned}
		& [\sigma_2] = \hspace{3pt}
		\begin{NiceArray}{>{\color{gray}}c @{\hspace{13pt}}c}
			abx & 2x \\
			bcx & 3x^3 \\
			adx^2 & 0 \\
			cdx^2 & 0 \\
			acx^3 & 0
		\CodeAfter
			\SubMatrix({1-2}{5-2})[left-xshift=3pt, right-xshift=3pt]
		\end{NiceArray}
	\\[30pt]
	\begin{aligned}
		\sigma_3 &= (cd-ad+ac)x^5 \\
		&= (adx^2)(-x^3) + (cdx^2)(x^3) + (acx^3)(x^2)
	\end{aligned}
		& [\sigma_3] = \hspace{3pt}
		\begin{NiceArray}{>{\color{gray}}c @{\hspace{13pt}}c}
			{abx}	& 0 \\
			{bcx}	& 0 \\
			{adx^2}	& -x^3 \\ 
			{cdx^2}	& x^3 \\
			{acx^3}	& x^2 
		\CodeAfter
			\SubMatrix({1-2}{5-2})[left-xshift=3pt, right-xshift=3pt]
		\end{NiceArray}
	\end{longtable}
	\vspace{-\baselineskip}
\end{example}

The functorial nature of the simplicial chain complex construction also allows us to extend the simplicial boundary maps to the case of persistence modules.
More specifically, the collection of simplicial boundary maps form a persistence morphism between the filtered chain groups.
We state this in more detail below.

\begin{corollary}\label{cor:filtered-boundary-map}
	Let $K_\bullet$ be a simplicial filtration and let $n \in \ints$.
	For each $t \in \nonnegints$, let $\boundary_n^{\hspace{0.5pt}t}: C_n(K_t;\field) \to C_{n-1}(K_t;\field)$ be the $n$\th simplicial boundary map of $K_t$.
	Then, the collection $\set{\boundary_n^{\hspace{0.5pt}t}}_{t \in \nonnegints}$ determines a persistence morphism 
	$\boundary_n^\sbullet: C_n(K_\bullet;\field) \to C_{n-1}(K_\bullet;\field)$.
\end{corollary}
\begin{proof}
	Fix $n \in \nonnegints$ and let $t,s \in \nonnegints$ with $t \leq s$.
	By \fref{prop:induced-maps-on-chains-respect-boundary}, 
		$i^{s,t}_{n-1,\hash} \circ \boundary_n^t = \boundary_n^s \circ i^{s,t}_{n,\hash}$ for all $t,s \in \nonnegints$ with $t \leq s$
	with $i^{s,t}_{n,\hash}: C_n(K_t;\field) \to C_n(K_s;\field)$ 
	and $i^{s,t}_{n-1,\hash}: C_{n-1}(K_t;\field) \to C_{n-1}(K_s;\field)$ referring to structure maps 
	of $C_n(K_\bullet;\field)$ and $C_{n-1}(K_\bullet;\field)$ respectively.
	That is, the following diagram commutes:
	\vspace{-5pt}
	\begin{displaymath}
	\begin{tikzcd}[row sep=2em]
		C_n(K_t;R) \arrow[r, "\boundary_n^{\hspace{0.5pt}t}"] 
					\arrow[d, "i_{n,\hash}^{s,t}" swap]
			&[2em] C_{n-1}(K;R) 
					\arrow[d, "i_{n-1,\hash}^{s,t}"]
		\\[1em]
		C_n(K_s;R) \arrow[r, "\boundary_n^{\hspace{0.5pt}s}"] & C_{n-1}(L;R)
	\end{tikzcd}
	\end{displaymath}
	Therefore, $\boundary_n^\sbullet = (\boundary_n^{\hspace{0.5pt}t})_{t \in \nonnegints}$ is a persistence morphism by \fref{defn:persmod-cat}. 
\end{proof}

We name these persistence morphisms, along with the corresponding graded homomorphism given by application of $\togrmod(-)$, below.

\begin{definition}\label{defn:filtered-n-boundary-morphism}
	For each $n \in \ints$,
		define the 
		$n$\th \textbf{filtered boundary morphism} 
		$\boundary_n^\sbullet$ of a simplicial filtration $K_\bullet$ 
		to be the persistence morphism 
		$\boundary_n^\sbullet: C_n(K_\bullet;\field) \to C_{n-1}(K_\bullet;\field)$
		given by $\boundary_n^\sbullet = (\boundary_n^{\hspace{0.5pt}t})_{t \in \nonnegints}$
		where $\boundary_n^{\hspace{0.5pt}t}: C_n(K_t;\field) \to C_{n-1}(K_t;\field)$ is the $n$\th simplicial boundary map of $K_t$.
\end{definition}
\negativespacer
\begin{definition}\label{defn:graded-n-boundary-map}
	Let $K_\bullet$ be a filtration of a simplicial complex.
	For each $n \in \ints$, define the 
	$n$\th \textbf{graded boundary morphism} $\boundary_n\graded$
	of $K_\bullet$ to be the graded $\field[x]$-module homomorphism 
	$\boundary_n\graded: C_n\graded(K_\bullet;\field) \to C_{n-1}\graded(K_\bullet;\field)$ given by 
	$\boundary_n\graded := \togrmod\bigl( \boundary_n^\sbullet \bigr)$.
	A \textbf{filtered} $n$-\textbf{cycle}
	is an element of $\ker(\boundary_n\graded)$ 
	and a \textbf{filtered} $n$-\textbf{boundary}
	that of $\im(\boundary_{n+1}\graded)$.
	If $K$ is oriented and finite,
		define the $n$\th \textbf{graded boundary matrix} $[\boundary_n\graded]$ to be the matrix of $\boundary_n\graded$ relative to the standard ordered bases $\basis{K}_n$ and $\basis{K}_{n-1}$.
\end{definition}

Below, we identify a useful characterization of filtered and graded boundary morphisms of a filtration $K_\bullet$ relative to the simplicial boundary morphism of a simplicial complex $K$.

\begin{corollary}
	Let $K_\bullet$ be a filtration of a simplicial complex $K$ and let $n \in \nonnegints$.
	\begin{enumerate}
		\item 
		For all $t \in \nonnegints$,
			$\boundary_n^{\w t}(\sigma) = \boundary_n(\sigma)$ for all $\sigma \in C_n(K_t;\field)$ and $t \in \nonnegints$.

		\item 
		The $n$\th graded boundary map satisfies 
			$
			\boundary_n\graded\bigl( \sigma \bigr)
			= \sum_{t=0}^\infty \boundary_n(\sigma_t) x^t
			$
		for all $\sigma = \sum_{t=0}^\infty \sigma_t x^t \in C_n(K_\bullet;\field)$.
	\end{enumerate}
\end{corollary}
\begin{proof}
	Fix $n \in \nonnegints$. 
	For (i):
	Since $\boundary_n^\sbullet: C_n(K_\bullet;\field) \to C_{n-1}(K_\bullet;\field)$ is a persistence morphism, 
		we have that 
		for all $t \in \nonnegints$
		and $\sigma \in C_n(K_t;\field)$,
	\begin{equation*}
		\boundary_n^{\w t}(\sigma)
		= \Bigl( \boundary_n^{\w t} \circ i_{n,\hash}^{[t]} \Bigr)(\sigma)
		= \Bigl( i_{n-1,\hash}^{[t]} \circ \boundary_n \Bigr)(\sigma)
		= \Bigl( \id_{C_{n-1}(K_t;\field)} \circ\ \boundary_n \Bigr)(\sigma)
		= \boundary_n(\sigma)
	\end{equation*}
	where $i^{[t]}: K_t \to K$ refers to the inclusion map.
	For (ii): For all $\sigma = \sum_{t=0}^\infty \sigma_t x^t \in C_n\graded(K_\bullet;\field)$,
	\begin{equation*}
		\boundary_n\graded(\sigma)
		= \boundary_n\left(
			\,\sum_{t=0}^\infty \sigma_t x^t
		\right)
		\stackrel{(\star)}{=} \sum_{t=0}^\infty \boundary_n^{\w t}(\sigma_t) x^t 
		= \sum_{t=0}^\infty \boundary_n(\sigma_t) x^t
	\end{equation*}
	with $(\star)$ given by definition of the morphism assignment of $\togrmod$.
\end{proof}

The corollary above tells us that, when evaluating filtered $n$-chains against graded boundary morphisms, 
	we can basically ignore the added $x^t$ in the notation and determine the boundary of an oriented $n$-simplex as usual.
We have an example below.

\begin{example}
	Let $K$ and $\filt{K}$ be as given in \fref{ex:pershom-one} and equip $K$ with the orientation by $(a,b,c,d)$.
	Let $\sigma_1 := (ab+bc-ac)x^4 \in C_1\graded(K_\bullet;\rationals)$. 
	Then, $\sigma_1$ is a filtered $1$-cycle by the following calculation:
	\begin{equation*}
		\boundary_1\graded(\sigma_1)
		= \boundary_1\graded\Bigl(
			(ab + bc - ac)x^4
		\Bigr)
		= \boundary_1(ab+bc-ac) x^4 
		= \big(
			(b-a) + (c-b) - (c-a)
		\big) x^4 = (0)x^4 = 0
	\end{equation*}
	Additionally, $\sigma_1$ is a filtered $1$-boundary since 
	$\sigma_1 \in \im(\boundary_2\graded)$ by the following calculation:
	\begin{equation*}
		\boundary_2\graded\Bigl(
			abc x^4
		\Bigr)
		= \boundary_{2}(abc)x^4 
		= \big(
			bc - ac + ab
		\big)x^4
		= \sigma_1
	\end{equation*}
	Let $\sigma_2 := (ab + bc - ac)x^3 \in C_1\graded(K_\bullet;\rationals)$. 
	Then, $\boundary_1\graded(\sigma_2) = \boundary_1(ab+bc-ac)x^3 = 0$ and $\sigma_2$ is a filtered $1$-cycle.
	However, $\sigma_2$ is not a filtered $1$-boundary since 
	$abc$ is not in $K_3$ and $abc x^3 \not\in C_2\graded(K_\bullet;\rationals)$.
\end{example}

We need to verify that the $\ints$-indexed collection of filtered chain modules $C_n(K_\bullet;\rationals)$ and filtered boundary morphisms $\boundary_n^\sbullet: C_n(K_\bullet;\rationals) \to C_{n-1}(K_\bullet;\rationals)$ does indeed determine a persistence complex. 
We also need to do the same for the graded chain modules and graded boundary morphisms.

\begin{proposition}
	Let $K_\bullet$ be a simplicial filtration.
	Then, 
		$\boundary_{n-1}^\sbullet \circ \boundary_{n}^\sbullet = 0_\sbullet$
		and 
		$\boundary_{n-1}\graded \circ \boundary_{n}\graded = 0$ for all $n \in \ints$.
\end{proposition}
\begin{proof}
	Fix $n \in \nonnegints$. 
	For the persistence modules case: 
	For all $t \in \nonnegints$, $
		(\boundary_{n-1}^\sbullet \circ \boundary_{n}^\sbullet)_t 
		= \boundary_{n-1}^{\w t} \circ \boundary_{n}^{\w t}
		= 0
	$ as $\field$-linear maps.
	Therefore, $\boundary_n^\sbullet \circ \boundary_{n-1}^\sbullet = 0_\sbullet$,
	with $0_\sbullet$ denoting the zero persistence morphism.
	For the graded module case:
	For all $\sigma = \sum_{t=0}^\infty \sigma_t x^{t} \in C_{n}\graded(K_\bullet;\field)$, we have the following:
	\begin{equation*}
		\Bigl( \boundary_{n-1}\graded \circ \boundary_{n}\graded \Bigr)(\sigma)
		= \boundary_{n-1}\graded\left(
			\,\sum_{t=0}^\infty \boundary_{n}(\sigma_t) x^t
		\right)
		= \sum_{t=0}^\infty 
			\Bigl( \boundary_{n-1} \circ \boundary_n \Bigr) x^t
		= \sum_{t=0}^\infty 0 \cdot x^t = 0
	\end{equation*}
	Therefore, $\boundary_{n-1}\graded \circ \boundary_n\graded = 0$.
\end{proof}

Now, we define the \textit{simplicial persistence complex} of a simplicial filtration, along with its corresponding chain complex of graded modules.

\begin{definition}\label{defn:simplicial-persistence-complex}
	Define the \textbf{simplicial persistence complex} $C_\ast(K_\bullet;\field)$ 
	and the \textbf{simplicial graded chain complex} $C_\ast(K_\bullet;\field)$
	of a simplicial filtration $K_\bullet$ with coefficients in $\field$
	as follows:
	\begin{equation*}
		C_\ast(K_\bullet;\field) := \Bigl(\,
			C_n(K_\bullet);\field, \boundary_n^\sbullet
		\,\Bigr)_{n \in \ints} 
		\quad\text{ and }\quad 
		C_\ast\graded(K_\bullet;\field) := \Bigl(\,
			C_n\graded(K_\bullet);\field, \boundary_n\graded
		\,\Bigr)_{n \in \ints}
	\end{equation*}
\end{definition}

The simplicial persistence complex $C_\ast(K_\bullet;\field)$ can be visualized as the following sequence of persistence modules and persistence morphisms:
\begin{equation*}
\begin{tikzcd}[column sep=5em]
	\cdots \arrow[r, "\boundary_{n+2}^\sbullet"]
	& \redmath{C_{n+1}(\filt{K};\field)} 
		\arrow[r, "\boundary_{n+1}^\sbullet"]
	& \bluemath{C_{n}(\filt{K};\field)} 
		\arrow[r, "\boundary_{n}^\sbullet"]
	& \greenmath{C_{n-1}(\filt{K};\field)} 
		\arrow[r, "\boundary_{n-1}^\sbullet"]
	& \cdots
\end{tikzcd}
\end{equation*}
The change in perspective from having the index $t \in \nonnegints$ of a filtration $K_\bullet$ take precedence to that for the index $n \in \nonnegints$ of the chain complex $C_\ast$ can be visualized using the following commutative diagram.
In the diagram below, each row fixes the index $t \in \nonnegints$ of the filtration and represents the simplicial chain complex $C_\ast(K_t;\field)$ of $K_t$.
Similarly, each column fixes the dimension $n \in \ints$ and corresponds to a filtered chain module $C_n(K_\bullet;\field)$.
For clarity, we suppressed the index $t$ in $i_\hash^{s,t}$ 
	and highlighted the vector spaces of 
	$\redmath{C_{n+1}(\filt{K};\field)}$ in \redtag,
	those of $\bluemath{C_{n}(\filt{K};\field)}$ in \bluetag,
	and $\greenmath{C_{n-1}(\filt{K};\field)}$ in \greentag.
\vspace{-10pt}
\begin{equation*}
\begin{tikzcd}[column sep=5em, row sep=2.1em]
	{}
		& \vdots
		& \vdots 
		& \vdots
		& {}
	\\
	\cdots \arrow[r, "\boundary_{n+2}"]
		& \redmath{C_{n+1}(K_{t+1};\field)} 
			\arrow[r, "\boundary_{n+1}"] \arrow[u, "i_\hash"]
		& \bluemath{C_{n}(K_{t+1};\field)} 
			\arrow[r, "\boundary_{n}"] \arrow[u, "i_\hash"]
		& \greenmath{C_{n-1}(K_{t+1};\field)} 
			\arrow[r, "\boundary_{n-1}"] \arrow[u, "i_\hash"]
		& \cdots
	\\
	\cdots \arrow[r, "\boundary_{n+2}"]
		& \redmath{C_{n+1}(K_t;\field)} 
			\arrow[r, "\boundary_{n+1}"] \arrow[u, "i_\hash"]
		& \bluemath{C_{n}(K_t;\field)}
			\arrow[r, "\boundary_n"] \arrow[u, "i_\hash"]
		& C_{n-1}(K_t;\field) 
			\arrow[r, "\boundary_{n-1}"] \arrow[u, "i_\hash"]
		& \cdots
	\\
	\cdots \arrow[r, "\boundary_{n+2}"]
		& \redmath{C_{n+1}(K_{t-1};\field)} 
			\arrow[r, "\boundary_{n+1}"] \arrow[u, "i_\hash"]
		& \bluemath{C_{n}(K_{t-1};\field)} 
			\arrow[r, "\boundary_n"] \arrow[u, "i_\hash"]
		& \greenmath{C_{n-1}(K_{t-1};\field)}
			\arrow[r, "\boundary_{n-1}"] \arrow[u, "i_\hash"]
		& \cdots
	\\[-10pt]
	{}
		& \vdots \arrow[u, "i_\hash"]
		& \vdots \arrow[u, "i_\hash"]
		& \vdots \arrow[u, "i_\hash"]
		& {}
\end{tikzcd}
\end{equation*}
\vspace{-15pt}

Finally, we present the result that allows us to use both the simplicial persistence complex and the simplicial graded complex in our calculations.

\begin{proposition}\label{prop:simp-pers-hom-graded-to-persmod}
	Let $K_\bullet$ be a simplicial filtration.
	For all $n \in \ints$, 
		we have that 
		\begin{equation*}
			H_n(K_\bullet;\field) 
			\cong 
			H_n\pers\Bigl( C_\ast(K_\bullet;\field) \Bigr)
			\cong 
			\Bigl(
				\topersmod \circ H_n\graded
			\Bigr) \Bigl( C_\ast\graded(K_\bullet;\field) \Bigr)
			\quad\text{ as persistence modules over $\field$ }
		\end{equation*}
	where $H_n\pers: \catchaincomplex{\catpersmod} \to \catpersmod$
	and $H_n\graded: \catchaincomplex{\catgradedmod{\field}} \to \catgradedmod{\field}$ refer to the chain homology functors on persistence complexes and graded chain complexes respectively.
\end{proposition}
\begin{proof}
	Recall that $H_n(-;\field) = H_n\pers \mathrel{\circ} C_\ast(-;\field)$ as functors $\catsimp \to \catvectspace$, where $H_n(-;\field)$ and $C_\ast(-;\field)$ are the simplicial homology and simplicial chain complex functors respectively.
	Since functor composition is associative when defined, we have the following:
	\begin{equation*}
		H_n(K_\bullet;\field)
		= H_n(-;\field)(K_\bullet) 
		= \Bigl( H_n\pers \circ C_\ast(-;\field) \Bigr)(K_\bullet)
		= H_n\pers\Bigl( C_\ast(K_\bullet;\field) \Bigr)
	\end{equation*}
	As stated in \fref{prop:catequiv-preserves-chains}, $\togrmod$ and $\topersmod$ preserve chain complexes and chain homology.
	In particular, both functors commute with the chain homology functors.
	Then, 
	\begin{align*}
		H_n\pers \circ C_\ast(K_\bullet;\field)
		&= \id_{\catpersmod} \mathrel{\circ} H_n\pers \mathrel{\circ} C_\ast(K_\bullet;\field)
		= \topersmod \mathrel{\circ} \togrmod \mathrel{\circ} H_n\pers \mathrel{\circ} C_\ast(K_\bullet;\field)
		\\[2pt]
		&= \topersmod \Bigl(
			\togrmod \circ H_n\pers \circ C_\ast(K_\bullet;\field)
		\Bigr)
		= \topersmod \Bigl(
			H_n\graded \circ \togrmod \circ C_\ast(K_\bullet;\field)
		\Bigr)
		\\[2pt]
		&= \Bigl(
			\topersmod \circ H_n\graded
		\Bigr)\Bigl( C_\ast\graded(K_\bullet;\field) \Bigr)
	\end{align*}

	\vspace{-3\baselineskip}
\end{proof}

The proposition above is a fundamental result behind the derivation of the \textit{matrix reduction algorithm for persistent homology}, as presented in \cite{matrixalg:zomorodian}, since it allows us to calculate at the level of graded modules and return back to that of persistence modules post-calculation.
For convenience, we name the chain homology of the graded persistence complex.

\begin{definition}\label{defn:graded-homology-module}
	For each $n \in \ints$,
		define the $n$\th \textbf{graded homology module}
		$H_n\graded(K_\bullet;\field)$ with coefficients in a field $\field$
		of a simplicial filtration $K_\bullet$ 
		as $H_n\graded(K_\bullet;\field) := H_n \circ C_\ast\graded(K_\bullet;\field)$,
		i.e.\ the $n$\th chain homology of the simplicial graded chain complex.
\end{definition} 
\clearpage


\onlyifstandalone{
	\setcounter{tocdepth}{1}
	\renewcommand*\contentsname{Table of Contents}
	\pdfbookmark{\contentsname}{toc}
	\begin{spacing}{1.2} 
		\tableofcontents 
	\end{spacing}
	\clearpage
}

\onlyifstandalone{\setcounter{chapter}{3}}  
\chapter{Calculation by Matrices}
\label{chapter:matrix-calculation}

The paper \textit{Calculating Persistent Homology} \cite{matrixalg:zomorodian} by Afra Zomorodian and Gunnar Carlsson 
describes how the \textit{matrix reduction algorithm for persistent homology} is based on a method of calculating {invariant factor decompositions} of finitely-generated $\field[x]$-modules using matrices over $\field[x]$.
Certain assumptions on the $\field[x]$-modules in question allow this calculation to be simulated using matrices over $\field$.
In this chapter, we frame the discussion in \cite{matrixalg:zomorodian} relative to the theoretical foundation established in the previous chapters.

Fix a field $\field$ and let $K_\bullet$ be a simplicial filtration of a finite simplicial complex $K$.
As established in \fref{section:simplicial-persistent-homology},
	the $n$\th persistent homology module $H_n(K_\bullet;\field)$ of $K_\bullet$ with coefficients in $\field$ 
	can be calculated
	using the following persistence isomorphism relation:
	\begin{equation*}
		H_n(K_\bullet;\field) \cong 
		\Bigl(
			\topersmod \circ H_n \circ \togrmod
		\Bigr)\!\Bigl(
			C_\ast(K_\bullet;\field)
		\Bigr)
	\end{equation*}
where 
$C_\ast(K_\bullet;\field) = (C_m(K_\bullet;\field); \boundary_m^\sbullet)_{m \in \ints}$ is the simplicial persistence complex of $K_\bullet$,
$H_n: \catchaincomplex{\catgradedmod{\field}} \to \catgradedmod{\field}$ is the $n$\th chain homology functor on graded chain complexes,
and $\togrmod: \catpersmod \to \catgradedmod{\field}$
and $\topersmod: \catgradedmod{\field} \to \catpersmod$ refer to the category equivalence discussed in \fref{section:cat-equiv-graded-modules}.

Following the arguments in \cite{matrixalg:zomorodian},
	the \textit{matrix reduction algorithm for persistent homology}
	calculates persistent homology at the level of graded modules.
In particular,
	the algorithm determines the \textit{graded invariant factor decomposition} of the following graded $\field[x]$-module by matrix reduction:
	\begin{equation*}
		H_n\graded(K_\bullet;\field) := \Bigl(
			H_n \circ \togrmod
		\Bigr)\!\Bigl(
			C_\ast(K_\bullet;\field)
		\Bigr)
	\end{equation*}
The resulting decomposition then determines the \textit{interval decomposition} of the persistence module $H_n(K_\bullet;\field)$, similarly as in \fref{lemma:intmods-gradedmods-corr} and \fref{cor:interval-decomp-from-structure-theorem}.

In this chapter, 
we take a more general view and explore how matrix reduction can be used to find graded invariant factor decompositions of the $n$\th chain homology of graded chain complexes.
Let $R$ be a PID and $\field$ be a field.
This chapter is structured as follows:
\begin{enumerate}[chapterdecompositionVTWO]
	\item[In \textbf{\fref{section:matrix-calculation-of-IFDs}.}]
	\textbf{The Structure Theorem and Smith Normal Decompositions}

	\noindent
	We consider the Structure Theorem for Finitely Generated Modules over a PID $R$ in the category $\catmod{R}$, 
	i.e.\ disregarding grading (if it exists), 
	and discuss how \textit{invariant factor decompositions} of said modules can be calculated using presentations and a matrix factorization called Smith Normal Decomposition (SND).
	Note that our examples in this section use $R = \ints$ for comparison to the graded case in the later sections.

	\item[In \textbf{\fref{section:snd-on-ungraded-chain-complexes}.}]
	\textbf{Matrix Calculation of Homology of Ungraded Chain Complexes}
	
	\noindent
	We consider chain complexes $C_\ast = (C_n, \boundary_n)_{n \in \ints}$ 
	in $\catchaincomplex{\catmod{R}}$
	such that for all $n \in \ints$, $C_n$ is a free $R$-module of finite-rank.
	For each $n \in \ints$, we present an existence result involving the existence of a decomposition of $C_{n}$ into three free direct summands:
	\begin{equation*}
		C_{n} \cong K_{n}^\text{free} \oplus K_{n}^\text{tor} \oplus \frac{C_n}{\ker(\boundary_n)}
	\end{equation*}
	such that the free component and torsion component of $H_n(C_\ast)$ are given by $\text{F}(H_n(C_\ast)) \cong K_{n}^\text{free}$ and $\text{T}(H_n(C_\ast)) \cong K_n^\text{tor} \bigmod \im(\boundary_{n+1})$ respectively.
	We also discuss how these components can be determined from specific SNDs of the matrices of the differentials $\boundary_{n+1}: C_{n+1} \to C_n$ and $\boundary_n: C_n \to C_{n-1}$.

	\item[In \textbf{\fref{section:calculation-graded-ifds}.}]
	\textbf{The Graded Structure Theorem and SNDs in the Graded Case}
	
	\noindent
	We present the
	Graded Structure Theorem for Finitely-Generated $\field[x]$-modules in the category $\catgradedmod{\field}$.
	We also discuss how this theorem can be considered a special case of the Structure Theorem in $\catmod{\field[x]}$
	and 
	how the method of calculating invariant factor decompositions by SNDs, limited to \textit{graded presentations},
	can be used to find \textit{graded} invariant factor decompositions of graded $\field[x]$-modules.
	For convenience, we use $\field = \rationals$ for our examples.

	\item[In \textbf{\fref{section:matrix-reduction-of-graded-matrices}.}]
	\textbf{Matrix Reduction of Graded Matrices}
	
	We consider matrices of homomorphisms in graded presentations, which we call \textit{graded matrices} for brevity.
	We discuss how specific matrix operations, i.e.\ elementary permutations, elementary dilations, and non-trivial elimination operations in matrices over $\field[x]$, preserve the homogeneity of graded matrices.

	We also consider the simplicial filtration $K_\bullet$ presented in \fref{ex:pershom-one}, i.e.\ that in \cite[Figure 1]{matrixalg:zomorodian}.
	An illustration of $K_\bullet$ is copied below for convenience.

	\begin{center}
		\baselineCenter{\includegraphics[height=1.0in]{zomfig1/zom-simp.png}}
		\quad
		\baselineCenter{\includegraphics[height=1.05in]{zomfig1/zom1-cleaned.png}}
	\end{center}
	\vspace{5pt}
	\noindent
	In particular, we perform matrix reduction on 
	$[\boundary_1\graded]$ and $[\boundary_2\graded]$, the matrices of the graded boundary maps 
	$\boundary_1\graded: C_1\graded(K_\bullet;\rationals) \to C_0\graded(K_\bullet;\rationals)$
	and 
	$\boundary_2\graded: C_2\graded(K_\bullet;\rationals) \to C_1\graded(K_\bullet;\rationals)$
	respectively,
	and discuss why the homogeneity of said matrices are preserved after the listed matrix operations.

	\item[In \textbf{\fref{section:snd-algorithm-in-the-graded-case}.}]
	\textbf{An Ungraded SND Algorithm in the Graded Case}
	
	\noindent
	We present a general algorithm for finding SNDs of graded matrices, adapted from an algorithm for finding SNDs of matrices over a PID $R$, 
	and discuss why the SNDs resulting from this algorithm can be used to determine \textit{graded} invariant factor decompositions.

	We also use this algorithm to determine the graded invariant factor decomposition of the graded $\rationals[x]$-module 
	$C_0\graded(K_\bullet;\rationals)$
	of $K_\bullet$ 
	and use said result, along with \fref{lemma:intmods-gradedmods-corr} and \fref{cor:interval-decomp-from-structure-theorem},
	to determine the interval decomposition of 
	the persistent homology module
	$C_0(K_\bullet;\rationals)$ of $K_\bullet$ in dimension $n=0$.

	\item[In \textbf{\fref{section:matrix-graded-chain-complex}.}]
	\textbf{Matrix Calculation of Homology of Graded Chain Complexes}
	
	We consider graded chain complexes $C_\ast = (C_n, \boundary_n)_{n \in \ints}$ in $\catchaincomplex{\catgradedmod{\field}}$ such that $C_n$ is a free graded $\field[x]$-module of finite rank.
	We briefly discuss why the decomposition of $C_n$ into the free summands $K_n^\text{free}$, $K_n^\text{tor}$, and $C_n \bigmod \ker\boundary_n$ and how the method of finding ungraded invariant factor decompositions of chain homology in $\catchaincomplex{\catmod{\field[x]}}$ extends to the graded case. 
\end{enumerate}
Some of the notation used involving matrices and matrix reduction are identified \fref{appendix:matrix-theory}.
In particular, we bring emphasis to the elementary matrices 
$\elswap{k_1, k_2}$, $\eldilate{k,\mu}$, and $\eladd{k_j, k_i\,;\alpha} \in \GL(n,R)$ described in \fref{defn:elementary-matrices}. 

\clearpage 


 
\section{The Structure Theorem and Smith Normal Decompositions}
\label{section:matrix-calculation-of-IFDs}

The Structure Theorem for Finitely Generated Modules over a PID, which we call the \textbf{Structure Theorem} in this paper for convenience,
proves the existence and uniqueness of invariant factor decompositions of certain modules.
In this section, we discuss one of the proofs of the Structure Theorem, as presented in \cite{algebra:dummit},
and how the method of calculating invariant factor decompositions can re-stated as a matrix calculation using a matrix factorization called a \textit{Smith Normal Decomposition}, defined later in this section in \fref{defn:smith-normal-decomposition}.
To start, we provide a statement of the Structure Theorem below.

\begin{statement}{Theorem}\label{thm:structure-theorem}
	\textbf{The Structure Theorem for Finitely Generated Modules over a PID.}
 
	Let $M$ be a finitely generated module over some PID $R$.
	There exists a set of \textbf{invariant factors} $\set{d_i}_{i=1}^n$ of non-invertible elements $d_i \in R$ with divisibility relation $d_1 \mid d_2 \mid \cdots \mid d_n$ such that $M$ is isomorphic to a direct sum of cyclic modules as follows:
	\begin{equation*}
		M 
			\cong \bigoplus_{i=1}^n R \bigmod (d_i)
			= R \bigmod (d_1) 
				\oplus R \bigmod (d_2) 
				\oplus \cdots 
				\oplus R \bigmod (d_n)
	\end{equation*}
	This direct sum is called the \textbf{invariant factor decomposition} of $M$ and is unique up to isomorphism.
	The invariant factors ${d_i}$ are unique up to multiplication by units. 

\end{statement}
\remark{
	We refer to \cite[Theorem 12.5]{algebra:dummit} and \cite[Theorem 12.9]{algebra:dummit} for the proofs of the existence and uniqueness claims respectively. 
	Note that a key characteristic of a PID is that all of its ideals can be generated by a single element.
	Given $d \in R$ with $R$ a PID, 
		we write $(d)$ to refer to the ideal generated by $d$ by 
		$(d) := Rd = \set{rd : r \in R}$.
}

The divisibility relation 
$d_1 \divides d_2 \divides \cdots \divides d_n$ 
on the invariant factors $\set{d_i}$ in the Structure Theorem
is sometimes stated in terms of proper ideals, wherein each $(d_i)$ must be a proper ideal of $R$, i.e.\ $(d_i) \neq R$, and we have the following decreasing sequence of ideals:
\begin{equation*}
	(d_1) \supseteq (d_2) \supseteq \cdots \supseteq (d_n)
\end{equation*}
Here, we use the term \textit{decreasing} relative to the subset relation.
Given $a, b \in R$ such that $a$ divides $b$, i.e.\ $a \divides b$, there must exist $q \in R$ such that $aq = b$.
Then, any element $rb \in (b)$ with $r \in R$ must also be in $(ra)$ since $rb = r(aq) = (rq)a$ and $rq \in R$.
Therefore, $(b) \subseteq (a)$ or equivalently, $(a) \supseteq (b)$.

Since the ideals $(d_i)$ must be proper ideals, this means that the invariant factors $d_i$ cannot be invertible.
More specifically, if $a \in R$ is invertible, then $(a) = (a\inv a) = (1) = R$ where $a\inv \in R$ refers to the multiplicative inverse of $a$ in $R$.
If $d_i \in R$ were invertible, then the summand $R \bigmod (d_i) = R \bigmod R$ would be the trivial module and, therefore, can be removed from the direct sum without invalidating the isomorphism.

We also want to emphasize that the Structure Theorem allows the invariant factors $d_i$ to be zero.
Note that, relative to the decreasing sequence of ideals, this means that $d_i = 0$ would occur at the end of the sequence since $(0) \subseteq (a)$ for any $a \in R$.
Consequently, some references prefer to only consider nonzero $d_i$'s for the invariant factors and state \fref{thm:structure-theorem} as follows:
\begin{equation*}
	M \cong R^f \oplus R \bigmod (d_1) 
	\oplus R \bigmod (d_2) 
	\oplus \cdots 
	\oplus R \bigmod (d_r)
\end{equation*}
with $f+r=n$ and $d_i \neq 0$ for all $i \in \set{1, \ldots, r}$.
In this case, we call $R^f$ the \textit{free component} of $M$ and $f \in \nonnegints$ the rank of $M$.
We call the remaining part of the decomposition $\bigoplus_{i=1}^r R \bigmod (d_i)$ the \textit{torsion} component of $M$.
We prefer the statement in \fref{thm:structure-theorem} since it is more suitable with the matrix calculation we present in this section.

\spacer 

The proof for the existence of invariant factor decompositions presented in \cite[Theorem 12.5]{algebra:dummit} relies on a system of generators and relations of a finitely generated module over a PID. 
We provide an alternate characterization of these systems below.

\begin{definition}\label{defn:presentation}
	A \textbf{presentation}
	of a module $M$ over a PID $R$ is an exact sequence 
	\begin{equation*}
		F_S \Xrightarrow{\,\phi\,} 
		F_G \Xrightarrow{\,\pi\,} 
		M \Xrightarrow{\,\,\,}\, 0
	\end{equation*}
	of free $R$-modules $F_S$ and $F_G$ with homomorphisms $\phi: F_S \to F_G$ and $\pi: F_G \to M$.
	We call $F_G$ and $F_S$ the \textbf{module of generators} and \textbf{module of relations} respectively.
	Given a basis 
	$A$ of $F_G$ and $S$ of $F_S$, 
	we call $\pi(A) \subseteq M$ and $(\phi \,\circ\, \pi)(S) \subseteq M$ a \textbf{system of generators and relations for $M$} respectively.
	A \textbf{finite presentation} of $M$ is a presentation wherein both $F_S$ and $F_G$ have finite rank. 

	When we say that a presentation of $M$ is given by $\phi: F_S \to F_G$ and $\pi: F_G \to M$, we refer to the exact sequence as given above.
\end{definition}
\remark{
	We usually indicate that a presentation $\phi:F_S \to F_G$ is finite by listing a finite basis $S = (\sigma_1, \ldots, \sigma_n)$ of $F_S$ and 
	$A = (\alpha_1, \ldots, \alpha_m)$ of $F_G$.
	Also, this notion of presentation is not generally compatible with the notion of presentation of groups, particularly in the case of non-abelian groups. 
}

In the proof presented in \cite[Theorem 12.5]{algebra:dummit}, 
a presentation of $M$ is constructed using a set of generators $\set{a_1, \ldots, a_m}$ of $M$, which exists by assumption of $M$ being finitely generated.
Let $\set{\alpha_1, \ldots, \alpha_m}$ be a set of indeterminates.
Define the module of generators to be $F_G = R\ket{\alpha_1, \ldots, \alpha_m}$ and the homomorphism $\pi: F_G \to M$ by $\alpha_j \mapsto a_j$ for $j \in \set{1, \ldots, m}$, i.e.\ $\alpha_j$ is essentially a relabeling of $a_j \in M$.
Note that we use different labels for $a_j$ and $\alpha_j$ since $a_j \in M$ may be a torsion element of $M$, i.e.\ there may exist $r \in R$ such that $r \cdot a_j = 0$, but $\alpha_j \in F_G$ cannot be since it is an element of the free $R$-module $F_G$.
Then, the following exact sequence is a finite presentation for $M$:
\begin{equation*}
	\equalsupto{F_S}{\ker(\pi)}
	\,\,\Xrightarrow{\,\phi\,}\,\,
		\equalsupto{F_G}{R\ket{\alpha_1, \ldots, \alpha_m}}
		\,\,\Xrightarrow{\,\pi\,}\,\,
		M 
		\,\,\Xrightarrow{\quad}\,\,
		0
\end{equation*}
where $\phi: \ker(\pi) \to R\ket{\alpha_1, \ldots, \alpha_m}$ is taken to be the inclusion map.
Since $\ker(\pi)$ is a submodule of a finitely generated module over a PID,
$\ker(\pi)$ must also be finitely generated. Note that this is not true in general if $R$ is not a PID.
Then, $\ker(\pi)$ is a free module with finite basis, i.e.\ of finite rank.
Then, the presentation above determines $M$ by the following isomorphism:
\begin{equation*}
	M 
		\,\overset{(1)}{=}\, 
			\im(\pi)
		\,\overset{(\star)}{\cong}\, 
			F_G \bigmod \ker(\pi)
		\,\overset{(2)}{=}\,
		 	F_G \bigmod \image(\phi)
\end{equation*}
where $(\star)$ is given by the first isomorphism theorem on $\pi: F_G \to M$, and 
$(1)$ and $(2)$ are both given by the exactness of the sequence $F_S \to F_G \to M \to 0$, i.e.\ $\im(\pi) = \ker(M \to 0) = M$ and $\ker(\pi) = \image(\phi)$ respectively.
Note that the homomorphism $\phi: F_S \to F_G$ determines $M$ up to isomorphism since $M \cong F_G \bigmod \im(\phi)$. 

Observe that this isomorphism holds true even if $F_S$ is not exactly $\ker(\pi)$, i.e.\ $F_S$ only needs to contain $\ker(\pi)$ as a submodule, 
or if $\phi:F_S \to F_G$ is not an inclusion map.
That is, the isomorphism holds for arbitrary presentations given by \fref{defn:presentation}.
We discuss state this in more detail later in \fref{prop:how-to-get-ifd}.

The next step in the proof of \cite[Theorem 12.5]{algebra:dummit} involves finding a basis on $\ker(\pi)$ and $F_G$, as denoted in \fref{defn:presentation}, of a presentation such that certain properties are fulfilled.
The existence of such a basis is guaranteed by the following theorem, taken from \cite[Theorem 12.5]{algebra:dummit}.

\begin{theorem}\label{thm:invariant-factor-theorem-for-submodules}
	\textbf{Invariant Factor Theorem for Submodules.} 

	Let $M$ be a free module over a PID $R$ with $\rank(M) = m$ and let $L$ be a submodule of $M$.
	Then, $L$ is a free submodule with $\rank(L) = r \leq m$
	and there exists a basis $B = \set{\beta_1, \ldots, \beta_m}$ of $M$ 
	and nonzero elements $d_1, \ldots, d_r \in R$
	with divisibility relation $d_1 \divides d_2 \divides \cdots \divides d_r$ 
	such that $\set{d_1\beta_1, \ldots, d_r\beta_r}$ is a basis of $L$.
	Furthermore, the elements $d_1, \ldots, d_r$ are unique up to multiplication by units.
\end{theorem}
\remark{
	For a proof of the existence and uniqueness claims, see \cite[Theorem 12.4]{algebra:dummit} and \cite[Theorem 12.9]{algebra:dummit} respectively.
}

With $M$ and $L$ as denoted in the theorem above, we have the following decompositions for $M$ and $L$ using the bases $\set{\beta_1, \ldots, \beta_m}$ and $\set{d_1\beta_1, \ldots, d_r\beta_r}$:
\begin{equation*}
	\begin{array}{cc cc cc cc cc cc cc cc}
		M &\cong
			&R\ket{\beta_1}
			&\oplus\,\cdots\,\oplus
			&R\ket{\beta_r}
			&\oplus
			&R\ket{\beta_{r+1}}
			&\oplus\,\cdots\,\oplus
			&R\ket{\beta_m}
		\\[2pt] 
		L &\cong
			&R\ket{d_1\beta_1}
			&\oplus\,\cdots\,\oplus
			&R\ket{d_r\beta_r}
	\end{array}
\end{equation*}
Note that the elements $d_1, \ldots, d_r \in R$ given by this theorem may still be invertible as elements of $R$, e.g. $d_j = 1$ for some $j \in \set{1, \ldots, r}$.
Observe that each basis element $d_j\beta_j$ of $L$ is then associated with a unique basis element $\beta_j$ of $M$.
Since for each $j \in \set{1, \ldots, r}$, $R\ket{d_j\beta_j}$ is a submodule of $R\ket{\beta_j}$,
we can characterize the quotient module $M \bigmod L$ by considering each pair of $R\ket{\beta_j}$ and $R\ket{d_j\beta_j}$ as a torsion summand of $M \bigmod L$.
We state the result that allows us to do this below.

\begin{lemma}\label{lemma:distribute-summands-over-quotients}
	Let $M$ and $N$ be modules over a ring $R$.
	Let $A$ be a submodule of $M$ and $B$ that of $N$.
	Then,
	\begin{equation*}
		\frac{M \oplus N}{A \oplus B}
		\cong 
		\paren{\frac{M}{A}}
		\oplus 
		\paren{\frac{N}{B}}
	\end{equation*}
\end{lemma}
\begin{proof}
	Let $\pi_1: M \to M \bigmod A$ and $\pi_2: N \to N \bigmod B$ be canonical quotient homomorphisms.
	Note that $\ker(\pi_1) = A$, $\ker(\pi_2) = B$ and 
	that both $\pi_1$ and $\pi_2$ are surjective, i.e.\
		$\im(\pi_1) = M \bigmod A$, $\im(\pi_2) = N \bigmod B$.
	The direct sum of modules induces a homomorphism $\pi = \pi_1 \oplus \pi_2$ as follows:
	\begin{align*}
		\pi: M \oplus N &\to \big( M \bigmod A \big) \oplus \big( N \bigmod B \big)
		\\
		(m,n) &\mapsto \Big( \pi_1(m), \pi_2(n) \Big) = (m+A, n+B)
	\end{align*}
	Then, $\ker(\pi) = \ker(\pi_1) \oplus \ker(\pi_2) = A \oplus B$ and 
	$\im(\pi) = \im(\pi_1) \oplus \im(\pi_2) = (M \bigmod A) \oplus (M \bigmod B)$.
	By the first isomorphism theorem on $\pi$,
	\begin{equation*}
		\im(\pi) = 
			\paren{\frac{M}{A}}
			\oplus 
			\paren{\frac{N}{B}}
		\cong 
		\frac{M \oplus N}{\ker(\pi)}
		= 
		\frac{M \oplus N}{A \oplus B}
		\,.
		\vspace{-\baselineskip}
	\end{equation*}
\end{proof}

We then apply \fref{thm:invariant-factor-theorem-for-submodules} on the image of a finite presentation of a finitely generated module over a PID. 
Note that the basis from this theorem allows us to apply \fref{lemma:distribute-summands-over-quotients}.
We state this in more detail below.

\begin{proposition}\label{prop:how-to-get-ifd}
	Let $M$ be a finitely generated module over a PID $R$.
	Let $\phi: F_S \to F_G$ and $\pi: F_G \to M$ correspond to a presentation of $M$ with $\rank(F_S) = n$ and $\rank(F_G) = m$.
	Then, there exists a basis $\set{\beta_1, \ldots, \beta_m}$ of $F_G$ and nonzero values 
	$d_1, \ldots, d_r \in R$ 
	with divisibility relation $d_1 \divides d_2 \divides \cdots \divides d_r$
	such that 
		$\set{d_1\beta_1, \ldots, d_r\beta_r}$ is a basis for $\im(\phi)$
		and 
	\begin{equation*}
		M 
		\cong 
		\frac{F_G}{\im(\phi)}
		\cong 
		\paren{\frac{R\ket{\beta_1}}{R\ket{d_1\beta_1}}}
		\oplus \cdots \oplus 
		\paren{\frac{R\ket{\beta_r}}{R\ket{d_r\beta_r}}}
		\oplus 
		R\ket{\beta_{r+1}}
		\oplus \cdots \oplus 
		R\ket{\beta_m}
	\end{equation*}
\end{proposition}
\begin{proof}
	Since images of module homomorphisms are submodules of the codomain,
	$\im(\phi)$ is a submodule of $F_G$.
	Let the basis $B = \set{\beta_1, \ldots, \beta_m}$ of $F_G$ and nonzero elements 
	$d_1, \ldots, d_r \in R$ be given by \fref{thm:invariant-factor-theorem-for-submodules} on $\im(\phi) =: L$.
	Then, $\im(\phi) \cong R\ket{d_1\beta_1} \oplus \cdots \oplus R\ket{d_r\beta_r}$
	and $F_G = R\ket{\beta_1} \oplus \cdots \oplus R\ket{\beta_m}$.
	Observe that for each $j \in \set{1, \ldots, r}$, $R\ket{d_1\beta_j}$ is a submodule of $R\ket{\beta_j}$.
	We apply \fref{lemma:distribute-summands-over-quotients} on $M \cong F_G \bigmod \im(\phi)$ as follows:
	\begin{align*}
		M
		&\cong 
			\frac{F_G}{\im(\phi)}
		\cong 
			\dfrac{
				R\ket{\beta_1} \oplus\cdots\oplus R\ket{\beta_m}
			}{
				R\ket{d_1\beta_1} \oplus\cdots\oplus R\ket{\beta_r}
			}
		\\[2pt]
		&\cong 
		\dfrac{
			R\ket{\beta_1} \oplus\cdots\oplus R\ket{\beta_m}
		}{
			R\ket{d_1\beta_1} \oplus\cdots\oplus R\ket{d_r\beta_r}
			\oplus \underbrace{0 \oplus \cdots \oplus 0}_{m-r \text{ times}}
		}
		\\[2pt]
		&\cong 
			\paren{
				\frac{R\ket{\beta_1}}{R\ket{d_1\beta_1}}
			}
			\oplus \cdots \oplus 
			\paren{
				\frac{R\ket{\beta_r}}{R\ket{d_r\beta_r}}
			}
			\oplus
			\paren{
				\frac{R\ket{\beta_{r+1}}}{0}
			}
			\oplus \cdots \oplus
			\paren{
				\frac{R\ket{\beta_{m}}}{0}
			}
		\\[2pt]
		&\cong 
			\paren{\frac{R\ket{\beta_1}}{R\ket{d_1\beta_1}}}
			\oplus \cdots \oplus 
			\paren{\frac{R\ket{\beta_r}}{R\ket{d_r\beta_r}}}
			\oplus 
			R\ket{\beta_{r+1}}
			\oplus \cdots \oplus 
			R\ket{\beta_m}
	\end{align*}
	where the trivial $R$-module is denoted by $0$.
\end{proof}

Note that, in the statement of the Structure Theorem (\fref{thm:structure-theorem}), the invariant factors $d_1, \ldots, d_n$ may include zero elements at the end of the sequence
while the elements $d_1, \ldots, d_r$, as denoted in \fref{prop:how-to-get-ifd}, are defined to be nonzero.
The notation does suggest that $d_1, \ldots, d_r$ are related to the invariant factors.
By \fref{prop:how-to-get-ifd}, the module of generators $F_G$ has $\rank(F_G) = m$. The zero invariant factors correspond to the free summands of $M$.
Since $R\ket{0\beta_j} = \set{r (0\beta_j) : r \in R} = \set{0}$ becomes trivial, we can define additional elements $d_{r+1}, \ldots, d_{m}$ to be zero and present the decomposition from \fref{prop:how-to-get-ifd} as follows:
\begin{equation*}
	M \cong 
	\bigoplus_{j=1}^m \paren{
		\frac{R\ket{\beta_j}}{d_j\beta_j}
	}
	= 
	\underbrace{\paren{
		\frac{R\ket{\beta_1}}{R\ket{d_1\beta_1}}
	}
	\oplus \cdots \oplus 
	\paren{
		\frac{R\ket{\beta_r}}{R\ket{d_r\beta_r}}
	}}_{\mathclap{\text{\small
		either torsion $R$-modules or trivial
	}}}
	\oplus
	\underbrace{\paren{
		\frac{R\ket{\beta_{r+1}}}{R\ket{d_{r+1}\beta_{r+1}}}
	}
	\oplus \cdots \oplus
	\paren{
		\frac{R\ket{\beta_{m}}}{R\ket{d_m\beta_m}}
	}}_{\mathclap{\text{\small these are free $R$-modules since $d_j = 0$}}}
\end{equation*}
This direct sum is then transformed into an invariant factor decomposition by replacing each summand into either a copy of $R$ or a cyclic ideal of $R$.
We state the required isomorphisms for these below.

\begin{lemma}\label{lemma:summands-of-ifd}
	Let $R\ket{a}$ be a free module over a PID $R$ with basis $\set{a}$ and let $d \in R$ be nonzero.
	Then, $R\ket{a} \cong R$.
	If $d$ is invertible, then $R\ket{a} \bigmod R\ket{da} \cong 0$, i.e.\ the trivial module. 
	Otherwise, $R\ket{a} \bigmod R\ket{da} \cong R \bigmod (d)$.
\end{lemma}
\begin{proof}
	Let $f: R\ket{a} \to R$ be given by $a \mapsto 1$ where $1 \in R$ refers to the identity element of $R$.
	Observe that $f$ is a homomorphism with inverse $r \mapsto ra \in R\ket{a}$ for all $r \in R$.
	Then, $f$ is an isomorphism and $R\ket{a} \cong R$.
	We then examine two cases below.
	\begin{enumerate}
		\item 
		Assume that $d \in R$ is invertible, i.e.\ there exists $d\inv \in R$ such that $d\inv d = 1$.
		Since $R\ket{da} \subseteq R\ket{a}$, it suffices to show $R\ket{da} \subseteq R\ket{a}$ to claim that $R\ket{da} = R\ket{a}$.
		Each element $ra \in R\ket{a}$ with $r \in R$ is generated by $s = d\inv r$ in $R\ket{da}$ as follows:
		\begin{equation*}
			(s)(da) = (d\inv r)(da) = d\inv d (ra) = ra 
			\,.
		\end{equation*}
		Note that PIDs are commutative by definition.
		Therefore, $R\ket{da} = R\ket{a}$ and $R\ket{a} \bigmod R\ket{da} = 0$.

		\item 
		Let $\pi: R \to R \bigmod (d)$ be the canonical quotient map.
		By the $1$\st isomorphism theorem on $(\pi \circ f)$:
		\begin{equation*}
			R \bigmod (d)
			=
			\im(\pi \circ f) 
			\cong 
			\frac{R\ket{a}}{\ker(\pi \circ f)}
			= 
			\frac{R\ket{a}}{
				\set{ra : r \in (d) }
			}
			=
			\frac{R\ket{a}}{R\ket{da}}
			\,.\vspace{-2\baselineskip}
		\end{equation*}
	\end{enumerate}
\end{proof}

Then, by applying \fref{lemma:summands-of-ifd} to the direct sum resulting from \fref{prop:how-to-get-ifd} and removing the trivial summands, we get an invariant factor decomposition of a module over a PID.
Since the nonzero values $d_1, \ldots, d_r$, as denoted in \fref{prop:how-to-get-ifd}, must satisfy the divisibility condition $d_1 \divides \cdots \divides d_r$, any invertible values must occur at the beginning of the list.
That is, we have some $k$ such that $1 \leq k \leq r$ and 
the elements $d_1, \ldots, d_k$ are invertible (contributing to trivial summands) and the elements $d_{k+1}, \ldots, d_r$ are not invertible.
Then, we get an invariant factor decomposition as follows:
\begin{equation*}
	\arraycolsep=1.5pt
	\begin{array}{
		*9{>{\displaystyle}c} *9{>{\displaystyle}c} 
	}
		{} && \overbrace{\hphantom{
			\paren{\frac{R\ket{\beta_1}}{R\ket{d_1\beta_1}}}
			\oplus \,\cdots\, \oplus 
			\paren{\frac{R\ket{\beta_k}}{R\ket{d_r\beta_k}}}
		}}^{\displaystyle\text{trivial summands}}
		&& \overbrace{\hphantom{
			\paren{\frac{R\ket{\beta_{k+1}}}{R\ket{d_{k+1}\beta_{k+1}}}}
			\oplus \,\cdots\, \oplus 
			\paren{\frac{R\ket{\beta_r}}{R\ket{d_r\beta_r}}}
		}}^{\displaystyle\text{torsion summands}}
		&& \overbrace{\hphantom{
			R\ket{\beta_{r+1}}
			\oplus \,\cdots\, \oplus 
			R\ket{\beta_m}
		}}^{\displaystyle\text{free summands}}
		\\
		M &\cong&
			\paren{\frac{R\ket{\beta_1}}{R\ket{d_1\beta_1}}}
			\oplus \,\cdots\, \oplus 
			\paren{\frac{R\ket{\beta_k}}{R\ket{d_k\beta_k}}}
		&\oplus& 
			\paren{\frac{R\ket{\beta_{k+1}}}{R\ket{d_{k+1}\beta_{k+1}}}}
			\oplus \,\cdots\, \oplus 
			\paren{\frac{R\ket{\beta_r}}{R\ket{d_r\beta_r}}}
		&\oplus& 
			R\ket{\beta_{r+1}}
			\oplus \,\cdots\, \oplus 
			R\ket{\beta_m}
		\\[14pt] 
		&\cong& 
			\text{\color{gray} (nothing in here)} 
		&\oplus&
			\frac{R}{(d_{k+1})}
			\oplus \,\cdots\, \oplus 
			\frac{R}{(d_r)}
		&\oplus& 
			R \oplus \,\cdots\, \oplus R
		\\[-12pt] 
			&& && && 
			\underbrace{\hphantom{R \oplus \,\cdots\, \oplus R}}_{
				\displaystyle\text{\small $f := m-r$ times}
			}
	\end{array}
\end{equation*}

\noindent 
The matrix calculation for invariant factor decompositions essentially comes from re-stating the Invariant Factor Theorem for Submodules (\fref{thm:invariant-factor-theorem-for-submodules}) in terms of matrices.
We state this in more detail below.

\begin{proposition}\label{prop:change-of-basis-on-homomorphism}
	Let $R$ be a PID.
	Let $\phi: N \to M$ be an $R$-module homomorphism between free $R$-modules $N$ and $N$ 
	with $\rank(N) = n$ and $\rank(M) = m$.
	Then, there exists a basis 
		$T = (\tau_1, \ldots, \tau_n)$ of $N$,
	a basis 
		$B = (\beta_1, \ldots, \beta_m)$ of $M$,
	and nonzero elements $d_1, \ldots, d_r \in R$ with $r = \rank(\im(\phi))\leq n$ 
	and divisibility relation $d_1 \divides d_2 \divides \cdots \divides d_r$
	such that 
	\begin{equation*}
		\phi(\tau_i) = \begin{cases}
			d_i \beta_i 
				&\text{ if } i \in \set{1, \ldots, r}
			\\ 
			0 
				&\text{ if } i \in \set{r+1, \ldots, n}
		\end{cases}
	\end{equation*}
	That is, the matrix $[\phi]_{B,T}$ of $\phi$ relative to $T$ and $B$ is given by the following block matrix
	\begin{equation*}
		[\phi]_{B,T} = \begin{pmatrix}
			D_r & 0 \\
			0 & 0
		\end{pmatrix} \in M_{m,n}(R)
		\qquad\text{ with }\qquad 
		D_r = \diag(d_1, \ldots, d_r)
		= \begin{pmatrix}
			d_1 & 0 & \cdots & 0 \\
			0 & d_2 & \cdots & 0 \\
			\vdots & \vdots & \ddots & 0 \\
			0 & 0 & \cdots & d_r
		\end{pmatrix}
	\end{equation*}
	Furthermore, the elements $d_1, \ldots, d_r$ are unique by multiplication of units.
\end{proposition}
\begin{proof}
	We provide an outline of the proof given for \cite[Proposition 4.3.20]{algebra:adkins} below.
	
	Since $\im(\phi)$ is a submodule of $M$, the Invariant Factor Theorem for Submodules 
	(\fref{thm:invariant-factor-theorem-for-submodules}) applies.
	Let $B = (\beta_1, \ldots, \beta_m)$ be the basis of $M$
	and $\set{d_1, \ldots, d_r}$ with $r \leq m$ be the set of nonzero elements given by 
	\fref{thm:invariant-factor-theorem-for-submodules}.
	Then, $\set{d_1\beta_1, \ldots, d_r\beta_r}$ is a basis for $\im(\phi)$ and the divisibility relation $d_1 \divides \cdots \divides d_r$ is satisfied.
	Note that $d_1, \ldots, d_r$ are also unique up to multiplication by units.

	For each $i \in \set{1, \ldots, r}$, choose $\tau_i \in N$ such that $\phi(\tau_i) = d_i\beta_i$.
	Note that $N \cong \ker(\phi) \oplus (N \bigmod \ker(\phi))$ and that 
	$\set{\tau_1, \ldots, \tau_r}$ is a basis of $N \bigmod \ker(\phi)$.
	Since $\ker(\phi)$ is a free submodule of $N$ with $n = \rank(\ker\phi) + r$, there exists a basis $\set{\tau_{r+1}, \ldots, \tau_n}$ of $\ker(\phi)$.
	One then proves that 
	$\set{\tau_1, \ldots, \tau_r, \tau_{r+1}, \ldots, \tau_n}$ is a basis of $N$ by showing linear independence.
	Then, by construction, $\phi(\tau_i) = d_i\beta_i$ if $i \in \set{1, \ldots, r}$ and $\phi(\tau_i) = 0$ otherwise.

	The matrix representation follows from the definition of matrices of homomorphisms. For reference, see \fref{defn:coordinate-matrices} in \fref{appendix:module-theory}.
	Let $i \in \set{1, \ldots, n}$ be a column index.
	By definition, $\col_i([\phi]_{B,T}) = [\phi(\tau_i)]_B$, i.e.\ the coordinate vector of $\phi(\tau_i)$ relative to $B$.
	If $i \in \set{1, \ldots, r}$:
	\begin{equation*}
		[\phi]_{B,T}(j,i) = [\phi(\tau_i)]_B(j) = \begin{cases}
			d_i	&\text{ if } j=i \\
			0 	&\text{ otherwise }
		\end{cases}
	\end{equation*}
	If $i \in \set{r+1, \ldots, n}$, then $\col_i([\phi]_{B,T}) = [\phi(\tau_i)]_B$ would be the zero column since $\phi(\tau_i) = 0$.	
\end{proof}

The proposition above is then used to prove the existence of Smith Normal Decompositions of matrices, defined below.
Note that this definition is taken from from \cite[Section 5.3]{algebra:adkins}.

\begin{definition}\label{defn:smith-normal-decomposition}
	The \textbf{Smith Normal Decomposition (SND)} of a matrix $A \in \M_{m,n}(R)$ is a triple $(U,D,V)$ of matrices $U \in \GL(m,R)$, $V \in \GL(n,R)$, $D \in \M_{m,n}(R)$ 
	such that 
	\begin{equation*}
		U\inv A V = D = \begin{pmatrix}
			D_r & 0 \\
			0 & 0
		\end{pmatrix}
		\quad\text{ and }\quad 
		D_r = \diag(d_1, \ldots, d_r)
	\end{equation*}	
	with $r = \rank(A)$ and nonzero elements $d_1, d_2, \ldots, d_r \in R$ that satisfy divisibility relation $d_1 \divides d_2 \divides \cdots \divides d_r$.
	We call $D \in \M_{m,n}(R)$ the \textbf{Smith Normal Form (SNF)} of $A \in M_{m,n}(R)$. 
\end{definition}

We do warn readers that, in most references, the matrix $U \in \GL(m,R)$ in the Smith Normal Decomposition $(U,D,V)$ of $A \in \M_{m,n}(R)$ 
is defined such that $UAV = D$, as opposed to our definition where $U$ corresponds to the factorization $U\inv AV = D$.
We have decided to change the definition for $U \in \GL(m,R)$ so that the factorization $U\inv AV = D$ is more compatible with our application, i.e.\ we should interpret $U\inv AV = D$ as corresponding to some $R$-module homomorphism equipped with change of bases on the domain and codomain, as we will see later in this section in \fref{prop:snd-on-homomorphisms}.

Next, we state the result for the existence and uniqueness of Smith Normal Decompositions.
The main idea here involves creating an $R$-module homomorphism from the given matrix and using \fref{prop:change-of-basis-on-homomorphism}.
Note that definitions for elementary matrices over $R$, including elementary dilations, are given in \fref{defn:elementary-matrices}.

\begin{theorem}\label{thm:smith-normal-existence-and-uniqueness}
	Let $A \in M_{m,n}(R)$ be a matrix over a PID $R$.
	Then, a Smith Normal Decomposition $(U,D,V)$ of $A$ exists and the Smith Normal Form $D$ of $A$ is unique up to elementary dilations over $R$,
		i.e.\ the diagonal elements of $D$ are unique up to multiplication by units in $R$.
\end{theorem}
\begin{proof}
	We refer to the proof in \cite[Theorem 5.3.1]{algebra:adkins} and provide an outline below.

	Let $\basis{S} = (\sigma_1, \ldots, \sigma_n)$ denote the standard ordered basis on $R^n$ and let $\basis{A} = (\alpha_1, \ldots, \alpha_m)$ be that on $R^m$.
	Define the $R$-module homomorphism $\phi: R^n \to R^m$ such that 
	$[\phi]_{\basis{A}, \basis{S}} = A$, i.e.\ $A$ is the matrix of $\phi$ relative to $\basis{A}$ and $\basis{S}$.
	By \fref{prop:change-of-basis-on-homomorphism},
		there exists a basis $\basis{T} = (\tau_1, \ldots, \tau_n)$ of $N$,
		a basis $\basis{B} = (\beta_1, \ldots, \beta_m)$ of $M$,
		and nonzero elements 
			$d_1, \ldots, d_r \in R$ with $r = \rank(A) \leq n$
	such that 
	\begin{equation*}
		[\phi]_{\basis{B}, \basis{T}} 
		= \begin{pmatrix}
			D_r & 0 \\
			0 & 0
		\end{pmatrix} \in M_{m,n}(R)
		\qquad\text{ with }\qquad 
		D_r = \diag(d_1, \ldots, d_r)
	\end{equation*}
	Let $[\id_M]_{\basis{B}, \basis{A}} \in \GL(m,R)$ be the matrix of the identity map $\id_M: M \to M$ relative to $\basis{A}$ and $\basis{B}$. Note that $[\id_M]_{\basis{B}, \basis{A}}$ is invertible with inverse $[\id_M]_{\basis{A}, \basis{B}}$.
	Similarly, let $[\id_N]_{\basis{T},\basis{S}} \in \GL(n,R)$ be the matrix of the identity map $\id_N: N \to N$ relative to $\basis{T}$ and $\basis{S}$.
	Then, we can relate $[\phi]_{\basis{A}, \basis{S}}$ 
	and $[\phi]_{\basis{B}, \basis{T}}$ by the following matrix equation:
	\begin{equation*}
		[\id_M]_{\basis{B}, \basis{A}}	\,
		[\phi]_{\basis{A}, \basis{S}}	\,
		[\id_N]_{\basis{S}, \basis{T}}
		= 
		[\phi]_{\basis{B}, \basis{T}}
	\end{equation*}
	Let $U \in \GL(m,R)$, $D \in \M_{m,n}(R)$, and $V \in \GL(n,R)$ be given by 
	\begin{equation*}
		U = [\id_M]_{\basis{A}, \basis{B}}
		\qquad 
		D = [\phi]_{\basis{B}, \basis{T}}
		\qquad 
		V = [\id_N]_{\basis{S}, \basis{T}}
	\end{equation*}
	Observe that the divisibility relation on the entries $d_1, \ldots, d_r$ is satisfied by \fref{prop:change-of-basis-on-homomorphism}.
	We can confirm that the matrix factorization is correct by the following calculation:
	\begin{equation*}
		U\inv A V 
		= 
		\Big(
			[\id_M]_{\basis{A}, \basis{B}}
		\Big)\inv
		\Big(
			[\phi]_{\basis{A}, \basis{S}}
		\Big)
		\Big(
			[\id_N]_{\basis{S}, \basis{T}}
		\Big)
		= 
		[\id_M]_{\basis{B}, \basis{A}}	\,
		[\phi]_{\basis{A}, \basis{S}}	\,
		[\id_N]_{\basis{S}, \basis{T}}
		= 
		[\phi]_{\basis{B}, \basis{T}}
		= 
		D
	\end{equation*}
	Since the elements $d_1, \ldots, d_r$ satisfy the divisibility relation $d_1 \divides \ldots \divides d_r$ and $U\inv AV = D$, $(U,D,V)$ is an SND of $A$.

	Note that the uniqueness of the Smith Normal Form $D$ of $A$ (up to multiplication by elementary dilations over $R$) is implied by the uniqueness of the invariant factors up to multiplication by units, as stated in the Invariant Factor Theorem for Submodules (\fref{thm:invariant-factor-theorem-for-submodules}).
\end{proof}

Below, we provide an example of two SNDs of a matrix over $\ints$.
Note that for both SNDs, the Smith Normal Forms on both SNDs only differ by an elementary dilation, e.g.\ multiplication by $(-1)$ since $(-1)$ is a unit in $\ints$, but the matrices $U$ and $V$ on an SND $(U,D,V)$ are not generally unique.

\begin{example}\def\arraystretch{0.9}
	Let $A \in \M_{4,3}(\ints)$ be given as follows:
	\begin{equation*}
		A = \begin{pmatrix}
			1 & 2 & 3 \\
			4 & 5 & 6 \\
			7 & 8 & 9 \\
			1 & 2 & 4
		\end{pmatrix}
	\end{equation*}
	By Mathematica (a software system), 
	$A$ admits the following Smith Normal Decomposition $(U_1,D_1,V_1)$ 
	where $(U_1)\inv A V_1 = D_1$ and the matrices $U_1 \in \GL(4,\ints)$, $D_1 = M_{4,3}(\ints)$, $V_1 = \GL(3,\ints)$ are given by 
	\begin{equation*}
		U_1 = \begin{pmatrix}
			1 & 3 & 0 & 0 \\
			4 & 3 & 1 & 0 \\
			7 & 3 & 2 & 1 \\
			1 & 4 & 0 & 0
		\end{pmatrix}
		\qquad 
		D_1 = \begin{pmatrix}
			1 & 0 & 0 \\
			0 & 1 & 0 \\
			0 & 0 & 3 \\
			0 & 0 & 0
		\end{pmatrix}
		\qquad 
		V_1 = \begin{pmatrix}
			1 & -2 & 2 \\
			0 & 1 & -1 \\
			0 & 1 & 0
		\end{pmatrix}
	\end{equation*}
	We can confirm this by doing the following calculation:
	\begin{equation*}
		(U_1)\inv A V_1 
		= 
		\begin{pmatrix}
			4 & 0 & 0 & -3 \\
			-1 & 0 & 0 & 1 \\
			-13 & 1 & 0 & 9 \\
			1 & -2 & 1 & 0
		\end{pmatrix}
		\begin{pmatrix}
			1 & 2 & 3 \\
			4 & 5 & 6 \\
			7 & 8 & 9 \\
			1 & 2 & 4
		\end{pmatrix}
		\begin{pmatrix}
			1 & -2 & 2 \\
			0 & 1 & -1 \\
			0 & 1 & 0
		\end{pmatrix}
		=
		\quad\cdots\quad
		= 
		\begin{pmatrix}
			1 & 0 & 0 \\
			0 & 1 & 0 \\
			0 & 0 & 3 \\
			0 & 0 & 0
		\end{pmatrix}
		= D_1
	\end{equation*}
	Another Smith Normal Decomposition of $A$ is $(U_2, D_2, V_2)$
	with 
		$U_2 \in \GL(4,\ints)$, 
		$D_2 = M_{4,3}(\ints)$, and 
		$V_2 \in \GL(3,\ints)$
	given as follows:
	\begin{equation*}
		U_2 = \begin{pmatrix}
			6 & 4 & -1 & 0	\\
			15 & 13 & -3 & 0	\\
			24 & 22 & -5 & 1 \\
			7 & 4 & -1 & 0 
		\end{pmatrix}
		\qquad
		D_2 = \begin{pmatrix}
			1 & 0 & 0 \\
			0 & 1 & 0 \\
			0 & 0 & -3 \\
			0 & 0 & 0
		\end{pmatrix}
		\qquad 
		V_2 = \begin{pmatrix}
			1 & 2 & 1 \\
			1 & 1 & 1 \\
			1 & 0 & 0
		\end{pmatrix}
	\end{equation*}
	As with $(U_1,D_1,V_1)$, we can confirm that $(U_2, D_2, V_2)$ is a valid factorization by doing the following calculation:
	\begin{equation*}
		(U_2)\inv 
		A 
		(V_2)
		= 
		\begin{pmatrix}
			-1 & 0 & 0 & 1 \\
			-6 & 1 & 0 & 3 \\
			-31 & 4 & 0 & 18 \\
			1 & -2 & 1 & 0
		\end{pmatrix}
		\begin{pmatrix}
			1 & 2 & 3 \\
			4 & 5 & 6 \\
			7 & 8 & 9 \\
			1 & 2 & 4
		\end{pmatrix}
		\begin{pmatrix}
			1 & 2 & 1 \\
			1 & 1 & 1 \\
			1 & 0 & 0
		\end{pmatrix}
		= 
		\quad\cdots\quad 
		= 
		\begin{pmatrix}
			1 & 0 & 0 \\
			0 & 1 & 0 \\
			0 & 0 & -3 \\
			0 & 0 & 0
		\end{pmatrix}
		= D_2
	\end{equation*}
	Observe that the diagonal elements $D_2(1,1)=1, D_2(2,2)=1, D_2(3,3) =-3$ of $D_2$ obey the divisibility rule where $d_1$ divides $d_2$ and $d_2$ divides $d_3$.
	Also, the first two diagonal elements of $D_1$ and $D_2$ match and for the $3$\rd diagonal element, we have $D_2(3,3) = -3 = (-1)D_1(3,3) = (-1)(3)$.
\end{example}

As implied in the proof of \fref{thm:smith-normal-existence-and-uniqueness}, 
given a homomorphism $\phi: N \to M$, the nonzero diagonal elements $d_1, \ldots, d_r$ of the Smith Normal Form $D$ of the matrix $[\phi]_{A,S}$ are exactly the nonzero elements guaranteed by the Invariant Factor Theorem on Submodules (\fref{thm:invariant-factor-theorem-for-submodules}).
We believe this is partly why the nonzero diagonal elements of the Smith Normal Form of a matrix are sometimes called the \textit{invariant factors} of the matrix $[\phi]_{A,S}$. 
To avoid confusion between the invariant factors of modules (which cannot be invertible elements) and that of matrices (which can be invertible elements), 
we will only use invariant factors in the context of modules, i.e.\ the nonzero elements of the SNF of matrices are not called invariant factors in this paper.

This means that the problem of finding invariant factor decompositions can now be expressed as a matrix factorization problem, i.e.\ that of finding an SND.
Below, we describe how we should interpret an SND of $[\phi]_{A,S}$, as denoted above, in order to apply \fref{prop:how-to-get-ifd} to calculate invariant factor decompositions.

\begin{proposition}\label{prop:snd-on-homomorphisms}
	Let $\phi: N \to M$ be a module homomorphism between free $R$-modules 
	$N$ and $M$ with ordered bases 
		$S = (\sigma_1, \ldots, \sigma_n)$ and 
		$A = (\alpha_1, \ldots, \alpha_m)$
	respectively.
	Let $(U,D,V)$ be an SND of $[\phi]_{A, S} \in \M_{m,n}(R)$
	and let $d_i \in R$ be given by $d_i = D(i,i) \neq 0$ for $i \in \set{1, \ldots, r}$ with $r = \rank([\phi]_{A,S})$.

	Then, 
	$V \in \GL(n,R)$ determines a basis $T = (\tau_1, \ldots, \tau_n)$ of $N$ given by $[\tau_i]_S = \col_i(V)$ for $i \in \set{1, \ldots, n}$.
	Similarly, 
		$U \in \GL(m,R)$ determines a basis $B = (\beta_1, \ldots, \beta_m)$ of $M$ given by $[\beta_j]_A = \col_i(U)$ for $j \in \set{1, \ldots, m}$.
	Furthermore, the divisibility relation $d_1 \divides d_2 \divides \cdots \divides d_r$ is satisfied and $D = [\phi]_{B,T}$, i.e.\ 
	\begin{equation*}
		\phi(\tau_i) = \begin{cases}
			d_i \beta_i 	&\text{ if } i \in \set{1, \ldots, r} \\
			0				&\text{ if } i \in \set{r+1, \ldots, n}
		\end{cases}
	\end{equation*}
	Note that $\set{\phi(\tau_1), \ldots, \phi(\tau_r)} = \set{d_1\beta_1, \ldots, d_r\beta_r}$ forms a basis of $\im(\phi) \subseteq M$
	and 
	$\set{\tau_{r+1}, \ldots, \tau_n}$ determines a basis for $\ker(\phi) \subseteq N$.
\end{proposition}
\begin{proof}
	Let $e_i^{[n]} \in \M_{n,1}(R)$ denote the $i$\th standard basis (column) vector of $R^n$ and let $e_j^{[m]} \in \M_{m,1}(R)$ be the $j$\th standard basis (column) vector of $R^m$.

	Since $V \in \GL(n,R)$ and $U \in \GL(m,R)$ are invertible matrices, they can be used as change of basis matrices.
	Define a basis $T = (\tau_1, \ldots, \tau_n)$ of $N$ by $V = [\id_N]_{S,T}$.
	Note that $[\id_N]_{T,S} = V\inv$. 
	Similarly, define a basis $B = (\beta_1, \ldots, \beta_m)$ of $M$ by 
	$[\id_M]_{A,B} = U$.
	Note that $[\id_M]_{B,A} = U\inv$.
	Then, the coordinate vectors of each $\tau_i$ and $\beta_j$ relative to $S$ and $T$ respectively are given below:
	\begin{equation*}\arraycolsep=5pt
		\begin{array}{c !{=} l !{=} l !{=} ll} 
			{[\tau_i]_S }
				& [\id_N]_{S,T} \, [\tau_i]_T
				& V e_i^{[n]}
				& \col_i(V) \in \M_{n,1}(R)
				&\text{for } i \in \set{1, \ldots, n}
			\\[2pt]
			{[\beta_j]_A}
				& [\id_M]_{A,B} \, [\beta_j]_B 
				& \big( [\id_M]_{B,A} \big)\inv e_j^{[m]}
					= U  e_j^{[m]}
				& \col_j(U) \in \M_{m,1}(R)
				& \text{for } j \in \set{1, \ldots, m}
		\end{array}
	\end{equation*}
	Since $(U,D,V)$ is an SND of $[\phi]_{A,S}$, $U\inv [\phi]_{A,S} V = D$.
	Then, $[\phi]_{B,T} = D$ by the following calculation:
	\begin{equation*}
		D = 
		U\inv [\phi]_{A,S} V 
		= 
		\Bigl(
			[\id_M]_{A,B}
		\Bigr)\inv 
		\Bigl(
			[\phi]_{A,S}
		\Bigr)
		\Bigl(
			[\id_N]_{S,T}
		\Bigr)
		= 
		[\id_M]_{B,A}
		[\phi]_{A,S}
		[\id_N]_{S,T}
		= 
		[\phi]_{B,T}
	\end{equation*}
	Then, for each $i \in \set{1, \ldots, n}$, the coordinate vector of $\phi(\tau_i)$ relative to $B$ is as follows: 
	\begin{equation*}
		[\phi(\tau_i)]_B 
		= [\phi]_{B,T} [\tau_i]_T 
		= D e_i^{[n]} 
		= \col_i(D)
		= \begin{cases}
			d_i e_i^{[m]} 	
				&\text{ if } i \in \set{1, \ldots, r} \\
			0 		
				&\text{ if } i \in \set{r+1, \ldots, n}
		\end{cases}
	\end{equation*}
	Therefore, $\phi(\tau_i) = d_i \beta_i$ if $\set{1, \ldots, r}$ and $\phi(\tau_i) = 0$ otherwise.
	Observe that $\set{\phi(\tau_1), \ldots, \phi(\tau_r)} = \set{d_1\beta_1, \ldots, d_r\beta_r}$ is a basis of $\im(\phi)$
	and $\set{\tau_{r+1}, \ldots, \tau_n}$ is a basis for $\ker(\phi)$.
	The divisibility relation $d_1 \divides d_2 \divides \cdots \divides d_r$ is satisfied by definition of SND.
\end{proof}

Since an SND $(U,D,V)$ of the matrix $[\phi]_{A,S}$ of a homomorphism $\phi: N \to M$, as denoted above, gives us the results of the Invariant Factor Theorem for Submodules (\fref{thm:invariant-factor-theorem-for-submodules}), i.e.\ 
we can identify the basis $B = (\beta_1, \ldots, \beta_m)$ of $M$ and the basis $\set{d_1\beta_1, \ldots, d_r\beta_r}$ of $\im(\phi)$, 
we can use \fref{prop:how-to-get-ifd} to calculate invariant factor decompositions.

Observe that given a finite presentation $\phi: F_S \to F_G$, $\pi: F_G \to M$ of an $R$-module $M$, the basis $B = (\beta_1, \ldots, \beta_m)$ and the nonzero elements $d_1, \ldots, d_r \in R$ given by the above proposition (\fref{prop:snd-on-homomorphisms}) are exactly the information required in the hypothesis of \fref{prop:how-to-get-ifd}.
In particular, if given ordered bases $S = (\sigma_1, \ldots, \sigma_n)$ of $F_S$
and $A = (\alpha_1, \ldots, \alpha_m)$, 
an SND $(U,D,V)$ of $[\phi]_{A,S}$ determines the basis $B = (\beta_1, \ldots, \beta_m)$ by $[\beta_j]_A = \col_j(U)$ and the nonzero elements $d_j = D(j,j)$
and yields the following isomorphism by application of \fref{prop:how-to-get-ifd}:
\begin{align*}
	M \cong 
	\frac{F_G}{\im(F_S)}
	&\cong 
	\frac{R\ket{\alpha_1, \ldots, \alpha_m}}{\im(F_S)}
		\\ 
	&\cong 
	\frac{
		R\ket{\beta_1, \ldots, \beta_m}
	}{
		R\ket{d_1\beta_1, \ldots, d_r\beta_r}
	} &&\text{by \fref{prop:snd-on-homomorphisms} }
	\\ 
	&\cong 
	\paren{\frac{R\ket{\beta_1}}{R\ket{d_1\beta_1}}}
	\oplus \cdots \oplus 
	\paren{\frac{R\ket{\beta_r}}{R\ket{d_1\beta_r}}}
	\oplus 
	R\ket{\beta_{r+1}}
	\oplus \cdots \oplus 
	R\ket{\beta_m}
		&&\text{by \fref{prop:how-to-get-ifd} }
	\\ 
	&\cong 
	\frac{R}{(d_1)}
	\oplus \cdots \oplus 
	\frac{R}{(d_k)}
	\oplus 
	R \oplus \cdots \oplus R
		&&\text{by \fref{lemma:summands-of-ifd} }
\end{align*}
where $k \in \set{1, \ldots, r}$ is chosen such that $d_{k+1}, \ldots, d_r$ are all the invertible elements of the set $\set{d_1, \ldots, d_r}$.
We provide an example of this calculation below, using a homomorphism between $\ints$-modules.

\begin{example}\def\arraystretch{0.9}
	Let $N = \ints\ket{a,b,c,d}$ and $M = \ket{x,y,z}$. 
	Define the module homomorphism $\phi: N \to M$ by
	\begin{equation*}
		a \mapsto x-y, \quad 
		b \mapsto y+z, \quad 
		c \mapsto 2x, \quad 
		d \mapsto z.
	\end{equation*}
	Then, the matrix of $\phi$ relative to the bases $A = (a,b,c,d)$ of $N$ and $X = (x,y,z)$ of $M$ is as follows:
	\begin{equation*}
		[\phi]_{X,A} =
		\begin{blockarray}{*{5}{>{\textstyle\color{gray}}c}}
			{} & a & b & c & d \\
			\begin{block}{>{\textstyle\color{gray}}c (*{4}{c})}
				{x} & 1		& 0 & 2 & 0 \\
				{y} & -1	& 1	& 0 & 0	\\
				{z} & 0		& 1	& 0 & 1	\\
			\end{block}
		\end{blockarray}\vspace{-\parskip}
	\end{equation*}
	Using Mathematica, the triple $(U,D,V)$ of matrices below corresponds to an SND of $[\phi]_{X,A}$:
	\begin{equation*}
		U = 
		\begin{pmatrix}
			1 & 0 & 1 \\
			-1 & 1 & 0 \\
			0 & 1 & 0
		\end{pmatrix}
		\qquad 
		D = \begin{pmatrix}
			1 & 0 & 0 & 0 \\
			0 & 1 & 0 & 0 \\
			0 & 0 & 1 & 0
		\end{pmatrix}
		\qquad 
		V = \begin{pmatrix}
			1 & 0 & -1 & -2 \\
			0 & 1 & -1 & -2 \\
			0 & 0 & 1 & 1 \\
			0 & 0 & 1 & 2
		\end{pmatrix}
	\end{equation*}
	The validity of this factorization can be confirmed by performing the following calculation: 
	\begin{equation*}
		U\inv [\phi]_{X,A} V 
		= 
		\begin{pmatrix}
			0 & -1 & 1 \\
			0 & 0 & 1 \\
			1 & 1 & -1
		\end{pmatrix}
		\begin{pmatrix}
			1		& 0 & 2 & 0 \\
			-1		& 1	& 0 & 0	\\
			0		& 1	& 0 & 1	\\
		\end{pmatrix}
		\begin{pmatrix}
			1 & 0 & -1 & -2 \\
			0 & 1 & -1 & -2 \\
			0 & 0 & 1 & 1 \\
			0 & 0 & 1 & 2
		\end{pmatrix}
		=
		\quad\cdots\quad 
		= 
		\begin{pmatrix}
			1 & 0 & 0 & 0 \\
			0 & 1 & 0 & 0 \\
			0 & 0 & 1 & 0
		\end{pmatrix}
		= 
		D
	\end{equation*}
	From $V \in \GL(4,\ints)$, we determine the basis $\basis{A} = (\alpha_1, \alpha_2, \alpha_3, \alpha_4)$ of $N$ by $[\alpha_i]_{A} = \col_i(V)$
	as follows:
	\begin{equation*}
		\begin{aligned}
			[\alpha_1]_A &= \col_1(V) \\
			\alpha_1 &= a
		\end{aligned}
		\qquad\quad
		\begin{aligned}
			[\alpha_2]_A &= \col_2(V) \\
			\alpha_2 &= b
		\end{aligned}
		\qquad\quad
		\begin{aligned}
			[\alpha_3]_A &= \col_3(V) \\
			\alpha_3 &= -a-b+c+d
		\end{aligned}
		\qquad\quad
		\begin{aligned}
			[\alpha_4]_A &= \col_4(V) \\
			\alpha_4 &= -2a-2b+c+2d
		\end{aligned}
	\end{equation*}
	From $U \in \GL(3,\ints)$, we get the basis $\basis{Y} = (\gamma_1, \gamma_2, \gamma_3)$ of $M$ by $[\gamma_i] = \col_i(U)$ as follows:
	\begin{equation*}
		\begin{aligned}
			[\gamma_1]_X &= \col_1(U) \\
			\gamma_1 &= x-y
		\end{aligned}
		\qquad\quad
		\begin{aligned}
			[\gamma_2]_X &= \col_2(U) \\
			\gamma_2 &= y+z
		\end{aligned}
		\qquad\quad
		\begin{aligned}
			[\gamma_3]_X &= \col_3(U) \\
			\gamma_3 &= x
		\end{aligned}
	\end{equation*}
	We can confirm that $D = [\phi]_{\basis{Y}, \basis{A}}$ by the following calculations:
	\begin{align*}
		\phi(\alpha_1) &= \phi(a) 
			= x-y 
			&= \gamma_1 
		\\
		\phi(\alpha_2) &= \phi(b)
			= y+z 
			&= \gamma_2 
		\\ 
		\phi(\alpha_3) &= \phi(-a-b+c+d) 
			= -(x-y)-(y+z)+2x+z  		
			= x 
			&= \gamma_3
		\\ 
		\phi(\alpha_4) &= \phi(-2a-2b+c+2d) 
			= -2(x-y)-2(y+z)+2x+2(z)  		
			&= 0
	\end{align*}
	Then, we can express the kernel and image of $\phi$ as follows:
	\begin{align*}
		\ker(\phi) &= \ints\ket{
			\alpha_4
		} = \ints\ket{-2a-2b+c+2d}
		\qquad&\text{ with } \rank(\ker\phi) &= 1  \\
		\image(\phi) &= \ints\ket{
			D(1,1)\gamma_1, D(2,2)\gamma_2, D(3,3)\gamma_3
		} = \ints\ket{
			x-y, y+z, x
		} 
		\qquad&\text{ with } \rank(\image\phi) &= 3 
	\end{align*}
\end{example}


While we have established that Smith Normal Decompositions always exist for arbitrary PIDs, we have not yet discussed an algorithm for their calculation.
However, if we restrict the PID $R$ to be a Euclidean domain, we can use the following result.

\begin{proposition}\label{prop:euclidean-domain-matrices}
	Any invertible matrix over a Euclidean domain $R$ can be expressed as a finite product of elementary matrices.
	That is, $\GL(n,R)$ is generated by the elementary matrices over $R$ of degree $n$.
\end{proposition}
\remarks{
	\item 
	For a proof, see \cite[Theorem 5.2.10]{algebra:adkins}.
	Note that definitions for elementary matrices (i.e.\ elementary dilations, elementary permutations, and elementary transvections) on arbitrary Euclidean domains
	are given in \fref{defn:elementary-matrices}. 

	\item 
	The paper \textit{Products of Elementary Matrices and Non-Euclidean Principal Ideal Domains} \cite{extra:elementary-matrices-non-euclidean} conjectures that this is not generally true if the ring $R$ is a PID but not a Euclidean domain.
}

Since elementary matrices correspond to row and column operations (as discussed in \fref{prop:row-operations} for row operations and in \fref{prop:column-operations} for column operations),
it is possible to calculate SNDs by doing matrix reduction assuming the matrices are over Euclidean domains.
We will not present a general algorithm for calculating SNDs of matrices over Euclidean domains in this paper and we refer to \cite[Remark 5.3.4]{algebra:adkins} for those interested.
Below, we provide an example of this reduction process in a matrix over $\ints$.

\begin{example}
	\label{ex:snd-by-reduction}
	Let $A \in \M_{3,4}(\ints)$ be the $(3 \times 4)$-matrix given by
	\begin{equation*}
		A = \begin{pmatrix}
			1 & 2 & 0 & 1 \\
			0 & 3 & 0 & 3\\ 
			0 & 0 & 1 & 1\\
		\end{pmatrix}
	\end{equation*}
	We can calculate an SND of $A$ by the following row and column reduction operation, represented by the sequence $(A_n)$ of matrices. Let $A_0 = A$.
	The elements of $A_n$ highlighted in \redtag are elements that are to be eliminated using an elementary transvection, in \bluetag is the pivot used for said transvection, and in \orangetag the pivot multiplier.
	\def\arraystretch{0.9}
	\begin{align*}
		A_1 &:= A_0 \eladd[4]{1,2 \,; -2}
			&=& \begin{pmatrix}
				\bluemath{1} & \redmath{2} & 0 & 1 \\
				0 & 3 & 0 & 3 \\ 
				0 & 0 & 1 & 1 \\
			\end{pmatrix}
			\begin{pmatrix}
				1 & \orangemath{-2} & 0 & 0 \\
				0 & 1 & 0 & 0 \\
				0 & 0 & 1 & 0 \\
				0 & 0 & 0 & 1
			\end{pmatrix}
			&=&
			\begin{pmatrix}
				1 & 0 & 0 & {1} \\
				0 & 3 & 0 & 3 \\ 
				0 & 0 & 1 & 1 \\
			\end{pmatrix}&
			\\
		A_2 &:= A_1 \eladd[4]{1,4 \,; -1}
			&=& \begin{pmatrix}
				\bluemath{1} & 0 & 0 & \redmath{1} \\
				0 & 3 & 0 & 3 \\ 
				0 & 0 & 1 & 1 \\
			\end{pmatrix}
			\begin{pmatrix}
				1 & 0 & 0 & \orangemath{-1} \\
				0 & 1 & 0 & 0 \\
				0 & 0 & 1 & 0 \\
				0 & 0 & 0 & 1
			\end{pmatrix}
			&=&
			\begin{pmatrix}
				1 & 0 & 0 & 0 \\
				0 & 3 & 0 & 3 \\ 
				0 & 0 & 1 & 1 \\
			\end{pmatrix}&
			\\
		A_3 &:= A_2 \eladd[4]{2,4 \,; -1}
			&=& \begin{pmatrix}
				1 & 0 & 0 & 0 \\
				0 & \bluemath{3} & 0 & \redmath{3} \\ 
				0 & 0 & 1 & 1 \\
			\end{pmatrix}
			\begin{pmatrix}
				1 & 0 & 0 & 0 \\
				0 & 1 & 0 & \orangemath{-1} \\
				0 & 0 & 1 & 0 \\
				0 & 0 & 0 & 1
			\end{pmatrix}
			&=& 
			\begin{pmatrix}
				1 & 0 & 0 & 0 \\
				0 & 3 & 0 & 0 \\ 
				0 & 0 & 1 & 1 \\
			\end{pmatrix}&
			\\
		A_4 &:= A_3 \eladd[4]{3,4 \,; -1}
			&=& \begin{pmatrix}
				1 & 0 & 0 & 0 \\
				0 & 3 & 0 & 0 \\ 
				0 & 0 & \bluemath{1} & \redmath{1} \\
			\end{pmatrix}
			\begin{pmatrix}
				1 & 0 & 0 & 0 \\
				0 & 1 & 0 & 0 \\
				0 & 0 & 1 & \orangemath{-1} \\
				0 & 0 & 0 & 1
			\end{pmatrix}
			&=& 
			\begin{pmatrix}
				1 & 0 & 0 & 0 \\
				0 & 3 & 0 & 0 \\ 
				0 & 0 & 1 & 0 \\
			\end{pmatrix}&
		\end{align*}
		The matrix $A_4$ is almost is Smith Normal Form. We need to apply row and column permutations. Highlighted in \greentag are the elements that are switched in each $A_n$.
		\begin{align*}
		A_5 &:= \elswap[3]{2,3} A_4
			&=& 
			\begin{pmatrix}
				1 & 0 & 0 \\ 
				0 & 0 & 1 \\
				0 & 1 & 0
			\end{pmatrix}
			\begin{pmatrix}
				1 & 0 & 0 & 0 \\
				0 & \greenmath{3} & 0 & 0 \\ 
				0 & 0 & \greenmath{1} & 0 \\
			\end{pmatrix}
			&=&
			\begin{pmatrix}
				1 & 0 & 0 & 0 \\
				0 & 0 & 1 & 0 \\ 
				0 & 3 & 0 & 0 \\
			\end{pmatrix}&
			\\
		A_6 &:= A_5 \elswap[4]{2,3}
			&=&
			\begin{pmatrix}
				1 & 0 & 0 & 0 \\
				0 & 0 & \greenmath{3} & 0 \\ 
				0 & \greenmath{1} & 0 & 0 \\
			\end{pmatrix}
			\begin{pmatrix}
				1 & 0 & 0 & 0 \\
				0 & 0 & 1 & 0 \\
				0 & 1 & 0 & 0 \\
				0 & 0 & 0 & 1
			\end{pmatrix}
			&=& 
			\begin{pmatrix}
				1 & 0 & 0 & 0 \\
				0 & 1 & 0 & 0 \\ 
				0 & 0 & 3 & 0 \\
			\end{pmatrix}
			=:
			D
	\end{align*}
	The matrices $U \in \GL(3,\ints)$ and $V \in \GL(4,\ints)$ are given by:
	\begin{align*}
		V &:= 
		\eladd[4]{1,2 \,; -2}\,
		\eladd[4]{1, 4 \,; -1}\,
		\eladd[4]{2, 4 \,; -1}\,
		\eladd[4]{2, 4 \,; -1}\,
		\elswap[4]{2,3} 
		&=&\,
		\begin{pmatrix}
			1 & 0 & -2 & 1 \\
			0 & 0 & 1 & -1 \\
			0 & 1 & 0 & -1 \\
			0 & 0 & 0 & 1
		\end{pmatrix}
		\\
		U &:= \paren{
			\elswap[3]{2,3}
		}\inv
		&=&\,
		\begin{pmatrix}
			1 & 0 & 0 \\ 
			0 & 0 & 1 \\
			0 & 1 & 0
		\end{pmatrix}
	\end{align*}
	Observe that $U$ and $V$ are invertible since they are both a finite product of elementary matrices, which are also invertible.
	Therefore, the triple $(U,D,V)$ represent an SND of $A$.
\end{example}

The method described in \cite[Remark 5.3.4]{algebra:adkins} mostly involves the elimination of specific entries in a matrix, like in the case of matrices over $\reals$.
However, since nonzero elements of an arbitrary Euclidean domain are not generally invertible, the entries used to eliminate other entries along the same row or column (usually called a \textit{pivot}) have to be chosen with care, i.e.\ we cannot use any nonzero element as in the case of matrices over $\reals$.
The main difference lies in the family of elementary dilations over $R$ which, as described in \fref{defn:elementary-matrices}, are defined on the group of units $R^\times$ of $R$.
In the case of $R = \reals$, any nonzero element of $\reals$ can be used to construct an elementary dilation.
So, any nonzero element $a \in \reals$ of a matrix over $\reals$ can be used to eliminate any entry in the same row or column since $a$ can be reduced by $1$ using an elementary dilation by $a\inv = \frac{1}{a}$.
In contrast, the only units of $\ints$ are $\set{-1,1}$ and an element such as $2 \in \ints$ cannot be reduced to $1$.

Consequently, there is no direct analog of row or column echelon form for matrices over $\ints$ and over $\field[x]$.
In particular, even if a matrix over $\ints$ or $\field[x]$ may appear to be in row or column echelon form (as we conventionally define those), 
the pivots of that matrix generally do not correspond to the diagonal elements of its Smith Normal Form.
Note that this applies more generally for Euclidean domains that are not fields.
We provide an example of this problem below.

\begin{example}\label{ex:illegal-reduction-on-ints}
	Let $A \in M_{3,4}(\ints)$ be the matrix given by
	\begin{equation*}
		A = \begin{pmatrix}
			\redmath{2} & 0 & 3 & 0 \\
			0 & \redmath{7} & 2 & 0 \\
			0 & 0 & 0 & \redmath{3}
		\end{pmatrix}
	\end{equation*}
	Note that the matrix $A$ is in row echelon form with \redtagged{pivots} $2,7,3$ (in the row echelon sense) above highlighted in \redtag\!.
	An SND $(U,D,V)$ of $A$, calculated by Mathematica, is given below:
	\begin{equation*}
		U = I_3 
		\qquad\quad 
		D = \begin{pmatrix}
			1 & 0 & 0 & 0 \\
			0 & 1 & 0 & 0 \\
			0 & 0 & 3 & 0
		\end{pmatrix}
		\qquad\quad 
		V = \begin{pmatrix}
			-10 & -6 & 0 & -21 \\
			-2 & -1 & 0 & -4 \\
			7 & 4 & 0 & 14 \\
			0 & 0 & 1 & 0
		\end{pmatrix}
	\end{equation*}
	Observe that the nonzero diagonal elements of $D$ are $1,1,3$, not $2,7,3$.
	Note that this is because we cannot do dilation operations
	on $A$ as a matrix over $\ints$ like the one below:
	\begin{equation*}
		A^\prime = \eldilate[3]{1,\frac{1}{2}} A 
		= \begin{pmatrix}
			\frac{1}{2} & 0 & 0 \\
			0 & 1 & 0 \\
			0 & 0 & 1
		\end{pmatrix}
		\begin{pmatrix}
			\redmath{2} & 0 & 3 & 0 \\
			0 & 7 & 2 & 0 \\
			0 & 0 & 0 & 3
		\end{pmatrix}
		=
		\begin{pmatrix}
			1 & 0 & \frac{2}{3} & 0 \\
			0 & 7 & 2 & 0 \\
			0 & 0 & 0 & 3
		\end{pmatrix}
	\end{equation*}
	Note that $\frac{1}{2}$ is not an element of $\ints$ and therefore, 
	$\eldilate[3]{1,\frac{1}{2}}$ is not an elementary matrix over $\ints$.
	However, if we consider $A$ as a matrix over $\reals$, 
	we can calculate an SND $(U_\reals, D_\reals, V_\reals)$ of $A$ as follows:
	\begin{equation*}
		U_\reals = I_3 
		\qquad\quad 
		D_\reals = \begin{pmatrix}
			1 & 0 & 0 & 0 \\
			0 & 1 & 0 & 0 \\
			0 & 0 & 1 & 0
		\end{pmatrix}
		\qquad\quad 
		V_\reals = \begin{pmatrix}
			\frac{1}{2} & 0 & 0 & -\frac{3}{2} \\
			0 & \frac{1}{7} & 0 & -\frac{2}{7} \\
			0 & 0 & 0 & 1 \\
			0 & 0 & \frac{1}{3} & 0
		\end{pmatrix}
	\end{equation*}
	Observe that $V_\reals$ is not a matrix over $\ints$.
	As a sidenote, since all nonzero field elements are invertible, any nonzero diagonal element of the Smith Normal Form of any matrix over a field $\field$ can always be made into $1 \in \field$ by multiplication of an appropriate elementary dilation over $\field$
\end{example}

\clearpage

\section{Matrix Calculation of Homology of Chain Complexes}
\label{section:snd-on-ungraded-chain-complexes}

Let $C_\ast = (C_n, \boundary_n)_{n \in \ints}$ be a chain complex of free modules $C_n$ over a PID $R$ with differentials $\boundary_n: C_n \to C_{n+1}$ and assume that each $C_n$ is of finite rank.
By the assumption on $\rank(C_n)$,
	the $n$\th homology $H_n(C_\ast)$ of $C_\ast$ is finitely generated
	and therefore admits an invariant factor decomposition by the Structure Theorem (\fref{thm:structure-theorem}).
In this section, 
	we discuss how SNDs of the matrices $[\boundary_{n+1}]$ and $[\boundary_n]$
		of $\boundary_{n+1}: C_{n+1} \to C_n$ and $\boundary_{n}: C_n \to C_{n-1}$ respectively (relative to some chosen ordered bases)
	can be used to calculate these decompositions for $H_n(C_\ast)$.

The method we discussed for calculating invariant factor decompositions in \fref{section:matrix-calculation-of-IFDs} starts with a finite presentation for the $R$-module in question.
The following proposition identifies a finite presentation for $H_n(C_\ast)$ that arises naturally from the definition $H_n(C_\ast) := \ker(\boundary_n) \bigmod \im(\boundary_{n+1})$.

\begin{proposition}\label{prop:chain-homology-presentation}
	Let $C_\ast = (C_n, \boundary_n)_{n \in \ints}$ be a chain complex of free modules $C_n$ over a PID $R$ and differentials $\boundary_n: C_n \to C_{n-1}$.
	Then, the $n$\th homology group $H_n(C_\ast)$ of $C_\ast$ admits the following presentation
	\begin{equation*}
		C_{n+1}
			\,\,\Xrightarrow{\boundary_{n+1}}\,\,
		\ker(\boundary_{n})
			\,\,\Xrightarrow{\,\pi\,}\,\, 
		H_n(C_\ast)
			\,\,\Xrightarrow{\quad}\,\,
		0
	\end{equation*}
	with $\pi$ being the canonical projection onto $H_n(C_\ast)$.
\end{proposition} 
\begin{proof}
	Note that since $C_\ast$ is a chain complex, i.e.\ $\image(\boundary_{n+1}) \subseteq \ker(\boundary_n)$,
		the restriction of the codomain of $\boundary_{n+1}: C_{n+1} \to C_{n}$ onto $\ker(\boundary_n)$ is well-defined. 
	Since $\ker(\boundary_n)$ is a submodule of a free module $C_n$ over a PID $R$, $\ker(\boundary_n)$ is a free $R$-module for all $n \in \ints$
	\cite[Theorem 12.1]{algebra:dummit}.
	By definition of homology of chain complexes,
		$H_n(C_\ast) = \ker(\boundary_n) \bigmod \im(\boundary_{n+1})$.
	Then, 
		$\im(\boundary_{n+1}) = \ker(\pi)$,
		and $\im(\pi) =  \ker(H_n(C_\ast) \to 0) = H_n(C_\ast)$.
	Therefore, Equation (1) is an exact sequence of free $R$-modules and is a presentation of $H_n(C_\ast)$.
\end{proof}

We want to emphasize that, in the proposition above, 
	the codomain of the $(n+1)$\th differential $\boundary_{n+1}$ is restricted to $\ker(\boundary_n)$.
	For clarity, we write 
	$\boundary_{n+1}^\text{\,ker}: C_{n+1} \to \ker(\boundary_n)$ 
	to refer to $\boundary_{n+1}: C_{n+1} \to C_n$ with this codomain restriction.
	Since it is rarely the case that $\ker(\boundary_n) = C_n$,
		the matrices $[\boundary_{n+1}^\text{\,ker}]$ (i.e.\ with the codomain restriction)
		and $[\boundary_{n+1}]$ (i.e.\ without the codomain restriction) are generally not the same.
		For example, since $\rank(\ker\boundary_n) \leq \rank(C_n)$,
			the number of rows of $[\boundary_{n+1}^\text{\,ker}]$ 
			is less than or equal to that of $[\boundary_n]$.
		It may also be the case that the basis of $\ker(\boundary_n)$ used for $[\boundary_{n+1}^\text{\,ker}]$ is not a subset of the basis of $C_n$ used for $[\boundary_{n+1}]$.
		
This means that the results of \fref{section:matrix-calculation-of-IFDs} only immediately apply to the matrix $[\boundary_{n+1}^\text{\,ker}]$.
More specifically, if we were to use an SND of $[\boundary_{n+1}]$, 
	then said calculation would correspond to the following sequence:
	\begin{equation*}
		C_{n+1}
			\,\,\Xrightarrow{\boundary_{n+1}}\,\,
		C_n
			\,\,\Xrightarrow{\,\pi\,}\,\,
		H_n(C_\ast)
			\,\,\Xrightarrow{\quad}\,\,
		0
	\end{equation*}
Note that this sequence is generally not exact and, therefore, is not a presentation of $H_n(C_\ast)$.
Consequently, in order for $H_n(C_\ast)$ to be computable from $[\boundary_{n+1}]$, we need to determine how SNDs of 
	$[\boundary_{n+1}]$ correspond to those of 
	$[\boundary_{n+1}^\text{\,ker}]$.

Before we discuss this further, 
	we provide an example calculation of the $0$\th homology group $H_0(K) = H_0(K;\ints)$ of a simplicial complex $K$ below.
	Note that, by definition, the $0$\th boundary map $\boundary_0: C_0(K) \to 0$ has trivial codomain.
	Therefore, $\ker(\boundary_0) = C_0(K) = C_0(K;\ints)$ and the method of calculation by \fref{section:matrix-calculation-of-IFDs} is directly applicable to the matrix 
	$[\boundary_1]$.
	Note that an orientation on $K$ induces a standard ordered basis on the chain groups of $K$. 

\begin{example}
	Define the simplicial complex $K$ as illustrated below and equip $K$ with the orientation given by $\Vertex(K) = (a_1, a_2, a_3, a_4, a_5, a_6)$.
	\begin{center}
		\includegraphics[width=0.3\textwidth]{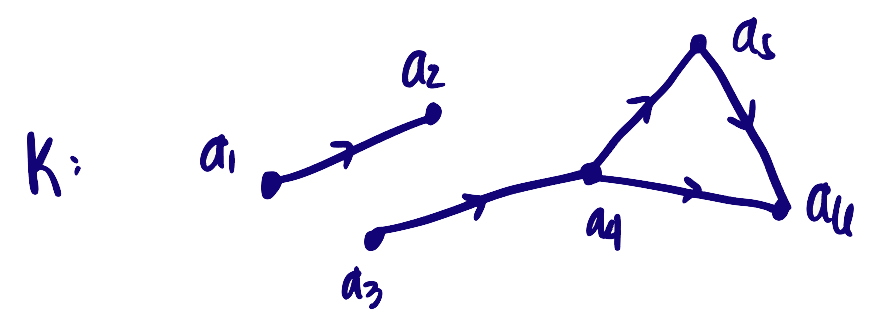}
	\end{center}
	The following sequence is a finite presentation for the $0$\th homology group $H_0(K) = H_0(K; \ints)$ of $K$ with coefficients in $\ints$:
	\begin{equation*}
		C_{1}
			\,\,\Xrightarrow{\quad\boundary_{1}\quad}\,\,
		\equalsupto{\mathclap{\ker(\boundary_0)}}{C_0}
			\,\,\Xrightarrow{\quad\pi\quad}\,\,
		H_0(K)
			\,\,\Xrightarrow{\qquad}\,\,
		0
	\end{equation*}
	Therefore, we can calculate $H_0(K)$ using an SND $(U_1, D_1, V_1)$ of the matrix
	$[\boundary_1]$ of $\boundary_1: C_1(K) \to C_0(K)$ relative to some ordered bases.
	The orientation on $K$ by $\Vertex(K) = (a_1, a_2, a_3, a_4, a_5, a_6)$ determines the following ordered bases of $C_0(K)$ and $C_1(K)$:
	\begin{align*}
		\text{ $0$-simplices: } (a_1,\, a_2,\, a_3,\, a_4,\, a_5,\, a_6) 
		\qquad\text{ and }\qquad
		\text{ $1$-simplices: } (a_1a_2,\, a_3a_4,\, a_4a_5,\, a_4a_6,\, a_5a_6)
	\end{align*}
	Then, the matrix $[\boundary_1] \in M_{6,5}(\ints)$ relative to the standard ordered bases can be calculated as follows:
	\vspace{-\parskip}
	\begin{equation*}\def\arraystretch{0.9}
		[\boundary_1] =
		\begin{NiceArray}{>{\color{gray}}c ccc cc}
		\RowStyle[]{\color{gray}}
			{} & a_1a_2 & a_3a_4 & a_4a_5 & a_4a_6 & a_5a_6 \\
			{a_1} 	& -1 & 0 & 0 & 0 & 0 \\
			{a_2} 	& 1 & 0 & 0 & 0 & 0 \\
			{a_3} 	& 0 & -1 & 0 & 0 & 0 \\
			{a_4} 	& 0 & 1 & -1 & -1 & 0 \\
			{a_5} 	& 0 & 0 & 1 & 0 & -1 \\
			{a_6} 	& 0 & 0 & 0 & 1 & 1
		\CodeAfter \SubMatrix({2-2}{7-6})
		\end{NiceArray}
	\end{equation*}
	Given below is an SND $(U_1, D_1, V_1)$ of $[\boundary_1]$ calculated using Mathematica. 
	Since $C_0 = \ker(\boundary_0)$, all six columns of $U_1$ correspond to a basis of $\ker(\boundary_0)$. 
	Highlighted in \greentag are the columns of $U_1 \in \GL(6,\ints)$ corresponding to \greentagged{torsion} summands of $H_0(K)$ and their corresponding diagonal elements in $D_1$
	and in \bluetag are the columns of $U_1$ corresponding to \bluetagged{free} summands of $H_0(K)$.
	\begin{equation*}\def\arraystretch{0.9}
		U_1 = 
		\begin{NiceArray}{ 
			>{\color{gray}}c ( ccc ccc ) 
		}
			\CodeBefore [create-cell-nodes]
				\tikz \node [green-cell = (1-2) (6-2)] {} ;
				\tikz \node [green-cell = (1-3) (3-3) (6-3)] {} ;
				\tikz \node [green-cell = (1-4) (4-4) (6-4)] {} ;
				\tikz \node [green-cell = (1-5) (4-5) (6-5)] {} ;
				\tikz \node [blue-cell = (1-6) (6-6)] {} ;
				\tikz \node [blue-cell = (1-7) (6-7)] {} ;
			\Body
			{a_1} 	& -1 & 0 & 0 & 0 & 0 & 0	\\
			{a_2} 	& 1 & 0 & 0 & 0 & 1 & 0		\\
			{a_3} 	& 0 & -1 & 0 & 0 & 0 & 0	\\
			{a_4} 	& 0 & 1 & -1 & -1 & 0 & 0	\\
			{a_5} 	& 0 & 0 & 1 & 0 & 0 & 0		\\
			{a_6} 	& 0 & 0 & 0 & 1 & 0 & 1
		\end{NiceArray} 
		\qquad 
		D_1 = \begin{NiceArray}{( ccc cc )}
			\CodeBefore [create-cell-nodes]
				\tikz \node [green-cell = (1-1)] {} ;
				\tikz \node [green-cell = (2-2)] {} ;
				\tikz \node [green-cell = (3-3)] {} ;
				\tikz \node [green-cell = (4-4)] {} ;
			\Body
			1 & 0 & 0 & 0 & 0 \\
			0 & 1 & 0 & 0 & 0 \\
			0 & 0 & 1 & 0 & 0 \\
			0 & 0 & 0 & 1 & 0 \\
			0 & 0 & 0 & 0 & 0 \\
			0 & 0 & 0 & 0 & 0 \\
		\end{NiceArray}
		\qquad 
		V_1 = \begin{NiceArray}{ >{\color{gray}}c ( ccc cc )  }
			{a_1 a_2} 	& 1 & 0 & 0 & 0 & 0		\\
			{a_3 a_4} 	& 0 & 1 & 0 & 0 & 0		\\
			{a_4 a_5} 	& 0 & 0 & 1 & 0 & 1		\\
			{a_4 a_6} 	& 0 & 0 & 0 & 1 & -1		\\
			{a_5 a_6} 	& 0 & 0 & 0 & 0 & 1		\\
		\end{NiceArray}
	\end{equation*}
	Then, by \fref{prop:how-to-get-ifd}, we can calculate $H_0(K)$ as follows:
	\begin{align*}
		H_0(K)
		= \frac{C_0(K)}{\image(\boundary_1)}
		&\cong \frac{
			\ints\big\lket
				\greenmath{a_2 - a_1},
				\greenmath{a_4 - a_5},
				\greenmath{a_5 - a_4},
				\greenmath{a_6 - a_4},
				\bluemath{a_2},
				\bluemath{a_6}
			\big\rket
		}{
			\ints\big\lket
				\greenmath{(1)(a_2 - a_1)},
				\greenmath{(1)(a_4 - a_3)},
				\greenmath{(1)(a_5 - a_4)},
				\greenmath{(1)(a_6 - a_4)}
			\big\rket
		}
		\\[1pt]
		&= 
		\greenmath{
			\frac{\ints\ket{a_2 - a_1}}{\ints\ket{a_2 - a_1}}
		} \oplus 
		\greenmath{
			\frac{\ints\ket{a_4-a_3}}{\ints\ket{a_4-a_3}}
		} \oplus 
		\greenmath{
			\frac{\ints\ket{a_5-a_4}}{\ints\ket{a_5 - a_4}}
		} \oplus 
		\greenmath{
			\frac{\ints\ket{a_6-a_4}}{\ints\ket{a_6-a_4}}
		} \oplus 
		\bluemath{ \ints\ket{a_2} } \oplus 
		\bluemath{ \ints\ket{a_6} } 
		\\[1pt]
		&= \bluemath{ \ints\ket{a_2} } \oplus 
			\bluemath{ \ints\ket{a_6} } 
		\cong \ints^2
	\end{align*}
	Observe that this result matches with the expected interpretation of the $0$\th homology group of $K$,
		wherein $K$ has two path components represented by the homology classes $[a_2]$ and $[a_6]$ as illustrated below:
	\begin{center}
		\includegraphics[width=0.3\textwidth]{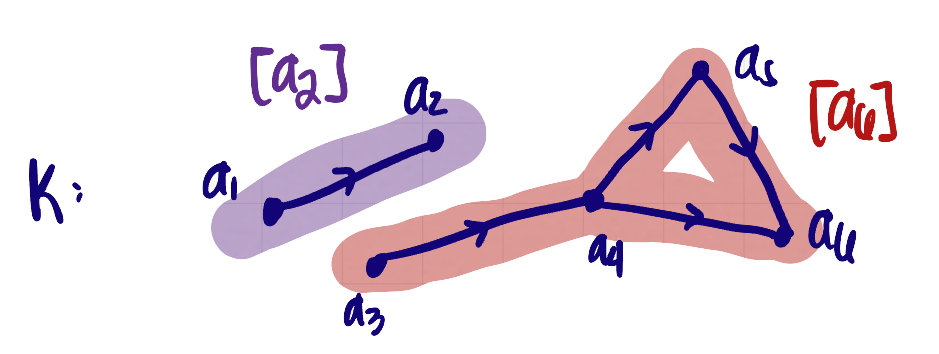}
	\end{center}
\end{example}

\spacer
 
We now
consider $H_n(C_\ast)$ and the matrix $[\boundary_{n+1}\upker]$ corresponding to the presentation of $H_n(C_\ast)$ by \fref{prop:chain-homology-presentation}, 
	i.e.\ $\boundary_{n+1}\upker$ refers to $\boundary_{n+1}$ with codomain restricted to  $\ker(\boundary_n)$.
Following the calculation presented in \fref{section:matrix-calculation-of-IFDs},
	an SND $(U_{n+1}\upker, D_{n+1}\upker, V_{n+1}\upker)$ of 
	$[\boundary_{n+1}^\text{\,ker}]$ determines a basis $B =( \beta_1, \ldots, \beta_{r+f})$ of $\ker(\boundary_n)$ with $\rank(\ker\boundary_n) = r+f$ by $[\beta_j] = \col_j(U_{n+1}\upker)$
	and nonzero elements $d_1, \ldots, d_r \in R$ by $d_j = D_{n+1}\upker(j,j)$ 
	such that $\set{d_1\beta_1, \ldots, d_r\beta_r}$ is a basis for $\im(\boundary_{n+1})$.

Observe that we can express the partition of the basis $B$ into $(\beta_1, \ldots, \beta_r)$ and $(\beta_{r+1}, \ldots, \beta_m)$ as a decomposition 
	$\ker(\boundary_n) = K_n^\text{tor} \oplus K_n^\text{free}$ 
of $\ker(\boundary_n)$ where $K_n^\text{tor}$ and $K_n^\text{free}$ are free $R$-submodules given by:
\begin{equation*}
	K_n^\text{tor} := R\ket{ \beta_1, \ldots, \beta_r }
	\qquad\text{ and }\qquad 
	K_n^\text{free} := R\ket{ \beta_{r+1}, \ldots, \beta_m }
\end{equation*}
As the notation suggests, 
	the submodule $K_n^\text{tor}$ accounts for all elements in 
	$\ker(\boundary_n)$ that map to torsion elements or trivial in $H_n(C_\ast)$
	and the submodule 
		$K_n^\text{free}$ accounts for the elements of $\ker(\boundary_n)$ that remain free in $H_n(C_\ast)$.
We illustrate this below 
	with the \greentagged{torsion} component of $H_n(C_\ast)$ highlighted in \greentag and \bluetagged{free} component of $H_n(C_\ast)$ in \bluetag\!.
\begin{align*}
	H_n(C_\ast)
	&\cong \frac{\ker(\boundary_n)}{\im(\boundary_{n+1})}
		\cong \frac{ 
				{K_n^\text{tor}} \oplus {K_n^\text{free}}
			}{ \im(\boundary_{n+1}) }
		\cong \greenmath{\frac{ K_n^\text{tor} }{ \im(\boundary_{n+1}) }}
				\oplus 
				\bluemath{K_n^\text{free}}
	\\ 
	&\cong 
		\greenmath{ \frac{R\ket{\beta_1}}{R\ket{d_1\beta_1}} } 
		\oplus \greenmath{\cdots} \oplus 
		\greenmath{ \frac{R\ket{\beta_r}}{R\ket{d_r\beta_r}} } 
		\oplus 
		\bluemath{ R\ket{\beta_{r+1}} } \oplus \bluemath{\cdots} \oplus 
			\bluemath{ R\ket{\beta_{m}} }
\end{align*}
In order for $H_n(C_\ast)$ to be computable from $[\boundary_{n+1}]$,
	we need to show that a basis $B = (\beta_1, \ldots, \beta_m)$ of $\ker(\boundary_n)$ and elements $d_1, \ldots, d_r \in R$ as described above 
	can be calculated from some SND $(U_{n+1}, D_{n+1}, V_{n+1})$ of $[\boundary_{n+1}]$
	with $U_{n+1}$ inducing a basis on $C_n$, not on $\ker(\boundary_n)$.
The existence of such an SND of $[\boundary_{n+1}]$ is given by the following result, 
taken from \cite[Theorem 11.4]{algtopo:munkres}.

\begin{theorem}\label{thm:chain-complex-decomposition}
	Let $C_\ast = (C_n, \boundary_n)$ be a chain complex of free 
	modules $C_n$ over a PID $R$ with differentials $\boundary_n: C_n \to C_{n-1}$ and assume each $C_n$ is of finite rank.
	For each $n \in \nonnegints$, 
		$C_n$ decomposes into the following direct sum of free $R$-modules:
		\begin{equation*}\tag{i}
			C_n = K_n^\text{tor} \oplus K_n^\text{free} \oplus 
				\frac{C_n}{\ker(\boundary_n)}
		\end{equation*}
		where $\ker(\boundary_n) = K_n^\text{free} \oplus K_n^\text{tor}$,  $\image(\boundary_{n+1}) \subseteq K_n^\text{tor}$,
		and $K_n^\text{tor} \bigmod \image(\boundary_{n+1})$ is the torsion component of $H_n(C_\ast)$.
		Furthermore,
			there exists an SND $(U_{n+1}, D_{n+1}, V_{n+1})$ of $[\boundary_{n+1}]$
			such that the basis $B = (\beta_1, \ldots, \beta_q)$ of $C_n$ 
			by $[\beta_j] = \col_j(U_{n+1})$ partitions into three sets of bases as follows:
			\begin{equation*}
				K_n^\text{tor} = R\ket{\, \beta_1, \ldots, \beta_r \,}
				\quad\text{ , }\quad
				K_n^\text{free} = R\ket{\, \beta_{r+1}, \ldots, \beta_{r+f} \,}
				\quad\text{ and }\quad
				\frac{C_n}{\ker(\boundary_n)}
					= R\ket{\, \beta_{r+f+1}, \ldots, \beta_q \,}
			\end{equation*}
		Note that $\im(\boundary_{n+1}) = R\ket{\, d_1\beta_1, \ldots, d_r\beta_r \,}$ with $d_1, \ldots, d_r \in R$ corresponding to the nonzero entries of $D_n$. 
\end{theorem}
\remark{
	For a proof, see \cite[Theorem 11.4]{algtopo:munkres} where they write $W_n$, $V_n$, and $U_n$ to refer to $K_n^\text{tor}$, $K_n^\text{free}$, and $C_n \bigmod \ker(\boundary_n)$ respectively.
	Note that since $K_n^\text{tor} \bigmod \im(\boundary_{n+1})$ has to be the torsion component of $H_n(C_\ast)$, we have the following by definition of torsion component:
	\begin{equation*}
		H_n(C_\ast) = \frac{
			\ker(\boundary_n)
		}{
			\image(\boundary_{n+1})
		}
		\cong 
		K_n^\text{free}
		\oplus
		\frac{
			K_n^\text{tor}
		}{
			\image(\boundary_{n+1})
		}\,.
	\end{equation*}
	with $K_n^\text{free} = R\ket{\beta_{r+1}, \ldots, \beta_{r+f}} \cong R^f$ consists of the torsion-free elements of $H_n(C_\ast)$.
}

The existence of this decomposition of $C_n$ into direct summands
guarantees that the information we get 
	about $\im(\boundary_{n+1})$ from an SND of 
		$[\boundary_{n+1}]$
	is compatible
	to that about $\ker(\boundary_n)$ from an SND of $[\boundary_{n}]$.
\newcommand{\phan}{\hphantom{v}}%
Let $\rank(C_n) = q$ and let the ranks of the direct summands of $C_n$ by \fref{thm:chain-complex-decomposition} be given as follows:
\begin{equation*}
	\greenmath{
		\rank(K_n^\text{tor}) = r
	}
	\quad;\quad
	\bluemath{
		\rank(K_n^\text{free}) = f
	}
	\quad;\quad
	\redmath{
		\rank\bigl( C_n \bigmod \ker(\boundary_n) \bigr) = q
	}
\end{equation*}
Let $(W_{n+1}, D_{n+1}, V_{n+1})$ be an arbitrary SND of $[\boundary_{n+1}]$.
The first $r$ columns of the matrix $W_{n+1}D_{n+1}$
	determines a basis $(d_1\alpha_1, \ldots, d_r\alpha_r)$ 
	of $\im(\boundary_{n+1})$
	(or equivalently, the first $r$ columns of $W_{n+1}$ and the first $r$ rows of $D_{n+1}$).
Let $A = (\alpha_1, \ldots, \alpha_m)$ be the basis of $C_n$ given by $[\alpha_j] = \col_j(W_{n+1})$.
Part of \fref{thm:chain-complex-decomposition}
	tells us that the basis $A$ of $C_n$ partitions into two sets:
	the subset 
	\greentagged{$(\alpha_1, \ldots, \alpha_r)$} is a basis of the direct summand 
	\greentagged{$K_n^\text{tor}$} of $C_n$	
	and the other subset 
	\purpletagged{$(\alpha_{r+1}, \ldots, \alpha_{m})$} form a basis of the direct sum \purpletagged{$K_n^\text{free} \oplus (C_n \bigmod \ker\boundary_n)$},
	as illustrated below:
	\vspace{30pt}
	\begin{equation*}
		\begin{NiceArray}{ ccc c ccc } 
		\CodeBefore [create-cell-nodes]
			\tikz \node [green-cell = (1-1) (2-1) (4-1), inner sep=2pt] {} ;
			\tikz \node [green-highlight = (2-2), inner sep=2pt] {} ;
			\tikz \node [green-cell = (1-3) (2-3) (4-3), inner sep=2pt] {} ;
			\tikz \node [purple-cell = (1-5) (2-5) (4-5), inner sep=2pt] {} ;
			\tikz \node [purple-highlight = (2-6), inner sep=2pt] {} ;
			\tikz \node [purple-cell = (1-7) (2-7) (4-7), inner sep=2pt] {} ;
		\Body
			\vdots & {} & \vdots && \vdots & {} & \vdots \\
			[\alpha_1] & \cdots & [\alpha_r]
				&& [\alpha_{r+1}] & \cdots & [\alpha_{m}] \\
			\vdots & {} & \vdots && \vdots & {} & \vdots \\[-7pt] 
			{\phan} & {\phan} & {\phan} && {\phan} & {\phan} & {\phan} \\
		\CodeAfter
			\SubMatrix({1-1}{3-7})[left-xshift=3pt, right-xshift=3pt]
			\OverBrace{1-1}{1-7}{
				W_{n+1} \in \GL(m, R)
			}[yshift=10pt] 
			\UnderBrace{4-1}{4-3}{\tikz \node [green-node-padded, yshift=-10pt] {
				$K_n^\text{tor}$
			};}[yshift=5pt, color=ForestGreen, shorten]
			\UnderBrace{4-5}{4-7}{\tikz \node [purple-node-padded, yshift=-10pt] {
				$K_n^\text{free} \oplus (C_n \bigmod \ker\boundary_n)$
			};}[yshift=5pt, color=violet, shorten]
		\end{NiceArray}
		\hspace{13pt}\raisebox{20pt}{$\scriptstyle -1$}\,
		{\color{gray} [\boundary_{n+1}]V_{n+1}} 
		= \hspace{13pt}
		\begin{NiceArray}{ ccc ccc }
		\CodeBefore [create-cell-nodes]
			\tikz \node [green-cell = (1-1)] {} ;
			\tikz \node [green-cell = (3-3)] {} ;
		\Body
			d_1 			\\[-4pt] 
			 & \ddots 		\\
			&& d_r 		\\[0pt]
			&&& 0 & \cdots & 0 			\\[-3pt]
			&&& \vdots & \ddots &\vdots 	\\
			&&& 0 & \cdots & 0					
		\CodeAfter
			\SubMatrix({1-1}{6-6})[left-xshift=3pt, right-xshift=3pt]
			\OverBrace{1-1}{1-6}{
				D_{n+1} \in \M_{q,(-)}(R)
			}[yshift=10pt] 
		\end{NiceArray} 
	\vspace{30pt}
	\end{equation*} 
Note that, in general, the basis $(\alpha_{r+1}, \ldots, \alpha_m)$ of $K_n^\text{free} \oplus (C_n \bigmod \ker(\boundary_n))$ generally does not partition into bases for each summand.

In contrast, 
let $(U_n, D_n, V_n)$ be an arbitrary SND of $[\boundary_n]$.
The matrix $V_n \in \GL(m,R)$ determines a basis 
$K = (\kappa_1, \ldots, \kappa_m)$ of $C_n$ by $[\kappa_i] = \col_i(V_n)$ such that 
	$(\kappa_1, \ldots, \kappa_q)$ is a basis of $C_n \bigmod \ker\boundary_n$ 
	and 
	$(\kappa_{q+1}, \ldots, \kappa_m)$ is a basis of $\ker(\boundary_n) = K_n^\text{tor} \oplus K_n^\text{free}$,
	as illustrated below.
	\vspace{30pt}
	\begin{equation*}
		{\color{gray} (U_n)[\boundary_n]} \hspace{13pt}
		\begin{NiceArray}{ ccc c ccc } 
		\CodeBefore [create-cell-nodes]
			\tikz \node [red-cell = (1-1) (2-1) (4-1), inner sep=2pt] {} ;
			\tikz \node [red-highlight = (2-2), inner sep=2pt] {} ;
			\tikz \node [red-cell = (1-3) (2-3) (4-3), inner sep=2pt] {} ;
			\tikz \node [orange-cell = (1-5) (2-5) (4-5), inner sep=2pt] {} ;
			\tikz \node [orange-highlight = (2-6), inner sep=2pt] {} ;
			\tikz \node [orange-cell = (1-7) (2-7) (4-7), inner sep=2pt] {} ;
		\Body
			\vdots & {} & \vdots && \vdots & {} & \vdots \\
			[\kappa_1] & \cdots & [\kappa_q]
				&& [\kappa_{q+1}] & \cdots & [\kappa_{m}] \\
			\vdots & {} & \vdots && \vdots & {} & \vdots \\[-7pt] 
			{\phan} & {\phan} & {\phan} && {\phan} & {\phan} & {\phan} \\
		\CodeAfter
			\SubMatrix({1-1}{3-7})[left-xshift=3pt, right-xshift=3pt]
			\OverBrace{1-1}{1-7}{
				V_n \in \GL(m, R)
			}[yshift=10pt] 
			\UnderBrace{4-1}{4-3}{\tikz \node [red-node-padded, yshift=-10pt] {
				$C_n \bigmod \ker(\boundary_n)$
			};}[yshift=5pt, color=BrickRed, shorten]
			\UnderBrace{4-5}{4-7}{\tikz \node [orange-node-padded, yshift=-10pt] {
				$\ker(\boundary_n) = K_n^\text{tor} \oplus K_n^\text{free}$
			};}[yshift=5pt, color=BurntOrange, shorten]
		\end{NiceArray}
		\hspace{13pt} = \hspace{13pt}
		\begin{NiceArray}{ ccc ccc }
		\CodeBefore [create-cell-nodes]
			\tikz \node [red-cell = (1-1)] {} ;
			\tikz \node [red-cell = (3-3)] {} ;
		\Body
			(*)_1 			\\[-4pt] 
			 & \ddots 		\\
			&& (*)_q 		\\[0pt]
			&&& 0 & \cdots & 0 			\\[-3pt]
			&&& \vdots & \ddots &\vdots 	\\
			&&& 0 & \cdots & 0					
		\CodeAfter
			\SubMatrix({1-1}{6-6})[left-xshift=3pt, right-xshift=3pt]
			\OverBrace{1-1}{1-6}{
				D_n \in \M_{(-),q}(R)
			}[yshift=10pt] 
		\end{NiceArray} 
	\vspace{25pt}
	\end{equation*} 
The problem here is that the basis 
	$(\kappa_{q+1}, \ldots, \kappa_m)$
	generally does not partition into a basis of $K_n^\text{tor}$ and a basis of $K_n^\text{free}$.

What \fref{thm:chain-complex-decomposition} guarantees is the existence of an SND $(U_{n+1}, D_{n+1}, V_{n+1})$ of $[\boundary_{n+1}]$ wherein a basis of $K_n^\text{free}$ can be identified, possibly in comparison with another SND of $[\boundary_n]$.
Then, rows of the SNF $D_{n+1}$ of $[\boundary_{n+1}]$ corresponding to basis elements of $C_n \bigmod \ker\boundary_n$ can be removed to get the SNF $D_{n+1}\upker$ of $[\boundary_{n+1}\upker]$, the matrix of $\boundary_{n+1}$ with the codomain restriction into $\ker(\boundary_n)$.
We provide an example of such an SND below,
		in a calculation of the $1$\st homology group of a simplicial complex.

\begin{example}\label{ex:matrix-firsthom-sadness}
	Let $K$ be the simplicial complex illustrated below with orientation $\Vertex(K) = (a_1, a_2, a_3, a_4)$.
	\begin{center}
		\includegraphics[width=0.15\textwidth]{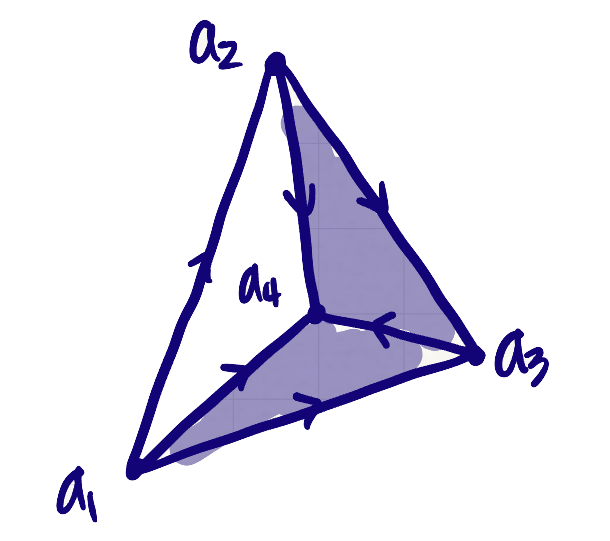}
	\end{center}
	The $1$\st homology group $H_1(K) = H_1(K;\ints)$ of $K$ can be calculated from a Smith Normal Decomposition of 
	$[\boundary_{2}]$ with $\boundary_2: C_2(K) \to C_1(K)$.
	The orientation on $K$ induces a standard ordered bases on the chain groups of $K$. We identify these bases for $C_0(K)$, $C_1(K)$, and $C_2(K)$ below.
	\begin{equation*}\def\arraystretch{1.2}
		\begin{array}{ccc}
			\text{four 0-simplices} &:&
				(a_1, a_2, a_3, a_4) \\
			\text{six 1-simplices} &:& 
				(a_1a_2, a_1a_3, a_1a_4, a_2a_3, a_2a_4 a_3a_4) \\
			\text{two 2-simplices} &:& 
				(a_1a_3a_4, a_2a_3a_4)
		\end{array}
	\end{equation*}
	The matrices $[\boundary_1]$ and $[\boundary_2]$ of the boundary maps $\boundary_1:C_1(K) \to C_0(K)$ and 
		$\boundary_2: C_2(K) \to C_1(K)$ respectively relative to these bases are given below.
	\begin{equation*}\arraycolsep=0.7\arraycolsep
		[\boundary_1] = \hspace{3pt}\begin{NiceArray}{
			>{\color{gray}}c<{\hspace{7pt}} 
		ccc ccc}
			\RowStyle{\color{gray}}
			{} & a_1a_2 & a_1a_3 & a_1a_4 & a_2a_3 & a_2a_4 & a_3a_4 \\
			{a_1}	& -1 & -1 & -1 & 0 & 0 & 0 \\
				{a_2}	& 1 & 0 & 0 & -1 & -1 & 0 \\
				{a_3}	& 0 & 1 & 0 & 1 & 0 & -1 \\
				{a_4}	& 0 & 0 & 1 & 0 & 1 & 1 
			\CodeAfter
			\SubMatrix({2-2}{5-7})[left-xshift=3pt, right-xshift=3pt]
		\end{NiceArray}
		\hspace{12pt}\qquad\text{ and }\qquad
		[\boundary_2] = \hspace{3pt}\begin{NiceArray}{
			>{\color{gray}}c<{\hspace{7pt}} cc
		}
			\RowStyle{\color{gray}}
			{} & a_1a_3a_4 & a_2a_3a_4 \\
			{a_1a_2} & 0 & 0 \\
			{a_1a_3} & 1 & 0 \\
			{a_1a_4} & -1 & 0 \\
			{a_2a_3} & 0 & 1 \\
			{a_2a_4} & 0 & -1 \\
			{a_3a_4} & 1 & 1 \\
			\CodeAfter
			\SubMatrix({2-2}{7-3})[left-xshift=3pt, right-xshift=3pt]
		\end{NiceArray}
	\end{equation*}
	Given below is an SND $(U_2, D_2, V_2)$ of $[\boundary_2]$,
		with the columns of $U_2$ corresponding to 
		\greentagged{$K_1^\text{tor}$} (following the notation in \fref{thm:chain-complex-decomposition})
		highlighted in \greentag.
	\begin{equation*}
		U_2 = \hspace{3pt}
		\begin{NiceArray}{>{\color{gray}}c<{\,\,} ccc ccc}
			\CodeBefore [create-cell-nodes]
				\tikz \node [green-cell = (1-2) (3-2) (6-2), inner xsep=2pt] {} ;
				\tikz \node [green-cell = (1-3) (5-3) (6-3), inner xsep=2pt] {} ;
			\Body
			{a_1 a_2} 	& 0 & 0  & 1  & 0 & 1 & 0		\\
			{a_1 a_3} 	& 1 & 0  & 0  & 0 & 0 & 0		\\
			{a_1 a_4} 	& -1 & 0 & -1 & 1 & 0 & 0		\\
			{a_2 a_3} 	& 0 & 1  & 0  & 0 & 0 & 0		\\
			{a_2 a_4} 	& 0 & -1 & 1  & 0 & 0 & 0		\\
			{a_3 a_4} 	& 1 & 1  & 0  & 0 & 0 & 1
			\CodeAfter
			\SubMatrix({1-2}{6-7})[left-xshift=3pt, right-xshift=3pt]
		\end{NiceArray}
		\qquad\quad 
		D_2 = \hspace{12pt}\begin{NiceArray}{cc}
			\CodeBefore [create-cell-nodes]
				\tikz \node [green-cell = (1-1)] {} ;
				\tikz \node [green-cell = (2-2)] {} ;
			\Body
			1 & 0 \\
			0 & 1 \\
			0 & 0 \\
			0 & 0 \\
			0 & 0 \\
			0 & 0
			\CodeAfter
			\SubMatrix({1-1}{6-2})[left-xshift=3pt, right-xshift=3pt]
		\end{NiceArray}
		\qquad\quad 
		V_2 = I_2
	\end{equation*}
	Columns 3 to 6 of $U_2$ determine a basis of $K_1^\text{free} \oplus (C_1 \bigmod \ker\boundary_1)$.
	We can confirm that the basis $B = (\beta_1, \ldots, \beta_6)$ of $C_1$ by $[\beta_i] = \col_i(U_2)$ partitions into bases for $K_1^\text{tor}$, $K_1^\text{free}$, and $C_1 \bigmod \ker\boundary_1$ as given in \fref{thm:chain-complex-decomposition} by considering the SND $(U_1, D_1, V_1)$ of $[\boundary_1]$ given as follows:
	\vspace{25pt}
	\begin{equation*}
		U_1 = \begin{NiceArray}{>{\color{gray}}c (cccc)}
			{a_1} 	& -1 & -1 & 0 & 0	\\
			{a_2} 	& 0 & 1 & 0 & 0		\\
			{a_3} 	& 0 & 1 & -1 & 0		\\
			{a_4} 	& 1 & 1 & 1 & 1
		\end{NiceArray}
		\qquad 
		D_1 = \hspace{9pt}\begin{NiceArray}{ccc ccc}
		\CodeBefore [create-cell-nodes]
			\tikz \node [red-cell = (1-1), inner xsep=4pt] {} ;
			\tikz \node [red-cell = (2-2), inner xsep=4pt] {} ;
			\tikz \node [red-cell = (3-3), inner xsep=4pt] {} ;
		\Body
			1 & 0 & 0 & 0 & 0 & 0 \\
			0 & 1 & 0 & 0 & 0 & 0 \\
			0 & 0 & 1 & 0 & 0 & 0 \\
			0 & 0 & 0 & 0 & 0 & 0 
		\CodeAfter
			\SubMatrix({1-1}{4-6})[left-xshift=3pt, right-xshift=3pt]
		\end{NiceArray}
		\hspace{9pt}
		\qquad 
		V_1 = \hspace{3pt}\begin{NiceArray}{>{\color{gray}}c !{\,\,}
			ccc ccc 
		}
		\CodeBefore [create-cell-nodes]
			\tikz \node [red-cell = (1-2) (6-2), inner xsep=4pt] {} ;
			\tikz \node [red-cell = (1-3) (6-3), inner xsep=4pt] {} ;
			\tikz \node [red-cell = (1-4) (6-4), inner xsep=4pt] {} ;
			\tikz \node [green-cell = (1-5) (3-5) (6-5), inner xsep=2pt] {} ;
			\tikz \node [green-cell = (1-6) (5-6) (6-6), inner xsep=2pt] {} ;
			\tikz \node [blue-cell = (1-7) (3-7) (6-7), inner xsep=2pt] {} ;
		\Body
			{a_1 a_2} 	& 0 & 1 & 0 & 0 & 0 & 1		\\
			{a_1 a_3} 	& 0 & 0 & 0 & 1 & 0 & 0		\\
			{a_1 a_4} 	& 1 & 0 & 0 & -1 & 0 & -1		\\
			{a_2 a_3} 	& 0 & 0 & 0 & 0 & 1 & 0		\\
			{a_2 a_4} 	& 0 & 0 & 0 & 0 & -1 & 1		\\
			{a_3 a_4} 	& 0 & 0 & 1 & 1 & 1 & 0
		\CodeAfter
			\SubMatrix({1-2}{6-7})[left-xshift=3pt, right-xshift=3pt]
			\UnderBrace{6-2}{6-4}{\tikz \node [red-node-padded, yshift=-10pt] {
					$C_1 \bigmod \ker\boundary_1$
				};}[yshift=5pt, shorten]
			\OverBrace{1-5}{1-6}{\tikz \node [green-node-padded, yshift=8pt] {
					$K_1^\text{tor}$
				};}[yshift=7pt, shorten]
			\UnderBrace{6-7}{6-7}{\tikz \node [blue-node-padded, yshift=-10pt] {
					$K_1^\text{free}$
				};}[yshift=5pt, shorten]
		\end{NiceArray}
		\vspace{30pt}
	\end{equation*} 
	Observe that the matrix $U_2$ is the matrix $V_1$ up to column permutation.
	The basis $B = (\beta_1, \ldots, \beta_6)$ of $C_1$ given by $[\beta_i] = \col_i(U_2)$ is described below, grouped relative to the decomposition of $C_1$.
	\begin{equation*}
		\begin{NiceArray}{rl}
			\Block{2-1}{ \greenmath{K_1^\text{tor}} :}
				& \sigma_1 = a_1 a_3 - a_1a_4 + a_3a_4 
					\text{ with } d_1 = D_2(1,1) = 1
				\\
				& \sigma_2 = a_2a_3 - a_2a_4 + a_3a_4 
					\text{ with } d_2 = D_2(2,2) = 1
				\\[2pt]
			\Block{}{ \bluemath{K_1^\text{free}} :}
				& \sigma_3 = a_1a_2 - a_1a_4 + a_2a_4 \\
			\Block{3-1}{ \redmath{C_1 \bigmod \ker\boundary_1} :}
				& \sigma_4 = a_1a_4 \\
				& \sigma_5 = a_1 a_2 \\
				& \sigma_6 = a_3a_4
		\end{NiceArray}
	\end{equation*}
	We can calculate $H_1(K)$ by disregarding the rows of $D_2$ corresponding to the basis elements of $C_1 \bigmod \ker\boundary_1$, i.e.\ by only considering the basis elements of $K_1^\text{free}$ and $K_1^\text{tor}$ as follows:
	\begin{align*}
		H_1(K)
		= 
		\frac{\ker\boundary_1}{\im\boundary_2}
		\cong 
			\greenmath{\frac{K_1^\text{tor}}{\im(\boundary_1)}}
			\oplus
			\bluemath{K_1^\text{free}}
		&\cong 
			\greenmath{\frac{
				\ints\ket{\sigma_1}
			}{
				\ints\ket{d_1\sigma_1}
			}}
			\oplus 
			\greenmath{\frac{
				\ints\ket{\sigma_2}
			}{
				\ints\ket{d_2\sigma_2}
			}}
			\oplus 
			\bluemath{\ints\ket{\sigma_3}}
		\\ 
		&\cong 
			\bluemath{\ints\ket{a_1a_2 - a_1a_4 + a_2a_4}}
			\cong \ints
	\end{align*}
	Note that both summands of $K_1^\text{tor}$ map to trivial elements in $H_1(K)$ since both $\sigma_1, \sigma_2 \in \im\boundary_2$.
	The cycle representative $\sigma_3$ of the homology class $[\sigma_3] \in H_1(K)$ of $K$ is illustrated below:
	\begin{center}
		\includegraphics[width=0.3\textwidth]{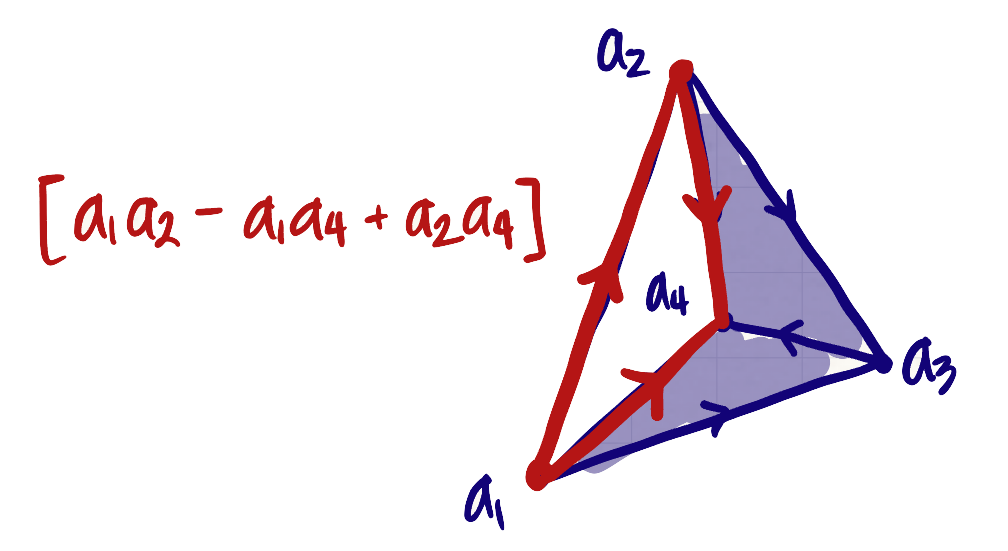}
	\end{center}
\end{example}

Once the decomposition of $C_n$ into the three direct summands
is established, 
	an SND $(U_{n+1}, D_{n+1}, V_{n+1})$ of $[\boundary_{n+1}]$ with the required properties can be generated from any SND $(W_{n+1}, D_{n+1}, V_{n+1})$ of $[\boundary_{n+1}]$.
The first $r$ columns of $W_{n+1}$, as discussed earlier, 
	determines a basis
	$(\alpha_1, \ldots, \alpha_r)$ of 
	\greentagged{$K_n^\text{tor}$} by $[\alpha_i] = \col_i(W_{n+1})$.
Since $K_n^\text{free}$ and $C_n \bigmod \ker(\boundary_n)$ are direct summands of $C_n$,
	there must exist bases 
	$(\gamma_1, \ldots, \gamma_f)$ of 
	\bluetagged{$K_n^\text{free}$} 
	and 
	$(\omega_1, \ldots, \omega_q)$ of 
		\redtagged{$C_n \bigmod \ker(\boundary_n)$}
	that is a subset of some basis of $C_n$.
Observe that $B$, as given below, is a basis of $C_n$.
\begin{equation*}
	B := \bigg\{\,
		\underbrace{
			\greenmath{\alpha_1}, \ldots, 
			\greenmath{\alpha_r}
		}_{\displaystyle\text{basis of } K_n^\text{tor}}
		\,,\, 
		\overbrace{
			\bluemath{\gamma_1}, \ldots, 
			\bluemath{\gamma_f}
		}^{\displaystyle\text{basis of } K_n^\text{free}}
		\,,\, 
		\underbrace{ 
			\redmath{\omega_1}, \ldots, \redmath{\omega_q}
		}_{\mathclap{\displaystyle\text{basis of } C_n \bigmod \ker\boundary_n}}
	\,\bigg\}
\end{equation*}
A matrix $U_{n+1} \in \GL(m,R)$ such that 
	$(U_{n+1}, D_{n+1}, V_{n+1})$ is an SND of $[\boundary_{n+1}]$ and is as given by \fref{thm:chain-complex-decomposition}
can be defined as follows:
\begin{equation*}
	U_{n+1} = \hspace{12pt}
	\begin{NiceArray}{ ccc c ccc c ccc } 
		\CodeBefore [create-cell-nodes]
			\tikz \node [green-cell = (1-1) (2-1) (4-1), inner sep=2pt] {} ;
			\tikz \node [green-highlight = (2-2), inner sep=2pt] {} ;
			\tikz \node [green-cell = (1-3) (2-3) (4-3), inner sep=2pt] {} ;
			\tikz \node [blue-cell = (1-5) (2-5) (4-5), inner sep=2pt] {} ;
			\tikz \node [blue-highlight = (2-6), inner sep=2pt] {} ;
			\tikz \node [blue-cell = (1-7) (2-7) (4-7), inner sep=2pt] {} ;
			\tikz \node [red-cell = (1-9) (2-9) (4-9), inner sep=2pt] {} ;
			\tikz \node [red-highlight = (2-10), inner sep=2pt] {} ;
			\tikz \node [red-cell = (1-11) (2-11) (4-11), inner sep=2pt] {} ;
		\Body
			\vdots & {} & \vdots &
				& \vdots & {} & \vdots & 
				& \vdots & {} & \vdots \\
			[\alpha_1] & \cdots & [\alpha_r] &
				& [\gamma_1] & \cdots & [\gamma_f] & 
				& [\omega_1] & \cdots & [\omega_q] \\
			\vdots & {} & \vdots &
				& \vdots & {} & \vdots & 
				& \vdots & {} & \vdots \\[-7pt]
			{\phan} & {\phan} & {\phan} &
				& {\phan} & {\phan} & {\phan} & 
				& {\phan} & {\phan} & {\phan} \\
		\CodeAfter
			\SubMatrix({1-1}{3-11})[left-xshift=3pt, right-xshift=3pt]
			\UnderBrace{4-1}{4-3}{\tikz \node [green-node-padded, yshift=-10pt] {
				$K_n^\text{tor}$
			};}[yshift=5pt, color=ForestGreen, shorten]
			\UnderBrace{4-5}{4-7}{\tikz \node [blue-node-padded, yshift=-10pt] {
				$K_n^\text{free}$
			};}[yshift=5pt, color=NavyBlue, shorten]
			\UnderBrace{4-9}{4-11}{\tikz \node [red-node-padded, yshift=-10pt] {
				$C_n \bigmod \ker\boundary_n$
			};}[yshift=5pt, color=red, shorten]
		\end{NiceArray}
	\vspace{25pt}
\end{equation*} 
Observe that since $\im(\boundary_{n+1}) \subseteq K_n^\text{tor}$,
	there should be no non-trivial elements of $C_{n+1}$ mapping into $K_n^\text{free}$ and $C_n \bigmod \ker\boundary_n$.
Then, the expression of $[\boundary_{n+1}]$ relative to the basis $B$ of $C_n$ (with the same initial basis of $C_{n+1}$)
	should be given exactly by $D_{n+1}$.
That is, replacing the basis elements $\alpha_{r+1}, \ldots, \alpha_{m}$ of $K_n^\text{free} \oplus (C_n \bigmod \ker\boundary_n)$ with any basis element leaves the matrix $D_{n+1}$ undisturbed.
We can confirm this by 
	showing that the equation 
	$(W_{n+1})\inv [\boundary_{n+1}] V_{n+1} = D_{n+1}$
from the SND $(W_{n+1}, D_{n+1}, V_{n+1})$
implies that 
$(U_{n+1})\inv [\boundary_{n+1}] V_{n+1} = D_{n+1}$.

\HIDE{ 
	We can confirm this by the following calculation:
	Let $Z$ represent the initial basis on $C_n$,
		i.e.\ we write $[\alpha]$ for $\alpha \in C_n$
		as shorthand for $[\alpha] := [\alpha]_Z$ and 
		$[\boundary_{n+1}]$ for 
		$[\boundary_{n+1}] := [\boundary_{n+1}]_{X,Y}$ 
			with $Y$ some initial basis on $C_{n+1}$.
	Then, 
		$W_{n+1} = [\id_{C_n}]_{Z,A}$, 
		$U_{n+1} = [\id_{C_n}]_{Z,B}$,
		and 
		$U_{n+1} = W_{n+1}[\id_{C_n}]_{A,B}$
		with $[\id_{C_n}]_{A,B}$ being the change of basis matrix on $C_n$ from the basis $B$ to $A$.
	Since 
	$(W_{n+1})\inv [\boundary_{n+1}] V_{n+1} 
	= 
	D_{n+1}$ by assumption of $(W_{n+1}, D_{n+1}, V_{n+1})$ being an SND,
	the following matrix equation is true:
	\begin{align*}
		\tag{S1}
		([\id_{C_n}]_{A,B})\inv
		(W_{n+1})\inv [\boundary_{n+1}] V_{n+1} 
		&= 
		([\id_{C_n}]_{A,B})\inv 
		D_{n+1}
	\end{align*}
	On the left-hand side of Equation \textbf{(S1)}, we have that:
	\begin{equation*}
		([\id_{C_n}]_{A,B})\inv (W_{n+1})\inv [\boundary_{n+1}] V_{n+1}
		= 
		\Big( W_{n+1}[\id_{C_n}]_{A,B} \Big)
		\raisebox{6pt}{$\scriptstyle -1$}\,
		[\boundary_{n+1}] V_{n+1}
		= 
		\big( U_{n+1} \big)\inv[\boundary_{n+1}] V_{n+1}
	\end{equation*}
	For the right-hand side of Equation \textbf{(S1)}:
	Observe that $[\id_{C_n}]_{A,B}$ is a block matrix of the following form:
	\begin{equation*}
		[\id_{C_n}]_{A,B} = \begin{pmatrix}
			I_r & 0 \\
			0 & P_{A,B}
		\end{pmatrix}
		\qquad\text{ with }\qquad 
		([\id_{C_n}]_{A,B})\inv 
		= \begin{pmatrix}
			I_r & 0 \\
			0 & (P_{A,B})\inv
		\end{pmatrix}
	\end{equation*}
	with $I_r \in \GL(r,R)$ the identity matrix 
	and $P_{A,B} \in \GL(f+q,R)$ the change of basis matrix on 
		$K_n^\text{free} \oplus (C_n \bigmod \ker\boundary_n)$ 
		from the basis 
		$(\gamma_1, \ldots, \gamma_f) \cup (\omega_1, \ldots, \omega_q)$
		to the basis $(\alpha_{r+1}, \ldots, \alpha_m)$.
	Then,
	\begin{equation*}\arraycolsep=0.5\arraycolsep
		([\id_{C_n}]_{A,B})\inv D_{n+1} 
		= 
		\begin{pmatrix}
			I_r & 0 \\
			0 & (P_{A,B})\inv
		\end{pmatrix}
		\begin{pmatrix}
			\diag(d_1, \ldots, d_r) & 0 \\
			0 & 0
		\end{pmatrix}
		= \begin{pmatrix}
			I_r \diag(d_1, \ldots, d_r) & 0 \\
			0 & 0 
		\end{pmatrix}
		= D_{n+1}
	\end{equation*}
	Therefore, $(U_{n+1})\inv [\boundary_{n+1}] V_{n+1} = D_{n+1}$.
	Since $D_{n+1}$ is in Smith Normal Form by assumption,
		the triple $(U_{n+1}, D_{n+1}, V_{n+1})$ is an SND of $[\boundary_{n+1}]$.
	By construction of $U_{n+1}$,
		the SND $(U_{n+1}, D_{n+1}, V_{n+1})$ is the SND of $[\boundary_{n+1}]$ guaranteed by \fref{thm:chain-complex-decomposition}.
}

\spacer 

Observe that if cycle representatives for the homology groups are not required,
	we can calculate the invariant factor decomposition of homology groups of free chain complexes using any SNDs of $[\boundary_{n+1}]$ (for the nonzero invariant factors) and of $[\boundary_n]$ (for $\rank(\boundary_n)$).
We state this as a result below.

\begin{corollary}\label{cor:chain-complex-snf-kernel}
	Let $C_\ast = (C_n, \boundary_n)$ be a chain complex of finitely-generated free modules $C_n$ over a PID $R$ with differentials $\boundary_n: C_n \to C_{n-1}$.
	For each $n \in \nonnegints$, 
	the $n$\th homology group $H_n(C_\ast)$ can be calculated using the SNF $D_{n+1}$ of $[\boundary_{n+1}]$ as follows:
	\begin{equation*}
		H_n(C_\ast) \cong 
		R^f \oplus 
		\frac{R}{R\ket{d_1}}
		\oplus 
		\cdots 
		\oplus 
		\frac{R}{R\ket{d_r}}
	\end{equation*}
	where 
	$\set{d_1, \ldots, d_r}$ with $r = \rank(\boundary_{n+1})$
	are exactly the nonzero diagonal elements of $D_{n+1}$
	and $f + r = \rank(\ker\boundary_n)$.
	Note that $\rank(\ker\boundary_n)$ can be determined from the SNF of $[\boundary_n]$.
\end{corollary}
\begin{proof}
	Let $(U_{n+1}, D_{n+1}, V_{n+1})$ be the SND of $[\boundary_{n+1}]$ as denoted in \fref{thm:chain-complex-decomposition},
	i.e.\ the basis $B = (\beta_1, \ldots, \beta_m)$ of $C_n$ induced by $[\beta_j] = \col_j(U_{n+1})$ with $U_{n+1} \in \GL(m,R)$ partitions into three bases:
		$(\beta_1, \ldots, \beta_r)$ of $K_n^\text{tor}$,
		$(\beta_{r+1}, \ldots, \beta_{r+f})$ of $K_n^\text{free}$,
		and 
		$(\beta_{r+f+1}, \ldots, \beta_m)$ of $C_n \bigmod \ker\boundary_n$.
	The nonzero entries of $D_{n+1}$ are given by $d_1, \ldots, d_r$ with $r = \rank(\boundary_{n+1})$ and $d_j = D_{n+1}(j,j)$.
	Since $\ker(\boundary_n) = K_n^\text{tor} \oplus K_n^\text{free}$, $\rank(\boundary_n) = r+f$.
	Then,
	\begin{align*}
		H_n(C_\ast) = \frac{\ker\boundary_n}{\im\boundary_{n+1}}
		= 
		\frac{K_n^\text{tor}}{\im(\boundary_{n+1})}
		\oplus K_n^\text{free} 
		&=
		\frac{
			R\ket{\beta_1, \ldots, \beta_r}
		}{
			R\ket{d_1\beta_1, \ldots, d_r\beta_r}
		}
		\oplus 
		R\ket{\beta_{r+1}, \ldots, \beta_{r+f}}
		\\ 
		&= 
		\frac{R}{R\ket{d_1}}
		\oplus 
		\frac{R}{R\ket{d_1}}
		\oplus 
		R^{f}
	\end{align*}
	Since the SNF $D_{n+1}$ of $[\boundary_{n+1}]$ is unique up to multiplication by elementary dilations over $R$,
		any SND of $[\boundary_{n+1}]$ will share the same SNF.
\end{proof}

This result is particularly powerful since the SNDs of $[\boundary_{n+1}]$ as guaranteed by \fref{thm:chain-complex-decomposition} can be cumbersome to find.
In addition, most software packages only provide the SNF of matrices over $\ints$ as output, not the full SND.

We re-calculate $H_1(K)$ of \fref{ex:matrix-firsthom-sadness} again using this corollary below.
Note that the SNDs we provide in the example below are the SNDs returned by the \href{https://reference.wolfram.com/language/ref/SmithDecomposition.html}{\texttt{SmithDecomposition[-]}} 
method with the matrices $[\boundary_1]$ and $[\boundary_2]$ as input.
That is, the SNDs presented in \fref{ex:matrix-firsthom-sadness} are specifically calculated to match \fref{thm:chain-complex-decomposition}.

\begin{example}
	Let the simplicial complex $K$ be as given in \fref{ex:matrix-firsthom-sadness}.
	An SND $(U_2\sprime, D_2\sprime, V_2\sprime)$ of $[\boundary_2]$ is given below
		with the columns of $U_2^\prime$ corresponding to \greentagged{$K_1^\text{tor}$} highlighted in \greentag and the columns corresponding to the direct sum \purpletagged{$K_1^\text{free} \oplus (C_1 \bigmod \ker\boundary_1)$} are highlighted in \purpletag.
	\vspace{30pt}
	\begin{equation*}
		U_2 = \hspace{3pt}
		\begin{NiceArray}{>{\color{gray}}c<{\,\,} ccc ccc}
			\CodeBefore [create-cell-nodes]
				\tikz \node [green-cell = (1-2) (3-2) (6-2)] {} ;
				\tikz \node [green-cell = (1-3) (5-3) (6-3)] {} ;
				\tikz \node [purple-cell = (1-4) (6-4)] {} ;
				\tikz \node [purple-cell = (1-5) (6-5)] {} ;
				\tikz \node [purple-cell = (1-6) (6-6)] {} ;
				\tikz \node [purple-cell = (1-7) (6-7)] {} ;
			\Body
			{a_1 a_2} 	& 0 & 0  & 0 & 1 & 0 & 0		\\
			{a_1 a_3} 	& 1 & 0  & 0 & 0 & 0 & 0		\\
			{a_1 a_4} 	& -1 & 0 & 1 & 0 & 0 & 0		\\
			{a_2 a_3} 	& 0 & 1  & 0 & 0 & 0 & 0		\\
			{a_2 a_4} 	& 0 & -1 & 0 & 0 & 1 & 0		\\
			{a_3 a_4} 	& 1 & 1  & 0 & 0 & 0 & 1		
			\CodeAfter
			\SubMatrix({1-2}{6-7})[left-xshift=3pt, right-xshift=3pt]
			\OverBrace{1-2}{1-3}{\tikz \node [green-node-padded, yshift=8pt] {
					$K_1^\text{tor}$
				};}[yshift=7pt, shorten]
			\UnderBrace{6-4}{6-7}{\tikz \node [purple-node-padded, yshift=-10pt] {
					$K_1^\text{free} \oplus (C_1 \bigmod \ker\boundary_1)$
				};}[yshift=5pt, shorten]
		\end{NiceArray}
		\qquad\quad 
		D_2^\prime = \hspace{12pt}\begin{NiceArray}{cc}
			\CodeBefore [create-cell-nodes]
				\tikz \node [green-cell = (1-1)] {} ;
				\tikz \node [green-cell = (2-2)] {} ;
			\Body
			1 & 0 \\
			0 & 1 \\
			0 & 0 \\
			0 & 0 \\
			0 & 0 \\
			0 & 0
			\CodeAfter
			\SubMatrix({1-1}{6-2})[left-xshift=3pt, right-xshift=3pt]
		\end{NiceArray}
		\qquad\quad 
		V_2^\prime = I_2
		\vspace{30pt}
	\end{equation*}
	An SND $(U_1^\prime, D_1^\prime, V_1^\prime)$ of $[\boundary_1]$ is also given below,
		with the columns of $V_1^\prime$ corresponding to \redtagged{$C_1 \bigmod \ker\boundary_1$} highlighted in \redtag and the columns corresponding to \orangetagged{$\ker\boundary_1$} in \orangetag.
	\vspace{30pt}
	\begin{equation*}
		U_1^\prime = \begin{NiceArray}{>{\color{gray}}c (cccc)}
			{a_1} 	& -1 & -1 & -1 & 0	\\
			{a_2} 	& 1 & 0 & 0 & 0		\\
			{a_3} 	& 0 & 1 & 0 & 0		\\
			{a_4} 	& 0 & 0 & 1 & 1
		\end{NiceArray}
		\qquad 
		D_1^\prime = \hspace{9pt}\begin{NiceArray}{ccc ccc}
		\CodeBefore [create-cell-nodes]
			\tikz \node [red-cell = (1-1), inner xsep=2pt] {} ;
			\tikz \node [red-cell = (2-2), inner xsep=2pt] {} ;
			\tikz \node [red-cell = (3-3), inner xsep=2pt] {} ;
		\Body
			1 & 0 & 0 & 0 & 0 & 0 \\
			0 & 1 & 0 & 0 & 0 & 0 \\
			0 & 0 & 1 & 0 & 0 & 0 \\
			0 & 0 & 0 & 0 & 0 & 0 
		\CodeAfter
			\SubMatrix({1-1}{4-6})[left-xshift=3pt, right-xshift=3pt]
		\end{NiceArray}
		\hspace{9pt}
		\qquad 
		V_1^\prime = \hspace{3pt}\begin{NiceArray}{>{\color{gray}}c !{\,\,}
			ccc ccc 
		}
		\CodeBefore [create-cell-nodes]
			\tikz \node [red-cell = (1-2) (6-2)] {} ;
			\tikz \node [red-cell = (1-3) (6-3)] {} ;
			\tikz \node [red-cell = (1-4) (6-4)] {} ;
			\tikz \node [orange-cell = (1-5) (2-5) (6-5)] {} ;
			\tikz \node [orange-cell = (1-6) (3-6) (6-6)] {} ;
			\tikz \node [orange-cell = (1-7) (3-7) (6-7)] {} ;
		\Body
			{a_1 a_2} 	& 1 & 0 & 0 & 1 & 1 & 0		\\
			{a_1 a_3} 	& 0 & 1 & 0 & -1 & 0 & 1		\\
			{a_1 a_4} 	& 0 & 0 & 1 & 0 & -1 & -1		\\
			{a_2 a_3} 	& 0 & 0 & 0 & 1 & 0 & 0		\\
			{a_2 a_4} 	& 0 & 0 & 0 & 0 & 1 & 0		\\
			{a_3 a_4} 	& 0 & 0 & 0 & 0 & 0 & 1 \\
		\CodeAfter
			\SubMatrix({1-2}{6-7})[left-xshift=3pt, right-xshift=3pt]
			\UnderBrace{6-2}{6-4}{\tikz \node [red-node-padded, yshift=-10pt] {
					$C_1 \bigmod \ker\boundary_1$
				};}[yshift=5pt, shorten]
			\OverBrace{1-5}{1-7}{\tikz \node [orange-node-padded, yshift=8pt] {
					$K_1^\text{tor} \oplus K_1^\text{free}$
				};}[yshift=7pt, shorten]
		\end{NiceArray}
		\vspace{30pt}
	\end{equation*} 
	We can determine the following facts from these two SNDs:
	\begin{enumerate}
		\item From $D_1^\prime$, $\rank(\ker\boundary_1) = 6-3 = 3$ and from $D_2^\prime$, $r = \rank(\boundary_2) = 2$.
		Therefore, $f = 1$. 

		\item 
		Let $A = (\alpha_1, \ldots, \alpha_6)$ be a basis of $C_1$ by $[\alpha_j] = \col_j(U_2^\prime)$.
		We can confirm that the basis $(\alpha_3, \ldots, \alpha_6)$ of $K_1^\text{free} \oplus (C_1 \bigmod \ker\boundary_1)$ does not partition into bases for each summand since $\boundary_1(\alpha_j) \neq 0$ for $j \in \set{3, \ldots, 6}$
		but $\rank(K_1^\text{free}) = f = 1$.

		\item 
		Let $K = (\kappa_1, \ldots, \kappa_6)$ be a basis on $C_1$ by $[\kappa_i] = \col_i(V_1^\prime)$.
		Since the basis elements $\kappa_4, \kappa_5, \kappa_6$ do not appear as basis elements in $A$, 
			i.e.\ the columns of $V_1^\prime$ are not columns of $U_2^\prime$.
		Note that it is not immediately clear if one of $\kappa_4, \kappa_5, \kappa_6$ generates $K_1^\text{free}$.
	\end{enumerate}
	However, if we do not need to calculate cycle representatives, we have enough information to calculate $H_1(K)$ up to isomorphism using \fref{cor:chain-complex-snf-kernel}:
	From $D_2^\prime$, we have $r = 2$, $d_1 = 1$, $d_2 = 1$.
	From $D_1^\prime$, we have $\rank(\ker\boundary_1) = 3$ and $f = 3-2 = 1$. Then,
	\begin{equation*}
		H_1(K) 
		\cong 
		\ints^{f} \oplus 
			\paren{\frac{\ints}{d_1\ints}} \oplus 
			\paren{\frac{\ints}{d_2\ints}}
		=
		\ints^{3-2} \oplus 
			\paren{\frac{\ints}{\ints}} \oplus 
			\paren{\frac{\ints}{\ints}}
		= \ints
		\,.
	\end{equation*}
	Observe that this matches with our calculation for $H_1(K)$ in \fref{ex:matrix-firsthom-sadness}.
\end{example}

\clearpage

%

\section{The Graded Structure Theorem and SNDs in the Graded Case} 
\label{section:calculation-graded-ifds} 

In \fref{section:matrix-calculation-of-IFDs}, we discussed a method of calculating the invariant factor decompositions of finitely generated modules over a PID $R$, i.e.\ finitely generated modules in the category $\catmod{\field[x]}$.
In this section, we extend this method to the category $\catgradedmod{\field}$
of \textit{graded} $\field[x]$-modules, wherein invariant factor decompositions need to be given by graded isomorphisms.
For reference, we discuss the category $\catgradedmod{\field}$ and the notation we use involving graded $\field[x]$-modules earlier in \fref{section:graded-mod-notation}.

To start, we state \textbf{Graded Structure Theorem} in $\catgradedmod{\field}$ for graded $\field[x]$-modules, which is the theorem corresponding to the Structure Theorem (\fref{thm:structure-theorem}) in the ungraded category $\catmod{\field[x]}$ of $\field[x]$-modules.

\begin{theorem}\textbf{The Graded Structure Theorem.}%
\label{thm:graded-structure-theorem}

	Let $M$ be a finitely generated graded $\field[x]$-module over $\field[x]$ for some field $\field$.
	Then, 
	there exists a finite direct sum of shifted cyclic graded submodules of $\field[x]$ that is graded isomorphic to $M$ as follows:
	\begin{equation*}
		M \,\upgraded\cong\,
			\mathlarger\Sigma^{s_1}\!\paren{\frac{\field[x]}{(x^{t_1})}}
			\oplus \cdots \oplus 
			\mathlarger\Sigma^{s_r}\!\paren{\frac{\field[x]}{(x^{t_r})}}
			\oplus 
			\Sigma^{s_{r+1}}\field[x]
			\oplus \cdots \oplus 
			\Sigma^{s_{m}}\field[x]
	\end{equation*}
	with indices $s_1, \ldots, s_r, \ldots, s_m \in \nonnegints$ and 
	non-zero, non-unit $x^{t_1}, x^{t_2}, \ldots, x^{t_r} \in \field[x]$
	such that the divisibility relation 
	$x^{t_1} \divides x^{t_2} \divides \cdots \divides x^{t_r}$ is satisfied.
	Furthermore, the collection 
	$
		\bigl\{
			(x^{t_1}, s_1), \ldots, (x^{t_r}, s_r),
			(0, s_{r+1}), \ldots, (0, s_m)
		\bigr\}
	$
	is uniquely determined by $M$ up to graded isomorphism.
\end{theorem}
\remark{
	This theorem is stated in \cite[Theorem 2.1]{matrixalg:zomorodian} without proof.
	The paper 
	\cite{persmod:clara} by Clara L\"{o}h provides a detailed proof, with this theorem listed as \cite[Theorem 5.1]{persmod:clara}.
}

Much like the case of the Structure Theorem in $\catmod{\field[x]}$, 
we give the decomposition guaranteed by the Graded Structure Theorem a special name.

\begin{definition}\label{defn:graded-invariant-factor-decomposition}
	Let $M$ be a finitely generated $\field[x]$-module and let the following direct sum decomposition of $M$ be as denoted in the Graded Structure Theorem 
	(\fref{thm:graded-structure-theorem}):
	\begin{equation*}
		M \,\upgraded\cong\,
			\mathlarger\Sigma^{s_1}\!\paren{\frac{\field[x]}{(x^{t_1})}}
			\oplus \cdots \oplus 
			\mathlarger\Sigma^{s_r}\!\paren{\frac{\field[x]}{(x^{t_r})}}
			\oplus 
			\Sigma^{s_{r+1}}\field[x]
			\oplus \cdots \oplus 
			\Sigma^{s_{m}}\field[x]
	\end{equation*}
	This decomposition is called the \textbf{graded invariant factor decomposition} of $M$.
	The \textbf{invariant factors} of $M$ are given by 
		$x^{t_1}, \ldots, x^{t_r} \in \field[x]$ (which are non-zero and non-unit)
	and the \textbf{grading shifts} of $M$ by $s_1, \ldots, s_m \in \nonnegints$.
\end{definition}
\remark{
	The term \textit{invariant factor} for graded $\field[x]$-modules is deliberately chosen since it corresponds to that of $\field[x]$-modules, 
	i.e.\ the invariant factors of a graded $\field[x]$-module $M$ are exactly the invariant factors of $M$ viewed as a $\field[x]$-module (disregarding grading).
	Note that some references, e.g.\ \cite{persmod:clara},
	define invariant factors to include zero elements.
	To avoid confusion, we sometimes use the term \textit{nonzero invariant factor} for clarity. 
}

In this section, we discuss how this graded invariant factor decomposition can be calculated, under certain restrictions, using a method similar to that for non-graded decompositions as presented in \fref{section:matrix-calculation-of-IFDs}.
Let $M$ be a finitely generated graded $\field[x]$-module.
Since $\field[x]$ is a PID,
	the non-graded Structure Theorem (\fref{thm:structure-theorem}) guarantees the existence of an invariant factor decomposition for $M$ as a $\field[x]$-module as follows:
	\begin{equation*}
		M \cong 
			\paren{\frac{\field[x]}{(f_1)}}
			\oplus \cdots \oplus 
			\paren{\frac{\field[x]}{(f_r)}}
			\oplus 
			\field[x]
			\oplus \cdots \oplus  
			\field[x]
	\end{equation*}
	with nonzero invariant factors $f_1, \ldots, f_r \in \field[x]$.
However, this isomorphism may not correspond to that between graded $\field[x]$-modules.
For example, if $f_j \in \field[x]$ is not a homogeneous element for some $j \in \set{1, \ldots, r}$, 
	then the ideal $(f_j)$ is not a graded ideal of $\field[x]$
	and the quotient $\field[x] \bigmod (f_j)$ is not necessarily a graded ring.
If we assume that each $f_j \in \field[x]$ is homogeneous of the form $f_j = x^{t_j}$ and each homogeneous element $q \in M$ is mapped to a homogeneous element in one of the summands, 
then we can shift each summand such that the degree of $q \in M$ matches that of its image in the decomposition as follows:
\begin{equation}\label{eqn:graded-starting-eqn}
	M \upgraded\cong 
		\Sigma^{\degh(q_1)}\!\paren{\frac{\field[x]}{(x^{t_1})}}
		\oplus \cdots \oplus 
		\Sigma^{\degh(q_r)}\!\paren{\frac{\field[x]}{(x^{t_r})}}
		\oplus 
		\Sigma^{\degh(q_{r+1})}\field[x]
		\oplus \cdots \oplus 
		\Sigma^{\degh(q_{r+1})}\field[x]
\end{equation}
where $q_j \in M$ is chosen such that $q_j$ corresponds to the generator $1_j$ of each summand. 
We claim that, assuming we start with a presentation that respects the graded structure, 
then a non-graded invariant factor decomposition can be transformed (for lack of a better word) to a graded invariant factor decomposition.
Below, we provide a definition for these kinds of presentations.

\begin{definition}\label{defn:graded-presentation}
	A \textbf{graded presentation} of a graded $\field[x]$-module $M$ 
	is a presentation
	\begin{equation*}
		F_S \,\,\Xrightarrow{\,\phi\,}\,\, 
		F_G \,\,\Xrightarrow{\,\pi\,}\,\,
		M \,\,\Xrightarrow{\,\,\,}\,\, 0
	\end{equation*}
	on $M$
	such that $F_S$ and $F_G$ are free graded $\field[x]$-modules and the homomorphisms $\phi: F_S \to F_G$ and $\pi: F_G \to M$ are graded homomorphisms.
\end{definition}


A graded presentation by $\phi: F_S \to F_G$ and $\pi: F_G \to M$ of a graded $\field[x]$-module $M$ determines $M$ up to non-graded isomorphism by $M \cong F_G \bigmod \im(\phi)$, as discussed under \fref{defn:presentation}.
Since $\im(\phi)$ is a graded submodule of $F_G$, the quotient $F_G \bigmod \im(\phi)$ is a graded $\field[x]$-module that inherits the grading of $F_G$ and 
	the graded presentation also determines $M$ up to graded isomorphism.
As with the non-graded case, the calculation of graded invariant factor decompositions start with graded presentations. Below, we provide an existence claim for these graded presentations.

\begin{lemma}
	There exists a graded presentation for any finitely generated graded $\field[x]$-module.
\end{lemma}
\begin{proof}
	Let $M$ be a finitely generated graded $\field[x]$-module.
	By assumption of $M$ being graded, 
		there must exist a homogeneous system of generators $\set{a_1 x^{s_1}, \ldots, a_m x^{s_m}}$ of $M$ with $\degh(a_j x^{s_j}) = s_j \in \nonnegints$ for all $j \in \set{1, \ldots, m}$.
	Note that we write $a_j x^{s_j} \in M$ such that $a_j$ is an element of the $\field$-vector space $M_{s_j}$ with $M_{s_j} x^{s_j}$ being the homogeneous component of $M$ of degree $s_j$.
	We can construct a graded presentation for $M$ as follows:

	Define the module of generators $F_G$ to be the free graded $\field[x]$-module with homogeneous basis given by 
	$A = (\alpha_1 x^{s_1}, \ldots, \alpha_m x^{s_m})$ with $\degh(\alpha_j x^{s_j}) = s_j$
	and define $\pi: F_G \to M$ by $\alpha_j x^{s_j} \mapsto a_j x^{s_j}$ for $j \in \set{1, \ldots, m}$.
	Then, $\pi$ is a graded homomorphism.
	Observe that if the set of generators of $M$ were not homogeneous, 
	then the requirement that $\degh(\alpha_j x^{s_j}) = \degh(a_j x^{s_j})$ may not be fulfilled since $\degh(a_j x^{s_j})$ may be undefined.
	For reference, $F_G$ can be expressed as a direct sum of shifted graded free $\field[x]$-modules, each with one basis element, as follows:
	\begin{equation}
		F_G 
		= 
		\field[x]\ket{\alpha_1 x^{s_1}} 
			\oplus \cdots \oplus 
			\field[x]\ket{\alpha_m x^{s_m}}
		= 
		\Sigma^{s_1}\fieldket{\alpha_1}
			\oplus \cdots \oplus 
			\Sigma^{s_m}\fieldket{\alpha_m}
	\end{equation}
	Define the module of relations $F_S$ by $F_S := \ker(\pi)$ and 
	define $\phi: F_S \to F_G$ to be the inclusion of $\ker(\pi)$ into $F_G$.
		Note that $\ker(\pi)$ is a graded submodule of $F_S$ as it is the kernel of the graded homomorphism $\pi: F_G \to M$, and $\phi$ is graded as the identity map. 
	As with the non-graded case, $\phi: F_S \to F_G$ and $\pi: F_G \to M$ correspond to a non-graded presentation of $M$.
	Since $F_S$, $F_G$ are graded $\field[x]$-modules and $\phi$, $\pi$ are graded homomorphisms,
		this presentation of $M$ is also a graded presentation of $M$.
\end{proof}

Let $M$ be a finitely generated graded $\field[x]$-module.
Let a graded presentation of $M$ be given by $\phi: F_S \to F_G$ and $\pi: F_G \to M$ with bases 
$S = (\sigma_1(x), \ldots, \sigma_n(x))$ of $F_S$ and 
$A = (\alpha_1(x), \ldots, \alpha_m(x))$ of $F_G$.
In the non-graded case, an SND $(U,D,V)$ of $[\phi]_{A,S}$, which exists since $\field[x]$ is a PID, gives us the following by \fref{prop:snd-on-homomorphisms}:
\begin{enumerate}
	\item A basis $T = (\tau_1(x), \ldots, \tau_n(x))$ of $F_G$ by $[\tau_i(x)]_S = \col_i(V)$.
	\item A basis $B = (\beta_1(x), \ldots, \beta_m(x))$ of $F_S$ by $[\beta_j(x)]_A = \col_j(U)$.
	\item The set of nonzero invariant factors $f_1, \ldots, f_r \in \field[x]$ by $f_j = D(j,j)$.
\end{enumerate}
Note that we write $\sigma_i(x)$ and $\alpha_j(x)$ here to emphasize that these elements are polynomials in $x$, i.e.\ the degree of $\sigma_i(x)$ and of $\alpha_j(x)$ are generally not zero. 

For the process described for \fref{eqn:graded-starting-eqn} to make sense, i.e.\ shifting the summands of the (non-graded) invariant factor decomposition,
	the basis $B = (\beta_j(x))$ must be a homogeneous basis.
Otherwise, $\degh(\beta_j(x))$ and $\Sigma^{\degh(\beta_j(x))}\field[x]$ would be undefined.
If we have that $B$ is homogeneous and that the nonzero invariant factors $\set{f_j}$ are also homogeneous, then we can calculate a graded decomposition for $M$ that would later correspond to a graded invariant factor decomposition.
We state this in more detail below.

\begin{proposition}\label{prop:how-to-get-graded-ifd}\label{prop:graded-SND-exists}
	Let $M$ be a finitely generated graded $\field[x]$-module. 
	Let $\phi: F_S \to F_G$, $\pi: F_G \to M$ correspond to a finite graded presentation on $M$ with $\rank(F_S) = n$ and $\rank(F_G) = m$.

	If there exists a homogeneous basis $\set{\beta_1 x^{s_1}, \ldots, \beta_m x^{s_m}}$ of $F_G$ 
	and nonzero elements $x^{t_1}, \ldots, x^{t_r} \in \field[x]$ with divisibility relation $x^{t_1} \divides \cdots \divides x^{t_r}$ such that 
	$\set{\beta_1 x^{s_1 + t_1}, \ldots, \beta_r x^{s_r + t_r}}$ is a basis for $\im(\phi)$, then we have the following graded isomorphism on $M$:
	\begin{align}
		M 
		&\upgraded\cong 
			\paren{\frac{
				\field[x]\ket{\beta_1 x^{s_1}}
			}{
				\field[x]\ket{\beta_1 x^{s_1 + t_1}}
			}}
			\oplus \cdots \oplus 
			\paren{\frac{
				\field[x]\ket{\beta_r x^{s_r}}
			}{
				\field[x]\ket{\beta_r x^{s_r + t_r}}
			}}
			\oplus 
			\field[x]\ket{\beta_{r+1}x^{s_{r+1}}}
			\oplus \cdots \oplus 
			\field[x]\ket{\beta_m x^{s_m}}
		\label{eqn:how-to-get-graded-ifd}
	\end{align}
	Note that we write the basis elements $\beta_j x^{s_j} \in F_G$ such that $\degh(\beta_j x^{s_j}) = s_j$.
\end{proposition}
\begin{proof}
	Assume there does exist a homogeneous basis $\set{\beta_j x^{s_j}}_{j=1}^m$ of $F_G$ and nonzero elements $\set{x^{t_j}}_{j=1}^r \in \field[x]$.
	For $j \in \set{1, \ldots, r}$:
		$\field[x]\ket{\beta_j x^{s_j + t_j}}$ is a graded submodule of $\field[x]\ket{\beta_j x^{s_j}}$ since it has a homogeneous basis.
	Therefore,
		the quotient module 
		$\field[x]\ket{\beta_j x^{s_j}} \bigmod \field[x]\ket{\beta_j x^{s_j + t_j}}$ 
		is graded with homogeneous components of degree $q \in \set{s_j, \ldots, t_j - 1}$ given by $[k \cdot \beta_j x^q]$ for $k \in \field$ nonzero.

	By \fref{prop:how-to-get-ifd}, \ref{eqn:how-to-get-graded-ifd} is a $\field[x]$-module isomorphism.
	For $j \in \set{1, \ldots, r}$:
		the isomorphism by (1) maps $[\beta_j x^q]
			\in \field[x]\ket{\beta_j x^{s_j}} \bigmod \field[x]\ket{\beta_j x^{s_j + t_j}}
		$ with $q \in \set{s_j, \ldots, t_j - 1}$ 
		to $\pi(\beta_j x^q) \in M$.
	Since $\pi$ is a graded homomorphism, we have that 
	\begin{equation*}
		\degh([\beta_j x^q]) = q = \degh(\beta_j x^q) = \degh(\pi(\beta_j x^q))
	\end{equation*}
	Similarly, for $j \in \set{r+1, \ldots, m}$,
		the isomorphism by (1) maps $\beta_j x^q \in \field[x]\ket{\beta_j x^{s_j}}$ with $q \geq s_j$ to $\pi(\beta_j x^q)$
		and $\degh(\beta_j x^q) = q = \degh(\pi(\beta_j x^q))$.
	Therefore, \ref{eqn:how-to-get-graded-ifd} is a graded isomorphism.
\end{proof}

We provide an example of this calculation below.
Note that we have not yet proven that a homogeneous basis $\set{\beta_j x^{s_j}}$ and a set of nonzero invariant factors $\set{x^{t_j}}$ with the desired properties generally exists.

\begin{example}
	\label{ex:first-graded-snd-calc}
	Let $M$ be a graded $\rationals[x]$-module and define a graded presentation on $M$ as follows:
	Let $S = \set{\sigma_1 x, \sigma_2 x^6, \sigma_3 x^8, \sigma_4 x^{10}}$ be a basis on the $\rationals[x]$-module $F_S$ of relations
	and let $A = \set{\alpha_1 x, \alpha_2 x^2, \alpha_3 x^3, \alpha_4 x^4}$ be a basis on the $\rationals[x]$-module $F_G$ of generators.
	Observe that since $S$ and $A$ are homogeneous bases, $F_S$ and $F_G$ are both graded $\rationals[x]$-modules.
	Let $\phi: F_S \to F_G$ be given by the following matrix:
	\begin{equation*}
		[\phi]_{A,S} = \begin{pmatrix}
			1 & 4x^5 & 0 & 3x^9 \\
			0 & x^4 & -x^6 & 0 \\
			0 & 0 & x^5 & -2x^7 \\
			0 & 0 & 0 & 0
		\end{pmatrix}
	\end{equation*}
	We can confirm by inspection that $\phi$ is a graded homomorphism.
	Assume that $M = F_G \bigmod \im(\phi)$, i.e.\ $\pi$ is the canonical quotient homomorphism.
	An SND $(U_1, D_1, V_1)$ of $[\phi]_{A,S}$ is given below:
	\begin{equation*}
		U_1 = \begin{pmatrix}
			1 & 2x & -x^2 & 0 \\
			0 & 1 & 0 & -x^2 \\
			0 & 0 & 1 & 2x \\
			0 & 0 & 0 & 1
		\end{pmatrix}
		\qquad
		D_1 = \begin{pmatrix}
			1 & 0 & 0 & 0 \\
			0 & x^4 & 0 & 0 \\
			0 & 0 & x^5 & 0 \\
			0 & 0 & 0 & 0
		\end{pmatrix}
		\qquad 
		V_1 = \begin{pmatrix}
			1 & -2x^5 & -5x^7 & -11x^9 \\
			0 & 1 & x^2 & 2x^4 \\
			0 & 0 & 1 & 2x^2 \\
			0 & 0 & 0 & 1
		\end{pmatrix}
	\end{equation*}
	We can confirm this matrix factorization is correct by the following calculation:
	\begin{equation*}
		(U_1)\inv 
		[\phi]_{A,S}
		V_1
		=
		\begin{pmatrix}
			1 & -2x & x^2 & -4x^3 \\
			0 & 1 & 0 & x^2 \\
			0 & 0 & 1 & -2x \\
			0 & 0 & 0 & 1
		\end{pmatrix}
		\begin{pmatrix}
			1 & 4x^5 & 0 & 3x^9 \\
			0 & x^4 & -x^6 & 0 \\
			0 & 0 & x^5 & -2x^7 \\
			0 & 0 & 0 & 0
		\end{pmatrix}
		\begin{pmatrix}
			1 & -2x^5 & -5x^7 & -11x^9 \\
			0 & 1 & x^2 & 2x^4 \\
			0 & 0 & 1 & 2x^2 \\
			0 & 0 & 0 & 1
		\end{pmatrix}
		= \cdots 
		= D_1
	\end{equation*}
	The matrix $U_1 \in \GL(4,\rationals[x])$ induces a new basis $B = (\beta_i x^{t_i})$ on $F_G$ by $[\beta_i x^{t_i}]_A = \col_i(U)$, given below,
	along with its corresponding diagonal element $d_i = D_1(i,i)$.
	Note that the basis elements of $B$ can be written as $\beta_i x^{t_i}$ with $\degh(\beta_i x^{t_i}) = t_i$ since the basis $B$ consists of homogeneous elements (confirmed post-calculation).
	\begin{equation*}\arraycolsep=3pt
		\begin{array}{r c c c r !{\qquad\text{ with }\qquad} ccc }
			\beta_1 x^{t_1} 
				&=& (1)(\alpha_1) x
				&=& (\alpha_1) x\phantom{^1}
				& d_1 &=& 1
				\\
			\beta_2 x^{t_2} 
				&=& (2x)(\alpha_1 x) + (1)(\alpha_2 x^2)
				&=& (2\alpha_1 + \alpha_2)x^2
				& d_2 &=& x^4
				\\ 
			\beta_3 x^{t_3}
				&=& (-x^2)(\alpha_1 x) + (1)(\alpha_3 x^3)
				&=& (-\alpha_1 + \alpha_3)x^3
				& d_3 &=& x^5
				\\
			\beta_4 x^{t_4}
				&=& (-x^2)(\alpha_2 x^2) + (2x)(\alpha_3 x^3) + (1)(\alpha_4 x^4)
				&=& (-\alpha_2 + 2\alpha_3 + \alpha_4)x^4
				& d_4 &=& 0
		\end{array}
	\end{equation*}
	By applying \fref{prop:how-to-get-ifd}, we have the following (non-graded) isomorphism on $M$. Note that since $B$ is a homogeneous basis and the nonzero elements $d_i \in \field[x]$ are also homogeneous, the following decomposition is also a graded isomorphism by \fref{prop:how-to-get-graded-ifd}:
	\begin{equation*}
		M \upgraded\cong 
			\paren{\frac{
				\rationals[x]\ket{\beta_1 x}
			}{
				\rationals[x]\ket{\beta_1 x}
			}}
			\oplus 
			\paren{\frac{
				\rationals[x]\ket{\beta_2 x^2}
			}{
				\rationals[x]\ket{x^4 \cdot \beta_2 x^2}
			}}
			\oplus
			\paren{\frac{
				\rationals[x]\ket{\beta_3 x^3}
			}{
				\rationals[x]\ket{x^5 \cdot \beta_3 x^3}
			}}
			\oplus 
			\rationals[x]\ket{\beta_4 x^4}
	\end{equation*}
	Since the only units of $\rationals[x]$ are nonzero field elements, the nonzero entries of the SNF $D_1$ of $[\phi]_{A,S}$ have to be homogeneous elements by uniqueness of SNFs.
	However, we can do row operations on the SNF $(U_1)\inv [\phi]_{A,S} V_1 = D_1$ without disturbing the nonzero elements of $D$, i.e.\ $1$, $x^4$ and $x^5$.
	In particular, we can add any $\rationals[x]$-multiple of the $4$\th row of $D_1$ to rows $1$, $2$, and $3$ of $D_1$ without changing $D_1$. However, this will affect the matrix that determines the basis for $F_G$.
	For example, we can do the following:
	\begin{equation*}
		\eladd[4]{1,4 \,; x^{10}}D_1 
		= 
		\begin{pmatrix}
			1 & 0 & 0 & x^{10} \\
			0 & 1 & 0 & 0 \\
			0 & 0 & 1 & 0 \\
			0 & 0 & 0 & 1
		\end{pmatrix}
		\begin{pmatrix}
			1 & 0 & 0 & 0 \\
			0 & x^4 & 0 & 0 \\
			0 & 0 & x^5 & 0 \\
			0 & 0 & 0 & 0
		\end{pmatrix}
		=
		\begin{pmatrix}
			1 & 0 & 0 & 0 \\
			0 & x^4 & 0 & 0 \\
			0 & 0 & x^5 & 0 \\
			0 & 0 & 0 & 0
		\end{pmatrix}
		= D_1
	\end{equation*}
	Then, we can get another SND $(U_2, D_1, V_1)$ of $[\phi]_{A,S}$ by considering the following factorization:
	\begin{equation*}
		D_1
		= \eladd[4]{1,4 \,; x^{10}}D_1 
		= \eladd[4]{1,4 \,; x^{10}}(U_1)\inv [\phi]_{A,S} V_1
		= \underbrace{\Big(
				U_1\eladd[4]{1,4 \,; -x^{10}}
				\Big)\inv
			}_{\textstyle\mathclap{\text{Set this as } (U_2)\inv}}
		[\phi]_{A,S} V_1
	\end{equation*}
	Note that the inverse of the elementary transvection $\eladd[4]{1,4 \,; x^{10}}$ is given by 
	$\eladd[4]{1,4 \,; -x^{10}}$.
	Then, the matrix $U_2 \in \GL(4,\rationals[x])$ determines a new basis $P = (p_1(x), \ldots, p_4(x))$ of $F_G$ as follows:
	\begin{equation*}
		U_2 = \begin{pmatrix}
			1 & 2x & -x^2 & -x^{10} \\
			0 & 1 & 0 & -x^2 \\
			0 & 0 & 1 & 2x \\
			0 & 0 & 0 & 1
		\end{pmatrix}
		\qquad\text{ and }\qquad 
		\begin{aligned}
			p_1(x) &= \beta_1 x \\
			p_2(x) &= \beta_2 x^2 \\
			p_3(x) &= \beta_3 x^3 \\ 
			p_4(x) &= (-x^{10})(\alpha_1 x) + \beta_4 x^4 = -\alpha_1 x^{11} + \beta_4 x^4
		\end{aligned}
	\end{equation*}
	Then, the isomorphism by \fref{prop:how-to-get-ifd} on the SND $(U_2, D_1, V_1)$ of $[\phi]_{A,S}$, given below, is \textbf{not} a graded isomorphism:
	\begin{equation*}
		M
		\upmod\cong 
			\paren{\frac{
				\field[x]\ket{\beta_1 x}
			}{
				\field[x]\ket{\beta_1 x}
			}}
			\oplus 
			\paren{\frac{
				\field[x]\ket{\beta_2 x^2}
			}{
				\field[x]\ket{x^4 \cdot \beta_2 x^4}
			}}
			\oplus
			\paren{\frac{
				\field[x]\ket{\beta_3 x^3}
			}{
				\field[x]\ket{x^5 \cdot \beta_2 x^5}
			}}
			\oplus 
			\field[x]\ket{-\alpha_1 x^{11} + \beta_4 x^4}
	\end{equation*}
	In particular, the summand $\field[x]\ket{-\alpha_1 x^{11} + \beta_4 x^4}$ fails to be a graded $\rationals[x]$-module since it cannot be generated by a homogeneous element. 
\end{example}

The decomposition given by \fref{prop:how-to-get-graded-ifd} is then transformed into a graded invariant factor decomposition by applying specific graded isomorphisms to each summand in the direct sum.
We state these isomorphisms in the following lemma.

\begin{lemma}\label{lemma:summands-on-graded-ifds}
	Let $\field[x]\ket{\alpha x^s}$ be the free graded $\field[x]$-module generated by the element $\alpha x^s$ with $\deg(\alpha x^s) = s \in \nonnegints$. Let $x^t \in \field[x]$ for some $t \in \nonnegints$ with $t \neq 0$.
	Then, we have the following graded isomorphisms:
	\begin{equation*}
		\field[x]\ket{\alpha x^s} 
			\upgraded\cong 
			\Sigma^{s}\field[x]
		\qquad\text{ and }\qquad 
			\frac{
				\field[x]\ket{\alpha x^s}
			}{
				\field[x]\ket{x^t \cdot \alpha x^s}
			}
			\upgraded\cong 
			\mathlarger\Sigma^{s}
			\Big( \field[x] \bigmod (x^t)
			\Big)
	\end{equation*}
	Note that if $t=0$, then $x^t = 1$ is a unit of $\field[x]$ and $\field[x]\ket{\alpha x^s} \bigmod \field[x]\ket{x^t \cdot \alpha x^s} = 0$.
\end{lemma}
\begin{proof}
	Let $\phi: \field[x]\ket{\alpha x^s} \to \Sigma^s \field[x]$ 
		be given by $\alpha x^s \mapsto (x^s)(1)$.
	This is an isomorphism with inverse $(x^s)(1) \mapsto \alpha x^s$. 
	Let $f \in \field[x]\ket{\alpha x^s}$ be a homogeneous element.
	Then, $f = k \cdot \alpha x^{s+r}$ for some $k \in \field$ nonzero and $r \in \nonnegints$
	and 
	\begin{align*}
		\degh(f) &= \degh(k \cdot \alpha x^{s+r}) = s+r 
			= \degh\!\big( (x^{s+r})(1) \big) = \degh(\phi(f)) 
	\end{align*}
	The quotients $\field[x] \bigmod (x^t)$ and $\field[x]\ket{\alpha x^{s}} \bigmod \field[x]\ket{\alpha x^{s+t}}$
	are graded $\field[x]$-modules with grading determined by their coset representatives, well-defined since $(x^t)$ and $\field[x]\ket{\alpha x^{s+t}}$ are graded submodules respectively.
	Then, the isomorphism 
	$\field[x]\ket{\alpha x^{s}} \bigmod \field[x]\ket{\alpha x^{s+t}}
	\to \Sigma^s(\field[x] \bigmod (x^t))$ given by $[\alpha x^s] \to (x^s)[1]$ is graded similarly as in the case of $\phi$.
\end{proof}

By this lemma, the graded decomposition given by \fref{prop:how-to-get-graded-ifd} can be transformed into a graded invariant factor decomposition 
as follows, with summands that have a unit for $x^{t_j}$, i.e.\ $t_j = 0$, removed from the decomposition.
\begin{align*}
	M 
	&\upgraded\cong\, 
		\paren{\frac{
			\field[x]\ket{\beta_1 x^{s_1}}
		}{
			\field[x]\ket{\beta_1 x^{s_1 + t_1}}
		}}
		\oplus \cdots \oplus 
		\paren{\frac{
			\field[x]\ket{\beta_r x^{s_r}}
		}{
			\field[x]\ket{\beta_r x^{s_r + t_r}}
		}}
		\oplus 
		\field[x]\ket{\beta_{r+1}x^{s_{r+1}}}
		\oplus \cdots \oplus 
		\field[x]\ket{\beta_m x^{s_m}}
	\\ 
	&\upgraded\cong\, 
		\Sigma^{s_1}\!\paren{\frac{
			\field[x]\ket{\beta_1}
		}{
			\field[x]\ket{\beta_1 x^{t_1}}
		}}
		\oplus \cdots \oplus 
		\Sigma^{s_r}\!\paren{\frac{
			\field[x]\ket{\beta_r }
		}{
			\field[x]\ket{\beta_r x^{t_r}}
		}}
		\oplus 
		\Sigma^{s_{r+1}}\field[x]\ket{\beta_{r+1}}
		\oplus \cdots \oplus 
		\Sigma^{s_m}\field[x]\ket{\beta_m}
	\\ 
	&\upgraded\cong\, 
		\Sigma^{s_1}\!\paren{\frac{
			\field[x]
		}{
			(x^{t_1})
		}}
		\oplus \cdots \oplus 
		\Sigma^{s_r}\!\paren{\frac{
			\field[x]
		}{
			(x^{t_r})
		}}
		\oplus 
		\Sigma^{s_{r+1}}\field[x]
		\oplus \cdots \oplus 
		\Sigma^{s_m}\field[x]
\end{align*}
Observe that the set of grading shifts $(s_1, \ldots, s_m)$ are identified from the degree of the homogeneous basis element $\deg(\beta_j x^{s_j})$.
We provide an example of calculating a graded invariant factor decomposition of a module below.

\begin{example}
	Let $M$ be the graded $\rationals[x]$-module as given in \fref{ex:first-graded-snd-calc}.
	We have calculated the following graded decomposition on $M$ by \fref{prop:how-to-get-graded-ifd}:
	\begin{equation*}
		M \upgraded\cong 
			\paren{\frac{
				\rationals[x]\ket{\beta_1 x}
			}{
				\rationals[x]\ket{\beta_1 x}
			}}
			\oplus 
			\paren{\frac{
				\rationals[x]\ket{\beta_2 x^2}
			}{
				\rationals[x]\ket{x^4 \cdot \beta_2 x^2}
			}}
			\oplus
			\paren{\frac{
				\rationals[x]\ket{\beta_3 x^3}
			}{
				\rationals[x]\ket{x^5 \cdot \beta_3 x^3}
			}}
			\oplus 
			\rationals[x]\ket{\beta_4 x^4}
	\end{equation*}
	Then, by application of \fref{lemma:summands-on-graded-ifds}, 
	the graded invariant factor decomposition on $M$ is given as follows:
	\begin{equation*}
		M \upgraded\cong 
			\Sigma^{2}\Bigl(
				\rationals[x] \bigmod (x^4)
			\Bigr)
			\oplus 
			\Sigma^{3}\Bigl(
				\rationals[x] \bigmod (x^5)
			\Bigr)
			\oplus 
			\Sigma^4 \rationals[x]
	\end{equation*}
	Note that the first summand with the basis element $\beta_1 x$ becomes a trivial module.
	Recall that in \fref{chapter:persistence-theory}, we discussed how graded $\rationals[x]$-modules correspond to persistence modules over $\rationals$.
	Then, by \fref{cor:interval-decomp-from-structure-theorem},
		we have the following isomorphism:
	\begin{equation*}
		\topersmod(M) \,\uppersmod\cong\, 
			\intmod{[2, 2+4)}
			\oplus 
			\intmod{[3, 3+5)}
			\oplus 
			\intmod{[4, \infty)}
		= 
			\intmod{[2, 6)}
			\oplus 
			\intmod{[3, 8)}
			\oplus 
			\intmod{[4, \infty)}
	\end{equation*}
	where $\uppersmod{\cong}$ refers to an isomorphism in the category $\catpersmod[\rationals]$ of persistence modules over $\rationals$.
\end{example}\clearpage
%

\section{Matrix Reduction of Graded Matrices}
\label{section:matrix-reduction-of-graded-matrices}

Let $M$ be a graded $\field[x]$-module with a graded presentation given by $\phi: F_S \to F_G$ and $\pi: F_G \to M$.
In \fref{prop:how-to-get-graded-ifd}, 
	we stated that 
	if there exists a homogeneous basis 
	$\set{\beta_j x^{s_j}}_{j=1}^m$ of $F_G$ and nonzero invariant factors $x^{t_1}, \ldots, x^{t_r} \in \field[x]$ 
	such that 
	$\set{\beta_1 x^{s_1 + t_1}, \ldots, \beta_r x^{s_r + t_r}}$ is a basis for $\im(\phi)$,
	then the graded invariant factor decomposition of $M$ can be calculated from said basis.
In this section and in \fref{section:snd-algorithm-in-the-graded-case},
	we argue that we can always calculate such a basis, assuming we know of homogeneous bases $A$ and $S$ of $F_S$ and $F_G$ respectively.

Our argument makes use of the general algorithm for the calculation of Smith Normal Decompositions described in \cite[Remark 5.3.4]{algebra:adkins} for matrices over arbitrary Euclidean domains,
	as discussed in \fref{section:matrix-calculation-of-IFDs} after \fref{ex:snd-by-reduction}.
This algorithm involves strategically performing matrix reduction on $[\phi]_{A,S}$, with $V$ being the product of elementary matrices corresponding to column reduction and $U\inv$ being that for row reduction.
Since $[\phi]_{A,S}$ is a matrix of a graded homomorphism,
	there are a number of restrictions on the entries of $[\phi]_{A,S}$ that makes matrix reduction more straightforward and allows a more simplified version of the algorithm by \cite[Remark 5.3.4]{algebra:adkins}.
	
Assume $A$ and $S$ are homogeneous bases.
In this section, we consider four types of elementary matrices used in matrix reduction and discuss how matrix reduction preserves the homogeneity of these bases.
Then, in \fref{section:snd-algorithm-in-the-graded-case}, we present 
\textit{\fref{alg:general-snf-graded-homs}.\ Matrix Reduction Algorithm for Graded SNDs}, which calculates the desired SND $(U,D,V)$ of $[\phi]_{A,S}$ using only said types of elementary matrices.

Our first result is about the entries of matrices over $\field[x]$ corresponding to graded presentations.

\begin{lemma}\label{lemma:matrices-of-graded-homs}
	Let $\phi: N \to M$ be a graded homomorphism between two free graded $\field[x]$-modules $N$ and $M$ with ordered homogeneous bases 
		$S = (\sigma_1(x), \ldots, \sigma_n(x))$
	and $A = (\alpha_1(x), \ldots, \alpha_m(x))$ 
	respectively.
	Let $[\phi]$ be the matrix of $\phi$ relative to $S$ and $A$, i.e.\ $[\phi] := [\phi]_{A,S} \in \M_{m,n}(\field[x])$.

	For all row indices $j \in \set{1, \ldots, m}$ and column indices $i \in \set{1, \ldots, n}$ such that $[\phi](j,i) \neq 0$,
		$[\phi](j,i)$ is a homogeneous element of $\field[x]$ and 
		\begin{equation*}
			\degh\!\Big(
				\alpha_j(x)
			\Big) + \degh\!\Big(
				[\phi](j,i)
			\Big) = \degh\!\Big(
				\sigma_i(x)
			\Big)
		\end{equation*}
		where $\degh(\alpha_j(x))$ and $\degh(\sigma_i(x))$ are defined by homogeneity of the bases $S$ and $A$.
\end{lemma}
\begin{proof}
	For brevity, write $\sigma_i := \sigma_i(x)$ for $i \in \set{1, \ldots, n}$ and 
	$\alpha_j := \alpha_j(x)$ for $j \in \set{1, \ldots, m}$.

	Let $i \in \set{1, \ldots, n}$ be a column index such that $\phi(\sigma_i) \neq 0$.
	Since $\phi: N \to M$ is a graded homomorphism by assumption,
		$\col_i([\phi])$ is not a zero column and 
		$\degh(\sigma_i) = \degh(\phi(\sigma_i))$.
	By definition of matrix of module homomorphisms, 
		we have that $\col_i([\phi]) = [\phi(\sigma_i)]_A$ and 
		\begin{equation*}
			\phi(\sigma_i)  = 
			\sum_{j=1}^m \Bigl(
				\col_i([\phi])(j) \, \alpha_j 
			\Bigr)
			= \sum_{j=1}^m [\phi](j,i) \, \alpha_j
		\end{equation*}
	Let $j \in \set{1, \ldots, m}$ be a row index such that $[\phi](j,i) \neq 0$.
	Since $\phi(\sigma_i)$ is homogeneous of degree $\degh(\phi(\sigma_i))$, 
		the summand $[\phi](j,i)\,\alpha_j$ must be homogeneous of the same degree.
	Moreover, $[\phi](j,i)$ must also be homogeneous. Otherwise, the product 
	$[\phi](j,i)\,\alpha_j$ would not be homogeneous.
	Then, we have that:
	\begin{equation*}
		\degh(\sigma_i)
		=
		\degh\!\big(\phi(\sigma_i) \big)
		= \degh\!\Big( [\phi](j,i)\,\alpha_j \Big)
		= \degh\!\Big( [\phi](j,i) \Big) + \degh(\alpha_j)
	\,.
	\vspace{-\baselineskip}
	\end{equation*}
\end{proof}

Below, we provide an example of a graded homomorphism involving the persistent homology of some filtration, as discussed in \fref{section:construction-of-persistent-homology}.

\begin{example}
	Let $\filt{K}$ and $K$ be given as in \fref{ex:pershom-one}
	and orient $K$ with the vertex order $\Vertex(K) = (a,b,c,d)$.
	For convenience, an illustration of $K$ and $K_\bullet$ (without orientation) is copied below:
	\begin{center}
		\baselineCenter{\includegraphics[height=1.05in]{zomfig1/zom-simp.png}}
		\quad
		\baselineCenter{\includegraphics[height=1.1in]{zomfig1/zom1-cleaned.png}}
	\end{center}
	\vspace{3pt}
	The $1$\st and $2$\nd graded chain modules of $\filt{K}$ with coefficients in $\rationals$ are given as follows, relative to the ordered bases induced by the orientation $\Vertex(K) = (a,b,c,d)$ on the simplicial complex $K$:
	\begin{equation*}
		C_1\graded(\filt{K}; \rationals)
		= \rationals[x]\ket{
			abx, bcx, adx^2, cdx^2, acx^3
		}
		\qquad\text{ and }\qquad 
		C_2\graded(\filt{K}; \rationals)
		= \rationals[x]\ket{
			abc x^4, acd x^5
		}
	\end{equation*}
	The matrix $[\boundary_2\graded]$ of the graded boundary map 
	$\boundary_2\graded: 
		C_2\graded(\filt{K};\rationals) \to C_{1}\graded(\filt{K};\rationals)$
	relative to the standard ordered bases is given as follows:
	\NiceMatrixOptions{
		code-for-first-row = \color{gray}\,\,,
		code-for-first-col = \color{gray}
	}
	\begin{equation*}
		[\boundary_2\graded] = 
		\begin{pNiceMatrix}[first-row, first-col]
			{} 	& abc x^4 & acd x^5 \\
			{abx}	& 	x^3 & 0 \\
			{bcx}	& 	x^3 & 0 \\
			{adx^2}	& 	0 & -x^3 \\
			{cdx^2}	& 	0 & x^3 \\
			{acx^3}	&	-x & x^2 \\
		\end{pNiceMatrix}
	\end{equation*}
	Observe that the nonzero entries of $[\boundary_2\graded]$ satisfy the degree relation stated in \fref{lemma:matrices-of-graded-homs}. 
	We list these below.
	\begin{equation*}
		\def\arraystretch{1.3}
		\begin{array}{
			c !{=} r !{\quad:\quad}  
			!{\degh(} c !{)} 
			!{\,+\,} 
			!{\degh(} c !{)}  
			!{\,=\,} 
			!{\degh(} c !{)}  
		}
			{}[\boundary_2\graded](1,1) & x^3 
				& x^3 & ab x & abc x^4 
			\\ 
			{}[\boundary_2\graded](2,1) & x^3 
				& x^3 & bc x & abc x^4
			\\
			{}[\boundary_2\graded](3,2) & -x^3 
				& -x^3 & ad x^2 & acd x^5
			\\
			{}[\boundary_2\graded](4,2) & x^3 
				& x^3 & cd x^2 & acd x^5
			\\ 
			{}[\boundary_2\graded](5,1) & -x^{\phantom{1}}
				& -x & ac x^3 & abc x^4
			\\
			{}[\boundary_2\graded](5,2) & x^2
				& x^2 & ac x^3 & acd x^5
		\end{array}
	\end{equation*}
	Note that \fref{lemma:matrices-of-graded-homs} does not apply for zero entries since $\degh(0)$ is undefined.
\end{example}

Since \fref{lemma:matrices-of-graded-homs} implies that the nonzero elements of the matrices over $\field[x]$ are homogeneous and that their degrees are fixed by the degrees of the initial homogeneous bases (with the added assumptions stated on the same lemma),
	we can show that elimination operations done on the matrix reduction algorithm for graded SNDs must also preserve the homogeneity of the initial bases, 
	i.e.\ the bases that they induce are also homogeneous assuming we start with homogeneous bases.

Column operations on a matrix correspond to multiplication of said matrix on the \textit{right} by an elementary matrix (see \fref{prop:column-operations} in \fref{appendix:matrix-theory}).
For clarity, we identify certain terminology involving column reduction of matrices.

\begin{definition}\label{defn:column-reduction-operation}
	A \textbf{column reduction operation} on a matrix $T \in \M_{m,n}(\field[x])$ consists of the following:
	\begin{enumerate}
		\item a fixed row index $r \in \set{1,\ldots,m}$, 
		\item a nonzero \textbf{target entry} $T(r,k)$ with \textbf{target column index} $k \in \set{1,\ldots,n}$ to be eliminated in $T$,
		\item a nonzero \textbf{pivot entry} $T(r,p)$ with \textbf{pivot column index} $p \in \set{1,\ldots,n}$, 
		\item a \textbf{pivot multiplier} $f \in \field[x]$ such that 
		$T(r,k) + f \cdot T(r,p) = 0$, and 
		\item an elementary transvection $V := \eladd[n]{p,k \,; f}$, multiplied to $T$ on the right. 
		Note that the pivot column index $p$ is the first argument and the target column index $k$ is the second.
	\end{enumerate}
	Then, the entry $T(r,k)$ is eliminated in the product 
	$TV = T\eladd[n]{p,k \,; f}$, i.e.\ $(TV)(r,k) = 0$. 
\end{definition}

Following the notation above, 
we can show that $(TV)(r,k) = 0$ by looking at the columns of the product $TV$ and using a column-wise description of $V = \eladd[n]{p,k \,; f}$: for all column indices $i \in \set{1, \ldots, n}$,
\begin{equation*}
	\col_i(TV) 
	=
	T \col_i(V)
	=
	T \col_i\paren{
		\eladd{p,k \,; f}
	} = \begin{cases}
		\col_k(T) + f \cdot \col_p(T) 
			&\text{ if } i = k \\
		\col_i(T) 	&\text{ otherwise }
	\end{cases}
\end{equation*}
Note that $T$ and $TV$ can differ only on the $k$\th column.
Then, the $(r,k)$\th entry of $TV$ would be zero, i.e.\ the target entry $T(r,k)$ is eliminated as follows:
\begin{equation*}
	\col_k(T)(r) + f \cdot \col_p(T)(r)
	= T(r,k) + f \cdot T(r,p)
	= 0
\end{equation*}
For arbitrary matrices over $\field[x]$, 
elimination operations as above generally cannot be done since $f$ exists if and only if the pivot entry $T(r,p)$ divides the target entry $T(r,k)$.
However, assuming $T = [\phi]$ is as given in \fref{lemma:matrices-of-graded-homs}, 
	the elements $[\phi](r,p)$ and $[\phi](r,k)$ are homogeneous
	and we have the following equivalence:
\begin{equation*}
	\left\lgroup\begin{gathered}
		[\phi](r,p) \,\divides\, [\phi](r,k) \\
		\text{ divisibility relation }
	\end{gathered}\right\rgroup
	\quad\text{ if and only if }\quad 
	\left\lgroup\,\begin{gathered}
		\degh\bigl( [\phi](r,p) \bigr) 
		\leq \degh\bigl( [\phi](r,k) \bigr) \\ 
		\text{ degree relation }
	\end{gathered}\,\right\rgroup
\end{equation*}
Then, $f \in \field[x]$ is given by $f = - [\phi](r,k) \bigmod [\phi](r,p)$ assuming that $\degh([\phi](r,p)) \leq \degh([\phi](r,k))$.
The question now is whether the basis induced by $V \in \GL(n,\field[x])$ and the remaining entries of $\col_k([\phi]V)$ correspond to homogeneous elements. The following proposition addresses this.

\begin{proposition}\label{prop:column-reduction-preserves-hom}
	\textbf{Column Reduction Preserves Homogeneity.}

	Let $\phi: N \to M$ be a graded homomorphism between two free graded $\field[x]$-modules $N$ and $M$ with ordered homogeneous bases 
		$S = (\sigma_1(x), \ldots, \sigma_n(x))$
	and $A = (\alpha_1(x), \ldots, \alpha_m(x))$ 
	respectively.
	Let $[\phi]$ be the matrix of $\phi$ relative to $S$ and $A$, i.e.\ $[\phi] := [\phi]_{A,S} \in \M_{m,n}(\field[x])$.

	Let $p,k \in \set{1, \ldots, n}$ be distinct column indices 
	and $r \in \set{1, \ldots, m}$ be a row index such that 
	the target $[\phi](r, k)$
	and 
	the pivot $[\phi](r, p)$
	are both nonzero 
	and $\degh([\phi](r,p)) \leq \degh([\phi](r,k))$.
	If $f \in \field[x]$ is nonzero and homogeneous such that
	\begin{equation*}
		\degh\Bigl( [\phi](r,k)
		\Bigr) = \degh(f) + 
		\degh\Bigl(
			[\phi](r,p)
		\Bigr)
		\,
	\end{equation*}
	then the basis $T = (\tau_1(x), \ldots, \tau_n(x))$ of $N$ induced by the elementary transvection $V := \eladd[n]{p,k \,; f}$ by $[\tau_i]_S = \col_i(V)$
	is homogeneous.
\end{proposition}
\begin{proof}
	Let $f \in \field[x]$ be nonzero and homogeneous such that
	$
		\degh([\phi](r,k)) = \degh(f) + \degh([\phi](r,p))
	$.
	For brevity, let $\sigma_i := \sigma_i(x)$, $\tau_i := \tau_i(x)$ and $\alpha_j := \alpha_j(x)$.

	Since $V := \eladd{p,k \,; f}$ differs from the identity matrix $I_n$ only in the $(p,k)$\th entry,
	we have that $\tau_i = \sigma_i$ and $\tau_i$ is homogeneous
	for all $i \in \set{1, \ldots, n}$ with $i \neq a$.
	We need to check that $\tau_k \in N$ is also homogeneous.
	Since $[\tau_k]_S = \col_k(V) = \col_k(I_n) + f \cdot \col_p(I_n)$, we have that $\tau_k = \sigma_k + f \cdot \sigma_p$.
	For $\tau_k$ to be homogeneous, we need to show that 
	$\degh(\sigma_k) = \degh(f \cdot \sigma_p)$:
	\begin{align*}
		\degh(f \cdot \sigma_p) 
		&= \degh(f) + \degh(\sigma_p) 
			&&\text{ by homogeneity of $f$ and $\sigma_p$ } \\
		&= \degh([\phi](r,k)) - \degh([\phi](r,p))
			&&\text{ by assumption on } f \in \field[x] 
			\\
			&\qquad + \degh(\sigma_p) 
			\\
		&=	\degh([\phi](r,k)) - \degh([\phi](r,p))  \\
			&\qquad + \degh([\phi](r,p)) + \degh(\alpha_r) 
			&& \text{ by \fref{lemma:matrices-of-graded-homs} on the $[\phi](r,p) \neq 0$}
			\\
		&= \degh([\phi](r,k)) + \degh(\alpha_r)
			\\ 
		&= \degh(\sigma_k)
			&&\text{ by \fref{lemma:matrices-of-graded-homs} on $[\phi](r,k)\neq 0$}
	\end{align*}
	Since $\degh(\sigma_k) = \degh(f \cdot \sigma_p)$, $\tau_k \in N$ is homogeneous with $\degh(\tau_k) = \degh(\sigma_k)$.
	Therefore, $T = (\tau_i)_{i=1}^n$ is a homogeneous basis of $N$.
\end{proof}

Note that the statement of \fref{prop:column-reduction-preserves-hom} applies even when $f \in \field[x]$, as denoted above, is not used to eliminate entries, i.e.\ it is not required that $[\phi](r,k) + f \cdot [\phi](r,p) = 0$,
only that $[\phi](r,k) + f \cdot [\phi](r,p)$ results in either zero or a homogeneous element of the same degree.
While this is not critical for the reduction algorithm, it does allow some steps that may be more preferable when working by hand.

Observe that the key part of the proof of \fref{prop:column-reduction-preserves-hom} is the existence of a row index $r \in \set{1, \ldots, m}$ such that both $[\phi](r,k)$ and $[\phi](r,p)$ are nonzero,
	in that they determine the degree of the element $f \in \field[x]$.
Note that we use the degree relation from \fref{lemma:matrices-of-graded-homs} on both $[\phi](r,k)$ and $[\phi](r,p)$, which require them to both be nonzero. 
If such a row index $r$ does not exist, then it is possible for the nonzero entries of $[\phi]V$ to be homogeneous but the basis $T = (\tau_i)$ induced by $V$ to not be a homogeneous basis.
A trivial example would involve a zero column on $[\phi]$, wherein if the nonzero entries of $[\phi]$ are homogeneous and $\col_p([\phi])$ is the zero column, 
	then $[\phi]\eladd[n]{p, k \; f}$ for any column index $k$ and any element $f \in \field[x]$ would also have entries that are either zero or homogeneous.
	We saw a similar issue happen with row reduction earlier in \fref{ex:first-graded-snd-calc} wherein the matrix $U_2 \in \GL(4,\field[x])$ produced a non-homogeneous basis since the $4$\th row is the zero row.
Fortunately, column operations of this type are generally not required for matrix reduction since these tend to create more nonzero entries, i.e.\ the opposite goal of matrix reduction.

As denoted above, $[\phi]V = [\phi]_{A,S} V = [\phi]_{A,T}$, i.e.\ 
	$[\phi]V$ corresponds to the matrix of $\phi$ relative to the bases $T = (\tau_i(x))_{i=1}^n$ of $N$ and $A = (\alpha_j(x))_{j=1}^m$ of $M$, 
Since we have shown that $T$ is a homogeneous basis of $N$, 
	\fref{lemma:matrices-of-graded-homs} also applies for $[\phi]V = [\phi]_{A,T}$, i.e.\ 
	for all nonzero entries $[\phi]_{A,T}(j,i)$, 
	\begin{equation*}
		\degh\!\Big(
			\alpha_j(x)
		\Big) + \degh\!\Big(
			[\phi]_{A,T}(j,i)
		\Big) = \degh\!\Big(
			\tau_i(x)
		\Big)
		\,.
	\end{equation*}
An immediate consequence of this is that we can do elimination by column reduction finitely many times and the resulting matrices from said operation will also preserve the homogeneity of the initial bases.
Below, we provide an example of this column reduction process in action.
Note that the following example does not use the algorithm by \cite[Remark 5.3.4]{algebra:adkins} (\fref{alg:general-snf-graded-homs}) since we have not established that row reduction preserves homogeneity.
Instead, this example more closely resembles the column reduction algorithm discussed in the next section.

\begin{example}\label{ex:column-reduction-of-graded-boundary-one}
	Let $\filt{K}$ and $K$ be given as in \fref{ex:pershom-one} and as illustrated below.
	\begin{center}
		\baselineCenter{\includegraphics[height=1.05in]{zomfig1/zom-simp.png}}
		\quad
		\baselineCenter{\includegraphics[height=1.1in]{zomfig1/zom1-cleaned.png}}
	\end{center}
	\vspace{3pt}
	The $0$\th and $1$\st graded chain modules of $\filt{K}$ with coefficients in $\rationals$ are given as follows, relative to ordered bases induced by the orientation $\Vertex(K) = (a,b,c,d)$ on the simplicial complex $K$:
	\begin{equation*}
		C_0\graded(\filt{K}; \rationals)
		= \rationals[x]\ket{
			a, b, cx, dx
		}
		\qquad\text{ and }\qquad 
		C_1\graded(\filt{K}; \rationals)
		= \rationals[x]\ket{
			abx, bcx, adx^2, cdx^2, acx^3
		}
	\end{equation*}
	Below, we identify the matrix $[\boundary_1\graded]$ of the $1$\st graded boundary morphism 
	$\boundary_1\graded: C_1\graded(\filt{K}; \rationals) \to C_0\graded(\filt{K}; \rationals)$
	relative to the same ordered bases:
	\begin{equation*}
		[\boundary_1\graded] 
		= \begin{NiceArray}{>{\color{black}}c ccccc}
			\RowStyle[color=black]{}
			& abx & bcx & adx^2 & cdx^2 & acx^3 \\
			{a}		&	-x & 0 & -x^2 & 0 & -x^3 \\
			{b}		&	x & -x & 0 & 0 & 0 \\ 
			{cx}	&	0 & 1 & 0 & x & x^2 \\
			{dx}	&	0 & 0 & x & -x & 0
		\CodeAfter \SubMatrix({2-2}{5-6})
		\end{NiceArray}
	\end{equation*}
	Given below is a sequence $Q_1, Q_2, \ldots$ of matrices $Q_n \in \M_{4,5}(\field[x])$ resulting from successive column reduction operations on $Q_0 := [\boundary_1\graded]$, with the following color scheme:
	\vspace{5pt}
	\begin{center}
		\bluetagged{chosen pivot entry}, 
		\redtagged{target entry for elimination},
		\greentagged{pivot multiplier},
		\orangetagged{affected column}
	\end{center} 
	\allowdisplaybreaks
	{\arraycolsep=3pt\begin{longtable}{L !{$=$} C !{$=$} C}
		Q_1 := Q_0 \eladd[5]{3,4 \,; 1}
			& \begin{pNiceMatrix}
				\CodeBefore [create-cell-nodes]
					\tikz \node [targetcell = (4-4)] {} ;
					\tikz \node [pivotcell = (4-3)] {} ;
				\Body
				-x & 0 & -x^2 & 0 & -x^3 \\
				x & -x & 0 & 0 & 0 \\
				0 & 1 & 0 & x & x^2 \\
				0 & 0 & x & -x & 0
			\end{pNiceMatrix} 
			\begin{pNiceMatrix}
				\CodeBefore [create-cell-nodes]
					\tikz \node [differentcell = (3-4)] {} ;
				\Body
				1 & 0 & 0 & 0 & 0 \\
				0 & 1 & 0 & 0 & 0 \\
				0 & 0 & 1 & 1 & 0 \\
				0 & 0 & 0 & 1 & 0 \\
				0 & 0 & 0 & 0 & 1
			\end{pNiceMatrix}
			& \begin{pNiceMatrix} 
				\CodeBefore [create-cell-nodes]
					\tikz \node [orangecell = (1-4) (4-4)] {} ;
				\Body
				-x & 0 & -x^2 & -x^2 & -x^3 \\
				x & -x & 0 & 0 & 0 \\
				0 & 1 & 0 & x & x^2 \\
				0 & 0 & x & 0 & 0
			\end{pNiceMatrix} 
		\\[30pt] 
		Q_2 := Q_1 \eladd[5]{2,4\,;-x}
			& \begin{pNiceMatrix}
				\CodeBefore [create-cell-nodes]
					\tikz \node [targetcell = (3-4)] {} ;
					\tikz \node [pivotcell = (3-2)] {} ;
				\Body
				-x & 0 & -x^2 & -x^2 & -x^3 \\
				x & -x & 0 & 0 & 0 \\
				0 & 1 & 0 & x & x^2 \\
				0 & 0 & x & 0 & 0
			\end{pNiceMatrix} 
			 \begin{pNiceMatrix}
				\CodeBefore [create-cell-nodes]
					\tikz \node [differentcell = (2-4)] {} ;
				\Body
				1 & 0 & 0 & 0 & 0 \\
				0 & 1 & 0 & -x & 0 \\
				0 & 0 & 1 & 0 & 0 \\
				0 & 0 & 0 & 1 & 0 \\
				0 & 0 & 0 & 0 & 1
			\end{pNiceMatrix}
			& \begin{pNiceMatrix}
				\CodeBefore [create-cell-nodes]
					\tikz \node [orangecell = (1-4) (4-4)] {} ;
				\Body
				-x & 0 & -x^2 & -x^2 & -x^3 \\
				x & -x & 0 & x^2 & 0 \\
				0 & 1 & 0 & 0 & x^2 \\
				0 & 0 & x & 0 & 0
			\end{pNiceMatrix} 
		\\[30pt]
		Q_3 := Q_2 \eladd[5]{1,4\,; -x}
			& \begin{pNiceMatrix}
				\CodeBefore [create-cell-nodes]
					\tikz \node [pivotcell  = (2-1)] {} ;
					\tikz \node [targetcell = (2-4)] {} ;
				\Body
				-x & 0 & -x^2 & -x^2 & -x^3 \\
				x & -x & 0 & x^2 & 0 \\
				0 & 1 & 0 & 0 & x^2 \\
				0 & 0 & x & 0 & 0
			\end{pNiceMatrix} 
			 \begin{pNiceMatrix}
				\CodeBefore [create-cell-nodes]
					\tikz \node [differentcell = (1-4)] {} ;
				\Body
				1 & 0 & 0 & -x & 0 \\
				0 & 1 & 0 & 0 & 0 \\
				0 & 0 & 1 & 0 & 0 \\
				0 & 0 & 0 & 1 & 0 \\
				0 & 0 & 0 & 0 & 1
			\end{pNiceMatrix}
			& \begin{pNiceMatrix} 
				\CodeBefore [create-cell-nodes]
					\tikz \node [orangecell = (1-4) (4-4)] {} ;
				\Body 
				-x & 0 & -x^2 & 0 & -x^3 \\
				x & -x & 0 & 0 & 0 \\
				0 & 1 & 0 & 0 & x^2 \\
				0 & 0 & x & 0 & 0
			\end{pNiceMatrix} 
		\\[30pt] 
		Q_4 := Q_3 \eladd[5]{2,5\,; -x^2}
			& \begin{pNiceMatrix}
				\CodeBefore [create-cell-nodes]
					\tikz \node [pivotcell  = (3-2)] {} ;
					\tikz \node [targetcell = (3-5)] {} ;
				\Body
				-x & 0 & -x^2 & 0 & -x^3 \\
				x & -x & 0 & 0 & 0 \\
				0 & 1 & 0 & 0 & x^2 \\
				0 & 0 & x & 0 & 0
			\end{pNiceMatrix} 
			 \begin{pNiceMatrix}
				\CodeBefore [create-cell-nodes]
					\tikz \node [differentcell = (2-5)] {} ;
				\Body
				1 & 0 & 0 & 0 & 0 \\
				0 & 1 & 0 & 0 & -x^2 \\
				0 & 0 & 1 & 0 & 0 \\
				0 & 0 & 0 & 1 & 0 \\
				0 & 0 & 0 & 0 & 1
			\end{pNiceMatrix}
			& \begin{pNiceMatrix}
				\CodeBefore [create-cell-nodes]
					\tikz \node [orangecell = (1-5) (4-5)] {} ;
				\Body
				-x & 0 & -x^2 & 0 & -x^3 \\
				x & -x & 0 & 0 & x^3 \\
				0 & 1 & 0 & 0 & 0 \\
				0 & 0 & x & 0 & 0
			\end{pNiceMatrix} 
		\\[30pt]
		Q_5 := Q_4 \eladd[5]{1,5\,;-x^2} 
			& \begin{pNiceMatrix}
				\CodeBefore [create-cell-nodes] 
					\tikz \node [pivotcell  = (2-1)] {} ;
					\tikz \node [targetcell = (2-5)] {} ;
				\Body
				-x & 0 & -x^2 & 0 & -x^3 \\
				x & -x & 0 & 0 & x^3 \\
				0 & 1 & 0 & 0 & 0 \\
				0 & 0 & x & 0 & 0
			\end{pNiceMatrix} 
			 \begin{pNiceMatrix}
				\CodeBefore [create-cell-nodes]
					\tikz \node [differentcell = (1-5)] {} ;
				\Body
				1 & 0 & 0 & 0 & -x^2 \\
				0 & 1 & 0 & 0 & 0 \\
				0 & 0 & 1 & 0 & 0 \\
				0 & 0 & 0 & 1 & 0 \\
				0 & 0 & 0 & 0 & 1
			\end{pNiceMatrix}
			& \begin{pNiceMatrix}
				\CodeBefore [create-cell-nodes]
					\tikz \node [orangecell = (1-5) (4-5)] {} ;
				\Body
				-x & 0 & -x^2 & 0 & 0 \\
				x & -x & 0 & 0 & 0 \\
				0 & 1 & 0 & 0 & 0 \\
				0 & 0 & x & 0 & 0
			\end{pNiceMatrix}
	\end{longtable}}\noindent
	Then, we have $Q_5 = [\boundary_1\graded] V$ with $V \in \GL(5,\field[x])$
	given as follows, with each successive elementary transvection matrix multiplied on the right of $Q_i$:
	\begin{equation*}
		V = \color{black}
			\overbrace{\color{black}\eladd[5]{3,4 \,; 1}}^{\mathclap{\text{
				first column elimination
			}}}
			\underbrace{\color{black}\eladd[5]{2,4 \,; -x}}_{\mathclap{\text{
				second column elimination
			}}}
			\cdots 
			\overbrace{\color{black}\eladd[5]{1,5 \,; -x^2}}^{\mathclap{\text{
				fifth and last column elim.
			}}}
		\color{black}
		=\vspace{5pt}
		\begin{NiceArray}{
				>{\color{black}}c (ccccc)
			}
			{abx} 	&	1 & 0 & 0 & -x & -x^2 \\
			{bcx} 	&	0 & 1 & 0 & -x & -x^2 \\
			{adx^2}	&	0 & 0 & 1 & 1 & 0 \\
			{cdx^2}	&	0 & 0 & 0 & 1 & 0 \\
			{acx^3}	&	0 & 0 & 0 & 0 & 1
		\end{NiceArray}
	\end{equation*}
	Then, $V \in \GL(5, \field[x])$ induces a new basis $T = (\tau_1(x), \ldots, \tau_5(x))$ by $[\tau_i(x)] = \col_i(V)$, i.e.\ 
	\begin{equation*}\def\arraystretch{1.1}
		\begin{array}{c !{=} c !{=} r}
			\tau_1(x) & (1)(abx) & (ab)x\phantom{^1} \\
			\tau_2(x) & (1)(bcx) & (bc)x\phantom{^1} \\
			\tau_3(x) & (1)(adx^2) & (ad)x^2 \\
			\tau_4(x) & (-x)(abx) + (-x)(bcx) + (1)(adx^2)
				& (-ab -bc + ad)x^2 \\
			\tau_5(x) & (-x^2)(abx) + (-x^2)(bcx) + (1)(acx^3)
				& (-ab - bc + ac)x^3
		\end{array}
	\end{equation*}
	Observe that $T$ is a homogeneous basis for $C_2\graded(\filt{K}; \rationals)$
	and that the nonzero entries of $Q_5$ are as expected by \fref{lemma:matrices-of-graded-homs}.
	For example, at the $(1,3)$\th entry of $Q_5$, we have that 
	\begin{equation*}
		\degh(a) + \degh\bigl( Q_5(1,3) \bigr)
		= 0 + \degh(-x^2) = 3 = \degh\bigl(\tau_3(x)\bigr) = \degh(ad x^2)
	\end{equation*}
\end{example}

We have a similar result for elimination by row reduction.
Row operations on a matrix correspond to multiplication of said matrix on the \textit{left} by an elementary matrix (\fref{prop:row-operations} in \fref{appendix:matrix-theory}).
For clarity, we identify certain terminology involving elimination by row reduction below.

\begin{definition}\label{defn:row-reduction-operation}
	A \textbf{row reduction operation} on a matrix $T \in \M_{m,n}(\field[x])$ consists of the following:

	\begin{enumerate}
		\item a fixed column index $c \in \set{1,\ldots,n}$, 
		\item a nonzero \textbf{target entry} $T(k,c)$ with \textbf{target row index} $k \in \set{1,\ldots,m}$ to be eliminated in $T$,
		\item a nonzero \textbf{pivot entry} $T(p,c)$ with \textbf{pivot row index} $p \in \set{1,\ldots,m}$, 
		\item a \textbf{pivot multiplier} $f \in \field[x]$ such that 
		$T(k,c) + f \cdot T(p,c) = 0$, and 
		\item an elementary transvection
		$U\inv := \eladd[m]{k,p \,; f}$, multiplied to $T$ on the left.
		Note that the pivot row index $p$ is the \textit{second} argument and the target row index $k$ is the \textit{first}, and that $U = \eladd[m]{k,p \,; -f}$.
	\end{enumerate}
	Then, the entry $T(k,c)$ is eliminated in the product 
	$U\inv T = \eladd[m]{k,p \,; f}T$, i.e.\ $(U\inv T)(k,c)=0$.
\end{definition}
\remark{
	We define $U\inv := \eladd[m]{k,p \,; f}$, as opposed to setting $U := \eladd[m]{k,p \,; f}$, since we want to interpret the resulting matrix as a change of basis on the codomain of homomorphisms.
}

Following the notation above, 
we can show that $(U\inv T)(k,c)=0$ by examining the 
the rows of the product $U\inv T \in \M_{m,n}(\field[x])$ and using a row-wise description of $U\inv = \eladd[m]{k,p \,; f}$: 
for all row indices $j \in \set{1, \ldots, m}$,
\begin{equation*}
	\row_j(U\inv T)
	= \row_j(U\inv) T 
	= \row_j\paren{
		\eladd[m]{k,p \,; f}
	} T 
	= \begin{cases}
		\row_k(T) + f \cdot \row_p(T) &\text{ if } j = k \\
		\row_j(T) &\text{ otherwise }
	\end{cases}
\end{equation*}
As with the case of column reduction, when $T = [\phi]$ is as given in \fref{lemma:matrices-of-graded-homs},
	the elimination process becomes more straightforward.
We have the following result involving row reduction of matrices given by \fref{lemma:matrices-of-graded-homs}.

\begin{proposition}\label{prop:row-reduction-preserves-hom}
	\textbf{Row Reduction Preserves Homogeneity.}

	Let $\phi: N \to M$ be a graded homomorphism between two free graded $\field[x]$-modules $N$ and $M$ with ordered homogeneous bases 
		$S = (\sigma_1, \ldots, \sigma_n)$
	and $A = (\alpha_1, \ldots, \alpha_m)$ 
	respectively.
	Let $[\phi]$ be the matrix of $\phi$ relative to $S$ and $A$, i.e.\ $[\phi] := [\phi]_{A,S} \in \M_{m,n}(\field[x])$.

	Let $p,k \in \set{1, \ldots, m}$ be distinct row indices and 
		$c \in \set{1, \ldots, n}$ be a column index 
		such that the pivot $[\phi](p,c)$ and the target $[\phi](k,c)$ are both nonzero and $\degh([\phi](p,c)) \leq \degh([\phi](k,c))$.
	If $f \in \field[x]$ is nonzero and homogeneous such that 
	\begin{equation*}
		\degh\Bigl( [\phi](k,c)
		\Bigr) = \degh(f) + 
		\degh\Bigl(
			[\phi](p,c)
		\Bigr)
		\,
	\end{equation*}
	then the basis $B = (\beta_1(x), \ldots, \beta_m(x))$ of $M$ induced by the elementary transvection $U := \eladd[n]{k,p \,; -f}$ by $[\beta_j]_A = \col_i(U)$ is homogeneous.
\end{proposition}
\begin{proof}
	Let $f \in \field[x]$ be nonzero and homogeneous such that
	$
		\degh([\phi](k,c)) = \degh(f) + \degh([\phi](p,c))
	$.
	For brevity, let $\sigma_i := \sigma_i(x)$, $\alpha_j := \alpha_j(x)$, 
	and $\beta_j := \beta_j(x)$.
	Note that while $U := \eladd{k,p \,; -f} \in \GL(m,\field[x])$ corresponds to a row operation, the basis $B = (\beta_1, \ldots, \beta_m)$ is determined by the columns of $U$.

	For $j \in \set{1, \ldots, m}$ with $j \neq p$, 
		$\col_j(U) = \col_j(I_m)$ and $\beta_j = \alpha_j$ is homogeneous.
	Consider the $p$\th column of $U$.
	Since $[\beta_p]_A = \col_p(U) = \col_p(I_m) - f \cdot \col_k(I_m)$,
		$\beta_p = \alpha_p - f \cdot \alpha_k$.
	For $\beta_p \in M$ to be homogeneous, we need to show that 
		$\degh(\alpha_p) = \degh(-f \cdot \alpha_k)$.
	Note that $\degh(-f) = \degh(f)$.
	\begin{align*}
		\degh(f \cdot \alpha_k)
		&= \degh(f) + \degh(\alpha_k) 
				&&\text{ by homogeneity of $f$ and $\alpha_k$ } \\
		&= \degh\bigl( [\phi](k,c) \bigr) - \bigl( [\phi](p,c) \bigr) 
				&&\text{ by assumption on } f \in \field[x] \\
			&\qquad 	+ \degh(\alpha_k) \\ 
		&= \degh\bigl( [\phi](k,c) \bigr) - \degh\bigl( [\phi](p,c) \bigr) \\ 
			&\qquad 	+ \degh(\sigma_c) 
					- \degh\bigl( [\phi](k,c) \bigr)
				&&\text{ by \fref{lemma:matrices-of-graded-homs} on } [\phi](k,c) \neq 0
		\\ 
		&= - \degh\bigl( [\phi](p,c) \bigr) + \degh(\sigma_c) \\
		&= \degh(\alpha_p) 
				&&\text{ by \fref{lemma:matrices-of-graded-homs} on } [\phi](p,c) \neq 0
	\end{align*}
	Since $\degh(\alpha_p) = \degh(-f \cdot \alpha_k)$,
		$\beta_p \in M$ is homogeneous with $\degh(\beta_p) = \degh(\alpha_p)$.
	Therefore, $B = (\beta_j)_{j=1}^m$ is a homogeneous basis.
\end{proof}

Observe that, as denoted above, $U\inv[\phi] = U\inv [\phi]_{A,S} = [\phi]_{B,S}$, 
	i.e.\ 
	$U\inv[\phi]$ is the matrix of $\phi$ relative to the bases 
	$S = (\sigma_i(x))_{i=1}^n$ of $N$ and 
	$B = (\beta_j(x))_{j=1}^m$ of $M$.
Since we have shown that $B$ is a homogeneous basis, then \fref{lemma:matrices-of-graded-homs} applies to $U\inv[\phi] = [\phi]_{B,S}$, i.e.\ 
for all nonzero entries $[\phi]_{B,T}(j,i)$:
\begin{equation*}
	\degh\!\Big(
		\beta_j(x)
	\Big) + \degh\!\Big(
		[\phi]_{B,S}(j,i)
	\Big) = \degh\!\Big(
		\sigma_i(x)
	\Big)
	\,.
\end{equation*}
As with the case for column reduction, 
we can do elimination by row reduction finitely many times and the resulting matrices will also preserve the homogeneity of the initial bases and that of its entries.
Below, we continue \fref{ex:column-reduction-of-graded-boundary-one} and perform row reduction.

\begin{example}\label{ex:row-red-of-graded-bdry-one}
	We continue from \fref{ex:column-reduction-of-graded-boundary-one}.
	We currently have that $Q_5 = [\boundary_1\graded] V$ with 
	\begin{equation*}
		Q_5 = \begin{pmatrix}
			-x & 0 & -x^2 & 0 & 0 \\
			x & -x & 0 & 0 & 0 \\
			0 & 1 & 0 & 0 & 0 \\
			0 & 0 & x & 0 & 0
		\end{pmatrix}
		\qquad 
		[\boundary_1\graded] =
			\begin{pmatrix}
				-x & 0 & -x^2 & 0 & -x^3 \\
				x & -x & 0 & 0 & 0 \\
				0 & 1 & 0 & x & x^2 \\
				0 & 0 & x & -x & 0
			\end{pmatrix}
		\qquad 
		V = \begin{pmatrix}
				1 & 0 & 0 & -x & -x^2 \\
				0 & 1 & 0 & -x & -x^2 \\
				0 & 0 & 1 & 1 & 0 \\
				0 & 0 & 0 & 1 & 0 \\
				0 & 0 & 0 & 0 & 1
		\end{pmatrix}
	\end{equation*}
	We perform the following row operations on $Q_5$, as written below with the following color scheme:
	\vspace{5pt}
	\begin{center}
		\bluetagged{chosen pivot entry}, 
		\redtagged{target entry for elimination},
		\greentagged{pivot multiplier},
		\orangetagged{affected row}
	\end{center} 
	\allowdisplaybreaks 
	{ \arraycolsep=3pt
	\begin{longtable}{L !{$=$} C !{$=$} C}
		Q_6 := \eladd[4]{1,4 \,; x} Q_5 
		& \begin{pNiceMatrix}
			\CodeBefore [create-cell-nodes]
				\tikz \node [differentcell = (1-4)] {} ;
			\Body
			1 & 0 & 0 & x \\
			0 & 1 & 0 & 0 \\
			0 & 0 & 1 & 0 \\
			0 & 0 & 0 & 1 \\
		\end{pNiceMatrix}
		\begin{pNiceMatrix}
			\CodeBefore [create-cell-nodes]
				\tikz \node [targetcell = (1-3)] {} ;
				\tikz \node [pivotcell = (4-3)] {} ;
			\Body
			-x & 0 & -x^2 & 0 & 0 \\
			x & -x & 0 & 0 & 0 \\
			0 & 1 & 0 & 0 & 0 \\
			0 & 0 & x & 0 & 0
		\end{pNiceMatrix} 
		& \begin{pNiceMatrix} 
			\CodeBefore [create-cell-nodes]
				\tikz \node [orangecell = (1-1) (1-5)] {} ;
			\Body
			-x & 0 & 0 & 0 & 0 \\
			x & -x & 0 & 0 & 0 \\
			0 & 1 & 0 & 0 & 0 \\
			0 & 0 & x & 0 & 0
		\end{pNiceMatrix} 
	\\[20pt] 
	Q_7 := \eladd[4]{2,3 \,; x} Q_6 
		& \begin{pNiceMatrix}
			\CodeBefore [create-cell-nodes]
				\tikz \node [differentcell = (2-3)] {} ;
			\Body
			1 & 0 & 0 & 0 \\
			0 & 1 & x & 0 \\
			0 & 0 & 1 & 0 \\
			0 & 0 & 0 & 1 \\
		\end{pNiceMatrix}
		\begin{pNiceMatrix}
			\CodeBefore [create-cell-nodes]
				\tikz \node [targetcell = (2-2)] {} ;
				\tikz \node [pivotcell = (3-2)] {} ;
			\Body
			-x & 0 & 0 & 0 & 0 \\
			x & -x & 0 & 0 & 0 \\
			0 & 1 & 0 & 0 & 0 \\
			0 & 0 & x & 0 & 0
		\end{pNiceMatrix} 
		& \begin{pNiceMatrix} 
			\CodeBefore [create-cell-nodes]
				\tikz \node [orangecell = (2-1) (2-5)] {} ;
			\Body
			-x & 0 & 0 & 0 & 0 \\
			x & 0 & 0 & 0 & 0 \\
			0 & 1 & 0 & 0 & 0 \\
			0 & 0 & x & 0 & 0
		\end{pNiceMatrix} 
	\\[20pt]
	Q_8 := \eladd[4]{1,2 \,; 1} Q_7 
		& \begin{pNiceMatrix}
			\CodeBefore [create-cell-nodes]
				\tikz \node [differentcell = (1-2)] {} ;
			\Body
			1 & 1 & 0 & 0 \\
			0 & 1 & 0 & 0 \\
			0 & 0 & 1 & 0 \\
			0 & 0 & 0 & 1 \\
		\end{pNiceMatrix}
		\begin{pNiceMatrix}
			\CodeBefore [create-cell-nodes]
				\tikz \node [targetcell = (1-1)] {} ;
				\tikz \node [pivotcell = (2-1)] {} ;
			\Body
			-x & 0 & 0 & 0 & 0 \\
			x & 0 & 0 & 0 & 0 \\
			0 & 1 & 0 & 0 & 0 \\
			0 & 0 & x & 0 & 0
		\end{pNiceMatrix} 
		& \begin{pNiceMatrix} 
			\CodeBefore [create-cell-nodes]
				\tikz \node [orangecell = (1-1) (1-5)] {} ;
			\Body
			0 & 0 & 0 & 0 & 0 \\
			x & 0 & 0 & 0 & 0 \\
			0 & 1 & 0 & 0 & 0 \\
			0 & 0 & x & 0 & 0
		\end{pNiceMatrix} 
	\end{longtable} }\noindent 
	Note that the elementary transvection for each row operation is multiplied on the left of $Q_i$.
	Let $U \in \GL(4, \field[x])$ such that 
	$U\inv Q_5 = U\inv [\boundary_1\graded] V$.
	Then, $U\inv$ is given as follows:
	\begin{equation*}
		U\inv =	
			\overbrace{
				\eladd[4]{1,2 \,; 1}
			}^{\mathclap{\text{
				third row elimination
			}}}
			\cdot 
			\underbrace{
				\eladd[4]{2,3 \,; x}
			}_{\mathclap{\text{
				second row elimination
			}}}
			\cdot 
			\overbrace{
				\eladd[4]{1,4 \,; x}
			}^{\mathclap{\text{
				first row elimination
			}}}
		= \begin{pNiceMatrix}
			1 & 1 & x & x \\
			0 & 1 & x & 0 \\
			0 & 0 & 1 & 0 \\
			0 & 0 & 0 & 1 \\
		\end{pNiceMatrix}
	\end{equation*}
	Since the inverse of $\eladd[m]{a,b \,; f}$ is given by $\eladd[m]{a,b \,; -f}$, we can calculate $U$ as follows:
	\begin{align*}
		U &= \Biggl(
			\eladd[4]{1,2 \,; 1}\,
			\eladd[4]{2,3 \,; x}\,
			\eladd[4]{1,4 \,; x}
		\Biggr)\inv 
		= 	\overbrace{
				\eladd[4]{1,4 \,; -x}
			}^{\mathclap{\substack{ 
				\text{ inverse corr. to } \\ \text{ first row elimination }
			}}}
			\underbrace{
				\eladd[4]{2,3 \,; -x}
			}_{\mathclap{\substack{ 
				\text{ inverse corr. to } \\ \text{ second row elimination }
			}}}
			\overbrace{
				\eladd[4]{1,2 \,; -1}
			}^{\mathclap{\substack{ 
				\text{ inverse corr. to } \\ \text{ third row elimination }
			}}}
		\\ 
		&= \begin{NiceArray}{>{\color{gray}}c (cccc)}
			{a} & 1 & -1 & 0 & -x \\
			{b} & 0 & 1 & -x & 0 \\
			{cx} & 0 & 0 & 1 & 0 \\
			{dx} & 0 & 0 & 0 & 1 \\
		\end{NiceArray}
	\end{align*}
	Then, the basis $B = (\beta_1(x), \ldots, \beta_4(x))$ of $C_0\graded(\filt{K}; \rationals)$
	given by $[\beta_j(x)] = \col_j(U)$ is as follows:
	\begin{equation*}\def\arraystretch{1.1}
		\begin{array}{c !{=} c !{=} l}
			\beta_1(x) & (1)(a) 
				& a \\
			\beta_2(x) & (-1)(a) + (1)(b) 
				& -a + b \\
			\beta_3(x) & (-x)(b) + (1)(cx)
				& (-b+c)x \\
			\beta_4(x) & (-x)(a) + (1)(dx)
				& (-a+d)x
		\end{array}
	\end{equation*}
	Observe that $B$ is a homogeneous basis and that the nonzero entries of $Q_8$ are as expected by \fref{lemma:matrices-of-graded-homs}. 
\end{example}

The last two types of elementary matrices involved in matrix reduction are elementary permutations and elementary dilations.
In the proposition below, we claim that these two kinds of matrices also preserve homogeneity of bases.
Note that, unlike the case for elementary transvections in column and row reduction, the elementary permutations and dilations are not restricted by the entries of the matrix being reduced.

\begin{proposition}\label{prop:perm-and-dilation-preserve-hom}
	Let $N$ be a free graded $\field[x]$-module with ordered homogeneous basis $S = (\sigma_1, \ldots, \sigma_n)$.
	If $V \in \GL(n, \field[x])$ is an elementary permutation or elementary dilation over $\field[x]$, then the basis 
	$T = (\tau_1, \ldots, \tau_n)$ of $N$ by $[\tau_i]_S = \col_i(V)$ is homogeneous.
\end{proposition}
\begin{proof}
	We examine two cases.
	Assume that $V = \elswap[n]{k_1, k_2}$ is an elementary permutation with distinct indices $k_1, k_2 \in \set{1, \ldots, n}$,
	Note that if $k_1 = k_2$, then $V = I_n$ and $T$ is trivially a homogeneous basis.
	For $i \in \set{1, \ldots, n}$ with $i \neq k_1$ and $i \neq k_2$, $\tau_i = \sigma_i$.
	Since $\col_{k_1}(V) = \col_{k_2}(I_n)$, $\tau_{k_1} = \sigma_{k_2}$ is homogeneous.
	Similarly, $\col_{k_2}(V) = \col_{k_1}(I_n)$ and $\tau_{k_2} = \sigma_{k_1}$ is homogeneous.
	Therefore, $T$ is a homogeneous basis.
	
	Assume that $V = \eldilate[n]{k, \mu}$ is an elementary dilation over $\field[x]$ with index $k \in \set{1, \ldots, n}$ and $\mu \in (\field[x])^\times = \field[x] \setminus \set{0}$.
	For $i \in \set{1, \ldots, n}$ with $i \neq k$, $\tau_i = \sigma_i$.
	Since $\col_k(V) = \mu \cdot \col_k(I_n)$ and $\mu$ is a nonzero homogeneous element of degree $0$,
		$\tau_k = \mu \sigma_k$ and $\degh(\tau_k) = \degh(\mu \sigma_k) = \degh(\sigma_k)$.
	Therefore, $T$ is a homogeneous basis.
\end{proof}

Then, multiplication on the left (for row reduction) and on the right (for column reduction) by an elementary dilation or elementary permutation over $\field[x]$ will preserve the homogeneity of the initial bases and \fref{lemma:matrices-of-graded-homs} applies to the resulting matrix. 

\spacer 

Putting it all together, 
	any reduction algorithm on a matrix over $\field[x]$ that satisfy \fref{lemma:matrices-of-graded-homs} that involves a finite product of elementary permutations, elementary transvections, and elementary transvections that satisfy \fref{prop:column-reduction-preserves-hom} or \fref{prop:row-reduction-preserves-hom} over $\field[x]$
	will produce a matrix that has zero or homogeneous elements and the bases they produce will also be homogeneous (assuming we start with homogeneous bases).\clearpage
%

\newcounter{ALGSTEPS}
\renewcommand{\theALGSTEPS}{\textbf{Step\ \Alph{ALGSTEPS}.}}
\newcommand{\NEWSTEP}[2]{\spacer\refstepcounter{ALGSTEPS}\noindent{\bfseries Step \Alph{ALGSTEPS}.\ #2}\label{#1}\par\noindent\ignorespaces}
\newcommand{\subDim}{\ensuremath [k\,\mathord{:}\,m, k\,\mathord{:}\,n]} 

%

\section{An Ungraded SND Algorithm in the Graded Case}
\label{section:snd-algorithm-in-the-graded-case}

In this section, 
	we present and discuss the algorithm for finding SNDs given \cite[Remark 5.3.4]{algebra:adkins} optimized for the case of matrices of graded homomorphisms relative to homogeneous bases.

\begin{algorithmenv}\label{alg:general-snf-graded-homs}
	\textbf{Matrix Reduction Algorithm for Graded SNDs}.
	\vspace{-3pt}\begin{spacing}{1.2} 
		\setlength\tabcolsep{0pt}
		\begin{tabularx}{\linewidth}{l !{\quad:\quad} X}
			\textbf{Input } & $[\phi] \in \M_{m,n}(\field[x])$ as given in the hypothesis of \fref{lemma:matrices-of-graded-homs}, i.e.\ $\phi$ is a graded $\field[x]$-module homomorphism relative to homogeneous bases.
			\\ 
			\textbf{Output} & A Smith Normal Decomposition $(U,D,V)$ of $[\phi]$
			with $U \in \GL(m, \field[x])$, $V \in \GL(n, \field[x])$ and $D \in \M_{m,n}(\field[x])$.
		\end{tabularx}
	\begin{algorithm}[H]
		\textbf{set} $D_0 := [\phi]$. \;
		\For{increasing indices $k := 1, \ldots, \min(m,n)$}{
			\If{$D_{k-1}(j,i) = 0$ for all $j \in \set{k, \ldots, m}$ and $i \in \set{k, \ldots, n}$}{
				\textbf{set} $D_{\min(m,n)} := D_{k-1}$. \;
				\lFor{indices $i := k, \ldots, \min(m,n)$}{%
					\textbf{set} $U_i := I_m$ and $V_i := I_n$.
				}
				\textbf{break.} \tcc*[f]{stop the for-loop, $D_{\min(m,n)}$ is in SNF.}
			}
			\textbf{find} a row index $r_k \in \set{k, \ldots, m}$ and a column index $c_k \in \set{k, \ldots, m}$ such that \newline%
				\null\hspace{5pt}$\displaystyle
					\degh\!\big(D_{k-1}(r_k, c_k)\big) 
					= \min\!\Big\{
						\degh\!\big( D_{k-1}(j,i) \big)
						\,:\, 
						j \in \set{k, \ldots, m},
						i \in \set{k, \ldots, n},
						D_{k-1}(j,i) \neq 0
					\,\Big\}
				$. \;
			\textbf{set} $U_{k,k} := \elswap[m]{k, r_k}$, $V_{k,k} := \elswap[n]{k, c_k}$ 
				and $W_k := U_{k,k} D_{k-1} V_{k,k}$. 
				\label{line:adkins-alg-permutations}
				\;
			\For{increasing column index $i := {k+1, \ldots, n}$}{
				\textbf{set} $f_{k,i} := - W_{k}(k,i) \bigmod W_{k}(k,k)$
				and $V_{k,i} := \eladd[n]{k, i \,; f_{k,i}}$.
			}
			\textbf{set} $V_k := V_{k,k} V_{k, k+1} \cdots V_{k,n}$. \;
			\For{increasing row index $j := {k+1, \ldots, m}$}{
				\textbf{set} $g_{k,j} := - W_{k}(j,k) \bigmod W_{k}(k,k)$
				and $U_{k,j} := \eladd[m]{j,k \,; -g_{k,j}}$.
				\tcc*[r]{ $(U_{k,j})\inv = \eladd[m]{j,k \,; f_{k,k}}$}
			}
			\textbf{set} $U_k := U_{k,k} U_{k,k+1} \cdots U_{k,m}$. 
				\tcc*[r]{$(U_k)\inv = (U_{k,m})\inv (U_{k,m-1})\inv \cdots (U_{k,k})\inv$} 
			\textbf{set} $D_k := (U_k)\inv D_{k-1} V_k$. \;
		}
		\textbf{set} $U := U_{1} U_2 \cdots U_{\min(m,n)}$, $V := V_1 V_2 \cdots V_{\min(m,n)}$. \;
		\textbf{set} $D := D_{\min(m,n)}$. \;
		\Return{$(U,D,V)$}
	\end{algorithm}
	\end{spacing}
\end{algorithmenv}
\remarks{
	\item 
	This algorithm is written such that no variable is re-defined, 
	e.g.\ for $[\phi] \in \M_{3,4}(\field[x])$, the matrix $D_2$ is defined once in Line 15 at the $k=2$ loop and not re-defined in later loops.
	This is done to help with exposition. 

	\item 
	We refer to \cite[Remark 5.3.4]{algebra:adkins} for the correctness of \fref{alg:general-snf-graded-homs}, i.e.\ that the result $(U,D,V)$ is indeed an SND of $[\phi]$ (as denoted above).
}

For the rest of this section, 
	we discuss how the assumption of $[\phi]$ being given as in \fref{lemma:matrices-of-graded-homs} allows certain optimizations in the algorithm 
	and how the matrix operations in the algorithm preserve the homogeneity of the initial matrices.
In particular, we examine 4 major steps that are done for each $k$\th iteration of the for-loop starting in Line 2 of \fref{alg:general-snf-graded-homs}:
\begin{enumerate}[left=0.2in, align=left]
	\item[\ref{step:one}] Check if the matrix $D_k$ is in Smith Normal Form. If not, continue.
	\item[\ref{step:two}] Choose an appropriate pivot and perform the appropriate row and column permutations.
	\item[\ref{step:column}] Eliminate the entries to the right of the pivot by column reduction.
	\item[\ref{step:row}] Eliminate the entries below the pivot by row reduction.
\end{enumerate}
We expand on what happens for each of these steps below.

\NEWSTEP{step:one}{Checking the Result of the Previous Iteration: Lines 3 to 7 of \fref{alg:general-snf-graded-homs}. }
	The algorithm in \cite[Remark 5.3.4]{algebra:adkins}
	examines the submatrix $D_{k-1}\subDim$ of $D_{k-1}$, 
		where $D_{k-1}$ is the result of the previous iteration, 
	and $D_{k-1}\subDim$ is 
	obtained by removing columns $1, \ldots, k-1$ and rows $1, \ldots, k-1$ of $D_{k-1}$, i.e.\ $D_{k-1}$ can be expressed as the following block matrix:
	\begin{equation*}
		D_{k-1} = \begin{pmatrix}
			\diag(d_1, \ldots, d_{k-1}) & 0 \\
			0 & D_{k-1}\subDim
		\end{pmatrix}
	\end{equation*}
	where $d_1, \ldots, d_{k-1}$ are to be the nonzero diagonal elements of the SNF $D$ of $[\phi]$.
	There are two cases:
	\begin{enumerate}[left=0.7in]
		\item[Case 1:] 
		If the submatrix $D_{k-1}\subDim$ is the zero matrix,
		then $D_{k-1}$ is already in Smith Normal Form and we do not need to do any more matrix operations.
		That is, $D_{k-1}$ is the Smith Normal Form of $[\phi]$.
		Note that we set $D_{\min(m,n)} := D_{k-1}$ in Line 4, as opposed to defining $D := D_{k-1}$ right there, for consistency.
		The \textbf{break} keyword in Line 6 tells us to stop the for-loop and proceed directly to Line 16.
		Then, at Line 17, we define $D := D_{\min(m,n)}$.

		\item[Case 2:]
		If the submatrix $D_{k-1}\subDim$ is not the zero matrix, 
		then $D_{k-1}$ is not in Smith Normal Form.
		The goal of the $k$\th iteration is then to find matrices $U_k \in \GL(m, \field[x])$ and $V_k \in \GL(n, \field[x])$ such that 
		\begin{equation*}
			(U_k)\inv D_{k-1} V_k = D_k = \begin{pmatrix}
				\diag(d_1, \ldots, d_{k-1}, d_k) & 0 \\
				0 & D_{k}[k{+}1\!:\!m, k{+}1\!:\!n]
			\end{pmatrix}
		\end{equation*}
		where $D_k$ is the end result of this iteration and passed to the next,
		i.e.\ the matrix $D_k$ determines the $k$\th diagonal element of the Smith Normal Form of $D$.
		We apply the appropriate elementary matrices on $D_{k-1}$ to achieve this.
	\end{enumerate}
	Below, we provide an example of a matrix on which \fref{alg:general-snf-graded-homs} stops at some $k$\th iteration and another where it does not and continues.

\begin{example}
	Let $A \in \M_{4,4}(\field[x])$ 
	and $B \in \M_{4,4}(\field[x])$ be given below:
	\begin{equation*}
		A = \begin{pmatrix}
			1 & 0 & 0 & 0 \\
			0 & 5x & 2 & 0 \\
			0 & x^2 & 3x & 3 \\
			0 & 4 & -2x & 2x^2
		\end{pmatrix}
		\qquad\text{ and }\qquad
		B = \begin{pmatrix}
			1 & 0 & 0 & 0 \\
			0 & x^2 & 0 & 0 \\
			0 & 0 & 0 & 0 \\
			0 & 0 & 0 & 0
		\end{pmatrix}
	\end{equation*}
	We can view the matrix $A$ to be the matrix $D_1$ relative to \fref{alg:general-snf-graded-homs}, i.e.\ this is the matrix that is processed in the iteration with $k=2$.
	Then, $A$ is not in Smith Normal Form and we consider the submatrix 
	$A[2\mathord{:}4, 2\mathord{:}4]$, highlighted in \redtag below:
	\begin{equation*}
		A = \begin{pNiceMatrix}
			\CodeBefore [create-cell-nodes]
				\tikz \node [red-cell = (2-2) (4-4)] {} ;
			\Body
			1 & 0 & 0 & 0 \\
			0 & 5x & 2 & 0 \\
			0 & x^2 & 3x & 3 \\
			0 & 4 & -2x & 2x^2
		\end{pNiceMatrix} 
		\qquad\text{ with }\qquad 
		A[2\mathord{:}4, 2\mathord{:}4] = \begin{pmatrix}
			5x & 2 & 0 \\
			x^2 & 3x & 3 \\
			4 & -2x & 2x^2
		\end{pmatrix}
	\end{equation*}
	Since $A[2\mathord{:}4, 2\mathord{:}4]$ is not the zero matrix, we proceed with the calculation.

	In contrast, we can assume that the matrix $B$ corresponds to the matrix $D_2$ relative to \fref{alg:general-snf-graded-homs}, i.e.\ the matrix processed in the iteration with $k=3$.
	Then, \fref{alg:general-snf-graded-homs} concludes that $B$ is in Smith Normal Form since the submatrix $B[3\mathord{:}4, 3\mathord{:}4]$, highlighted in \redtag below, is the zero matrix.
	\begin{equation*}
		B = D_2 = \begin{pNiceMatrix}
			\CodeBefore [create-cell-nodes]
				\tikz \node [red-cell = (3-3) (4-4)] {} ;
			\Body
			1 & 0 & 0 & 0 \\
			0 & x^2 & 0 & 0 \\
			0 & 0 & 0 & 0 \\
			0 & 0 & 0 & 0
		\end{pNiceMatrix}
	\end{equation*}
	Then, the for-loop stops here and sets $D_{\min(m,n)} = D_2 = B$ in Line 4 and $D = B$ in Line 17 of \fref{alg:general-snf-graded-homs}. 
\end{example}

\NEWSTEP{step:two}{On Selection of Pivot Elements: Lines 7-8 of \fref{alg:general-snf-graded-homs} }
	The next step involves finding an element in the submatrix $D_{k-1}\subDim$ that can be used to eliminate any nonzero entry in $D_{k-1}\subDim$ as necessary.
	For convenience, we call this element the \textit{pivot}, 
		consistent with the terminology for our reduction operations by \fref{defn:column-reduction-operation} and \fref{defn:row-reduction-operation}.

		This pivot, given in Line 7 of \fref{alg:general-snf-graded-homs} by $D_{k-1}(r_k, c_k)$,
			is an element that has minimal degree across the nonzero entries of $D_{k-1}$,
			where the term \textit{degree} refers to the Euclidean function $\deg(-)$ on a Euclidean domain.
		Since the pivot is chosen to be of minimal degree, 
			it can be shown that any nonzero element of $D_{k-1}\subDim$ can be eliminated by a finite number of elementary row/column operations.
			Note that if the matrix $[\phi]$ is over a field $\field$, then any nonzero element can serve as the pivot since $\deg(\mu) = 1$ for all nonzero $\mu \in \field$. 

	Observe that, by \fref{lemma:matrices-of-graded-homs}, the nonzero entries of $[\phi]$ are homogeneous and are of the form $kx^t$ for some $k \in \field$ and $t \in \nonnegints$.
	Then, the degree $\deg(-)$ on $\field[x]$ as a Euclidean domain and the degree $\degh(-)$ on $\field[x]$ as a graded ring agree on the entries of the matrices involved.
	Moreover, if the chosen pivot is $k x^s \in \field[x]$ for some $k \in \field$ and $s \in \nonnegints$,
	the pivot multiplier required to eliminate a nonzero entry $k' x^t$ can immeidately be determined by $k' x^t \bigmod k x^s = (k' \bigmod k)x^{t-s}$ since $t \leq s$ by minimality of the degree of the chosen pivot. 
	This is reflected in the Step C involving column reduction and Step D involving row reduction. 

	We then apply appropriate permutations on $D_{k-1}$ such that the pivot $D_{k-1}(r_k, c_k)$ is found on the $(k,k)$\th entry.
	Relative to the submatrix $D_{k-1}\subDim$, we want the pivot to be the $(1,1)$\th entry post-permutation.
	The row permutation is given by $U_{k,k} := \elswap[m]{k, r_k}$, which swaps rows $k$ and $r_k$ of $D_{k-1}$ in the product $(U_{k,k})\inv D_{k-1}$.
	Note that $(\elswap[m]{k,r_k})\inv = \elswap[m]{k,r_k}$.
	The column permutation is given by $V_{k,k} := \elswap[n]{k, c_k}$, which swaps columns $k$ and $c_k$ of $U_{k,k}D_{k-1}$ in the product $U_{k,k}D_{k-1}V_{k,k}$.
	We save this permuted matrix as $W_k$ and perform the matrix operations on $W_k$.
	
	Note that, by \fref{prop:perm-and-dilation-preserve-hom} (permutations preserve homogeneity), 
		the homogeneity of the entries of $W_k$ 
		and of the bases induced by the multiplication of the permutation matrices $U_{k,k} = \elswap[m]{k, r_k}$ and $V_{k,k} = \elswap[n]{k, c_k}$ are preserved.
	We provide an example of this calculation below.

\begin{example}\label{ex:zom1-gradedboundary1-permutation}
	Let $[\boundary_1\graded] \in \M_{4,5}(\field[x])$ be as given in \fref{ex:column-reduction-of-graded-boundary-one}.
	For convenience, the description of $[\boundary_1\graded]$ is provided below.
	We apply \fref{alg:general-snf-graded-homs} on $[\boundary_1\graded]$.
	Following the notation in the algorithm,
		we have that $D_0 := [\boundary_1\graded]$.
	At the $k=1$ iteration, the entry $D_0(3,2)$, highlighted in \bluetag below, can serve as the pivot.
	\begin{equation*}
		D_0 := 
		[\boundary_1\graded] 
		= \begin{NiceArray}{>{\color{gray}}c ccccc}
			\CodeBefore [create-cell-nodes]
				\tikz \node [blue-cell = (4-3)] {} ;
			\Body
			\RowStyle[color=gray]{}
			& abx & bcx & adx^2 & cdx^2 & acx^3 \\
			{a}		&	-x & 0 & -x^2 & 0 & -x^3 \\
			{b}		&	x & -x & 0 & 0 & 0 \\ 
			{cx}	&	0 & 1 & 0 & x & x^2 \\
			{dx}	&	0 & 0 & x & -x & 0
		\CodeAfter \SubMatrix({2-2}{5-6})
		\end{NiceArray}
	\end{equation*}
	Then, we permute the rows of $D_0$ by $U_{1,1} = \elswap[4]{1,r_1}$ with $r_1 := 3$ and columns of $D_0$ by $V_{1,1} = \elswap[5]{1, c_1}$ with $c_1 := 2$ as follows.
	Highlighted in \greentag are the rows and columns of the \greentagged{pivot} $D_0(r_1, c_1) = 1$ and in \orangetag are those of the \orangetagged{entry} $D_0(1,1) = -x$.
	{\arraycolsep=3pt\begin{longtable}{L !{$=$} L}
		U_{1,1} D_0 V_{1,1} 
		= \elswap[4]{1,3} D_0 \elswap[5]{1,2} 
		&\begin{pNiceMatrix}
			\CodeBefore [create-cell-nodes]
				\tikz \node [greencell = (1-3)] {} ;
				\tikz \node [orangecell = (3-1)] {} ;
			\Body
			0 & 0 & 1 & 0 \\
			0 & 1 & 0 & 0 \\
			1 & 0 & 0 & 0 \\
			0 & 0 & 0 & 1
		\end{pNiceMatrix}
		\begin{pNiceMatrix}
			\CodeBefore [create-cell-nodes]
				\tikz \node [orangestrip = (1-1) (1-5)] {} ;
				\tikz \node [orangestrip = (1-1) (4-1)] {} ;
				\tikz \node [greencell = (3-1) (3-5)] {} ;
				\tikz \node [greencell = (1-2) (4-2)] {} ;
				\tikz \node [pivotcell = (3-2)] {} ;
			\Body
			-x & 0 & -x^2 & 0 & -x^3 \\
			x & -x & 0 & 0 & 0 \\
			0 & 1 & 0 & x & x^2 \\
			0 & 0 & x & -x & 0
		\end{pNiceMatrix} 
		\begin{pNiceMatrix}
			\CodeBefore [create-cell-nodes]
				\tikz \node [greencell = (1-2)] {} ;
				\tikz \node [orangecell = (2-1)] {} ;
			\Body
			0 & 1 & 0 & 0 & 0 \\
			1 & 0 & 0 & 0 & 0 \\
			0 & 0 & 1 & 0 & 0 \\
			0 & 0 & 0 & 1 & 0 \\
			0 & 0 & 0 & 0 & 1
		\end{pNiceMatrix}
		\\ 
		&\begin{pNiceMatrix}
			\CodeBefore [create-cell-nodes]
				\tikz \node [orangestrip = (3-1) (3-5)] {} ;
				\tikz \node [orangestrip = (1-2) (4-2)] {} ;
				\tikz \node [greencell = (1-1) (1-5)] {} ;
				\tikz \node [greencell = (1-1) (4-1)] {} ;
				\tikz \node [pivotcell = (1-1)] {} ;
			\Body
			1	& 0	& 0	& x	& x^2 \\
			-x	& x	& 0	& 0	& 0 \\
			0	& -x &	-x^2 &	0 & -x^3 \\
			0	& 0	& x	& -x & 0
		\end{pNiceMatrix} 
		=: W_1
	\end{longtable}} 
	\noindent 
	with $W_1 \in \M_{4,5}(\field[x])$ as denoted in \fref{alg:general-snf-graded-homs}. 
\end{example}

\NEWSTEP{step:column}{On Elimination by Column Reduction: Lines 9-11 of \fref{alg:general-snf-graded-homs}}
The matrices $V_{k,k+1}, V_{k, k+2}, \ldots, V_{k, n} \in \GL(n, \field[x])$ correspond to column operations that eliminate the entries to the left of the pivot $W_k(k,k)$.
Recall that 
	any nonzero homogeneous element of $\field[x]$ divides any nonzero homogeneous element of equal or greater degree.
Since the pivot $W_k(k,k)$ is of minimal degree by construction, 
	we can eliminate any nonzero entry $W_k(k,i)$ to the left of $W_k(k,k)$ with $i \in \set{k+1, \ldots, n}$ by addition of $f_{k,i} \cdot W_k(k,k)$ with $f_{k,i} \in \field[x]$ homogeneous.

Note that since none of the transvections $V_{k,i} := \eladd[n]{k,i \,; f_{k,i}}$ act on the $k$\th columns of $W_k$ and all of them add a $f_{k,i}$-multiple of the $k$\th column of $W_k$,
	the order in which $V_{k,i}$ is multiplied to $W_k$ does not matter for the definition of $V_{k,i}$.
	In particular, $\col_k(W_k) = \col_k(W_k V_{k,i})$ for any $i \in \set{k+1, \ldots, n}$.
	Therefore, we can define $f_{k,i}$ for each transvection $V_{k,i}$ relative to $W_k$, as done in Line 10, as opposed to a matrix product like $W_k V_{k,k+1} \cdots V_{k,i-1}$ that accounts for matrix operations that are already done.

Observe that if $W_k(k,i) = 0$, then there is nothing to eliminate.
In that case, $f_{k,i} = 0$, $\eladd[n]{k,i \,; 0} = I_n$, and homogeneity is trivially preserved.
However, if $W_k(k,i) \neq 0$, then $f_{k,i}$ must be homogeneous with 
\begin{equation*}
	\degh\bigl( f_{k,i} \bigr) 
		+ \degh\bigl( W_k(k,k) \bigr)
		= \degh\bigl( W_{k}(k,i) \bigr)
	\,.
\end{equation*}
Then, \fref{prop:column-reduction-preserves-hom} (column reduction preserves homogeneity) applies and these column reduction operations preserve the homogeneity of the bases and of the entries.
We provide an example below.

\begin{example}\label{ex:zom1-gradedboundary1-columnred}
	Continue \fref{ex:zom1-gradedboundary1-permutation}. 
	For reference, a description of $W_1 \in \M_{4,5}(\rationals[x])$ is given below, with the \bluetagged{pivot entry $W_1(1,1)$} highlighted in \bluetag 
	and the \redtagged{target elements}, i.e.\ the elements to the left of the pivot, highlighted in \redtag\!.
	\begin{equation*}
		W_1 = 
		\begin{pNiceMatrix}
			\CodeBefore [create-cell-nodes]
				\tikz \node [blue-cell = (1-1)] {} ;
				\tikz \node [red-cell = (1-2)] {} ;
				\tikz \node [red-cell = (1-3)] {} ;
				\tikz \node [red-cell = (1-4)] {} ;
				\tikz \node [red-cell = (1-5)] {} ;
			\Body
			1	& 0	& 0	& x	& x^2 \\
			-x	& x	& 0	& 0	& 0 \\
			0	& -x &	-x^2 &	0 & -x^3 \\
			0	& 0	& x	& -x & 0
		\end{pNiceMatrix}
	\end{equation*} 
	Below, we define the elementary transvections $\eladd[5]{1,i \,; f_{1,i}}$ for $i \in \set{2, \ldots, 5}$
	used to eliminate the entry $W_1(1,i)$ with the \orangetagged{pivot multiplier $f_{1,i}$} highlighted in \orangetag\!.
	\begin{longtable}{
		L !{:} >{\displaystyle}L !{ and } 
		L !{=} >{\def\arraystretch{0.9}}L 
	} 
		i = 2 
			& \orangemath{f_{1,2} = -\frac{W_1(1,2)}{W_1(1,1)} = 0} 
			& V_{1,2} := \eladd[5]{1,2 \,; \orangemath{0}} = I_5 
			& \begin{pNiceMatrix}
				\CodeBefore [create-cell-nodes]
					\tikz \node [orange-cell = (1-2)] {} ;
				\Body
				1 & 0 & 0 & 0 & 0 \\
				0 & 1 & 0 & 0 & 0 \\
				0 & 0 & 1 & 0 & 0 \\
				0 & 0 & 0 & 1 & 0 \\
				0 & 0 & 0 & 0 & 1 \\
			\end{pNiceMatrix}
		\\[30pt] 
		i=3
			& \orangemath{f_{1,3} = -\frac{W_1(1,3)}{W_1(1,1)} = 0}
			& V_{1,3} := \eladd[5]{1,3 \,; \orangemath{0}} = I_5 
			& \begin{pNiceMatrix}
				\CodeBefore [create-cell-nodes]
					\tikz \node [orange-cell = (1-3)] {} ;
				\Body
				1 & 0 & 0 & 0 & 0 \\
				0 & 1 & 0 & 0 & 0 \\
				0 & 0 & 1 & 0 & 0 \\
				0 & 0 & 0 & 1 & 0 \\
				0 & 0 & 0 & 0 & 1 \\
			\end{pNiceMatrix}
		\\[30pt] 
		i=4 
			& \orangemath{f_{1,4} = -\frac{W_1(1,4)}{W_1(1,1)} = -x} 
			& V_{1,4} := \eladd[5]{1,3 \,; \orangemath{-x}} 
			& \begin{pNiceMatrix}
				\CodeBefore [create-cell-nodes]
					\tikz \node [orange-cell = (1-4)] {} ;
				\Body
				1 & 0 & 0 & -x & 0 \\
				0 & 1 & 0 & 0 & 0 \\
				0 & 0 & 1 & 0 & 0 \\
				0 & 0 & 0 & 1 & 0 \\
				0 & 0 & 0 & 0 & 1 \\
			\end{pNiceMatrix} 
		\\[30pt] 
		i=5 
			& \orangemath{f_{1,5} = -\frac{W_1(1,5)}{W_1(1,1)} = -x^2} 
			& V_{1,4} := \eladd[5]{1,3 \,; \orangemath{-x^2}} 
			& \begin{pNiceMatrix}
				\CodeBefore [create-cell-nodes]
					\tikz \node [orange-cell = (1-5)] {} ;
				\Body
				1 & 0 & 0 & 0 & -x^2 \\
				0 & 1 & 0 & 0 & 0 \\
				0 & 0 & 1 & 0 & 0 \\
				0 & 0 & 0 & 1 & 0 \\
				0 & 0 & 0 & 0 & 1 \\
			\end{pNiceMatrix} 
	\end{longtable}
	\noindent 
	As a sidenote, we can calculate the matrix product 
		$(V_{1,2}) (V_{1,3}) (V_{1,4}) (V_{1,5})$ 
	or any re-ordering of the product without explicitly identifying each $\eladd[5]{1,i \,; f_{1,i}}$ as follows:
	\begin{equation*}\def\arraystretch{0.9}
		(V_{1,2}) (V_{1,3}) (V_{1,4}) (V_{1,5})
		=
		\begin{pNiceMatrix}
			\CodeBefore [create-cell-nodes]
				\tikz \node [orange-highlight = (1-2)] {} ;
				\tikz \node [orange-highlight = (1-3)] {} ;
				\tikz \node [orange-highlight = (1-4)] {} ;
				\tikz \node [orange-highlight = (1-5)] {} ;
			\Body
			1 & f_{1,2} & f_{1,3} & f_{1,4} & f_{1,5} \\
			0 & 1 & 0 & 0 & 0 \\
			0 & 0 & 1 & 0 & 0 \\
			0 & 0 & 0 & 1 & 0 \\
			0 & 0 & 0 & 0 & 1 \\
		\end{pNiceMatrix} 
		=
		\begin{pNiceMatrix}
			\CodeBefore [create-cell-nodes]
				\tikz \node [orange-cell = (1-2)] {} ;
				\tikz \node [orange-cell = (1-3)] {} ;
				\tikz \node [orange-cell = (1-4)] {} ;
				\tikz \node [orange-cell = (1-5)] {} ;
			\Body
			1 & 0 & 0 & -x & -x^2 \\
			0 & 1 & 0 & 0 & 0 \\
			0 & 0 & 1 & 0 & 0 \\
			0 & 0 & 0 & 1 & 0 \\
			0 & 0 & 0 & 0 & 1 \\
		\end{pNiceMatrix} 
	\end{equation*} 
	Then, the column reduction step in the $k=1$ loop produces the following matrix, with the columns highlighted in \orangetag being the columns affected by the column reduction.
	\begin{longtable}{
		C !{=} L
	}\def\arraystretch{0.9}
		(W_1) (V_{1,2}) \cdots (V_{1,5})
		& {\arraycolsep=3pt\begin{pNiceMatrix}
			\CodeBefore [create-cell-nodes]
				\tikz \node [blue-cell = (1-1)] {} ;
				\tikz \node [red-cell = (1-2)] {} ;
				\tikz \node [red-cell = (1-3)] {} ;
				\tikz \node [red-cell = (1-4)] {} ;
				\tikz \node [red-cell = (1-5)] {} ;
			\Body
			1	& 0	& 0	& x	& x^2 \\
			-x	& x	& 0	& 0	& 0 \\
			0	& -x &	-x^2 &	0 & -x^3 \\
			0	& 0	& x	& -x & 0
		\end{pNiceMatrix}}
		\begin{pNiceMatrix}
			\CodeBefore [create-cell-nodes]
				\tikz \node [orange-cell = (1-2)] {} ;
				\tikz \node [orange-cell = (1-3)] {} ;
				\tikz \node [orange-cell = (1-4)] {} ;
				\tikz \node [orange-cell = (1-5)] {} ;
			\Body
			1 & 0 & 0 & -x & -x^2 \\
			0 & 1 & 0 & 0 & 0 \\
			0 & 0 & 1 & 0 & 0 \\
			0 & 0 & 0 & 1 & 0 \\
			0 & 0 & 0 & 0 & 1 \\
		\end{pNiceMatrix} 
		= {\arraycolsep=3pt\begin{pNiceMatrix}
			\CodeBefore [create-cell-nodes]
				\tikz \node [blue-cell = (1-1)] {} ;
				\tikz \node [orange-highlight = (1-2) (4-2)] {} ;
				\tikz \node [orange-highlight = (1-3) (4-3)] {} ;
				\tikz \node [orange-highlight = (1-4) (4-4)] {} ;
				\tikz \node [orange-highlight = (1-5) (4-5)] {} ;
			\Body
			1 & 0 & 0 & 0 & 0	\\
			-x & x & 0 & x^2 & x^3	\\
			0 & -x & -x^2 & 0 & -x^3	\\
			0 & 0 & x & -x & 0
		\end{pNiceMatrix}}
	\end{longtable} 
	\noindent 
	Note that the definition of $V_1$ by Line 11 of \fref{alg:general-snf-graded-homs} includes the column permutation in its definition. Note that since $V_{1,1} = \elswap[5]{1,2}$ is multiplied on the left of $(V_{1,2})\cdots(V_{1,5})$, it permutes the rows of $(V_{1,2})\cdots(V_{1,5})$ instead.
	Highlighted below are \orangetagged{Row $1$} and \greentagged{Row $2$} of $(V_{1,2})\cdots(V_{1,5})$ and where they are mapped to after the row permutation.
	\begin{longtable}{ 
		C
	}\def\arraystretch{0.9}
		V_1 = (V_{1,1})(V_{1,2})(V_{1,3})(V_{1,4})(V_{1,5})
		= \elswap[5]{1,2}
		\begin{pNiceMatrix}
				\CodeBefore [create-cell-nodes]
					\tikz \node [green-highlight = (2-1) (2-5)] {} ;
					\tikz \node [orange-highlight = (1-1) (1-5)] {} ;
				\Body
				1 & 0 & 0 & -x & -x^2 \\
				0 & 1 & 0 & 0 & 0 \\
				0 & 0 & 1 & 0 & 0 \\
				0 & 0 & 0 & 1 & 0 \\
				0 & 0 & 0 & 0 & 1
			\end{pNiceMatrix} 
		= \begin{NiceArray}{>{\color{gray}}c(ccccc)}
			\CodeBefore [create-cell-nodes]
				\tikz \node [green-highlight = (1-2) (1-6)] {} ;
				\tikz \node [orange-highlight = (2-2) (2-6)] {} ;
			\Body
			{abx} & 0 & 1 & 0 & 0 & 0 \\
			{bcx} & 1 & 0 & 0 & -x & -x^2 \\
			{adx^2} & 0 & 0 & 1 & 0 & 0 \\
			{cdx^2} & 0 & 0 & 0 & 1 & 0 \\
			{acx^3} & 0 & 0 & 0 & 0 & 1 
		\end{NiceArray} 
	\end{longtable} 
	\noindent 
	The matrix $V_1 \in \GL(5,\rationals[x])$ induces a basis $T = (\tau_1(x), \ldots, \tau_5(x))$ of $C_1\graded(\filt{K}; \rationals)$ given as follows:
	\begin{equation*}\def\arraystretch{1.1}
		\begin{array}{c !{=} c !{=} l}
			\tau_1(x) & (1)(bcx) & (bc)x \\
			\tau_2(x) & (1)(abx) & (ab)x \\
			\tau_3(x) & (1)(adx^2) & (ad)x^2 \\
			\tau_4(x) & (-x)(bcx) + (1)(cdx^2) 
				& (-bc + cd)x^2 \\
			\tau_5(x) & (-x^2)(bcx) + (acx^3) 
				& (-bc + ac)x^3
		\end{array}
	\end{equation*}
	Observe that the basis $T = (\tau_i(x))$ is a homogeneous basis.
\end{example}

\NEWSTEP{step:row}{On Elimination by Row Reduction: Lines 12-14 of \fref{alg:general-snf-graded-homs}}
The matrices $U_{k,k+1}, U_{k,k+2}, \ldots, U_{k,m} \in \GL(m, \field[x])$ correspond to row operations that eliminate the entries under the pivot entry $W_k(k,k)$.
As with column reduction, 
	since the elements of $W_k$ are homogeneous by \fref{lemma:matrices-of-graded-homs}
	and the pivot $W_k(k,k)$ is chosen to be of minimal degree,
	each entry $W_k(j,k)$ under the pivot can be eliminated by the addition of $g_{k,j} W_k(k,k)$ with $g_{k,j} \in \field[x]$ homogeneous and $j \in \set{k+1, \ldots, m}$.

As mentioned on Item 4, the column operations in Lines 9-11 of \fref{alg:general-snf-graded-homs} do not affect the $k$\th column of $W_k$.
Then, the entries below the pivot $W_k(k,k)$, 
	i.e.\ entries indexed by $(j,k)$ with $j \in \set{k+1, \ldots, m}$, 
	are unaffected by the column reduction.
	Therefore, we can define the transvections $U_{k,j}$ that eliminate the $(j,k)$\th entry of $W_k$ relative to $W_k$, as opposed to the matrix $W_k V_k$ post-column reduction.
Furthermore, all of the row operations by the transvection $(U_{k,j})\inv := \eladd[m]{j,k \,; g_{k,j}}$
	also do not affect the $k$\th row of $W_k$.
Therefore, 
	we can define the elements $g_{k,j}$ such that $W_k(j,k) + g_{k,j} \cdot W_k(k,k) = 0$ relative to $W_k$, as done in Line 13 of \fref{alg:general-snf-graded-homs}.
As a sidenote, the order in which the row operations $\eladd[m]{j,k \,; g_{k,j}}$ with $j \in \set{k+1,\ldots,m}$ \textit{and} the column operations $\eladd[n]{k,i \,; f_{k,i}}$ with $i \in \set{k+1, \ldots, n}$ do not affect the definition of $g_{k,j}$ and $f_{k,i}$.

Since we have labeled the matrix $U_k \in \GL(m,\field)$ in Item (1) such that $(U_k)\inv D_{k-1} V_k = D_k$, i.e.\ $(U_k)\inv$ more closely describes the row operations done on $W_k = U_{k,k} D_{k-1} V_{k,k}$,
	we will have to take the inverse of the matrix product 
	$(U_{k,m})\inv (U_{k,m-1}) \cdots (U_{k,k})\inv$
	resulting from the row operations done on $D_{k-1}$
	if we want an expression for $U_k$.
To simplify this calculation, we can apply the property that $(AB)\inv = B\inv A\inv$ for any pair of invertible matrices $A,B \in \GL(m,\field[x])$.
Then, the matrix reduction on $D_{k-1}$ is as follows:

\begin{equation*}\arraycolsep=2pt
	\begin{NiceArray}{c !{\,=\,} c ccc c}
		D_k
		& (U_{k,m})\inv (U_{k,m-1})\inv \cdots (U_{k,k+1})\inv
		& (U_{k,k})\inv & (D_{k-1}) & (V_{k,k})
		& (V_{k,k+1}) (V_{k,k+2}) \cdots (V_{k,n})
	\\[4pt] 
		&\span\Big(\, 
			\underbrace{(U_{k,k})(U_{k+1}) \cdots (U_{k,m-1})(U_{k,m})}_{\raisebox{-5pt}{$U_k$}} 
		\,\Big)^{-1}
		& (D_{k-1})
		&\span
			\underbrace{(V_{k,k})(V_{k,k+1}) (V_{k,k+2}) \cdots (V_{k,n})}_{\raisebox{-5pt}{$V_k$}}
	\CodeAfter
		\OverBrace[shorten,yshift=3pt]{1-2}{1-2}{\text{\small
			row reduction on $W_k$
		}}
		\OverBrace[shorten,yshift=3pt]{1-3}{1-5}{W_k}
		\OverBrace[shorten,yshift=3pt]{1-6}{1-6}{\text{\small
			column reduction on $W_k$
		}}
	\end{NiceArray}
\end{equation*}

Observe that if $W_k(j,k) = 0$, then there is nothing to eliminate and $\eladd[m]{j,k \,; 0} = I_m$ trivially preserves homogeneity.
If $W_k(j,k) \neq 0$, $g_{k,j} \in \field[x]$ must be homogeneous and satisfy the following equation:
\begin{equation*}
	\degh\bigl( g_{k,j} \bigr) 
		+ \degh\bigl( W_k(k,k) \bigr)
		= \degh\bigl( W_k(j,k) \bigr)
\end{equation*}
Then, \fref{prop:row-reduction-preserves-hom} (row reduction preserves homogeneity) applies and these row operations preserve the homogeneity of the bases and of the entries.
We continue the calculation of the SNF of $[\boundary_1\graded]$ below.

\begin{example}\label{ex:zom1-gradedboundary1-rowred}
	Continue \fref{ex:zom1-gradedboundary1-columnred}.
	Listed below is the matrix $W_1$, 
		with the \bluetagged{pivot entry $W_1(1,1)$} highlighted in \bluetag
		and the \redtagged{target entries}, i.e.\ the entries below the chosen pivot, in \redtag\!.
	\begin{equation*}
		W_1 = 
		\begin{pNiceMatrix}
			\CodeBefore [create-cell-nodes]
				\tikz \node [blue-cell = (1-1)] {} ;
				\tikz \node [red-cell = (2-1)] {} ;
				\tikz \node [red-cell = (3-1)] {} ;
				\tikz \node [red-cell = (4-1)] {} ;
			\Body
			1	& 0	& 0	& x	& x^2 \\
			-x	& x	& 0	& 0	& 0 \\
			0	& -x &	-x^2 &	0 & -x^3 \\
			0	& 0	& x	& -x & 0
		\end{pNiceMatrix}
	\end{equation*} 
	Below, we define the elementary transvections $\eladd[4]{j,i \,; f_{1,i}}$ for $j \in \set{2, \ldots, 4}$ used to eliminate the entry $W_1(j,1)$,
	with the \orangetagged{pivot multiplier $g_{1,j}$} highlighted in \orangetag\!.
	\begin{longtable}{
		L !{:} >{\displaystyle}L !{ and } 
		L !{=} >{\def\arraystretch{0.9}}L
	} 
		j = 2 
			& \orangemath{g_{1,2} = -\frac{W_1(2,1)}{W_1(1,1)} = x}
			& (U_{1,2})\inv := \eladd[4]{2,1 \,; \orangemath{x\vphantom{0}}}
			& \begin{pNiceMatrix}
				\CodeBefore [create-cell-nodes]
					\tikz \node [orange-cell = (2-1)] {} ;
				\Body
				1 & 0 & 0 & 0 \\
				x & 1 & 0 & 0 \\
				0 & 0 & 1 & 0 \\
				0 & 0 & 0 & 1 \\
			\end{pNiceMatrix} 
		\\[15pt]
		j = 3 
			& \orangemath{g_{1,3} = -\frac{W_1(3,1)}{W_1(1,1)} = 0}
			& (U_{1,3})\inv := \eladd[4]{3,1 \,; \orangemath{0}} = I_4
			& \begin{pNiceMatrix}
				\CodeBefore [create-cell-nodes]
					\tikz \node [orange-cell = (3-1)] {} ;
				\Body
				1 & 0 & 0 & 0 \\
				0 & 1 & 0 & 0 \\
				0 & 0 & 1 & 0 \\
				0 & 0 & 0 & 1 \\
			\end{pNiceMatrix} 
		\\[15pt]
		j = 4 
			& \orangemath{g_{1,4} = -\frac{W_1(4,1)}{W_1(1,1)} = 0}
			& (U_{1,4})\inv := \eladd[4]{4,1 \,; \orangemath{0}} = I_4
			& \begin{pNiceMatrix}
				\CodeBefore [create-cell-nodes]
					\tikz \node [orange-cell = (4-1)] {} ;
				\Body
				1 & 0 & 0 & 0 \\
				0 & 1 & 0 & 0 \\
				0 & 0 & 1 & 0 \\
				0 & 0 & 0 & 1 \\
			\end{pNiceMatrix} 
	\end{longtable}
	\noindent 
	Like in the case of $(V_{1,2})\cdots(V_{1,5})$, 
	the matrix product 
	$(U_{1,4})\inv (U_{1,3})\inv (U_{1,2})\inv = (U_{1,2} U_{1,3} U_{1,2})\inv$ can be identified explicitly without identifying each matrix $U_{1,j}$ as follows:
	\begin{equation*}\def\arraystretch{0.9}
		(U_{1,4})\inv (U_{1,3})\inv (U_{1,2})\inv
		=
		\begin{pNiceMatrix}
			\CodeBefore [create-cell-nodes]
				\tikz \node [orange-highlight = (2-1)] {} ;
				\tikz \node [orange-highlight = (3-1)] {} ;
				\tikz \node [orange-highlight = (4-1)] {} ;
			\Body
			1 & 0 & 0 & 0 \\
			g_{1,2} & 1 & 0 & 0 \\
			g_{1,3} & 0 & 1 & 0 \\
			g_{1,4} & 0 & 0 & 1 \\
		\end{pNiceMatrix} 
		=
		\begin{pNiceMatrix}
			\CodeBefore [create-cell-nodes]
				\tikz \node [orange-highlight = (2-1)] {} ;
				\tikz \node [orange-highlight = (3-1)] {} ;
				\tikz \node [orange-highlight = (4-1)] {} ;
			\Body
			1 & 0 & 0 & 0 \\
			x & 1 & 0 & 0 \\
			0 & 0 & 1 & 0 \\
			0 & 0 & 0 & 1 \\
		\end{pNiceMatrix} 
	\end{equation*} 
	Since the inverses of elementary transvections are known,
		we can also calculate the product $(U_{1,2})(U_{1,3})(U_{1,4})$, i.e.\ the inverse of the matrix calculated above, directly as follows:
	\begin{equation*}\def\arraystretch{0.9}
		(U_{1,2})(U_{1,3})(U_{1,4})
		=
		\begin{pNiceMatrix}
			\CodeBefore [create-cell-nodes]
				\tikz \node [orange-highlight = (2-1)] {} ;
				\tikz \node [orange-highlight = (3-1)] {} ;
				\tikz \node [orange-highlight = (4-1)] {} ;
			\Body
			1 & 0 & 0 & 0 \\
			-g_{1,2} & 1 & 0 & 0 \\
			-g_{1,3} & 0 & 1 & 0 \\
			-g_{1,4} & 0 & 0 & 1 \\
		\end{pNiceMatrix} 
		=
		\begin{pNiceMatrix}
			\CodeBefore [create-cell-nodes]
				\tikz \node [orange-highlight = (2-1)] {} ;
				\tikz \node [orange-highlight = (3-1)] {} ;
				\tikz \node [orange-highlight = (4-1)] {} ;
			\Body
			1 & 0 & 0 & 0 \\
			-x & 1 & 0 & 0 \\
			0 & 0 & 1 & 0 \\
			0 & 0 & 0 & 1 \\
		\end{pNiceMatrix} 
	\end{equation*} 
	To calculate $U_1 := (U_{1,1})(U_{1,2})(U_{1,3})(U_{1,4})$ with 
		$U_{1,1} = \elswap[4]{1,3}$ from \fref{ex:zom1-gradedboundary1-rowred}, 
		we can permute the rows of the product $(U_{1,2})(U_{1,3})(U_{1,4})$ as follows:
	\begin{equation*}
		U_1 = (U_{1,1})(U_{1,2})(U_{1,3})(U_{1,4})
		= \elswap[4]{1,3}\begin{pNiceMatrix}
			\CodeBefore [create-cell-nodes]
				\tikz \node [green-highlight = (1-1) (1-4)] {} ;
				\tikz \node [orange-highlight = (3-1) (3-4)] {} ;
			\Body
			1 & 0 & 0 & 0 \\
			-x & 1 & 0 & 0 \\
			0 & 0 & 1 & 0 \\
			0 & 0 & 0 & 1 \\
		\end{pNiceMatrix} 
		= \begin{NiceArray}{>{\color{gray}}c(cccc)}
			\CodeBefore [create-cell-nodes]
				\tikz \node [orange-highlight = (1-2) (1-5)] {} ;
				\tikz \node [green-highlight = (3-2) (3-5)] {} ;
			\Body
			{a} 	& 0 & 0 & 1 & 0 \\
			{b}		& -x & 1 & 0 & 0 \\
			{cx}	& 1 & 0 & 0 & 0 \\
			{dx}	& 0 & 0 & 0 & 1 \\
		\end{NiceArray} 
	\end{equation*}
	The matrix $U_1 \in \GL(4,\rationals[x])$ 
	determines a basis $B = (\beta_1(x), \ldots, \beta_4(x))$ of $C_0\graded(\filt{K};\rationals[x])$
	by $[\beta_j(x)] = \col_j(U_1)$ as given below:
	\begin{equation*}\def\arraystretch{1.1}
		\begin{array}{c !{=} c !{=} l}
			\beta_1(x) & (-x)(b) + (1)(cx) & (-b+c)x \\
			\beta_2(x) & (1)(b) & b \\
			\beta_3(x) & (1)(a) & a \\
			\beta_4(x) & (1)(dx) & (d)x
		\end{array}
	\end{equation*}
	We can confirm by direct calculation that the \redtagged{entries} to the right and below the pivot $W_1(1,1)$ in the product $(U_1)\inv D_0 V_1 =: D_1$ are zero, highlighted in \redtag below:
	{\def\arraystretch{0.9}\begin{longtable}{L !{$=$} L }
		(U_1)\inv D_0 V_1 
		& \begin{pmatrix}
			0 & 0 & 1 & 0 \\
			x & 1 & 0 & 0 \\
			1 & 0 & 0 & 0 \\
			0 & 0 & 0 & 1 \\
		\end{pmatrix}
		\begin{pmatrix}
			1	& 0	& 0	& x	& x^2 \\
			-x	& x	& 0	& 0	& 0 \\
			0	& -x &	-x^2 &	0 & -x^3 \\
			0	& 0	& x	& -x & 0
		\end{pmatrix}
		\begin{pmatrix}
			0 & 1 & 0 & 0 & 0 \\
			1 & 0 & 0 & -x & -x^2 \\
			0 & 0 & 1 & 0 & 0 \\
			0 & 0 & 0 & 1 & 0 \\
			0 & 0 & 0 & 0 & 1 
		\end{pmatrix}
		\\ 
		& \begin{pNiceMatrix}
			\CodeBefore [create-cell-nodes]
				\tikz \node [red-highlight = (1-2) (1-5)] {} ;
				\tikz \node [red-highlight = (2-1) (4-1)] {} ;
				\tikz \node [blue-cell = (1-1)] {} ;
			\Body
			1 &	0 &	0 &	0 &	0 \\
			0 &	x &	0 &	x^2 & x^3 \\
			0 & -x & -x^2 & 0 & x^3 \\
			0 & 0 & x & -x & 0
		\end{pNiceMatrix} 
		=: D_1
	\end{longtable}}
\end{example}

\spacer\noindent
{\bfseries Final Step. Aggregating the Results: Lines 16-17 of \fref{alg:general-snf-graded-homs}}
\par\noindent 
Observe that in the $k$\th iteration, the zero entries found below and to the right of the entries $d_1, \ldots, d_{k-1}$ of $D_{k-1}$ with $d_i = D_{k-1}(i,i)$ remain zero on $D_k$, i.e.\ the permutation, column reduction, and row reduction operations done on $D_{k-1}$ leave columns $1, \ldots, k-1$ and rows $1, \ldots, k-1$ of $D_{k-1}$ undisturbed. 
Then, at each iteration with $k \in \set{1, \ldots, \min(m,n)}$, the matrix $[\phi]$ gets cleared row by row and column by column like so:

\begin{longtable}{L}
	\stackrel{\raisebox{5pt}{$D_0 := [\phi]$}}
	{\scalebox{0.8}{$\begin{NiceArray}{(cc c cc c c)}
		* & * & \cdots & * & * & \cdots & * 	\\
		* & * & \cdots & * & * & \cdots & * 	\\
		\vdots & \vdots & {\ddots} 
			& \vdots & \vdots & {} & \vdots \\
		* & * & \cdots & * & * & \cdots & *		\\
		* & * & \cdots & * & * & \cdots & *		\\
		\vdots & \vdots & {} 
			& \vdots & \vdots & {\ddots} & \vdots \\ 
		* & * & \cdots & * & * & \cdots & * 
	\end{NiceArray}$}}

	\,\scalebox{1.5}{$\rightsquigarrow$}\,
	\stackrel{\raisebox{5pt}{$D_1 := (U_1)\inv D_0 V_1$}}
	{\scalebox{0.8}{$\begin{NiceArray}{(cc c cc c c)}
		\CodeBefore [create-cell-nodes]
			\tikz \node [orangecell = (1-2) (1-7)] {} ;
			\tikz \node [orangecell = (2-1) (7-1)] {} ;
			\tikz \node [pivotcell = (1-1)] {} ;
		\Body
		d_1 & 0 & \cdots & 0 & 0 & \cdots & 0 	\\
		0 & * & \cdots & * & * & \cdots & * 	\\
		\vdots & \vdots & {\ddots} 
			& \vdots & \vdots & {} & \vdots \\
		0 & * & \cdots & * & * & \cdots & *		\\
		0 & * & \cdots & * & * & \cdots & *		\\
		\vdots & \vdots & {} 
			& \vdots & \vdots & {\ddots} & \vdots \\ 
		0 & * & \cdots & * & * & \cdots & * 
	\end{NiceArray}$}} 

	\,\scalebox{1.5}{$\rightsquigarrow$}\,
	\stackrel{\raisebox{5pt}{$D_2 := (U_2)\inv D_1 V_2$}}
	{\scalebox{0.8}{$\begin{NiceArray}{(cc c cc c c)}
		\CodeBefore [create-cell-nodes]
			\tikz \node [orangecell = (2-3) (2-7)] {} ;
			\tikz \node [orangecell = (3-2) (7-2)] {} ;
			\tikz \node [pivotcell = (1-1)] {} ;
			\tikz \node [pivotcell = (2-2)] {} ;
		\Body
		d_1 & 0 & \cdots & 0 & 0 & \cdots & 0 	\\
		0 & d_2 & \cdots & 0 & 0 & \cdots & 0 	\\
		\vdots & \vdots & {\ddots} 
			& \vdots & \vdots & {} & \vdots \\
		0 & 0 & \cdots & * & * & \cdots & *		\\
		0 & 0 & \cdots & * & * & \cdots & *		\\
		\vdots & \vdots & {} 
			& \vdots & \vdots & {\ddots} & \vdots \\ 
		0 & 0 & \cdots & * & * & \cdots & * 
	\end{NiceArray}$}}

	\,\scalebox{1.5}{$\rightsquigarrow \cdots$}\,

	\\[50pt] 

	\qquad\quad
	\,\scalebox{1.5}{$\cdots\, \rightsquigarrow$}\,

	\stackrel{\raisebox{5pt}{$D_{k-1} := (U_{k-1})\inv D_{k-2} V_{k-1}$}}
	{\scalebox{0.8}{$\begin{NiceArray}{(cc c cc c c)}
		\CodeBefore [create-cell-nodes]
			\tikz \node [orangecell = (4-5) (4-7)] {} ;
			\tikz \node [orangecell = (5-4) (7-4)] {} ;
			\tikz \node [pivotcell = (1-1)] {} ;
			\tikz \node [pivotcell = (2-2)] {} ;
			\tikz \node [pivotcell = (4-4)] {} ;
		\Body
		d_1 & 0 & \cdots & 0 & 0 & \cdots & 0 	\\
		0 & d_2 & \cdots & 0 & 0 & \cdots & 0 	\\
		\vdots & \vdots & {\color{NavyBlue}\ddots} 
			& \vdots & \vdots & {} & \vdots \\
		0 & 0 & \cdots & d_{k-1} & 0 & \cdots & 0		\\
		0 & 0 & \cdots & 0 & * & \cdots & *		\\
		\vdots & \vdots & {} 
			& \vdots & \vdots & {\ddots} & \vdots \\ 
		0 & 0 & \cdots & 0 & * & \cdots & * 
	\end{NiceArray}$}}

	\,\scalebox{1.5}{$\rightsquigarrow$}\,
	\stackrel{\raisebox{5pt}{$D_{k} := (U_{k})\inv D_{k-1} V_{k}$}}
	{\scalebox{0.8}{$\begin{NiceArray}{(cc c cc c c)}
		\CodeBefore [create-cell-nodes]
			\tikz \node [orangecell = (5-6) (5-7)] {} ;
			\tikz \node [orangecell = (6-5) (7-5)] {} ;
			\tikz \node [pivotcell = (1-1)] {} ;
			\tikz \node [pivotcell = (2-2)] {} ;
			\tikz \node [pivotcell = (4-4)] {} ;
			\tikz \node [pivotcell = (5-5)] {} ;
		\Body
		d_1 & 0 & \cdots & 0 & 0 & \cdots & 0 	\\
		0 & d_2 & \cdots & 0 & 0 & \cdots & 0 	\\
		\vdots & \vdots & {\color{NavyBlue}\ddots} 
			& \vdots & \vdots & {} & \vdots \\
		0 & 0 & \cdots & d_{k-1} & 0 & \cdots & 0		\\
		0 & 0 & \cdots & 0 & d_k & \cdots & *		\\
		\vdots & \vdots & {} 
			& \vdots & \vdots & {\ddots} & \vdots \\ 
		0 & 0 & \cdots & 0 & * & \cdots & * 
	\end{NiceArray}$}}

	\,\scalebox{1.5}{$\rightsquigarrow$}\,\,

	\Bigg(\parbox{1in}{\small\centering
		Repeat until SNF is achieved
	}\Bigg)
\end{longtable}

\noindent 
We also get the matrices $U \in \GL(m,\field[x])$ and $V \in \GL(n, \field)$ by unfolding the matrix reduction as follows:

{\setlength{\tabcolsep}{2pt}\begin{longtable}{C R RCL}
	D &:=& & D_{\min(m,n)} & \\[3pt]
	&=& 
		\bigl( U_{\min(m,n)} \bigr)\raisebox{4pt}{$\scriptstyle -1$}
		& \bigl( D_{\min(m,n)-1} \bigr)
		& \bigl( V_{\min(m,n)} \bigr)
	\\[3pt] 
	&& & \vdots \\[2pt]
	&=& 
		\bigl( U_{\min(m,n)} \bigr)\raisebox{4pt}{$\scriptstyle -1$}
			\cdots 
			(U_{k+1})\raisebox{4pt}{$\scriptstyle -1$}
			(U_{k})\raisebox{4pt}{$\scriptstyle -1$}
		& (D_{k-1})
		& (V_k) (V_{k+1})
			\cdots 
			\bigl( V_{\min(m,n)} \bigr) 
		\\ 
	&&& \vdots \\[3pt]
	&=& 
		\bigl( U_{\min(m,n)} \bigr)\raisebox{4pt}{$\scriptstyle -1$}
		\,\cdots\, 
		(U_k)\raisebox{4pt}{$\scriptstyle -1$}
		\,\cdots\,
		(U_2)\raisebox{4pt}{$\scriptstyle -1$}
		(U_1)\raisebox{4pt}{$\scriptstyle -1$}
		& D_0 
		& (V_1) (V_2) \,\cdots\, 
			(V_k) \,\cdots\, \bigl( V_{\min(m,n)} \bigr) 
	\\[5pt] 
	&=& 
		\biggl(\, 
			\underbrace{(U_1)(U_2) \,\cdots\, 
			(U_k) \,\cdots\, \bigl( U_{\min(m,n)} \bigr)}_{
				\raisebox{-5pt}{$U$}
			} 
		\,\biggr)\raisebox{4pt}{$\scriptstyle -1$}
		& [\phi]
		& \underbrace{(V_1) (V_2) \,\cdots\, 
			(V_k) \,\cdots\, \bigl( V_{\min(m,n)} \bigr)}_{
				\raisebox{-5pt}{$V$}
			}
\end{longtable}}

Below, we finish the calculation on $[\boundary_1\graded]$ started on \fref{ex:zom1-gradedboundary1-permutation}.

\begin{example}\label{ex:zom1-gradedboundary1-finishSND}
	Continue from \fref{ex:zom1-gradedboundary1-rowred}.
	In this example, we finish the calculation of an SND $(U,D,V)$ of 
	$[\boundary_1\graded]$ by \fref{alg:general-snf-graded-homs}.
	For reference, $[\phi] =: D_0 \in \M_{4,5}(\rationals[x])$ is given below.
	\begin{equation*}
		D_0 := [\boundary_1\graded]
		= \begin{NiceArray}{>{\color{gray}}c ccccc}
			\RowStyle[color=gray]{}
			& abx & bcx & adx^2 & cdx^2 & acx^3 \\
			{a}		&	-x & 0 & -x^2 & 0 & -x^3 \\
			{b}		&	x & -x & 0 & 0 & 0 \\ 
			{cx}	&	0 & 1 & 0 & x & x^2 \\
			{dx}	&	0 & 0 & x & -x & 0
		\CodeAfter \SubMatrix({2-2}{5-6})
		\end{NiceArray}
	\end{equation*}
	From \fref{ex:zom1-gradedboundary1-rowred}, 
	we have calculated the matrices $U_1 \in \GL(4,\rationals[x])$,
	$V_1 \in \GL(5,\rationals[x])$,
	and $D_1 \in \M_{4,5}(\rationals[x])$ as labeled in \fref{alg:general-snf-graded-homs}.
	For convenience, we copied the results below.
	\begin{equation*}
		U_1 =
		\begin{pmatrix}
			0 & 0 & 1 & 0 \\
			-x & 1 & 0 & 0 \\
			1 & 0 & 0 & 0 \\
			0 & 0 & 0 & 1 \\
		\end{pmatrix} 
		\qquad 
		D_1 = 
		\begin{pmatrix}
			1 &	0 &	0 &	0 &	0 \\
			0 &	x &	0 &	x^2 & x^3 \\
			0 & -x & -x^2 & 0 & x^3 \\
			0 & 0 & x & -x & 0
		\end{pmatrix}
		\qquad 
		V_1 = \begin{pmatrix}
			0 & 1 & 0 & 0 & 0 \\
			1 & 0 & 0 & -x & -x^2 \\
			0 & 0 & 1 & 0 & 0 \\
			0 & 0 & 0 & 1 & 0 \\
			0 & 0 & 0 & 0 & 1 
		\end{pmatrix}
	\end{equation*}

	\noindent\textbf{For the $k=2$ loop of \fref{alg:general-snf-graded-homs}:}

	We perform matrix operations on $D_1$.
	For brevity, we provide the matrices $U_2 \in \GL(4,\rationals[x])$ and $V_2 \in \GL(5,\rationals[x])$ that account for all row and column operations on $D_1$ respectively.
	Note that since $D_1(2,2) = x$ is already of minimal degree across the nonzero entries $D_1(j,i)$ with $j \in \set{2, \ldots, 4}$ and $i \in \set{2, \ldots, 5}$, 
	$D_1(2,2)$ can serve as the pivot and we do not need to permute $D_1$.
	Highlighted 
		in \bluetag is the \bluetagged{pivot $D_1(2,2)=x$}, 
		in \redtag are the entries to be eliminated by \redtagged{row reduction}, and 
		in \orangetag are those to be eliminated by \orangetagged{column operations}.
	\begin{equation*}
		D_1 = 
		\begin{pNiceMatrix}
			\CodeBefore [create-cell-nodes]
				\tikz \node [orange-highlight = (3-2) (4-2)] {} ;
				\tikz \node [red-highlight = (2-3) (2-5)] {} ;
				\tikz \node [blue-cell = (2-2)] {} ;
			\Body
			1 &	0 &	0 &	0 &	0 \\
			0 &	x &	0 &	x^2 & x^3 \\
			0 & -x & -x^2 & 0 & x^3 \\
			0 & 0 & x & -x & 0
		\end{pNiceMatrix}
	\end{equation*}
	We set $U_2 \in \GL(4,\rationals[x])$ for the row operations 
	and $V_2 \in \GL(5,\rationals)$ for the column operations as follows,
		with their respective pivot multipliers highlighted in 
		\orangetag and \redtag respectively.
	\begin{equation*}
		(U_2)\inv = \begin{pNiceMatrix}
			\CodeBefore [create-cell-nodes]
				\tikz \node [orange-cell = (3-2)] {} ;
				\tikz \node [orange-cell = (4-2)] {} ;
			\Body
			1 &	0 & 0 & 0 \\
			0 &	1 & 0 & 0 \\
			0 &	1 & 1 & 0 \\
			0 &	0 & 0 & 1 
		\end{pNiceMatrix}
		\quad\text{ with }\quad 
		U_2 = \begin{pNiceMatrix}
			\CodeBefore [create-cell-nodes]
				\tikz \node [orange-cell = (3-2)] {} ;
				\tikz \node [orange-cell = (4-2)] {} ;
			\Body
			1 &	0 & 0 & 0 \\
			0 &	1 & 0 & 0 \\
			0 &	-1 & 1 & 0 \\
			0 &	0 & 0 & 1 
		\end{pNiceMatrix}
		\quad\text{ and }\quad
		V_2 = \begin{pNiceMatrix}
			\CodeBefore [create-cell-nodes]
				\tikz \node [red-cell = (2-3)] {} ;
				\tikz \node [red-cell = (2-4)] {} ;
				\tikz \node [red-cell = (2-5)] {} ;
			\Body
			1 &	0 & 0 & 0 & 0 \\
			0 &	1 & 0 & -x & -x^2 \\
			0 &	0 & 1 & 0 & 0 \\
			0 &	0 & 0 & 1 & 0 \\
			0 &	0 & 0 & 0 & 1 
		\end{pNiceMatrix}
	\end{equation*}

	\noindent\textbf{For the $k=3$ loop of \fref{alg:general-snf-graded-homs}:}

	Given below is $D_2 := (U_2)\inv D_1 V_2$, with the chosen \bluetagged{pivot entry $D_1(4,3)=x$} highlighted in \bluetag\!.
	Since the pivot is not given by the $(3,3)$\th entry of $D_2$, we need to permute the rows and columns of $D_2$. 
	\begin{equation*}
		D_2 := (U_2)\inv D_1 V_2
		= {\arraycolsep=0.7\arraycolsep\def\arraystretch{0.9}
		\begin{pmatrix}
			1 &	0 & 0 & 0 \\
			0 &	1 & 0 & 0 \\
			0 &	1 & 1 & 0 \\
			0 &	0 & 0 & 1 
		\end{pmatrix}
		\begin{pmatrix}
			1 &	0 &	0 &	0 &	0 \\
			0 &	x &	0 &	x^2 & x^3 \\
			0 & -x & -x^2 & 0 & x^3 \\
			0 & 0 & x & -x & 0
		\end{pmatrix}
		\begin{pmatrix}
			1 &	0 & 0 & 0 & 0 \\
			0 &	1 & 0 & -x & -x^2 \\
			0 &	0 & 1 & 0 & 0 \\
			0 &	0 & 0 & 1 & 0 \\
			0 &	0 & 0 & 0 & 1 
		\end{pmatrix}}
		= \begin{pNiceMatrix}
			\CodeBefore [create-cell-nodes]
				\tikz \node [pivotcell = (4-3)] {} ;
			\Body
			1 &	0 &	0 &	0 &	0 \\
			0 &	x &	0 &	0 & 0 \\
			0 & 0 & -x^2 & x^2 & 0 \\
			0 & 0 & x & -x & 0
		\end{pNiceMatrix}
	\end{equation*}
	We have the following calculations for 
		$U_3 \in \GL(4,\rationals[x])$ and $V_3 \in \GL(5, \rationals[x])$ below.
	Highlighted in \greentag are the entries and rows related to \greentagged{row permutations}, those in \redtag for \redtagged{column reduction}, and those in \orangetag for \orangetagged{row reduction}.

	{\setlength{\tabcolsep}{2pt}\def\arraystretch{0.9}\arraycolsep=0.7\arraycolsep
	\begin{longtable}{C R L C C}
		W_3
		&:=& \elswap[4]{3,4} D_2
		= \begin{pNiceMatrix}
			\CodeBefore [create-cell-nodes]
				\tikz \node [green-highlight = (3-4)] {} ;
				\tikz \node [green-highlight = (4-3)] {} ;
			\Body
			1 &	0 & 0 & 0 \\
			0 &	1 & 0 & 0 \\
			0 &	0 & 0 & 1 \\
			0 &	0 & 1 & 0 
		\end{pNiceMatrix}
		\begin{pNiceMatrix}
			\CodeBefore [create-cell-nodes]
				\tikz \node [green-highlight = (3-1) (3-5)] {} ;
				\tikz \node [green-highlight = (4-1) (4-5)] {} ;
				\tikz \node [blue-cell = (4-3)] {} ;
			\Body
			1 &	0 &	0 &	0 &	0 \\
			0 &	x &	0 &	0 & 0 \\
			0 & 0 & -x^2 & x^2 & 0 \\
			0 & 0 & x & -x & 0
		\end{pNiceMatrix}
		= 
		\begin{pNiceMatrix}
			\CodeBefore [create-cell-nodes]
				\tikz \node [pivotcell = (3-3)] {} ;
			\Body
			1 &	0 &	0 &	0 &	0 \\
			0 &	x &	0 &	0 & 0 \\
			0 & 0 & x & -x & 0 \\
			0 & 0 & -x^2 & x^2 & 0 \\
		\end{pNiceMatrix}

		\\[25pt] 

		D_3
		&:=& \begin{pNiceMatrix}
			\CodeBefore [create-cell-nodes]
				\tikz \node [orange-highlight = (4-3)] {} ;
			\Body
			1 &	0 & 0 & 0 \\
			0 &	1 & 0 & 0 \\
			0 &	0 & 1 & 0 \\
			0 &	0 & x & 1 
		\end{pNiceMatrix}
		\begin{pNiceMatrix}
			\CodeBefore [create-cell-nodes]
				\tikz \node [orange-highlight = (4-3)] {} ;
				\tikz \node [red-highlight = (3-4)] {} ;
				\tikz \node [red-highlight = (3-5)] {} ;
				\tikz \node [blue-cell = (3-3)] {} ;
			\Body
			1 &	0 &	0 &	0 &	0 \\
			0 &	x &	0 &	0 & 0 \\
			0 & 0 & x & -x & 0 \\
			0 & 0 & -x^2 & x^2 & 0 \\
		\end{pNiceMatrix}
		\begin{pNiceMatrix}
			\CodeBefore [create-cell-nodes]
				\tikz \node [red-highlight = (3-4)] {} ;
				\tikz \node [red-highlight = (3-5)] {} ;
			\Body
			1 &	0 & 0 & 0 & 0 \\
			0 &	1 & 0 & 0 & 0 \\
			0 &	0 & 1 & 1 & 0 \\
			0 &	0 & 0 & 1 & 0 \\
			0 &	0 & 0 & 0 & 1 
		\end{pNiceMatrix}
		= 
		\begin{pNiceMatrix}
			\CodeBefore [create-cell-nodes]
				\tikz \node [orange-highlight = (4-3)] {} ;
				\tikz \node [red-highlight = (3-4)] {} ;
				\tikz \node [red-highlight = (3-5)] {} ;
				\tikz \node [blue-cell = (3-3)] {} ;
			\Body
			1 &	0 &	0 &	0 &	0 \\
			0 &	x &	0 &	0 & 0 \\
			0 & 0 & x & 0 & 0 \\
			0 & 0 & 0 & 0 & 0 
		\end{pNiceMatrix}

		\\[25pt] 

		(U_3)\inv 
		&:=& 
		\begin{pNiceMatrix}
			\CodeBefore [create-cell-nodes]
				\tikz \node [orange-highlight = (4-3)] {} ;
			\Body
			1 &	0 & 0 & 0 \\
			0 &	1 & 0 & 0 \\
			0 &	0 & 1 & 0 \\
			0 &	0 & x & 1 
		\end{pNiceMatrix}
		\begin{pNiceMatrix}
			\CodeBefore [create-cell-nodes]
				\tikz \node [green-highlight = (3-4)] {} ;
				\tikz \node [green-highlight = (4-3)] {} ;
			\Body
			1 &	0 & 0 & 0 \\
			0 &	1 & 0 & 0 \\
			0 &	0 & 0 & 1 \\
			0 &	0 & 1 & 0 
		\end{pNiceMatrix}
		= 
		\begin{pmatrix}
			1 &	0 & 0 & 0 \\
			0 &	1 & 0 & 0 \\
			0 &	0 & 0 & 1 \\
			0 &	0 & 1 & x 
		\end{pmatrix}
		
		\\[25pt]

		U_3
		&=&
		\begin{pNiceMatrix}
			\CodeBefore [create-cell-nodes]
				\tikz \node [green-highlight = (3-4)] {} ;
				\tikz \node [green-highlight = (4-3)] {} ;
			\Body
			1 &	0 & 0 & 0 \\
			0 &	1 & 0 & 0 \\
			0 &	0 & 0 & 1 \\
			0 &	0 & 1 & 0 
		\end{pNiceMatrix}\raisebox{6pt}{$\scriptstyle -1$}
		\begin{pNiceMatrix}
			\CodeBefore [create-cell-nodes]
				\tikz \node [orange-highlight = (4-3)] {} ;
			\Body
			1 &	0 & 0 & 0 \\
			0 &	1 & 0 & 0 \\
			0 &	0 & 1 & 0 \\
			0 &	0 & x & 1 
		\end{pNiceMatrix}\raisebox{6pt}{$\scriptstyle -1$}
		=
		\begin{pNiceMatrix}
			\CodeBefore [create-cell-nodes]
				\tikz \node [green-highlight = (3-4)] {} ;
				\tikz \node [green-highlight = (4-3)] {} ;
			\Body
			1 &	0 & 0 & 0 \\
			0 &	1 & 0 & 0 \\
			0 &	0 & 0 & 1 \\
			0 &	0 & 1 & 0 
		\end{pNiceMatrix}
		\begin{pNiceMatrix}
			\CodeBefore [create-cell-nodes]
				\tikz \node [orangecell = (4-3)] {} ;
			\Body
			1 &	0 & 0 & 0 \\
			0 &	1 & 0 & 0 \\
			0 &	0 & 1 & 0 \\
			0 &	0 & -x & 1 
		\end{pNiceMatrix}
		= 
		\begin{pmatrix}
			1 &	0 & 0 & 0 \\
			0 &	1 & 0 & 0 \\
			0 &	0 & -x & 1 \\
			0 &	0 & 1 & 0 
		\end{pmatrix}

		\\[25pt]

		V_3 &:=& 
		\begin{pNiceMatrix}
			\CodeBefore [create-cell-nodes]
				\tikz \node [red-highlight = (3-4)] {} ;
				\tikz \node [red-highlight = (3-5)] {} ;
			\Body
			1 &	0 & 0 & 0 & 0 \\
			0 &	1 & 0 & 0 & 0 \\
			0 &	0 & 1 & 1 & 0 \\
			0 &	0 & 0 & 1 & 0 \\
			0 &	0 & 0 & 0 & 1 
		\end{pNiceMatrix}
	\end{longtable}}

	\noindent\textbf{For the $k=4$ loop of \fref{alg:general-snf-graded-homs}:}

	Observe that $D_3$ is already in Smith Normal Form. 
	Then, \fref{alg:general-snf-graded-homs} stops and returns the matrices $U$ and $V$ as denoted below and the matrix $D := D_3$ as the Smith Normal Form of $D_0 = [\boundary_1\graded]$.
	
	{\setlength{\tabcolsep}{2pt}\def\arraystretch{0.9}\arraycolsep=0.7\arraycolsep
	\begin{longtable}{C C C L C}
		U &=& (U_1) (U_2) (U_3 )
		= 
			\begin{pmatrix}
				0 & 0 & 1 & 0 \\
				-x & 1 & 0 & 0 \\
				1 & 0 & 0 & 0 \\
				0 & 0 & 0 & 1 \\
			\end{pmatrix} 
			\begin{pmatrix}
				1 &	0 & 0 & 0 \\
				0 &	1 & 0 & 0 \\
				0 &	-1 & 1 & 0 \\
				0 &	0 & 0 & 1 
			\end{pmatrix}
			\begin{pmatrix}
				1 &	0 & 0 & 0 \\
				0 &	1 & 0 & 0 \\
				0 &	0 & -x & 1 \\
				0 &	0 & 1 & 0 
			\end{pmatrix}
		&=& 
		\begin{NiceArray}{>{\color{gray}}c !{\,} (cccc)}
			{a} & 0 & -1 & -x & 1 \\
			{b} & -x & 1 & 0 & 0 \\
			{cx} & 1 & 0 & 0 & 0 \\
			{dx} & 0 & 0 & 1 & 0
		\end{NiceArray}

		\\[25pt]

		V &=& (V_1) (V_2) (V_3)
		=
		\begin{pmatrix}
			0 &	1 & 0 & 0 & 0 \\
			1 &	0 & 0 & -x & -x^2 \\
			0 &	0 & 1 & 0 & 0 \\
			0 &	0 & 0 & 1 & 0 \\
			0 &	0 & 0 & 0 & 1 
		\end{pmatrix}
		\begin{pmatrix}
			1 &	0 & 0 & 0 & 0 \\
			0 &	1 & 0 & -x & -x^2 \\
			0 &	0 & 1 & 0 & 0 \\
			0 &	0 & 0 & 1 & 0 \\
			0 &	0 & 0 & 0 & 1 
		\end{pmatrix}
		\begin{pmatrix}
			1 &	0 & 0 & 0 & 0 \\
			0 &	1 & 0 & 0 & 0 \\
			0 &	0 & 1 & 1 & 0 \\
			0 &	0 & 0 & 1 & 0 \\
			0 &	0 & 0 & 0 & 1 
		\end{pmatrix}
		&=&
		\begin{NiceArray}{>{\color{gray}}c !{\,} (ccc cc)}
			{abx} & 0 &	1 & 0 & -x & -x^2 \\
			{bcx} & 1 &	0 & 0 & -x & -x^2 \\
			{adx^2} & 0 &	0 & 1 & 1 & 0 \\
			{cdx^2} & 0 &	0 & 0 & 1 & 0 \\
			{acx^3} & 0 &	0 & 0 & 0 & 1 
		\end{NiceArray}
	\end{longtable}}
	\noindent
	We can confirm that $U\inv [\boundary_1\graded] V = D$ by the following calculation:
	\begin{equation*}\def\arraystretch{0.9}\arraycolsep=0.7\arraycolsep
		U\inv [\boundary_1\graded] V 
		= 
		\begin{pmatrix}
			0 & 0 & 1 & 0 \\
			0 & 1 & x & 0 \\
			0 & 0 & 0 & 1 \\
			1 & 1 & x & 0 \\
		\end{pmatrix}
		\begin{pmatrix}
			-x & 0 & -x^2 & 0 & -x^3 \\
			x & -x & 0 & 0 & 0 \\ 
			0 & 1 & 0 & x & x^2 \\
			0 & 0 & x & -x & 0
		\end{pmatrix}
		\begin{pmatrix}
			0 &	1 & 0 & -x & -x^2 \\
			1 &	0 & 0 & -x & -x^2 \\
			0 &	0 & 1 & 1 & 0 \\
			0 &	0 & 0 & 1 & 0 \\
			0 &	0 & 0 & 0 & 1 
		\end{pmatrix}
		= \cdots
		= \begin{pmatrix}
			1 &	0 &	0 &	0 &	0 \\
			0 &	x &	0 &	0 & 0 \\
			0 & 0 & x & 0 & 0 \\
			0 & 0 & 0 & 0 & 0 
		\end{pmatrix}
		= D
	\end{equation*}
	Then, the matrices $U \in \GL(4,\rationals[x])$ and $V \in \GL(5,\rationals[x])$ induce bases 
		$B = (\beta_1(x), \ldots, \beta_4(x))$ of $C_0\graded(\filt{K};\rationals)$
		by $[\beta_j] = \col_j(U)$
		and 
		$T = (\tau_1(x), \ldots, \tau_5(x))$ of $C_1\graded(\filt{K};\rationals)$
		by $[\tau_i] = \col_i(V)$
	as follows:
	\begin{equation*}
		\begin{array}{c !{=} c !{=} r}
			\beta_1(x) 
				& (-x)(b) + (1)(cx) & (-b+c)x \\
			\beta_2(x)
				& (-1)(a) + (1)(b) & (-a+b)\phantom{x} \\
			\beta_3(x) 
				& (-x)(a) + (1)(dx) & (-a+d)x \\
			\beta_4(x) 
				& (1)(a) & (a)\phantom{x}
		\end{array}
	\end{equation*}
	\begin{equation*}
		\begin{array}{c !{=} c !{=} r}
			\tau_1(x) 
				& (1)(bc x) & (bc)x\phantom{^2} \\
			\tau_2(x) 
				& (1)(ab x) & (ab)x\phantom{^2} \\
			\tau_3(x) 
				& (1)(adx^2) & (ad)x^2 \\
			\tau_4(x) 
				& (-x)(abx) + (-x)(bcx) + (1)(adx^2) + (1)(cdx^2)
				& (-ab+bc+ad+cd)x^2 \\
			\tau_5(x) 
				& (-x^2)(abx) + (-x^2)(bcx) + (1)(acx^3)
				& (-ab + bc + ac)x^3
		\end{array}
	\end{equation*}
	Observe that both $B$ and $T$ are homogeneous bases and that the nonzero elements $1,x,x \in \rationals[x]$ of the SNF $D$ are also homogeneous.
\end{example}

With \fref{alg:general-snf-graded-homs} established, 
	we can now state the existence claim for the Graded Structure Theorem (\fref{thm:graded-structure-theorem}) in relation to the SNDs compatible with \fref{prop:how-to-get-graded-ifd}.

\begin{proposition}
	Let $M$ be a finitely generated graded $\field[x]$-module
	with graded presentation given by $\phi: F_S \to F_G$ and $\pi: F_G \to M$ 
	with homogeneous bases $S = (\sigma_1(x), \ldots, \sigma_n(x))$ of $F_S$ 
	and $A = (\alpha_1(x), \ldots, \alpha_m(x))$ of $F_G$.

	Then, there exists an SND $(U,D,V)$ of $[\phi]_{A,S}$ such that 
	the basis $B = (\beta_1(x), \ldots, \beta_m(x))$ of $F_G$ by $[\beta_j(x)] = \col_j(U)$ is a homogeneous basis
	and the nonzero elements $d_1(x), \ldots, d_r(x)$ of the SNF $D$ by $d_j(x) = D(j,j)$ are also homogeneous.
	Furthermore, 
		the basis $B = \set{\beta_j(x)}_{j=1}^m$ of $F_G$ and 
		the nonzero elements $\set{d_j(x)}_{j=1}^r$ 
		satisfy the hypothesis of \fref{prop:how-to-get-graded-ifd} 
		and can be used to calculate the graded invariant factor decomposition of $M$.
\end{proposition}
\begin{proof}
	Since $[\phi]_{A,S}$ satisfies \fref{lemma:matrices-of-graded-homs},
		\fref{alg:general-snf-graded-homs} will yield an SND $(U,D,V)$ with the required properties.
	Note that \fref{prop:column-reduction-preserves-hom} states that the properties given by \fref{lemma:matrices-of-graded-homs} are preserved in the column reduction operations done by \fref{alg:general-snf-graded-homs},
		\fref{prop:row-reduction-preserves-hom} states this for row reduction,
		and \fref{prop:perm-and-dilation-preserve-hom} states this for the permutations and dilations.

	The divisibility relation $d_1(x) \divides \cdots \divides d_r(x)$ is guaranteed by how the pivots in each iteration in \fref{alg:general-snf-graded-homs} are chosen, 
		i.e.\ $\degh(d_1(x)) \leq \degh(d_2(x)) \leq \cdots \degh(d_r(x))$ since each $d_r(x)$ is chosen to be of minimal degree.
	Note that while $d_j(x)$ is generally not in the form $x^{t_j}$ with $\degh(d_j(x)) = t_j$ as required by \fref{prop:how-to-get-graded-ifd}, there exists a nonzero $\mu_j \in \field$ such that $\mu \cdot d_j(x) = x^{t_j}$.
	Then, the SNF $D$ of $[\phi]$ can be multiplied by appropriate elementary dilations over $\field[x]$ to have $d_j(x) = x^{t_j}$ for all $j \in \set{1,\ldots, r}$.
\end{proof}

While we have not strictly established this for arbitrary graded chain complexes of free graded $\field[x]$-modules, we can continue the calculation on $[\boundary_1\graded]$ on \fref{ex:zom1-gradedboundary1-finishSND} and determine the $0$\th persistent homology module of the filtration $\filt{K}$.

\begin{example}\label{ex:zom1-gradedboundary1-finalcalc}
	Let $\filt{K}$ and $K$ be given as in \fref{ex:pershom-one}
	and orient $K$ with the vertex order $\Vertex(K) = (a,b,c,d)$.
	For convenience, an illustration of $K$ and $K_\bullet$ (without orientation) is copied below:
	\begin{center}
		\baselineCenter{\includegraphics[height=1.05in]{zomfig1/zom-simp.png}}
		\quad
		\baselineCenter{\includegraphics[height=1.1in]{zomfig1/zom1-cleaned.png}}
	\end{center}
	\vspace{3pt}
	We claim that the following sequence is a graded presentation of the $0$\th graded persistent homology module $H_0\graded(\filt{K};\rationals)$ of $\filt{K}$ with coefficients in $\rationals$:
	\begin{equation*}
		C_1\graded(\filt{K};\rationals)
			\,\,\Xrightarrow{\quad\boundary_1\graded\quad}\,\,
			C_0\graded(\filt{K};\rationals)
			\,\,\Xrightarrow{\quad\pi\quad}\,\, 
			H_0\graded(\filt{K};\rationals)
			\,\,\Xrightarrow{\quad\quad}\,\,
			0
	\end{equation*}
	where $\pi$ refers to the canonical quotient map.
	The justification for this is similar to the non-graded case by \fref{prop:chain-homology-presentation} with $\ker(\boundary_0\graded) = C_0\graded(\filt{K};\rationals)$. 

	An SND $(U,D,V)$ of $[\boundary_1\graded]$ in the graded presentation of 
	$H_0\graded(\filt{K};\rationals)$,
		as calculated in \fref{ex:zom1-gradedboundary1-finishSND},
	is given below.
	The columns of $U$ and the diagonal entries of $D$ highlighted in \bluetag correspond to a basis of \bluetagged{$\im(\phi)$},
		and the columns of $U$ highlighted in \redtag correspond to the basis elements that map to the \redtagged{free component} of 
		$H_0\graded(\filt{K};\rationals)$.
	\begin{equation*}\def\arraystretch{0.9}\arraycolsep=\arraycolsep
		U = \begin{NiceArray}{
			>{\color{gray}}c !{}( !{\,} cccc !{\,\,} )
		}
			\CodeBefore [create-cell-nodes]
				\tikz \node [blue-cell = (1-4) (2-4) (4-4)] {} ;
				\tikz \node [blue-cell = (1-2) (2-2) (4-2)] {} ;
				\tikz \node [blue-cell = (1-3) (4-3)] {} ;
				\tikz \node [red-cell = (1-5) (4-5)] {} ;
			\Body
			{a} & 0 & -1 & -x & 1 \\
			{b} & -x & 1 & 0 & 0 \\
			{cx} & 1 & 0 & 0 & 0 \\
			{dx} & 0 & 0 & 1 & 0
		\end{NiceArray}
		\qquad 
		D = \begin{NiceArray}{( !{\,\,} ccccc !{\,} )}
			\CodeBefore [create-cell-nodes]
				\tikz \node [blue-cell = (1-1)] {} ;
				\tikz \node [blue-cell = (2-2)] {} ;
				\tikz \node [blue-cell = (3-3)] {} ;
			\Body
			1 &	0 &	0 &	0 &	0 \\
			0 &	x &	0 &	0 & 0 \\
			0 & 0 & x & 0 & 0 \\
			0 & 0 & 0 & 0 & 0 
		\end{NiceArray}
		\qquad 
		V = \begin{NiceArray}{>{\color{gray}}c !{} (ccc cc)}
			{abx} & 0 &	1 & 0 & -x & -x^2 \\
			{bcx} & 1 &	0 & 0 & -x & -x^2 \\
			{adx^2} & 0 &	0 & 1 & 1 & 0 \\
			{cdx^2} & 0 &	0 & 0 & 1 & 0 \\
			{acx^3} & 0 &	0 & 0 & 0 & 1 
		\end{NiceArray}
	\end{equation*}
	Then, $H_0\graded(\filt{K};\rationals)$ is calculated using 
	\fref{prop:how-to-get-graded-ifd} as follows, with the \bluetagged{trivial and torsion summands} highlighted in \bluetag and the \redtagged{free summands} in \redtag\!.
	\vspace{-10pt}\begin{equation*}\arraycolsep=2pt\def\arraystretch{2}
	\begin{NiceArray}{
		ccc c ccccc 
	}
	\CodeBefore [create-cell-nodes]
		\tikz \node [blue-highlight = (1-3)] {} ;
		\tikz \node [blue-highlight = (1-5)] {} ;
		\tikz \node [blue-highlight = (1-7)] {} ;
		\tikz \node [red-highlight = (1-9)] {} ;
	\Body
		\RowStyle[]{\color{black}}
			&& \text{\footnotesize (trivial)} 
			&& \text{\footnotesize (torsion)}
			&& \text{\footnotesize (torsion)}
			&& \text{\footnotesize (free)}
		\\
		H_0\graded(\filt{K};\rationals)
			&\upgraded\cong 
			& \dfrac{ \rationals[x]\ket{ (-b+c)x } }
				{ \rationals[x]\ket{ (1)(-b+c)x } }
			&\oplus& \dfrac{ \rationals[x]\ket{ -a+b } }
				{ \rationals[x]\ket{ (x)(-a+b) } }
			&\oplus& \dfrac{ \rationals[x]\ket{ (-a+d)x } }
				{ \rationals[x]\ket{(x)(-a+d)x} }
			&\oplus& \rationals[x]\ket{a}
		\\ 
			&\upgraded\cong
			& {\color{gray} (0)}
			& {\color{gray}\oplus}
			& \paren{\dfrac{ \rationals[x]\ket{-a+b} }{
				\rationals[x]\ket{(x)(-a+b)}
			}}
			&\oplus& \Sigma^{1}\!\paren{\dfrac{ \rationals[x]\ket{-a+d} }{
				\rationals[x]\ket{(x)(-a+d)}
			}}
			&\oplus& \rationals[x]\ket{a}
		\\ 
			&\upgraded\cong
			&&
			& \dfrac{\rationals[x]}{(x)} 
			&\oplus& \Sigma^1\!\paren{ \dfrac{\rationals[x]}{(x)} }
			&\oplus& \rationals[x] 
	\end{NiceArray} 
	\end{equation*}
	We can also determine the $0$\th persistence barcode $\barcode(\filt{K};\rationals)$ of the filtration $\filt{K}$ by calculating the interval decomposition of the $0$\th persistent homology module $H_0(\filt{K};\rationals)$ by \fref{cor:interval-decomp-from-structure-theorem} as follows: 
	\vspace{-12pt}\begin{equation*}\arraycolsep=5pt\def\arraystretch{2}
	\begin{NiceArray}{
		ccc c ccc >{\hspace{10pt}\color{black}}l
	}
	\CodeBefore [create-cell-nodes]
		\tikz \node [bluecell = (1-3)] {} ;
		\tikz \node [bluecell = (1-5)] {} ;
		\tikz \node [bluecell = (1-7)] {} ;
		\tikz \node [redcell = (1-7)] {} ;
	\Body
		\RowStyle[]{\color{black}}
			&& \text{\footnotesize (torsion)}
			&& \text{\footnotesize (torsion)}
			&& \text{\footnotesize (free)}
		\\[-1pt]
		H_0\graded(\filt{K};\rationals)
			&\upgraded\cong& 
			\dfrac{\rationals[x]}{(x)} 
			&\oplus& 
			\Sigma^1\!\paren{ \dfrac{\rationals[x]}{(x)} }
			&\oplus& 
			\rationals[x]
			& \text{ as graded $\rationals[x]$-modules }
		\\[-2pt] 
		H_0(\filt{K};\rationals) 
			&\uppersmod\cong& 
				\intmod{[0,1)} 
			&\oplus& 
				\intmod{[1,2)}
			&\oplus& 
				\intmod{[0, \infty)}
			& \text{ as persistence modules over $\rationals$ }
		\\[-7pt] 
		\barcode(\filt{K};\rationals)
			&\ni&
			[0,1) &,& [1,2) &,& [0,\infty)
			& \text{ as intervals in a persistence barcode }
	\end{NiceArray}
	\end{equation*}
\end{example}\clearpage

\section{Matrix Calculation of Homology of Graded Chain Complexes}
\label{section:matrix-graded-chain-complex}

Let $C_\ast = (C_n, \boundary_n)_{n \in \ints}$ be a chain complex of free graded $\field[x]$-modules $C_n$ of finite rank with graded differentials $\boundary_n: C_{n-1} \to C_n$.
As with the case of chain complexes of free $R$-modules in \fref{section:snd-on-ungraded-chain-complexes}, 
	we can calculate the $n$\th homology of $C_\ast$ using SNDs of $[\boundary_{n+1}]$.
Recall that this calculation relies on the existence of a decomposition of $C_n$ by \fref{thm:chain-complex-decomposition} into three direct summands as denoted below:
\begin{equation*}
	C_n \upmod\cong K_n^\text{tor} \oplus K_n^\text{free} \oplus 
		(C_n \bigmod \ker\boundary_n)
\end{equation*}
We discussed that this decomposition can be represented by a specific SND 
$(U_{n+1}, D_{n+1}, V_{n+1})$ of $[\boundary_{n+1}]$ 
and we have also established that this SND can be calculated using an arbitrary SND $(W_{n+1}, D_{n+1}, V_{n+1})$ of $[\boundary_{n+1}]$.
Then, using the same arguments for Smith Normal Decompositions of matrices over $\field[x]$ given by \fref{lemma:matrices-of-graded-homs} in \fref{section:calculation-graded-ifds},
	we claim that the results of \fref{section:snd-on-ungraded-chain-complexes} also extend to the graded case.
That is, we have the following decomposition on $C_n$ as a \textit{graded} $\field[x]$-module:
\begin{equation*}
	C_n \,\upgraded\cong\, K_n^\text{tor} \oplus K_n^\text{free} \oplus 
		(C_n \bigmod \ker\boundary_n)
\end{equation*}
with graded $\field[x]$-modules $K_n^\text{tor}$, $K_n^\text{free}$, and $C_n \bigmod \ker\boundary_n$ defined similarly as in \fref{thm:chain-complex-decomposition}.

In \fref{section:construction-of-persistent-homology},
	we discussed how the persistent homology of filtrations can be expressed as the homology of graded chain complexes.
That is, given a filtration $\filt{K}$, we can determine the $n$\th persistent homology module $H_n(\filt{K};\field)$ of $\filt{K}$ with coefficients in $\field$ by the calculating the $n$\th simplicial persistent homology module:
\begin{equation*}
	H_n\graded(\filt{K};\field)
	\subgraded\cong 
	\frac{\ker(\boundary_n\graded)}{\im(\boundary_{n+1}\graded)} 
\end{equation*}
where, for all $n \in \ints$,
	the filtered $n$\th chain module $C_n\graded(\filt{K};\field)$ of $\filt{K}$ is a graded $\field[x]$-module
	and graded $n$\th boundary map 
	$\boundary_n\graded: C_n\graded(\filt{K};\field)
	\to C_{n-1}\graded(\filt{K};\field)$
	is a graded homomorphism.
Then, $H_n\graded(\filt{K};\field)$ admits the following graded presentation:
\begin{equation*}
	C_{n+1}\graded(\filt{K};\field)
		\,\,\Xrightarrow{\,\boundary_{n+1}\graded \,}\,\,
	\ker(\boundary_n\graded)
		\,\,\Xrightarrow{\,\pi\,}\,\,
	H_n\graded(\filt{K}; \field)
		\,\,\Xrightarrow{\,\,}\,\,
			0
\end{equation*}
Since the arguments for the calculation of graded chain complexes follow exactly as that of \fref{section:snd-on-ungraded-chain-complexes} (with the addition of the modifier \textit{graded} for $R = \field[x]$),
	we only present examples and some comments about calculation in this section.

Earlier in \fref{ex:zom1-gradedboundary1-finalcalc}, 
we calculated the $0$\th persistent homology module $H_0(\filt{K};\rationals)$ of some filtration $\filt{K}$ where $\ker(\boundary_0\graded) = C_0\graded(\filt{K};\rationals)$.
In the example below,
	we calculate the $1$\st persistent homology module 
	$H_1(\filt{K};\rationals)$ of the same filtration.

\begin{example}\label{ex:zom1-h1-matrix-calculation}
	Let $K$ be the simplicial complex given below with orientation $\Vertex(K) = (a,b,c,d)$.
	Let $\filt{K}$ be a filtration on $K$ with $K_t$ given as follows:
	\begin{center}
		\baselineCenter{\includegraphics[height=1.1in]{zomfig1/zom-simp.png}}
		\quad\quad 
		\baselineCenter{\includegraphics[height=1.2in]{zomfig1/zom1-cleaned.png}}
	\end{center} 
	To calculate $H_1(\filt{K};\rationals)$, we consider the matrices of the graded boundary maps 
	$\boundary_2\graded: C_2\graded(\filt{K};\rationals)$
	and 
	$\boundary_1\graded: C_1\graded(\filt{K};\rationals)$,
	which determine the torsion and free component of $H_1\graded(\filt{K};\rationals)$ respectively.
	Provided below are $[\boundary_1\graded]$ and $[\boundary_2\graded]$ relative to the standard bases.
	\begin{equation*}
		[\boundary_1\graded] 
		= \begin{NiceArray}{>{\color{gray}}c ccccc}
			\RowStyle[color=gray]{}
			& abx & bcx & adx^2 & cdx^2 & acx^3 \\
			{a}		&	-x & 0 & -x^2 & 0 & -x^3 \\
			{b}		&	x & -x & 0 & 0 & 0 \\ 
			{cx}	&	0 & 1 & 0 & -x & x^2 \\
			{dx}	&	0 & 0 & x & x & 0
		\CodeAfter \SubMatrix({2-2}{5-6})
		\end{NiceArray}
		\quad\quad\text{ and }\quad
		[\boundary_2\graded] = 
		\begin{NiceArray}{>{\color{gray}}c cc}
			\RowStyle[color=gray]{}
			{} 		& abc x^4 & acd x^5 \\
			{abx}	& 	x^3 & 0 \\
			{bcx}	& 	x^3 & 0 \\
			{adx^2}	& 	0 & -x^3 \\
			{cdx^2}	& 	0 & x^3 \\
			{acx^3}	&	-x & x^2 \\
		\CodeAfter \SubMatrix({2-2}{6-3})
		\end{NiceArray}
	\end{equation*}
	We want to find an SNDs $(U_1, D_1, V_1)$ and $(U_2, D_2, V_2)$ of $[\boundary_1\graded]$
	and $[\boundary_2\graded]$ respectively 
	such that a basis compatible 
		with the following decomposition of $C_1\graded(\filt{K};\rationals)$
	can be identified from said SNDs:
	\begin{equation*}
		C_1\graded(\filt{K}; \rationals)
		\,\upgraded\cong\, 
		K_1^\text{tor} \oplus K_1^\text{free} 
		\oplus \frac{
			C_1\graded(\filt{K}; \rationals)
		}{
			\ker(\boundary_1\graded)
		}
	\end{equation*}
	with $K_1^\text{tor}$ and $K_1^\text{free}$ as denoted in \fref{thm:chain-complex-decomposition}.
	Recall that 
	\begin{equation*}
		\ker(\boundary_1\graded) \cong K_1^\text{tor} \oplus K_1^\text{free} 
		\quad\text{ and }\quad 
		\text{T}\Bigl( H_1\graded(K_\bullet; \rationals) \Bigr)
		\cong 
		K_1^\text{tor} \bigmod \im\bigl( \boundary_2\graded \bigr)
	\end{equation*}
	where $\text{T}(-)$ denotes the torsion submodule of a graded $\field[x]$-module.
	Note that if we want to identify cycle representatives, the basis we get for $K_1^\text{tor}$ from the SND of $[\boundary_2\graded]$ must be identifiable in the basis for $\ker(\boundary_1\graded)$ we get from the SND of $[\boundary_1\graded]$.
	To get SNDs with this property, we do matrix reduction differently than described by \fref{alg:general-snf-graded-homs}.
	
	For the SND $(U_2, D_2, V_2)$ of $[\boundary_2\graded]$:
	Define $T_0 := [\boundary_2\graded]$ and let $T_k$ denote matrix $[\boundary_2\graded]$ after $k$ column or row elimination operations.
	We use the following color scheme in the calculation below:
	\begin{center}
		\bluetagged{chosen pivot}, \redtagged{entry to be eliminated},
		\greentagged{pivot multiplier},
		\orangetagged{affected row or column}
	\end{center} 
	\vspace{5pt}
	
	{\def\arraystretch{0.9}\begin{longtable}{L !{$=$} C !{$=$} C}
		T_1 := T_0 \cdot \eladd[2]{1,2 \,; x}
			& 
			\left(\,\,\begin{NiceMatrix}
				\CodeBefore [create-cell-nodes]
					\tikz \node [blue-cell = (5-1)] {} ;
					\tikz \node [red-cell = (5-2)] {} ;
				\Body
				x^3 & 0 \\
				x^3 & 0 \\
				0 & -x^3 \\
				0 & x^3 \\
				-x & x^2 \\
			\end{NiceMatrix}\,\,\right)
			\left(\,\,\begin{NiceMatrix}
				\CodeBefore [create-cell-nodes]
					\tikz \node [green-cell = (1-2)] {} ;
				\Body
				1 & x \\
				0 & 1
			\end{NiceMatrix}\,\,\right)
			& 
			\left(\begin{NiceMatrix}
				\CodeBefore [create-cell-nodes]
					\tikz \node [orange-highlight = (1-2) (5-2)] {} ;
				\Body
				x^3 & x^4 \\
				x^3 & x^4 \\
				0 & -x^3 \\
				0 & x^3 \\
				-x & 0 \\
			\end{NiceMatrix}\right)
		\\[30pt]
		
		T_2 := \eladd[5]{2,5 \,; x^2} \cdot T_1 
			& 
			\left(\,\,\begin{NiceMatrix}
				\CodeBefore [create-cell-nodes]
					\tikz \node [green-cell = (2-5)] {} ;
				\Body
				1 & 0 & 0 & 0 & 0 \\
				0 & 1 & 0 & 0 & x^2 \\
				0 & 0 & 1 & 0 & 0 \\
				0 & 0 & 0 & 1 & 0 \\
				0 & 0 & 0 & 0 & 1 
			\end{NiceMatrix}\,\,\right)
			\left(\,\,\begin{NiceMatrix}
				\CodeBefore [create-cell-nodes]
					\tikz \node [blue-cell = (5-1)] {} ;
					\tikz \node [red-cell = (2-1)] {} ;
				\Body
				x^3 & x^4 \\
				x^3 & x^4 \\
				0 & -x^3 \\
				0 & x^3 \\
				-x & 0 \\
			\end{NiceMatrix}\,\,\right)
			& 
			\left(\begin{NiceMatrix}
				\CodeBefore [create-cell-nodes]
					\tikz \node [orange-highlight = (2-1) (2-2)] {} ;
				\Body
				x^3 & x^4 \\
				0 & x^4 \\
				0 & -x^3 \\
				0 & x^3 \\
				-x & 0 \\
			\end{NiceMatrix}\right)
			\\[30pt]

		T_3 := \eladd[5]{1,5 \,; x^2} \cdot T_2 
			& 
			\left(\,\,\begin{NiceMatrix}
				\CodeBefore [create-cell-nodes]
					\tikz \node [green-cell = (2-5)] {} ;
				\Body
				1 & 0 & 0 & 0 & x^2 \\
				0 & 1 & 0 & 0 & 0 \\
				0 & 0 & 1 & 0 & 0 \\
				0 & 0 & 0 & 1 & 0 \\
				0 & 0 & 0 & 0 & 1 
			\end{NiceMatrix}\,\,\right)
			\left(\,\,\begin{NiceMatrix}
				\CodeBefore [create-cell-nodes]
					\tikz \node [blue-cell = (5-1)] {} ;
					\tikz \node [red-cell = (1-1)] {} ;
				\Body
				x^3 & x^4 \\
				0 & x^4 \\
				0 & -x^3 \\
				0 & x^3 \\
				-x & 0 \\
			\end{NiceMatrix}\,\,\right)
			& 
			\left(\begin{NiceMatrix}
				\CodeBefore [create-cell-nodes]
					\tikz \node [orange-highlight = (1-1) (1-2)] {} ;
				\Body
				0 & x^4 \\
				0 & x^4 \\
				0 & -x^3 \\
				0 & x^3 \\
				-x & 0 \\
			\end{NiceMatrix}\right)
			\\[30pt]

		T_4 := \eladd[5]{3,4 \,; 1} \cdot T_3
			& 
			\left(\,\,\begin{NiceMatrix}
				\CodeBefore [create-cell-nodes]
					\tikz \node [green-cell = (3-4)] {} ;
				\Body
				1 & 0 & 0 & 0 & 0 \\
				0 & 1 & 0 & 0 & 0 \\
				0 & 0 & 1 & 1 & 0 \\
				0 & 0 & 0 & 1 & 0 \\
				0 & 0 & 0 & 0 & 1 
			\end{NiceMatrix}\,\,\right)
			\left(\,\,\begin{NiceMatrix}
				\CodeBefore [create-cell-nodes]
					\tikz \node [blue-cell = (4-2)] {} ;
					\tikz \node [red-cell = (3-2)] {} ;
				\Body
				0 & x^4 \\
				0 & x^4 \\
				0 & -x^3 \\
				0 & x^3 \\
				-x & 0 \\
			\end{NiceMatrix}\,\,\right)
			& 
			\left(\begin{NiceMatrix}
				\CodeBefore [create-cell-nodes]
					\tikz \node [orange-highlight = (3-1) (3-2)] {} ;
				\Body
				0 & x^4 \\
				0 & x^4 \\
				0 & 0 \\
				0 & x^3 \\
				-x & 0 \\
			\end{NiceMatrix}\right)
			\\[30pt]

		T_5 := \eladd[5]{2,4 \,; -x} \cdot T_4 
			& 
			\left(\,\,\begin{NiceMatrix}
				\CodeBefore [create-cell-nodes]
					\tikz \node [green-cell = (2-4)] {} ;
				\Body
				1 & 0 & 0 & 0 & 0 \\
				0 & 1 & 0 & -x & 0 \\
				0 & 0 & 1 & 0 & 0 \\
				0 & 0 & 0 & 1 & 0 \\
				0 & 0 & 0 & 0 & 1 
			\end{NiceMatrix}\,\,\right)
			\left(\begin{NiceMatrix}
				\CodeBefore [create-cell-nodes]
					\tikz \node [blue-cell = (2-2)] {} ;
					\tikz \node [red-cell = (4-2)] {} ;
				\Body
				0 & x^4 \\
				0 & x^4 \\
				0 & 0 \\
				0 & x^3 \\
				-x & 0 \\
			\end{NiceMatrix}\right)
			& 
			\left(\begin{NiceMatrix}
				\CodeBefore [create-cell-nodes]
					\tikz \node [orange-highlight = (2-1) (2-2)] {} ;
				\Body
				0 & x^4 \\
				0 & 0 \\
				0 & 0 \\
				0 & x^3 \\
				-x & 0 \\
			\end{NiceMatrix}\right)
			\\[30pt]

		T_6 := \eladd[5]{1,4 \,; -x} \cdot T_5
			& 
			\left(\,\,\begin{NiceMatrix}
				\CodeBefore [create-cell-nodes]
					\tikz \node [green-cell = (1-4)] {} ;
				\Body
				1 & 0 & 0 & -x & 0 \\
				0 & 1 & 0 & 0 & 0 \\
				0 & 0 & 1 & 0 & 0 \\
				0 & 0 & 0 & 1 & 0 \\
				0 & 0 & 0 & 0 & 1 
			\end{NiceMatrix}\,\,\right)
			\left(\,\,\begin{NiceMatrix}
				\CodeBefore [create-cell-nodes]
					\tikz \node [blue-cell = (4-2)] {} ;
					\tikz \node [red-cell = (1-2)] {} ;
				\Body
				0 & x^4 \\
				0 & 0 \\
				0 & 0 \\
				0 & x^3 \\
				-x & 0 \\
			\end{NiceMatrix}\,\,\right)
			& 
			\left(\begin{NiceMatrix}
				\CodeBefore [create-cell-nodes]
					\tikz \node [orange-highlight = (1-1) (1-2)] {} ;
				\Body
				0 & 0 \\
				0 & 0 \\
				0 & 0 \\
				0 & x^3 \\
				-x & 0 \\
			\end{NiceMatrix}\right)
	\end{longtable}}\noindent 
	Observe that the matrix $T_6$ is only a few elementary permutations away from being in Smith Normal Form.
	Instead of doing row and column permutations on $T_6$, 
	we stop here and work with the SND $(U_2, D_2, V_2)$ of $[\boundary_2\graded]$ as given below,
		with the columns of $U_2$ corresponding to a basis of $K_1^\text{tor}$ and the entries in $D_2$ corresponding to the invariant factors of $H_1\graded(K_\bullet;\rationals)$ highlighted in \redtag and \purpletag\!.
	\vspace{-5pt}
	\begin{longtable}{L@{\ }C@{\ }L}
		V_2 
		&=&
		\eladd[2]{1,2 \,; x} 
		= \left(\begin{NiceMatrix}
			1 & x \\
			0 & 1
		\end{NiceMatrix}\right),
		\qquad\text{ and }\qquad 
		D_2
			= T_6
			= \left(\,\,\begin{NiceMatrix}
				\CodeBefore [create-cell-nodes]
					\tikz \node [red-cell = (4-2), inner xsep=4pt] {} ;
					\tikz \node [purple-cell = (5-1), inner xsep=4pt] {} ;
				\Body
				0 & 0 \\
				0 & 0 \\
				0 & 0 \\
				0 & x^3 \\
				-x & 0 \\
			\end{NiceMatrix}\,\,\right),
		\\[30pt]
		U_2 
		&=& \biggl( 
			\eladd[5]{1,4 \,; -x} \cdot
			\eladd[5]{2,4 \,; -x} \cdot
			\eladd[5]{3,4 \,; 1}  \cdot
			\eladd[5]{1,5 \,; x^2} \cdot
			\eladd[5]{2,5 \,; x^2}
		\biggr)\raisebox{5pt}{$-1$}
		\\[15pt]
		&=& 
			\eladd[5]{2,5 \,; -x^2} \cdot
			\eladd[5]{1,5 \,; -x^2} \cdot
			\eladd[5]{3,4 \,; -1} \cdot 
			\eladd[5]{2,4 \,; x} \cdot
			\eladd[5]{1,4 \,; x}
		\\[15pt]
		&=& 
		\hspace{3pt}\begin{NiceArray}{>{\color{gray}}c @{\hspace{13pt}}
			ccc cc 
		}
		\CodeBefore [create-cell-nodes]
			\tikz \node [red-cell = (1-5) (3-5) (5-5), inner xsep=4pt] {} ;
			\tikz \node [purple-cell = (1-6) (5-6), inner xsep=4pt] {} ;
		\Body
			{abx} 	& 1 & 0 & 0 & x & -x^2 \\
			{bcx} 	& 0 & 1 & 0 & x & -x^2 \\
			{adx^2} 	& 0 & 0 & 1 & -1 & 0 \\
			{cdx^2} 	& 0 & 0 & 0 & 1 & 0 \\
			{acx^3} 	& 0 & 0 & 0 & 0 & 1 
		\CodeAfter
			\SubMatrix({1-2}{5-6})[left-xshift=3pt, right-xshift=3pt]
			\UnderBrace{5-5}{5-6}{\tikz \node [red-node-padded, yshift=-10pt] {
					$K_1^\text{tor}$
				};}[yshift=5pt, shorten]
		\end{NiceArray}
	\end{longtable} 
	\vspace{25pt}

	\noindent 
	Let $\redmath{\beta_1 x^{s_1}}, \purplemath{\beta_2 x^{s_2}} \in C_n\graded(K_\bullet;\rationals)$
	be given by $\redmath{[\beta_1 x^{s_1}] = \col_4(U_2)}$ and $\purplemath{[\beta_2 x^{s_2}] = \col_5(U_2)}$, following the color scheme above.
	Then, we have the following information for the graded invariant factor decomposition of $K_1^\text{tor} \bigmod \im(\boundary_2\graded)$:
	\begin{longtable}{ccccc}
		Filtered Cycle &
		Cycle in $C_1(K;\rationals)$ & Degree of Invariant Factor & Grading Shift
		\\[5pt]
		$\redmath{\beta_1 x^{s_1}} = (ab+bc-ad+cd)x^2$ &
		$\beta_1 = ab+bc-ad+cd$ & $t_1 = \degh(x^3) = 3$ & $s_1 = 2$ 
		\\[5pt]
		$\purplemath{\beta_2 x^{s_2}} = (-ab-bc+ac)x^3$ &
		$\beta_2 = -ab-bc+ac$ & $t_1 = \degh(-x) = 1$
		& $s_2 = 3$
	\end{longtable}
	\noindent 
	Then, the torsion component 
	$\text{T}( H_1\graded(K_\bullet; \rationals) )$
	of $H_1\graded(K_\bullet;\rationals)$ is given by:
	\begin{align*}
		\text{T}\Bigl( H_1\graded(K_\bullet; \rationals) \Bigr)
		\upgraded\cong 
		\frac{ K_1^\text{tor} }{ \im(\boundary_2\graded) }
		&\cong 
			\frac{ \rationals[x]\ket{\beta_1 x^2} }{ \rationals[x]\ket{x^3 \cdot \beta_1 x^2 } }
			\oplus 
			\frac{ \rationals[x]\ket{\beta_2 x^3} }{ \rationals[x]\ket{x \cdot \beta_2 x^3} }
		\cong 
		\Sigma^2 \left(\frac{\rationals[x]}{(x^3)}\right)
		\oplus 
		\Sigma^3 \left(\frac{\rationals[x]}{(x)}\right)
	\end{align*}

	\spacer 
	To determine the free component 
		$\text{F}( H_1\graded(K_\bullet; \rationals) ) \cong K_1^\text{free}$
	of $H_1\graded(K_\bullet;\rationals)$,
	we need information about $\ker(\boundary_1\graded) = K_1^\text{free} \oplus K_1^\text{tor}$, which we can get from an SND of $[\boundary_1\graded]$.
	Given below is the SND $(U_1, D_1, V_1)$ of $[\boundary_1\graded]$ calculated in \fref{ex:column-reduction-of-graded-boundary-one} (column reduction)
	and \fref{ex:row-red-of-graded-bdry-one} (row reduction):
	\begin{equation*}
		U_1 = 
		\left(\begin{NiceMatrix}
			1 & 1 & x & x \\
			0 & 1 & x & 0 \\
			0 & 0 & 1 & 0 \\
			0 & 0 & 0 & 1 \\
		\end{NiceMatrix}\right)
		\qquad
		D_1 = 
		\left(\,\,\begin{NiceMatrix} 
			\CodeBefore [create-cell-nodes]
				\tikz \node [orange-highlight = (1-4) (4-4)] {} ;
				\tikz \node [orange-highlight = (1-5) (4-5)] {} ;
			\Body
			0 & 0 & 0 & 0 & 0 \\
			x & 0 & 0 & 0 & 0 \\
			0 & 1 & 0 & 0 & 0 \\
			0 & 0 & x & 0 & 0
		\end{NiceMatrix}\,\,\right) 
		\qquad 
		V_1 = 
		\hspace{3pt}\begin{NiceArray}{>{\color{gray}}c @{\hspace{13pt}}
			ccc cc 
		}
		\CodeBefore [create-cell-nodes]
			\tikz \node [orange-cell = (1-5) (3-5) (5-5), inner xsep=4pt] {} ;
			\tikz \node [orange-cell = (1-6) (5-6), inner xsep=4pt] {} ;
		\Body
			{abx} 	& 1 & 0 & 0 & x & -x^2 \\
			{bcx} 	& 0 & 1 & 0 & x & -x^2 \\
			{adx^2} 	& 0 & 0 & 1 & -1 & 0 \\
			{cdx^2} 	& 0 & 0 & 0 & 1 & 0 \\
			{acx^3} 	& 0 & 0 & 0 & 0 & 1 
		\CodeAfter
			\SubMatrix({1-2}{5-6})[left-xshift=3pt, right-xshift=3pt]
			\UnderBrace{5-5}{5-6}{\tikz \node [orange-node-padded, yshift=-10pt] {
					$\ker(\boundary_1\graded)$
				};}[yshift=5pt, shorten]
		\end{NiceArray} 
	\end{equation*} 
	\vspace{17pt}

	\noindent 
	Since the $4$\th and $5$\th columns of the SNF $D_1$ of $[\boundary_1\graded]$ are \orangetagged{zero columns}\!, 
		the $4$\th and $5$\th columns of $V_1$ determine a basis of $\ker(\boundary_1\graded)$.
	Observe that these two columns are exactly the columns corresponding to 
		$\redmath{\beta_1 x^{s_1}}$ and $\purplemath{\beta_2 x^{s_2}}$.
	Then, $K_1^\text{free} = 0$ and the free component of $H_1\graded(K_\bullet, \rationals)$ is also trivial.
	Therefore,
	\begin{equation*}
		H_1\graded(K_\bullet;\rationals)
		\upgraded\cong 
		\text{F}\Bigl( H_1\graded(K_\bullet; \rationals) \Bigr)
		\oplus 
		\text{T}\Bigl( H_1\graded(K_\bullet; \rationals) \Bigr)
		\cong 
		0 \oplus \Sigma^2 \left(\frac{\rationals[x]}{(x^3)}\right)
		\oplus 
		\Sigma^3 \left(\frac{\rationals[x]}{(x)}\right)
	\end{equation*}
	By application of $\topersmod(-)$ on the graded invariant factor decomposition of $H_1\graded(K_\bullet;\rationals)$,
		we get the following for the interval decomposition and persistence barcode $\barcode_1(K_\bullet;\rationals)$ of $H_1(K_\bullet;\rationals)$:
	\begin{gather*}
		H_1(K_\bullet; \rationals) 
		\,\uppersmod\cong\,
		\topersmod\Biggl(\,
			\Sigma^2 \left(\frac{\rationals[x]}{(x^3)}\right)
			\oplus 
			\Sigma^3 \left(\frac{\rationals[x]}{(x)}\right)
		\,\Biggr)
		\cong 
		\intmod{[2,2+3)} \oplus \intmod{[3, 3+1)}
		\\[2pt]
		\text{ and }
		\\[2pt]
		\barcode_1(K_\bullet; \rationals) = \Bigl\{
			[2,5), [3,4)
		\Bigr\}
	\end{gather*}
	The intervals in the $1$\st persistent barcode $\barcode_1(K_\bullet;\rationals)$ of $K_\bullet$ over $\rationals$ 
	correspond to the $1$\st filtered homology classes illustrated below:
	\begin{longtable}{c@{\qquad}c}
		Interval Module 
		& Illustrated in $K_\bullet$
		\\[5pt]
		$\begin{gathered}
			\displaystyle\intmod{[2,2+3)} \cong \topersmod\left(
			\frac{ \rationals[x]\ket{\beta_1 x^2} }{ \rationals[x]\ket{x^3 \cdot \beta_1 x^2 } }
			\right)
			\\
			\text{ with } \beta_1 = ab+bc-ad+cd
		\end{gathered}$
			& 
			\baselineCenter{\includegraphics[width=4in]{zomfig1/zom1-h1-beta1.png}}
		\\[15pt]
		$\begin{gathered}
			\displaystyle\intmod{[3,3+1)} \cong \topersmod\left(
				\frac{ \rationals[x]\ket{\beta_2 x^3} }{ \rationals[x]\ket{x \cdot \beta_2 x^3} }
			\right)
			\\ 
		\text{ with } \beta_2 = -ab-bc+ac
		\end{gathered}$
			& \baselineCenter{\includegraphics[width=4in]{zomfig1/zom1-h1-beta2.png}}
	\end{longtable}
\end{example} 
\clearpage

\onlyifstandalone{\input{"../BX. Reference Page.tex"}}


\addcontentsline{toc}{chapter}{Acknowledgements}
\chapter*{Acknowledgements}

This expository paper is the culmination of my two and a half year's stay at Oregon State University.
I will be the first to admit that this paper is far from perfect and there is so much more that can be done to improve it.
The weaknesses of this paper are reflective of my journey as a graduate student and a learner of mathematics.
Before I started my program at OSU, 
I had very little to no knowledge of topics like category theory, graded module theory, algebraic topology, and homological algebra -- all of which I had to learn in sufficient depth to write this paper.
Mistakes were a huge part of my journey but I like to see my recognition of these mistakes as part of growth and as a testament to how much I have improved from where I started.

Despite this paper's weaknesses, 
	I'd like to think that this paper has some meaningful contribution to mathematics.
Persistent homology theory is a new and young field of mathematics, one that emerged only a few decades ago.
A huge part of the writing process involved gathering all these seemingly disparate papers and results, and trying to make sense of the differences in definitions, convention, and notation.
This paper is my attempt at reconciling these differences and making the connections between all these algebraic concepts clearer and more concrete.
I'd like to think that, in the future, what I had written would help another student progress further in their learning journey.

The strengths of this paper are also reflective of my support system, without whom I could not have completed this paper and my Master's degree. 
I'd like to express my sincere gratitude to the following people:
\begin{enumerate}
	\item 
	My advisor, Dr.\ Christine Escher, 
	for her unwavering support and guidance throughout my career at OSU, undeterred by my failings and weaknesses.

	\item 
	Dr.\ Chad Giusti and Dr.\ Ren Guo, for their encouragement, patience, and understanding as part of my graduate committee, especially during the times I had faced hardships and roadblocks.

	\item 
	The mathematics community at OSU, including my fellow graduate students, the graduate and undergraduate math faculty, and the math support staff, for making me feel welcome at OSU.

	\item 
	My family and friends, for their support throughout the challenges I've faced, whether they be academic, financial, emotional, or medical in nature.
\end{enumerate}

\onlyifstandalone{\input{"../BX. Reference Page.tex"}}

\clearpage
 
	\onlyifstandalone{
		\renewcommand*\contentsname{Table of Contents}
		\pdfbookmark{\contentsname}{toc}
		\begin{spacing}{1.2} 
			\tableofcontents 
		\end{spacing} 
		\clearpage
	}

	\appendix
	\titleformat{\section}[hang]{\large\bfseries}{Appendix \thesection.}{1em}{}[\vspace{-10pt}\grayline] 
	\renewcommand{\thesection}{A\arabic{section}}%
	\addcontentsline{toc}{chapter}{Appendices}

{ 
\section{List of Symbols}
\label{appendix:miscellaneous}

\newlist{symbollist}{itemize}{2}
\setlist[symbollist]{
	label={}, 
	labelwidth=40pt, labelsep=5pt, itemsep=-0pt,
	leftmargin=45pt, 
	itemindent=0pt, align=left
}

\raggedright
\setlength{\columnsep}{20pt}

\vspace{0pt}

\noindent
\textbf{Symbols involving Sets, Modules, and Graded Modules.}
\vspace{-7pt}
\begin{multicols}{2}\begin{symbollist}
	\item[$\nonnegints$] 
	the nonnegative integers, 
	i.e.\ $\nonnegints = \set{n \in \ints: n \geq 0} = \set{0, 1, 2, \ldots}$.

	\item[$R$]
	unless otherwise specified, typically used to refer to an arbitrary principal ideal domain (PID).

	\item[$R^\times$] 
	the group of units (invertible elements) of some ring $R$.

	\item[$\field$]
	typically used to refer to an arbitrary field.

	\item[$\text{F}(M)$] 
	the free component of an $R$-module $M$. 

	\item[$\text{T}(M)$] 
	the torsion component of an $R$-module $M$.

	\item[$\frac{M}{N}$] 
	usu.\ the quotient module with some $R$-module $M$ and submodule $N$ of $M$; equiv.\ to $M \bigmod N$.

	\item[$R\ket{-}$]
	refers to the free $R$-module formed by $R$-formal sums of elements in the argument,
	see \fref{defn:formal-sums}.

	\columnbreak

	\item[${R[x]}$] 
	the polynomial ring with indeterminate $x$,
	equipped with the standard grading when grading is considered, see \fref{defn:graded-polynomial-ring}.

	\item[$Mx^t$]
	the $R$-module isomorphic to $M$ by
	$mx^t \mapsto m \in M$;
	may be used to identify the degree of a homogeneous component of a graded $R[x]$-module, see \fref{remark:assumption-on-graded-modules}.

	\item[$\deg(-)$]
	the degree of a polynomial with indeterminate $x$,
	not necessarily as an element of a graded $\field[x]$-module.

	\item[$\degh(-)$]
	the degree of a homogeneous element of a graded $\field[x]$-module, usually used when a relation is only valid for homogeneous elements,
	agrees with $\deg(-)$ for nonzero homogeneous elements and is undefined for non-homogeneous elements and zero.

	\item[$\id_A$] 
	usu.\ the identity function on some set $A$.
\end{symbollist}\end{multicols}

\noindent
\textbf{Symbols for Relevant Categories.}
\vspace{-7pt}
\begin{multicols}{2}\begin{symbollist}
	\item[$\catname{C}$] 
	usu.\ used to refer to an arbitrary category.

	\item[$\cattop$]
	the category of topological spaces and continuous maps.

	\item[$\poset(R,\leq)$]
	the category induced by a poset $(R,\leq)$,
	see \fref{defn:poset-cat},
	conventionally denoted as $(R,\leq)$ in category theory;

	not to be confused with the category $\catname{Poset}$ of posets and order-preserving maps (not used in this paper).

	\item[$\catsimp$]
	the category of (abstract) simplicial complexes and simplicial maps, 
	see \fref{defn:cat-simp}.

	\item[$\catname{A}$]
	usu.\ for an arbitrary abelian category.

	\item[$\catname{AbGrp}$]
	the category of abelian groups and group homomorphisms.

	\item[$\catmod{R}$]
	the category of right $R$-modules and $R$-module homomorphisms for some ring $R$,
	the modifier ``right'' is usu.\ dropped when $R$ is commutative.

	\item[$\catmod{\ints}$] 
	category of $\ints$-modules, equiv.\ to $\catname{AbGrp}$.

	\columnbreak

	\item[$\catmod{R[x]}$]
	the category of ungraded $R[x]$-modules, i.e.\ disregarding any graded structure (if such exists).

	\item[$\catgradedmod{R}$]
	the category of graded modules over $R[x]$ and graded $R[x]$-module homomorphisms, 
	see \fref{thm:category-of-graded-modules}.

	\item[$\catpersmod$]
	the category of persistence modules over $\field$ and persistence morphisms,
	see \fref{defn:persmod-cat}.

	\item[$\catchaincomplex{\catname{A}}$]
	the category of chain complexes of an abelian category $\catname{A}$ and chain maps;
	in this paper, $\ast$ is used for the index $n \in \ints$, e.g.\ $C_\ast = (C_n, \boundary_n)_{n \in \ints}$ is a chain complex.

	\item[$\catchaincomplex{\catgradedmod{\field}}$]
	the category of graded chain complexes,
	see \fref{defn:cat-graded-chain-complexes}.

	\item[$\catchaincomplex{\catpersmod}$]
	the category of persistence complexes, see \fref{defn:category-of-persistence-chains}.

	\item[$H_n(-)$] 
	the $n$\th chain homology functor 
	$\catchaincomplex{\catname{A}} \to \catname{A}$ in some abelian category $\catname{A}$. \\
	For $\catname{A} = \catgradedmod{\field}$, see \fref{defn:cat-graded-chain-complexes}. 
	For $\catname{A} = \catpersmod$, see \fref{defn:category-of-persistence-chains}(iii).

	\item[$\id_{\catname{C}}$] 
	the identity functor on a category $\catname{C}$.
\end{symbollist}\end{multicols}

\noindent
\textbf{Symbols and Shorthands for Specific Binary Relations.}
\vspace{-7pt}
\begin{multicols}{2}\begin{symbollist}
	\item[$\oplus$]
	refers to the direct sum operation on an abelian category $\catname{A}$;
	when the category $\catname{A}$ is ambiguous or emphasized, an accompanying relation $\cong_\catname{A}$ is identified in the relevant line.
	
	For $\catname{A} = \catmod{\field[x]}$ or $\catname{A} = \catgradedmod{\field}$, see \fref{remark:shorthand-for-graded-relations}.
	For $\catname{A} = \catpersmod$, see \fref{defn:persmod-directsum}.

	\item[$\cong_{\catname{C}}$]
	denotes an isomorphism relation on a category $\catname{C}$, usu.\ used when the category in which the relation is considered is ambiguous or emphasized.

	\item[$\smash{\upabelian=}$] 
	denotes an equality at the level of abelian groups or $\field[x]$-vector spaces,
	see \fref{remark:shorthand-for-graded-relations}.

	\item[$\smash{\upmod\cong}$]
	shorthand for an isomorphism relation in $\catmod{R[x]}$, disregarding any graded structure (if such exists), 
	see \fref{remark:shorthand-for-graded-relations}.

	\item[$\smash{\upgraded\cong}$]
	shorthand for an isomorphism relation in $\catgradedmod{\field}$,
	see \fref{remark:shorthand-for-graded-relations}.

	\item[$\smash{\uppersmod\cong}$]
	shorthand for an isomorphism relation in $\catpersmod$, see remark under \fref{defn:persistence-isomorphism}.
\end{symbollist}\end{multicols}

\noindent
\textbf{Notation involving Simplicial Complexes and Simplicial Homology.} See \fref{chapter:simplicial-homology}.
\vspace{-7pt}
\begin{multicols}{2}\begin{symbollist}
	\item[$K$] 
	unless otherwise specified, usually used to refer to an arbitrary simplicial complex.

	\item[$\Vertex(K)$]
	the vertex set of a simplicial complex $K$, 
	see \fref{defn:simp-complex-abstract}.

	\item[$\basis{K}_n$]
	the standard ordered basis of $C_n(K;R)$ 
	of a simplicial complex $K$ relative to some given orientation on $\Vertex(K)$, 
	see \fref{defn:standard-basis-on-chain-groups}.

	\item[$C_n(K;R)$] 
	the $n$\th simplicial chain group of a simplicial complex $K$, see \fref{defn:chain-groups}.

	\item[$C_n(K)$] 
	shorthand for $C_n(K;\ints)$.

	\item[$C_n(-;R)$]
	the $n$\th simplicial chain group functor 
	$\catsimp \to \catmod{R}$ with coefficients in $R$, see \fref{defn:simplicial-chain-group-functor}.

	\item[$f_{n,\hash}$] 
	usu.\ the map $C_n(K;R) \to C_n(L;R)$ on the $n$\th simplicial chain groups induced by the simplicial map $f: K \to L$, see \fref{defn:induced-maps-on-chain-groups}.

	\item[$f_\hash$] 
	shorthand for $f_{n,\hash}$, used when the dimension $n$ is arbitrary or unambiguous.

	\item[${[v_0, \ldots, v_n]}$]
	an oriented $n$-simplex, 
	i.e.\ a basis element of an $n$\th chain group $C_n(K;R)$ with ordering $(v_0, \ldots, v_n)$,
	see \fref{defn:standard-basis-on-chain-groups}.

	\item[${v_0 \cdots v_n}$] 
	string shorthand for the oriented $n$-simplex $[v_0, \ldots, v_n]$,
	see \fref{remark:shorthand-for-oriented-simplices}. 

	\columnbreak

	\item[$\boundary_n$]
	usu.\ refers to a simplicial boundary map $\boundary_n: C_n(K;R) \to C_{n-1}(K;R)$, see \fref{defn:boundary-map};
	also used to denote differentials of an arbitrary chain complex $C_\ast = (C_n, \boundary_n)_{n \in \ints}$.

	\item[$C_\ast(K;R)$] 
	the simplicial chain complex $C_\ast(K;R) = (C_n(K;R), \boundary_n)_{n \in \ints}$ with boundary maps $\boundary_n:C_n(K;R) \to $\\$C_{n-1}(K;R)$, see \fref{defn:simplicial-chain-complex}.

	\item[$C_\ast(K)$]
	shorthand for $C_\ast(K;\ints)$.

	\item[$C_\ast(-;R)$] 
	the simplicial chain complex functor 
	$\catsimp \to \catchaincomplex{\catmod{R}}$, see \fref{defn:simplicial-chain-complex-functor}.

	\item[$H_n(K;R)$] 
	the $n$\th simplicial homology group of a simplicial complex $K$ with coeff.\ in $R$, see \fref{defn:simplicial-chain-complex}.

	\item[$H_n(K)$] 
	shorthand for $H_n(K;\ints)$.

	\item[$H_n(-;R)$] 
	usu.\ refers to the $n$\th simplicial homology functor 
	$\catsimp \to \catmod{R}$ with coeff.\ in $R$, 
	see \fref{defn:simplicial-homology-functor}.

	\item[$\betti_n(K;R)$]
	the $n$\th Betti number of a simplicial complex $K$ with coefficients in a PID $R$, i.e.\ 
	$\betti_n(K;R) = \rank(H_n(K;R))$.

	\item[$\betti_n(K)$] 
	shorthand for $\betti_n(K;\ints)$.

	\item[${[\sigma]}$] 
	usu.\ refers to a homology class in $H_n(K;R)$ 
	with the $n$-chain $\sigma$ an $R$-formal sum of oriented $n$-simplices, typically written in string notation.

\end{symbollist}\end{multicols}
\vspace{-14pt}

\clearpage

\noindent
\textbf{Notation involving Persistence Modules.}
See \fref{chapter:persistence-theory}.
\vspace{-7pt}
\begin{multicols}{2}\begin{symbollist}
	\item[$(V_\bullet, \alpha_\bullet)$]
	a persistence module with vector spaces $V_t$ and structure maps $\alpha_{s,t}: V_t \to V_s$, 
	see \fref{defn:persmod}.

	\item[$\persmod{V}$]
	shorthand for $(\persmod{V}, \alpha_\bullet)$,
	see \fref{defn:persmod}.

	\item[$\phi_\bullet$] 
	usu.\ a persistence morphism $\phi_\bullet: V_\bullet \to W_\bullet$ with $\phi_\bullet = (\phi_t)_{t \in \nonnegints}$, see \fref{defn:persmod-cat}(ii).

	\item[$(V_\ast^\bullet, \alpha_\ast^\bullet, \boundary_\ast^\sbullet)$] 
	a persistence complex, i.e.\ a chain complex of persistence modules $(V_n^\bullet, \alpha_n^\bullet)$ with differentials 
	$\boundary_n^\sbullet: V_n^\bullet \to V_{n-1}^\bullet$ over $n \in \ints$;
	see \fref{defn:category-of-persistence-chains}.

	\item[$(V_\ast^\bullet, \boundary_\ast^\sbullet)$]
	shorthand for $(V_\ast^\bullet, \alpha_\ast^\sbullet, \boundary_\ast^\sbullet)$.

	\item[$\intmod{J}$]
	a $J$-interval (persistence) module with $J \subseteq \nonnegints$ an interval,
	see \fref{defn:persistence-interval-modules}.

	\item[$\barcode(V_\bullet)$]
	the persistence barcode of a persistence module $V_\bullet$, 
	see \fref{defn:persistence-barcode}.

	\item[$\togrmod$]
	the functor $\catpersmod \to \catgradedmod{\field}$
	in the isomorphism of categories bet.\ $\catpersmod$ and $\catgradedmod{\field}$,
	see \fref{defn:togrmod}.

	\item[$\topersmod$]
	the functor $\catgradedmod{\field} \to \catpersmod$
	in the isomorphism of categories bet.\ $\catpersmod$ and $\catgradedmod{\field}$,
	see \fref{defn:topersmod}.
\end{symbollist}\end{multicols}

\vspace{3pt}
\noindent
\textbf{Notation involving Filtrations and Persistent Homology.}
See \fref{chapter:filtrations-and-pershoms}.
\vspace{-7pt}
\begin{multicols}{2}\begin{symbollist}
	\item[$\filt{K}$]
	a filtration $\filt{K}: \posetN \to \catsimp$
	of a simplicial complex $K$, 
	see \fref{defn:filtration}.

	\item[$C_n(\filt{K}; \field)$]
	the $n$\th filtered chain module of a filtration $\filt{K}$ with coefficients in $\field$, 
	see \fref{defn:filtered-n-chain-module}.

	\item[$C_n\graded(\filt{K}; \field)$]
	the $n$\th graded chain module of a filtration $\filt{K}$ with coefficients in $\field$, 
	see \fref{defn:graded-chain-module}.

	\item[${\basis{K}_n\graded}$]
	the standard ordered basis of $C_n\graded(K_\bullet;\field)$ induced by the orientation on $K$,
	see \fref{defn:standard-basis-graded-chain}.

	\item[$\boundary_n^{\sbullet}$]
	the $n$\th filtered boundary morphism 
	$\boundary_n^\sbullet: C_n(\filt{K};\field) \to C_{n-1}(\filt{K})$, 
	see \fref{defn:filtered-n-boundary-morphism}.

	\item[${\boundary_n\graded}$]
	the $n$\th graded boundary map
	$\boundary_n\graded: C_n\graded(K_\bullet;\field) \to C_{n-1}\graded(K_\bullet;\field)$, 
	see \fref{defn:graded-n-boundary-map}.

	\columnbreak 

	\item[$C_\ast(K_\bullet;\field)$] 
	the simplicial persistence complex of $K_\bullet$
	with $C_\ast(K_\bullet;\field) = (C_n(K_\bullet;\field), \boundary_n^\sbullet)_{n \in \ints}$, see \fref{defn:simplicial-persistence-complex}.

	\item[$C_\ast\graded(K_\bullet;\field)$] 
	the simp.\ graded chain complex of $K_\bullet$ 
	with $C_\ast\graded(K_\bullet;\field) = (C_n\graded(K_\bullet;\field), \boundary_n\graded)_{n \in \ints}$, 
	see \fref{defn:simplicial-persistence-complex}

	\item[$H_n(\filt{K}; \field)$]
	the $n$\th persistent homology module of $\filt{K}$ with coeff.\ in $\field$, 
	see \fref{defn:n-persistent-homology}.

	\item[$H_n\graded(K_\bullet;\field)$] 
	the $n$\th graded homology module of a filtration $\filt{K}$, see \fref{defn:graded-homology-module}.

	\item[${[-]_t}$]
	a homology class in a persistent homology module at 
	scale $t \in \nonnegints$,
	see Defn.\ \ref{defn:n-persistent-homology}.
	
	\item[$\barcode_n(\filt{K}; \field)$]
	the $n$\th persistence barcode of a filtration $\filt{K}$ with coefficients in $\field$, 
	see Defn.\ \ref{defn:n-persistent-homology}.

	\item[$H_n(K_t;p;\field)$]
	the $p$-persistent $n$\th homology group 
	with coefficients in $\field$, 
	see Defn.\ \ref{defn:p-persistent-bullshit}.

	\item[$\betti_n(K_t; p; \field)$]
	the $p$-persistent $n$\th Betti number with coefficients in $\field$, 
	see Defn.\ \ref{defn:p-persistent-bullshit}.

\end{symbollist}\end{multicols}

\noindent
\textbf{Notation involving Matrices and Smith Normal Decompositions (SNDs).} See \fref{appendix:matrix-theory}.
\vspace{-7pt}
\begin{multicols}{2}\begin{symbollist}
	\item[${[v]_S}$] 
	the coordinate matrix of $v \in M$ relative to a basis $S$ of some free $R$-module $M$, see \fref{defn:coordinate-matrices}. 

	\item[${[\phi]_{A,S}}$]
	the matrix of a module homomorphism $\phi: M \to N$ relative to bases $S$ and $A$ of $R$-modules $M$ and $N$ respectively.

	\item[$(U,D,V)$] 
	usu.\ refers to an SND of some matrix $A$ such that 
	$U\inv AV = D$, see Defn.\ \ref{defn:smith-normal-decomposition}.

	\columnbreak

	\item[$I_n$] 
	the identity matrix in $\GL(n,R)$.

	\item[${\elswap[n]{-}}$]
	an elementary permutation of degree $n$ over some ring $R$,
	see \fref{defn:elementary-matrices}.

	\item[${\eldilate[n]{-}}$]
	an elementary dilation of degree $n$ over some ring $R$,
	see \fref{defn:elementary-matrices}.

	\item[${\eladd{-}}$]
	an elementary transvection of degree $n$ over some ring $R$,
	see \fref{defn:elementary-matrices}. 

\end{symbollist}\end{multicols}
\vspace{-14pt}

} 
\clearpage

\section{Matrices over Euclidean Domains}
\label{appendix:matrix-theory}

In this paper, we discuss algorithms over matrices over some Euclidean domain $R$, i.e.\ matrices whose entries are elements of $R$. Note that Euclidean domains have multiplicative identities, typically denoted by $1$ or $1_R$, which are unique.
Euclidean domains are also commutative rings, i.e.\ the multiplication operation is commutative.

Note that most definitions involving matrices over $R$ are generally the same as that over $\reals$.
For example, since $R$ has identity $1$, identity matrices and invertible matrices are defined in the same way as that of $\GL(n,\reals)$.
Below, we identify notation involving families of matrices used in this paper.

\begin{definition}
	Let $R$ be a Euclidean domain.
	Let $\M_{m,n}(R)$ refer to the \textbf{collection of all $(m \times n)$-matrices}, i.e.\ having $m$ rows and $n$ columns, \textbf{over $R$}.
	If $m=n$, we write $M_n(R)$ to refer to $M_{n,n}(R)$.
	Let $\GL(n,R)$ refer to the \textbf{general linear group of degree $n$ over $R$}, i.e.\ 
	the {collection of all invertible $(n\times n)$-matrices over $R$}.
\end{definition}

Next, we identify notation for certain components of some given matrix.

\begin{definition}
	Let $A \in M_{m,n}(R)$.
	Let $A(j,i)$ denote the \textbf{entry} of $A$ at the $j$\th row and $i$\th column.

		If $n=1$, 
		we call $A \in M_{m,1}(R)$ a \textbf{column vector} and write $A(j) := A(j,1)$ for the $j$\th entry of $A$.
		If $m=1$, we call $A \in M_{1,n}(R)$ a 
		\textbf{row vector} and write $A(i) := A(1,i)$ for the $i$\th entry of $A$.

		Let $\col_i(A) \in M_{m,1}(R)$ refer to the column vector corresponding to the 
		\textbf{$i$\th column of $A \in M_{m,n}(R)$}, 
		i.e.\ for all $j \in \set{1, \ldots, m}$, 
			$\col_i(A)(j) = A(j,i)$ for fixed $i$.
		Let $\row_j(A) \in M_{1,n}(R)$ refer to the row vector corresponding to the 
		\textbf{$j$\th row of $A \in M_{m,n}(R)$}, 
		i.e.\ 
			for all $i \in \set{1, \ldots, n}$, 
			$\row_j(A)(i) = A(j,i)$ for fixed $j$.
\end{definition}

Note that indices of the rows and columns of matrices always start at $1$. We emphasize this since, for some objects referenced in this paper, e.g.\ persistence modules, the indexing starts at $0$.
In this paper, we make an effort to use $j$ and $i$ to refer to an index of some row and column respectively but this is sometimes not possible, e.g.\ in cases where $i$ or $j$ have been defined beforehand.

We also talk about diagonal matrices in this paper. For brevity, we identify notation for describing diagonal matrices by the elements on their diagonal.

\begin{definition}\label{defn:diagonal-matrix}
	Let $D = \diag(d_1, d_2, \ldots, d_n)$ with $d_i \in R$ for all $i \in \set{1, \ldots, n}$ refer to matrix $D \in M_n(R)$ given by 
	$D(i,j) = d_i$ if $i = j$ and $D(i,j) = 0$ otherwise.
\end{definition}

We also talk about block matrices in this paper. We identify the notation used for these below.

\begin{definition}
	Let $A \in M_{m_1, n_1}(R)$, 
		$B \in M_{m_1, n_2}(R)$,
		$C \in M_{m_2, n_1}(R)$,
		$D \in M_{m_2, n_2}(R)$.
	We write  
	\begin{equation*}
		X = \begin{pmatrix}
			A & B \\
			C & D
		\end{pmatrix}
	\end{equation*}
	to describe the matrix $X \in M_{m, n}(R)$ with $m = m_1 + m_2$ and $n = n_1 + n_2$.
	We call $X$ a \textbf{block matrix} and the matrices $A,B,C,D$ the \textbf{blocks} of $X$.
	When a block is given by a zero matrix, we usually do not write the dimensions of the zero matrix and assume its dimensions are appropriately defined.
\end{definition}


Listed below are three families of matrices that are generalizations of the \textit{elementary matrices} over $\reals$, i.e.\ those used in matrix reduction and are usually called Type (I), Type (II), and Type (III) matrices.
Below, we state a definition for these families of matrices, taken from~\cite[Definition 4.1.8]{algebra:adkins}. 
Note that interpretation of the left or right multiplication of these matrices are given later in \fref{prop:row-operations} (row operations) and \fref{prop:column-operations} (column operations).

\begin{definition}\label{defn:elementary-matrices}
	The \textbf{elementary matrices over $R$ of degree $n$} consists of 
	elementary permutations, elementary dilations, and elementary transvections, defined below.
	Note that all indices are elements of $\set{1, \ldots, n}$.
	\begin{enumerate}
		\item 
			An \textbf{elementary permutation} $\elswap{k_1, k_2} \in \GL(n,R)$ on two indices $k_1,k_2 \in \set{1, \ldots, n}$ is the matrix obtained by interchanging rows $k_1$ and row $k_2$ 
			(or equivalently, columns $k_1$ and $k_2$) of the identity matrix $I_n$, i.e.\ 
			\begin{equation*}
				\row_j\paren{\elswap{k_1, k_2}} = \begin{cases}
					\row_{k_2}({I_n})	&\text{ if } j = k_1 \\
					\row_{k_1}({I_n})	&\text{ if } j = k_2 \\
					\row_{j}({I_n}) 	&\text{ otherwise } \\
				\end{cases}
				\,\,\text{ or }\,\,
				\col_{i}\paren{\elswap{k_1, k_2}} = \begin{cases}
					\col_{k_2}({I_n})	&\text{ if } i = k_1 \\
					\col_{k_1}({I_n})	&\text{ if } i = k_2 \\
					\col_{i}({I_n}) 	&\text{ otherwise } \\
				\end{cases}
			\end{equation*}
			Elementary permutations are involutions, i.e.\ the inverse of $\elswap{k_1, k_2}$ is $\elswap{k_1, k_2}$.
			
		\item 
			An \textbf{elementary transvection}
			$\eladd{k_j, k_i \,; \alpha} \in \GL(n,R)$ by $\alpha \in R$ on the $(k_j, k_i)$\th entry with $k_j,k_i \in \set{1, \ldots, n}$ and $k_j \neq k_i$ is the matrix obtained by taking the identity matrix and 
			replacing the $(k_j, k_i)$\th entry with $\alpha$, i.e.\ 
			\begin{equation*}
				\eladd{k_j, k_i \,; \alpha}(j,i) = \begin{cases}
					\alpha 		&\text{ if } j = k_j \text{ and } i = k_i \\
					I_n(j,i)	&\text{ otherwise }
				\end{cases}
			\end{equation*}
			We may call $\alpha$ the \textbf{transvection multiplier} of $\eladd{k_j, k_i \,; \alpha}$.
			The inverse of an elementary transvection 
			$\eladd{k_j, k_i; \alpha}$ is the transvection 
			$\eladd{k_j, k_i, -\alpha}$ where $-\alpha$ is the additive inverse of $\alpha$ in $R$.

		\item 
			An \textbf{elementary dilation} $\eldilate{k, \mu} \in \GL(n,R)$ on index $k \in \set{1, \ldots, n}$ by a unit $\mu \in R\units$ is the matrix obtained by replacing the $k$\th diagonal element of $I_n$ with $\mu$, i.e.\ 
			\begin{equation*}
				\eldilate{k, \mu}(j,i) = \begin{cases}
					\mu 		&\text{ if } j=k \text{ and } i=k \\
					I_n(j,i)	&\text{ otherwise }
				\end{cases}
			\end{equation*}
			We sometimes call $\mu$ the \textbf{dilation multiplier} of $\eldilate{k,\mu}$.
			The inverse of an elementary dilation $\eldilate{k, \mu}$ 
			is the dilation $\eldilate{k, \mu\inv}$ with $\mu\inv$ the multiplicative inverse of $\mu$ in $R$.

	\end{enumerate}
\end{definition}

One reason why elementary matrices are considered significant is due to the following result for matrices over Euclidean domains, stated below.

\begin{proposition}
	Any invertible matrix over a Euclidean domain $R$ can be expressed as a finite product of elementary matrices.
	That is, $\GL(n,R)$ is generated by the elementary matrices over $R$ of degree $n$.
\end{proposition}
\remark{
	For a proof, see under \cite[Theorem 5.2.10]{algebra:adkins}.
	As a sidenote, 
}

Observe that, given two Euclidean domains $R_1$ and $R_2$, the family of elementary permutations on $R_1$ and $R_2$ are defined very similarly
since only the multiplicative identities are used to generate them and these identities are unique for each ring.
For example, the elementary permutations of matrices over $\ints$ are exactly the elementary permutations of matrices of $\reals$ since they share identity elements.

We provide a visual description of these elementary permutations below.
Observe that the order of the arguments $k_1$ and $k_2$ in $\elswap{k_1,k_2}$ does not matter, i.e.\ $\elswap{k_1,k_2} = \elswap{k_2,k_1}$.

\vspace{-5pt}
\begin{equation*}
	\elswap{\greenmath{k_1\!}, \greenmath{k_2\!}}
	= \hspace{5pt}
	\begin{NiceArray}{>{\color{black}\scriptstyle}c @{\hspace{13pt}}
		ccc c ccc @{\hspace{5pt}} l
	}
	\CodeBefore [create-cell-nodes]
		\tikz \node [green-highlight = (1-4), inner xsep=2pt] {} ;
		\tikz \node [green-highlight = (1-6), inner xsep=2pt] {} ;
		\tikz \node [green-highlight = (4-1), inner xsep=2pt] {} ;
		\tikz \node [green-highlight = (6-1), inner xsep=2pt] {} ;
		\tikz \node [red-cell = (4-6), inner xsep=4pt] {} ;
		\tikz \node [red-cell = (6-4), inner xsep=4pt] {} ;
		\tikz \node [orange-highlight = (4-4), inner xsep=4pt] {} ;
		\tikz \node [orange-highlight = (6-6), inner xsep=4pt] {} ;
	\Body
		\RowStyle[nb-rows=1]{\color{black}\scriptstyle}
		{} & {\mathllap{i=}1} & {\cdots} & {k_1\!} & {\cdots} & {k_2\!} & {\cdots} & {n} 
			& \mathrlap{\triangleleft \text{ column indices}}
			\\[2pt]
		{j=1} 	& 1 & \cdots & 0 & \cdots & 0 & \cdots & 0 \\[-3pt]
		{\vdots} 	& \vdots & \ddots & \vdots &  & \vdots &  & \vdots \\[1pt]
		{k_1\!} 	& 0 	 &  & 0 &  & 1 &  & 0 \\[-1pt]
		{\vdots} 	& \vdots &  & \vdots & \ddots & \vdots &  & \vdots \\[1pt]
		{k_2\!} 	& 0 	 &  & {1} &  & 0 &  & 0 \\
		{\vdots} 	& \vdots &  & \vdots &  & \vdots & \ddots & \vdots \\
		{n} 		& 0 	 & \cdots & 0 & \cdots & 0 & \cdots & 1 \\
		\mathclap{\overset{\triangle}{\text{row indices}}}
	\CodeAfter
		\SubMatrix({2-2}{8-8})[left-xshift=3pt, right-xshift=3pt]
	\end{NiceArray}
	\hspace{30pt}
	\begin{gathered}
		\text{\small Color Scheme:} \\
		\greentagged{\small indices $k_1$ and $k_2$} \\
		\redtagged{\small entries in $I_n$ set to $1$} \\
		\orangetagged{\small entries in $I_n$ set to $0$}
		\\[0pt]
		\text{\small with uncolored entries } \\[-5pt]
		\text{\small agreeing with $I_n$ }
	\end{gathered}
\end{equation*} 
\vspace{-2pt}

Elementary transvections are relatively straightforward to work with since they allow any ring element to be the transvection multiplier.
Note that if $\alpha = 0$ is the additive identity of $R$, then $\eladd{k_j, k_i \,; \alpha} = I_n \in \GL(n,R)$.
Also, note that the order of arguments in the notation $\eladd{k_j, k_i \,; \alpha}$ is important, unlike for elementary permutations, i.e.\ $\eladd{k_j, k_i \,; \alpha} \neq \eladd{k_i, k_j \,; \alpha}$.

Provided below is a visual description of $\eladd{k_\text{row}, k_\text{col} \,; \alpha}$. 
As the labels suggest, the first two arguments $k_\text{row}, k_\text{col}$ determine that the $(k_\text{row}, k_\text{col})$\th entry of $\eladd{k_\text{row}, k_\text{col} \,; \alpha}$ is multiplier $\alpha$.

\vspace{-15pt}
\begin{equation*}
	\eladd{
		\greenmath{k_\text{row}}, 
		\bluemath{k_\text{col}} \,;
		\redmath{\rule[-1pt]{0pt}{9pt}\alpha}
	}
	= \hspace{5pt}
	\begin{NiceArray}{>{\scriptstyle}c @{\hspace{13pt}}
		ccc c ccc @{\hspace{5pt}} l
	}
	\CodeBefore [create-cell-nodes]
		\tikz \node [green-highlight = (1-4), inner xsep=2pt] {} ;
		\tikz \node [blue-highlight = (1-6), inner xsep=2pt] {} ;
		\tikz \node [green-highlight = (4-1), inner xsep=2pt] {} ;
		\tikz \node [blue-highlight = (6-1), inner xsep=2pt] {} ;
		\tikz \node [red-cell = (4-6), inner xsep=4pt] {} ; 
	\Body
		\RowStyle[nb-rows=1]{\color{black}\scriptstyle}
		{} & {\mathllap{i=}1} & {\cdots} & {k_\text{row}} & {\cdots} & {k_\text{col}} & {\cdots} & {n} 
			& \mathrlap{\triangleleft \text{ column indices}}
			\\[2pt]
		{j=1} 	& 1 	 & \cdots & 0 & \cdots & 0 & \cdots & 0 \\[-2pt]
		{\vdots} 			& \vdots & \ddots & \vdots &  & \vdots &  & \vdots \\[1pt]
		{k_\text{row}} 				& 0 	 &  & 1 &  & \alpha &  & 0 \\
		{\vdots} 			& \vdots &  & \vdots & \ddots & \vdots &  & \vdots \\[1pt]
		{k_\text{col}} 				& 0 	 &  & {0} &  & 1 &  & 0 
			\\
		{\vdots} 			& \vdots &  & \vdots &  & \vdots & \ddots & \vdots \\
		{n} 				& 0 	 & \cdots & 0 & \cdots & 0 & \cdots & 1 \\
		\mathclap{\overset{\triangle}{\text{row indices}}}
	\CodeAfter
		\SubMatrix({2-2}{8-8})[left-xshift=3pt, right-xshift=3pt]
	\end{NiceArray}
	\hspace{30pt}
	\begin{gathered}
		\text{\small Color Scheme:} \\
		\greentagged{\small row index $k_\text{row}$} \\
		\bluetagged{\small column index $k_\text{column}$} \\
		\redtagged{\small transvection multiplier $\alpha$} 
		\\[0pt]
		\text{\small with uncolored entries } \\[-5pt]
		\text{\small agreeing with $I_n$ }
	\end{gathered}
\end{equation*} 
\vspace{-7pt}

An elementary transvection $\eladd{k_\text{row}, k_\text{col} \,; \alpha}$ can also be described row-wise as follows. Note that only row $k_\text{row}$ of $\eladd{k_\text{row}, k_\text{col} \,; \alpha}$ differs from the identity matrix $I_n$.

\begin{equation*}
	\row_{j}\paren{\eladd{k_\text{row}, k_\text{col} \,; \alpha}} = \begin{cases}
		\row_{j}{(I_n)} + \alpha \cdot \row_{k_\text{col}}{(I_n)}
			&\text{ if } j = k_\text{row} \\
		\row_{j}({I_n})	&\text{ otherwise }
	\end{cases}
\end{equation*}
It also has the following column-wise description. 
Note that only column $k_\text{col}$ of $E := \eladd{k_\text{row}, k_\text{col} \,; \alpha}$ differs from $I_n$.
\begin{equation*}
	\col_{i}\paren{\eladd{k_\text{row}, k_\text{col} \,; \alpha}} = \begin{cases}
		\col_{i}(I_n) + \alpha \cdot \col_{k_\text{row}}(I_n)
			&\text{ if } i = k_\text{col} \\
		\col_{i}{(I_n)}	&\text{ otherwise }
	\end{cases}
\end{equation*}

\spacer

Let $R$ be a commutative ring with identity and $S$ be a subring of $R$.
Since $S$ and $R$ must share an identity element, the elementary dilations over $R$ are exactly those over $S$.
Given an elementary dilation $\eladd{k_j, k_i \,; r}$ over $R$ with $r \in R$,
	it should be clear that if $r \in S$, then $\eladd{k_j, k_i \,; r}$ is also an elementary dilation over $S$.
	More specifically, $\eladd{k_j,k_i \,; r} \in \GL(n,S)$.

However, we do not have this nice relationship for elementary dilations. That is,
an elementary dilation $\eldilate{k,r}$ over $R$ with $r \in R^\times$ is not generally an elementary dilation over $S$ even if $r \in S$.
In particular, $\eldilate{k,r}$ is that over $S$ only if $r$ is invertible in $S$.
This distinction becomes extremely significant when considering which elementary matrices we can use for matrix reduction.
We consider some relevant cases below.
\begin{enumerate}
	\item 
	If $R = \field$ is a field, then any nonzero element of $\field$ can be used as the dilation multiplier $\mu$ in $\eldilate{k,\mu}$.
	For example, row reduction of a matrix with $\reals$ entries allows the multiplication of any entry by a nonzero constant.
	This is because $\eldilate{k,a}$ is a valid elementary dilation over $\reals$ for any nonzero $a \in \reals$.

	\item 
	In the case of $R = \ints$, 
	the family of elementary dilations over $\ints$ is much smaller than that of $\reals$, despite $\ints \subseteq \reals$, since the only units of $\ints$ are $1$ and $-1$.
	So, when performing row reduction on matrices over $\ints$, we can only multiply each row by either $1$ (which does nothing) or $-1$.
	
	As another example, let $M = \eldilate[3]{1,2}$ with multiplier $\mu = 2$. 
	Observe that $M$ is a valid elementary dilation over $\reals$ since $\mu=2$ has an inverse of $\frac{1}{2} \in \reals$.
	Given below is a more visual description of $M$ and $M\inv = \eldilate[3]{1,\frac{1}{2}}$.
	\begin{equation*}
		M = \eldilate[3]{1,2} = \begin{pmatrix}
			2 & 0 & 0 \\
			0 & 1 & 0 \\
			0 & 0 & 1
		\end{pmatrix}
		\qquad\text{ and }\qquad 
		M\inv = \eldilate[3]{1,\frac{1}{2}} = \begin{pmatrix}
			\frac{1}{2} & 0 & 0 \\
			0 & 1 & 0 \\
			0 & 0 & 1
		\end{pmatrix}
	\end{equation*}
	Since $\frac{1}{2} \not\in \ints$, $M\inv \not\in \M_{3,3}(\ints)$ and $M\inv$ is not a valid elementary dilation over $\ints$.

	\item 
	In the case of $R = \rationals[x]$,
	the elementary dilations over $\rationals$ and those of $\rationals[x]$ are exactly the same since they share the same group of units, i.e.\ $\rationals^\times = (\rationals[x])^\times = \rationals \setminus \set{0}$.

	So, when performing row reduction on a matrix over $\rationals[x]$, we cannot use an elementary dilation to reduce the powers since $x^t \not\in \rationals[x]$ for $t \in \ints$ with $t < 0$.
	For example, let $M \in \M_{3,2}(\rationals[x])$ be as given below.
	\begin{equation*}
		M = \begin{pmatrix}
			x^3 & 0 & 0 \\
			0 & x & 1 \\
		\end{pmatrix}
	\end{equation*}
	The following row operation on $M$ is not valid as a matrix operation over $\rationals[x]$:
	\begin{equation*}
		\eldilate{1,x^{-3}}M 
		=
		\begin{pmatrix}
			x^{-3} & 0 \\
			0 & 1
		\end{pmatrix}\begin{pmatrix}
			x^3 & 0 & 0 \\
			0 & x & 1 \\
		\end{pmatrix}
		= 
		\begin{pmatrix}
			1 & 0 & 0 \\
			0 & x & 1
		\end{pmatrix}
	\end{equation*}
	In particular, $\eldilate{1,x^{-3}} \not\in \M_{2,2}(\rationals[x])$
	and we cannot multiply the first row of $M$ by $x^{-3}$.
\end{enumerate}
We provide a visual description of elementary dilations below.

\vspace{-5pt}
\begin{equation*}
	\eldilate{\greenmath{k}, \redmath{\rule{0pt}{6pt}\mu}}
	= \hspace{5pt}
	\begin{NiceArray}{>{\scriptstyle}c @{\hspace{13pt}}
		ccc c ccc @{\hspace{5pt}} l
	}
	\CodeBefore [create-cell-nodes]
		\tikz \node [red-cell = (5-5), inner xsep=4pt] {} ;
		\tikz \node [green-highlight = (1-5), inner xsep=2pt] {} ;
		\tikz \node [green-highlight = (5-1), inner xsep=2pt] {} ;
	\Body
		\RowStyle[nb-rows=1]{\scriptstyle}
		{} & {\mathllap{i=}1} & {\cdots} & {k-1} & {\color{black} k} & {k+1} & {\cdots} & {n} 
			& \mathrlap{\triangleleft \text{ column indices}}
			\\[2pt]
		{j=1} & 1 & \cdots & 0 & 0 & 0 & \cdots & 0 \\[-3pt]
		{\vdots} & \vdots & \ddots & \vdots & \vdots & \vdots & {} & \vdots \\
		{k-1} & 0 & {} & 1 & 0 & 0 & {} & 0 \\
		{\color{black} k} & 0 & {} & 0 & {\rule{0pt}{5pt}\mu} & 0 & {} & 0 \\[2pt]
		{k+1} & 0 & {} & 0 & 0 & 1 & {} & 0 \\
		{\vdots} & \vdots & {} & \vdots & \vdots & \vdots & \ddots & \vdots \\
		{n} & 0 & \cdots & 0 & 0 & 0 & \cdots & 1 \\
		\mathclap{\overset{\triangle}{\text{row indices}}}
	\CodeAfter
		\SubMatrix({2-2}{8-8})[left-xshift=3pt, right-xshift=3pt]
	\end{NiceArray}
	\hspace{30pt}
	\begin{gathered}
		\text{\small Color Scheme:} \\
		\greentagged{\small row or column index} \\
		\redtagged{\small dilation multiplier $\mu$} 
		\\[0pt]
		\text{\small with uncolored entries } \\[-5pt]
		\text{\small agreeing with $I_n$ }
	\end{gathered}
\end{equation*} 
\vspace{5pt}


\noindent
Below, we provide a characterization of row and column operations on matrices relative to these elementary matrices.

\begin{proposition}\label{prop:row-operations}
	\textbf{Elementary Row Operations}
	on a matrix ${A} \in \M_{m,n}(R)$ correspond to multiplying $A$ by some elementary matrix in $\GL(m,R)$ on the left. In particular:
	\begin{enumerate}[itemsep=-\parsep]
		\item 
		\textbf{Row Permutation or Row Swapping.} 
		The product $\elswap[m]{k_1, k_2}{A}$ corresponds to swapping the $k_1$\th and $k_2$\th rows of ${A}$, i.e.\ 
		\begin{equation*}
			\row_j\paren{\elswap{k_1, k_2}A} = 
			\begin{cases}
				\row_{k_2}(A)		&\text{ if } j = k_1 \\
				\row_{k_1}(A)		&\text{ if } j = k_2 \\
				\row_j(A)			&\text{ otherwise }
			\end{cases}
		\end{equation*}

		\item  
		\textbf{Row Dilation or Row Multiplication.}
		The product $\eldilate[m]{k,\mu}{A}$ corresponds to multiplying the $k$\th row by a unit $\mu \in R^\times$, 
		i.e.\ 
		\begin{equation*}
			\row_j\paren{\eldilate{k,\mu}{A}} = 
			\begin{cases}
				\mu \row_j({A})	&\text{ if } j = k \\
				\row_j({A})				&\text{ otherwise }
			\end{cases}
		\end{equation*}

		\item 
		\textbf{Row Addition.}
		With $\alpha \in R$ and column indices $k \neq p$,
		the product $\eladd[m]{k, p, \alpha}{A}$ corresponds to adding an $\alpha$-multiple of the $p$\th row of ${A}$ to the $k$\th row of ${A}$,
		i.e.\ 
		\begin{equation*}
			\row_j\paren{\eladd{k,p \,; \alpha}{A}} = 
			\begin{cases}
				\row_k({A})	+ \alpha \row_p(A)
					&\text{ if } j=k \\
				\row_j({A})	&\text{ otherwise }
			\end{cases}
		\end{equation*}
		Note that the first argument $k$ in $\eladd[m]{k, p, \alpha}{A}$ is the target row.

	\end{enumerate}
\end{proposition}





\negativespacer

\begin{proposition}\label{prop:column-operations}
	\textbf{Elementary Column Operations} on a matrix $A \in M_{m,n}(R)$ correspond multiplying ${A}$ by some elementary matrix in $\GL(n,R)$ on the \textit{left}. In particular:
	\begin{enumerate}
		\item \textbf{Column Permutation or Column Swapping.} 
		The product ${A}\elswap{k_1, k_2}$ corresponds to swapping the $k_1$\th and $k_2$\th column of ${A}$, i.e.\ 
		\begin{equation*}
			\col_i\paren{A\elswap{k_1, k_2}} = 
			\begin{cases}
				\col_{k_2}(A)		&\text{ if } i = k_1 \\
				\col_{k_1}(A)		&\text{ if } i = k_2 \\
				\col_i(A)			&\text{ otherwise }
			\end{cases}
		\end{equation*}

		\item \textbf{Column Dilation or Column Multiplication.} 
		The product ${A}\eldilate{k, \mu}$ corresponds to multiplying the $k$\th column of ${A}$ by a unit $\mu \in R^\times$,
		i.e.\ 
		\begin{equation*}
			\col_i\paren{A\eldilate{k,\mu}} = 
			\begin{cases}
				\mu \col_k({A})			&\text{ if } i = k \\
				\col_i({A})				&\text{ otherwise }
			\end{cases}
		\end{equation*}

		\item \textbf{Column Addition.} 
		With $\alpha \in R$ and column indices $k \neq p$,
		the product $\eladd{p,k \,; \alpha}$ corresponds to adding an $\alpha$-multiple of the $p$\th column of $A$ to the $k$\th column of $A$, i.e.\ 
		\begin{equation*}
			\col_i\paren{A\eladd{p,k \,; \alpha}} = 
			\begin{cases}
				\col_k({A})	
					+ \alpha \col_p(A)
					&\text{ if } i = k \\
				\col_i({A})			&\text{ otherwise }
			\end{cases}
		\end{equation*}
		Note that the second argument $k$ in $\eladd{p,k \,; \alpha}$ is the target column.
	\end{enumerate}
\end{proposition}

\clearpage

\section{Notes on Ring and Module Theory}
\label{appendix:module-theory}


For reference, we have listed some relevant basic definitions and notation involving rings and modules.
These are mostly taken from the texts \textit{Algebra: An Approach to Module Theory} \cite{algebra:adkins} by William Adkins and Steven Weintraub, and 
\textit{Abstract Algebra} \cite{algebra:dummit} by David Dummit and Richard Foote.

\begin{definition}
	Let $R$ be a commutative ring.
	An \textbf{ideal} $S$ of a ring $R$ is a subring of $R$ such that $RS \subseteq S$ and $SR \subseteq S$, i.e.\ the subring $S$ is closed under left and right multiplication of elements of $R$. 
	An {ideal} $S \subseteq R$ is called \textbf{principal} if $S$ can be generated by a single element, i.e.\ there exists $a \in R$ such that $S = Ra = \set{ra : r \in R}$. In this case, we write $S = (a)$ and say $S$ is \textbf{generated} by $a \in R$.
\end{definition}
\remark{
	In some references, $S = (a)$ is written as $S = \ket{a}$, i.e.\ with angle brackets. We avoid using this notation in this paper for clarity, e.g.\ in cases where torsion is present.
	Instead, we use angle brackets as described in \fref{defn:formal-sums}, i.e.\ $\ket{a}$ is usually interpreted as $\ints\ket{a}$.
}

Below, we provide a definition for PIDs.
Note that, in this paper, we almost always require that the ring $R$ be at least a PID.

\begin{definition}
	An \textbf{integral domain} is a nonzero commutative ring $R$ with (multiplicative) identity $1_R \in R$ 
	such that product of any two nonzero elements of $R$ is nonzero.
	A \textbf{principal ideal domain} $R$ is an integral domain such that every ideal of $R$ is principal.
\end{definition}
\remark{
	Note that some references do not necessarily require integral domains to have identities.
	For contrast, an example of a commutative ring without identity is $2\ints$.
}

We also have references to \textit{Euclidean domains} in this paper but we have determined that an exact definition for such is not relevant.
It should suffice to know that all Euclidean domains are PIDs and are equipped with some sort of division algorithm, in which quotients and remainders are well-defined and unique.
Examples of Euclidean domains (and therefore, PIDs) include the integers $\ints$, all fields $\field$, and polynomial rings $\field[x]$ for any field $\field$.
Note that the polynomial ring $\ints[x]$ is not a PID. For example, the ideal $(2,x) = \set{2f + x g: f,g \in \ints[x]} \subseteq \ints[x]$ is not principal.

\spacer 

Next, we identify terminology involving properties of modules.

\begin{definition}
	An $R$-module $M$ is called \textbf{cyclic} if it can be generated by a single element $m \in M$, i.e.\ $Rm = \set{rm : r \in R} = M$.
\end{definition}

Note that $R$ can be viewed as an cyclic $R$-module generated by its identity element $1 \in R$.
Similarly, given $d \in R$, the quotient $R \bigmod (d)$ is a cyclic $R$-module with underlying abelian group $R \bigmod (d)$ and action $s \cdot [r] = [rs]$ generated by the coset $[1]$.

We also talk about notions of torsion and free involving modules over a PID. We provide definitions involving these below.

\begin{definition}
	Let $M$ be a module over a commutative ring $R$.
	A subset $B \subseteq M$ is called a \textbf{basis} of $M$ if each element of $M$ can be written as a unique $R$-linear combination of elements of $B$.
	An $R$-module $M$ is called \textbf{free} if it has a basis.
\end{definition}

\negativespacer

\HIDE{\remark{
	In \fref{section:graded-mod-notation}, we discussed graded modules. Since we only consider graded modules over $\field[x]$, the term ``\textit{free graded}'' simply means both graded and free in this paper. More generally, there is a notion of ``graded-free'' in \cite{algebra:nastasescu}, which is a stronger condition than just both free and graded.
}}

\begin{definition}
	Let $M$ be a module over a commutative ring $R$.
	We say that $m \in M$ is a \textbf{torsion element} of $M$ if there exists some nonzero $r \in R$ such that $rm = 0$.
	If $M$ has no torsion elements, then $M$ is called \textbf{torsion-free}.
	If every element of $M$ is torsion, then $M$ is called a \textbf{torsion module}.
\end{definition}

Modules over PIDs are considered ``well-behaved'' in that torsion and free can considered distinct notions in such modules. 
A rigorous discussion of how this works is outside the scope of this paper and we refer to \cite[Chapter 12]{algebra:dummit} and \cite[Chapter 3]{algebra:adkins} for more details. 
Listed below are a number of relevant results for reference.
\begin{enumerate}
	\item 
	Let $R$ be an integral domain and let $M$ be an $R$-module. 
	Then, the subset 
	$
		\text{T}(M) := \{ m \in M : rm = 0 \text{ for some nonzero } r \in R \} \subseteq M
	$ 
	is a submodule of $M$
	and the quotient $M \bigmod \text{T}(M)$ is torsion-free
	\cite[Proposition 3.2.18]{algebra:adkins}.

	\item 
	Let $R$ be an integral domain and let $M$ be an $R$-module.  
	If $M$ is free, then $M$ is torsion-free \cite[Proposition 3.4.8]{algebra:adkins}.
	Conversely:
	If $M$ is a finitely generated torsion-free module over a PID $R$, then $M$ is free \cite[Theorem 3.6.6]{algebra:adkins}.
	Therefore, a module over a PID are free if and only if it is torsion-free. 


	\item 
	Let $M$ be a module over a PID $R$.
	If $M$ finitely generated, then 
	$M \cong \text{T}(M) \oplus M \bigmod \text{T}(M)$ 
	\cite[Corollary 3.4.17]{algebra:adkins}.
	Since $M \bigmod \text{T}(M)$ is torsion-free and therefore free, 
		each element $m \in M$ decomposes uniquely into a ``free'' component and ``torsion'' component.

	\item 
	Let $R$ be a commutative ring with identity.
	If $M$ is a free $R$-module with a finite basis, then every basis of $M$ has the same number of elements
	\cite[Corollary 3.6.18]{algebra:adkins}. 
	This allows the notion of ``rank'' to be well-defined for free $R$-modules
	and implies that all free $R$-modules of the same rank are isomorphic.

	\item 
	If $M$ and $N$ are finitely generated modules over a PID $R$, then $M \cong N$ if and only if 
	$\text{T}(M) \cong \text{T}(N)$ and 
	$\rank(M \bigmod \text{T}(M)) = \rank(N \bigmod \text{T}(N))$ 
	\cite[Corollary 3.6.20]{algebra:adkins}.
	This determines that the ``free part'' of $M$, i.e.\ the summand $M \bigmod \text{T}(M)$ is unique up to isomorphism.
\end{enumerate}
This motivates the following terminology for finitely generated modules over a PID.

\begin{definition}
	Let $M$ be a finitely generated module over a PID $R$. 
	Define the \textbf{torsion component} $\text{T}(M)$ (also \textbf{torsion submodule})
	and the \textbf{free component} $\text{F}(M)$ of $M$ 
	to be submodules of $M$ given as follows:
	\begin{equation*}
		\text{T}(M) := \Bigl\{ m \in M : rm = 0 \text{ for some nonzero } r \in R \Bigr\} 
		\quad\text{ and }\quad 
		\text{F}(M) := M \bigmod \text{T}(M)
	\end{equation*}
	Let the \textbf{rank} of $M$, denoted $\rank(M)$, be the cardinality of any basis of $\text{F}(M)$.
\end{definition}

In \fref{section:matrix-calculation-of-IFDs} and \fref{section:snd-on-ungraded-chain-complexes}, we examine a proof of the Structure Theorem (\fref{thm:structure-theorem}) in $\catmod{R}$ and discuss how we can calculate invariant factor decompositions for finitely generated modules over a PID using presentations and matrix reduction.

In this expository paper, we usually define $R$-modules (particularly in examples) using formal sums of some set of indeterminates for convenience. 
We include relevant definitions below.

\begin{statement}{Definition}\label{defn:formal-sums}
	Let $R$ be a PID.
	Let $A = \set{a_1, a_2, \ldots, a_n}$ be some set of indeterminates.
	\begin{enumerate}
		\item 
		An $R$-\textbf{formal sum in $A$} is an expression in the form 
		$a = \sum_{i=1}^n r_ia_i = r_1 a_1 + r_2 a_2 + \cdots + r_n a^n$ for some $r_1, \ldots, r_n \in R$. When $R$ is unambiguous, we may refer to $\sum_{i=1}^n r_ia_i$ as a formal sum in $A$.
		We may also say that \textbf{$a$ is a formal sum of elements of $A$}.

		\item 
		The \textbf{$R$-module generated by $A$} is the free $R$-module with basis in correspondence with $A$.
		We usually write elements of $R\ket{A}$ as formal sums in $A$ with coefficients in $R$, i.e.\ 
		\begin{equation*}
			R\ket{A} := \set{\, 
				\sum_{i=1}^n r_i a_i
				\,:\,
				r_i \in R, a_i \in R \text{ for all } i =
				1, \ldots, n
			\,}
		\end{equation*}
		When an $R$-module $M$ is generated this way, we may write $M = R\ket{A}$ or $M = R\ket{a_1, \ldots, a_n}$.
	\end{enumerate}
\end{statement}
\remark{
	We somewhat abuse notation and also use $R\ket{-}$ to generate free submodules of $R\ket{A}$. For example, let $M = R\ket{a,b}$. Then, we may write $R\ket{2b}$ to refer to the submodule $R\ket{2b} = \set{r \cdot 2b: r \in R} = (2b)$, treating $2b$ as an element of $R\ket{b}$ rather than an indeterminate distinct from $b$. 
	In this case, we prefer to write $R\ket{2b}$ to emphasize that $R\ket{2b}$ is free.
}

Below, we identify notation for matrices related to free modules and homomorphisms between free modules, adapted from \cite[Section 4.3]{algebra:adkins}.

\begin{statement}{Definition}\label{defn:coordinate-matrices}
	Let $M$ and $N$ be modules over a commutative ring $R$ with ordered bases 
	$\basis{A} = (\alpha_1, \ldots, \alpha_m)$ and
	$\basis{S} = (\sigma_1, \ldots, \sigma_n)$ respectively.
	Let $\phi: N \to M$ be a module homomorphism.
	\begin{enumerate}
		\item 
		The \textbf{coordinate vector 
			$[\alpha]_\basis{A} \in \M_{m,1}(R)$ 
		of $\alpha = \sum_{j=1}^m r_j \alpha_j \in M$ 
		relative to $\basis{A}$} is the column vector given by 
		$[\alpha]_\basis{A}(j) = r_j$ for $j \in \set{1, \ldots, m}$, i.e. 
		$
			[\alpha]_\basis{A} = (
				r_1, r_2, \cdots, r_m
			)^\top
		$.

		\item 
		The 
		\textbf{matrix $[\phi]_{\basis{A}, \basis{S}} \in \M_{m,n}(R)$
		of $\phi$ relative to $\basis{A}$ and $\basis{S}$}
		is the matrix given by 
		\begin{equation*}
			\col_i[\phi]_{\basis{A}, \basis{S}} = [\phi(\sigma_i)]_\basis{A}
			\quad\text{ for all }\quad
			i \in \set{1, \ldots, n}.
		\end{equation*}
	\end{enumerate}
	If $\basis{A}$ and $\basis{S}$ are the standard bases for $M$ and $N$ respectively (if such are defined), then we may suppress the subscripts and write $[\alpha]$ and $[\phi]$ for $[\alpha]_\basis{A}$
	and $[\phi]_{\basis{S}, \basis{A}}$ respectively.
\end{statement}

\clearpage

\section{A Brief Review of Categories and Functors}
\label{appendix:cat-theory}

Below, we list a number of basic definitions and results relevant in this expository paper, mostly taken from the text \textit{Category Theory in Context} \cite[Chapter 1]{cattheory:rhiel} by Emily Riehl.

Note that there are some topics in category theory that are relevant to what is discussed in this paper (albeit some tangentially) that are not included in this appendix, e.g.\ abelian categories and homological algebra.
For said topics, we recommend \textit{Introduction to Homological Algebra} \cite[Chapters 1 and 5]{cattheory:rotman} by Joseph Rotman
and \textit{Introduction to Homological Algebra} \cite[Chapter 1 and Appendix A]{hom-algebra:weibel} for introductory reading.

\begin{definition}\cite[Definition 1.1.1]{cattheory:rhiel}
	\label{defn:category}
	A \textbf{category} $\catname{C}$ consists of a class $\objectset{\catname{C}}$ of \textbf{objects} and a class $\homset{\catname{C}}$ of \textbf{morphisms} such that all the following conditions are satisfied:
	\begin{enumerate}
		\item Each morphism has specified \textbf{domain} and \textbf{codomain} objects. 
		This relationship is typically denoted by $f: X \to Y$, where $f \in \homset{\catname{C}}$ with $X$ and $Y$ as its domain and codomain objects respectively.

		\item 
		Each object $X \in \objectset{\catname{C}}$ has a designated \textbf{identity morphism} $1_X: X \to X$.

		\item For any pairs for morphisms $f,g \in \homset{\catname{C}}$ such that $f: X \to Y$ and $g: Y \to Z$ (i.e.\ $f$ and $g$ are \textbf{composable}), there exists a specified \textbf{composite morphism} $g \circ f \in \homset{\catname{C}}$ such that $gf: X \to Z$. 
		The collection of these assignments is usually referred as the \textbf{composition law} of $\catname{C}$.

		\item 
		Composition is unital with identity morphisms, i.e.\ 
		for any morphism $f \in \homset{\catname{C}}$ with $f: X \to Y$, $f \circ 1_X = f$ and $f = 1_Y \circ f$.

		\item 
		Composition is associative, i.e.\ 
		for any $f,g,h \in \homset{\catname{C}}$ with $f: X \to Y$, $g: Y \to Z$, and $h: Y \to Z$, the composites $h \circ (g \circ f)$ and $(h \circ g) \circ f$ must be equal. 
	\end{enumerate}
	A \textbf{subcategory} $\catname{D}$ of some category $\catname{C}$ is another category whose classes of objects and of morphisms are subclasses of those of $\catname{C}$. 
\end{definition}
\remark{ 
	For our purposes, it suffices to know that classes act similarly to sets, in that they are collections of objects. They are distinct from sets to avoid paradoxes such as Russell's paradox.
}

The class of morphisms and the composition law of a category determines the notion of similarity, i.e.\ isomorphism relations, in said category.
Below, we include a category-level definition of isomorphisms.

\begin{definition}
	A morphism $f: X \to Y$ in a category $\catname{C}$ is called an \textbf{isomorphism} if there exists another morphism $g: X \to Y$ such that $g \circ f = \id_X$ and $f \circ g = \id_Y$. 
	The objects $X$ and $Y$ are \textbf{isomorphic} in $\catname{C}$ if there exists a isomorphism $f: X \to Y$. 
	In this case, we write $X \cong_\catname{C} Y$.
	If the category $\catname{C}$ is clear from context, we usually write $X \cong Y$ for brevity.
\end{definition}

Note that, in general, equality between objects and morphisms in a category $\catname{C}$ are considered in the abstract sense.
This is because not all categories are concrete categories.
Loosely speaking, a category $\catname{C}$ is a \textit{concrete category} if each object $X \in \objectset{\catname{C}}$ are sets (usually with some additional structure).
If a category is not concrete, then notions such as \textit{injective maps} and \textit{surjective maps}, as defined in set theory, are ambiguous.
For morphisms that act similarly, new category-level terms are introduced. For example, injective maps are \textit{monomorphisms} and surjective maps are \textit{epimorphisms}. 

With the exception of $\poset(I,\leq)$ (see \fref{defn:poset-cat} below), 
all relevant categories discussed in this expository paper are concrete categories. 
For clarity, we prefer the set-level terminology when applicable. 
Listed below are some of these categories, taken from \cite[Example 1.1.3]{cattheory:rhiel} and \cite{persmod:bubenik-categorification}.
Note that the composition law on the categories below are given by the usual function composition.

\begin{enumerate}
	\item $\catset$ 
		denotes the category of sets (as objects) and functions (as morphisms). Isomorphisms in this category are \textit{bijections}, as is usually defined in set theory.

	\item $\cattop$ 
		denotes the category of topological spaces (as objects) and continuous maps (as morphisms). Isomorphisms in this category are homeomorphisms.
		Recall that composition of continuous maps are continuous.

	\item
		Let $\field$ be a field.
		$\catvectspace$ 
		denotes the category of $\field$-vector spaces (as objects) and $\field$-linear maps (as morphisms).
		Isomorphisms in $\catvectspace$ are usually called isomorphisms, although we may call these vector space isomorphisms for clarity,
	$\catvectspace$ has a subcategory denoted $\catfinvectspace$ composed of finite-dimensional $\field$-vector spaces, with the class of morphisms appropriately restricted.
		
	\item 
		Let $R$ be a ring.
		$\catmod{R}$
		(or \!$\prescript{}{R}{\catname{Mod}}$)
		refers to the category of right (or left) $R$-modules as objects and $R$-module homomorphisms as morphisms.
		Note that if $R$ is a commutative ring, then $\catmod{R} = \prescript{}{R}{\catname{Mod}}$.

		Since we generally require $R$ to be a PID in this paper, $R$ is commutative and we refer to $\catmod{R}$ as simply the category of $R$-modules.
		We also prefer using the symbol $\catmod{R}$ for formatting purposes.
		As with $\catvectspace$, we may refer to isomorphisms in $\catmod{R}$ $R$-module isomorphisms for clarity.

	\item 
		Let $R$ be a PID.
		$\catchaincomplex{\catmod{R}}$ denotes the category of chain complexes of $R$-modules as objects and chain maps as morphisms.
		In this paper, we use the following definition for chain complexes:
		A chain complex $C_\ast$ of $R$-modules is a $\ints$-indexed collection of $R$-modules $C_n$ and $R$-module homomorphisms $\boundary_n: C_n \to C_{n-1}$. In this case, we write $C_\ast = (C_n, \boundary_n)_{n \in \ints}$ and call $\boundary_n$ the \textit{differentials} of $C_\ast$. 
\end{enumerate}
In this paper, for brevity, we usually call categories only by their objects, e.g.\ we say $\catmod{R}$ is the category of $R$-modules, without reference to $R$-module homomorphisms. 

\spacer 

Poset categories play a huge role in this expository paper. We provide a definition taken from \cite[Example 1.1.4]{cattheory:rhiel}, but with notation slightly changed for clarity.

\begin{definition}\label{defn:poset-cat}
	The \textbf{poset category} $\poset(I, \leq)$, also called \textbf{indexing category}, induced by a poset $(I, \leq)$
	is the category constructed as follows:
	\begin{enumerate}
		\item The elements of $I$ are exactly the objects of $\poset(I, \leq)$.
		\item For all $a,b \in I$, 
			there exists a unique morphism $a \to b$ if and only if $a \leq b$.
		\item Given all $a \in I$, the morphism $a \to a$ is the identity morphism of the object $a$.
		\item For all $a,b,c \in I$ with $a \leq b \leq c$, the composition law on the morphism is given by $(c \leftarrow a) = (c \leftarrow b) \circ (b \leftarrow c)$.
	\end{enumerate}
\end{definition}

	Posets can usually be described and defined using simple directed graphs.
	In particular, a poset $(I,\leq)$ corresponds to a simple directed graph where the vertices represent the elements of $I$ and the arrows represent the relations in $\leq$.
	For example, the poset $(\nonnegints, \leq)$ has the corresponding directed graph:
	\begin{center}
		\includegraphics[width=0.7\linewidth]{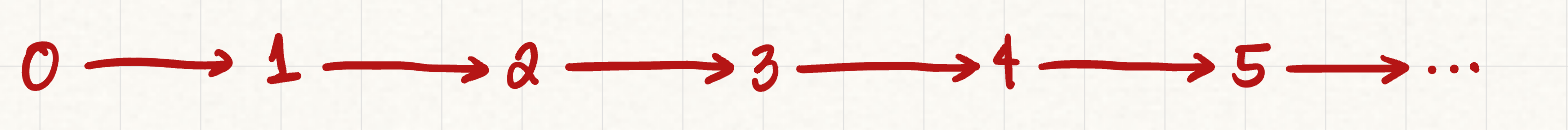}
	\end{center}
	Note that $\nonnegints$ is a totally ordered set and therefore, is also a poset.
	This interpretation also extends to the induced category $\poset(I,\leq)$, wherein the vertices and arrows of the graph represent the objects and morphisms of $\poset(I,\leq)$ respectively.

	Since posets have to be antisymmetric, the only isomorphisms on $\poset(I, \leq)$ are the identity morphisms. 
	Given any $a,b \in I$ with $a \leq b$, there cannot be a relation $b \leq a$ unless $a = b$. In that case, there does not exist a morphism $b \to a$.
	In other words, the notions of equality and isomorphisms are the same in $\poset(I, \leq)$.

\spacer 

Functors are also present in a significant amount in this expository paper.
Roughly speaking, functors are functions between categories, but instead of just having assignments on the objects, there are also assignments on the morphisms.
We provide a definition of functors below taken from~\cite[Definition 1.3.1]{cattheory:rhiel}

\begin{definition}\label{defn:functor}
	A \textbf{functor} $F: \catname{C} \to \catname{D}$ between categories $\catname{C}$ and $\catname{D}$ consists of the following:
	\begin{enumerate}
		\item 
			An (object) assignment of each object $X$ in $\catname{C}$ to some object $A$ in $\catname{D}$, denoted $F(X) = A$.
		\item 
			A (morphism) assignment of each morphism $f: X \to Y$ in $\catname{C}$ to a morphism $h: F(X) \to F(Y)$ in $\catname{D}$, denoted $F(f) = h$. 
			Note that the domain and codomain of the morphism assignment are determined by the object assignment.
		\item 
			The functor must respect composition, i.e.\ 
				for any composable pair $f$ and $g$ in $\catname{C}$, 
				$F(f)$ and $F(g)$ must also be composable in $\catname{D}$, i.e.\ $F( f \circ g ) = F(f) \circ F(g)$.
		\item 
			The functor must respect identity maps, i.e.\ 
				the functor assigns the identity morphism $\id_X$ of any object $X$ in $\catname{C}$ to the identity morphism $\id_{F(X)}$ of $F(X)$ in $\catname{D}$, i.e.\ $F(\id_X) = \id_{F(X)}$.
	\end{enumerate}
	We call $\catname{C}$ the \textbf{domain category} and $\catname{D}$ as the \textbf{codomain category}.
	The last two conditions listed above are often called \textbf{functorial properties}.
\end{definition}
\remark{
	More generally, the statement above defines a \textit{covariant functor}. 
	There is another type of functor called a \textit{contravariant functor} but it is not relevant to this paper.
}

A number of constructions in algebraic topology can be considered functors, which is appropriate given that category theory is said to be motivated by observations on algebraic topology theory. 

\begin{example}\label{ex:functor-examples}
	Listed below are some functors related to homology and simplicial homology, taken from \cite[Example 1.3.2]{cattheory:rhiel} and \cite[Section 2.3]{algtopo:hatcher}.
	\begin{enumerate}
		\item 
		The \textbf{$n$\th singular homology functor} $H_n(-; R): \cattop \to \catmod{R}$ sends a topological space $X$ to its $n$\th homology group $H_n(X; R)$ with coefficients in a PID $R$ and $n \in \nonnegints$. When $R = \ints$, the $n$\th homology functor is usually written as $H_n(-)$, suppressing the reference to coefficient ring $\ints$.
		The morphism assignment sends an inclusion map $i: X \to Y$ to an induced homomorphism $i_\ast: H_n(X; R) \to H_n(Y; R)$ as denoted in \cite{algtopo:hatcher}.
		
		\item 
		Abstractly, homology can be calculated from any chain complex.
		For each $n \in \ints$, 
		we have the \textbf{$n$\th chain homology functor} $H_n(-): \catchaincomplex{\catmod{R}} \to \catmod{R}$ that sends a chain complex $C_\ast = (C_n, \boundary_n)_{n \in \ints}$ of $R$-modules to its $n$\th homology group by the following construction:
		\begin{equation*}
			H_n(C_\ast) = \ker(\boundary_n) \,/\, \image(\boundary_{n+1})
		\end{equation*}
		The morphism assignment sends the chain map $\phi_\ast: C_\ast \to A_\ast$ with $\phi_\ast = (\phi_n: C_n \to A_n)_{n \in \ints}$ 
		to the homomorphism $H_n(C_\ast) \to H_n(A_\ast)$
		induced by the cokernel operation on $\phi_n: C_n \to A_n$ and $\phi_{n+1}: C_{n+1} \to A_{n+1}$.

		\item 
		The ``free'' functor $F: \catset \to \catmod{R}$ 
		sends a set $X$ to the free $R$-module $R\ket{X}$, with $X$ treated as a set of indeterminates (see \fref{defn:formal-sums})
		and sends a function $f: X \to Y$ 
		to the unique homomorphism $R\ket{X} \to R\ket{Y}$ defined by mapping basis elements to basis elements.  
	\end{enumerate}
\end{example}

As mentioned earlier, poset categories $\poset(I, \leq)$ can be described using graphs.
In the same vein, we can illustrate commutative diagrams as functors with some poset category as the domain category.
We provide a definition below taken from~\cite[Definition 1.6.4]{cattheory:rhiel}.

\begin{definition}\label{defn:cat-diagrams}
	A \textbf{diagram} is a functor $F: \poset(I, \leq) \to \catname{C}$ from some poset category $\poset(I, \leq)$ to some other category $\catname{C}$.
\end{definition}

Given a diagram $F: \poset(I,\leq) \to \catname{C}$,
	the vertices of the directed graph corresponding to $\poset(I,\leq)$ are replaced by objects and the arrows by morphisms.
The transitivity of $(I,\leq)$ and the composition axiom for functors require the resulting diagram (in the non-category theory-sense) be commutative.
We provide a simple example involving squares:

\begin{example}
	Let the poset $(S, \leq)$ be given by the following directed graph:
	\begin{center}
		\includegraphics[height=0.7in]{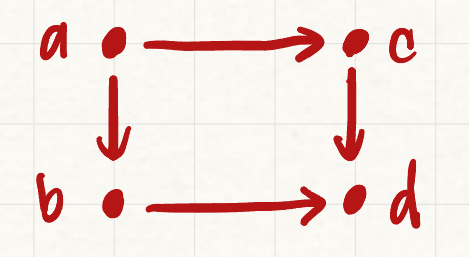}
	\end{center}
	Let $V_1, V_2, W_1, W_2$ be $\reals$-vector spaces.
	Let $\alpha_1: V_1 \to V_2, \alpha_2: W_1 \to W_2$ and 
		$\phi_1: V_1 \to W_1, \phi_2: V_2 \to W_2$ be linear maps.
	Then, the requirement that $\phi_2 \circ \alpha_1 = \alpha_2 \circ \phi_1$ can be restated as follows:
	Let $F: \poset(S, \leq) \to \catvectspace$ be an assignment on the objects and morphisms of $\poset(S, \leq)$ as follows:
	\begin{align*}
		F(a) &= V_1
		&	F(c) &= V_2
		&	F(b) &= W_1
		&	F(d) &= W_2
		\\
		F(a \to b) 		&= \phi_1
		&	F(a \to c) 	&= \alpha_1
		& 	F(b \to d) 	&= \alpha_2
		&	F(c \to d)	&=	\phi_2
	\end{align*}
	We illustrate the assignment $F$ as a diagram below:
	\begin{center}
		\includegraphics[height=1in]{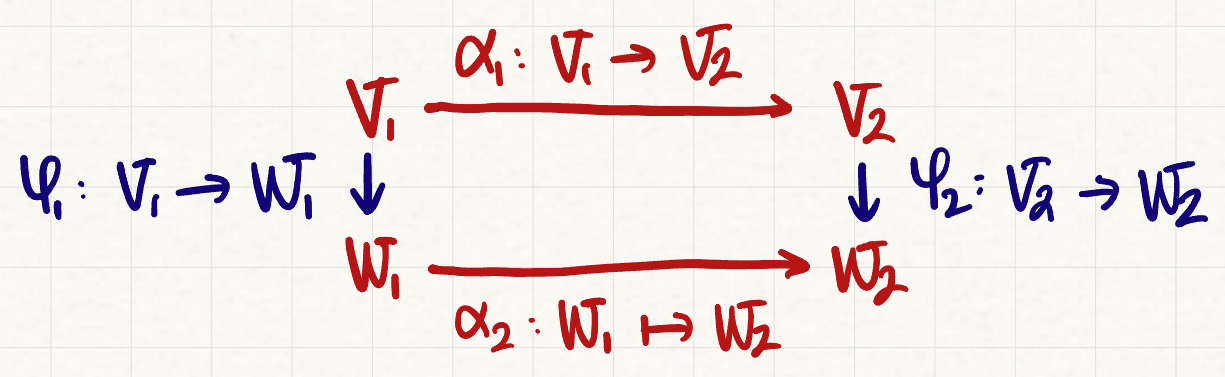}
	\end{center}
	Then, $\phi_2 \circ \alpha_1 = \alpha_2 \circ \phi_1$ if and only if $F$ represents a functor. 
	By definition of functor, the morphisms in the form $x \to x$ are mapped to identity linear maps.
	In terms of compositions, we have the following chain of equalities:
	\begin{align*}
		\phi_2 \circ \alpha_1
		= F(d \leftarrow c) \circ F(c \leftarrow a)
		\stackrel{\star}{=} F(d \leftarrow a)
		\stackrel{\star}{=} F(d \leftarrow b) \circ F(b \leftarrow a)
		= \alpha_2 \circ \phi_1
	\end{align*}
	The equalities labeled by $\star$ are due to the uniqueness of morphisms in $\poset(S, \leq)$.
\end{example}



\setlength\bibitemsep{\parskip}
\printbibliography[
	heading=bibintoc,
	title={References},
] 

\end{document}